\documentclass[12pt,english,refpage,intoc,bibliography=totoc,index=totoc,BCOR7.5mm,captions=tableheading]{amsart}
\usepackage{lmodern}

\usepackage[T1]{fontenc}
\usepackage[latin9]{inputenc}
\usepackage{geometry}
\usepackage{color}
\usepackage{babel}
\usepackage{textcomp}
\usepackage{mathrsfs}
\usepackage{enumitem}
\usepackage{amstext}
\usepackage{amsthm}
\usepackage{amssymb}
\usepackage{cancel}
\PassOptionsToPackage{normalem}{ulem}
\usepackage{ulem}
\usepackage[unicode=true,
 bookmarks=true,bookmarksnumbered=true,bookmarksopen=false,
 breaklinks=false,pdfborder={0 0 1},backref=false,colorlinks=true]
 {hyperref}
\hypersetup{pdftitle={Algebras of $p$-Adic Distributions Induced by Pointwise Products
of F-Series},
 pdfauthor={LyX Team},
 pdfsubject={LyX},
 pdfkeywords={LyX},
 linkcolor=black, citecolor=black, urlcolor=blue, filecolor=blue, pdfpagelayout=OneColumn, pdfnewwindow=true, pdfstartview=XYZ, plainpages=false}

\makeatletter
\numberwithin{equation}{section}
\numberwithin{figure}{section}
\theoremstyle{remark}
    \ifx\thechapter\undefined
      \newtheorem{note}{\protect\notename}
    \else
      \newtheorem{note}{\protect\notename}[chapter]
    \fi
\theoremstyle{plain}
    \ifx\thechapter\undefined
      \newtheorem{assumption}{\protect\assumptionname}
    \else
      \newtheorem{assumption}{\protect\assumptionname}[chapter]
    \fi
\theoremstyle{definition}
    \ifx\thechapter\undefined
      \newtheorem{defn}{\protect\definitionname}
    \else
      \newtheorem{defn}{\protect\definitionname}[chapter]
    \fi
\theoremstyle{plain}
    \ifx\thechapter\undefined
      \newtheorem{prop}{\protect\propositionname}
    \else
      \newtheorem{prop}{\protect\propositionname}[chapter]
    \fi
\theoremstyle{remark}
    \ifx\thechapter\undefined
      \newtheorem{rem}{\protect\remarkname}
    \else
      \newtheorem{rem}{\protect\remarkname}[chapter]
    \fi
\theoremstyle{plain}
    \ifx\thechapter\undefined
	    \newtheorem{thm}{\protect\theoremname}
	  \else
      \newtheorem{thm}{\protect\theoremname}[chapter]
    \fi
\theoremstyle{plain}
    \ifx\thechapter\undefined
  \newtheorem{cor}{\protect\corollaryname}
\else
      \newtheorem{cor}{\protect\corollaryname}[chapter]
    \fi
\theoremstyle{definition}
    \ifx\thechapter\undefined
      \newtheorem{example}{\protect\examplename}
    \else
      \newtheorem{example}{\protect\examplename}[chapter]
    \fi
\theoremstyle{remark}
    \ifx\thechapter\undefined
      \newtheorem{notation}{\protect\notationname}
    \else
      \newtheorem{notation}{\protect\notationname}[chapter]
    \fi
\theoremstyle{plain}
    \ifx\thechapter\undefined
      \newtheorem{question}{\protect\questionname}
    \else
      \newtheorem{question}{\protect\questionname}[chapter]
    \fi
\theoremstyle{remark}
    \ifx\thechapter\undefined
      \newtheorem{claim}{\protect\claimname}
    \else
      \newtheorem{claim}{\protect\claimname}[chapter]
    \fi
\theoremstyle{plain}
    \ifx\thechapter\undefined
      \newtheorem{lem}{\protect\lemmaname}
    \else
      \newtheorem{lem}{\protect\lemmaname}[chapter]
    \fi
\theoremstyle{plain}
    \ifx\thechapter\undefined
      \newtheorem{fact}{\protect\factname}
    \else
      \newtheorem{fact}{\protect\factname}[chapter]
    \fi
\theoremstyle{plain}
\newtheorem*{cor*}{\protect\corollaryname}

\usepackage[figure]{hypcap}

\usepackage{fancyhdr}

\@ifundefined{extrarowheight}
 {\usepackage{array}}{}
\makeatother

\providecommand{\assumptionname}{Assumption}
\providecommand{\claimname}{Claim}
\providecommand{\corollaryname}{Corollary}
\providecommand{\definitionname}{Definition}
\providecommand{\examplename}{Example}
\providecommand{\factname}{Fact}
\providecommand{\lemmaname}{Lemma}
\providecommand{\notationname}{Notation}
\providecommand{\notename}{Note}
\providecommand{\propositionname}{Proposition}
\providecommand{\questionname}{Question}
\providecommand{\remarkname}{Remark}
\providecommand{\theoremname}{Theorem}

\begin{document}
\title{Algebras of $p$-Adic Distributions Induced by Pointwise Products
of F-Series}
\author{Maxwell C. Siegel}
\date{17 June 2025}
\begin{abstract}
Let $p$ be an integer $\geq2$ and let $K$ be a global field. A
foliated $p$-adic F-series is a function $X$ of a $p$-adic integer
variable $\mathfrak{z}$ satisfying the functional equations $X\left(p\mathfrak{z}+j\right)=a_{j}X\left(\mathfrak{z}\right)+b_{j}$
for all $\mathfrak{z}\in\mathbb{Z}_{p}$ and all $j\in\left\{ 0,\ldots,p-1\right\} $,
where the $a_{j}$s and $b_{j}$s are indeterminates. Treating $X$
as taking value in a certain ring of formal power series over $K$,
this paper establishes a universal/functorial Fourier theory for F-series:
we show that $X$ has a Fourier transform, and that, for nearly any
ideal $I\subseteq R$, where of $R=\mathcal{O}_{K}\left[a_{0},\ldots,a_{p-1},b_{0},\ldots,b_{p-1}\right]$,
this Fourier transform descends through the quotient mod $I$ which
imposes on $X$ the relations encoded by $I$. Furthermore, we show
that the pointwise product of $X$ with itself $n$ times also has
a Fourier transform compatible with descent. These results generalize
products $X_{1}^{e_{1}}\cdots X_{d}^{e_{d}}$ of any $d$ distinct
F-series $X_{1},\ldots,X_{d}$ with integer exponents $e_{1},\ldots,e_{d}\geq0$.
Using these Fourier transforms, F-series and their products can be
identified with distributions on $\mathbb{Z}_{p}$ in a manner compatible
with descent mod $I$, forming algebras under pointwise multiplication,
in striking contrast to the general principle that the pointwise product
of two distributions will not be a distribution. Also, to any given
F-series or product thereof, one can associate an affine algebraic
variety over $K$ which I call the breakdown variety. The distributions
induced by a product of F-series under descent mod $I$ exhibit sensitivity
to $I$'s containment of the ideal corresponding to the distributions'
breakdown varieties. This presents a novel method of encoding given
affine algebraic varieties via the breakdown varieties of F-series
distributions in a manner compatible with pointwise products, convolutions,
and tensor products.
\end{abstract}

\address{1626 Thayer Ave., Los Angeles, CA, 90024}
\email{siegelmaxwellc@ucla.edu}

\maketitle
\tableofcontents{}

\section*{Preface}

This is a long, unorthodox paper, filled with lengthy, detailed, formal
computations. It is well-known that such computations often obscure
the underlying ideas, and make arguments more difficult to read. I
would like to think that this paper is the exception which proves
the rule. In this paper, \emph{the computations }\emph{\uline{are}}\emph{
the central object of study}. I first discovered them in the course
of doing my PhD Dissertation \cite{My Dissertation}. In writing this
paper, my overriding objective was to clothe these computations in
a rigorous theoretical framework that not only clarifies what is going
on, but allows us to begin interfacing my discoveries with extant
mathematics. I invite the reader to treat this paper more as the records
of a scientific expedition, rather than a more traditional research
article. In my studies, I have stumbled onto a strange new world;
my present research is a chronicle of my observations of this world
and its inhabitants, and of my attempts to describe the order that
appears to govern them.

 To begin with, even in the introduction (which properly begins on
page \ref{subsec:Introduction}), there is a great deal of non-standard
terminology. I'm not pulling out this terminology just for novelty's
sake, or to show off; rather, I use it because I \emph{must}. To my
knowledge, there is no established vocabulary for describing many
of the objects and phenomena that I have discovered, and though I
have subsequently discovered that at least one of the key concepts\textemdash that
of a frame\textemdash can (and will) be given an elegant presentation
in terms of locally convex topological rings, using that language
obscures the unusual nature of what is actually going on.

Regarding my use of novel terminology, the most important terms are
\textbf{F-series},\textbf{ frames}, \textbf{quasi-integrability},
\textbf{rising-continuity}, \textbf{reffinite},\textbf{ }and \textbf{breakdown
varieties}. Variants of these, such as ``digital frames'', or ``foliated
F-series'' are just specific instances of more general constructions,
as is common in mathematics. The reader can rest assured that \emph{all
}of these concepts will be defined in full at the appropriate time.
While these concepts can be studied in isolation, for the purposes
of this paper, I invoke them only because they are necessary to make
rigorous sense of what is going on. To that end, I implore the reader
to \emph{suspend their disbelief} for the time being. The notation
and conventions set up immediately below this message are ten-thousand
times more important than the meaning of quasi-integrability or the
nature of what I call a frame.

If the paper's introduction seems opaque or difficult to follow, I
recommend taking a look at the examples on pages \pageref{exa:choosing absolute values},
\pageref{exa:Let-,-where}. Being able to follow my computations there
and verify them for yourself is far, far more important to understanding
this paper than immediate knowledge of what a digital reffinite frame
could ever be. As long as you can do the computations, you can discover
for yourself the phenomena I have encountered, and thereby understand
why we need these new concepts to make sense of them.

Finally, I need to speak about the paper's intended audience, or rather,
its \emph{lack} of one. I would love to be able to say that I wrote
this for a particular group of specialists who could immediately gobble
it up, but I can't, and for the simple fact that, by mathematical
standards, the paper's contents are fundamentally interdisciplinary.
The main computational tools are elementary Fourier-theoretic computations
for locally constant/Schwartz-Bruhat functions of a $p$-adic integer
variable, but with the twist that instead of being real or complex-valued,
our functions will take values in a wide array of different rings.
Topologies will not always be needed, let alone metrics or absolute
values, and our analysis can and will change from working in ordinary
metric spaces to non-archimedean/ultrametric spaces at the drop of
a hat. Dealing with this variability requires drawing from functional
analysis\textemdash particularly, the ideas of locally convex topological
vector spaces/algebras/rings\textemdash both of the archimedean and
non-archimedean variants. Because we will be examining absolute values
and seminorms on rings of number-theoretic importance, we will also
need algebraic number theory and many of the basic concepts surrounding
Dedekind domains, such as global fields, places, rings of integers,
and completions thereof with respect to an absolute value. Finally,
fulfilling the abstract's promise about a universal form of Fourier
analysis and relating it to algebraic varieties demands that we draw
from basic constructions in commutative algebra and algebraic geometry,
such as localizations of rings (both at and away from a given prime
ideal), and the correspondence between algebraic sets and ideals of
polynomial rings.

I recognize that this is a lot. Thankfully, there's a saving grace:
although a great amount of material is involved, we won't need to
get particularly deep in any of the subject areas we will draw from.
You won't need more than an advanced undergraduate- or early graduate-level
understanding of any of the subjects involved. The difficulty is horizontal,
not vertical. I happen to find a great deal of satisfaction in that;
it is beautiful to see all these different ideas working in unison,
and a pleasant surprise that their coöperation results in unexpected
structure, rather than an incoherent mess.

\section{Entr'acte\label{sec:Introduction}}

\subsection{\label{subsec:Preliminaries-(Notation,-etc.)}Preliminaries (Notation,
etc.)}
\begin{note}[\textbf{¡WARNING!}]
Though \textbf{Section \ref{subsec:Quasi-Integrability-=000026-Degenerate}}
begins with a brief overview of $\left(p,q\right)$-adic Fourier analysis,
for considerations of length, this paper DOES NOT include a detailed
exposition of the relevant material. My paper \cite{2nd blog paper}
gives all of the necessary details; the content can be found for free
\href{https://arxiv.org/abs/2208.11082}{on arXiv} at <https://arxiv.org/abs/2208.11082>.
That being said, the present paper \textbf{should be readable }provided
the reader knows the details of Pontryagin duality and integration
of real- or complex-valued functions of a $p$-adic integer variable,
such as they occur in any of \cite{Ramakrishnan,Tate's thesis,Automorphic Representations,Folland - harmonic analysis}.
Ironically enough, the computational formalism (Schwartz-Bruhat functions,
unitary characters, Fourier integrals, Haar measures) for the case
of functions $\mathbb{Z}_{p}\rightarrow\mathbb{Q}_{q}$ for distinct
primes $p$ and $q$ is \emph{identical }to the case of functions
$\mathbb{Z}_{p}\rightarrow\mathbb{C}$. It is only in certain finer
matters such as series convergence and Banach spaces of well-behaved
functions (ex: continuous; integrable; admitting absolutely convergent
Fourier series; etc.) that the meaningful differences between $\left(p,q\right)$-adic
and archimedean-valued Fourier analysis emerge.
\end{note}
We write $\mathbb{N}_{0}$ to denote the set of all non-negative integers;
$\mathbb{N}_{1}$ denotes the set of all integers $\geq1$; $\mathbb{N}_{\pi}$
denotes the set of all integers $\geq\pi$, and so on. We fix an integer
$p\geq2$ and a global field $K$; in the case where $K$ has positive
characteristic, we require $\textrm{char}K$ to be co-prime to $p$.
We write $K\left(\zeta_{p^{\infty}}\right)$ to denote the direct
limit $\lim_{n\rightarrow\infty}K\left(\zeta_{p^{n}}\right)$, the
maximal cyclotomic extension of $K$ obtained by adjoining to $K$
roots of unity of order $p^{n}$ for all $n\geq1$. More generally,
for any ring $R$, we write $R\left(\zeta_{p^{\infty}}\right)$ to
denote $\lim_{n\rightarrow\infty}R\left(\zeta_{p^{n}}\right)$.

When $p$ is composite, we define $\mathbb{Z}_{p}$, the ring of $p$-adic
integers, as:
\begin{equation}
\mathbb{Z}_{p}\overset{\textrm{def}}{=}\prod_{\ell\mid p}\mathbb{Z}_{\ell}^{v_{\ell}\left(p\right)}
\end{equation}
where the product is taken over all prime divisors $\ell$ of $p$,
where $v_{\ell}$ is the $\ell$-adic valution (valuation, \emph{not
}absolute value!), and where $\overset{\textrm{def}}{=}$ means ``by
definition''. Likewise, we write:
\begin{equation}
\mathbb{Q}_{p}\overset{\textrm{def}}{=}\prod_{\ell\mid p}\mathbb{Q}_{\ell}^{v_{\ell}\left(p\right)}
\end{equation}

We write $\mathbb{Z}_{p}^{\prime}$ to denote $\mathbb{Z}_{p}\backslash\mathbb{N}_{0}$.
We write $\hat{\mathbb{Z}}_{p}$ to denote $\mathbb{Z}\left[1/p\right]/\mathbb{Z}$,
the Pontryagin dual of $\mathbb{Z}_{p}$. We identify $\hat{\mathbb{Z}}_{p}$
with the group of rational numbers in $\left[0,1\right)$ whose denominators
are powers of $p$. For any $t\in\hat{\mathbb{Z}}_{p}\backslash\left\{ 0\right\} $,
observe that $t$ can be written in irreducible form as $k/p^{n}$
for some integers $n\geq1$ and some $k\in\left\{ 0,\ldots,p^{n}-1\right\} $.
Consequently, the $p$-adic absolute value of $t$ is $\left|t\right|_{p}=p^{n}$,
and the numerator of $t$ is given by $t\left|t\right|_{p}=k$, with
both values being $0$ when $t=0$.

Also, for aesthetic reasons, we write $p$-adic variables in lower-case
$\mathfrak{fraktur}$ font.

We use \textbf{Iverson bracket notation} extensively, writing $\left[S\right]$
for a statement $S$ to denote a quantity which is $1$ when $S$
is true and $0$ otherwise. For example, given integers $n\geq0$
and $k\in\left\{ 0,\ldots,p^{n}-1\right\} $:
\begin{equation}
\mathfrak{z}\in\mathbb{Z}_{p}\mapsto\left[\mathfrak{z}\overset{p^{n}}{\equiv}k\right]\in\left\{ 0,1\right\} 
\end{equation}
is the function which is equal to $1$ when the $p$-adic integer
$\mathfrak{z}$ is congruent to $k$ mod $p^{n}$ and which is equal
to $0$ otherwise, where $\overset{p^{n}}{\equiv}$ denotes congruence
mod $p^{n}$. We write $\left[\cdot\right]_{p^{n}}:\mathbb{Z}_{p}\rightarrow\left\{ 0,\ldots,p^{n}-1\right\} $
to denote the map which, for any $\mathfrak{z}\in\mathbb{Z}_{p}$,
outputs the unique integer $\left[\mathfrak{z}\right]_{p^{n}}\in\left\{ 0,\ldots,p^{n}-1\right\} $
congruent to $\mathfrak{z}$ mod $p^{n}$.

When a limit is being taken in a metric space, we write the metric
space in which the convergence occurs over the equals sign. Hence:
\begin{equation}
\lim_{n\rightarrow\infty}f\left(n\right)\overset{\mathbb{R}}{=}0
\end{equation}
means that $f$ converges to $0$ in the reals, whereas:
\begin{equation}
\lim_{n\rightarrow\infty}f\left(n\right)\overset{\mathbb{Q}_{5}}{=}0
\end{equation}
means that $f$ converges to $0$ in the $5$-adics.

Letting $\left(M,d\right)$ be a metric space, letting $L$ be a point
of $M$, and considering a function $\varphi:\hat{\mathbb{Z}}_{p}\rightarrow M$
(where $\hat{\mathbb{Z}}_{p}$ is the pontryagin dual of $\mathbb{Z}_{p}$)
we say that \textbf{$\varphi\left(t\right)$ converges (in $M$) to
$L$ as $\left|t\right|_{p}\rightarrow\infty$}, written:
\begin{equation}
\lim_{\left|t\right|_{p}\rightarrow\infty}\varphi\left(t\right)\overset{M}{=}L
\end{equation}
when, for all $\epsilon>0$, there exists an $N\geq1$ which makes
$d\left(\varphi\left(t\right),L\right)<\epsilon$ for all $t\in\hat{\mathbb{Z}}_{p}$
with $\left|t\right|_{p}\geq p^{N}$. Thus, for example, given $\hat{\chi}:\hat{\mathbb{Z}}_{p}\rightarrow\mathbb{Q}_{q}$,
for a prime number $q$, we say $\lim_{\left|t\right|_{p}\rightarrow\infty}\left|\hat{\chi}\left(t\right)\right|_{q}\overset{\mathbb{R}}{=}0$
to mean that, for all $\epsilon>0$, there exists an $N\geq1$ so
that the real number $\left|\hat{\chi}\left(t\right)\right|_{q}$
is less than $\epsilon$ for all $t\in\hat{\mathbb{Z}}_{p}$ with
$\left|t\right|_{p}\geq p^{N}$.

We write $\left\{ \cdot\right\} _{p}$ to denote the \textbf{$p$-adic
fractional part}, with:
\begin{equation}
\left\{ \mathfrak{y}\right\} _{p}\overset{\textrm{def}}{=}\sum_{\ell\mid p}\left\{ \mathfrak{y}\right\} _{\ell}
\end{equation}
where $\left\{ \cdot\right\} _{\ell}:\mathbb{Q}_{\ell}\rightarrow\hat{\mathbb{Z}}_{\ell}$
is the $\ell$-adic fractional part defined by:
\begin{equation}
\left\{ \mathfrak{x}\right\} _{\ell}=\begin{cases}
0 & \textrm{if }n_{0}\geq0\\
\sum_{n=n_{0}}^{-1}c_{n}\ell^{n} & \textrm{if }n_{0}\leq-1
\end{cases}
\end{equation}
where:
\begin{equation}
\mathfrak{x}=\sum_{n=n_{0}}^{\infty}c_{n}\ell^{n}
\end{equation}
is the Hensel series representation of $\mathfrak{x}\in\mathbb{Q}_{\ell}$.
Thus, for example, given $t\in\hat{\mathbb{Z}}_{p}$ and $\mathfrak{z}\in\mathbb{Z}_{p}$:
\begin{equation}
e^{2\pi i\left\{ t\mathfrak{z}\right\} _{p}}=\begin{cases}
1 & \textrm{if }t=0\\
e^{2\pi it\left[\mathfrak{z}\right]_{\left|t\right|_{p}}} & \textrm{else}
\end{cases}
\end{equation}
where $\left[\mathfrak{z}\right]_{\left|t\right|_{p}}$ is the projection
of $\mathfrak{z}$ mod the power of $p$ in $t$'s denominator. Thus,
if $p$ is odd:
\begin{equation}
e^{2\pi i\left\{ \frac{2\mathfrak{z}}{p^{3}}\right\} _{p}}=e^{2\pi i\left(2\left[\mathfrak{z}\right]_{p^{3}}/p^{3}\right)}
\end{equation}

Given an abelian topological group $G$ and a function $\hat{\chi}:\hat{\mathbb{Z}}_{p}\rightarrow G$,
we write $\sum_{t\in\hat{\mathbb{Z}}_{p}}\hat{\chi}\left(t\right)$
to denote the limit:
\begin{equation}
\sum_{t\in\hat{\mathbb{Z}}_{p}}\hat{\chi}\left(t\right)\overset{\textrm{def}}{=}\lim_{N\rightarrow\infty}\sum_{\left|t\right|_{p}\leq p^{N}}\hat{\chi}\left(t\right)\overset{\textrm{def}}{=}\lim_{N\rightarrow\infty}\sum_{k=0}^{p^{N}-1}\hat{\chi}\left(\frac{k}{p^{N}}\right)
\end{equation}
where the limits are taken in the topology of $G$, provided that
the limit exists.

Given $s\in\hat{\mathbb{Z}}_{p}$, we write $\mathbf{1}_{s}:\hat{\mathbb{Z}}_{p}\rightarrow\left\{ 0,1\right\} $
to denote the indicator function of $\left\{ s\right\} $:
\begin{equation}
\mathbf{1}_{s}\left(t\right)\overset{\textrm{def}}{=}\left[t\overset{1}{\equiv}s\right]
\end{equation}
Unless stated otherwise, sums of the form $\sum_{n=n_{1}}^{n_{2}}f\left(n\right)$
are defined to be $0$ when $n_{2}<n_{1}$; likewise, products of
the form $\prod_{n=n_{1}}^{n_{2}}f\left(n\right)$ are defined to
be $1$ when $n_{2}<n_{1}$.

\begin{note}[Convention on Complex Exponentials \& Roots of Unity]
\label{note:In-our-approach}In our approach to Pontryagin duality,
a unitary character $G\rightarrow K$ of a locally compact abelian
group $G$ is an expression of the form:
\begin{equation}
g\in G\mapsto e^{2\pi i\left\langle \gamma,g\right\rangle }\in K
\end{equation}
where $\gamma$ is an element of the Pontryagin dual $\hat{G}$ of
$G$ and $\left\langle \cdot,\cdot\right\rangle :\hat{G}\times G\rightarrow\mathbb{R}/\mathbb{Z}$
is the duality bracket/duality pairing (a continuous, $\mathbb{Z}$-bilinear
homomorphism of locally compact abelian groups). To that end, \textbf{we
fix once and for all a choice of an embedding of the field of all
roots of unity in $\overline{K}$}, so that for any integer $d\geq1$,
we can write $e^{2\pi i/d}$ to denote our favorite primitive $d$th
root of unity in $\overline{K}$, and then exploit the familiar algebraic
properties of complex exponential notation, such as:
\begin{equation}
e^{2\pi ik/d}=\left(e^{2\pi i/d}\right)^{k}
\end{equation}
and so that, for any prime $p$ and any $n\geq0$:
\begin{equation}
\left(e^{2\pi i/p^{n+1}}\right)^{p}=e^{2\pi i/p^{n}}
\end{equation}
In this notation, every $\overline{K}$-valued unitary character on
$\mathbb{Z}_{p}$ is expressible as a map of the form:
\begin{align}
\mathbb{Z}_{p} & \rightarrow\overline{K}\\
\mathfrak{z} & \mapsto e^{2\pi i\left\{ t\mathfrak{z}\right\} _{p}}
\end{align}
for some fixed $t\in\hat{\mathbb{Z}}_{p}$; note that multiplication
of $t$ and $\mathfrak{z}$ makes sense here, because both quantities
lie in $\mathbb{Q}_{p}$.

In terms of our abuse-of-notation, in having identified $\hat{\mathbb{Z}}_{p}$
with $\mathbb{Z}\left[1/p\right]/\mathbb{Z}$, for any given $\mathfrak{y}\in\mathbb{\mathbb{Q}}_{p}$,
the fractional part of $\mathfrak{y}$ will be uniquely expressible
as an irreducible $p$-power fraction in $\left[0,1\right)$ of the
form $k/p^{n}$ for some integers $n\geq0$ and $k\in\left\{ 0,\ldots,p^{n}-1\right\} $.
As such, we write $e^{2\pi i\left\{ \mathfrak{y}\right\} _{p}}$ to
denote the $k$th power of the primitive $p^{n}$th root of unity
chosen by our embedding. Likewise, $e^{-2\pi i\left\{ \mathfrak{y}\right\} _{p}}$
is the reciprocal of the root of unity denoted by $e^{2\pi i\left\{ \mathfrak{y}\right\} _{p}}$.
\end{note}
\begin{assumption}
ALL RINGS ARE ASSUMED TO BE COMMUTATIVE AND UNITAL UNLESS STATED OTHERWISE.
\end{assumption}
\begin{defn}
For any $\mathfrak{z}\in\mathbb{Z}_{p}$ and any integer $n\geq0$,
we write $\varepsilon_{n}\left(\mathfrak{z}\right)$ to denote:

\begin{equation}
\varepsilon_{n}\left(\mathfrak{z}\right)\overset{\textrm{def}}{=}e^{\frac{2\pi i}{p^{n+1}}\left(\left[\mathfrak{z}\right]_{p^{n+1}}-\left[\mathfrak{z}\right]_{p^{n}}\right)}
\end{equation}
the output being a root of unity in whatever ring (or extension thereof)
we happen to be working in, as per \textbf{Note \ref{note:In-our-approach}}.
\end{defn}
\begin{defn}
Given an integer $n\geq0$, and given $k\in\left\{ 1,\ldots,p-1\right\} $,
we write $\#_{p:k}\left(n\right)$ to denote the number of $j$s among
the base $p$ digits of $n$. We have that:
\begin{equation}
\#_{p:k}\left(pn+j\right)=\#_{p:k}\left(n\right)+\left[j=k\right],\textrm{ }\forall n\geq0,\textrm{ }\forall j,k\in\left\{ 1,\ldots,p-1\right\} 
\end{equation}
where $\left[k=j\right]$ is an Iverson bracket. Letting $\lambda_{p}\left(n\right)$
denote the total number of digits in $n$'s base $p$-adic representation,
we then define:
\begin{equation}
\#_{p:0}\left(n\right)\overset{\textrm{def}}{=}\lambda_{p}\left(n\right)-\sum_{k=1}^{p-1}\#_{p:k}\left(n\right)
\end{equation}
so that $\#_{p:0}\left(0\right)=0$. Note that $\lambda_{p}$ can
be written in closed form as:
\begin{equation}
\lambda_{p}\left(n\right)\overset{\textrm{def}}{=}\left\lceil \log_{p}\left(n+1\right)\right\rceil =\left\lceil \frac{\ln\left(n+1\right)}{\ln p}\right\rceil ,\textrm{ }\forall n\geq0
\end{equation}

\textbf{Note}:\textbf{ }When there is no chance of confusion, we write
$\#_{j}$ instead of $\#_{p:j}$, to save space.
\end{defn}
\begin{defn}
The \textbf{$p$-adic shift map }$\theta_{p}:\mathbb{Z}_{p}\rightarrow\mathbb{Z}_{p}$
is defined by:
\begin{equation}
\theta_{p}\left(\mathfrak{z}\right)\overset{\textrm{def}}{=}\frac{\mathfrak{z}-\left[\mathfrak{z}\right]_{p}}{p},\textrm{ }\forall\mathfrak{z}\in\mathbb{Z}_{p}
\end{equation}
If we think of $p$-adic integers as base $p$ numbers with infinitely
many digits (viz. their Hensel series representation):
\begin{equation}
\mathfrak{z}=\sum_{n=0}^{\infty}d_{n}p^{n}
\end{equation}
then $\theta_{p}$ has the effect of deleting $\mathfrak{z}$'s smallest
$p$-adic digit (the value of $\mathfrak{z}$ mod $p$) and shifting
the remaining digit sequence over by one:
\begin{equation}
\theta_{p}\left(\mathfrak{z}\right)=\sum_{n=0}^{\infty}d_{n+1}p^{n}
\end{equation}
We write $\theta_{p}^{\circ n}$ to denote the composition of $n$
copies of $\theta_{p}$, with $\theta_{p}^{\circ0}$ being defined
as the identity map on $\mathbb{Z}_{p}$.
\end{defn}
The arithmetic functions defined above satisfy the following functional
equations. These will be used throughout this paper.
\begin{prop}
\label{prop:fundamental functional equations-1}

I.
\begin{equation}
\varepsilon_{n}\left(p\mathfrak{z}+k\right)=\begin{cases}
\varepsilon_{0}\left(k\right) & \textrm{if }n=0\\
\varepsilon_{n-1}\left(\mathfrak{z}\right) & \textrm{if }n\geq1
\end{cases},\textrm{ }\forall\mathfrak{z}\in\mathbb{Z}_{p},\textrm{ }\forall k\in\left\{ 0,\ldots,p-1\right\} ,\textrm{ }\forall n\geq0\label{eq:epsilon identity-1}
\end{equation}

II.
\begin{equation}
\#_{p:k}\left(\left[p\mathfrak{z}+j\right]_{p^{n}}\right)=\begin{cases}
0 & \textrm{if }n=0\\
\#_{p:k}\left(\left[\mathfrak{z}\right]_{p^{n-1}}\right)+\left[j=k\right] & \textrm{if }n\geq1
\end{cases},\textrm{ }\forall\mathfrak{z}\in\mathbb{Z}_{p},\textrm{ }\forall j,k\in\left\{ 0,\ldots,p-1\right\} ,\textrm{ }\forall n\geq0\label{eq:number of js identity-1}
\end{equation}
More generally: 
\begin{equation}
\#_{p:k}\left(\left[\theta_{p}^{\circ m}\left(\mathfrak{z}\right)\right]_{p^{n}}\right)=\#_{p:k}\left(\left[\mathfrak{z}\right]_{p^{n+m}}\right)-\#_{p:k}\left(\left[\mathfrak{z}\right]_{p^{m}}\right),\textrm{ }\forall m,n\geq0\label{eq:digit shift}
\end{equation}
\end{prop}
We need the notion of \textbf{(semi)norms }and \textbf{absolute values}
on rings.
\begin{defn}
\label{def:A-seminorm-on}A \textbf{seminorm }on $R$ is a function
$\left\Vert \cdot\right\Vert :R\rightarrow\left[0,\infty\right)$
so that:
\begin{align}
\left\Vert 0_{R}\right\Vert  & =0\\
\left\Vert 1_{R}\right\Vert  & =1\\
\left\Vert r+s\right\Vert  & \leq\left\Vert r\right\Vert +\left\Vert s\right\Vert ,\textrm{ }\forall r,s\in R\\
\left\Vert rs\right\Vert  & \leq\left\Vert r\right\Vert \left\Vert s\right\Vert ,\textrm{ }\forall r,s\in R
\end{align}
If, in addition, $\left\Vert \cdot\right\Vert $ satisfies the ultrametric
inequality:
\begin{equation}
\left\Vert f+g\right\Vert \leq\max\left\{ \left\Vert f\right\Vert ,\left\Vert g\right\Vert \right\} 
\end{equation}
with equality whenever $\left\Vert f\right\Vert \neq\left\Vert g\right\Vert $,
we say $\left\Vert \cdot\right\Vert $ is \textbf{non-archimedean}.
We say $\left\Vert \cdot\right\Vert $ is \textbf{multiplicative }if
$\left\Vert rs\right\Vert =\left\Vert r\right\Vert \left\Vert s\right\Vert $
for all $r,s\in R$.

We call $\left\Vert \cdot\right\Vert $ a \textbf{ring} \textbf{norm}
if $\left\Vert r\right\Vert =0$ if and only if $r=0_{R}$. A \textbf{ring
absolute value }is a ring norm which is multiplicative.

A \textbf{valued ring} is a ring with an absolute value. A \textbf{normed
ring }is a ring with a norm. A \textbf{completed valued ring }is a
valued ring so that the metric induced by the absolute value is complete.
A \textbf{Banach ring }is a normed ring so that the metric induced
by the norm is complete.

A \textbf{normed ($R$-)algebra }is an $R$-algebra $\mathcal{A}$
which, as a ring, is a normed ring, and which satisfies $\left\Vert r\alpha\right\Vert =\left\Vert r\right\Vert \left\Vert \alpha\right\Vert $
for all $r\in R$ and all $\alpha\in\mathcal{A}$. A \textbf{Banach
($R$-)algebra }is a normed algebra so that the metric induced by
the norm is complete.

Finally, by the \textbf{quality }of an absolute value, field, norm,
normed, space, Banach algebra, or the like, I refer to the property
of whether the thing under consideration is archimedean or non-archimedean.
Thus, we can have a seminorm of archimedean quality, a valued ring
of non-archimedean quality, and so on.
\end{defn}
We will also need Schwartz-Bruhat functions and related notions.
\begin{defn}
Let $R$ be a set and let $N\in\mathbb{N}_{0}$. Recall we say a function
$\phi:\mathbb{Z}_{p}\rightarrow R$ is \textbf{locally constant mod
$p^{N}$ }whenever:
\begin{equation}
\phi\left(\mathfrak{z}\right)=\phi\left(\left[\mathfrak{z}\right]_{p^{N}}\right),\textrm{ }\forall\mathfrak{z}\in\mathbb{Z}_{p}
\end{equation}
That is, $\phi$'s values depend only on the value of $\mathfrak{z}$
mod $p^{N}$. We say $\phi:\mathbb{Z}_{p}\rightarrow R$ is \textbf{locally
constant }whenever it is locally constant mod $p^{N}$ for some $N\in\mathbb{N}_{0}$.

Now, give $R$ a ring structure. Thanks to the profinite topology
of $\mathbb{Z}_{p}$, every locally constant function $\mathbb{Z}_{p}\rightarrow R$
is of the form:
\begin{equation}
\mathfrak{z}\in\mathbb{Z}_{p}\mapsto\sum_{n=1}^{N}r_{n}\left[\mathfrak{z}\overset{p^{e_{n}}}{\equiv}k_{n}\right]\in R
\end{equation}
for some constants $r_{n}\in R$, and integers $e_{n}\geq0$ and $k_{n}\in\left\{ 0,\ldots,p^{e_{n}}-1\right\} $.
These form an $R$-algebra under pointwise addition and multiplication,
with scalar multiplication by elements of $R$. We denote this $R$-algebra
by $\mathcal{S}\left(\mathbb{Z}_{p},R\right)$, and refer to its elements
as \textbf{$\left(p,R\right)$-adic Schwartz-Bruhat functions}, or
\textbf{SB functions}, for short. We then write $\mathcal{S}\left(\mathbb{Z}_{p},R\right)^{\prime}$
to denote the dual of $\mathcal{S}\left(\mathbb{Z}_{p},R\right)$.
This is the $R$-module of $R$-linear maps $\mathcal{S}\left(\mathbb{Z}_{p},R\right)\rightarrow R$
under pointwise addition and scalar multiplication. Elements of $\mathcal{S}\left(\mathbb{Z}_{p},R\right)^{\prime}$
are called \textbf{$\left(p,R\right)$-adic} \textbf{distributions}.
Given $d\mu\in\mathcal{S}\left(\mathbb{Z}_{p},R\right)^{\prime}$,
we write:
\begin{equation}
\int_{\mathbb{Z}_{p}}\phi\left(\mathfrak{z}\right)d\mu\left(\mathfrak{z}\right)
\end{equation}
to denote the image of an SB function $\phi$ under $d\mu$. Note
that $\mathcal{S}\left(\mathbb{Z}_{p},R\right)^{\prime}$ is an $\mathcal{S}\left(\mathbb{Z}_{p},R\right)$-module,
where the action of a given $\phi\in\mathcal{S}\left(\mathbb{Z}_{p},R\right)$
on a given $d\mu\in\mathcal{S}\left(\mathbb{Z}_{p},R\right)^{\prime}$
is to send $d\mu$ to $\phi\left(\mathfrak{z}\right)d\mu\left(\mathfrak{z}\right)$,
the distribution defined by:
\begin{equation}
\int_{\mathbb{Z}_{p}}\psi\left(\mathfrak{z}\right)\left(\phi\left(\mathfrak{z}\right)d\mu\left(\mathfrak{z}\right)\right)\overset{\textrm{def}}{=}\int_{\mathbb{Z}_{p}}\left(\psi\left(\mathfrak{z}\right)\phi\left(\mathfrak{z}\right)\right)d\mu\left(\mathfrak{z}\right),\textrm{ }\forall\psi\in\mathcal{S}\left(\mathbb{Z}_{p},R\right)
\end{equation}
We say a distribution $d\mu$ is \textbf{translation-invariant }whenever:
\begin{equation}
\int_{\mathbb{Z}_{p}}\phi\left(\mathfrak{z}+\mathfrak{a}\right)d\mu\left(\mathfrak{z}\right)=\int_{\mathbb{Z}_{p}}\phi\left(\mathfrak{z}\right)d\mu\left(\mathfrak{z}\right),\textrm{ }\forall\mathfrak{a}\in\mathbb{Z}_{p},\textrm{ }\forall\phi\in\mathcal{S}\left(\mathbb{Z}_{p},R\right)
\end{equation}
When $R$ has either characteristic zero or positive characteristic
co-prime to $p$, there exists a unique \textbf{Haar (probability)
distribution}; this is the unique translation-invariant $\left(p,R\right)$-adic
distribution with:
\begin{equation}
\int_{\mathbb{Z}_{p}}d\mu\left(\mathfrak{z}\right)=1
\end{equation}
We denote the Haar distribution by $d\mathfrak{z}$. It satisfies:
\begin{equation}
\int_{\mathbb{Z}_{p}}\left[\mathfrak{z}\overset{p^{n}}{\equiv}k\right]d\mathfrak{z}=\frac{1}{p},\textrm{ }\forall n\geq0,\textrm{ }\forall k\in\left\{ 0,\ldots,p^{n}-1\right\} 
\end{equation}
\end{defn}
\begin{rem}
When $R=\mathbb{R}$ or $\mathbb{C}$, the definitions of Schwartz-Bruhat
functions and distributions coincide with their classical counterparts,
such as in Tate's Thesis (\cite{Tate's thesis}) or Chapter 12 of
Washington's book on cyclotomic fields (\cite{washington cyclo}).
\end{rem}
\begin{rem}
Unless mentioned otherwise, we will only be considering the case where
$R$ has either characteristic $0$, or a positive characteristic
co-prime to $p$. This is essential; translation-invariant distributions
cannot exist without it.
\end{rem}
\begin{rem}
The two most important identities you need to know are:
\begin{equation}
\left[\mathfrak{z}\overset{p^{n}}{\equiv}k\right]=\frac{1}{p}\sum_{\left|t\right|_{p}\leq p^{n}}e^{2\pi i\left\{ t\left(\mathfrak{z}-k\right)\right\} _{p}},\textrm{ }\forall\mathfrak{z}\in\mathbb{Z}_{p},\textrm{ }\forall n\in\mathbb{N}_{0},\textrm{ }\forall k\in\left\{ 0,\ldots,p^{n}-1\right\} 
\end{equation}
\begin{equation}
\sum_{\left|t\right|_{p}=p^{n}}e^{2\pi i\left\{ t\left(\mathfrak{z}-k\right)\right\} _{p}}=p^{n}\left[\mathfrak{z}\overset{p^{n}}{\equiv}k\right]-p^{n-1}\left[\mathfrak{z}\overset{p^{n-1}}{\equiv}k\right],\textrm{ }\forall\mathfrak{z}\in\mathbb{Z}_{p},\textrm{ }\forall n\in\mathbb{N}_{1},\textrm{ }\forall k\in\left\{ 0,\ldots,p^{n}-1\right\} 
\end{equation}
Note that these equalities hold in any ring $R$ which is of characteristic
co-prime to $p$ and contains all $p^{n}$th roots of unity.
\end{rem}
For us, the essential difference between measures and distributions
is that measures are distributions which are continuous with respect
to a non-trivial topology. For technical reasons, all of our measures
will be defined using Wiener algebras.
\begin{defn}
Let $R$ be an algebraically closed ring with a norm $\left\Vert \cdot\right\Vert _{R}$
with respect to which $R$ is complete as a metric space. If $R$
is non-archimedean, we require the characteristic of its residue field
to be co-prime to $p$. Under these hypotheses: we define the \textbf{$\left(p,R\right)$-adic
Wiener algebra}, denoted $W\left(\mathbb{Z}_{p},R\right)$ like so:

\textbullet{} If $R$ is non-archimedean, $W\left(\mathbb{Z}_{p},R\right)$
is the space of all functions $f:\mathbb{Z}_{p}\rightarrow R$ possessing
a Fourier series representation which converges in $R$ uniformly
over $\mathbb{Z}_{p}$.

\textbullet{} If $R$ is archimedean, $W\left(\mathbb{Z}_{p},R\right)$
is the space of all functions $f:\mathbb{Z}_{p}\rightarrow R$ possessing
a Fourier series representation which is absolutely convergent with
respect to $\left\Vert \cdot\right\Vert _{R}$. (Note that $f$'s
Fourier series is then also convergent in $R$ uniformly over $\mathbb{Z}_{p}$.)

Regardless of the quality of $R$, $W\left(\mathbb{Z}_{p},R\right)$
is a Banach ring with respect to the norm:
\begin{equation}
\left\Vert f\right\Vert _{W\left(\mathbb{Z}_{p},R\right)}\overset{\textrm{def}}{=}\sup_{t\in\hat{\mathbb{Z}}_{p}}\left\Vert \hat{f}\left(t\right)\right\Vert _{R}
\end{equation}
where $\hat{f}$ is the Fourier transform of $f$ and $\left|\cdot\right|_{R}$
is the absolute value on $R$. The multiplication operation is pointwise
multiplication of functions. \textbf{Young's convolution inequality
}of abstract harmonic analysis then shows that this is indeed a Banach
ring. When $\left\Vert \cdot\right\Vert _{R}$ is an absolute value,
$W\left(\mathbb{Z}_{p},R\right)$ is a Banach $R$-algebra.

We write $W\left(\mathbb{Z}_{p},R\right)^{\prime}$ to denote the
continuous dual of $W\left(\mathbb{Z}_{p},R\right)$. We call elements
of $W\left(\mathbb{Z}_{p},R\right)^{\prime}$ \textbf{$\left(p,R\right)$-adic
measures}. If the absolute value (resp. norm) on $R$ is written,
say, $\left|\cdot\right|_{q}$ (resp. $\left\Vert \cdot\right\Vert _{q}$),
we call the measures \textbf{$\left(p,q\right)$-adic measures}. I
also use \textbf{$\infty$-adic }to mean that archimedean absolute
values are in play; for example, a ``$\left(3,\infty\right)$-adic''
measure means an element of $W\left(\mathbb{Z}_{3},\mathbb{C}\right)^{\prime}$.
\end{defn}
\begin{rem}
When $R$ is non-archimedean, $W\left(\mathbb{Z}_{p},R\right)$ is
\emph{equal to} the Banach algebra of continuous functions $\mathbb{Z}_{p}\rightarrow R$,
in stark contrast to the archimedean case. Regardless of the quality
of $R$, however, the Fourier-Stieltjes transform is an isometric
isomorphism of Banach rings $W\left(\mathbb{Z}_{p},R\right)^{\prime}\rightarrow B\left(\hat{\mathbb{Z}}_{p},R\right)$,
where $B\left(\hat{\mathbb{Z}}_{p},R\right)$ is the Banach ring of
bounded functions $\hat{\mathbb{Z}}_{p}\rightarrow R$ under pointwise
multiplication. See \cite{2nd blog paper} for details.
\end{rem}

\subsection{\label{subsec:Introduction}Introduction}

The main subject of this paper are F-series. These are power series
expressions (sometimes formal, sometimes not) that arise as functions
out of $\mathbb{Z}_{p}$.
\begin{defn}
\label{def:F-series}Let $r_{0},\ldots,r_{p-1}$ be indeterminates,
and let $R$ be an integral domain. A \textbf{$p$-adic F-series over
$R$} is a function $X:\mathbb{Z}_{p}\rightarrow R\left[\left[r_{0},\ldots,r_{p-1}\right]\right]$
of the shape:
\begin{equation}
X\left(\mathfrak{z}\right)\overset{\textrm{def}}{=}\sum_{n=0}^{\infty}f_{n}\left(\mathfrak{z}\right)r_{0}^{n}\prod_{k=1}^{p-1}\left(r_{k}/r_{0}\right)^{\#_{p:k}\left(\left[\mathfrak{z}\right]_{p^{n}}\right)}\label{eq:general p-adic F-series}
\end{equation}
where $\left\{ f_{n}\right\} _{n\geq0}$ is a sequence of functions
in $\mathcal{S}\left(\mathbb{Z}_{p},R\right)$. The $r_{k}$s are
referred to as the \textbf{parameters }of the F-series. Other quantities
may also be referred to as the parameters, depending on the values
of $f$.
\end{defn}
\begin{rem}
Even though we have an $r_{k}/r_{0}$ in (\ref{eq:general p-adic F-series}),
it can be shown that for any $n$, the exponent of $a_{0}$ in:
\begin{equation}
r_{0}^{n}\prod_{k=1}^{p-1}\left(r_{k}/r_{0}\right)^{\#_{p:k}\left(\left[\mathfrak{z}\right]_{p^{n}}\right)}
\end{equation}
will be non-negative. Indeed, by definition:
\begin{equation}
\#_{p:0}\left(\left[\mathfrak{z}\right]_{p^{n}}\right)=\lambda_{p}\left(\left[\mathfrak{z}\right]_{p^{n}}\right)-\sum_{k=1}^{p-1}\#_{p:k}\left(\left[\mathfrak{z}\right]_{p^{n}}\right)
\end{equation}
and so:
\begin{equation}
r_{0}^{n}\prod_{k=1}^{p-1}\left(r_{k}/r_{0}\right)^{\#_{p:k}\left(\left[\mathfrak{z}\right]_{p^{n}}\right)}=r_{0}^{n-\lambda_{p}\left(\left[\mathfrak{z}\right]_{p^{n}}\right)}\prod_{k=0}^{p-1}r_{k}^{\#_{p:k}\left(\left[\mathfrak{z}\right]_{p^{n}}\right)}\label{eq:rewritten product in F-series}
\end{equation}
which is a non-negative exponent for $r_{0}$, since $\lambda_{p}\left(\left[\mathfrak{z}\right]_{p^{n}}\right)\leq n$
for all $\mathfrak{z}\in\mathbb{Z}_{p}$ and all $n\in\mathbb{N}_{0}$.
\end{rem}
That being said, using the functional equations from \textbf{Proposition
\ref{prop:fundamental functional equations-1}}, it is a simple, straightforward
task to verify that certain kinds of F-series can be characterized
as formal solution of systems of affine linear functional equations.
\begin{thm}
The F-series:
\begin{equation}
X\left(\mathfrak{z}\right)=\sum_{n=0}^{\infty}r_{0}^{n}\prod_{k=1}^{p-1}\left(r_{k}/r_{0}\right)^{\#_{p:k}\left(\left[\mathfrak{z}\right]_{p^{n}}\right)}
\end{equation}
is the unique formal solution of the functional equations:
\begin{equation}
X\left(p\mathfrak{z}+j\right)=r_{j}X\left(\mathfrak{z}\right)+1,\textrm{ }\forall\mathfrak{z}\in\mathbb{Z}_{p},\textrm{ }\forall j\in\left\{ 0,\ldots,p-1\right\} 
\end{equation}
Meanwhile, the F-series:
\begin{equation}
X\left(\mathfrak{z}\right)=\sum_{n=0}^{\infty}b_{\left[\theta_{p}^{\circ n}\left(\mathfrak{z}\right)\right]_{p}}a_{0}^{n}\prod_{k=1}^{p-1}\left(a_{k}/a_{0}\right)^{\#_{p:k}\left(\left[\mathfrak{z}\right]_{p^{n}}\right)}\label{eq:foliated F-series}
\end{equation}
is the unique formal solution of the functional equations:
\begin{equation}
X\left(p\mathfrak{z}+j\right)=a_{j}X\left(\mathfrak{z}\right)+b_{j},\textrm{ }\forall\mathfrak{z}\in\mathbb{Z}_{p},\textrm{ }\forall j\in\left\{ 0,\ldots,p-1\right\} \label{eq:a_j X plus b_j}
\end{equation}
\end{thm}
\begin{rem}
In (\ref{eq:foliated F-series}), note that we can write:
\begin{equation}
b_{\left[\theta_{p}^{\circ n}\left(\mathfrak{z}\right)\right]_{p}}=\sum_{j=0}^{p-1}\hat{b}_{j}\left(\varepsilon_{n}\left(\mathfrak{z}\right)\right)^{j}
\end{equation}
where the $b_{j}$s and $\hat{b}_{j}$s are finite Fourier transforms
of one another:
\begin{equation}
\hat{b}_{j}=\frac{1}{p}\sum_{k=0}^{p-1}b_{k}e^{-2\pi ijk/p}
\end{equation}
\begin{equation}
b_{j}=\sum_{k=0}^{p-1}\hat{b}_{k}e^{2\pi ijk/p}
\end{equation}
\end{rem}
\begin{rem}
Henceforth, unless stated otherwise, all F-series are foliated and
$p$-adic. We then refer to the $a_{j}$s and $b_{j}$s as the \textbf{parameters
}of the F-series.
\end{rem}
\begin{defn}
\label{def:foliated}Letting $a_{0},\ldots,a_{p-1}$ and $b_{0},\ldots,b_{p-1}$
be indeterminates, and letting $R$ be an integral domain, a \textbf{$p$-adic}
\textbf{foliated F-series over $R$ of degree $1$ }is a function
$X:\mathbb{Z}_{p}\rightarrow R\left[\left[a_{0},\ldots,a_{p-1}\right]\right]\left[b_{0},\ldots,b_{p-1}\right]$
satisfying the functional equations: 
\begin{equation}
X\left(p\mathfrak{z}+k\right)=a_{k}X\left(\mathfrak{z}\right)+b_{k},\textrm{ }\forall\mathfrak{z}\in\mathbb{Z}_{p},\textrm{ }\forall k\in\left\{ 0,\ldots,p-1\right\} \label{eq:degree 1}
\end{equation}
More generally, given an integer $n\geq1$, a \textbf{$p$-adic foliated
F-series over $R$ of degree $n$ }is the pointwise product of $n$
foliated F-series of degree $1$. Thus, if $X$ and $Y$ are degree
$1$, then $X^{3}$ and $X^{2}Y$ are both of degree $3$. A \textbf{$p$-adic}
\textbf{foliated F-series over $R$ of degree $0$ }is defined to
be a constant function $\mathbb{Z}_{p}\rightarrow\left\{ c\right\} \subset R$.
\end{defn}
\begin{rem}
With the help of $\theta_{p}$, the functional equations (\ref{eq:degree 1})
can be written in a single line:
\begin{equation}
X\left(\mathfrak{z}\right)=a_{\left[\mathfrak{z}\right]_{p}}X\left(\theta_{p}\left(\mathfrak{z}\right)\right)+b_{\left[\mathfrak{z}\right]_{p}},\textrm{ }\forall\mathfrak{z}\in\mathbb{Z}_{p}
\end{equation}
\end{rem}
\begin{assumption}
HENCEFORTH, UNLESS STATED OTHERWISE, ALL F-SERIES ARE $p$-ADIC AND
FOLIATED.
\end{assumption}
\begin{rem}
It's important to note that the definition given in (\ref{def:foliated})
is primarily a linguistic one, meant to give us the words we need
to talk about F-series. Realizing the kind of ``universal'' Fourier
analysis promised in the paper's abstract will require us to be more
particular in how we define the codomain of an F-series. The details
are covered in \textbf{Section \ref{subsec:Ascent-from-Fractional}}.
\end{rem}
In my PhD dissertation, I considered F-series where the parameters
were all constants in a global field and showed that, modulo certain
technical conditions, F-series had well-defined Fourier transforms
\cite{My Dissertation}. For our purposes, a function $\hat{X}$ is
a Fourier transform of the F-series $X$ when the limit of the Fourier
series generated by $\hat{X}$:
\begin{equation}
\lim_{N\rightarrow\infty}\sum_{\left|t\right|_{p}\leq p^{N}}\hat{X}\left(t\right)e^{2\pi i\left\{ t\mathfrak{z}\right\} _{p}}\label{eq:sense of Fourier transform}
\end{equation}
converges to $X$ in a suitable sense, to be made precise in \textbf{Section
\ref{subsec:Quasi-Integrability-=000026-Degenerate}}. The task of
making rigorous sense of (\ref{eq:sense of Fourier transform}) in
the general case compelled me to develop what I now call a \textbf{frame}
\cite{My Dissertation,my frames paper}. This paper will give an exposition
of frames, but because of the concept's technical nature, the explanation
will be withheld until \textbf{Section \ref{subsec:Frame-Theory-=000026}},
so as to not get in the way of the main ideas presented in \textbf{Section
\ref{sec:The-Big-Idea}}. For now, all you need to know is that, in
general, the limit (\ref{eq:sense of Fourier transform}) will converge
at many points $\mathfrak{z}\in\mathbb{Z}_{p}$, however, the topology
of convergence (that is, the metric being used to define the limit)
will vary from point to point. \textbf{Section \ref{subsec:Frame-Theory-=000026}}
begins with a concrete example of how this can happen.

Because of these technical details, the way in which $X$'s Fourier
transform exists in the sense of (\ref{eq:sense of Fourier transform}\textbf{)
}involves enough of a departure from classical Fourier analysis that
I had to coin the notion of \textbf{quasi-integrability }to deal with
it. I say that $X$ is \textbf{quasi-integrable }when it has a Fourier
transform in the sense of (\ref{eq:sense of Fourier transform}\textbf{)}.
\textbf{Section \ref{subsec:Quasi-Integrability-=000026-Degenerate}}
explains these features in full. Finally, using $\hat{X}$, I showed
that F-series can be realized as measures on $\mathbb{Z}_{p}$, and
more generally, as $p$-adic distributions.

The following theorem summarizes these findings.
\begin{thm}[M.C. Siegel. (2022 - 2024) \cite{My Dissertation,first blog paper}]
\label{thm:The-breakdown-variety}Let $a_{0},\ldots,a_{p-1},b_{0},\ldots,b_{p-1}$
be constants in $K$, with none of the $a_{j}$s being $0$, and with
$a_{0}\neq1$. Then, there exists a unique function $X:\mathbb{N}_{0}\rightarrow K$
satisfying the system of functional equations:
\begin{equation}
X\left(pn+j\right)=a_{j}X\left(n\right)+b_{j},\textrm{ }\forall n\in\mathbb{N}_{0},\forall j\in\left\{ 0,\ldots,p-1\right\} \label{eq:N_0 functional equations for a foilated F-series}
\end{equation}
If, on the other hand, $a_{0}=1$ and $b_{0}=0$, there are infinitely
many solutions, one for each initial value $X\left(0\right)\in K$.

In either case, given a solution $X:\mathbb{N}_{0}\rightarrow K$
of (\ref{eq:N_0 functional equations for a foilated F-series}), suppose
that there exists a reffinite $K$-frame $\mathcal{F}$ on $\mathbb{Z}_{p}$
so that the limit:
\begin{equation}
\lim_{n\rightarrow\infty}X\left(\left[\mathfrak{z}\right]_{p^{n}}\right)
\end{equation}
is $\mathcal{F}$-convergent for all $\mathfrak{z}\in\mathbb{Z}_{p}$.
Then, this limit defines an $\mathcal{F}$-compatible, rising-continuous
function, which we shall also denote by $X$. Moreover, $X$ is the
unique $\mathcal{F}$-compatible, rising-continuous function satisfying
the functional equations:
\begin{equation}
X\left(p\mathfrak{z}+j\right)=a_{j}X\left(\mathfrak{z}\right)+b_{j},\textrm{ }\forall\mathfrak{z}\in\mathbb{Z}_{p},\forall j\in\left\{ 0,\ldots,p-1\right\} \label{prop:fundamental functional equations}
\end{equation}
which is equal to $X\left(0\right)$ at $\mathfrak{z}=0$.

Next, defining the functions $\alpha_{X},\beta_{X},\gamma_{X},\hat{A}_{X}:\hat{\mathbb{Z}}_{p}\rightarrow K\left(\zeta_{p^{\infty}}\right)$
by:
\begin{equation}
\alpha_{X}\left(t\right)\overset{\textrm{def}}{=}\frac{1}{p}\sum_{j=0}^{p-1}a_{j}e^{-2\pi ijt}
\end{equation}
\begin{equation}
\beta_{X}\left(t\right)\overset{\textrm{def}}{=}\frac{1}{p}\sum_{j=0}^{p-1}b_{j}e^{-2\pi ijt}
\end{equation}
\begin{equation}
\gamma_{X}\left(t\right)\overset{\textrm{def}}{=}\frac{\beta_{X}\left(t\right)}{\alpha_{X}\left(t\right)}
\end{equation}
\begin{equation}
\hat{A}_{X}\left(t\right)\overset{\textrm{def}}{=}\prod_{m=0}^{-v_{p}\left(t\right)-1}\alpha_{X}\left(p^{m}t\right)
\end{equation}
$X$ then admits an $\mathcal{F}$-convergent Fourier series representation:
\begin{equation}
X\left(\mathfrak{z}\right)\overset{\mathcal{F}}{=}\sum_{t\in\hat{\mathbb{Z}}_{p}}\hat{X}\left(t\right)e^{2\pi i\left\{ t\mathfrak{z}\right\} _{p}}\label{eq:X's Fourier series}
\end{equation}
where $\hat{X}:\hat{\mathbb{Z}}_{p}\rightarrow K\left(\zeta_{p^{\infty}}\right)$
is given by:
\begin{equation}
\hat{X}\left(t\right)=\begin{cases}
0 & \textrm{if }t=0\\
\left(\beta_{X}\left(0\right)v_{p}\left(t\right)+\gamma_{X}\left(\frac{t\left|t\right|_{p}}{p}\right)\right)\hat{A}_{X}\left(t\right) & \textrm{if }t\neq0
\end{cases},\textrm{ }\forall t\in\hat{\mathbb{Z}}_{p}
\end{equation}
when $\alpha_{X}\left(0\right)=1$, and is given by:
\begin{equation}
\hat{X}\left(t\right)=\begin{cases}
\frac{\beta_{X}\left(0\right)\hat{A}_{X}\left(0\right)}{1-\alpha_{X}\left(0\right)} & \textrm{if }t=0\\
\frac{\beta_{X}\left(0\right)\hat{A}_{X}\left(t\right)}{1-\alpha_{X}\left(0\right)}+\gamma_{X}\left(\frac{t\left|t\right|_{p}}{p}\right)\hat{A}_{X}\left(t\right) & \textrm{if }t\neq0
\end{cases},\textrm{ }\forall t\in\hat{\mathbb{Z}}_{p}
\end{equation}
when $\alpha_{X}\left(0\right)\neq1$. Summing (\ref{eq:X's Fourier series})
yields an $\mathcal{F}$-convergent representation of $X$ as a foliated
$p$-adic F-series:
\begin{equation}
X\left(\mathfrak{z}\right)\overset{\mathcal{F}}{=}\sum_{n=0}^{\infty}\left(\sum_{j=0}^{p-1}\hat{b}_{j}\left(\varepsilon_{n}\left(\mathfrak{z}\right)\right)^{j}\right)a_{0}^{n}\prod_{k=1}^{p-1}\left(a_{k}/a_{0}\right)^{\#_{k}\left(\left[\mathfrak{z}\right]_{p^{n}}\right)}
\end{equation}
Finally, in the spirit of the Parseval-Plancherel identity, the formula:
\begin{equation}
\int_{\mathbb{Z}_{p}}\phi\left(\mathfrak{z}\right)X\left(\mathfrak{z}\right)d\mathfrak{z}\overset{\textrm{def}}{=}\sum_{t\in\hat{\mathbb{Z}}_{p}}\hat{\phi}\left(t\right)\hat{X}\left(-t\right),\textrm{ }\forall\phi\in\mathcal{S}\left(\mathbb{Z}_{p},K\right)\label{eq:Par-Plan formula}
\end{equation}
defines a $K$-valued $p$-adic distribution $X\left(\mathfrak{z}\right)d\mathfrak{z}$,
which we identify with $X$. Furthermore, if there is a place $\ell$
of $K$ so that $\sup_{t\in\hat{\mathbb{Z}}_{p}}\left|\hat{X}\left(t\right)\right|_{\ell}<\infty$,
the above formula holds for all $\phi\in W\left(\mathbb{Z}_{p},\overline{\overline{K_{\ell}}}\right)$,
where $\overline{\overline{K_{\ell}}}$ denotes the metric completion
of the algebraic closure of the metric completion of $K$ with respect
to $\ell$, in which case (\ref{eq:Par-Plan formula}) defines $X\left(\mathfrak{z}\right)d\mathfrak{z}$
as $\left(p,\ell\right)$-adic measure.
\end{thm}
\begin{rem}
When $a_{0}=a_{1}=\cdots=a_{p-1}$, note that: 
\begin{equation}
\alpha_{X}\left(t\right)=\frac{a_{0}}{p}\sum_{k=0}^{p-1}e^{-2\pi ikt}
\end{equation}
By the geometric series formula, this will vanish for all $\left|t\right|_{p}=p$.
In this case, as written, $\hat{A}_{X}\left(t\right)=0$ occurs for
all $t\neq0$, seeing as:
\begin{equation}
\hat{A}_{X}\left(t\right)=\prod_{n=0}^{-v_{p}\left(t\right)-1}\alpha_{X}\left(p^{n}t\right)
\end{equation}
and so, if $\left|t\right|_{p}\geq p$, the product will contain the
term $\alpha_{X}\left(p^{-v_{p}\left(t\right)-1}t\right)=\alpha_{X}\left(t\left|t\right|_{p}/p\right)$.
Since $t\left|t\right|_{p}/p$ has $p$-adic absolute value $p$ for
all $t\neq0$, this makes $\alpha_{X}\left(t\left|t\right|_{p}/p\right)=0$.
\emph{Nevertheless}, this is not a problem, because, regardless of
the value of $\alpha_{X}\left(0\right)$, both formulae for $\hat{X}\left(t\right)$
given in this theorem have:
\begin{equation}
\gamma_{X}\left(\frac{t\left|t\right|_{p}}{p}\right)\hat{A}_{X}\left(t\right)=\frac{\beta_{X}\left(t\left|t\right|_{p}/p\right)}{\alpha_{X}\left(t\left|t\right|_{p}/p\right)}\hat{A}_{X}\left(t\right)
\end{equation}
As such, the $\alpha_{X}$ in the denominator of $\gamma_{X}$ will
cancel out the vanishing $\alpha_{X}$ in $\hat{A}_{X}\left(t\right)$.
We could deal with this issue by defining $\hat{A}_{X}$ as:
\begin{equation}
\hat{A}_{X}\left(t\right)\overset{\textrm{alt. def}}{=}\prod_{n=0}^{-v_{p}\left(t\right)-2}\hat{\alpha}_{X}\left(t\right)
\end{equation}
but this needlessly complicates the computations, so we shall leave
things alone. Indeed, in the universal case we consider, we will never
need to worry about this, because the $a_{j}$s will be unrelated
formal indeterminates, and thus, $\alpha_{X}\left(t\right)$ will
be non-vanishing for all $t$.
\end{rem}
Using \textbf{Theorem \ref{thm:The-breakdown-variety}} as a springboard,
this paper extends its ideas to yield three-and-a-half new findings.
These are stated below with full rigor, however, the statements of
the results contain many technical conditions which, for clarity's
sake, will not be defined until \textbf{Section \ref{sec:The-Small-Ideas}}.
For this reason, it's worth giving the three-and-a-half main results
an informal presentation before confronting them in full.

The first of our main results is the title of this paper: when viewed
as distributions, F-series form algebras under \emph{pointwise multiplication}.
(Moreover, in the case when the F-series in question can all be realized
as measures taking values in the same space, these measures form an
algebra under pointwise multiplication.) That is, if $X$ and $Y$
are F-series, both of which have Fourier transforms in the sense of
(\ref{eq:sense of Fourier transform}), then, provided a simple condition
on $X$'s and $Y$'s parameters holds true, the pointwise product
$XY$ also has a Fourier transform, in the sense that there is a function
$\hat{Z}$ so that:
\begin{equation}
\lim_{N\rightarrow\infty}\sum_{\left|t\right|_{p}\leq p^{N}}\hat{Z}\left(t\right)e^{2\pi i\left\{ t\mathfrak{z}\right\} _{p}}=X\left(\mathfrak{z}\right)Y\left(\mathfrak{z}\right)
\end{equation}
pointwise in a suitable fashion. Moreover, much like $\hat{X}$ and
$\hat{Y}$ realize $X$ and $Y$ as distributions/measures, $\hat{Z}$
does the same for the product $XY$. Furthermore, the condition on
$X$ and $Y$ needed for $\hat{Z}$ to exist then guarantees the existence
of a Fourier transform for $X^{m}Y^{n}$ for all integers $m,n\geq0$.
This is highly unusual.

It is a classic result of the theory of distributions that while distributions
(and, more generally, measures) can be \emph{convolved} with one another
to obtain more distributions (resp., measures), the \emph{pointwise
products} of distributions/measures generally do not exist. Because
distributions/measures have Fourier transforms that need not decay
to $0$ at $\infty$, the convolution of these Fourier transforms
will tend to be divergent or non-convergent, and as that convolution
corresponds to the Fourier transform of whatever the  distributions'/measures'
pointwise product would be, if it existed, the convolution's non-existence
prevents the pointwise product from existing, either. What is noteworthy
about our result is that it does not overcome this obstacle so much
as it steps \emph{around }it. In this respect, the computations given
in \textbf{Section \ref{subsec:The-Central-Computation}} can be viewed
as a method of ``renormalizing'' the divergent/non-convergent convolutions
to produce a consistent answer.

Stated rigorously, this main result is as follows:
\begin{cor}
\label{cor:main result}Let $d$ and $p$ be positive integers, with
$d\geq1$ and $p\geq2$, let $K$ be a global field (if $\textrm{char}K>0$,
we require that $\textrm{char}K$ be co-prime to $p$). For each $j\in\left\{ 1,\ldots,d\right\} $
and $k\in\left\{ 0,\ldots,p-1\right\} $ let $a_{j,k}$ and $b_{j,k}$
be indeterminates, and let $R_{d}\left(K\right)$ be the ring of polynomials
in the $a_{j,k}$s and $b_{j,k}$s with coefficients in $\mathcal{O}_{K}$,
so that $R_{d}\left(K\right)$ is then a free $\mathcal{O}_{K}$-algebra
in $2dp$ indeterminates. Then, let:
\begin{equation}
\mathcal{R}_{d}\left(K\right)\overset{\textrm{def}}{=}\frac{R\left[x_{1},\ldots,x_{d}\right]}{\left\langle \left(1-a_{1,0}\right)x_{1}-1,\ldots,\left(1-a_{d,0}\right)x_{d}-1\right\rangle }
\end{equation}
be the localization of $R$ away from the prime ideal $\left\langle 1-a_{1,0},\ldots,1-a_{d,0}\right\rangle $.

Writing:
\begin{equation}
\kappa_{j}\left(n\right)\overset{\textrm{def}}{=}\prod_{k=1}^{p-1}\left(\frac{a_{j,k}}{a_{j,0}}\right)^{\#_{p:k}\left(n\right)}
\end{equation}
let $\mathcal{F}$ be a digital $\mathcal{R}_{d}\left(K\right)$-frame
on $\mathbb{Z}_{p}$ so that:
\begin{equation}
\left|a_{j,0}\right|_{\mathcal{F}\left(\mathfrak{z}\right)}\limsup_{n\rightarrow\infty}\prod_{k=1}^{p-1}\left|\frac{a_{j,k}}{a_{j,0}}\right|_{\mathcal{F}\left(\mathfrak{z}\right)}^{\#_{p:k}\left(\left[\mathfrak{z}\right]_{p^{n}}\right)}<1,\textrm{ }\forall\mathfrak{z}\in D\left(\mathcal{F}\right),\textrm{ }\forall j\in\left\{ 1,\ldots,d\right\} \label{eq:root condition for individual j-1-1}
\end{equation}
Finally, let $I$ be any ideal of $R_{d}\left(K\right)$ so that $\left\langle 1-a_{j,0}\right\rangle \nsubseteq I$
for all $j\in\left\{ 1,\ldots,d\right\} $. Then:

\vphantom{}I. For each $j\in\left\{ 1,\ldots,d\right\} $, there
exists a unique rising-continuous function $X_{j}\in C\left(\mathcal{F}/I\mathcal{R}_{d}\left(K\right)\right)$
so that:
\begin{equation}
X_{j}\left(\mathfrak{z}\right)\overset{\left(\mathcal{F}/I\mathcal{R}_{d}\left(K\right)\right)\left(\mathfrak{z}\right)}{=}a_{j,\left[\mathfrak{z}\right]_{p}}X_{j}\left(\theta_{p}\left(\mathfrak{z}\right)\right)+b_{j,\left[\mathfrak{z}\right]_{p}},\textrm{ }\forall\mathfrak{z}\in D\left(\mathcal{F}/I\mathcal{R}_{d}\left(K\right)\right)\label{eq:functional equations-2}
\end{equation}
where $\left(\mathcal{F}/I\mathcal{R}_{d}\left(K\right)\right)\left(\mathfrak{z}\right)$
is the completion of $\mathcal{R}_{d}\left(K\right)/I\mathcal{R}_{d}\left(K\right)$
with respect to the absolute value that $\mathcal{F}/I\mathcal{R}_{d}\left(K\right)$
associates to $\mathfrak{z}$, and where $\mathcal{F}/I\mathcal{R}_{d}\left(K\right)$
is the quotient frame induced by $\mathcal{F}$ and the ideal $I\mathcal{R}_{d}\left(K\right)\subseteq\mathcal{R}_{d}\left(K\right)$.

\vphantom{}

II. For any $\mathbf{n}\in\mathbb{N}_{0}^{d}$, the function $X_{\mathbf{n}}=\prod_{j=1}^{d}X_{j}^{n_{j}}\in C\left(\mathcal{F}/I\mathcal{R}_{d}\left(K\right)\right)$
is $\mathcal{F}/I\mathcal{R}_{d}\left(K\right)$-quasi-integrable.
Moreover, we denote a Fourier transform of this function by:
\begin{equation}
\hat{X}_{\mathbf{n}}=\bigstar_{j=1}^{d}\hat{X}_{j}^{*n_{j}}\left(t\right)
\end{equation}
where $\bigstar$ denotes formal convolution, with $\hat{X}_{j}^{*n_{j}}\left(t\right)$
being defined as the convolution identity element $\mathbf{1}_{0}$
whenever $n_{j}=0$. Here, $\bigstar_{j=1}^{d}\hat{X}_{j}^{*n_{j}}$
is a function $\hat{\mathbb{Z}}_{p}\rightarrow\textrm{Frac}\left(\mathcal{R}_{d}\left(K\right)/I\mathcal{R}_{d}\left(K\right)\right)\left(\zeta_{p^{\infty}}\right)$.

Using \textbf{Theorem \ref{thm:powers of X}} (see page \pageref{thm:powers of X})
to compute $\hat{X}_{j}$ for all $j\in\left\{ 1,\ldots,d\right\} $,
we then have the following formula for $\hat{X}_{\mathbf{n}}$:

\begin{equation}
\hat{X}_{\mathbf{n}}\left(t\right)=\hat{f}_{\mathbf{n}}\left(t\right)-\hat{g}_{\mathbf{n}}\left(t\right)
\end{equation}
where:
\begin{equation}
\hat{f}_{\mathbf{n}}\left(t\right)=\sum_{\mathbf{m}<\mathbf{n}}\sum_{n=0}^{-v_{p}\left(t\right)-2}\left(\prod_{m=0}^{n-1}\alpha_{\mathbf{n}}\left(p^{m}t\right)\right)\alpha_{\mathbf{m},\mathbf{n}}\left(p^{n}t\right)\hat{X}_{\mathbf{m}}\left(p^{n+1}t\right)\label{eq:f_script J hat-1-1}
\end{equation}
where the sum is taken over all $\mathbf{m}=\left(m_{1},\ldots,m_{d}\right)$
so that $m_{j}\leq n_{j}$ for all $j\in\left\{ 1,\ldots,d\right\} $
and $m_{j}<n_{j}$ for at least one $j\in\left\{ 1,\ldots,d\right\} $,
and where:
\begin{equation}
\hat{g}_{\mathbf{n}}\left(t\right)=\begin{cases}
0 & \textrm{if }t=0\\
\left(\beta_{\mathbf{n}}\left(0\right)v_{p}\left(t\right)+\gamma_{\mathbf{n}}\left(\frac{t\left|t\right|_{p}}{p}\right)\right)\hat{A}_{\mathbf{n}}\left(t\right) & \textrm{if }t\neq0
\end{cases},\textrm{ }\forall t\in\hat{\mathbb{Z}}_{p}\label{eq:g_script J hat , alpha equals 1-1-2}
\end{equation}
if $\alpha_{\mathbf{n}}\left(0\right)=1$ and:
\begin{equation}
\hat{g}_{\mathbf{n}}\left(t\right)=\begin{cases}
\frac{\beta_{\mathbf{n}}\left(0\right)}{1-\alpha_{\mathbf{n}}\left(0\right)} & \textrm{if }t=0\\
\left(\frac{\beta_{\mathbf{n}}\left(0\right)}{1-\alpha_{\mathbf{n}}\left(0\right)}+\gamma_{\mathbf{n}}\left(\frac{t\left|t\right|_{p}}{p}\right)\right)\hat{A}_{\mathbf{n}}\left(t\right) & \textrm{if }t\neq0
\end{cases},\textrm{ }\forall t\in\hat{\mathbb{Z}}_{p}\label{eq:g_script J hat , alpha not equal to 1-1-2}
\end{equation}
if $\alpha_{\mathbf{n}}\left(0\right)\neq1$, where, for all $\mathbf{m},\mathbf{n}\in\mathbb{N}_{0}^{d}$
with $m_{j}\leq n_{j}$ for all $j\in\left\{ 1,\ldots,d\right\} $,
$\alpha_{\mathbf{m},\mathbf{n}},\beta_{\mathbf{n}},\gamma_{\mathbf{n}},\hat{A}_{\mathbf{n}}:\hat{\mathbb{Z}}_{p}\rightarrow\textrm{Frac}\left(\mathcal{R}_{d}\left(K\right)/I\mathcal{R}_{d}\left(K\right)\right)\left(\zeta_{p^{\infty}}\right)$
are given by:
\begin{equation}
\alpha_{\mathbf{m},\mathbf{n}}\left(t\right)\overset{\textrm{def}}{=}\frac{1}{p}\sum_{k=0}^{p-1}r_{\mathbf{m},\mathbf{n},k}e^{-2\pi ikt}
\end{equation}
\begin{equation}
\beta_{\mathbf{n}}\left(t\right)\overset{\textrm{def}}{=}\frac{1}{p}\sum_{k=0}^{p-1}c_{\mathbf{n},k}e^{-2\pi ikt}
\end{equation}
\begin{equation}
\gamma_{\mathbf{n}}\left(t\right)\overset{\textrm{def}}{=}\frac{\beta_{\mathbf{n}}\left(t\right)}{\alpha_{\mathbf{n}}\left(t\right)}
\end{equation}
\begin{equation}
\hat{A}_{\mathbf{n}}\left(t\right)\overset{\textrm{def}}{=}\prod_{n=0}^{-v_{p}\left(t\right)-1}\alpha_{\mathbf{n}}\left(p^{n}t\right)
\end{equation}
where:
\begin{equation}
c_{\mathbf{n},k}=-\sum_{\mathbf{m}<\mathbf{n}}r_{\mathbf{m},k}\hat{X}_{\mathbf{m}}\left(0\right)
\end{equation}
\begin{equation}
r_{\mathbf{m},\mathbf{n},k}\overset{\textrm{def}}{=}\binom{\mathbf{n}}{\mathbf{m}}a_{\mathbf{m},k}b_{\mathbf{n}-\mathbf{m},k}
\end{equation}
\begin{equation}
r_{\mathbf{n},k}\overset{\textrm{def}}{=}r_{\mathbf{n},\mathbf{n},k}
\end{equation}
\begin{equation}
a_{\mathbf{m},k}\overset{\textrm{def}}{=}\prod_{j=1}^{d}a_{j,k}^{m_{j}}
\end{equation}
\begin{equation}
b_{\mathbf{n}-\mathbf{m},k}\overset{\textrm{def}}{=}\prod_{j=1}^{d}b_{j,k}^{n_{j}-m_{j}}
\end{equation}
\begin{equation}
\alpha_{\mathbf{n}}\left(t\right)\overset{\textrm{def}}{=}\alpha_{\mathbf{n},\mathbf{n}}\left(t\right)
\end{equation}

\vphantom{}III. For all $j\in\left\{ 1,\ldots,d\right\} $, $\hat{X}_{j}$
will be the Fourier-Stieltjes transform of a distribution in $\mathcal{S}\left(\mathbb{Z}_{p},\textrm{Frac}\left(\mathcal{R}_{d}\left(K\right)/I\mathcal{R}_{d}\left(K\right)\right)\right)^{\prime}$,
which we denote by $X_{j}\left(\mathfrak{z}\right)d\mathfrak{z}$.
Then, for any $\mathbf{n}\in\mathbb{N}_{0}^{d}$, we define a pointwise
product:
\begin{equation}
X_{\mathbf{n}}\left(\mathfrak{z}\right)d\mathfrak{z}\overset{\textrm{def}}{=}\prod_{j=1}^{d}\left(X_{j}\left(\mathfrak{z}\right)d\mathfrak{z}\right)^{n_{j}}\overset{\textrm{def}}{=}\left(\prod_{j=1}^{d}X_{j}^{n_{j}}\left(\mathfrak{z}\right)\right)d\mathfrak{z}
\end{equation}
where $\left(\prod_{j=1}^{d}X_{j}^{n_{j}}\left(\mathfrak{z}\right)\right)d\mathfrak{z}$
is the distribution whose Fourier-Stieltjes transform is $\bigstar_{j=1}^{d}\hat{X}_{j}^{*n_{j}}$,
as given in (II). Consequently, we get a polynomial algebra $K\left[X_{1}\left(\mathfrak{z}\right)d\mathfrak{z},\ldots,X_{d}\left(\mathfrak{z}\right)d\mathfrak{z}\right]$.
Moreover, the Fourier transform induces an isomorphism from this algebra
to the $K$-algebra generated by $\hat{X}_{1},\ldots,\hat{X}_{d}$
under convolution.

\textbf{WARNING}: In the convolution algebra generated by $\hat{X}_{1},\ldots,\hat{X}_{d}$,
it will generally \uline{not} be the case that $\bigstar_{j=1}^{d}\hat{X}_{j}^{*n_{j}}$
will be given by an actual convolution; that is, it will not be given
by iterated applications of the $*$ operator defined by:
\begin{equation}
\left(\hat{f}*\hat{g}\right)\left(t\right)\overset{\textrm{def}}{=}\sum_{s\in\hat{\mathbb{Z}}_{p}}\hat{f}\left(s\right)\hat{g}\left(t-s\right)
\end{equation}
because that sum will be either divergent or non-convergent. Rather,
we \emph{define }$\bigstar_{j=1}^{d}\hat{X}_{j}^{*n_{j}}$ by the
formula from (II).
\end{cor}
\begin{rem}
The reason \textbf{Corollary \ref{cor:main result} }is a Corollary
rather than a Theorem is because it is deduced from \textbf{Theorem
\ref{thm:powers of X}}, which considers the \textbf{power case},
where $d=1$. 
\end{rem}
\textbf{Corollary \ref{cor:main result}} also contains the paper's
second main result; this lies in our use of the rings $R_{d}\left(K\right)$
and $\mathcal{R}_{d}\left(K\right)$. In \textbf{Theorem \ref{thm:The-breakdown-variety}},
the parameters of $X$ are scalars, as was the case in my dissertation.
However, both while I was writing my dissertation and in the years
that have followed since its completion in 2022 and 2025, I had a
habit of treating those scalars as variables/formal indeterminates.
I find it's easier to show how this works than it is to definite it
rigorously. As such, consider the following example:
\begin{example}
\label{exa:choosing absolute values}Beginning with the elementary
generating function identity:
\begin{equation}
\prod_{n=0}^{N-1}\left(1+az^{2^{n}}\right)=\sum_{n=0}^{2^{N}-1}a^{\#_{2:1}\left(n\right)}z^{n}
\end{equation}
and letting $q$ and $p$ be, say, any non-zero rational numbers,
it is easy to show that, for any $n\geq0$, the locally constant function
$\mathfrak{z}\in\mathbb{Z}_{2}\mapsto q^{\#_{2:1}\left(\left[\mathfrak{z}\right]_{2^{n}}\right)}\in\mathbb{Q}$
has the Fourier transform:
\begin{equation}
\int_{\mathbb{Z}_{2}}q^{\#_{1}\left(\left[\mathfrak{z}\right]_{2^{n}}\right)}e^{-2\pi i\left\{ t\mathfrak{z}\right\} _{2}}d\mathfrak{z}=\left(q+1\right)^{n-\max\left\{ 0,-v_{2}\left(t\right)\right\} }\frac{\mathbf{1}_{0}\left(2^{n}t\right)}{2^{n}}\hat{\Pi}_{q}\left(t\right)\label{eq:fourier transform of q to the number of 1s}
\end{equation}
where the function $\hat{\Pi}_{q}:\hat{\mathbb{Z}}_{2}\rightarrow\mathbb{Q}\left(\zeta_{2^{\infty}}\right)$
is defined by:
\begin{equation}
\hat{\Pi}_{q}\left(t\right)\overset{\textrm{def}}{=}\begin{cases}
1 & \textrm{if }t=0\\
\prod_{n=0}^{-v_{2}\left(t\right)-1}\left(1+qe^{-2\pi i\left(2^{n}t\right)}\right) & \textrm{else}
\end{cases}\label{eq:Definition of Big Pi without N}
\end{equation}
Here, we are using the fact that for any function $f$ which is locally
constant mod $p^{n}$:
\[
\int_{\mathbb{Z}_{p}}f\left(\mathfrak{z}\right)e^{-2\pi i\left\{ t\mathfrak{z}\right\} _{p}}d\mathfrak{z}=\frac{1}{p^{n}}\sum_{m=0}^{p^{n}-1}f\left(m\right)e^{-2\pi imt},\textrm{ }\forall t\in\hat{\mathbb{Z}}_{p}
\]
Summing (\ref{eq:fourier transform of q to the number of 1s}) from
$n=0$ to $n=N-1$, we obtain:
\begin{equation}
\int_{\mathbb{Z}_{2}}\sum_{n=0}^{N-1}\frac{q^{\#_{2:1}\left(\left[\mathfrak{z}\right]_{2^{n}}\right)}}{p^{n}}e^{-2\pi i\left\{ t\mathfrak{z}\right\} _{2}}d\mathfrak{z}=\begin{cases}
\frac{N+\min\left\{ 0,v_{2}\left(t\right)\right\} }{\left(2p\right)^{\max\left\{ 0,-v_{2}\left(t\right)\right\} }}\hat{\Pi}_{q}\left(t\right) & \textrm{if }\frac{q+1}{2p}=1\\
\frac{1-\left(\frac{q+1}{2p}\right)^{N+\min\left\{ 0,v_{2}\left(t\right)\right\} }}{1-\frac{q+1}{2p}}\frac{\hat{\Pi}_{q}\left(t\right)}{\left(2p\right)^{\max\left\{ 0,-v_{2}\left(t\right)\right\} }} & \textrm{else}
\end{cases}\label{eq:transform of sum of q p}
\end{equation}

Now, suppose $\left(q+1\right)/2p\notin\left\{ 1,-1\right\} $. Then,
there is a place $\ell$ of $\mathbb{Q}$ so that: 
\begin{equation}
\left|\frac{q+1}{2p}\right|_{\ell}<1
\end{equation}
Consequently, the right-hand side of (\ref{eq:transform of sum of q p})
converges $\ell$-adically as $N\rightarrow\infty$:
\begin{equation}
\lim_{N\rightarrow\infty}\frac{1-\left(\frac{q+1}{2p}\right)^{N+\min\left\{ 0,v_{2}\left(t\right)\right\} }}{1-\frac{q+1}{2p}}\frac{\hat{\Pi}_{q}\left(t\right)}{\left(2p\right)^{\max\left\{ 0,-v_{2}\left(t\right)\right\} }}\overset{\mathbb{C}_{\ell}}{=}\frac{1}{1-\frac{q+1}{2p}}\frac{\hat{\Pi}_{q}\left(t\right)}{\left(2p\right)^{\max\left\{ 0,-v_{2}\left(t\right)\right\} }}\label{eq:limit with pi hat}
\end{equation}
This takes the value $1/\left(1-\frac{q+1}{2p}\right)$ when $t=0$.
If $\ell$ is non-archimedean and $\left|2p\right|_{\ell}=1$, the
ultrametric inequality gives us:
\begin{equation}
\left|\frac{\hat{\Pi}_{q}\left(t\right)}{\left(2p\right)^{-v_{2}\left(t\right)}}\right|_{\ell}=\prod_{n=0}^{-v_{2}\left(t\right)-1}\left|\frac{1+qe^{-2\pi i2^{n}t}}{2p}\right|_{\ell}\leq\prod_{n=0}^{-v_{2}\left(t\right)-1}\max\left\{ 1,\left|qe^{-2\pi i2^{n}t}\right|_{\ell}\right\} \leq1
\end{equation}
which shows that right-hand side of (\ref{eq:limit with pi hat})
is then $\ell$-adically bounded. This then shows that:
\begin{equation}
\lim_{N\rightarrow\infty}\sum_{n=0}^{N-1}\frac{q^{\#_{2:1}\left(\left[\mathfrak{z}\right]_{2^{n}}\right)}}{p^{n}}d\mathfrak{z}
\end{equation}
converges in $W\left(\mathbb{Z}_{2},\mathbb{C}_{\ell}\right)^{\prime}$,
and hence, that we may identify the F-series:
\begin{equation}
X\left(\mathfrak{z}\right)=\sum_{n=0}^{\infty}\frac{q^{\#_{2:1}\left(\left[\mathfrak{z}\right]_{2^{n}}\right)}}{p^{n}}
\end{equation}
with the $\left(2,\ell\right)$-adic measure $X\left(\mathfrak{z}\right)d\mathfrak{z}$
which acts on functions $\phi\in W\left(\mathbb{Z}_{2},\mathbb{C}_{\ell}\right)$
by the formula:
\begin{equation}
\int_{\mathbb{Z}_{2}}\phi\left(\mathfrak{z}\right)X\left(\mathfrak{z}\right)d\mathfrak{z}\overset{\mathbb{C}_{\ell}}{=}\frac{1}{1-\frac{q+1}{2p}}\sum_{t\in\hat{\mathbb{Z}}_{2}}\frac{\hat{\phi}\left(t\right)\hat{\Pi}_{q}\left(-t\right)}{\left(2p\right)^{\max\left\{ 0,-v_{2}\left(-t\right)\right\} }}
\end{equation}
Note that this same analysis holds when $q$ and $p$ are \emph{algebraic}
numbers, not just rational numbers. Moreover, there is no reason we
must limit ourselves to working with non-archimedean places. Indeed,
given and algebraic integers $p,q$ so that $\left|q\right|_{\infty}<\left|p\right|_{\infty}$,
where $\left|\cdot\right|_{\infty}$ is the complex absolute value,
observing that $\#_{2:1}\left(\left[\mathfrak{z}\right]_{2^{n}}\right)$
is maximized when $\mathfrak{z}\overset{2^{n}}{\equiv}2^{n}-1$, with
$\#_{2:1}\left(\left[\mathfrak{z}\right]_{2^{n}}\right)=n$. This
gives us the bound:
\begin{equation}
\left|\sum_{n=0}^{\infty}\frac{q^{\#_{2:1}\left(\left[\mathfrak{z}\right]_{2^{n}}\right)}}{p^{n}}\right|_{\infty}\leq\sum_{n=0}^{\infty}\left|\frac{q}{p}\right|_{\infty}^{n}<\infty
\end{equation}
By the M-test, we see that $X$ converges absolutely in $\mathbb{C}$
uniformly with respect to $\mathfrak{z}\in\mathbb{Z}_{2}$. Not only
does this mean that $X:\mathbb{Z}_{2}\rightarrow\mathbb{C}$ is continuous
and absolutely integrable with respect to the real-valued $2$-adic
Haar probability measure, but also that, for any $n\geq2$:
\begin{equation}
\int_{\mathbb{Z}_{2}}\left|X^{n}\left(\mathfrak{z}\right)\right|_{\infty}d\mathfrak{z}\leq\int_{\mathbb{Z}_{2}}\left(\sum_{m=0}^{\infty}\left|\frac{q}{p}\right|_{\infty}^{m}\right)^{n}d\mathfrak{z}<\frac{1}{\left(1-\left|\frac{q}{p}\right|_{\infty}\right)^{n}}<\infty
\end{equation}
and thus, that the pointwise product $X^{n}$ is continuous and absolutely
integrable as well, with the Fourier integral:
\begin{equation}
\int_{\mathbb{Z}_{2}}\left(\sum_{m=0}^{\infty}\frac{q^{\#_{2:1}\left(\left[\mathfrak{z}\right]_{2^{m}}\right)}}{p^{m}}\right)^{n}e^{-2\pi i\left\{ t\mathfrak{z}\right\} _{2}}d\mathfrak{z}
\end{equation}
being absolutely convergent in $\mathbb{C}$ for all $n\geq1$, uniformly
with respect to $t$. \textbf{Corollary \ref{cor:main result}} can
be viewed as a generalization of this to the case where the pointwise
product of our F-series is no longer integrable in the classical sense,
but rather must be viewed as a measure.
\end{example}
As a rule, in classical analysis, both $p$-adic and otherwise, the
ambient field in which the analysis occurs is fixed, as is the absolute
value used to induce a metric on said field. One doesn't usually stop
to consider what might happen if our absolute value was changed to
an inequivalent one. However, as \textbf{Example \ref{exa:choosing absolute values}}
shows, when working with F-series, it is not only easy, but \emph{natural}
to want to consider working with potentially \emph{any} absolute value
on the background field $K$ for which expressions happen to converge
in some sense. As I explained in my dissertation, instead of restricting
ourselves to working in a single metric completion of the background
field, we can\textemdash and, in many cases, \emph{must}\textemdash be
flexible enough to choose whichever absolute values happen to work
best at any given moment. Much of the impetus behind my developed
of the concept of frames in my dissertation was implement this flexibility
in a rigorous way. However, even there, an obstacle remained: the
frames I used were ad-hoc construction generally only worked for a
given choice of values for the F-series' parameters.

For example, on the one hand:
\begin{equation}
\text{\ensuremath{\sum_{n=0}^{\infty}\frac{3^{\#_{2:1}\left(\left[\mathfrak{z}\right]_{2^{n}}\right)}}{4^{n}}}}
\end{equation}
is a uniformly continuous real-valued function on $\mathbb{Z}_{2}$,
while: 
\begin{equation}
\text{\ensuremath{\sum_{n=0}^{\infty}\frac{5^{\#_{2:1}\left(\left[\mathfrak{z}\right]_{2^{n}}\right)}}{4^{n}}}}\label{eq:4 5}
\end{equation}
is a $\left(2,3\right)$-adic and $\left(2,\infty\right)$-adic measure
(yes, $3$, not $5$; this is because $\ell=3$ and $\ell=\infty$
are the two places of $\mathbb{Q}$ for which $\left(q+1\right)/2p=6/8=3/4$
has $\ell$-adic absolute value $<1$). Moreover, it's clear that
both of these F-series are merely special cases of:
\begin{equation}
\text{\ensuremath{\sum_{n=0}^{\infty}\frac{q^{\#_{2:1}\left(\left[\mathfrak{z}\right]_{2^{n}}\right)}}{4^{n}}}}
\end{equation}
where $q$ is an indeterminate. However, interpreting these three
F-series as functions with $\mathbb{Z}_{2}$ as their domain requires
working in completely different spaces for their codomains. For the
first, we need to deal with continuous functions $\mathbb{Z}_{2}\rightarrow\mathbb{R}$;
for the second, we can work with functions $\mathbb{Z}_{2}\rightarrow\mathbb{Z}_{5}$,
provided we sum (\ref{eq:4 5}) using the topology of the reals when
$\mathfrak{z}\in\mathbb{N}_{0}$ and using the topology of the $5$-adics
when $\mathfrak{z}\in\mathbb{Z}_{2}^{\prime}$. For the third, we
can treat the series as taking values in various rings of formal power
series, such as $\mathbb{Z}\left[\frac{1}{4}\right]\left[\left[q\right]\right]$
or $\mathbb{Q}\left[\left[q\right]\right]$. (Though, even there,
subtleties arise. If we wished to evaluate $q$ at some non-zero rational
number $q_{0}=a_{0}/b_{0}$ by quotienting our ring by the ideal $\left\langle b_{0}q-a_{0}\right\rangle $,
note that $\mathbb{Q}\left[\left[q\right]\right]/\left\langle b_{0}q-a_{0}\right\rangle $
would be trivial, as $b_{0}q-a_{0}$ is a unit of $\mathbb{Q}\left[\left[q\right]\right]$,
while $\mathbb{Z}\left[\frac{1}{4}\right]\left[\left[q\right]\right]/\left\langle b_{0}q-a_{0}\right\rangle $
would vary, depending on the value of $b_{0}$. We'll see more of
this in \textbf{Section \ref{subsec:Ascent-from-Fractional}}, specifically
\textbf{Fact \ref{fact:Let--be}} on page \pageref{fact:Let--be}.)

To that end, it would be very nice to have a formalism that puts all
of these different cases on equal footing. This is precisely the second
main accomplishment of this paper: the creation of such a formalism.
As explained in \textbf{Section \pageref{subsec:Frame-Theory-=000026}},
frames are naturally (and, likely, functorially) compatible with ring-theoretic
quotients, thanks to the algebraic properties of ring seminorms. In
this way, we can take a page from classical algebraic geometry and
view specific F-series like (\ref{eq:4 5}) as being obtained from
more general F-series taking values in rings of formal power series
by quotients by appropriately chosen ideals. \textbf{Corollary \pageref{cor:main result}}
shows that Fourier transforms of products of F-series exist in this
generalized context, and that the existence of those Fourier transforms
is naturally compatible with quotients of the background rings. This
furnishes the ``universal'', purely formal Fourier theory promised
in the abstract.

Much like the second main result, the paper's third main result also
comes from my desire to give a firm foundation to an informal computational
procedure that I've been using for years without rigorous justification.
The most general version of the computation is presented in \textbf{Example
\ref{exa:proceeding formally}} on page \pageref{exa:proceeding formally}
of \textbf{Section \ref{subsec:The-Central-Computation}}. The idea
comes from the observation that the functional equations:
\begin{equation}
X\left(p\mathfrak{z}+k\right)=a_{k}X\left(\mathfrak{z}\right)+b_{k}
\end{equation}
satisfied by the F-series $X$ should translate to functional equations
satisfied by $\hat{X}$. In that vein, if we \emph{assume }that $X$
is ``integrable'' (for the appropriate sense of ``integrable''),
we can obtain a functional equation for $\hat{X}$ by using a standard
$p$-adic change of variables in conjunction with $X$'s functional
equation. Namely:
\begin{align*}
\hat{X}\left(t\right) & =\int_{\mathbb{Z}_{p}}X\left(\mathfrak{z}\right)e^{-2\pi i\left\{ t\mathfrak{z}\right\} _{p}}d\mathfrak{z}\\
 & =\sum_{k=0}^{p-1}\int_{p\mathbb{Z}_{p}+k}X\left(\mathfrak{z}\right)e^{-2\pi i\left\{ t\mathfrak{z}\right\} _{p}}d\mathfrak{z}\\
\left(\mathfrak{x}=\frac{\mathfrak{z}-k}{p}\right); & =\frac{1}{p}\int_{\mathbb{Z}_{p}}X\left(p\mathfrak{x}+k\right)e^{-2\pi i\left\{ t\left(p\mathfrak{x}+k\right)\right\} _{p}}d\mathfrak{x}\\
 & =\frac{1}{p}\int_{\mathbb{Z}_{p}}\left(a_{k}X\left(\mathfrak{x}\right)+b_{k}\right)e^{-2\pi ikt}e^{-2\pi i\left\{ pt\mathfrak{x}\right\} _{p}}d\mathfrak{x}\\
 & =\frac{1}{p}\int_{\mathbb{Z}_{p}}\left(a_{k}X\left(\mathfrak{x}\right)+b_{k}\right)e^{-2\pi ikt}e^{-2\pi i\left\{ pt\mathfrak{x}\right\} _{p}}d\mathfrak{x}\\
 & =\frac{1}{p}\sum_{k=0}^{p-1}a_{k}e^{-2\pi ikt}\int_{\mathbb{Z}_{p}}X\left(\mathfrak{x}\right)e^{-2\pi i\left\{ pt\mathfrak{x}\right\} _{p}}d\mathfrak{x}+\frac{1}{p}\sum_{k=0}^{p-1}b_{k}e^{-2\pi ikt}\int_{\mathbb{Z}_{p}}e^{-2\pi i\left\{ pt\mathfrak{x}\right\} _{p}}d\mathfrak{x}\\
 & =\underbrace{\left(\frac{1}{p}\sum_{k=0}^{p-1}a_{k}e^{-2\pi ikt}\right)}_{\alpha_{X}\left(t\right)}\hat{X}\left(pt\right)+\underbrace{\left(\frac{1}{p}\sum_{k=0}^{p-1}b_{k}e^{-2\pi ikt}\right)}_{\beta_{X}\left(t\right)}\mathbf{1}_{0}\left(pt\right)
\end{align*}
By this \emph{purely formal }argument (no limits or topologies involved),
this yields the functional equation:
\begin{equation}
\hat{X}\left(t\right)=\alpha_{X}\left(t\right)\hat{X}\left(pt\right)+\beta_{X}\left(t\right)\mathbf{1}_{0}\left(pt\right)\label{eq:functional equation in t}
\end{equation}
The question is: assuming this equation even has a solution, are its
solution(s) Fourier transforms of $X$ in the sense of (\ref{eq:sense of Fourier transform})?
Moreover, can we use the natural generalization of this computation
to obtain Fourier transforms of powers of $X$, or products of $X$
and other F-series? The answer is presented in \textbf{Theorem \ref{thm:formal solutions}}
on page \pageref{thm:formal solutions} of \textbf{Section \ref{sec:encoding}}.
In short, the answer is that as long as (\ref{eq:functional equation in t})
has a \emph{unique }solution, that solution will then be a Fourier
transform of $X$. However, when the solution fails to be unique,
or even fails to exist altogether, we have to resort to other means
to compute $X$'s Fourier transform. Not only is this a useful fact
for computations and further exploration, it resonates with the ``universal''
Fourier theory by showing that, save for important singular cases,
one can, in fact, derive the Fourier transform of an F-series through
purely formal considerations like the ones above.

The paper's final, third-and-a-halfth accomplishment, is intimately
intertwined with its third main result. By some strange mechanism
whose nature I've yet to fully pin down, there is an interplay between
the solvability of (\ref{eq:functional equation in t}), the properties
of the distributions/measure obtained by solving it, and a certain
algebraic variety associated to $X$ which I call the \textbf{breakdown
variety}.\textbf{ }So far, this behavior is most comparable to the
relationship between the eigenvalues and eigenspaces of linear operators,
with the breakdown variety playing the role of the eigenvalue. For
example, given an affine algebraic variety $V$, we can construct
a distribution $dM$, viewed here as a linear functional $\mathcal{S}\left(\mathbb{Z}_{p},R\right)\rightarrow R$
for some ring $R$, so that ``evaluating'' $dM$ at any given affine
point $\mathbf{p}$ (i.e., considering the map induced by $dM$ on
the quotient $\mathcal{S}\left(\mathbb{Z}_{p},R/I\right)\rightarrow R/I$
by the maximal ideal $I$ corresponding to $\mathbf{p}$) yields a
distribution which has the property of being what I call \textbf{degenerate
}if and only if $\mathbf{p}\in V$. This evinces a novel means of
``encoding'' affine algebraic varieties through distributions, one
which interacts with the pointwise multiplicative structure on the
algebra generated by these distributions. There are also interactions
with the tensor product; see \textbf{Proposition \ref{prop:tensors}
}on page \ref{prop:tensors} for more. I count this finding as only
half of an accomplishment because, as of the writing of this paper,
I've yet to comprehensively work out the details of this dictionary
between distributions and varieties, least of all its robustness.

\subsection{An Outline of this Paper}

\vphantom{}

\textbf{\uline{Section \mbox{\ref{sec:The-Big-Idea}}}}

The methods of mathematical analysis needed to rigorously prove this
paper's main results are unorthodox. I cannot point you to prior examples
of them outside of my own work. As such, I have chosen to first present
the central computations in a purely formal matter. Nearly all of
the opaque or unexplained concepts brought up in what I've presented
so far arise as part of the labor of justify the computations we are
about to do. Because of this, I have made a \emph{deliberate decision
}to invert the typical order of mathematical exposition and withhold
a full explanation of concepts like frames and quasi-integrability
until you have seen the material that forced me to devise these concepts
in the first place. Explicating these concepts up front would be dry,
unmotivated, and unenlightening; it would be an ungrounded solution
in search of problem to which to anchor itself.

\textbf{Section \ref{subsec:M-functions}} introduces what I call
\textbf{M-functions}. M-functions are the building blocks of F-series,
and this subsection introduces elementary identities and definitions
that will be invoked throughout everything that follows.

\textbf{Section \ref{subsec:The-Central-Computation}} turns to our
main object of study: a family of F-series $X_{j}$ indexed by $j\in\left\{ 1,\ldots,d\right\} $,
where $d$ is an integer $\geq1$. This subsection is one long computation,
devoted almost entirely to showing that a given formula formally sums
to the Fourier series of a given product of F-series $\prod_{j=1}^{d}X_{j}^{n_{j}}$,
where $n_{1},\ldots,n_{d}$ are any non-negative integers. The subsection
begins with some important notation that will be used throughout the
paper; for example, we write $X_{\mathbf{n}}$ to denote the product
$\prod_{j=1}^{d}X_{j}^{n_{j}}$, where $\mathbf{n}=\left(n_{1},\ldots,n_{d}\right)\in\mathbb{N}_{0}^{d}$.
The most important details left unaddressed in this subsection are
contained in \textbf{Assumption \ref{assu:main limit lemma}} on page
\pageref{assu:main limit lemma}. \textbf{Theorem \ref{thm:formal quasi-integrability of X script J}}
on page \pageref{thm:formal quasi-integrability of X script J} is
the main result of this section, the quasi-integrability of $X_{\mathbf{n}}$,
where an explicit formula for the function's Fourier transform is
also given. The rest of the subsection is devoted to computing various
identities that will be used later on in the paper.

Overall, \textbf{Section \ref{sec:The-Big-Idea} }is almost entirely
elementary, consisting primarily of a variety of detailed but straightforward
Fourier analytic computations.

\vphantom{}

\textbf{\uline{Section \mbox{\ref{sec:The-Small-Ideas}}}}

Before we can even begin to show that \textbf{Assumption \ref{assu:main limit lemma}}
is justified, we need to have the proper conceptual frameworks to
deal with the questions therein.

\textbf{Section \ref{subsec:Frame-Theory-=000026}} is a streamlined
exposition of the latest iteration of the theory of \textbf{frames}
that I initially developed in my dissertation \cite{My Dissertation}
to deal with the novel forms of convergence that occurred with the
sequences of functions I happened to be working with. The main tools
used are the notions of seminorms, norms, and absolute values on rings,
in particular, the Berkovich spectrum. Implicitly, we will also be
dealing with locally convex topologies on algebras and rings. I must
stress that knowledge of either Berkovich spaces or locally convex
topological vector spaces/rings/algebras is in no way needed to understand
what's going on, though some familiarity with locally convex topological
vector spaces will make things much easier to appreciate.

\textbf{Section \ref{subsec:Ascent-from-Fractional}} is dedicated
to giving a rigorous formulation of what has, up until this point,
been a merely informal practice on my part: treating the parameters
of an F-series as ``variables''. In essence, what we want is a means
of treating the parameters of one or more F-series as if we were working
with the indeterminates of an equation for defining an algebraic variety
along classical lines in terms of coordinate rings (the quotient of
a ring of polynomials over a field by an ideal corresponding to the
locus of points on a given algebraic surface or curve). The main technical
issue at hand is to build a formalism that is compatible with the
theory of frames presented in \textbf{Section \ref{subsec:Frame-Theory-=000026}}.
This can be done in the ``natural''/universal way by using the notion
of a quotient frame from \textbf{Section \ref{subsec:Frame-Theory-=000026}}.
In this context, we arrive at a rigorous treatment of F-series by
viewing them as elements of a frame-theoretic completion of a localization
of a ring of formal polynomials over a base field. By the quotient
frame construction, this completion and the associated frame naturally
descend through the quotient of the ring of polynomials by any given
ideal. In doing so, we will then be able to rigorously justify \textbf{Section
\ref{subsec:The-Central-Computation}}'s computations in this satisfyingly
general context.

\textbf{Section \ref{subsec:Quasi-Integrability-=000026-Degenerate}}
begins with a brief overview of the essential ideas behind the Fourier
theory used throughout this paper. I repeat that all of the essential
information is available in the paper \cite{2nd blog paper}, which
is freely accessible on arXiv. Anyhow, the main content of this subsection
is in giving a rigorous definition of the concept of \textbf{quasi-integrability},
and an exploration of how the Fourier transform of quasi-integrable
function is not unique as in classical Fourier analysis, but rather
is only unique modulo a vector space (or module) of what I call \textbf{degenerate
measures}.

\vphantom{}

\textbf{\uline{Section \mbox{\ref{sec:Making-Things-Precise}}}}

This is the core of the paper, where the computations from \textbf{Section
\ref{subsec:The-Central-Computation}} are given their proper foundation
and fully justified.

\textbf{Section \ref{subsec:M-functions-=000026-Asymptotics}} revisits
M-functions in the context of frames. To justify \textbf{Section \ref{subsec:The-Central-Computation}}'s
computations, we must use the decay rate of the tail of the $X_{j}$'s
Fourier series (that is, the decay of the difference between each
$X_{j}$ and a partial sum of $X_{j}$'s Fourier series) to show that
the same is true for the tails of the Fourier series of the products
$X_{\mathbf{n}}$ for any $\mathbf{n}\in\mathbb{N}_{0}^{d}\backslash\left\{ \mathbf{0}\right\} $.
Doing this requires arguing by induction/recursion along the $\mathbf{n}$s,
and noticing that the decay happens to be bounded by an M-function.
\textbf{Section \ref{subsec:M-functions-=000026-Asymptotics}} introduces
the technical tools that will be used to show that sums and other
combinations/transformations of M-functions are once again bounded
by M-functions. Our guiding light here is the classic \textbf{Root
Test }for series convergence from elementary analysis. In brief, the
main condition we will need is that an M-function $\left\{ M_{n}\right\} _{n\geq0}$
is summable in $n$ at a given $\mathfrak{z}\in\mathbb{Z}_{p}$ with
respect to a given absolute value $\left|\cdot\right|_{\ell}$. Provided
that $\lim_{n\rightarrow\infty}\left|M_{n}\left(\mathfrak{z}\right)\right|_{\ell}^{1/n}$
exists in the reals and is $<1$, we can then bounded $M_{n}\left(\mathfrak{z}\right)$
by a sequence which decays geometrically with respect to $n$ and
then construct another, slightly larger M-function, which nevertheless
decays to $0$ as $n\rightarrow\infty$. For technical reasons, in
order to guarantee the existence of the root test limit, we will need
to restrict ourselves to those $p$-adic integers $\mathfrak{z}$
such that the densities of $0$s, $1$s, $\ldots$, $\left(p-1\right)$s
among the $p$-adic digits of $\mathfrak{z}$ are well-defined.

\textbf{Section \ref{subsec:Frame-Theoretic-Asymptotics}} begins
with a presentation of classic asymptotic notation (Big-Oh, Vinogradov,
etc.) modified to be compatible with frames. The focus of this section
is in proving \textbf{Lemma \ref{lem:main limit lemma}} (page \pageref{lem:main limit lemma}),
the paper's main lemma, which is used to justify \textbf{Assumption
\ref{assu:main limit lemma} }from \textbf{Section \ref{subsec:The-Central-Computation}}.
This lemma establishes the decay conditions on the tail of an F-series'
Fourier series representation needed to guarantee that limits of the
type from \textbf{Assumption \ref{assu:main limit lemma}} actually
hold true. This reduces the task to making an inductive argument that
shows that the satisfaction of those decay conditions for the individual
$X_{j}$s is essentially enough to guarantee that \textbf{Lemma \ref{lem:main limit lemma}}
can be applied to rigorously justify \textbf{Assumption \ref{assu:main limit lemma}}.

\textbf{Section \ref{subsec:Assumption--=000026}} is the inductive
step where everything comes together. Given $\mathbf{n}$ and an $\mathbf{m}$
so that all of $\mathbf{m}$'s entries are less than or equal to the
entries in the corresponding positions of $\mathbf{n}$ (denoted $\mathbf{m}\leq\mathbf{n}$),
by using the formulae from \textbf{Section \ref{subsec:The-Central-Computation}},
we get explicit estimates for the Fourier tail decay of $X_{\mathbf{n}}$
in terms of the tail decay of $X_{\mathbf{m}}$ for all $\mathbf{m}<\mathbf{n}$
(where the strict inequality means that at least one entry of $\mathbf{m}$
is strictly less than the entry in the corresponding position of $\mathbf{n}$).
After breaking up the formula for the tail decay into several terms,
\textbf{Lemma \ref{lem:main limit lemma}} in conjunction with the
results of \textbf{Section \ref{subsec:M-functions-=000026-Asymptotics}}
are used to show that the necessary induction happens: the tail decay
of the $X_{\mathbf{m}}$s implies that $X_{\mathbf{n}}$ also has
the necessary tail decay. This is done in \textbf{Theorem \ref{thm:inductive step}}
on page \pageref{thm:inductive step}. \textbf{Theorem \ref{thm:mini thesis}}
(page \pageref{thm:mini thesis}) serves as the base case, among other
things. As a result of this, we get then get our main result, \textbf{Corollary
\ref{cor:main result}}, after first proving the special case (the
\textbf{Power Case})\textbf{ }where $X_{i}=X_{j}$ for all $i,j$;
the power case is dealt with in \textbf{Theorem \ref{thm:powers of X}},
on page \pageref{thm:powers of X}.

\vphantom{}

\textbf{\uline{Section \mbox{\ref{sec:encoding}}}}

This section is in two subsections, the second of which is of a primarily
speculative nature.

\textbf{Section \ref{subsec:Formal-Solutions} }deals with the third
main result of the paper: justifying the situations in which one can
compute Fourier transforms of F-series using purely formal calculation
in the manner described near the end of the paper's Introduction (\textbf{Section
\ref{subsec:Introduction}}). After introducing some terminology to
deal with the recursive nature of the functional equations being dealt
with, we proceed to see how the distributions/measures induced by
a formal solution of the functional equations in $t$ can be computed
in closed form provided a certain linear operator on the space of
SB functions ends up being invertible. Both the invertibility of this
operator and the existence of unique solutions of the functional equations
in $t$ hinge on whether or not the parameters of the F-series under
consideration satisfy a certain algebraic relation, the locus of solutions
of which I call the \textbf{breakdown varieties} of the F-series involved.
By the end, we establish \textbf{Theorem \ref{thm:formal solutions}}
(see page \pageref{thm:formal solutions}), where it is shown that
the functional equations in $t$ have a unique solution if and only
if the F-series' parameters do not lie in the breakdown variety. Moreover,
whenever this unique solution exists, it is then necessarily the Fourier
transform of the F-series in question. In particular, given $\mathbf{n}$,
if we know Fourier transforms for $X_{\mathbf{m}}$ for all $\mathbf{m}<\mathbf{n}$,
then we can use the functional equations in $t$ to solve for a unique
formula for a Fourier transform of $X_{\mathbf{n}}$ if and only if
the parameters of the $X_{j}$s do not lie in the breakdown variety
of $X_{\mathbf{n}}$.

As the computations in this paper show, given an algebraic variety,
we can construct an F-series/measure which will be degenerate for
a given choice of parameters if and only if those parameters lie in
the F-series' breakdown variety. This, coupled with the results of
the current paper, hint at a greater structure present in the kinds
of measures we have been examining, one that can be used to detect
points in algebraic varieties and, more generally, to encode algebraic
varieties. \textbf{Section \ref{subsec:Varieties-=000026-Measures}}
deals with this issue as best as I currently can. After some examples
of the encoding in practice, w\textbf{ }show that, for F-series induced
distributions whose degeneracy/non-degeneracy encodes points of algebraic
varieties, the degeneracy is equivalent to the containment of the
distribution's kernel in the image of a certain linear operator. Furthermore,
while the effects of pointwise multiplication on the varieties encoded
by two or more F-series still need to be explored in greater detail,
it is an easy and immediate consequence of \textbf{Proposition \ref{prop:L_A and dA}}
that, when considering the vector space generated by finitely many
F-series distributions/measures, the natural tensor algebra structure
possessed by distributions/measures allows us to encode unions of
algebraic varieties. This is shown in \textbf{Proposition \ref{prop:tensors}}
on page \ref{prop:tensors}. There, we have that if $V_{1},\ldots,V_{d}$
are the breakdown varieties of the F-series distributions $dM_{1},\ldots,dM_{d}$,
then $V_{1}\cup\cdots\cup V_{d}$ is the breakdown variety of the
distribution $dM_{1}\otimes\cdots\otimes dM_{d}$.

\textbf{Section \ref{subsec:Arithmetic-Dynamics,-Varieties,}}, the
paper's final section, is primarily speculative. The original impetus
behind my PhD dissertation was to explore certain arithmetic dynamical
systems. One of the principal results of my dissertation was that,
in many cases, given the arithmetical dynamical system generated by
what I call a \textbf{hydra map}, one can construct an F-series, which
I call the \textbf{numen},\textbf{ }whose value distribution properties
control the dynamics of the associated hydra map, such as the map's
periodic points and cycles. \textbf{Theorem \ref{thm:The-breakdown-variety}}
shows that the numen of a hydra map can be realized as a measure,
and open the door for the application of Fourier-theoretic techniques
to study the dynamics of hydra maps by way of their numens. Together,
these two observations suggest it might be possible to improve our
knowledge of hydra maps' dynamics by fleshing out the dictionary between
algebraic varieties and F-series and then using this dictionary to
study the breakdown varieties encoded by their numens. To that end,
\textbf{Section \ref{subsec:Arithmetic-Dynamics,-Varieties,} }explores
a handful of examples of how algebraic varieties appear to be intertwined
with hydra maps' numens, in the hopes inspiring new avenues of future
work.
\begin{note}
Throughout this paper, I have strived to make my exposition as accessible
as possible. An advanced undergraduate mathematics education and/or
early graduate-level education should be more than sufficient preparation.
\end{note}

\section{\label{sec:The-Big-Idea}The Big Idea}

After many failed attempts at alternative structures, I've found that
the cleanest way to go about the proof is to do everything formally,
without regard to delicate issues of convergence, so as to make the
argument as clear as possible. In this section, I will present the
computational arc of the proof, while marking those steps where more
thorough justification is needed.

For psychological reasons, I recommend that readers turn off their
brains for a while and \emph{focus on the formal, symbol-pushing computation}.
The nature of this work is best understood as coming up with analysis-based
justifications for results that should arguably be provable in a purely
algebraic context.

\vphantom{}
\begin{assumption}
THROUGHOUT THIS SECTION, WE WORK WITH A FIXED GLOBAL FIELD $K$ AND
AN INTEGER $p\geq2$. IF $K$ HAS POSITIVE CHARACTERISTIC, WE REQUIRE
$\textrm{char}K$ TO BE CO-PRIME TO $p$.
\end{assumption}

\subsection{\label{subsec:M-functions}M-functions}

In the work to come, we will often find ourselves dealing with a product
of the form:
\begin{equation}
\prod_{k=0}^{n-1}r_{\left[\theta_{p}^{\circ k}\left(\mathfrak{z}\right)\right]_{p}}\label{eq:prototypical M-function}
\end{equation}
where $r_{0},\ldots,r_{p-1}\in K$ are constants in a global field,
and where $n\in\mathbb{N}_{0}$ and $\mathfrak{z}\in\mathbb{Z}_{p}$.
The notation is as follows: for each $\mathfrak{z}$ and each $k$,
compute the $p$-adic integer $\theta_{p}^{\circ k}\left(\mathfrak{z}\right)$,
and then set $j_{k}$ to be the projection of $\theta_{p}^{\circ k}\left(\mathfrak{z}\right)$
mod $p$. Then:
\begin{equation}
\prod_{k=0}^{n-1}r_{\left[\theta_{p}^{\circ k}\left(\mathfrak{z}\right)\right]_{p}}=\prod_{k=0}^{n-1}r_{j_{k}}
\end{equation}
We can compute this product in closed form by scrutinizing the $p$-adic
digits of $\mathfrak{z}$.
\begin{prop}
\label{prop:Product identity for kappa_X}Let $n\in\mathbb{N}_{0}$
and $\mathfrak{z}\in\mathbb{Z}_{p}$. Then:
\begin{equation}
\prod_{k=0}^{n-1}r_{\left[\theta_{p}^{\circ k}\left(\mathfrak{z}\right)\right]_{p}}=r_{0}^{n}\prod_{j=1}^{p-1}\left(r_{j}/r_{0}\right)^{\#_{j}\left(\left[\mathfrak{z}\right]_{p^{n}}\right)}
\end{equation}
\end{prop}
Proof: We begin by writing $\mathfrak{z}$ as a Hensel series:
\begin{equation}
\mathfrak{z}=\sum_{j=0}^{\infty}d_{j}p^{j}
\end{equation}
Here: 
\begin{equation}
\left[\theta_{p}^{\circ k}\left(\mathfrak{z}\right)\right]_{p}=\left[\sum_{j=0}^{\infty}d_{j+k}p^{j}\right]_{p}=d_{k}
\end{equation}
and so:
\begin{equation}
\prod_{k=0}^{n-1}r_{\left[\theta_{p}^{\circ k}\left(\mathfrak{z}\right)\right]_{p}}=\prod_{k=0}^{n-1}r_{d_{k}}
\end{equation}
Here, we note that for any $j\in\left\{ 1,\ldots,p-1\right\} $, the
number of $k$ in $\left\{ 0,\ldots,n-1\right\} $ for which $d_{k}=j$
is equal to $\#_{j}\left(\left[\mathfrak{z}\right]_{p^{n}}\right)$.
Thus:
\begin{equation}
\prod_{k=0}^{n-1}r_{\left[\theta_{p}^{\circ k}\left(\mathfrak{z}\right)\right]_{p}}=r_{0}^{\left|\left\{ k\in\left\{ 0,\ldots,n-1\right\} :d_{k}=0\right\} \right|}\prod_{j=1}^{p-1}r_{j}^{\#_{j}\left(\left[\mathfrak{z}\right]_{p^{n}}\right)}
\end{equation}

Now, let $k_{n}$ be the integer so that $n-k_{n}-1$ is the largest
index $\ell\in\left\{ 0,\ldots,n-1\right\} $ for which $d_{\ell}\neq0$.
That is:
\begin{align*}
d_{n-1}=d_{n-2}=\cdots=d_{n-k_{n}} & =0\\
d_{n-k_{n}-1} & \neq0
\end{align*}
Observe that $\lambda_{p}\left(\left[\mathfrak{z}\right]_{p^{n}}\right)=n-k_{n}$,
while:
\begin{align*}
\left|\left\{ k\in\left\{ 0,\ldots,n-1\right\} :d_{k}=0\right\} \right| & =k_{n}+\lambda_{p}\left(\left[\mathfrak{z}\right]_{p^{n}}\right)-\sum_{j=1}^{p-1}\#_{j}\left(\left[\mathfrak{z}\right]_{p^{n}}\right)\\
 & =n-\sum_{j=1}^{p-1}\#_{j}\left(\left[\mathfrak{z}\right]_{p^{n}}\right)
\end{align*}
 Thus:
\begin{equation}
\prod_{k=0}^{n-1}r_{\left[\theta_{p}^{\circ k}\left(\mathfrak{z}\right)\right]_{p}}=r_{0}^{n-\sum_{j=1}^{p-1}\#_{j}\left(\left[\mathfrak{z}\right]_{p^{n}}\right)}\prod_{j=1}^{p-1}r_{j}^{\#_{j}\left(\left[\mathfrak{z}\right]_{p^{n}}\right)}=r_{0}^{n}\prod_{j=1}^{p-1}\left(r_{j}/r_{0}\right)^{\#_{j}\left(\left[\mathfrak{z}\right]_{p^{n}}\right)}
\end{equation}

Q.E.D.

\vphantom{}We can characterize expressions of the form (\ref{eq:prototypical M-function})
as follows:
\begin{prop}
\label{prop:M function characterization}Let $R$ be an integral domain,
and let $\left\{ M_{n}\right\} _{n\geq0}$ be a sequence in $\mathcal{S}\left(\mathbb{Z}_{p},R\right)$.
Then, the following are equivalent:

I. There exist $r_{j}\in R\backslash\left\{ 0\right\} $ for $j\in\left\{ 0,\ldots,p-1\right\} $
so that:
\begin{equation}
M_{n}\left(\mathfrak{z}\right)=\prod_{k=0}^{n-1}r_{\left[\theta_{p}^{\circ k}\left(\mathfrak{z}\right)\right]_{p}}
\end{equation}

II. The following two properties hold:

i. For each $n$, $M_{n}\in\mathcal{S}\left(\mathbb{Z}_{p},R\right)$
is locally constant mod $p^{n}$ and has no zeroes.

ii. 
\begin{equation}
M_{n}\left(\theta_{p}^{\circ m}\left(\mathfrak{z}\right)\right)=\frac{M_{m+n}\left(\mathfrak{z}\right)}{M_{m}\left(\mathfrak{z}\right)},\textrm{ }\forall\mathfrak{z}\in\mathbb{Z}_{p},\textrm{ }\forall m,n\in\mathbb{N}_{0}\label{eq:kappa shift equation}
\end{equation}
where both sides are contained in $R$.

III. There exist $r_{j}\in R\backslash\left\{ 0\right\} $ for $j\in\left\{ 0,\ldots,p-1\right\} $
so that:

\begin{equation}
M_{n}\left(\mathfrak{z}\right)=r_{\left[\mathfrak{z}\right]_{p}}M_{n-1}\left(\theta_{p}\left(\mathfrak{z}\right)\right),\textrm{ }\forall\mathfrak{z}\in\mathbb{Z}_{p},\textrm{ }\forall n\geq1
\end{equation}
\end{prop}
\begin{rem}
The integral domain condition on $R$ is needed to ensure that $r_{0},\ldots,r_{p-1}\in R\backslash\left\{ 0\right\} $
guarantees that $M_{n}\left(\mathfrak{z}\right)\neq0$ for all $\mathfrak{z}\in\mathbb{Z}_{p}$.
\end{rem}
Proof:
\textbullet{} (II $\Rightarrow$ I) Suppose (II) holds. For each $j\in\left\{ 0,\ldots,p-1\right\} $,
let $r_{j}$ denote $M_{1}\left(j\right)$. Then, for $j\in\left\{ 0,\ldots,p-1\right\} $
and $\mathfrak{z}\in\mathbb{Z}_{p}$, note that $\theta_{p}\left(p\mathfrak{z}+j\right)=\mathfrak{z}$.
Thus (ii) for $m=1$ gives:
\begin{equation}
M_{n}\left(\mathfrak{z}\right)=\frac{M_{n+1}\left(p\mathfrak{z}+j\right)}{M_{1}\left(p\mathfrak{z}+j\right)}
\end{equation}
and so:
\begin{equation}
M_{n+1}\left(p\mathfrak{z}+j\right)=M_{1}\left(p\mathfrak{z}+j\right)M_{n}\left(\mathfrak{z}\right),\textrm{ }\forall j\in\left\{ 0,\ldots,p-1\right\} ,\textrm{ }\forall\mathfrak{z}\in\mathbb{Z}_{p},\textrm{ }\forall n\geq0
\end{equation}
Since $M_{1}$ is locally constant mod $p$, we have $M_{1}\left(p\mathfrak{z}+j\right)=M_{1}\left(j\right)=r_{j}$,
and thus, that:
\begin{equation}
M_{n+1}\left(p\mathfrak{z}+j\right)=r_{j}M_{n}\left(\mathfrak{z}\right),\textrm{ }\forall j\in\left\{ 0,\ldots,p-1\right\} ,\textrm{ }\forall\mathfrak{z}\in\mathbb{Z}_{p},\textrm{ }\forall n\geq0
\end{equation}
Setting $\mathfrak{y}=p\mathfrak{z}+j$, we get:
\begin{equation}
M_{n+1}\left(\mathfrak{y}\right)=r_{\left[\mathfrak{y}\right]_{p}}M_{n}\left(\theta_{p}\left(\mathfrak{y}\right)\right),\textrm{ }\forall\mathfrak{y}\in\mathbb{Z}_{p},\textrm{ }\forall n\geq0
\end{equation}
Setting $n=\mathfrak{y}=0$, we get:
\begin{equation}
\underbrace{M_{1}\left(0\right)}_{r_{0}}=r_{0}M_{0}\left(0\right)
\end{equation}
Since $r_{0}\neq0$, this forces $M_{0}\left(0\right)=1$. Since $M_{0}$
is locally constant mod $p^{0}$, it is a constant function, and thus
$M_{0}=1$. Consequently, by induction:
\begin{align*}
M_{n}\left(\mathfrak{y}\right) & =r_{\left[\mathfrak{y}\right]_{p}}M_{n-1}\left(\theta_{p}\left(\mathfrak{y}\right)\right)\\
 & =r_{\left[\mathfrak{y}\right]_{p}}r_{\left[\theta_{p}\left(\mathfrak{y}\right)\right]_{p}}M_{n-2}\left(\theta_{p}^{\circ2}\left(\mathfrak{y}\right)\right)\\
 & \vdots\\
 & =\underbrace{M_{0}\left(\theta_{p}^{\circ n}\left(\mathfrak{y}\right)\right)}_{1}\prod_{k=0}^{n-1}r_{\left[\theta_{p}^{\circ k}\left(\mathfrak{y}\right)\right]_{p}}\\
 & =\prod_{k=0}^{n-1}r_{\left[\theta_{p}^{\circ k}\left(\mathfrak{y}\right)\right]_{p}}
\end{align*}
Thus, (I) is satisfied. $\checkmark$

\textbullet{} (I $\Rightarrow$ II) Suppose (I) holds. Then, by \textbf{Proposition
\ref{prop:Product identity for kappa_X}}, we have:
\begin{equation}
M_{n}\left(\mathfrak{z}\right)=r_{0}^{n}\prod_{j=1}^{p-1}\left(r_{j}/r_{0}\right)^{\#_{j}\left(\left[\mathfrak{z}\right]_{p^{n}}\right)}
\end{equation}
This shows that $M_{n}$ is non-vanishing and locally constant mod
$p^{n}$, proving (i). As for (ii), \textbf{Proposition \ref{prop:fundamental functional equations-1}}
can be written as:
\begin{equation}
\#_{j}\left(\left[\mathfrak{z}\right]_{p^{n}}\right)=\begin{cases}
0 & \textrm{if }n=0\\
\#_{j}\left(\left[\theta_{p}\left(\mathfrak{z}\right)\right]_{p^{n-1}}\right)+\left[\left[\mathfrak{z}\right]_{p}=j\right] & \textrm{if }n\geq1
\end{cases},\textrm{ }\forall\mathfrak{z}\in\mathbb{Z}_{p},\textrm{ }\forall j\in\left\{ 0,\ldots,p-1\right\} 
\end{equation}
Hence, by induction:
\begin{equation}
\#_{j}\left(\left[\theta_{p}^{\circ m}\left(\mathfrak{z}\right)\right]_{p^{n-m}}\right)=\#_{j}\left(\left[\mathfrak{z}\right]_{p^{n}}\right)-\sum_{k=0}^{m-1}\left[\left[\theta_{p}^{\circ k}\left(\mathfrak{z}\right)\right]_{p}=j\right],\textrm{ }\forall m\in\left\{ 0,\ldots,n\right\} 
\end{equation}
So, adding $m$ to $n$ gives:
\begin{align*}
\#_{j}\left(\left[\theta_{p}^{\circ m}\left(\mathfrak{z}\right)\right]_{p^{n}}\right) & =\#_{j}\left(\left[\mathfrak{z}\right]_{p^{n+m}}\right)-\sum_{k=0}^{m-1}\left[\left[\theta_{p}^{\circ k}\left(\mathfrak{z}\right)\right]_{p}=j\right],\textrm{ }\forall m\in\left\{ 0,\ldots,n+m\right\} \\
 & =\#_{j}\left(\left[\mathfrak{z}\right]_{p^{n+m}}\right)-\#_{p:j}\left(\left[\mathfrak{z}\right]_{p^{m}}\right),\textrm{ }\forall m\in\left\{ 0,\ldots,n+m\right\} 
\end{align*}
Hence:
\begin{align*}
M_{n}\left(\theta_{p}^{\circ m}\left(\mathfrak{z}\right)\right) & =r_{0}^{n}\prod_{j=1}^{p-1}\left(r_{j}/r_{0}\right)^{\#_{j}\left(\left[\theta_{p}^{\circ m}\left(\mathfrak{z}\right)\right]_{p^{n}}\right)}\\
 & =r_{0}^{n}\prod_{j=1}^{p-1}\left(r_{j}/r_{0}\right)^{\#_{j}\left(\left[\mathfrak{z}\right]_{p^{n+m}}\right)-\#_{p:j}\left(\left[\mathfrak{z}\right]_{p^{m}}\right)}\\
 & =\frac{r_{0}^{n+m}\prod_{j=1}^{p-1}\left(r_{j}/r_{0}\right)^{\#_{j}\left(\left[\mathfrak{z}\right]_{p^{n+m}}\right)}}{r_{0}^{m}\prod_{j=1}^{p-1}\left(r_{j}/r_{0}\right)^{\#_{p:j}\left(\left[\mathfrak{z}\right]_{p^{m}}\right)}}\\
 & =\frac{M_{m+n}\left(\mathfrak{z}\right)}{M_{m}\left(\mathfrak{z}\right)}
\end{align*}
and so, (II) is satisfied.

\textbullet{} (II $\Leftrightarrow$ III) This equivalence is obvious,
and follows immediately from the arguments given in (II $\Rightarrow$
I).

Q.E.D.

\vphantom{}

This justifies the following definition:
\begin{defn}
Given an integral domain $R$ and an integer $p\geq2$, a \textbf{$p$-adic
$R$-valued M-function }is a sequence of functions $\left\{ M_{n}\right\} _{n\geq0}$
in $\mathcal{S}\left(\mathbb{Z}_{p},R\right)$ satisfying either of
the two equivalent conditions of \textbf{Proposition \ref{prop:M function characterization}}.
The constants $r_{0},\ldots,r_{p-1}\in R\backslash\left\{ 0\right\} $
given in \textbf{Proposition \ref{prop:M function characterization}}
are called the \textbf{multipliers }of the M-function. $r_{0}$ is
called the \textbf{zero multiplier}.

The \textbf{$\kappa$-factor }of an M-function $\left\{ M_{n}\right\} _{n\geq0}$
is the function $\kappa:\mathbb{N}_{0}\rightarrow\textrm{Frac}R$
defined by:
\begin{equation}
\kappa\left(n\right)\overset{\textrm{def}}{=}\prod_{j=1}^{p-1}\left(\frac{r_{j}}{r_{0}}\right)^{\#_{p:j}\left(n\right)}
\end{equation}
where the $r_{j}$s are the multipliers of the M-function.

Because M-functions are uniquely determined by their zero multipliers
and $\kappa$-factors, we will sometimes give M-functions as a pair
$\left(r,\kappa\right)$, where $\kappa$ is a $\kappa$-factor, and
$r$ is the zero factor.

Furthermore, we write $\textrm{MFunc}\left(\mathbb{Z}_{p},R\right)$
to denote the \textbf{set of all $R$-valued $p$-adic M-functions}.
An element of $\textrm{MFunc}\left(\mathbb{Z}_{p},R\right)$ is written
as a sequence $\left\{ M_{n}\right\} _{n\geq0}$.
\end{defn}
\begin{rem}
\textbf{Proposition \ref{prop:M function characterization}} then
implies that every M-function can be uniquely written as:
\begin{equation}
M_{n}\left(\mathfrak{z}\right)=r_{0}^{n}\kappa\left(\left[\mathfrak{z}\right]_{p^{n}}\right),\textrm{ }\forall n\geq0,\textrm{ }\forall\mathfrak{z}\in\mathbb{Z}_{p}
\end{equation}
where $r_{0}$ is the zero multiplier and $\kappa$ is the $\kappa$-factor.
\end{rem}
\begin{prop}
\label{prop:Kappa shift equation}Let $\kappa$ be the $\kappa$-factor
of an M-function. Then:
\begin{equation}
\kappa\left(\left[\theta_{p}^{\circ m}\left(\mathfrak{z}\right)\right]_{p^{n}}\right)=\frac{\kappa\left(\left[\mathfrak{z}\right]_{p^{m+n}}\right)}{\kappa\left(\left[\mathfrak{z}\right]_{p^{m}}\right)}
\end{equation}
\end{prop}
Proof: Equivalent to the proof of (\ref{eq:kappa shift equation}).
Q.E.D.
\begin{rem}
Equation (\ref{eq:kappa shift equation}) from \textbf{Proposition
\ref{prop:Kappa shift equation}} is incredibly important, and will
be used repeatedly in nearly everything that follows.
\end{rem}

\subsection{\label{subsec:The-Central-Computation}The Central Computation}

We fix once and for all integers $d\geq1$, $p\geq2$, and a global
field $K$. If $K$ has positive characteristic, we require $\textrm{char}K$
to be co-prime to $p$.

Now we introduce the notation to be used for the rest of this section
and for a good deal of the rest of the paper.
\begin{notation}
For each $j\in\left\{ 1,\ldots,d\right\} $ and each $k\in\left\{ 0,\ldots,p-1\right\} $,
let $a_{j,k}$ and $b_{j,k}$ be indeterminates, completely unrelated
to one another. Fix a $K$-algebra $\mathcal{A}$ containing the $a_{j,k}$s
and $b_{j,k}$s so that, for each $j\in\left\{ 1,\ldots,d\right\} $,
there is a degree $1$ F-series $X_{j}:\mathbb{Z}_{p}\rightarrow\mathcal{A}$
satisfying:
\begin{equation}
X_{j}\left(p\mathfrak{z}+k\right)=a_{j,k}X_{j}\left(\mathfrak{z}\right)+b_{j,k},\textrm{ }\forall\mathfrak{z}\in\mathbb{Z}_{p},\textrm{ }\forall k\in\left\{ 0,\ldots,p-1\right\} \label{eq:X_j functional equations}
\end{equation}
Finally, given any $d$-tuple $\mathbf{n}=\left(n_{1},\ldots,n_{d}\right)\in\mathbb{N}_{0}^{d}$
of non-negative integers, we write:
\begin{equation}
X_{\mathbf{n}}\left(\mathfrak{z}\right)\overset{\textrm{def}}{=}\prod_{j=1}^{d}X_{j}^{n_{j}}\left(\mathfrak{z}\right)
\end{equation}
where $X_{j}^{n_{j}}$ is defined to be the constant function $1$
for any $j$ for which $n_{j}=0$. We then write:
\begin{equation}
a_{\mathbf{n},k}\overset{\textrm{def}}{=}\prod_{j=1}^{d}a_{j,k}^{n_{j}}
\end{equation}
\begin{equation}
b_{\mathbf{n},k}\overset{\textrm{def}}{=}\prod_{j=1}^{d}b_{j,k}^{n_{j}}
\end{equation}
We write $\mathbf{m}\leq\mathbf{n}$ to mean the $\mathbf{m}=\left(m_{1},\ldots,m_{d}\right)$
is a $d$-tuple of non-negative integers so that $m_{j}\leq n_{j}$
for all $j$; we write $\mathbf{m}<\mathbf{n}$ to mean $\mathbf{m}\leq\mathbf{n}$
and there exists at least one $j$ so that $m_{j}<n_{j}$. We also
write:
\begin{equation}
\binom{\mathbf{n}}{\mathbf{m}}\overset{\textrm{def}}{=}\prod_{j=1}^{d}\binom{n_{j}}{m_{j}}
\end{equation}
and:
\begin{equation}
\Sigma\left(\mathbf{n}\right)\overset{\textrm{def}}{=}\sum_{j=1}^{d}n_{j}
\end{equation}
\end{notation}
\begin{rem}
$\mathcal{A}$ is our stand-in for whatever the correct choice of
codomain will turn out to be for our F-series.
\end{rem}
Our goal is to show that if each $X_{j}$ has a Fourier transform,
then so does $X_{\mathbf{n}}$ for all $\mathbf{n}\in\mathbb{N}_{0}^{d}$.
To do this, we begin by establishing $X_{\mathbf{n}}$'s functional
equations.
\begin{prop}
\label{prop:Z functional equation}Let $\mathbf{n}\in\mathbb{N}_{0}^{d}$.
Then, $X_{\mathbf{n}}$ satisfies functional equations:
\begin{equation}
X_{\mathbf{n}}\left(p\mathfrak{z}+k\right)=\sum_{\mathbf{m}\leq\mathbf{n}}r_{\mathbf{m},\mathbf{n},k}X_{\mathbf{m}}\left(\mathfrak{z}\right),\textrm{ }\forall\mathfrak{z}\in\mathbb{Z}_{p},\textrm{ }\forall k\in\left\{ 0,\ldots,p-1\right\} \label{eq:Z functional equation}
\end{equation}
where:
\begin{equation}
r_{\mathbf{m},\mathbf{n},k}\overset{\textrm{def}}{=}\binom{\mathbf{n}}{\mathbf{m}}a_{\mathbf{m},k}b_{\mathbf{n}-\mathbf{m},k}\label{eq:def of r e m k}
\end{equation}
Using the $p$-adic shift operator, (\ref{eq:Z functional equation})
can be written as:
\begin{equation}
X_{\mathbf{n}}\left(\mathfrak{z}\right)=\sum_{\mathbf{m}\leq\mathbf{n}}r_{\mathbf{n},\mathbf{m},\left[\mathfrak{z}\right]_{p}}X_{\mathbf{m}}\left(\theta_{p}\left(\mathfrak{z}\right)\right),\textrm{ }\forall\mathfrak{z}\in\mathbb{Z}_{p}\label{prop:shift reformulation}
\end{equation}
\end{prop}
\begin{rem}
Here, $\mathbf{n}-\mathbf{m}$ denotes the $d$-tuple $\left(n_{1}-m_{1},\ldots,n-m_{d}\right)$.
\end{rem}
Proof: Letting $\mathfrak{z}\in\mathbb{Z}_{p}$ and $k\in\left\{ 0,\ldots,p-1\right\} $,
we have:
\begin{align*}
X_{\mathbf{n}}\left(p\mathfrak{z}+k\right) & =\prod_{j=1}^{d}\left(X_{j}\left(p\mathfrak{z}+k\right)\right)^{n_{j}}\\
\left(\textrm{use }(\ref{eq:X_j functional equations})\right); & =\prod_{j=1}^{d}\left(a_{j,k}X_{j}\left(\mathfrak{z}\right)+b_{j,k}\right)^{n_{j}}\\
\left(\textrm{Binomial theorem}\right); & =\prod_{j=1}^{d}\sum_{m=0}^{n_{j}}\binom{n_{j}}{m}a_{j,k}^{m}b_{j,k}^{n_{j}-m}X_{j}^{m}\left(\mathfrak{z}\right)\\
 & =\sum_{\mathbf{m}\leq\mathbf{n}}\binom{\mathbf{n}}{\mathbf{m}}a_{\mathbf{m},k}b_{\mathbf{n}-\mathbf{m},k}X_{\mathbf{m}}\left(\mathfrak{z}\right)
\end{align*}

Q.E.D.

\vphantom{}

Next, we investigate Fourier transforms. For this, we need some more
notation:
\begin{defn}
Given $\mathbf{m},\mathbf{n}\in\mathbb{N}_{0}^{d}$ with $\mathbf{m}\leq\mathbf{n}$,
we write $\alpha_{\mathbf{m},\mathbf{n}}:\hat{\mathbb{Z}}_{p}\rightarrow\mathcal{A}\otimes_{K}K\left(\zeta_{p^{\infty}}\right)$
to denote the function:
\begin{equation}
\alpha_{\mathbf{m},\mathbf{n}}\left(t\right)\overset{\textrm{def}}{=}\frac{1}{p}\sum_{k=0}^{p-1}r_{\mathbf{m},\mathbf{n},k}e^{-2\pi ikt},\textrm{ }\forall t\in\hat{\mathbb{Z}}_{p}
\end{equation}
We adopt a shorthand and write $\alpha_{\mathbf{n}}$ to denote $\alpha_{\mathbf{n},\mathbf{n}}$,
and write $r_{\mathbf{n},k}$ to denote $r_{\mathbf{n},\mathbf{n},k}$
, so that:
\begin{equation}
r_{\mathbf{n},k}=\binom{\mathbf{n}}{\mathbf{n}}a_{\mathbf{n},k}b_{\mathbf{n}-\mathbf{n},k}=a_{\mathbf{n},k}=\prod_{j=1}^{d}a_{j,k}^{n_{j}}
\end{equation}
and:
\begin{equation}
\alpha_{\mathbf{n}}\left(t\right)=\frac{1}{p}\sum_{k=0}^{p-1}r_{\mathbf{n},k}e^{-2\pi ikt}
\end{equation}
\end{defn}
\begin{rem}
$\alpha_{\mathbf{m},\mathbf{n}}$ takes values in $\mathcal{A}\otimes_{K}K\left(\zeta_{p^{\infty}}\right)$,
rather than $\mathcal{A}$ because of those roots of unity given by
the $e^{-2\pi ikt}$s.
\end{rem}
\begin{example}
\label{exa:proceeding formally}Much of what we will do is motivated
by two desires: to make formal computations rigorous, and to understand
how to compensate for when those formal computations break down.

Proceeding formally, let's suppose that, for all $\mathbf{n}$, we
can realize $X_{\mathbf{n}}$ as a measure $X_{\mathbf{n}}\left(\mathfrak{z}\right)d\mathfrak{z}$.
Its Fourier-Stieltjes transform would then be given by:
\begin{equation}
\hat{X}_{\mathbf{n}}\left(\mathfrak{z}\right)=\int_{\mathbb{Z}_{p}}e^{-2\pi i\left\{ t\mathfrak{z}\right\} _{p}}X_{\mathbf{n}}\left(\mathfrak{z}\right)d\mathfrak{z}
\end{equation}
Using the change of variable formula for $p$-adic integration:
\begin{equation}
\int_{\mathfrak{a}\mathbb{Z}_{p}+\mathfrak{b}}f\left(\mathfrak{z}\right)d\mathfrak{z}=\left|\mathfrak{a}\right|_{p}\int_{\mathbb{Z}_{p}}f\left(\mathfrak{a}\mathfrak{z}+\mathfrak{b}\right)d\mathfrak{z}
\end{equation}
for all $\mathfrak{a},\mathfrak{b}\in\mathbb{Z}_{p}$ with $\mathfrak{a}\neq0$
and all nice $f$, along with the partition of unity:
\begin{equation}
1=\sum_{k=0}^{p-1}\left[\mathfrak{z}\overset{p}{\equiv}k\right],\textrm{ }\forall\mathfrak{z}\in\mathbb{Z}_{p}
\end{equation}
and the fact that:
\begin{equation}
\int_{\mathbb{Z}_{p}}\left[\mathfrak{z}\overset{p}{\equiv}k\right]f\left(\mathfrak{z}\right)d\mathfrak{z}=\int_{p\mathbb{Z}_{p}+k}f\left(\mathfrak{z}\right)d\mathfrak{z}
\end{equation}
and:
\begin{equation}
\int_{\mathbb{Z}_{p}}d\mathfrak{z}=1
\end{equation}
we have, for any $\mathbf{n}$:
\begin{align*}
\hat{X}_{\mathbf{n}}\left(t\right) & =\int_{\mathbb{Z}_{p}}e^{-2\pi i\left\{ t\mathfrak{z}\right\} _{p}}X_{\mathbf{n}}\left(\mathfrak{z}\right)d\mathfrak{z}\\
 & =\sum_{k=0}^{p-1}\int_{\mathbb{Z}_{p}}\left[\mathfrak{z}\overset{p}{\equiv}k\right]e^{-2\pi i\left\{ t\mathfrak{z}\right\} _{p}}X_{\mathbf{n}}\left(\mathfrak{z}\right)d\mathfrak{z}\\
 & =\sum_{k=0}^{p-1}\int_{p\mathbb{Z}_{p}+k}e^{-2\pi i\left\{ t\mathfrak{z}\right\} _{p}}X_{n}\left(\mathfrak{z}\right)d\mathfrak{z}\\
\left(\textrm{change of variable}\right); & =\frac{1}{p}\sum_{k=0}^{p-1}\int_{\mathbb{Z}_{p}}e^{-2\pi i\left\{ t\left(p\mathfrak{z}+k\right)\right\} _{p}}X_{\mathbf{n}}\left(p\mathfrak{z}+k\right)d\mathfrak{z}\\
\left(\textrm{apply }(\ref{eq:Z functional equation})\right); & =\frac{1}{p}\sum_{k=0}^{p-1}\int_{\mathbb{Z}_{p}}e^{-2\pi i\left\{ pt\mathfrak{z}\right\} _{p}}e^{-2\pi ikt}\left(\sum_{\mathbf{m}\leq\mathbf{n}}r_{\mathbf{m},\mathbf{n},k}X_{\mathbf{m}}\left(\mathfrak{z}\right)\right)d\mathfrak{z}\\
 & =\sum_{\mathbf{m}\leq\mathbf{e}}\underbrace{\left(\frac{1}{p}\sum_{k=0}^{p-1}r_{\mathbf{m},\mathbf{n},k}e^{-2\pi ikt}\right)}_{\alpha_{\mathbf{m},\mathbf{n}}\left(t\right)}\underbrace{\int_{\mathbb{Z}_{p}}e^{-2\pi i\left\{ pt\mathfrak{z}\right\} _{p}}X_{\mathbf{m}}\left(\mathfrak{z}\right)d\mathfrak{z}}_{\hat{X}_{\mathbf{m}}\left(pt\right)}\\
 & =\sum_{\mathbf{m}\leq\mathbf{n}}\alpha_{\mathbf{m},\mathbf{n}}\left(t\right)\hat{X}_{\mathbf{m}}\left(pt\right)
\end{align*}
This suggests the functional equation:
\begin{equation}
\hat{X}_{\mathbf{n}}\left(t\right)=\sum_{\mathbf{m}\leq\mathbf{n}}\alpha_{\mathbf{m},\mathbf{n}}\left(t\right)\hat{X}_{\mathbf{m}}\left(pt\right)\label{eq:formal fourier functional equation}
\end{equation}
ought to be a relation between the Fourier transforms of the $X_{\mathbf{n}}$.
Setting $t=0$, pulling out the $\mathbf{m}=\mathbf{n}$ term gives
us:
\begin{equation}
\hat{X}_{\mathbf{n}}\left(0\right)=\alpha_{\mathbf{n}}\left(0\right)\hat{X}_{\mathbf{n}}\left(0\right)+\sum_{\mathbf{m}<\mathbf{n}}\alpha_{\mathbf{m},\mathbf{n}}\left(0\right)\hat{X}_{\mathbf{m}}\left(0\right)
\end{equation}
If $\alpha_{\mathbf{n}}\left(0\right)\neq1$, we can then solve for
$\hat{X}_{\mathbf{m}}\left(0\right)$ in terms of $\alpha_{\mathbf{m},\mathbf{n}}$s
and $\hat{X}_{\mathbf{m}}$ for all $\mathbf{m}<\mathbf{n}$. Moreover,
given any $t\in\hat{\mathbb{Z}}_{p}\backslash\left\{ 0\right\} $,
(\ref{eq:formal fourier functional equation}) shows that the value
of $\hat{X}_{\mathbf{n}}$ at $t$ depends on the values of $\hat{X}_{\mathbf{m}}$
at $pt$ for all $\mathbf{m}\leq\mathbf{n}$. Since $\left|pt\right|_{p}=\left|t\right|_{p}/p$,
it then follows that we can recursively solve (\ref{eq:formal fourier functional equation})
for $\hat{X}_{\mathbf{n}}\left(t\right)$ for all $t\in\hat{\mathbb{Z}}_{p}\backslash\left\{ 0\right\} $.
\end{example}
This motivates two key questions:
\begin{question}
\label{que:key question (easy part)}Fix $\mathbf{n}$, and suppose
the functions $\hat{X}_{\mathbf{m}}$ are known for all $\mathbf{m}<\mathbf{n}$,
and that they are Fourier transforms of $\hat{X}_{\mathbf{m}}$ in
the sense that:
\begin{equation}
\lim_{N\rightarrow\infty}\sum_{\left|t\right|_{p}\leq p^{N}}\hat{X}_{\mathbf{m}}\left(t\right)e^{2\pi i\left\{ t\mathfrak{z}\right\} _{p}}=X_{\mathbf{m}}\left(\mathfrak{z}\right)
\end{equation}
Then, when does there exist $\hat{\chi}:\hat{\mathbb{Z}}_{p}\rightarrow\mathcal{A}\otimes_{K}K\left(\zeta_{p^{\infty}}\right)$
satisfying:
\begin{equation}
\hat{\chi}\left(t\right)=\alpha_{\mathbf{n}}\left(t\right)\hat{\chi}\left(pt\right)+\sum_{\mathbf{m}<\mathbf{n}}\alpha_{\mathbf{m},\mathbf{n}}\left(t\right)\hat{X}_{\mathbf{m}}\left(pt\right),\textrm{ }\forall t\in\hat{\mathbb{Z}}_{p}\textrm{ ?}\label{eq:formal solution}
\end{equation}
Moreover, if $\hat{\chi}$ exists, is it unique?
\end{question}
\begin{question}
\label{que:key question (hard part)}Let everything be as given in
\textbf{Question \ref{que:key question (easy part)}}. What, if any,
solutions $\hat{\chi}$ of (\ref{eq:formal solution}) are Fourier
transforms of $X_{\mathbf{n}}$ in the sense that:
\begin{equation}
\lim_{N\rightarrow\infty}\sum_{\left|t\right|_{p}\leq p^{N}}\hat{\chi}\left(t\right)e^{2\pi i\left\{ t\mathfrak{z}\right\} _{p}}=X_{\mathbf{n}}\left(\mathfrak{z}\right)
\end{equation}
\end{question}
\textbf{Proposition} \textbf{\ref{prop:X_3-hat formal solution}},
given below, answers\textbf{ Question \ref{que:key question (easy part)}}
in full.
\begin{prop}
\label{prop:X_3-hat formal solution}If:

I. $\hat{X}_{\mathbf{m}}$ is given for all $\mathbf{m}<\mathbf{n}$;

II. $\alpha_{\mathbf{n}}\left(0\right)\neq1$;

then (\ref{eq:formal solution}) has a unique solution $\hat{\chi}:\hat{\mathbb{Z}}_{p}\rightarrow\mathcal{A}\otimes_{K}K\left(\zeta_{p^{\infty}}\right)$
given by:
\begin{align}
\hat{\chi}\left(t\right) & =\frac{\sum_{\mathbf{m}<\mathbf{n}}\alpha_{\mathbf{m},\mathbf{n}}\left(0\right)\hat{X}_{\mathbf{m}}\left(0\right)}{1-\alpha_{\mathbf{n}}\left(0\right)}\prod_{n=0}^{-v_{p}\left(t\right)-1}\alpha_{\mathbf{n}}\left(p^{n}t\right)\\
 & +\sum_{\mathbf{m}<\mathbf{n}}\sum_{n=0}^{-v_{p}\left(t\right)-1}\left(\prod_{m=0}^{n-1}\alpha_{\mathbf{n}}\left(p^{m}t\right)\right)\alpha_{\mathbf{m},\mathbf{n}}\left(p^{n}t\right)\hat{X}_{\mathbf{m}}\left(p^{n+1}t\right)\nonumber 
\end{align}

On the other hand, if $\alpha_{\mathbf{n}}\left(0\right)=1$, there
are two possibilities:

\vphantom{}

I. $\sum_{\mathbf{m}<\mathbf{n}}\alpha_{\mathbf{m},\mathbf{n}}\left(0\right)\hat{X}_{\mathbf{m}}\left(0\right)\neq0$,
in which case (\ref{eq:formal solution}) has no solutions.

\vphantom{}

II. $\sum_{\mathbf{m}<\mathbf{n}}\alpha_{\mathbf{m},\mathbf{n}}\left(0\right)\hat{X}_{\mathbf{m}}\left(0\right)=0$,
in which case (\ref{eq:formal solution}) has infinitely many solutions.
In particular, the set of $\hat{\chi}:\hat{\mathbb{Z}}_{p}\rightarrow\mathcal{A}\otimes_{K}K\left(\zeta_{p^{\infty}}\right)$
satisfying (\ref{eq:formal solution}) is a $1$-dimensional affine
$K\left(\zeta_{p^{\infty}}\right)$-vector subspace of the space of
functions $\hat{\mathbb{Z}}_{p}\rightarrow\mathcal{A}\otimes_{K}K\left(\zeta_{p^{\infty}}\right)$.
\end{prop}
Proof: First, suppose $\alpha_{\mathbf{n}}\left(0\right)\neq1$. Then:
\begin{align*}
\hat{\chi}\left(t\right) & =\alpha_{\mathbf{n}}\left(t\right)\overbrace{\hat{\chi}\left(pt\right)}^{\textrm{expand}}+\sum_{\mathbf{m}<\mathbf{n}}\alpha_{\mathbf{m},\mathbf{n}}\left(t\right)\hat{X}_{\mathbf{m}}\left(pt\right)\\
 & =\alpha_{\mathbf{n}}\left(t\right)\left(\alpha_{\mathbf{n}}\left(pt\right)\hat{\chi}\left(p^{2}t\right)+\sum_{\mathbf{m}<\mathbf{n}}\alpha_{\mathbf{m},\mathbf{n}}\left(pt\right)\hat{X}_{\mathbf{m}}\left(p^{2}t\right)\right)+\sum_{\mathbf{m}<\mathbf{n}}\alpha_{\mathbf{m},\mathbf{n}}\left(t\right)\hat{X}_{\mathbf{m}}\left(pt\right)\\
 & \vdots\\
 & =\hat{\chi}\left(p^{N}t\right)\prod_{n=0}^{N-1}\alpha_{\mathbf{n}}\left(p^{n}t\right)+\sum_{\mathbf{m}<\mathbf{n}}\sum_{n=0}^{N-1}\left(\prod_{m=0}^{n-1}\alpha_{\mathbf{n}}\left(p^{m}t\right)\right)\alpha_{\mathbf{m},\mathbf{n}}\left(p^{n}t\right)\hat{X}_{\mathbf{m}}\left(p^{n+1}t\right)
\end{align*}
For $t\neq0$, choose $N=-v_{p}\left(t\right)$, so that $p^{N}t\overset{1}{\equiv}0$.
Then:
\begin{equation}
\hat{\chi}\left(t\right)=\hat{\chi}\left(0\right)\prod_{n=0}^{-v_{p}\left(t\right)-1}\alpha_{\mathbf{n}}\left(p^{n}t\right)+\sum_{\mathbf{m}<\mathbf{n}}\sum_{n=0}^{-v_{p}\left(t\right)-1}\left(\prod_{m=0}^{n-1}\alpha_{\mathbf{n}}\left(p^{m}t\right)\right)\alpha_{\mathbf{m},\mathbf{n}}\left(p^{n}t\right)\hat{X}_{\mathbf{m}}\left(p^{n+1}t\right)
\end{equation}
Using (\ref{eq:formal solution}), since $\alpha_{\mathbf{n}}\left(0\right)\neq1$,
we can set $t=0$ and solve for $\hat{\chi}\left(0\right)$ directly
to obtain:
\begin{equation}
\hat{\chi}\left(0\right)=\frac{\sum_{\mathbf{m}<\mathbf{n}}\alpha_{\mathbf{m},\mathbf{n}}\left(0\right)\hat{X}_{\mathbf{m}}\left(0\right)}{1-\alpha_{\mathbf{n}}\left(0\right)}
\end{equation}

Next, suppose $\alpha_{\mathbf{n}}\left(0\right)=1$. Then setting
$t=0$ in (\ref{eq:formal solution}) gives us:
\begin{equation}
\cancel{\hat{\chi}\left(0\right)}=\cancel{\hat{\chi}\left(0\right)}+\sum_{\mathbf{m}<\mathbf{n}}\alpha_{\mathbf{m},\mathbf{n}}\left(0\right)\hat{X}_{\mathbf{m}}\left(0\right),\textrm{ }\forall t\in\hat{\mathbb{Z}}_{p}\textrm{ }?
\end{equation}
and so, we see that the vanishing of the sum on the right is a necessary
condition for a solution to exist. So, supposing the sum on the right
vanishes, fix a constant $c\in\mathcal{A}\otimes_{K}K\left(\zeta_{p^{\infty}}\right)$,
and let us look for solutions $\hat{\chi}$ with $\hat{\chi}\left(0\right)=c$.
Then, for $\left|t\right|_{p}=p$, (\ref{eq:formal solution}) becomes:
\begin{equation}
\hat{\chi}\left(t\right)=c\alpha_{\mathbf{n}}\left(t\right)+\sum_{\mathbf{m}<\mathbf{n}}\alpha_{\mathbf{m},\mathbf{n}}\left(t\right)\hat{X}_{\mathbf{m}}\left(0\right)
\end{equation}
Since the $\hat{X}_{\mathbf{m}}$s are given functions, this then
determines the values of $\hat{\chi}\left(t\right)$ for all $\left|t\right|_{p}=p$.
Proceeding recursively, we can solve (\ref{eq:formal solution}) for
$\hat{\chi}\left(t\right)$ for all $t\in\hat{\mathbb{Z}}_{p}\backslash\left\{ 0\right\} $.

Now, picking $c^{\prime}\in\mathcal{A}\otimes_{K}K\left(\zeta_{p^{\infty}}\right)$,
let $\hat{\chi}^{\prime}$ be the solution of (\ref{eq:formal solution})
constructed from the initial condition $\hat{\chi}^{\prime}\left(0\right)=c^{\prime}$.
Then, writing $\hat{\mu}\left(t\right)$ to denote $\hat{\chi}\left(t\right)-\hat{\chi}^{\prime}\left(t\right)$,
the functional equations satisfied by $\hat{\chi}$ and $\hat{\chi}^{\prime}$
tell us that:
\begin{equation}
\hat{\mu}\left(t\right)=\alpha_{\mathbf{n}}\left(t\right)\hat{\mu}\left(pt\right),\textrm{ }\forall t\in\hat{\mathbb{Z}}_{p}
\end{equation}
By induction, for all $t\neq0$, it follows that:
\begin{align*}
\hat{\mu}\left(t\right) & =\alpha_{\mathbf{n}}\left(t\right)\hat{\mu}\left(pt\right)\\
 & =\alpha_{\mathbf{n}}\left(t\right)\alpha_{\mathbf{n}}\left(pt\right)\hat{\mu}\left(p^{2}t\right)\\
 & \vdots\\
\left(p^{-v_{p}\left(t\right)}t\overset{1}{\equiv}0\right); & =\hat{\mu}\left(0\right)\prod_{n=0}^{-v_{p}\left(t\right)-1}\alpha_{\mathbf{n}}\left(p^{n}t\right)
\end{align*}
which shows that the difference of any two solutions of (\ref{eq:formal solution})
lie in the $1$-dimensional $K\left(\zeta_{p^{\infty}}\right)$-vector
space generated by the function:
\begin{equation}
t\in\hat{\mathbb{Z}}_{p}\mapsto\prod_{n=0}^{-v_{p}\left(t\right)-1}\alpha_{\mathbf{n}}\left(p^{n}t\right)\in\mathcal{A}\otimes_{K}K\left(\zeta_{p^{\infty}}\right)
\end{equation}
with the right-hand side being defined as $1$ for $t=0$. This shows
that the solution space of (\ref{eq:formal solution}) is a $1$-dimensional
affine subspace of the vector space of functions $\hat{\mathbb{Z}}_{p}\rightarrow\mathcal{A}\otimes_{K}K\left(\zeta_{p^{\infty}}\right)$.

Q.E.D.

\vphantom{}

Answering \textbf{Question \ref{que:key question (hard part)}}, on
the other hand, is significantly more involved.
\begin{example}
\label{exa:J has one element}To start with, let's consider the simplest
possible case, where $\Sigma\left(\mathbf{n}\right)=1$. That is,
$X_{\mathbf{n}}$ is just a single F-series, characterized by the
equations:
\begin{equation}
X_{\mathbf{n}}\left(p\mathfrak{z}+k\right)=a_{\mathbf{n},k}X_{\mathbf{n}}\left(\mathfrak{z}\right)+b_{\mathbf{n},k},\textrm{ }\forall\mathfrak{z}\in\mathbb{Z}_{p},\textrm{ }\forall k\in\left\{ 0,\ldots,p-1\right\} 
\end{equation}
The only tuple $\mathbf{m}$ with $\mathbf{m}<\mathbf{n}$ is the
zero tuple, where $X_{\mathbf{0}}$, recall, is defined as the constant
function $1$. As such:
\begin{equation}
r_{\mathbf{m},\mathbf{n},k}=\begin{cases}
a_{\mathbf{n},k} & \textrm{if}\textrm{ }\mathbf{m}=\mathbf{n}\\
b_{\mathbf{n},k} & \textrm{if }\mathbf{m}=\mathbf{0}
\end{cases}
\end{equation}
Turning to (\ref{eq:formal solution}), our candidate $\hat{\chi}$
for $\hat{X}_{\mathbf{m}}$ must satisfy:
\begin{equation}
\hat{\chi}\left(t\right)=\alpha_{\mathbf{n}}\left(t\right)\hat{\chi}\left(pt\right)+\alpha_{\mathbf{0}}\left(t\right)\hat{X}_{\mathbf{0}}\left(pt\right)
\end{equation}
Here, we have $\hat{X}_{\mathbf{0}}\left(pt\right)=\mathbf{1}_{0}\left(pt\right)$,
which vanishes when $\left|t\right|_{p}\geq p^{2}$. Meanwhile:
\begin{align}
\alpha_{\mathbf{n}}\left(t\right) & =\frac{1}{p}\sum_{k=0}^{p-1}r_{\mathbf{n},k}e^{-2\pi ikt}=\frac{1}{p}\sum_{k=0}^{p-1}a_{\mathbf{n},k}e^{-2\pi ikt}\\
\nonumber \\
\alpha_{\mathbf{0}}\left(t\right) & =\frac{1}{p}\sum_{k=0}^{p-1}r_{\mathbf{0},\mathbf{n},k}e^{-2\pi ikt}=\frac{1}{p}\sum_{k=0}^{p-1}b_{\mathbf{n},k}e^{-2\pi ikt}
\end{align}
and so:
\begin{equation}
\hat{\chi}\left(t\right)=\left(\frac{1}{p}\sum_{k=0}^{p-1}a_{\mathbf{n},k}e^{-2\pi ikt}\right)\hat{\chi}\left(pt\right)+\left(\frac{1}{p}\sum_{k=0}^{p-1}b_{\mathbf{n},k}e^{-2\pi ikt}\right)\mathbf{1}_{0}\left(pt\right)
\end{equation}
When $\alpha_{\mathbf{n}}\left(0\right)=1$, setting $t=0$ gives
us:
\begin{equation}
\hat{\chi}\left(0\right)=\underbrace{\left(\frac{1}{p}\sum_{k=0}^{p-1}a_{\mathbf{n},k}\right)}_{1}\hat{\chi}\left(0\right)+\frac{1}{p}\sum_{k=0}^{p-1}b_{\mathbf{n},k}
\end{equation}
So, there will be infinitely many solutions of (\ref{eq:formal solution})
if and only if:
\begin{equation}
\frac{1}{p}\sum_{k=0}^{p-1}b_{\mathbf{n},k}=0\label{eq:average}
\end{equation}
and will be \emph{no solutions }otherwise. However, the formula for
$\hat{X}_{\mathbf{n}}$ given by \textbf{Theorem \ref{thm:The-breakdown-variety}}
works \emph{even when }(\ref{eq:average}) fails to hold.

This shows that the answer to \textbf{Question \ref{que:key question (hard part)}}
is, in general, a resounding \textbf{\emph{no}}.
\end{example}
Thankfully, we can deal with both of these cases by using the same
general approach I took in my dissertation, described below.
\begin{defn}
Given a function $f$ on $\mathbb{Z}_{p}$, and given an integer $N\geq0$,
we write $f_{N}$ to denote the \textbf{$N$th truncation of $f$},
defined by:
\begin{equation}
f_{N}\left(\mathfrak{z}\right)\overset{\textrm{def}}{=}f\left(\left[\mathfrak{z}\right]_{p^{N}}\right)
\end{equation}
As defined, $f_{N}$ is locally constant mod $p^{N}$. Because of
this, the Fourier transform of $f_{N}$ exists, and is given by:
\begin{equation}
\hat{f}_{N}\left(t\right)\overset{\textrm{def}}{=}\frac{\mathbf{1}_{0}\left(p^{N}t\right)}{p^{N}}\sum_{n=0}^{p^{N}-1}f\left(n\right)e^{-2\pi int}
\end{equation}
where $t\in\hat{\mathbb{Z}}_{p}$. Note that $\hat{f}_{N}\left(t\right)$
then vanishes for all $\left|t\right|_{p}>p^{N}$.

If $f$ itself has subscripts, the Fourier transform of $f$'s $N$th
truncation will be given by putting a hat on $f$ and appending an
$N$ as the right-most subscript of $f$. For example:
\begin{equation}
\hat{f}_{p_{\sqrt{7}},\&^{2},N}\left(t\right)
\end{equation}
denotes the Fourier transform of the $N$th truncation of the function
$f_{p_{\sqrt{7}},\&^{2},N}$.
\end{defn}
\begin{rem}
It is important to remember that even if $f$ itself might not have
a well-defined Fourier transform, the locally constant function $f_{N}$
obtained by truncating $f$ \emph{always} has a Fourier transform.
So, even though we might not yet know whether or not $\hat{X}_{\mathbf{n}}$
exists, we can unambiguously write $\hat{X}_{\mathbf{n},N}$, because
this is the Fourier transform of the locally constant function $X_{\mathbf{n},N}$.
\end{rem}
\begin{prop}
\label{prop:truncated Fourier recursion}Let $N\in\mathbb{N}_{0}$.
For all $\mathbf{m},\mathbf{n}\in\mathbb{N}_{0}^{d}$ with $\mathbf{m}\leq\mathbf{n}$:
\begin{equation}
\hat{X}_{\mathbf{n},N}\left(t\right)=\sum_{\mathbf{m}\leq\mathbf{n}}\alpha_{\mathbf{m},\mathbf{n}}\left(t\right)\hat{X}_{\mathbf{m},N-1}\left(pt\right),\textrm{ }\forall\left|t\right|_{p}\leq p^{N}\label{eq:truncated Fourier transform recursion for Z}
\end{equation}
Furthermore:
\begin{equation}
\hat{X}_{\mathbf{n},N}\left(t\right)=X_{\mathbf{n}}\left(0\right)\prod_{n=0}^{N-1}\alpha_{\mathbf{n}}\left(p^{n}t\right)+\sum_{\mathbf{m}\leq\mathbf{n}}\sum_{n=0}^{N-1}\left(\prod_{m=0}^{n-1}\alpha_{\mathbf{n}}\left(p^{m}t\right)\right)\alpha_{\mathbf{m},\mathbf{n}}\left(p^{n}t\right)\hat{X}_{\mathbf{m},N-1-n}\left(p^{n+1}t\right)
\end{equation}
for all $\left|t\right|_{p}\leq p^{N}$ and all $N\in\mathbb{N}_{1}$.
\end{prop}
Proof: We apply the functional equations of $X_{\mathbf{n}}$ to the
definition of $\hat{X}_{\mathbf{n},N}$. Letting $\left|t\right|_{p}\leq p^{N}$,
we have:

\begin{align*}
\hat{X}_{\mathbf{n},N}\left(t\right) & =\frac{1}{p^{N}}\sum_{n=0}^{p^{N}-1}X_{\mathbf{n}}\left(n\right)e^{-2\pi int}\\
 & =\frac{1}{p^{N}}\sum_{k=0}^{p-1}\sum_{n=0}^{p^{N-1}-1}X_{\mathbf{n}}\left(pn+k\right)e^{-2\pi i\left(pn+k\right)t}\\
 & =\frac{1}{p^{N}}\sum_{k=0}^{p-1}\sum_{n=0}^{p^{N-1}-1}\left(\sum_{\mathbf{m}\leq\mathbf{n}}r_{\mathbf{m},\mathbf{n},k}X_{\mathbf{m}}\left(\mathfrak{z}\right)e^{-2\pi i\left(pn+k\right)t}\right)\\
 & =\sum_{\mathbf{m}\leq\mathbf{n}}\underbrace{\frac{1}{p}\sum_{k=0}^{p-1}r_{\mathbf{m},\mathbf{n},k}e^{-2\pi ikt}}_{\alpha_{\mathbf{m},\mathbf{n}}\left(t\right)}\underbrace{\frac{1}{p^{N-1}}\sum_{n=0}^{p^{N-1}-1}X_{\mathbf{m}}\left(\mathfrak{z}\right)e^{-2\pi in\left(pt\right)}}_{\hat{X}_{\mathbf{m},N-1}\left(pt\right)}\\
 & =\sum_{\mathbf{m}\leq\mathbf{n}}\alpha_{\mathbf{m},\mathbf{n}}\left(t\right)\hat{X}_{\mathbf{m},N-1}\left(pt\right)
\end{align*}
Thus, in particular:
\begin{equation}
\hat{X}_{\mathbf{n},N}\left(t\right)=\alpha_{\mathbf{n}}\left(t\right)\hat{X}_{\mathbf{n},N-1}\left(pt\right)+\sum_{\mathbf{m}<\mathbf{n}}\alpha_{\mathbf{m},\mathbf{n}}\left(t\right)\hat{X}_{\mathbf{m},N-1}\left(pt\right)
\end{equation}
As such, we can use the left-hand side to substitute for $\hat{X}_{\mathbf{n},N-1}\left(pt\right)$
on the right-hand side: 
\begin{align*}
\hat{X}_{\mathbf{n},N}\left(t\right) & =\alpha_{\mathbf{n}}\left(t\right)\hat{X}_{\mathbf{n},N-1}\left(pt\right)+\sum_{\mathbf{m}<\mathbf{n}}\alpha_{\mathbf{m},\mathbf{n}}\left(t\right)\hat{X}_{\mathbf{m},N-1}\left(pt\right)\\
 & =\alpha_{\mathbf{n}}\left(t\right)\alpha_{\mathbf{n}}\left(pt\right)\hat{X}_{\mathbf{n},N-2}\left(p^{2}t\right)+\alpha_{\mathbf{n}}\left(t\right)\sum_{\mathbf{m}<\mathbf{n}}\alpha_{\mathbf{m},\mathbf{n}}\left(pt\right)\hat{X}_{\mathbf{m},N-2}\left(p^{2}t\right)\\
 & +\sum_{\mathbf{m}<\mathbf{n}}\alpha_{\mathbf{m},\mathbf{n}}\left(t\right)\hat{X}_{\mathbf{m},N-1}\left(pt\right)\\
 & =\cdots
\end{align*}
Proceeding to the $N$th step, we see that the pattern is:
\begin{equation}
\hat{X}_{\mathbf{n},N}\left(t\right)=\hat{X}_{\mathbf{n},0}\left(p^{N}t\right)\prod_{n=0}^{N-1}\alpha_{\mathbf{n}}\left(p^{n}t\right)+\sum_{\mathbf{m}<\mathbf{n}}\sum_{n=0}^{N-1}\left(\prod_{m=0}^{n-1}\alpha_{\mathbf{n}}\left(p^{m}t\right)\right)\alpha_{\mathbf{m},\mathbf{n}}\left(p^{n}t\right)\hat{X}_{\mathbf{m},N-1-n}\left(p^{n+1}t\right)
\end{equation}
as desired. Finally, note that: 
\begin{equation}
X_{\mathcal{\mathbf{n}},0}\left(\mathfrak{z}\right)=X_{\mathcal{\mathbf{n}}}\left(\left[\mathfrak{z}\right]_{p^{0}}\right)=X_{\mathbf{n}}\left(0\right)
\end{equation}
is a constant. As such, its Fourier transform is: 
\begin{equation}
\hat{X}_{\mathbf{n},0}\left(t\right)=X_{\mathcal{\mathbf{n}}}\left(0\right)\mathbf{1}_{0}\left(t\right)
\end{equation}
and so:
\begin{equation}
\hat{X}_{\mathbf{n},0}\left(p^{N}t\right)=X_{\mathbf{n}}\left(0\right)\mathbf{1}_{0}\left(p^{N}t\right)
\end{equation}
Since $\left|t\right|_{p}\leq p^{N}$, $p^{N}t\overset{1}{\equiv}0$,
which means:
\begin{equation}
\hat{X}_{\mathbf{n},0}\left(p^{N}t\right)=X_{\mathbf{n}}\left(0\right)\mathbf{1}_{0}\left(p^{N}t\right)=\hat{X}_{\mathbf{n},0}\left(p^{N}t\right)=X_{\mathbf{n}}\left(0\right)
\end{equation}
and so:
\begin{equation}
\hat{X}_{\mathbf{n},N}\left(t\right)=X_{\mathbf{n}}\left(0\right)\prod_{n=0}^{N-1}\alpha_{\mathbf{n}}\left(p^{n}t\right)+\sum_{\mathbf{m}<\mathbf{n}}\sum_{n=0}^{N-1}\left(\prod_{m=0}^{n-1}\alpha_{\mathbf{n}}\left(p^{m}t\right)\right)\alpha_{\mathbf{m},\mathbf{n}}\left(p^{n}t\right)\hat{X}_{\mathbf{m},N-1-n}\left(p^{n+1}t\right)
\end{equation}

Q.E.D.

\vphantom{}

For reference, if we were to apply this procedure to a single $p$-adic
F-series $X\left(\mathfrak{z}\right)$ characterized by the equations:
\begin{equation}
X\left(p\mathfrak{z}+j\right)=a_{j}X\left(\mathfrak{z}\right)+b_{j},\textrm{ }\forall j\in\left\{ 0,\ldots,p-1\right\} ,\textrm{ }\forall\mathfrak{z}\in\mathbb{Z}_{p}
\end{equation}
we would find that:
\begin{equation}
\hat{X}_{N}\left(t\right)=\alpha_{X}\left(t\right)\hat{X}_{N-1}\left(pt\right)+\beta_{X}\left(t\right)\mathbf{1}_{0}\left(pt\right)
\end{equation}
for certain functions $\alpha_{X}$ and $\beta_{X}$ of $t\in\hat{\mathbb{Z}}_{p}$.
Like with $\hat{X}_{3,N}$, this formula can be nested to obtain a
closed-form expression for $\hat{X}_{N}\left(t\right)$:
\begin{equation}
\hat{X}_{N}\left(t\right)=\sum_{n=0}^{N-1}\mathbf{1}_{0}\left(p^{n+1}t\right)\beta_{X}\left(p^{n}t\right)\prod_{m=0}^{n-1}\alpha_{X}\left(p^{m}t\right),\textrm{ }\forall N\geq1,\textrm{ }\forall t\in\hat{\mathbb{Z}}_{p}
\end{equation}
Observe that since $\alpha_{X}$ is a function on $\hat{\mathbb{Z}}_{p}=\mathbb{Z}\left[1/p\right]/\mathbb{Z}$,
it is $1$-periodic, and thus:
\begin{equation}
\alpha_{X}\left(p^{m}t\right)=\alpha_{X}\left(0\right),\textrm{ }\forall m\geq-v_{p}\left(t\right)
\end{equation}
Consequently, for all $t\in\hat{\mathbb{Z}}_{p}\backslash\left\{ 0\right\} $
and all $n\geq0$:
\begin{equation}
\prod_{m=0}^{n-1}\alpha_{X}\left(p^{m}t\right)=\left(\alpha_{X}\left(0\right)\right)^{\max\left\{ n+v_{p}\left(t\right),0\right\} }\prod_{m=0}^{-v_{p}\left(t\right)-1}\alpha_{X}\left(p^{m}t\right)
\end{equation}
In \cite{My Dissertation}, I used this to derive an exact formula
for $\hat{X}_{N}\left(t\right)$ for fixed $t$. This was of the form:
\begin{equation}
\hat{X}_{N}\left(t\right)=\hat{F}\left(t\right)+\hat{G}_{N}\left(t\right)\label{eq:drop the divergent terms}
\end{equation}
where $\hat{F},\hat{G}_{N}:\hat{\mathbb{Z}}_{p}\rightarrow\overline{\mathbb{Q}}$
were functions. $\hat{G}_{N}\left(t\right)$ can be thought of as
the ``divergent term'' which failed to converge to a limit as $N\rightarrow\infty$,
while $\hat{F}\left(t\right)$ was the fine structure hidden beneath
the divergent tumult. Miraculously, it turned out that the Fourier
series generated by $\hat{F}$ converged pointwise everywhere to $X$:
\begin{equation}
\lim_{N\rightarrow\infty}\sum_{\left|t\right|_{p}\leq p^{N}}\hat{F}\left(t\right)e^{2\pi i\left\{ t\mathfrak{z}\right\} _{p}}\overset{\mathcal{F}}{=}X\left(\mathfrak{z}\right),\textrm{ }\forall\mathfrak{z}\in\mathbb{Z}_{p}\label{eq:summation step}
\end{equation}
albeit with the caveat that the limit had to be evaluated with respect
to some frame $\mathcal{F}$. The convergence of (\ref{eq:summation step})
then justified calling $\hat{F}$ the Fourier transform of $X$, and
denoting it by $\hat{X}$.

In summary, the method is to obtain an asymptotic decomposition of
the form (\ref{eq:drop the divergent terms}), discard the divergent
terms, and then hope that leftovers generate a Fourier series that
sums to $X$. The most important step is the summation of the Fourier
series (\ref{eq:summation step}), where we verify that our candidate
for the Fourier transform of our function does, indeed, generate our
function's Fourier series.

So, beginning with the formula:

\begin{equation}
\hat{X}_{\mathbf{n},N}\left(t\right)=X_{\mathbf{n}}\left(0\right)\prod_{n=0}^{N-1}\alpha_{\mathbf{n}}\left(p^{n}t\right)+\sum_{\mathbf{m}<\mathbf{n}}\sum_{n=0}^{N-1}\left(\prod_{m=0}^{n-1}\alpha_{\mathbf{n}}\left(p^{m}t\right)\right)\alpha_{\mathbf{m},\mathbf{n}}\left(p^{n}t\right)\hat{X}_{\mathbf{m},N-1-n}\left(p^{n+1}t\right)\label{eq:X_3,N recursion}
\end{equation}
obtained in \textbf{Proposition \ref{prop:truncated Fourier recursion}},
we make the following educated guess as to what the Fourier transform
of $\hat{X}_{\mathbf{n}}$ might be.
\begin{defn}
Now, for all $\mathbf{m},\mathbf{n}\in\mathbb{N}_{0}^{d}$ with $\mathbf{m}\leq\mathbf{n}$,
set:
\begin{equation}
\hat{f}_{\mathbf{m},\mathbf{n}}\left(t\right)\overset{\textrm{def}}{=}\sum_{n=0}^{-v_{p}\left(t\right)-2}\left(\prod_{m=0}^{n-1}\alpha_{\mathbf{n}}\left(p^{m}t\right)\right)\alpha_{\mathbf{m},\mathbf{n}}\left(p^{n}t\right)\hat{X}_{\mathbf{m}}\left(p^{n+1}t\right)\label{eq:def of f_j}
\end{equation}
with $n$-sum being defined as $0$ when $\left|t\right|_{p}\leq p$.
Then, let:
\begin{equation}
\hat{f}_{\mathbf{n}}\left(t\right)\overset{\textrm{def}}{=}\sum_{\mathbf{m}<\mathbf{n}}\hat{f}_{\mathbf{m},\mathbf{n}}\left(t\right)
\end{equation}
$\hat{f}_{\mathbf{n}}$ is going to be our guess for the Fourier transform
of $X_{\mathbf{n}}$. To see how accurate the guess is, we'll need
to sum the Fourier series it generates, which we denote by:
\begin{equation}
\tilde{f}_{\mathbf{m},N}\left(\mathfrak{z}\right)\overset{\textrm{def}}{=}\sum_{\left|t\right|_{p}\leq p^{N}}\hat{f}_{\mathbf{m}}\left(t\right)e^{2\pi i\left\{ t\mathfrak{z}\right\} _{p}}
\end{equation}
\end{defn}
\begin{rem}
I would like to be able to explain rigorously why these completely
ignore the $X_{\mathbf{n}}\left(0\right)\prod_{n=0}^{N-1}\alpha_{\mathbf{n}}\left(p^{n}t\right)$
term of (\ref{eq:X_3,N recursion}), but I can't, so I won't. Suffice
it to say, we'll deal with that term in a bit.
\end{rem}
\begin{rem}
The attentive reader will note that this guess does not agree with
the solution of \textbf{Proposition \ref{prop:X_3-hat formal solution}}
in the case where $\alpha_{\mathbf{n}}\left(0\right)\neq1$. We will
discuss this soon enough.
\end{rem}
Computing $\tilde{f}_{\mathbf{m},N}\left(\mathfrak{z}\right)$ will
be relatively easy, thanks to a Fourier transform identity given below.
Before we can state and prove it, however, we need an auxiliary identity:

For this, we will need the following summation identity:
\begin{prop}
\label{prop:adjoint}Let $A$ be an abelian group, let $\hat{f},\hat{g}:\hat{\mathbb{Z}}_{p}\rightarrow A$,
and let $r\in\mathbb{N}_{1}$. Then, for any $N\geq r$:
\begin{equation}
\sum_{\left|t\right|_{p}\leq p^{N}}\hat{f}\left(t\right)\hat{g}\left(p^{r}t\right)=\sum_{\left|t\right|_{p}\leq p^{N-r}}\sum_{k=0}^{p^{r}-1}\hat{f}\left(\frac{t+k}{p^{r}}\right)\hat{g}\left(t\right)=\sum_{\left|t\right|_{p}\leq p^{N-r}}\sum_{\left|s\right|_{p}\leq p^{r}}\hat{f}\left(\frac{t}{p^{r}}+s\right)\hat{g}\left(t\right)
\end{equation}
where the operator $\sum_{\left|t\right|_{p}\leq p^{N-r}}$ is evaluation
at $t=0$ in the case where $N=r$.
\end{prop}
Proof: Direct computation.

Q.E.D.

\vphantom{}

The main Fourier transform identity we need is:
\begin{prop}
\label{prop:Fourier transforms and shifts}Let $\hat{\chi}:\hat{\mathbb{Z}}_{p}\rightarrow\mathcal{A}\otimes_{K}K\left(\zeta_{p^{\infty}}\right)$.
Then, for any constants $q_{0},\ldots,q_{p-1}$, all $n\in\mathbb{N}_{0}$,
and all $N\geq n$:
\begin{equation}
\sum_{\left|t\right|_{p}\leq p^{N}}\left(\hat{\chi}\left(p^{n}t\right)\prod_{k=0}^{n-1}\left(\frac{1}{p}\sum_{j=0}^{p-1}q_{j}e^{-2\pi ijp^{k}t}\right)\right)e^{2\pi i\left\{ t\mathfrak{z}\right\} _{p}}=\left(\prod_{k=0}^{n-1}q_{\left[\theta_{p}^{\circ k}\left(\mathfrak{z}\right)\right]_{p}}\right)\tilde{\chi}_{N-n}\left(\theta_{p}^{\circ n}\left(\mathfrak{z}\right)\right),\forall\mathfrak{z}\in\mathbb{Z}_{p}
\end{equation}
where:
\begin{equation}
\tilde{\chi}_{L}\left(\mathfrak{z}\right)=\sum_{\left|t\right|_{p}\leq p^{L}}\hat{\chi}\left(t\right)e^{2\pi i\left\{ t\mathfrak{z}\right\} _{p}}
\end{equation}
In particular, this shows that $\left(\prod_{k=0}^{n-1}q_{\left[\theta_{p}^{\circ k}\left(\mathfrak{z}\right)\right]_{p}}\right)\tilde{\chi}_{L}\left(\theta_{p}^{\circ n}\left(\mathfrak{z}\right)\right)$
has:
\begin{equation}
\hat{\chi}\left(p^{n}t\right)\prod_{k=0}^{n-1}\left(\frac{1}{p}\sum_{j=0}^{p-1}q_{j}e^{-2\pi ijp^{k}t}\right)
\end{equation}
as a Fourier transform.
\end{prop}
Proof:
\begin{claim}
Let $g:\mathbb{Z}_{p}\rightarrow K$ be locally constant. Then, we
have the Fourier transform pair:
\begin{equation}
q_{\left[\mathfrak{z}\right]_{p}}g\left(\theta_{p}\left(\mathfrak{z}\right)\right)\sim\left(\frac{1}{p}\sum_{j=0}^{p-1}q_{j}e^{-2\pi ijt}\right)\hat{g}\left(pt\right)\label{eq:shift claim}
\end{equation}

Proof of claim: Let $g$ be locally constant mod $p^{N_{g}}$ for
some $N_{g}\geq0$. Then, note that for $n\geq0$:
\begin{equation}
f\left(\mathfrak{z}\right)\overset{\textrm{def}}{=}q_{\left[\mathfrak{z}\right]_{p}}g\left(\theta_{p}\left(\mathfrak{z}\right)\right)
\end{equation}
is locally constant as well. In particular, letting $n\geq1$ and
$\mathfrak{y}\in\mathbb{Z}_{p}$ be arbitrary:
\begin{align*}
f\left(\mathfrak{z}+p^{n}\mathfrak{y}\right) & =q_{\left[\mathfrak{z}+p^{n}\mathfrak{y}\right]_{p}}g\left(\frac{\mathfrak{z}+p^{n}\mathfrak{y}-\left[\mathfrak{z}+p^{n}\mathfrak{y}\right]_{p}}{p}\right)\\
 & =q_{\left[\mathfrak{z}\right]_{p}}g\left(\frac{\mathfrak{z}-\left[\mathfrak{z}\right]_{p}}{p}+p^{n-1}\mathfrak{y}\right)\\
 & =q_{\left[\mathfrak{z}\right]_{p}}g\left(\theta_{p}\left(\mathfrak{z}\right)+p^{n-1}\mathfrak{y}\right)
\end{align*}
So, choosing $n=N_{g}+1$, we get: 
\begin{align*}
f\left(\mathfrak{z}+p^{N_{g}+1}\mathfrak{y}\right) & =q_{\left[\mathfrak{z}\right]_{p}}g\left(\theta_{p}\left(\mathfrak{z}\right)+p^{N_{g}}\mathfrak{y}\right)\\
\left(g\textrm{ is loc. const. mod }p^{N_{g}}\right); & =q_{\left[\mathfrak{z}\right]_{p}}g\left(\theta_{p}\left(\mathfrak{z}\right)\right)\\
 & =f\left(\mathfrak{z}\right)
\end{align*}
Thus, we see that $f$ is locally constant mod $p^{N_{g}+1}$. As
such, $f$'s  Fourier transform is given by:
\begin{align*}
\hat{f}\left(t\right) & =\frac{1}{p^{N_{g}+1}}\sum_{n=0}^{p^{N_{g}+1}-1}f\left(n\right)e^{-2\pi itn}\\
 & =\frac{1}{p^{N_{g}+1}}\sum_{n=0}^{p^{N_{g}+1}-1}q_{\left[n\right]_{p}}g\left(\theta_{p}\left(n\right)\right)e^{-2\pi itn}\\
\left(\textrm{split sum mod }n\right); & =\frac{1}{p^{N_{g}+1}}\sum_{n=0}^{p^{N_{g}}-1}\sum_{j=0}^{p-1}q_{\left[pn+j\right]_{p}}g\left(\theta_{p}\left(pn+j\right)\right)e^{-2\pi it\left(pn+j\right)}\\
 & =\frac{1}{p^{N_{g}+1}}\sum_{n=0}^{p^{N_{g}}-1}\sum_{j=0}^{p-1}q_{\left[j\right]_{p}}g\left(n\right)e^{-2\pi it\left(pn+j\right)}\\
 & =\left(\frac{1}{p}\sum_{j=0}^{p-1}q_{j}e^{-2\pi ijt}\right)\underbrace{\frac{1}{p^{N_{g}}}\sum_{n=0}^{p^{N_{g}}-1}g\left(n\right)e^{-2\pi iptn}}_{\hat{g}\left(pt\right),\textrm{ since }g\textrm{ is loc. const. mod }p^{N_{g}}}
\end{align*}
This proves the claim. $\checkmark$
\end{claim}
So, letting $\chi$ be quasi-integrable, observe that $\tilde{\chi}_{L}$
is then a locally constant function, and $\mathbf{1}_{0}\left(p^{L}t\right)\hat{\chi}\left(t\right)$
is its Fourier transform. Hence, the above gives:
\begin{equation}
\sum_{t\in\hat{\mathbb{Z}}_{p}}\left(\mathbf{1}_{0}\left(p^{L+1}t\right)\hat{\chi}\left(pt\right)\left(\frac{1}{p}\sum_{j=0}^{p-1}q_{j}e^{-2\pi ijt}\right)\right)e^{2\pi i\left\{ t\mathfrak{z}\right\} _{p}}=\underbrace{q_{\left[\mathfrak{z}\right]_{p}}\tilde{\chi}_{L}\left(\theta_{p}\left(\mathfrak{z}\right)\right)}_{g_{1}\left(\mathfrak{z}\right)},\forall\mathfrak{z}\in\mathbb{Z}_{p}
\end{equation}
which holds for all $L$.

Next, fixing $L$, let $g_{n}\left(\mathfrak{z}\right)=q_{\left[\theta_{p}^{\circ n}\left(\mathfrak{z}\right)\right]_{p}}\tilde{\chi}_{L}\left(\theta_{p}^{\circ n+1}\left(\mathfrak{z}\right)\right)$.
Since $g_{1}$ is locally constant, we have the Fourier transform
pair:
\begin{align*}
q_{\left[\mathfrak{z}\right]_{p}}g_{1}\left(\theta_{p}\left(\mathfrak{z}\right)\right) & \sim\hat{g}_{1}\left(pt\right)\left(\frac{1}{p}\sum_{j=0}^{p-1}q_{j}e^{-2\pi ijt}\right)\\
 & =\left(\mathbf{1}_{0}\left(p^{L+2}t\right)\hat{\chi}\left(p^{2}t\right)\left(\frac{1}{p}\sum_{j=0}^{p-1}q_{j}e^{-2\pi ijpt}\right)\right)\left(\frac{1}{p}\sum_{j=0}^{p-1}q_{j}e^{-2\pi ijt}\right)\\
 & =\mathbf{1}_{0}\left(p^{L+2}t\right)\hat{\chi}\left(p^{2}t\right)\prod_{k=0}^{2-1}\left(\frac{1}{p}\sum_{j=0}^{p-1}q_{j}e^{-2\pi ijp^{k}t}\right)
\end{align*}
Since:
\begin{equation}
q_{\left[\mathfrak{z}\right]_{p}}g_{1}\left(\theta_{p}\left(\mathfrak{z}\right)\right)=q_{\left[\mathfrak{z}\right]_{p}}q_{\left[\theta_{p}\left(\mathfrak{z}\right)\right]_{p}}\tilde{\chi}_{L}\left(\theta_{p}^{\circ2}\left(\mathfrak{z}\right)\right)
\end{equation}
by induction, we have the Fourier transform pair:
\begin{equation}
\left(\prod_{k=0}^{n-1}q_{\left[\theta_{p}^{\circ k}\left(\mathfrak{z}\right)\right]_{p}}\right)\tilde{\chi}_{L}\left(\theta_{p}^{\circ n}\left(\mathfrak{z}\right)\right)\sim\mathbf{1}_{0}\left(p^{L+n}t\right)\hat{\chi}\left(p^{n}t\right)\prod_{k=0}^{n-1}\left(\frac{1}{p}\sum_{j=0}^{p-1}q_{j}e^{-2\pi ijp^{k}t}\right)
\end{equation}
and hence:
\begin{equation}
\sum_{\left|t\right|_{p}\leq p^{L+n}}\left(\hat{\chi}\left(p^{n}t\right)\prod_{k=0}^{n-1}\left(\frac{1}{p}\sum_{j=0}^{p-1}q_{j}e^{-2\pi ijp^{k}t}\right)\right)e^{2\pi i\left\{ t\mathfrak{z}\right\} _{p}}=\left(\prod_{k=0}^{n-1}q_{\left[\theta_{p}^{\circ k}\left(\mathfrak{z}\right)\right]_{p}}\right)\tilde{\chi}_{L}\left(\theta_{p}^{\circ n}\left(\mathfrak{z}\right)\right),\forall\mathfrak{z}\in\mathbb{Z}_{p}
\end{equation}
Finally, letting $N\geq n$, set $L=N-n$.

Q.E.D.

\vphantom{}Using this, we can compute the Fourier series generated
by the $\hat{f}_{\mathbf{m}}$s.
\begin{notation}
For all $\mathbf{m},\mathbf{n}\in\mathbb{N}_{0}^{d}$ with $\mathbf{m}\leq\mathbf{n}$,
define $\kappa_{\mathbf{m},\mathbf{n}}:\mathbb{N}_{0}\rightarrow\textrm{Frac}\mathcal{A}$
be defined by:
\begin{equation}
\kappa_{\mathbf{m},\mathbf{n}}\left(n\right)\overset{\textrm{def}}{=}\prod_{k=1}^{p-1}\left(\frac{r_{\mathbf{m},\mathbf{n},k}}{r_{\mathbf{m},\mathbf{n},0}}\right)^{\#_{p:k}\left(n\right)}
\end{equation}
and let $\kappa_{\mathbf{n}}$ denote $\kappa_{\mathbf{n},\mathbf{n}}$.
Also, write:
\begin{equation}
\hat{Y}_{\mathbf{m},\mathbf{n}}\left(t\right)\overset{\textrm{def}}{=}\alpha_{\mathbf{m},\mathbf{n}}\left(t\right)\hat{X}_{\mathbf{m}}\left(pt\right)
\end{equation}
\end{notation}
We break the computation into several steps:
\begin{prop}
\label{prop:1st Y_j formula}For all $N\geq0$, $n\geq1$, $\mathbf{m}\leq\mathbf{n}$,
and $\mathfrak{z}\in\mathbb{Z}_{p}$:
\begin{equation}
\sum_{\left|t\right|_{p}\leq p^{N}}\left(\prod_{m=0}^{n-1}\alpha_{\mathbf{n}}\left(p^{m}t\right)\right)\hat{Y}_{\mathbf{m},\mathbf{n}}\left(p^{n}t\right)e^{2\pi i\left\{ t\mathfrak{z}\right\} _{p}}=r_{\mathbf{m},\mathbf{n},0}^{n}\kappa_{\mathbf{n}}\left(\left[\mathfrak{z}\right]_{p^{n}}\right)r_{\mathbf{m},\mathbf{n},\left[\theta_{p}^{\circ n}\left(\mathfrak{z}\right)\right]_{p}}\tilde{X}_{\mathbf{m},N-n-1}\left(\theta_{p}^{\circ n+1}\left(\mathfrak{z}\right)\right)
\end{equation}
where:
\begin{equation}
\tilde{X}_{\mathbf{m},N-n-1}\left(\mathfrak{y}\right)=\sum_{\left|t\right|_{p}\leq p^{N-n-1}}\hat{X}_{\mathbf{m}}\left(t\right)e^{2\pi i\left\{ t\mathfrak{y}\right\} _{p}},\textrm{ }\forall\mathfrak{y}\in\mathbb{Z}_{p}
\end{equation}
\end{prop}
Proof: By \textbf{Proposition \ref{prop:Fourier transforms and shifts}},
we can write:
\begin{align*}
\sum_{\left|t\right|_{p}\leq p^{N}}\left(\prod_{m=0}^{n-1}\alpha_{\mathbf{n}}\left(p^{m}t\right)\right)\hat{Y}_{\mathbf{m},\mathbf{n}}\left(p^{n}t\right)e^{2\pi i\left\{ t\mathfrak{z}\right\} _{p}} & =\left(\prod_{k=0}^{n-1}r_{\mathbf{n},\left[\theta_{p}^{\circ k}\left(\mathfrak{z}\right)\right]_{p}}\right)\tilde{Y}_{\mathbf{m},\mathbf{n},N-n}\left(\theta_{p}^{\circ n}\left(\mathfrak{z}\right)\right)\\
\left(\textrm{Use \textbf{Proposition} \textbf{\ref{prop:Product identity for kappa_X}}}\right); & =r_{\mathbf{n},0}^{n}\kappa_{\mathbf{n}}\left(\left[\mathfrak{z}\right]_{p^{n}}\right)\tilde{Y}_{\mathbf{m},\mathbf{n},N-n}\left(\theta_{p}^{\circ n}\left(\mathfrak{z}\right)\right)
\end{align*}
Next, letting $\mathfrak{y}\in\mathbb{Z}_{p}$, we have:
\begin{align*}
\tilde{Y}_{\mathbf{m},\mathbf{n},N-n}\left(\mathfrak{y}\right) & =\sum_{\left|t\right|_{p}\leq p^{N-n}}\hat{Y}_{\mathbf{m},\mathbf{n}}\left(t\right)e^{2\pi i\left\{ t\mathfrak{y}\right\} _{p}}\\
 & =\sum_{\left|t\right|_{p}\leq p^{N-n}}\alpha_{\mathbf{m},\mathbf{n}}\left(t\right)\hat{X}_{\mathbf{m}}\left(pt\right)e^{2\pi i\left\{ t\mathfrak{y}\right\} _{p}}\\
\left(\textrm{Use \textbf{Proposition} \textbf{\ref{prop:Fourier transforms and shifts}}}\right); & =r_{\mathbf{m},\mathbf{n},\left[\mathfrak{y}\right]_{p}}\tilde{X}_{\mathbf{m},N-n-1}\left(\theta_{p}\left(\mathfrak{y}\right)\right)
\end{align*}
Finally, set $\mathfrak{y}=\theta_{p}^{\circ n}\left(\mathfrak{z}\right)$.

Q.E.D.
\begin{prop}
\label{prop:partial sum of f-J-hat fourier series}For all $\mathbf{m},\mathbf{n}\in\mathbb{N}_{0}^{d}$
with $\mathbf{m}\leq\mathbf{n}$, and for all $N\geq2$, we have:

\begin{equation}
\sum_{\left|t\right|_{p}\leq p^{N}}\hat{f}_{\mathbf{m},\mathbf{n}}\left(t\right)e^{2\pi i\left\{ t\mathfrak{z}\right\} _{p}}=\sum_{n=0}^{N-2}r_{\mathbf{n},0}^{n}\kappa_{\mathbf{n}}\left(\left[\mathfrak{z}\right]_{p^{n}}\right)r_{\mathbf{m},\mathbf{n},\left[\theta_{p}^{\circ n}\left(\mathfrak{z}\right)\right]_{p}}\left(\tilde{X}_{\mathbf{m},N-n-1}\left(\theta_{p}^{\circ n+1}\left(\mathfrak{z}\right)\right)-\hat{X}_{\mathbf{m}}\left(0\right)\right)\label{eq:ready for twiddle limit}
\end{equation}
for all $\mathfrak{z}\in\mathbb{Z}_{p}$, so that:
\begin{align}
\tilde{f}_{\mathbf{n},N}\left(\mathfrak{z}\right) & =\sum_{\mathbf{m}<\mathbf{n}}\sum_{\left|t\right|_{p}\leq p^{N}}\hat{f}_{\mathbf{m},\mathbf{n}}\left(t\right)e^{2\pi i\left\{ t\mathfrak{z}\right\} _{p}}\nonumber \\
 & =\sum_{n=0}^{N-2}r_{\mathbf{n},0}^{n}\kappa_{\mathbf{n}}\left(\left[\mathfrak{z}\right]_{p^{n}}\right)\sum_{\mathbf{m}<\mathbf{n}}r_{\mathbf{m},\mathbf{n},\left[\theta_{p}^{\circ n}\left(\mathfrak{z}\right)\right]_{p}}\left(\tilde{X}_{\mathbf{m},N-n-1}\left(\theta_{p}^{\circ n+1}\left(\mathfrak{z}\right)\right)-\hat{X}_{\mathbf{m}}\left(0\right)\right)\label{eq:f_script J N twiddle}
\end{align}
\end{prop}
Proof: Since $\hat{f}_{\mathbf{m},\mathbf{n}}\left(t\right)$ vanishes
for $\left|t\right|_{p}\leq p$, we can write:
\begin{align*}
\sum_{\left|t\right|_{p}\leq p^{N}}\hat{f}_{\mathbf{m},\mathbf{n}}\left(t\right)e^{2\pi i\left\{ t\mathfrak{z}\right\} _{p}} & =\sum_{p^{2}\leq\left|t\right|_{p}\leq p^{N}}\sum_{n=0}^{-v_{p}\left(t\right)-2}\left(\prod_{m=0}^{n-1}\alpha_{\mathbf{n}}\left(p^{m}t\right)\right)\hat{Y}_{\mathbf{m},\mathbf{n}}\left(p^{n}t\right)e^{2\pi i\left\{ t\mathfrak{z}\right\} _{p}}\\
 & =\sum_{h=2}^{N}\sum_{\left|t\right|_{p}=p^{h}}\sum_{n=0}^{h-2}\left(\prod_{m=0}^{n-1}\alpha_{\mathbf{n}}\left(p^{m}t\right)\right)\hat{Y}_{\mathbf{m},\mathbf{n}}\left(p^{n}t\right)e^{2\pi i\left\{ t\mathfrak{z}\right\} _{p}}\\
\left(\textrm{re-arrange \ensuremath{\sum}s}\right); & =\sum_{n=0}^{N-2}\sum_{h=n+2}^{N}\sum_{\left|t\right|_{p}=p^{h}}\left(\prod_{m=0}^{n-1}\alpha_{\mathbf{n}}\left(p^{m}t\right)\right)\hat{Y}_{\mathbf{m},\mathbf{n}}\left(p^{n}t\right)e^{2\pi i\left\{ t\mathfrak{z}\right\} _{p}}\\
 & =\sum_{n=0}^{N-2}\sum_{p^{n+2}\leq\left|t\right|_{p}\leq p^{N}}\left(\prod_{m=0}^{n-1}\alpha_{\mathbf{n}}\left(p^{m}t\right)\right)\hat{Y}_{\mathbf{m},\mathbf{n}}\left(p^{n}t\right)e^{2\pi i\left\{ t\mathfrak{z}\right\} _{p}}
\end{align*}
Here, we write the $t$-sum on the bottom line as the difference:
\begin{equation}
\sum_{\left|t\right|_{p}\leq p^{N}}\left(\prod_{m=0}^{n-1}\alpha_{\mathbf{n}}\left(p^{m}t\right)\right)\hat{Y}_{\mathbf{m},\mathbf{n}}\left(p^{n}t\right)e^{2\pi i\left\{ t\mathfrak{z}\right\} _{p}}-\sum_{\left|t\right|_{p}\leq p^{n+1}}\left(\prod_{m=0}^{n-1}\alpha_{\mathbf{n}}\left(p^{m}t\right)\right)\hat{Y}_{\mathbf{m},\mathbf{n}}\left(p^{n}t\right)e^{2\pi i\left\{ t\mathfrak{z}\right\} _{p}}
\end{equation}
Applying \textbf{Proposition} \textbf{\ref{prop:1st Y_j formula}}
to this yields:

\begin{eqnarray*}
 & \sum_{\left|t\right|_{p}\leq p^{N}}\hat{f}_{\mathbf{m},\mathbf{n}}\left(t\right)e^{2\pi i\left\{ t\mathfrak{z}\right\} _{p}}\\
 & =\\
 & \sum_{n=0}^{N-2}r_{\mathbf{n},0}^{n}\kappa_{\mathbf{n}}\left(\left[\mathfrak{z}\right]_{p^{n}}\right)r_{\mathbf{m},\mathbf{n},\left[\theta_{p}^{\circ n}\left(\mathfrak{z}\right)\right]_{p}}\left(\tilde{X}_{\mathbf{m},N-n-1}\left(\theta_{p}^{\circ n+1}\left(\mathfrak{z}\right)\right)-\tilde{X}_{\mathbf{m},0}\left(\theta_{p}^{\circ n+1}\left(\mathfrak{z}\right)\right)\right)
\end{eqnarray*}
Finally, noting: 
\begin{equation}
\tilde{X}_{\mathbf{m},0}\left(\mathfrak{y}\right)=\sum_{\left|t\right|_{p}\leq0}\hat{X}_{\mathbf{m}}\left(t\right)e^{2\pi i\left\{ t\mathfrak{y}\right\} _{p}}=\hat{X}_{\mathbf{m}}\left(0\right),\textrm{ }\forall\mathfrak{y}\in\mathbb{Z}_{p}
\end{equation}
we can then write:
\begin{equation}
\sum_{n=0}^{N-2}r_{\mathbf{n},0}^{n}\kappa_{\mathbf{n}}\left(\left[\mathfrak{z}\right]_{p^{n}}\right)r_{\mathbf{m},\mathbf{n},\left[\theta_{p}^{\circ n}\left(\mathfrak{z}\right)\right]_{p}}\left(\tilde{X}_{\mathbf{m},N-n-1}\left(\theta_{p}^{\circ n+1}\left(\mathfrak{z}\right)\right)-\hat{X}_{\mathbf{m}}\left(0\right)\right)
\end{equation}

Q.E.D.

\vphantom{}

So, summing over all $\mathbf{m}<\mathbf{n}$, we have:
\begin{align*}
\tilde{f}_{\mathbf{n},N}\left(\mathfrak{z}\right) & =\sum_{\left|t\right|_{p}\leq p^{N}}\hat{f}_{\mathbf{n}}\left(t\right)e^{2\pi i\left\{ t\mathfrak{z}\right\} _{p}}\\
 & =\sum_{\mathbf{m}<\mathbf{n}}\sum_{\left|t\right|_{p}\leq p^{N}}\hat{f}_{\mathbf{m},\mathbf{n}}\left(t\right)e^{2\pi i\left\{ t\mathfrak{z}\right\} _{p}}\\
 & =\sum_{\mathbf{m}<\mathbf{n}}\sum_{n=0}^{N-2}r_{\mathbf{n},0}^{n}\kappa_{\mathbf{n}}\left(\left[\mathfrak{z}\right]_{p^{n}}\right)r_{\mathbf{m},\mathbf{n},\left[\theta_{p}^{\circ n}\left(\mathfrak{z}\right)\right]_{p}}\left(\tilde{X}_{\mathbf{m},N-n-1}\left(\theta_{p}^{\circ n+1}\left(\mathfrak{z}\right)\right)-\hat{X}_{\mathbf{m}}\left(0\right)\right)
\end{align*}
Now, we skip over one of the biggest details:
\begin{assumption}
\label{assu:main limit lemma}For all $\mathbf{m}<\mathbf{n}$:
\begin{equation}
\lim_{N\rightarrow\infty}\tilde{f}_{\mathbf{m},\mathbf{n},N}\left(\mathfrak{z}\right)=\lim_{N\rightarrow\infty}\sum_{n=0}^{N-2}r_{\mathbf{n},0}^{n}\kappa_{\mathbf{n}}\left(\left[\mathfrak{z}\right]_{p^{n}}\right)r_{\mathbf{m},\mathbf{n},\left[\theta_{p}^{\circ n}\left(\mathfrak{z}\right)\right]_{p}}\left(\tilde{X}_{\mathbf{m},N-n-1}\left(\theta_{p}^{\circ n+1}\left(\mathfrak{z}\right)\right)-\hat{X}_{\mathbf{m}}\left(0\right)\right)\label{eq:assumption}
\end{equation}
is equal to:
\begin{equation}
f_{\mathbf{m},\mathbf{n}}\left(\mathfrak{z}\right)\overset{\textrm{def}}{=}\sum_{n=0}^{\infty}r_{\mathbf{n},0}^{n}\kappa_{\mathbf{n}}\left(\left[\mathfrak{z}\right]_{p^{n}}\right)r_{\mathbf{m},\mathbf{n},\left[\theta_{p}^{\circ n}\left(\mathfrak{z}\right)\right]_{p}}\left(X_{\mathbf{m}}\left(\theta_{p}^{\circ n+1}\left(\mathfrak{z}\right)\right)-\hat{X}_{\mathbf{m}}\left(0\right)\right)
\end{equation}
Hence, setting:
\begin{equation}
f_{\mathbf{n}}\left(\mathfrak{z}\right)\overset{\textrm{def}}{=}\lim_{N\rightarrow\infty}\tilde{f}_{\mathbf{n},N}\left(\mathfrak{z}\right)=\sum_{\mathbf{m}<\mathbf{n}}\lim_{N\rightarrow\infty}\tilde{f}_{\mathbf{m},\mathbf{n},N}\left(\mathfrak{z}\right)
\end{equation}
we have that:
\begin{equation}
f_{\mathbf{n}}\left(\mathfrak{z}\right)=\sum_{n=0}^{\infty}r_{\mathbf{n},0}^{n}\kappa_{\mathbf{n}}\left(\left[\mathfrak{z}\right]_{p^{n}}\right)\sum_{\mathbf{m}<\mathbf{n}}r_{\mathbf{m},\mathbf{n},\left[\theta_{p}^{\circ n}\left(\mathfrak{z}\right)\right]_{p}}\left(X_{\mathbf{m}}\left(\theta_{p}^{\circ n+1}\left(\mathfrak{z}\right)\right)-\hat{X}_{\mathbf{m}}\left(0\right)\right)
\end{equation}
\end{assumption}
\begin{rem}
In this, note that we are avoiding any discussion of either the mode
of convergence and the set of $\mathfrak{z}\in\mathbb{Z}_{p}$ for
which the convergence occurs. Furthermore, we are assuming the quasi-integrability
of $X_{\mathbf{m}}$ for all $\mathbf{m}<\mathbf{n}$. These details
will be dealt with in \textbf{Sections \ref{sec:The-Small-Ideas}}
and \textbf{\ref{sec:Making-Things-Precise}}. The notion of frames
presented in \textbf{Section \ref{subsec:Frame-Theory-=000026}} was
originally tailor-made to give a rigorous foundation to the notion
of convergence occurring in (\ref{eq:assumption}). The nature of
quasi-integrability, meanwhile, is chronicled in \textbf{Section \ref{subsec:Quasi-Integrability-=000026-Degenerate}}.
Moreover, while the work of \textbf{Section \ref{sec:The-Small-Ideas}
}will give us the language needed to \emph{define} the sense in which
(\ref{eq:assumption}) converges, we will still need to formulate
an argument to justify that the limit in (\ref{eq:assumption}) actually
converges in the sense defined in \textbf{Section \ref{sec:The-Small-Ideas}}.
This argument, as well as the terminology and preparatory work needed
to see it through, are covered in \textbf{Section \ref{sec:Making-Things-Precise}},
with \textbf{Lemma \ref{lem:main limit lemma}} on page \pageref{lem:main limit lemma}
giving sufficient conditions for which a general form of the limit
in (\ref{eq:assumption}) holds in the sense of \textbf{Section \ref{subsec:Frame-Theory-=000026}}.
\end{rem}
\begin{lem}
\label{lem:Assuming-Assumption-}If \textbf{Assumption \ref{assu:main limit lemma}
}holds, $f_{\mathbf{n}}$ is quasi-integrable, with:
\begin{equation}
\hat{f}_{\mathbf{n}}\left(t\right)=\sum_{\mathbf{m}<\mathbf{n}}\hat{f}_{\mathbf{m},\mathbf{n}}\left(t\right)
\end{equation}
as a Fourier transform.
\end{lem}
Provided everything up to this point is justified, this will be enough
to prove the quasi-integrability of $X_{\mathbf{n}}$.
\begin{thm}
\label{thm:formal quasi-integrability of X script J}$X_{\mathbf{n}}$
is quasi-integrable. In particular, by recursion-induction, for all
$\mathbf{n}$ with $\Sigma\left(\mathbf{n}\right)\geq2$, we have:
\begin{equation}
\hat{X}_{\mathbf{n}}\left(t\right)=\hat{f}_{\mathbf{n}}\left(t\right)-\hat{g}_{\mathbf{n}}\left(t\right)
\end{equation}
where:
\begin{equation}
\hat{f}_{\mathbf{n}}\left(t\right)=\sum_{\mathbf{m}<\mathbf{n}}\sum_{n=0}^{-v_{p}\left(t\right)-2}\left(\prod_{m=0}^{n-1}\alpha_{\mathbf{n}}\left(p^{m}t\right)\right)\alpha_{\mathbf{m},\mathbf{n}}\left(p^{n}t\right)\hat{X}_{\mathbf{m}}\left(p^{n+1}t\right)\label{eq:f_script J hat}
\end{equation}
and where:
\begin{equation}
\hat{g}_{\mathbf{n}}\left(t\right)=\begin{cases}
0 & \textrm{if }t=0\\
\left(\beta_{\mathbf{n}}\left(0\right)v_{p}\left(t\right)+\gamma_{\mathbf{n}}\left(\frac{t\left|t\right|_{p}}{p}\right)\right)\hat{A}_{\mathbf{n}}\left(t\right) & \textrm{if }t\neq0
\end{cases},\textrm{ }\forall t\in\hat{\mathbb{Z}}_{p}\label{eq:g_script J hat , alpha equals 1}
\end{equation}
if $\alpha_{\mathbf{n}}\left(0\right)=1$ and:
\begin{equation}
\hat{g}_{\mathbf{n}}\left(t\right)=\begin{cases}
\frac{\beta_{\mathbf{n}}\left(0\right)}{1-\alpha_{\mathbf{n}}\left(0\right)} & \textrm{if }t=0\\
\left(\frac{\beta_{\mathbf{n}}\left(0\right)}{1-\alpha_{\mathbf{n}}\left(0\right)}+\gamma_{\mathbf{n}}\left(\frac{t\left|t\right|_{p}}{p}\right)\right)\hat{A}_{\mathbf{n}}\left(t\right) & \textrm{if }t\neq0
\end{cases},\textrm{ }\forall t\in\hat{\mathbb{Z}}_{p}\label{eq:g_script J hat , alpha not equal to 1}
\end{equation}
if $\alpha_{\mathbf{n}}\left(0\right)\neq1$, where $\alpha_{\mathbf{n}},\beta_{\mathbf{n}},\gamma_{\mathbf{n}},\hat{A}_{\mathbf{n}}:\hat{\mathbb{Z}}_{p}\rightarrow\textrm{Frac}\left(\mathcal{A}\right)\otimes_{K}K\left(\zeta_{p^{\infty}}\right)$
are given by:
\begin{equation}
\alpha_{\mathbf{n}}\left(t\right)\overset{\textrm{def}}{=}\frac{1}{p}\sum_{k=0}^{p-1}r_{\mathbf{n},k}e^{-2\pi ikt}
\end{equation}
\begin{equation}
\beta_{\mathbf{n}}\left(t\right)\overset{\textrm{def}}{=}\frac{1}{p}\sum_{k=0}^{p-1}c_{\mathbf{n},k}e^{-2\pi ikt}
\end{equation}
\begin{equation}
\gamma_{\mathbf{n}}\left(t\right)\overset{\textrm{def}}{=}\frac{\beta_{\mathbf{n}}\left(t\right)}{\alpha_{\mathbf{n}}\left(t\right)}
\end{equation}
\begin{equation}
\hat{A}_{\mathbf{n}}\left(t\right)\overset{\textrm{def}}{=}\prod_{n=0}^{-v_{p}\left(t\right)-1}\alpha_{\mathbf{n}}\left(p^{n}t\right)
\end{equation}
where:
\begin{equation}
c_{\mathbf{n},k}=-\sum_{\mathbf{m}<\mathbf{n}}r_{\mathbf{m},k}\hat{X}_{\mathbf{m}}\left(0\right)
\end{equation}
\end{thm}
Proof: Let $\mathfrak{z}$ and $\mathbf{m}<\mathbf{n}$ be arbitrary,
let $k\in\left\{ 0,\ldots,p-1\right\} $ and pull out the $n=0$ term
of $f_{\mathbf{m},\mathbf{n}}$ (note: $\kappa_{\mathbf{n}}\left(\left[\mathfrak{z}\right]_{p^{0}}\right)=\kappa_{\mathbf{n}}\left(0\right)=1$).
This gives:
\begin{align*}
f_{\mathbf{m},\mathbf{n}}\left(p\mathfrak{z}+k\right) & =r_{\mathbf{m},\mathbf{n},k}\left(X_{\mathbf{m}}\left(\theta_{p}\left(p\mathfrak{z}+k\right)\right)-\hat{X}_{\mathbf{m}}\left(0\right)\right)\\
 & +\sum_{n=1}^{\infty}r_{\mathbf{n},0}^{n}\kappa_{\mathbf{n}}\left(\left[p\mathfrak{z}+k\right]_{p^{n}}\right)r_{\mathbf{m},\mathbf{n},\left[\theta_{p}^{\circ n}\left(p\mathfrak{z}+k\right)\right]_{p}}\left(X_{\mathbf{m}}\left(\theta_{p}^{\circ n+1}\left(p\mathfrak{z}+k\right)\right)-\hat{X}_{\mathbf{m}}\left(0\right)\right)\\
\left(\mathbf{Proposition}\textrm{ }\mathbf{\ref{prop:fundamental functional equations-1}}\right); & =r_{\mathbf{m},\mathbf{n},k}\left(X_{\mathbf{m}}\left(\mathfrak{z}\right)-\hat{X}_{\mathbf{m}}\left(0\right)\right)\\
 & +\underbrace{r_{\mathbf{n},0}\kappa_{\mathbf{n}}\left(k\right)}_{r_{\mathbf{n},k}}\underbrace{\sum_{n=1}^{\infty}r_{\mathbf{n},0}^{n-1}\kappa_{\mathbf{n}}\left(\left[\mathfrak{z}\right]_{p^{n-1}}\right)r_{\mathbf{m},\left[\theta_{p}^{\circ n-1}\left(\mathfrak{z}\right)\right]_{p}}\left(X_{\mathbf{m}}\left(\theta_{p}^{\circ n}\left(\mathfrak{z}\right)\right)-\hat{X}_{\mathbf{m}}\left(0\right)\right)}_{f_{\mathbf{m},\mathbf{n}}\left(\mathfrak{z}\right)}\\
 & =r_{\mathbf{m},\mathbf{n},k}\left(X_{\mathbf{m}}\left(\mathfrak{z}\right)-\hat{X}_{\mathbf{m}}\left(0\right)\right)+r_{\mathbf{n},k}f_{\mathbf{m},\mathbf{n}}\left(\mathfrak{z}\right)
\end{align*}

\begin{align*}
f_{\mathbf{n}}\left(p\mathfrak{z}+k\right) & =\sum_{\mathbf{m}<\mathbf{n}}f_{\mathbf{m}}\left(p\mathfrak{z}+k\right)\\
 & =\sum_{\mathbf{m}<\mathbf{n}}\left(r_{\mathbf{m},\mathbf{n},k}\left(X_{\mathbf{m}}\left(\mathfrak{z}\right)-\hat{X}_{\mathbf{m}}\left(0\right)\right)+r_{\mathbf{n},k}f_{\mathbf{m},\mathbf{n}}\left(\mathfrak{z}\right)\right)\\
 & =\sum_{\mathbf{m}<\mathbf{n}}r_{\mathbf{m},\mathbf{n},k}\left(X_{\mathbf{m}}\left(\mathfrak{z}\right)-\hat{X}_{\mathbf{m}}\left(0\right)\right)+r_{\mathbf{n},k}\underbrace{\sum_{\mathbf{m}<\mathbf{n}}f_{\mathbf{m},\mathbf{n}}\left(\mathfrak{z}\right)}_{f_{\mathbf{n}}\left(\mathfrak{z}\right)}
\end{align*}
Using \textbf{Proposition \ref{prop:Z functional equation}}, we have:
\begin{equation}
\sum_{\mathbf{m}<\mathbf{n}}r_{\mathbf{m},\mathbf{n},k}X_{\mathbf{m}}\left(\mathfrak{z}\right)=X_{\mathbf{n}}\left(p\mathfrak{z}+k\right)-r_{\mathbf{n},k}X_{\mathbf{n}}\left(\mathfrak{z}\right)
\end{equation}
Thus:
\begin{align*}
f_{\mathbf{n}}\left(p\mathfrak{z}+k\right) & =\sum_{\mathbf{m}<\mathbf{n}}r_{\mathbf{m},\mathbf{n},k}\left(X_{\mathbf{m}}\left(\mathfrak{z}\right)-\hat{X}_{\mathbf{m}}\left(0\right)\right)+r_{\mathbf{n},k}f_{\mathbf{n}}\left(\mathfrak{z}\right)\\
 & =X_{\mathbf{n}}\left(p\mathfrak{z}+k\right)-r_{\mathbf{n},k}X_{\mathbf{n}}\left(\mathfrak{z}\right)-\sum_{\mathbf{m}<\mathbf{n}}r_{\mathbf{m},\mathbf{n},k}\hat{X}_{\mathbf{m}}\left(0\right)+r_{\mathbf{n},k}f_{\mathbf{n}}\left(\mathfrak{z}\right)\\
 & \Updownarrow\\
f_{\mathbf{n}}\left(p\mathfrak{z}+k\right)-X_{\mathbf{n}}\left(p\mathfrak{z}+k\right) & =r_{\mathbf{n},k}\left(f_{\mathbf{n}}\left(\mathfrak{z}\right)-X_{\mathbf{n}}\left(\mathfrak{z}\right)\right)-\sum_{\mathbf{m}<\mathbf{n}}r_{\mathbf{m},k}\hat{X}_{\mathbf{m}}\left(0\right)
\end{align*}
Now, set $g_{\mathbf{n}}\left(\mathfrak{z}\right)\overset{\textrm{def}}{=}f_{\mathbf{n}}\left(\mathfrak{z}\right)-X_{\mathbf{n}}\left(\mathfrak{z}\right)$,
and let:
\begin{equation}
c_{\mathbf{n},k}\overset{\textrm{def}}{=}-\sum_{\mathbf{m}<\mathbf{n}}r_{\mathbf{m},k}\hat{X}_{\mathbf{m}}\left(0\right),\textrm{ }\forall k\in\left\{ 0,\ldots,p-1\right\} 
\end{equation}
Then, we have:
\begin{equation}
g_{\mathbf{n}}\left(p\mathfrak{z}+k\right)=r_{\mathbf{n},k}g\left(\mathfrak{z}\right)+c_{\mathbf{n},k},\textrm{ }\forall\mathfrak{z},\textrm{ }\forall k\in\left\{ 0,\ldots,p-1\right\} 
\end{equation}
Ah, but this shows that $g_{\mathbf{n}}$ is a foliated F-series,
and we already know how to compute the Fourier transform of such a
function by using the formulae from \textbf{Theorem \ref{thm:The-breakdown-variety}}.
This shows that $g_{\mathbf{n}}$ is quasi-integrable. Finally, since
$X_{\mathbf{n}}\left(\mathfrak{z}\right)=f_{\mathbf{n}}\left(\mathfrak{z}\right)-g_{\mathbf{n}}\left(\mathfrak{z}\right)$,
we see that $\hat{f}_{\mathbf{n}}-\hat{g}_{\mathbf{n}}$ is then a
Fourier transform of $X_{\mathbf{n}}$, which proves the quasi-integrability
of $X_{\mathbf{n}}$, and gives us a formula for it, as well!

Q.E.D.
\begin{note}
There are frame-theoretic concerns that need to be addressed here.
Fortunately, these are already ones that we've set aside for discussion
later on.
\end{note}
What remains is to deal with the inductive-recursive aspect of our
argument. Proving the quasi-integrability of a degree $d$ F-series
will require leveraging the quasi-integrability of F-series of degrees
$0$ through $d-1$. The actual balancing point is in \textbf{Assumption
\ref{assu:main limit lemma}}. Proving that this assumption is justified
requires knowing an \emph{explicit} upper bound on the difference:

\begin{equation}
\Delta_{N}^{\left(n\right)}\left(\mathfrak{z}\right)=\left(\tilde{X}_{N}\circ\theta_{p}^{\circ n}\right)\left(\mathfrak{z}\right)-\left(X\circ\theta_{p}^{\circ n}\right)\left(\mathfrak{z}\right)\label{eq:delta notation}
\end{equation}
\begin{equation}
\Delta_{\mathbf{m},N}^{\left(0\right)}\left(\mathfrak{z}\right)\overset{\textrm{def}}{=}\tilde{X}_{\mathbf{m},N}\left(\mathfrak{z}\right)-X_{\mathbf{m}}\left(\mathfrak{z}\right)\label{eq:initial delta}
\end{equation}
between $X_{\mathbf{m}}$ and the $N$th partial sum of $X_{\mathbf{m}}$'s
Fourier series. Obtaining this bound is a tedious but ultimately straightforward
matter of computing $\tilde{X}_{\mathbf{m},N}\left(\mathfrak{z}\right)$.
This is quite nice, because it tells us that we can deal with all
of the necessary computations in one fell swoop.

We'll deal with $\hat{g}_{\mathbf{n}}$ first. This is a simple matter
of computing (\ref{eq:initial delta}) when $\Sigma\left(\mathbf{n}\right)=1$.
\begin{prop}
\label{prop:A_X sum}Set:
\begin{equation}
\kappa_{X}\left(m\right)\overset{\textrm{def}}{=}\prod_{j=1}^{p-1}\left(\frac{a_{j}}{a_{0}}\right)^{\#_{p:j}\left(m\right)}
\end{equation}
Then:
\begin{equation}
\sum_{\left|t\right|_{p}\leq p^{N}}\hat{A}_{X}\left(t\right)e^{2\pi i\left\{ t\mathfrak{z}\right\} _{p}}=a_{0}^{N}\kappa_{X}\left(\left[\mathfrak{z}\right]_{p^{N}}\right)+\left(1-\alpha_{X}\left(0\right)\right)\sum_{n=0}^{N-1}a_{0}^{n}\kappa_{X}\left(\left[\mathfrak{z}\right]_{p^{n}}\right)\label{eq:A_X sum}
\end{equation}
for all $\mathfrak{z}\in\mathbb{Z}_{p}$ and all $N\in\mathbb{N}_{0}$.

Next, letting:
\begin{equation}
\alpha_{X}\left(t\right)\overset{\textrm{def}}{=}\frac{1}{p}\sum_{j=0}^{p-1}a_{j}e^{-2\pi ijt}
\end{equation}
\begin{equation}
\beta_{X}\left(t\right)\overset{\textrm{def}}{=}\frac{1}{p}\sum_{j=0}^{p-1}b_{j}e^{-2\pi ijt}
\end{equation}
\begin{equation}
\gamma_{X}\left(t\right)\overset{\textrm{def}}{=}\frac{\beta_{X}\left(t\right)}{\alpha_{X}\left(t\right)}
\end{equation}
we have:

\begin{equation}
\sum_{0<\left|t\right|_{p}\leq p^{N}}\gamma_{X}\left(\frac{t\left|t\right|_{p}}{p}\right)\hat{A}_{X}\left(t\right)e^{2\pi i\left\{ t\mathfrak{z}\right\} _{p}}=\sum_{n=0}^{N-1}\left(\sum_{j=1}^{p-1}\beta_{X}\left(\frac{j}{p}\right)\varepsilon_{n}^{j}\left(\mathfrak{z}\right)\right)a_{0}^{n}\kappa_{X}\left(\left[\mathfrak{z}\right]_{p^{n}}\right)\label{eq:gamma_X A_X sum}
\end{equation}
\end{prop}
Proof: First, we write:
\begin{align*}
\hat{A}_{X}\left(t\right) & =\prod_{m=0}^{-v_{p}\left(t\right)-1}\left(\frac{1}{p}\sum_{j=0}^{p-1}a_{j}e^{-2\pi ijp^{m}t}\right)\\
 & =\left(\frac{a_{0}}{p}\right)^{-v_{p}\left(t\right)}\prod_{m=0}^{-v_{p}\left(t\right)-1}\left(1+\sum_{j=1}^{p-1}\frac{a_{j}}{a_{0}}e^{-2\pi ijp^{m}t}\right)\\
 & \overset{!}{=}\left(\frac{a_{0}}{p}\right)^{-v_{p}\left(t\right)}\sum_{m=0}^{\left|t\right|_{p}-1}\left(\prod_{j=1}^{p-1}\left(\frac{a_{j}}{a_{0}}\right)^{\#_{p:j}\left(m\right)}\right)e^{-2\pi imt}
\end{align*}
The step marked (!) is an elementary digit-counting combinatorial
argument. Consequently:
\begin{align*}
\sum_{0<\left|t\right|_{p}\leq p^{N}}\hat{A}_{X}\left(t\right)e^{2\pi i\left\{ t\mathfrak{z}\right\} _{p}} & =\sum_{n=1}^{N}\sum_{\left|t\right|_{p}=p^{n}}\hat{A}_{X}\left(t\right)e^{2\pi i\left\{ t\mathfrak{z}\right\} _{p}}\\
 & =\sum_{n=1}^{N}\sum_{\left|t\right|_{p}=p^{n}}\left(\left(\frac{a_{0}}{p}\right)^{-v_{p}\left(t\right)}\sum_{m=0}^{\left|t\right|_{p}-1}\left(\prod_{j=1}^{p-1}\left(\frac{a_{j}}{a_{0}}\right)^{\#_{p:j}\left(m\right)}\right)e^{-2\pi imt}\right)e^{2\pi i\left\{ t\mathfrak{z}\right\} _{p}}\\
 & =\sum_{n=1}^{N}\left(\frac{a_{0}}{p}\right)^{n}\sum_{m=0}^{p^{n}-1}\left(\prod_{j=1}^{p-1}\left(\frac{a_{j}}{a_{0}}\right)^{\#_{p:j}\left(m\right)}\right)\sum_{\left|t\right|_{p}=p^{n}}e^{2\pi i\left\{ t\left(\mathfrak{z}-m\right)\right\} _{p}}\\
 & =\sum_{n=1}^{N}\left(\frac{a_{0}}{p}\right)^{n}\sum_{m=0}^{p^{n}-1}\left(\prod_{j=1}^{p-1}\left(\frac{a_{j}}{a_{0}}\right)^{\#_{p:j}\left(m\right)}\right)\left(p^{n}\left[\mathfrak{z}\overset{p^{n}}{\equiv}m\right]-p^{n-1}\left[\mathfrak{z}\overset{p^{n-1}}{\equiv}m\right]\right)
\end{align*}
Here, we have:
\begin{equation}
\sum_{m=0}^{p^{n}-1}\left(\prod_{j=1}^{p-1}\left(\frac{a_{j}}{a_{0}}\right)^{\#_{p:j}\left(m\right)}\right)\left[\mathfrak{z}\overset{p^{n}}{\equiv}m\right]=\prod_{j=1}^{p-1}\left(\frac{a_{j}}{a_{0}}\right)^{\#_{p:j}\left(\left[\mathfrak{z}\right]_{p^{n}}\right)}=\kappa_{X}\left(\left[\mathfrak{z}\right]_{p^{n}}\right)
\end{equation}
while:
\begin{align*}
\sum_{m=0}^{p^{n}-1}\left(\prod_{j=1}^{p-1}\left(\frac{a_{j}}{a_{0}}\right)^{\#_{p:j}\left(m\right)}\right)\left[\mathfrak{z}\overset{p^{n-1}}{\equiv}m\right] & =\sum_{k=0}^{p-1}\sum_{m=kp^{n-1}}^{\left(k+1\right)p^{n-1}-1}\left(\prod_{j=1}^{p-1}\left(\frac{a_{j}}{a_{0}}\right)^{\#_{p:j}\left(m\right)}\right)\left[\mathfrak{z}\overset{p^{n-1}}{\equiv}m\right]\\
 & =\sum_{k=0}^{p-1}\sum_{m=0}^{p^{n-1}-1}\left(\prod_{j=1}^{p-1}\left(\frac{a_{j}}{a_{0}}\right)^{\#_{p:j}\left(m+kp^{n-1}\right)}\right)\left[\mathfrak{z}\overset{p^{n-1}}{\equiv}m+kp^{n-1}\right]\\
 & =\sum_{k=0}^{p-1}\sum_{m=0}^{p^{n-1}-1}\left(\prod_{j=1}^{p-1}\left(\frac{a_{j}}{a_{0}}\right)^{\left[j=k\right]+\#_{p:j}\left(m\right)}\right)\left[\mathfrak{z}\overset{p^{n-1}}{\equiv}m\right]\\
 & =\sum_{k=0}^{p-1}\prod_{j=1}^{p-1}\left(\frac{a_{j}}{a_{0}}\right)^{\left[j=k\right]+\#_{p:j}\left(\left[\mathfrak{z}\right]_{p^{n-1}}\right)}\\
 & =\frac{1}{a_{0}}\left(\sum_{k=0}^{p-1}a_{k}\right)\kappa_{X}\left(\left[\mathfrak{z}\right]_{p^{n-1}}\right)
\end{align*}
and so:
\begin{align*}
\sum_{0<\left|t\right|_{p}\leq p^{N}}\hat{A}_{X}\left(t\right)e^{2\pi i\left\{ t\mathfrak{z}\right\} _{p}} & =\sum_{n=1}^{N}\left(\left(\frac{a_{0}}{p}\right)^{n}p^{n}\kappa_{X}\left(\left[\mathfrak{z}\right]_{p^{n}}\right)-\left(\frac{a_{0}}{p}\right)^{n}p^{n-1}\left(\frac{1}{a_{0}}\sum_{k=0}^{p-1}a_{k}\right)\kappa_{X}\left(\left[\mathfrak{z}\right]_{p^{n-1}}\right)\right)\\
 & =\sum_{n=1}^{N}a_{0}^{n}\kappa_{X}\left(\left[\mathfrak{z}\right]_{p^{n}}\right)-\sum_{n=1}^{N}a_{0}^{n-1}\underbrace{\left(\frac{1}{p}\sum_{k=0}^{p-1}a_{k}\right)}_{\alpha_{X}\left(0\right)}\kappa_{X}\left(\left[\mathfrak{z}\right]_{p^{n-1}}\right)\\
 & =\sum_{n=1}^{N}a_{0}^{n}\kappa_{X}\left(\left[\mathfrak{z}\right]_{p^{n}}\right)-\alpha_{X}\left(0\right)\sum_{n=0}^{N-1}a_{0}^{n}\kappa_{X}\left(\left[\mathfrak{z}\right]_{p^{n}}\right)\\
\left(\kappa_{X}\left(\left[\mathfrak{z}\right]_{p^{0}}\right)=1\right); & =a_{0}^{N}\kappa_{X}\left(\left[\mathfrak{z}\right]_{p^{N}}\right)-\alpha_{X}\left(0\right)+\left(1-\alpha_{X}\left(0\right)\right)\sum_{n=1}^{N-1}a_{0}^{n}\kappa_{X}\left(\left[\mathfrak{z}\right]_{p^{n}}\right)
\end{align*}
Thus:
\begin{equation}
\sum_{0<\left|t\right|_{p}\leq p^{N}}\hat{A}_{X}\left(t\right)e^{2\pi i\left\{ t\mathfrak{z}\right\} _{p}}=a_{0}^{N}\kappa_{X}\left(\left[\mathfrak{z}\right]_{p^{N}}\right)-\alpha_{X}\left(0\right)+\left(1-\alpha_{X}\left(0\right)\right)\sum_{n=1}^{N-1}a_{0}^{n}\kappa_{X}\left(\left[\mathfrak{z}\right]_{p^{n}}\right)
\end{equation}
Since $\hat{A}_{X}\left(0\right)=1$, we have:
\begin{align*}
\sum_{\left|t\right|_{p}\leq p^{N}}\hat{A}_{X}\left(t\right)e^{2\pi i\left\{ t\mathfrak{z}\right\} _{p}} & =a_{0}^{N}\kappa_{X}\left(\left[\mathfrak{z}\right]_{p^{N}}\right)+\left(1-\alpha_{X}\left(0\right)\right)+\left(1-\alpha_{X}\left(0\right)\right)\sum_{n=1}^{N-1}a_{0}^{n}\kappa_{X}\left(\left[\mathfrak{z}\right]_{p^{n}}\right)\\
 & =a_{0}^{N}\kappa_{X}\left(\left[\mathfrak{z}\right]_{p^{N}}\right)+\left(1-\alpha_{X}\left(0\right)\right)\sum_{n=0}^{N-1}a_{0}^{n}\kappa_{X}\left(\left[\mathfrak{z}\right]_{p^{n}}\right)
\end{align*}

The proof of (\ref{eq:gamma_X A_X sum}) is similar, only more involved.
One must deal with it by writing:
\begin{align}
\sum_{\left|t\right|_{p}\leq p^{N}}\gamma_{X}\left(\frac{t\left|t\right|_{p}}{p}\right)\hat{A}_{X}\left(t\right)e^{2\pi i\left\{ t\mathfrak{z}\right\} _{p}} & =\sum_{k=0}^{p^{N}-1}\gamma_{X}\left(\frac{k}{p}\right)\hat{A}_{X}\left(\frac{k}{p^{N}}\right)e^{2\pi i\left\{ k\mathfrak{z}/p^{N}\right\} _{p}}
\end{align}
and then using the sum decomposition formula:
\begin{equation}
\sum_{k=0}^{p^{N}-1}f\left(k\right)=\sum_{k=0}^{p^{N-1}-1}\sum_{j=0}^{p-1}f\left(pk+j\right)
\end{equation}
\cite{My Dissertation} contains the proof in \textbf{Lemma 4.4 }of
Section 4, on page 259.

Q.E.D.
\begin{prop}
\label{prop:X}Given a $p$-adic F-series $X$ characterized by:
\begin{equation}
X\left(p\mathfrak{z}+j\right)=a_{j}X\left(\mathfrak{z}\right)+b_{j},\textrm{ }\forall\mathfrak{z}\in\mathbb{Z}_{p},\textrm{ }\forall j\in\left\{ 0,\ldots,p-1\right\} 
\end{equation}
that is compatible with a frame $\mathcal{F}$, let $\hat{X}$ be
as given by \textbf{Theorem \ref{thm:The-breakdown-variety}}, though
without assuming that said theorem is true.

Under these conditions, we have that:
\begin{equation}
\hat{A}_{X}\left(pt\right)=\mathbf{1}_{0}\left(pt\right)\hat{A}_{X}\left(0\right)+\left(1-\mathbf{1}_{0}\left(pt\right)\right)\frac{\hat{A}_{X}\left(t\right)}{\alpha_{X}\left(t\right)}\label{eq:A_X-hat recurrence}
\end{equation}
and:
\begin{equation}
\alpha_{X}\left(t\right)\hat{X}\left(pt\right)=\begin{cases}
\mathbf{1}_{0}\left(pt\right)\hat{X}\left(0\right)\alpha_{X}\left(t\right)+\left(1-\mathbf{1}_{0}\left(pt\right)\right)\hat{X}\left(t\right) & \textrm{if }\alpha_{X}\left(0\right)\neq1\\
\left(1-\mathbf{1}_{0}\left(pt\right)\right)\left(\hat{X}\left(t\right)+\beta_{X}\left(0\right)\hat{A}_{X}\left(t\right)\right) & \textrm{if }\alpha_{X}\left(0\right)=1
\end{cases},\textrm{ }\forall t\in\hat{\mathbb{Z}}_{p}\label{eq:X-hat recurrence}
\end{equation}
Consequently:
\begin{equation}
\tilde{X}_{N}\left(\mathfrak{z}\right)=a_{\left[\mathfrak{z}\right]_{p}}\tilde{X}_{N-1}\left(\theta_{p}\left(\mathfrak{z}\right)\right)+b_{\left[\mathfrak{z}\right]_{p}}-\left[\alpha_{X}\left(0\right)=1\right]\beta_{X}\left(0\right)a_{0}^{N}\kappa_{X}\left(\left[\mathfrak{z}\right]_{p^{N}}\right)\label{eq:X_N twiddle recurrence formula}
\end{equation}
for all $\mathfrak{z}\in\mathbb{Z}_{p}$ and all $N\in\mathbb{N}_{0}$.
\end{prop}
Proof: The proof of (\ref{eq:A_X-hat recurrence}) follows by simple
re-indexing. When $\left|t\right|_{p}\geq p^{2}$:
\begin{align*}
\hat{A}_{X}\left(pt\right) & =\prod_{m=0}^{-v_{p}\left(pt\right)-1}\alpha_{X}\left(p^{m}pt\right)\\
 & =\prod_{m=0}^{-v_{p}\left(t\right)-2}\alpha_{X}\left(p^{m+1}t\right)\\
\left(\textrm{re-index}\right); & =\prod_{m=1}^{-v_{p}\left(t\right)-1}\alpha_{X}\left(p^{m}t\right)\\
\left(\times\frac{\alpha_{X}\left(t\right)}{\alpha_{X}\left(t\right)}\right) & =\frac{1}{\alpha_{X}\left(t\right)}\underbrace{\prod_{m=0}^{-v_{p}\left(t\right)-1}\alpha_{X}\left(p^{m}t\right)}_{\hat{A}_{X}\left(t\right)}
\end{align*}
When $\left|t\right|_{p}\leq p$, on the other hand, $pt$ is an integer,
and so, the $1$-periodicity of $\hat{A}_{X}$ forces $\hat{A}_{X}\left(pt\right)=\hat{A}_{X}\left(0\right)=1$.
Thus:
\begin{equation}
\hat{A}_{X}\left(pt\right)=\mathbf{1}_{0}\left(pt\right)\hat{A}_{X}\left(0\right)+\left(1-\mathbf{1}_{0}\left(pt\right)\right)\frac{\hat{A}_{X}\left(t\right)}{\alpha_{X}\left(t\right)}
\end{equation}
where $\mathbf{1}_{0}\left(pt\right)$ is $1$ if $\left|t\right|_{p}\leq p$
and is $0$ otherwise.

Using this, the proof of (\ref{eq:X-hat recurrence}) follows from
replacing $t$ in \textbf{Theorem \ref{thm:The-breakdown-variety}}'s
formulae for $\hat{X}$ with $pt$ and then applying (\ref{eq:A_X-hat recurrence}).
We have two cases, based on whether or not $\alpha_{X}\left(0\right)$
is equal to $1$.

\vphantom{}

I. Suppose $\alpha_{X}\left(0\right)\neq1$. Since $\left|pt\right|_{p}=0$
all $\left|t\right|_{p}\leq p$, we have
\begin{align*}
\hat{X}\left(pt\right) & =\begin{cases}
\frac{\beta_{X}\left(0\right)\hat{A}_{X}\left(0\right)}{1-\alpha_{X}\left(0\right)} & \textrm{if }\left|t\right|_{p}\leq p\\
\left(\frac{\beta_{X}\left(0\right)}{1-\alpha_{X}\left(0\right)}+\gamma_{X}\left(\frac{pt\left|pt\right|_{p}}{p}\right)\right)\hat{A}_{X}\left(pt\right) & \textrm{if }\left|t\right|_{p}\geq p^{2}
\end{cases}\\
\left(\textrm{Use }(\ref{eq:A_X-hat recurrence})\right); & =\begin{cases}
\frac{\beta_{X}\left(0\right)\hat{A}_{X}\left(0\right)}{1-\alpha_{X}\left(0\right)} & \textrm{if }\left|t\right|_{p}\leq p\\
\left(\frac{\beta_{X}\left(0\right)}{1-\alpha_{X}\left(0\right)}+\gamma_{X}\left(\frac{t\left|t\right|_{p}}{p}\right)\right)\frac{\hat{A}_{X}\left(t\right)}{\alpha_{X}\left(t\right)} & \textrm{if }\left|t\right|_{p}\geq p^{2}
\end{cases}\\
 & =\begin{cases}
\hat{X}\left(0\right) & \textrm{if }\left|t\right|_{p}\leq p\\
\frac{\hat{X}\left(t\right)}{\alpha_{X}\left(t\right)} & \textrm{if }\left|t\right|_{p}\geq p^{2}
\end{cases}
\end{align*}
Multiplying through by $\alpha_{X}\left(t\right)$ gives us:
\begin{equation}
\alpha_{X}\left(t\right)\hat{X}\left(pt\right)=\begin{cases}
\alpha_{X}\left(t\right)\hat{X}\left(0\right) & \textrm{if }\left|t\right|_{p}\leq p\\
\hat{X}\left(t\right) & \textrm{if }\left|t\right|_{p}\geq p^{2}
\end{cases}=\mathbf{1}_{0}\left(pt\right)\hat{X}\left(0\right)\alpha_{X}\left(t\right)+\left(1-\mathbf{1}_{0}\left(pt\right)\right)\hat{X}\left(t\right)
\end{equation}
as desired.

\vphantom{}II. Suppose $\alpha_{X}\left(0\right)=1$. Then, since
$v_{p}\left(pt\right)=v_{p}\left(t\right)+1$:

\begin{align*}
\hat{X}\left(pt\right) & =\begin{cases}
0 & \textrm{if }\left|t\right|_{p}\leq p\\
\left(\beta_{X}\left(0\right)v_{p}\left(pt\right)+\gamma_{X}\left(\frac{pt\left|pt\right|_{p}}{p}\right)\right)\hat{A}_{X}\left(pt\right) & \textrm{if }\left|t\right|_{p}\geq p^{2}
\end{cases}\\
\left(\textrm{Use }(\ref{eq:A_X-hat recurrence})\right); & =\begin{cases}
0 & \textrm{if }\left|t\right|_{p}\leq p\\
\beta_{X}\left(0\right)\frac{\hat{A}_{X}\left(t\right)}{\alpha_{X}\left(t\right)}+\underbrace{\left(\beta_{X}\left(0\right)v_{p}\left(t\right)+\gamma_{X}\left(\frac{t\left|t\right|_{p}}{p}\right)\right)\frac{\hat{A}_{X}\left(t\right)}{\alpha_{X}\left(t\right)}}_{\hat{X}\left(t\right)/\alpha_{X}\left(t\right)} & \textrm{if }\left|t\right|_{p}\geq p^{2}
\end{cases}\\
 & =\begin{cases}
0 & \textrm{if }\left|t\right|_{p}\leq p\\
\frac{\beta_{X}\left(0\right)\hat{A}_{X}\left(t\right)+\hat{X}\left(t\right)}{\alpha_{X}\left(t\right)} & \textrm{if }\left|t\right|_{p}\geq p^{2}
\end{cases}
\end{align*}
Multiplying through by $\alpha_{X}\left(t\right)$ gives the desired
result:
\begin{equation}
\alpha_{X}\left(t\right)\hat{X}\left(pt\right)=\left(1-\mathbf{1}_{0}\left(pt\right)\right)\left(\hat{X}\left(t\right)+\beta_{X}\left(0\right)\hat{A}_{X}\left(t\right)\right)
\end{equation}

Finally, we take (\ref{eq:X-hat recurrence}) and sum the generated
Fourier series. This necessitates two cases, based on whether or not
$\alpha_{X}\left(0\right)=1$.

\vphantom{}

I. When $\alpha_{X}\left(0\right)\neq1$ and $N\geq1$:

\begin{align*}
\sum_{\left|t\right|_{p}\leq p^{N}}\alpha_{X}\left(t\right)\hat{X}\left(pt\right)e^{2\pi i\left\{ t\mathfrak{z}\right\} _{p}} & =\hat{X}\left(0\right)\sum_{\left|t\right|_{p}\leq p}\alpha_{X}\left(t\right)e^{2\pi i\left\{ t\mathfrak{z}\right\} _{p}}+\sum_{p^{2}\leq\left|t\right|_{p}\leq p^{N}}\hat{X}\left(t\right)e^{2\pi i\left\{ t\mathfrak{z}\right\} _{p}}\\
 & =\hat{X}\left(0\right)\sum_{\left|t\right|_{p}\leq p}\alpha_{X}\left(t\right)e^{2\pi i\left\{ t\mathfrak{z}\right\} _{p}}+\underbrace{\sum_{\left|t\right|_{p}\leq p^{N}}\hat{X}\left(t\right)e^{2\pi i\left\{ t\mathfrak{z}\right\} _{p}}}_{\tilde{X}_{N}\left(\mathfrak{z}\right)}\\
 & -\sum_{\left|t\right|_{p}\leq p}\hat{X}\left(t\right)e^{2\pi i\left\{ t\mathfrak{z}\right\} _{p}}\\
 & =\tilde{X}_{N}\left(\mathfrak{z}\right)+\sum_{\left|t\right|_{p}\leq p}\left(\hat{X}\left(0\right)\alpha_{X}\left(t\right)-\hat{X}\left(t\right)\right)e^{2\pi i\left\{ t\mathfrak{z}\right\} _{p}}
\end{align*}
Here:
\begin{align}
\sum_{\left|t\right|_{p}\leq p}\left(\hat{X}\left(0\right)\alpha_{X}\left(t\right)-\hat{X}\left(t\right)\right)e^{2\pi i\left\{ t\mathfrak{z}\right\} _{p}} & =\hat{X}\left(0\right)\left(\alpha_{X}\left(0\right)-1\right)+\sum_{\left|t\right|_{p}=p}\left(\hat{X}\left(0\right)\alpha_{X}\left(t\right)-\hat{X}\left(t\right)\right)e^{2\pi i\left\{ t\mathfrak{z}\right\} _{p}}
\end{align}
Since $\hat{A}_{X}\left(t\right)=\alpha_{X}\left(t\right)$ for $\left|t\right|_{p}=p$,
using \textbf{\ref{thm:The-breakdown-variety}}'s formula for $\hat{X}$,
and noting that $t\left|t\right|_{p}/p=t$ when $\left|t\right|_{p}$,
we have:
\begin{align*}
\sum_{\left|t\right|_{p}=p}\left(\hat{X}\left(0\right)\alpha_{X}\left(t\right)-\hat{X}\left(t\right)\right)e^{2\pi i\left\{ t\mathfrak{z}\right\} _{p}} & =\sum_{\left|t\right|_{p}=p}\left(\hat{X}\left(0\right)-\frac{\beta_{X}\left(0\right)}{1-\alpha_{X}\left(0\right)}-\gamma_{X}\left(t\right)\right)\alpha_{X}\left(t\right)e^{2\pi i\left\{ t\mathfrak{z}\right\} _{p}}\\
\left(\gamma_{X}\left(t\right)=\frac{\beta_{X}\left(t\right)}{\alpha_{X}\left(t\right)}\right); & =\sum_{\left|t\right|_{p}=p}\left(\underbrace{\hat{X}\left(0\right)-\frac{\beta_{X}\left(0\right)}{1-\alpha_{X}\left(0\right)}}_{0}-\frac{\beta_{X}\left(t\right)}{\alpha_{X}\left(t\right)}\right)\alpha_{X}\left(t\right)e^{2\pi i\left\{ t\mathfrak{z}\right\} _{p}}\\
 & =-\sum_{\left|t\right|_{p}=p}\beta_{X}\left(t\right)e^{2\pi i\left\{ t\mathfrak{z}\right\} _{p}}\\
 & =-\frac{1}{p}\sum_{j=0}^{p-1}b_{j}\sum_{\left|t\right|_{p}=p}e^{2\pi i\left\{ t\left(\mathfrak{z}-j\right)\right\} _{p}}\\
 & =-\frac{1}{p}\sum_{j=0}^{p-1}b_{j}\left(p\left[\mathfrak{z}\overset{p}{\equiv}j\right]-1\right)\\
 & =-\sum_{j=0}^{p-1}b_{j}\left[\mathfrak{z}\overset{p}{\equiv}j\right]+\frac{1}{p}\sum_{j=0}^{p-1}b_{j}\\
 & =-b_{\left[\mathfrak{z}\right]_{p}}+\beta_{X}\left(0\right)
\end{align*}
and so:
\begin{align}
\sum_{\left|t\right|_{p}\leq p^{N}}\alpha_{X}\left(t\right)\hat{X}\left(pt\right)e^{2\pi i\left\{ t\mathfrak{z}\right\} _{p}} & =\tilde{X}_{N}\left(\mathfrak{z}\right)+\hat{X}\left(0\right)\left(\alpha_{X}\left(0\right)-1\right)-b_{\left[\mathfrak{z}\right]_{p}}+\beta_{X}\left(0\right)\label{eq:half}
\end{align}
Here, we note that:
\begin{align*}
\sum_{\left|t\right|_{p}\leq p^{N}}\alpha_{X}\left(t\right)\hat{X}\left(pt\right)e^{2\pi i\left\{ t\mathfrak{z}\right\} _{p}} & =\sum_{\left|t\right|_{p}\leq p^{N}}\left(\frac{1}{p}\sum_{j=0}^{p-1}a_{j}e^{-2\pi ijt}\right)\hat{X}\left(pt\right)e^{2\pi i\left\{ t\mathfrak{z}\right\} _{p}}\\
 & =\sum_{t\in\hat{\mathbb{Z}}_{p}}\left(\frac{1}{p}\sum_{j=0}^{p-1}a_{j}e^{-2\pi ijt}\right)\mathbf{1}_{0}\left(p^{N}t\right)\hat{X}\left(pt\right)e^{2\pi i\left\{ t\mathfrak{z}\right\} _{p}}\\
\left(\textrm{Use \textbf{Proposition }\textbf{\ref{prop:Fourier transforms and shifts}}}\right); & =a_{\left[\mathfrak{z}\right]_{p}}\sum_{t\in\hat{\mathbb{Z}}_{p}}\mathbf{1}_{0}\left(p^{N-1}t\right)\hat{X}\left(t\right)e^{2\pi i\left\{ t\theta_{p}\left(\mathfrak{z}\right)\right\} _{p}}\\
 & =a_{\left[\mathfrak{z}\right]_{p}}\tilde{X}_{N-1}\left(\theta_{p}\left(\mathfrak{z}\right)\right)
\end{align*}
Hence, (\ref{eq:half}) is:
\begin{align*}
a_{\left[\mathfrak{z}\right]_{p}}\tilde{X}_{N-1}\left(\theta_{p}\left(\mathfrak{z}\right)\right) & =\tilde{X}_{N}\left(\mathfrak{z}\right)+\hat{X}\left(0\right)\left(\alpha_{X}\left(0\right)-1\right)-b_{\left[\mathfrak{z}\right]_{p}}+\beta_{X}\left(0\right)\\
\left(\hat{X}\left(0\right)=\frac{\beta_{X}\left(0\right)\hat{A}_{X}\left(0\right)}{1-\alpha_{X}\left(0\right)}=\frac{\beta_{X}\left(0\right)}{1-\alpha_{X}\left(0\right)}\right); & =\tilde{X}_{N}\left(\mathfrak{z}\right)-\beta_{X}\left(0\right)-b_{\left[\mathfrak{z}\right]_{p}}+\beta_{X}\left(0\right)\\
 & =\tilde{X}_{N}\left(\mathfrak{z}\right)-b_{\left[\mathfrak{z}\right]_{p}}
\end{align*}
Adding $b_{\left[\mathfrak{z}\right]_{p}}$ to both sides gives the
desired result.

\vphantom{}II. When $\alpha_{X}\left(0\right)=1$: 
\[
\sum_{\left|t\right|_{p}\leq p^{N}}\alpha_{X}\left(t\right)\hat{X}\left(pt\right)e^{2\pi i\left\{ t\mathfrak{z}\right\} _{p}}=\sum_{p^{2}\leq\left|t\right|_{p}\leq p^{N}}\left(\hat{X}\left(t\right)+\beta_{X}\left(0\right)\hat{A}_{X}\left(t\right)\right)e^{2\pi i\left\{ t\mathfrak{z}\right\} _{p}}
\]
As we already know what the left-hand side is, we just need to deal
with the right-hand side. For this, \textbf{Proposition \ref{prop:A_X sum}}
gives us:
\[
\sum_{\left|t\right|_{p}\leq p^{N}}\beta_{X}\left(0\right)\hat{A}_{X}\left(t\right)e^{2\pi i\left\{ t\mathfrak{z}\right\} _{p}}=\beta_{X}\left(0\right)\left(a_{0}^{N}\kappa_{X}\left(\left[\mathfrak{z}\right]_{p^{N}}\right)+\underbrace{\left(1-\alpha_{X}\left(0\right)\right)\sum_{n=0}^{N-1}a_{0}^{n}\kappa_{X}\left(\left[\mathfrak{z}\right]_{p^{n}}\right)}_{0,\textrm{ since }\alpha_{X}\left(0\right)=1}\right)
\]
while:
\begin{equation}
\sum_{\left|t\right|_{p}\leq p}\beta_{X}\left(0\right)\hat{A}_{X}\left(t\right)e^{2\pi i\left\{ t\mathfrak{z}\right\} _{p}}=\beta_{X}\left(0\right)\left(a_{0}\kappa_{X}\left(\left[\mathfrak{z}\right]_{p}\right)+\underbrace{\left(1-\alpha_{X}\left(0\right)\right)}_{0,\textrm{ since }\alpha_{X}\left(0\right)=1}\right)
\end{equation}
Subtracting these two gives:
\begin{equation}
\sum_{p^{2}\leq\left|t\right|_{p}\leq p^{N}}\beta_{X}\left(0\right)\hat{A}_{X}\left(t\right)e^{2\pi i\left\{ t\mathfrak{z}\right\} _{p}}=\beta_{X}\left(0\right)\left(a_{0}^{N}\kappa_{X}\left(\left[\mathfrak{z}\right]_{p^{N}}\right)-a_{0}\kappa_{X}\left(\left[\mathfrak{z}\right]_{p}\right)\right)
\end{equation}

Meanwhile:
\begin{align}
\sum_{p^{2}\leq\left|t\right|_{p}\leq p^{N}}\hat{X}\left(t\right)e^{2\pi i\left\{ t\mathfrak{z}\right\} _{p}} & =\tilde{X}_{N}\left(\mathfrak{z}\right)-\tilde{X}_{1}\left(\mathfrak{z}\right)
\end{align}
Hence:
\begin{align*}
\sum_{\left|t\right|_{p}\leq p^{N}}\alpha_{X}\left(t\right)\hat{X}\left(pt\right)e^{2\pi i\left\{ t\mathfrak{z}\right\} _{p}} & =\sum_{p^{2}\leq\left|t\right|_{p}\leq p^{N}}\left(\hat{X}\left(t\right)+\beta_{X}\left(0\right)\hat{A}_{X}\left(t\right)\right)e^{2\pi i\left\{ t\mathfrak{z}\right\} _{p}}\\
 & \Updownarrow\\
a_{\left[\mathfrak{z}\right]_{p}}\tilde{X}_{N-1}\left(\theta_{p}\left(\mathfrak{z}\right)\right) & =\beta_{X}\left(0\right)\left(a_{0}^{N}\kappa_{X}\left(\left[\mathfrak{z}\right]_{p^{N}}\right)-a_{0}\kappa_{X}\left(\left[\mathfrak{z}\right]_{p}\right)\right)+\tilde{X}_{N}\left(\mathfrak{z}\right)-\tilde{X}_{1}\left(\mathfrak{z}\right)
\end{align*}
Here, we note that:
\begin{equation}
\kappa_{X}\left(\left[\mathfrak{z}\right]_{p}\right)=\frac{a_{\left[\mathfrak{z}\right]_{p}}}{a_{0}}
\end{equation}
which leaves us with:
\begin{equation}
a_{\left[\mathfrak{z}\right]_{p}}\tilde{X}_{N-1}\left(\theta_{p}\left(\mathfrak{z}\right)\right)=\beta_{X}\left(0\right)\left(a_{0}^{N}\kappa_{X}\left(\left[\mathfrak{z}\right]_{p^{N}}\right)-a_{\left[\mathfrak{z}\right]_{p}}\right)+\tilde{X}_{N}\left(\mathfrak{z}\right)-\tilde{X}_{1}\left(\mathfrak{z}\right)\label{eq:ready for solving}
\end{equation}

Finally:
\begin{align*}
\tilde{X}_{1}\left(\mathfrak{z}\right) & =\overbrace{\hat{X}\left(0\right)}^{0}+\sum_{\left|t\right|_{p}=p}\hat{X}\left(t\right)e^{2\pi i\left\{ t\mathfrak{z}\right\} _{p}}\\
 & =\sum_{\left|t\right|_{p}=p}\left(\beta_{X}\left(0\right)\overbrace{v_{p}\left(t\right)}^{-1}+\gamma_{X}\left(\frac{t\left|t\right|_{p}}{p}\right)\right)\hat{A}_{X}\left(t\right)e^{2\pi i\left\{ t\mathfrak{z}\right\} _{p}}\\
\left(\hat{A}_{X}\left(t\right)=\alpha_{X}\left(t\right),\textrm{ }\forall\left|t\right|_{p}=p\right); & =\sum_{\left|t\right|_{p}=p}\left(\gamma_{X}\left(t\right)-\beta_{X}\left(0\right)\right)\alpha_{X}\left(t\right)e^{2\pi i\left\{ t\mathfrak{z}\right\} _{p}}\\
 & =\sum_{\left|t\right|_{p}=p}\left(\beta_{X}\left(t\right)-\beta_{X}\left(0\right)\alpha_{X}\left(t\right)\right)e^{2\pi i\left\{ t\mathfrak{z}\right\} _{p}}\\
 & =-\beta_{X}\left(0\right)\overbrace{\left(1-\alpha_{X}\left(0\right)\right)}^{0}+\sum_{\left|t\right|_{p}\leq p}\left(\beta_{X}\left(t\right)-\beta_{X}\left(0\right)\alpha_{X}\left(t\right)\right)e^{2\pi i\left\{ t\mathfrak{z}\right\} _{p}}\\
\left(\alpha_{X}\left(0\right)=1\right); & =b_{\left[\mathfrak{z}\right]_{p}}-\beta_{X}\left(0\right)a_{\left[\mathfrak{z}\right]_{p}}
\end{align*}
Hence (\ref{eq:ready for solving}) becomes:
\begin{align*}
a_{\left[\mathfrak{z}\right]_{p}}\tilde{X}_{N-1}\left(\theta_{p}\left(\mathfrak{z}\right)\right) & =\beta_{X}\left(0\right)\left(a_{0}^{N}\kappa_{X}\left(\left[\mathfrak{z}\right]_{p^{N}}\right)-a_{\left[\mathfrak{z}\right]_{p}}\right)+\tilde{X}_{N}\left(\mathfrak{z}\right)-\left(b_{\left[\mathfrak{z}\right]_{p}}-\beta_{X}\left(0\right)a_{\left[\mathfrak{z}\right]_{p}}\right)\\
 & =\beta_{X}\left(0\right)a_{0}^{N}\kappa_{X}\left(\left[\mathfrak{z}\right]_{p^{N}}\right)+\tilde{X}_{N}\left(\mathfrak{z}\right)-b_{\left[\mathfrak{z}\right]_{p}}
\end{align*}
and so:
\begin{equation}
\tilde{X}_{N}\left(\mathfrak{z}\right)=a_{\left[\mathfrak{z}\right]_{p}}\tilde{X}_{N-1}\left(\theta_{p}\left(\mathfrak{z}\right)\right)+b_{\left[\mathfrak{z}\right]_{p}}-\beta_{X}\left(0\right)a_{0}^{N}\kappa_{X}\left(\left[\mathfrak{z}\right]_{p^{N}}\right)
\end{equation}
which is (\ref{eq:X_N twiddle recurrence formula}).

Q.E.D.
\begin{prop}
\label{prop:iterative triangle}Letting everything be as given in
\textbf{Proposition \ref{prop:X}}, we have:
\begin{equation}
\Delta_{N}^{\left(0\right)}\left\{ X\right\} \left(\mathfrak{z}\right)=\left(\Delta_{0}^{\left(N\right)}\left\{ X\right\} \left(\mathfrak{z}\right)+\left[\alpha_{X}\left(0\right)=1\right]N\beta_{X}\left(0\right)\right)a_{0}^{N}\kappa_{X}\left(\left[\mathfrak{z}\right]_{p^{N}}\right)\label{eq:iterative triangle}
\end{equation}
where:
\begin{equation}
\Delta_{N}^{\left(n\right)}\left\{ X\right\} \left(\mathfrak{z}\right)\overset{\textrm{def}}{=}X\left(\theta_{p}^{\circ n}\left(\mathfrak{z}\right)\right)-\tilde{X}_{N}\left(\theta_{p}^{\circ n}\left(\mathfrak{z}\right)\right)
\end{equation}
for all $n,N\geq0$.
\end{prop}
Proof: We subtract $X\left(\mathfrak{z}\right)$ from both sides of
equation (\ref{eq:X_N twiddle recurrence formula})\textbf{ Proposition
\ref{prop:X}}, and then use the identity:
\begin{equation}
X\left(\mathfrak{z}\right)=a_{\left[\mathfrak{z}\right]_{p}}X\left(\theta_{p}\left(\mathfrak{z}\right)\right)+b_{\left[\mathfrak{z}\right]_{p}}
\end{equation}
from \textbf{Proposition \ref{prop:shift reformulation}} to get:
\[
\underbrace{X\left(\mathfrak{z}\right)-\tilde{X}_{N}\left(\mathfrak{z}\right)}_{\Delta_{N}^{\left(0\right)}\left\{ X\right\} \left(\mathfrak{z}\right)}=a_{\left[\mathfrak{z}\right]_{p}}\underbrace{\left(X\left(\theta_{p}\left(\mathfrak{z}\right)\right)-\tilde{X}_{N-1}\left(\theta_{p}\left(\mathfrak{z}\right)\right)\right)}_{\Delta_{N-1}^{\left(1\right)}\left\{ X\right\} \left(\mathfrak{z}\right)}+\beta_{X}\left(0\right)\left[\alpha_{X}\left(0\right)=1\right]a_{0}^{N}\kappa_{X}\left(\left[\mathfrak{z}\right]_{p^{N}}\right)
\]
which is:
\begin{equation}
\Delta_{N}^{\left(0\right)}\left\{ X\right\} \left(\mathfrak{z}\right)=a_{\left[\mathfrak{z}\right]_{p}}\Delta_{N-1}^{\left(1\right)}\left\{ X\right\} \left(\mathfrak{z}\right)+h_{N}\left(\mathfrak{z}\right)\label{eq:induction}
\end{equation}
where:
\begin{equation}
h_{N}\left(\mathfrak{z}\right)\overset{\textrm{def}}{=}\beta_{X}\left(0\right)\left[\alpha_{X}\left(0\right)=1\right]a_{0}^{N}\kappa_{X}\left(\left[\mathfrak{z}\right]_{p^{N}}\right)
\end{equation}

Now, by induction on (\ref{eq:induction}), we have:
\begin{equation}
\Delta_{N}^{\left(0\right)}\left\{ X\right\} \left(\mathfrak{z}\right)=\Delta_{0}^{\left(N\right)}\left\{ X\right\} \left(\mathfrak{z}\right)\prod_{n=0}^{N-1}a_{\left[\theta_{p}^{\circ n}\left(\mathfrak{z}\right)\right]_{p}}+\sum_{n=0}^{N-1}h_{N-n}\left(\theta_{p}^{\circ n}\left(\mathfrak{z}\right)\right)\prod_{m=0}^{n-1}a_{\left[\theta_{p}^{\circ m}\left(\mathfrak{z}\right)\right]_{p}}\label{eq:induction of induction}
\end{equation}
By \textbf{Proposition \ref{prop:Product identity for kappa_X}},
we can write the products in (\ref{eq:induction of induction}) as:
\begin{equation}
\Delta_{N}^{\left(0\right)}\left\{ X\right\} \left(\mathfrak{z}\right)=\Delta_{0}^{\left(N\right)}\left\{ X\right\} \left(\mathfrak{z}\right)a_{0}^{N}\kappa_{X}\left(\left[\mathfrak{z}\right]_{p^{N}}\right)+\sum_{n=0}^{N-1}h_{N-n}\left(\theta_{p}^{\circ n}\left(\mathfrak{z}\right)\right)a_{0}^{n}\kappa_{X}\left(\left[\mathfrak{z}\right]_{p^{n}}\right)\label{eq:induction of induction simplified}
\end{equation}
Here, applying \textbf{Proposition \ref{prop:Kappa shift equation}},
we have that:
\begin{align*}
h_{N-n}\left(\theta_{p}^{\circ n}\left(\mathfrak{z}\right)\right) & =\beta_{X}\left(0\right)\left[\alpha_{X}\left(0\right)=1\right]a_{0}^{N-n}\kappa_{X}\left(\left[\theta_{p}^{\circ n}\left(\mathfrak{z}\right)\right]_{p^{N-n}}\right)\\
\left(\kappa_{X}\left(\left[\theta_{p}^{\circ n}\left(\mathfrak{z}\right)\right]_{p^{k}}\right)=\frac{\kappa_{X}\left(\left[\mathfrak{z}\right]_{p^{k+n}}\right)}{\kappa_{X}\left(\left[\mathfrak{z}\right]_{p^{n}}\right)}\right); & =\beta_{X}\left(0\right)\left[\alpha_{X}\left(0\right)=1\right]a_{0}^{N-n}\frac{\kappa_{X}\left(\left[\mathfrak{z}\right]_{p^{N}}\right)}{\kappa_{X}\left(\left[\mathfrak{z}\right]_{p^{n}}\right)}\\
 & =\beta_{X}\left(0\right)\left[\alpha_{X}\left(0\right)=1\right]\frac{a_{0}^{N}\kappa_{X}\left(\left[\mathfrak{z}\right]_{p^{N}}\right)}{a_{0}^{n}\kappa_{X}\left(\left[\mathfrak{z}\right]_{p^{n}}\right)}
\end{align*}
and so, (\ref{eq:induction of induction simplified}) becomes:
\begin{align*}
\Delta_{N}^{\left(0\right)}\left\{ X\right\} \left(\mathfrak{z}\right) & =\Delta_{0}^{\left(N\right)}\left\{ X\right\} \left(\mathfrak{z}\right)a_{0}^{N}\kappa_{X}\left(\left[\mathfrak{z}\right]_{p^{N}}\right)+\sum_{n=0}^{N-1}\beta_{X}\left(0\right)\left[\alpha_{X}\left(0\right)=1\right]\frac{a_{0}^{N}\kappa_{X}\left(\left[\mathfrak{z}\right]_{p^{N}}\right)}{a_{0}^{n}\kappa_{X}\left(\left[\mathfrak{z}\right]_{p^{n}}\right)}a_{0}^{n}\kappa_{X}\left(\left[\mathfrak{z}\right]_{p^{n}}\right)\\
 & =\Delta_{0}^{\left(N\right)}\left\{ X\right\} \left(\mathfrak{z}\right)a_{0}^{N}\kappa_{X}\left(\left[\mathfrak{z}\right]_{p^{N}}\right)+\left[\alpha_{X}\left(0\right)=1\right]N\beta_{X}\left(0\right)a_{0}^{N}\kappa_{X}\left(\left[\mathfrak{z}\right]_{p^{N}}\right)\\
 & =\left(\Delta_{0}^{\left(N\right)}\left\{ X\right\} \left(\mathfrak{z}\right)+\left[\alpha_{X}\left(0\right)=1\right]N\beta_{X}\left(0\right)\right)a_{0}^{N}\kappa_{X}\left(\left[\mathfrak{z}\right]_{p^{N}}\right)
\end{align*}

Q.E.D.

\vphantom{}Next, we deal with $\hat{f}_{\mathbf{n}}$ and all of
$X_{\mathbf{n}}$.
\begin{prop}
For all $N\geq1$ and all $\mathfrak{z}\in\mathbb{Z}_{p}$:
\begin{equation}
\tilde{f}_{\mathbf{n},N}\left(\mathfrak{z}\right)=r_{\mathbf{n},\left[\mathfrak{z}\right]_{p}}\tilde{f}_{\mathbf{n},N-1}\left(\theta_{p}\left(\mathfrak{z}\right)\right)+\sum_{\mathbf{m}<\mathbf{n}}r_{\mathbf{m},\mathbf{n},\left[\mathfrak{z}\right]_{p}}\left(\tilde{X}_{\mathbf{m},N-1}\left(\theta_{p}\left(\mathfrak{z}\right)\right)-\hat{X}_{\mathbf{m}}\left(0\right)\right)\label{eq:f_script J N twiddle formula}
\end{equation}
Hence, if \textbf{Assumption \ref{assu:main limit lemma}} holds:
\begin{align}
\Delta_{N}^{\left(0\right)}\left\{ f_{\mathbf{n}}\right\} \left(\mathfrak{z}\right) & =r_{\mathbf{n},0}^{N}\kappa_{\mathbf{n}}\left(\left[\mathfrak{z}\right]_{p^{N}}\right)\Delta_{0}^{\left(N\right)}\left\{ f_{\mathbf{n}}\right\} \left(\mathfrak{z}\right)\label{eq:f_script J N delta formula}\\
 & +\sum_{n=0}^{N-1}r_{\mathbf{n},0}^{n}\kappa_{\mathbf{n}}\left(\left[\mathfrak{z}\right]_{p^{n}}\right)\sum_{\mathbf{m}<\mathbf{n}}r_{\mathbf{m},\mathbf{n},\left[\theta_{p}^{\circ n}\left(\mathfrak{z}\right)\right]_{p}}\Delta_{N-1-n}^{\left(n+1\right)}\left\{ X_{\mathbf{m}}\right\} \left(\mathfrak{z}\right)\nonumber 
\end{align}
\end{prop}
Proof: By \textbf{Proposition \ref{prop:partial sum of f-J-hat fourier series}},
we have the formula:
\begin{equation}
\tilde{f}_{\mathbf{n},N}\left(\mathfrak{z}\right)=\sum_{n=0}^{N-2}r_{\mathbf{n},0}^{n}\kappa_{\mathbf{n}}\left(\left[\mathfrak{z}\right]_{p^{n}}\right)\sum_{\mathbf{m}<\mathbf{n}}r_{\mathbf{m},\mathbf{n},\left[\theta_{p}^{\circ n}\left(\mathfrak{z}\right)\right]_{p}}\left(\tilde{X}_{\mathbf{m},N-n-1}\left(\theta_{p}^{\circ n+1}\left(\mathfrak{z}\right)\right)-\hat{X}_{\mathbf{m}}\left(0\right)\right)
\end{equation}
Direct computation shows:
\begin{align*}
\tilde{f}_{\mathbf{n},N}\left(p\mathfrak{z}+k\right) & =r_{\mathbf{n},k}\underbrace{\sum_{n=1}^{N-3}r_{\mathbf{n},0}^{n}\kappa_{\mathbf{n}}\left(\left[\mathfrak{z}\right]_{p^{n}}\right)\sum_{\mathbf{m}<\mathbf{n}}r_{\mathbf{m},\mathbf{n},\left[\theta_{p}^{\circ n}\left(\mathfrak{z}\right)\right]_{p}}\left(\tilde{X}_{\mathbf{m},N-1-n-1}\left(\theta_{p}^{\circ n+1}\left(\mathfrak{z}\right)\right)-\hat{X}_{\mathbf{m}}\left(0\right)\right)}_{\tilde{f}_{\mathbf{n},N-1}\left(\mathfrak{z}\right)}\\
 & +\sum_{\mathbf{m}<\mathbf{n}}r_{\mathbf{m},\mathbf{n},\left[k\right]_{p}}\left(\tilde{X}_{\mathbf{m},N-1}\left(\mathfrak{z}\right)-\hat{X}_{\mathbf{m}}\left(0\right)\right)
\end{align*}
and hence:
\begin{equation}
\tilde{f}_{\mathbf{n},N}\left(\mathfrak{z}\right)=r_{\mathbf{n},\left[\mathfrak{z}\right]_{p}}\tilde{f}_{\mathbf{n},N-1}\left(\theta_{p}\left(\mathfrak{z}\right)\right)+\sum_{\mathbf{m}<\mathbf{n}}r_{\mathbf{m},\mathbf{n},\left[\mathfrak{z}\right]_{p}}\left(\tilde{X}_{\mathbf{m},N-1}\left(\theta_{p}\left(\mathfrak{z}\right)\right)-\hat{X}_{\mathbf{m}}\left(0\right)\right)
\end{equation}
which proves (\ref{eq:f_script J N twiddle formula}).

Assuming \textbf{Assumption \ref{assu:main limit lemma}}, we can
take $N\rightarrow\infty$ and get:
\begin{equation}
f_{\mathbf{n}}\left(\mathfrak{z}\right)=r_{\mathbf{n},\left[\mathfrak{z}\right]_{p}}f_{\mathbf{n}}\left(\theta_{p}\left(\mathfrak{z}\right)\right)+\sum_{\mathbf{m}<\mathbf{n}}r_{\mathbf{m},\mathbf{n},\left[\mathfrak{z}\right]_{p}}\left(X_{\mathbf{m}}\left(\theta_{p}\left(\mathfrak{z}\right)\right)-\hat{X}_{\mathbf{m}}\left(0\right)\right)
\end{equation}
Subtracting this from (\ref{eq:f_script J N twiddle formula}) gives:
\begin{align*}
\Delta_{N}^{\left(0\right)}\left\{ f_{\mathbf{n}}\right\} \left(\mathfrak{z}\right) & =f_{\mathbf{n}}\left(\mathfrak{z}\right)-\tilde{f}_{\mathbf{n},N}\left(\mathfrak{z}\right)\\
 & =\left(r_{\mathbf{n},\left[\mathfrak{z}\right]_{p}}f_{\mathbf{n}}\left(\theta_{p}\left(\mathfrak{z}\right)\right)+\sum_{\mathbf{m}<\mathbf{n}}r_{\mathbf{m},\mathbf{n},\left[\mathfrak{z}\right]_{p}}\left(X_{\mathbf{m}}\left(\theta_{p}\left(\mathfrak{z}\right)\right)-\hat{X}_{\mathbf{m}}\left(0\right)\right)\right)\\
 & -\left(r_{\mathbf{n},\left[\mathfrak{z}\right]_{p}}\tilde{f}_{\mathbf{n},N-1}\left(\theta_{p}\left(\mathfrak{z}\right)\right)+\sum_{\mathbf{m}<\mathbf{n}}r_{\mathbf{m},\mathbf{n},\left[\mathfrak{z}\right]_{p}}\left(\tilde{X}_{\mathbf{m},N-1}\left(\theta_{p}\left(\mathfrak{z}\right)\right)-\hat{X}_{\mathbf{m}}\left(0\right)\right)\right)\\
 & =r_{\mathbf{n},\left[\mathfrak{z}\right]_{p}}\left(f_{\mathbf{n}}\left(\theta_{p}\left(\mathfrak{z}\right)\right)-\tilde{f}_{\mathbf{n},N-1}\left(\theta_{p}\left(\mathfrak{z}\right)\right)\right)\\
 & +\sum_{\mathbf{m}<\mathbf{n}}r_{\mathbf{m},\mathbf{n},\left[\mathfrak{z}\right]_{p}}\left(X_{\mathbf{m}}\left(\theta_{p}\left(\mathfrak{z}\right)\right)-\tilde{X}_{\mathbf{m},N-1}\left(\theta_{p}\left(\mathfrak{z}\right)\right)\right)\\
 & =r_{\mathbf{n},\left[\mathfrak{z}\right]_{p}}\Delta_{N-1}^{\left(1\right)}\left\{ f_{\mathbf{n}}\right\} \left(\mathfrak{z}\right)+\sum_{\mathbf{m}<\mathbf{n}}r_{\mathbf{m},\mathbf{n},\left[\mathfrak{z}\right]_{p}}\Delta_{N-1}^{\left(1\right)}\left\{ X_{\mathbf{m}}\right\} \left(\mathfrak{z}\right)
\end{align*}
So, by induction:
\begin{align*}
\Delta_{N}^{\left(0\right)}\left\{ f_{\mathbf{n}}\right\} \left(\mathfrak{z}\right) & =r_{\mathbf{n},\left[\mathfrak{z}\right]_{p}}\Delta_{N-1}^{\left(1\right)}\left\{ f_{\mathbf{n}}\right\} \left(\mathfrak{z}\right)+\sum_{\mathbf{m}<\mathbf{n}}r_{\mathbf{m},\mathbf{n},\left[\mathfrak{z}\right]_{p}}\Delta_{N-1}^{\left(1\right)}\left\{ X_{\mathbf{m}}\right\} \left(\left[\mathfrak{z}\right]_{p}\right)\\
\left(\textrm{iterate once}\right); & =r_{\mathbf{n},\left[\mathfrak{z}\right]_{p}}r_{\mathbf{n},\left[\theta_{p}\left(\mathfrak{z}\right)\right]_{p}}\Delta_{N-2}^{\left(2\right)}\left\{ f_{\mathbf{n}}\right\} \left(\mathfrak{z}\right)\\
 & +r_{\mathbf{n},\left[\mathfrak{z}\right]_{p}}\sum_{\mathbf{m}<\mathbf{n}}r_{\mathbf{m},\mathbf{n},\left[\theta_{p}\left(\mathfrak{z}\right)\right]_{p}}\Delta_{N-2}^{\left(2\right)}\left\{ X_{\mathbf{m}}\right\} \left(\mathfrak{z}\right)+\sum_{\mathbf{m}<\mathbf{n}}r_{\mathbf{m},\mathbf{n},\left[\mathfrak{z}\right]_{p}}\Delta_{N-1}^{\left(1\right)}\left\{ X_{\mathbf{m}}\right\} \left(\mathfrak{z}\right)\\
 & \vdots\\
\left(\textrm{iterate }M\textrm{ times}\right); & =\Delta_{N-M}^{\left(M\right)}\left\{ f_{\mathbf{n}}\right\} \left(\mathfrak{z}\right)\prod_{m=0}^{M-1}r_{\mathbf{n},\left[\theta_{p}^{\circ m}\left(\mathfrak{z}\right)\right]_{p}}\\
 & +\sum_{m=0}^{M-1}\left(\prod_{n=0}^{m-1}r_{\mathbf{n},\left[\theta_{p}^{\circ n}\left(\mathfrak{z}\right)\right]_{p}}\right)\sum_{\mathbf{m}<\mathbf{n}}r_{\mathbf{m},\mathbf{n},\left[\theta_{p}^{\circ m}\left(\mathfrak{z}\right)\right]_{p}}\Delta_{N-1-m}^{\left(m+1\right)}\left\{ X_{\mathbf{m}}\right\} \left(\mathfrak{z}\right)
\end{align*}
Using \textbf{Proposition \ref{prop:Product identity for kappa_X}},
we can simplify this to:
\begin{align*}
\Delta_{N}^{\left(0\right)}\left\{ f_{\mathbf{n}}\right\} \left(\mathfrak{z}\right) & =\Delta_{N-M}^{\left(M\right)}\left\{ f_{\mathbf{n}}\right\} \left(\mathfrak{z}\right)r_{\mathbf{n},0}^{M}\kappa_{\mathbf{n}}\left(\left[\mathfrak{z}\right]_{p^{M}}\right)\\
 & +\sum_{m=0}^{M-1}r_{\mathbf{n},0}^{m}\kappa_{\mathbf{n}}\left(\left[\mathfrak{z}\right]_{p^{m}}\right)\sum_{\mathbf{m}<\mathbf{n}}r_{\mathbf{m},\mathbf{n},\left[\theta_{p}^{\circ m}\left(\mathfrak{z}\right)\right]_{p}}\Delta_{N-1-m}^{\left(m+1\right)}\left\{ X_{\mathbf{m}}\right\} \left(\mathfrak{z}\right)
\end{align*}
Setting $M=N$ gives:
\begin{align*}
\Delta_{N}^{\left(0\right)}\left\{ f_{\mathbf{n}}\right\} \left(\mathfrak{z}\right) & =r_{\mathbf{n},0}^{N}\kappa_{\mathbf{n}}\left(\left[\mathfrak{z}\right]_{p^{N}}\right)\Delta_{0}^{\left(N\right)}\left\{ f_{\mathbf{n}}\right\} \left(\mathfrak{z}\right)\\
 & +\sum_{n=0}^{N-1}r_{\mathbf{n},0}^{n}\kappa_{\mathbf{n}}\left(\left[\mathfrak{z}\right]_{p^{n}}\right)\sum_{\mathbf{m}<\mathbf{n}}r_{\mathbf{m},\mathbf{n},\left[\theta_{p}^{\circ n}\left(\mathfrak{z}\right)\right]_{p}}\Delta_{N-1-n}^{\left(n+1\right)}\left\{ X_{\mathbf{m}}\right\} \left(\mathfrak{z}\right)
\end{align*}

Q.E.D.
\begin{rem}
\label{rem:end of big idea remark}Using the symbology of \textbf{Proposition
\ref{prop:X}}, the $\Sigma\left(\mathbf{n}\right)=1$ case of the
functional equation (\ref{eq:formal solution}) dealt with in \textbf{Example
\ref{exa:J has one element}} can be written as:
\begin{equation}
\hat{X}\left(t\right)=\alpha_{X}\left(t\right)\hat{X}\left(pt\right)+\beta_{X}\left(t\right)\mathbf{1}_{0}\left(pt\right)\label{eq:formal solution J equals 1 case with alpha and beta}
\end{equation}
where $\hat{X}$ is used in place of $\hat{X}_{\mathbf{n}}$, where
$X_{\mathbf{0}}=\mathbf{1}_{0}$, where $\alpha_{X}=\alpha_{\mathbf{n}}$
and $\beta_{X}=\alpha_{\mathbf{0},\mathbf{n}}$. Solving for $\alpha_{X}\left(t\right)\hat{X}\left(pt\right)$
gives:
\begin{equation}
\alpha_{X}\left(t\right)\hat{X}\left(pt\right)=\begin{cases}
\hat{X}\left(t\right)-\beta_{X}\left(t\right) & \textrm{if }\left|t\right|_{p}\leq p\\
\hat{X}\left(t\right) & \textrm{if }\left|t\right|_{p}>p
\end{cases}\label{eq:formal solution J equals 1 case with alpha and beta rewritten}
\end{equation}
Meanwhile, the result of \textbf{Proposition \ref{que:key question (hard part)}
}is:
\begin{equation}
\alpha_{X}\left(t\right)\hat{X}\left(pt\right)=\begin{cases}
\begin{cases}
\alpha_{X}\left(t\right)\hat{X}\left(0\right) & \textrm{if }\left|t\right|_{p}\leq p\\
\hat{X}\left(t\right) & \textrm{if }\left|t\right|_{p}>p
\end{cases} & \textrm{if }\alpha_{X}\left(0\right)\neq1\\
\\
\begin{cases}
0 & \textrm{if }\left|t\right|_{p}\leq p\\
\hat{X}\left(t\right)+\beta_{X}\left(0\right)\hat{A}_{X}\left(t\right) & \textrm{if }\left|t\right|_{p}>p
\end{cases} & \textrm{if }\alpha_{X}\left(0\right)=1
\end{cases}\label{eq:Corrected equation}
\end{equation}
We will call (\ref{eq:Corrected equation}) the\textbf{ corrected
equation}, while (\ref{eq:formal solution J equals 1 case with alpha and beta rewritten})
will be called the \textbf{uncorrected equation}.\textbf{ }Thus, we
have the following analysis:

\textbullet{} Suppose $\alpha_{X}\left(0\right)\neq1$. When $t=0$,
the uncorrected equation is: 
\begin{equation}
\alpha_{X}\left(0\right)\hat{X}\left(0\right)=\hat{X}\left(0\right)-\beta_{X}\left(0\right)
\end{equation}
Hence:
\begin{equation}
\frac{\beta_{X}\left(0\right)}{1-\alpha_{X}\left(0\right)}=X\left(0\right)
\end{equation}
On the other hand, for $t=0$, the corrected equation is $\alpha_{X}\left(0\right)\hat{X}\left(0\right)=\alpha_{X}\left(0\right)\hat{X}\left(0\right)$,
which is a tautology. Thus, the two equations are \emph{not in conflict
}for $t=0$. That is, solutions of one are not automatically barred
from being solutions of the other.

When $\left|t\right|_{p}=p$, the uncorrected equation is:
\begin{equation}
\alpha_{X}\left(t\right)\hat{X}\left(0\right)=\hat{X}\left(t\right)-\beta_{X}\left(t\right)
\end{equation}
while the corrected equation is once again a tautology.
\begin{equation}
\alpha_{X}\left(t\right)\hat{X}\left(0\right)=\alpha_{X}\left(t\right)\hat{X}\left(0\right)
\end{equation}
Thus, the two equations \emph{not in conflict }for $\left|t\right|_{p}=p$,
either.

Finally, when $\left|t\right|_{p}>p$, both equations are identical:
$\alpha_{X}\left(t\right)\hat{X}\left(pt\right)=\hat{X}\left(t\right)$.

Consequently, \emph{the corrected and uncorrected equations are not
in conflict with one another when $\alpha_{X}\left(0\right)\neq1$.
}In particular, this shows that, in the $\alpha_{X}\left(0\right)\neq1$
case, \textbf{Theorem \ref{thm:The-breakdown-variety}}'s formula
for $\hat{X}$ is precisely the  unique solution of the uncorrected
equation (\ref{eq:formal solution J equals 1 case with alpha and beta}).
As such, the formal procedure for computing $\hat{X}$ described in
\textbf{Example \ref{exa:proceeding formally}} is justified in the
case where $\Sigma\left(\mathbf{n}\right)=1$ and $\alpha_{X}\left(0\right)\neq1$.

\textbullet{} Suppose $\alpha_{X}\left(0\right)=1$. When $\left|t\right|_{p}>p$,
the uncorrected equation is:
\begin{equation}
\alpha_{X}\left(t\right)\hat{X}\left(pt\right)=\hat{X}\left(t\right)
\end{equation}
while the corrected equation is:
\begin{equation}
\alpha_{X}\left(t\right)\hat{X}\left(pt\right)=\hat{X}\left(t\right)+\beta_{X}\left(0\right)\hat{A}_{X}\left(t\right)
\end{equation}
Since $\hat{A}_{X}$ is not identically $0$, these two are consistent
if and only if $\beta_{X}\left(0\right)=0$. When $\beta_{X}\left(0\right)\neq0$,
on the other hand, the corrected equation will have solutions, while\textemdash as
we saw in \textbf{Proposition \ref{prop:X_3-hat formal solution}},
the uncorrected equation will have no solutions!

When $t=0$, the uncorrected equation forces $\beta_{X}\left(0\right)=0$,
while the corrected one forces $\hat{X}\left(0\right)=0$ and leaves
$\beta_{X}\left(0\right)$ unaffected. On the other hand, as we saw
in \textbf{Proposition \ref{prop:X_3-hat formal solution}}, when
$\beta_{X}\left(0\right)=0$, the uncorrected equation has infinitely
many solutions, each being uniquely determined by the initial value
$\hat{X}\left(0\right)$, which may be chosen freely.

This gives the complete answer to \textbf{Question \ref{que:key question (hard part)}}
in the case where $\Sigma\left(\mathbf{n}\right)=1$:

I. If $\alpha_{X}\left(0\right)\neq1$, then (\ref{eq:formal solution J equals 1 case with alpha and beta})
(our simplified version of equation (\ref{eq:formal solution}) for
the $\Sigma\left(\mathbf{n}\right)=1$ case) has a unique solution,
and that solution is a Fourier transform of $X$.

II. If $\alpha_{X}\left(0\right)=1$ and $\beta_{X}\left(0\right)=0$,
of the infinitely many solutions that (\ref{eq:formal solution J equals 1 case with alpha and beta})
possesses, only the one with the initial condition $\hat{X}\left(0\right)=0$
is a Fourier transform of $X$.

III. If $\alpha_{X}\left(0\right)=1$ and $\beta_{X}\left(0\right)\neq0$,
(\ref{eq:formal solution J equals 1 case with alpha and beta}) has
no solutions, and all of these non-existent solutions are not Fourier
transform of $X$.
\end{rem}

\section{\label{sec:The-Small-Ideas}The Smaller Ideas}

With the main computation out of the way, we can now turn to introducing
the concepts and set-up needed to rigorously justify our computations.

\subsection{\label{subsec:Frame-Theory-=000026}Frame Theory}

Intuitively, a frame device for recording unusual recipes for interpreting
sequences of functions as ``convergent''. The following is an archetypical
example of the situation that motivated them.
\begin{example}
\label{exa:Let-,-where}Let $p,q,d\in\mathbb{Z}\backslash\left\{ 0,1,-1\right\} $,
where $p$ and $q$ are distinct primes and $d$ is co-prime to both
$p$ and $q$. Then, consider the $3$-adic F-series:
\begin{equation}
X\left(\mathfrak{z}\right)=\sum_{n=0}^{\infty}\frac{p^{\#_{3:1}\left(\left[\mathfrak{z}\right]_{3^{n}}\right)}q^{\#_{3:2}\left(\left[\mathfrak{z}\right]_{3^{n}}\right)}}{d^{n}}\label{eq:3-adic example}
\end{equation}
Observe that for all non-negative integer values of $\mathfrak{z}$,
$X\left(\mathfrak{z}\right)$ converges in $\mathbb{R}$ to a rational
number. Meanwhile, for any $\mathfrak{z}\in\mathbb{Z}_{3}\backslash\left\{ 0,1,2,\ldots\right\} $,
the series will converge:

\textbullet{} $p$-adically if $\mathfrak{z}$ has infinitely many
$1$s digits;

\textbullet{} $q$-adically if $\mathfrak{z}$ has infinitely many
$2$s digits;

\textbullet{} both $p$-adically \emph{and} $q$-adically if $\mathfrak{z}$
has infinitely many $1$s and $2$s digits.

\textbullet{} $\infty$-adically, depending on the value of $\mathfrak{z}$.
Indeed, by the \textbf{root test}, we have that:
\begin{equation}
\limsup_{n\rightarrow\infty}\left(\frac{\#_{3:1}\left(\left[\mathfrak{z}\right]_{3^{n}}\right)}{n}\frac{\ln\left|p\right|}{\ln\left|d\right|}+\frac{\#_{3:2}\left(\left[\mathfrak{z}\right]_{3^{n}}\right)}{n}\frac{\ln\left|q\right|}{\ln\left|d\right|}\right)<1
\end{equation}
is a sufficient, though not necessary, condition for $X\left(\mathfrak{z}\right)$
to converge in $\mathbb{R}$.

Thus, for example, we could realize $X$ as a function $\mathbb{Z}_{3}\rightarrow\mathbb{Q}_{p}\cup\mathbb{Q}_{q}\cup\mathbb{Q}$
by defining the sum of (\ref{eq:3-adic example}) using the real topology
if $\mathfrak{z}\in\left\{ 0,1,2,\ldots\right\} $, using the $p$-adic
topology if $\mathfrak{z}$ has infinitely many $1$s digits but finitely
many $2$s digits, and using the $q$-adic topology if $\mathfrak{z}$
has infinitely many $2$s digits. This recipe for summing (\ref{eq:3-adic example})
is an example of a frame. Obviously, there are infinitely many possible
choices of a frame that makes $X$ summable at every $3$-adic integer
$\mathfrak{z}$, and the flexibility of these choices encodes arithmetic
information about (\ref{eq:3-adic example}), despite their seemingly
ad-hoc nature. Indeed, as we will see, the everything that we will
do with F-series will work for any choice of frame satisfying the
necessary regularity hypotheses.
\end{example}
As the theory of frames is in its infancy, there is still much to
be done to develop it. One long term goal is to remove the ad-hoc
aspects of the construction and turn it into a means of capturing
arithmetic data.

To define frames in general, we need to use norms and seminorms on
rings. These, recall, were defined in \textbf{Definition \ref{def:A-seminorm-on}}
on page \pageref{def:A-seminorm-on} of \textbf{Section \ref{subsec:Preliminaries-(Notation,-etc.)}}.
For brevity, we use the following notation:
\begin{defn}
Let $R$ be a ring. I write:

\textbullet{} $\mathscr{S}\left(R\right)$, to denote the set of all
seminorms on $R$.

\textbullet{} $\mathscr{N}\left(R\right)$, to denote the set of all
norms on $R$.

\textbullet{} $\mathscr{V}\left(R\right)$, to denote the set of all
absolute values on $R$.

\textbullet{} $\mathscr{M}\left(R\right)$, to denote the set of all
multiplicative seminorms on $R$; this is known as the \textbf{Berkovich
spectrum of $R$}.
\end{defn}
Here, then, is the definition of a frame:
\begin{defn}
Let $\mathcal{X}$ be a set and let $R$ be a ring. Then, an \textbf{$R$-frame
on $\mathcal{X}$, }denoted $\mathcal{F}$, consists of the following
data:

I. A set $D\left(\mathcal{F}\right)\subseteq\mathcal{X}$, called
the \textbf{domain} of $\mathcal{F}$. (Many times, we will have $D\left(\mathcal{F}\right)=\mathcal{X}$,
but there will be cases where $\mathcal{X}$ is a measure space and
$\mathcal{X}\backslash D\left(\mathcal{F}\right)$ is a set of measure
zero.)

II. A map $\mathcal{F}:D\left(\mathcal{F}\right)\rightarrow\mathscr{S}\left(R^{D\left(\mathcal{F}\right)}\right)$
which sends each $x\in D\left(\mathcal{F}\right)$ to a seminorm $\mathcal{F}_{x}:R^{D\left(\mathcal{F}\right)}\rightarrow\left[0,\infty\right)$,
where $R^{D\left(\mathcal{F}\right)}$ is the ring of $R$-valued
functions on $D\left(\mathcal{F}\right)$.

We say a frame is \textbf{multiplicative} when $\mathcal{F}_{x}\in\mathscr{M}\left(R^{D\left(\mathcal{F}\right)}\right)$
for all $x\in D\left(\mathcal{F}\right)$. We say $\mathcal{F}$ is
\textbf{evaluative }whenever, for each $x\in D\left(\mathcal{F}\right)$,
there is an absolute value $\left|\cdot\right|_{\mathcal{F}\left(x\right)}\in\mathscr{V}\left(R\right)$
so that:
\begin{equation}
\mathcal{F}_{x}\left(f\right)=\left|f\left(x\right)\right|_{\mathcal{F}\left(x\right)},\textrm{ }\forall f\in R^{D\left(\mathcal{F}\right)}
\end{equation}
\end{defn}
\begin{assumption}
HENCEFORTH, UNLESS STATED OTHERWISE, ALL FRAMES WILL BE MULTIPLICATIVE
AND EVALUATIVE.
\end{assumption}
\begin{rem}
Aside from making things easier overall, the multiplicativity of $\mathcal{F}_{x}$
guarantees that the only elements $f\in R^{D\left(\mathcal{F}\right)}$
with $\mathcal{F}_{x}\left(f\right)=1$ are the units of $R^{D\left(\mathcal{F}\right)}$.
Moreover, viewing elements of $R$ as constant $R$-valued functions
on $D\left(\mathcal{F}\right)$, we then have that, for any $r\in R$,
$\mathcal{F}_{x}\left(r\right)\neq1$ implies $r\notin R^{\times}$.
We'll discuss the assumption of the evaluativity of $\mathcal{F}$
momentarily.
\end{rem}
As $\mathcal{F}_{x}$ is a seminorm, it is not guaranteed that $f$
will be the zero element of $R^{D\left(\mathcal{F}\right)}$ simply
when $\mathcal{F}_{x}\left(f\right)=0$; that is, $\textrm{Ker}\mathcal{F}_{x}$
is not guaranteed to be trivial. However, since $\mathcal{F}$ is
evaluative, we see that $\textrm{Ker}\mathcal{F}_{x}$ is precisely
the set of $f\in R^{D\left(\mathcal{F}\right)}$ which vanish at $x$.
Because of this, the ring quotient $R^{D\left(\mathcal{F}\right)}/\textrm{Ker}\mathcal{F}_{x}$
is isomorphic to $R$ itself.
\begin{defn}
Let $\mathcal{F}$ be an evaluative $R$-frame on $\mathcal{X}$.

I. For each $x\in D\left(\mathcal{F}\right)$, let $R_{\mathcal{F}\left(x\right)}$
denote the metric completion of $R$ with respect to the absolute
value $\left|\cdot\right|_{\mathcal{F}\left(x\right)}$. \emph{\uline{Note}}:
for compactness of notation, we will often write $\mathcal{F}\left(x\right)$
to denote $R_{\mathcal{F}\left(x\right)}$.

II. We write $C\left(\mathcal{F}\right)$\textbf{ }to denote the \textbf{space
of functions} \textbf{compatible with $\mathcal{F}$}.\textbf{ }This
is defined by the direct product:
\begin{equation}
C\left(\mathcal{F}\right)\overset{\textrm{def}}{=}\prod_{x\in D\left(\mathcal{F}\right)}R_{\mathcal{F}\left(x\right)}
\end{equation}
We view elements of $C\left(\mathcal{F}\right)$ as functions $f$
on $D\left(\mathcal{F}\right)$ with the property that $f\left(x\right)\in R_{\mathcal{F}\left(x\right)}$
for all $x\in X$. We can make this view rigorous by defining $I\left(\mathcal{F}\right)$,
called the \textbf{image of $\mathcal{F}$}, as:
\begin{equation}
I\left(\mathcal{F}\right)\overset{\textrm{def}}{=}\bigcup_{x\in D\left(\mathcal{F}\right)}R_{\mathcal{F}\left(x\right)}
\end{equation}
Then, we have: 
\begin{equation}
C\left(\mathcal{F}\right)=\left\{ f:D\left(\mathcal{F}\right)\rightarrow I\left(\mathcal{F}\right)\mid f\left(x\right)\in R_{\mathcal{F}\left(x\right)},\textrm{ }\forall x\in D\left(\mathcal{F}\right)\right\} 
\end{equation}
Observe that $C\left(\mathcal{F}\right)$ is then an $R$-algebra
under scalar multiplication, with pointwise addition and multiplication.

Note that because $\mathcal{F}$ is assumed to be evaluative, we can
make sense of $\mathcal{F}_{x}\left(f\right)$ for all $f\in C\left(\mathcal{F}\right)$;
this is simply $\left|f\left(x\right)\right|_{\mathcal{F}\left(x\right)}$,
where, in an acceptable abuse of notation, $\left|\cdot\right|_{\mathcal{F}\left(x\right)}$
here denotes the extension of the absolute value of that name on $R$
to one on the completion $R_{\mathcal{F}\left(x\right)}$.
\end{defn}
\begin{rem}
In this paper, we use the hypothesis that $\mathcal{F}$ is evaluative
in order to guarantee that $R^{D\left(\mathcal{F}\right)}/\textrm{Ker}\mathcal{F}_{x}\cong R$.
Without assuming the evaluativity of $\mathcal{F}$, all of the definitions
given above hold verbatim, with the exception that instead of denoting
an absolute value on $R$, the symbol $\left|\cdot\right|_{\mathcal{F}\left(x\right)}$
instead denotes the absolute value induced by $\mathcal{F}_{x}$ on
the quotient $R^{D\left(\mathcal{F}\right)}/\textrm{Ker}\mathcal{F}_{x}$,
which in the general case need not be isomorphic to $R$. Consequently,
$R_{\mathcal{F}\left(x\right)}$ is not the completion of $R$ with
respect to $\left|\cdot\right|_{\mathcal{F}\left(x\right)}$, but
the completion of $R^{D\left(\mathcal{F}\right)}/\textrm{Ker}\mathcal{F}_{x}$.
Provided that modification is made, everything else in the above follows
without change.

In terms of how it affects the results of this paper, if we remove
the assumption that $\mathcal{F}$ is evaluative, everything we are
doing will hold provided that we make the assumption that for every
$f\in C\left(\mathcal{F}\right)$ and every $\epsilon>0$, we have
a non-empty intersection:
\begin{equation}
\bigcap_{x\in D\left(\mathcal{F}\right)}\pi_{x}^{-1}\left(U_{x}\left(f;\epsilon\right)\right)\neq\varnothing
\end{equation}
Here, $\pi_{x}:R^{D\left(\mathcal{F}\right)}\rightarrow R^{D\left(\mathcal{F}\right)}/\textrm{Ker}\mathcal{F}_{x}$
is the canonical projection, while:
\begin{equation}
U_{x}\left(f;\epsilon\right)\overset{\textrm{def}}{=}\left\{ r\in R^{D\left(\mathcal{F}\right)}/\textrm{Ker}\mathcal{F}_{x}:\left|\pi_{x}\left(r\right)-f\left(x\right)\right|_{\mathcal{F}\left(x\right)}<\epsilon\right\} 
\end{equation}
where, as mentioned above, $\left|\cdot\right|_{\mathcal{F}\left(x\right)}$
is the absolute value induced by $\mathcal{F}_{x}$ on $R^{D\left(\mathcal{F}\right)}/\textrm{Ker}\mathcal{F}_{x}$.
Also, some of the notation given below will need to be modified slightly.
I intend to publish a paper on these and other details of frame theory
in the future.
\end{rem}
\begin{defn}
Let $\mathcal{F}$ be an evaluative frame. We say a sequence $\left\{ f_{n}\right\} _{n\geq0}\in C\left(\mathcal{F}\right)$
\textbf{converges to $f\in C\left(\mathcal{F}\right)$ with respect
to $\mathcal{F}$} (a.k.a., the \textbf{$f_{n}$s $\mathcal{F}$-converge
to $\mathcal{F}$}; the \textbf{$f_{n}$s converge frame-wise to $f$})
whenever:
\begin{equation}
\lim_{n\rightarrow\infty}\mathcal{F}_{x}\left(f_{n}-f\right)\overset{\mathbb{R}}{=}0,\textrm{ }\forall x\in D\left(\mathcal{F}\right)
\end{equation}
We denote this by writing:
\begin{equation}
\lim_{n\rightarrow\infty}f_{n}\left(x\right)\overset{\mathcal{F}}{=}f\left(x\right)
\end{equation}

Also, we say a series of compatible functions $\sum_{n=0}^{N-1}f_{n}$
is \textbf{absolutely convergent with respect to $\mathcal{F}$ }/
is \textbf{absolutely $\mathcal{F}$-convergent} whenever:

\begin{equation}
\lim_{N\rightarrow\infty}\sum_{n=0}^{N-1}\mathcal{F}_{x}\left(f_{n}\right)
\end{equation}
converges in $\mathbb{R}$ for all $x\in D\left(\mathcal{F}\right)$.
\end{defn}
When working with $\mathbb{Z}_{p}$, we will want our frames to be
compatible with both the profinite structure of $\mathbb{Z}_{p}$
and the action of the shift operator $\theta_{p}$ on $\mathbb{Z}_{p}$.
\begin{defn}
Let $\mathcal{F}$ be an $R$-frame on $\mathbb{Z}_{p}$. We say $\mathcal{F}$
is \textbf{reffinite }whenever:

\vphantom{}I. $\mathbb{N}_{0}\subseteq D\left(\mathcal{F}\right)$.

\vphantom{}II. $\theta_{p}\left(D\left(\mathcal{F}\right)\right)\subseteq D\left(\mathcal{F}\right)$.

\vphantom{}III. For any $\mathfrak{z}\in D\left(\mathcal{F}\right)$,
the absolute values $\left|\cdot\right|_{\mathcal{F}\left(\mathfrak{z}\right)}$
and $\left|\cdot\right|_{\mathcal{F}\left(\theta_{p}\left(\mathfrak{z}\right)\right)}$
are the same.
\end{defn}
\begin{assumption}
HENCEFORTH, UNLESS STATED OTHERWISE, ALL FRAMES ARE ASSUMED TO BE
REFFINITE.
\end{assumption}
We will also need to explicate how frames interact with quotients
of the underlying ring.
\begin{defn}
Let $\mathcal{F}$ be an $R$-frame on $\mathcal{X}$, and let $I$
be an ideal of $R$. Then, $\mathcal{F}$ induces a canonical $R/I$-frame
on $\mathcal{X}$, denoted $\mathcal{F}/I$, by the rule that, for
each $x\in D\left(\mathcal{F}\right)$, $\left(\mathcal{F}/I\right)_{x}$
is the multiplicative seminorm on $\left(R/I\right)^{D\left(\mathcal{F}\right)}$defined
by: 
\begin{equation}
\left(\mathcal{F}/I\right)_{x}\left(f\right)\overset{\textrm{def}}{=}\inf_{g\in I^{D\left(\mathcal{F}\right)}}\mathcal{F}_{x}\left(f+g\right)
\end{equation}
We call this the \textbf{quotient of $\mathcal{F}$ by $I$}, and
refer to $\mathcal{F}/I$ as a \textbf{quotient frame}.
\end{defn}
Next, we have the notion of rising-continuity.
\begin{defn}
Let $\mathcal{F}$ be a $R$-frame on $\mathbb{Z}_{p}$ so that $\mathbb{N}_{0}\subseteq D\left(\mathcal{F}\right)$.
Then, we say a function $f\in C\left(\mathcal{F}\right)$ is \textbf{$\mathcal{F}$-rising-continuous}
/ \textbf{rising-continuous with respect to $\mathcal{F}$} whenever:
\begin{equation}
f\left(\mathfrak{z}\right)\overset{\mathcal{F}}{=}\lim_{n\rightarrow\infty}f\left(\left[\mathfrak{z}\right]_{p^{n}}\right)
\end{equation}
This generalizes the notion of rising-continuity as defined in my
dissertation, where, for a metrically complete valued field $\mathbb{F}$,
a function $f:\mathbb{Z}_{p}\rightarrow\mathbb{F}$ was said to be
rising-continuous whenever $\lim_{n\rightarrow\infty}f\left(\left[\mathfrak{z}\right]_{p^{n}}\right)$
converged to $f\left(\mathfrak{z}\right)$ in the topology of $\mathbb{F}$
for all $\mathfrak{z}\in\mathbb{Z}_{p}$.
\end{defn}
I originally introduced the notion of rising-continuity in my dissertation,
where I devoted a subsection (Section 3.2) to covering its basic properties
\cite{My Dissertation}, though traces of the idea seem to have been
known since at least Schikhof's time. A notion of rising-continuity
occurs, though not by that name, in Exercises 62A and 62B in Part
3 of Chapter 6 of Schikhof's \emph{Ultrametric Calculus}, as part
of his discussion of the van der Put basis \cite{Ultrametric Calculus}.

In between the writing of my dissertation and the writing of this
paper, I have succeeded in giving a categorical description of the
property of rising-continuity. The details are given below for the
interested reader, though they will not be needed for this paper,
except at the end of \textbf{Section \ref{subsec:Arithmetic-Dynamics,-Varieties,}},
where they become relevant to some of the speculation there.
\begin{defn}
A \textbf{profinite set}, also known as a \textbf{profinite space
}or \textbf{Stone space }is a topological space which is homeomorphic
to the projective/inverse limit in the category of topological spaces
of a projective/inverse system of finite topological spaces, each
of which is equipped with the discrete topology. The index set of
this projective system is required to be \textbf{directed from below},
meaning that for every $i,j\in I$, there exists $k\in I$ so that
$k\leq i$ and $k\leq j$.

We write \textbf{$\mathbf{ProS}$ }to denote the category of profinite
spaces. The objects of this category are profinite spaces; the morphisms
are continuous maps.

A \textbf{light profinite space/set }is a profinite space whose projective
system's index set is $\mathbb{N}_{0}^{\textrm{op}}$, meaning that
the projective system looks like:
\begin{equation}
A_{0}\overset{f_{0,1}}{\leftarrow}A_{1}\overset{f_{1,2}}{\leftarrow}A_{2}\overset{f_{2,3}}{\leftarrow}\cdots
\end{equation}
We write \textbf{$\mathbf{LProS}$} to denote the category of light
profinite spaces. The objects of this category are light profinite
spaces; the morphisms are continuous maps.
\end{defn}
Despite the simplicity of the rising continuity condition:
\begin{equation}
\lim_{n\rightarrow\infty}f\left(\left[\mathfrak{z}\right]_{p^{n}}\right)=f\left(\mathfrak{z}\right)
\end{equation}
it does not seem to have been given systematic treatment in the literature.
In hindsight, this is not entirely surprising. I have has shown that
the ``natural'' context for rising-continuity appears to be profinite
sets that possess an additional direct limit structure, one that is
compatible with their profinite structure.
\begin{defn}
Let $\mathcal{C}$ be a category. A \textbf{bidirectional system }in
$\mathcal{C}$ consists of:

\vphantom{}

I. A poset $I$ with a partial order, denoted $\leq$.

\vphantom{}

II. A collection of objects $\left\{ A_{i}\right\} _{i\in I}$ in
$\mathcal{C}$.

\vphantom{}

III. For every $i,j\in I$ with $i\leq j$, a morphism (called a\textbf{
projection})\textbf{ }$\pi_{i,j}:A_{j}\rightarrow A_{i}$ so that:

i. In the case where $i=j$, the morphism $\pi_{i,i}$ is $\textrm{Id}_{i}$,
the identity morphism on $A_{i}$ (i.e. $\textrm{Id}_{i}:A_{i}\rightarrow A_{i}$).

ii. For all $i,j,k\in I$ with $i\leq j\leq k$:
\begin{equation}
\pi_{i,k}=\pi_{i,j}\circ\pi_{j,k}
\end{equation}

That is, we have a commutative diagram:
\begin{equation}
\begin{array}{cccc}
 & \pi_{j,k}\\
A_{k} & \rightarrow & A_{j}\\
 & \searrow & \downarrow & \pi_{i,j}\\
\pi_{i,k} &  & A_{i}
\end{array}
\end{equation}

\vphantom{}

IV. For every $i,j\in I$ with $i\leq j$, a morphism (called a\textbf{
injection})\textbf{ }$\iota_{i,j}:A_{i}\rightarrow A_{j}$ so that:

i. In the case where $i=j$, $\iota_{i,i}=\textrm{Id}_{i}$

ii. For all $i,j,k\in I$ with $i\leq j\leq k$:
\begin{equation}
\iota_{i,k}=\iota_{j,k}\circ\iota_{i,j}
\end{equation}

That is, we have a commutative diagram:
\begin{equation}
\begin{array}{cccc}
 & \iota_{j,k}\\
A_{k} & \leftarrow & A_{j}\\
 & \nwarrow & \uparrow & \iota_{i,j}\\
\iota_{i,k} &  & A_{i}
\end{array}
\end{equation}

\vphantom{}V. The projections and injections satisfy the relations:
\begin{equation}
\pi_{i,j}\circ\iota_{i,j}=\textrm{Id}_{i}
\end{equation}
for all $i,j\in I$ with $i\leq j$. Moreover:
\begin{equation}
\pi_{i,\infty}\circ\iota_{i,\infty}=\textrm{Id}_{i}
\end{equation}
\[
\]
\end{defn}
I then define the following:
\begin{defn}
A \textbf{reffinite set / reffinite space }(or \textbf{ref-finite},
or \textbf{split profinite})\textbf{ }$S$ is a profinite space obtained
by taking the projective limit of a bidirectional system of finite
discrete topological spaces whose index set is directed from below.
we write $\mathbb{N}_{0}\left(S\right)$ to denote the injective limit
of the bidirectional system associated to $S$, and we write $S^{\prime}$
to denote the set of all elements of $S$ which are not in $\mathbb{N}_{0}\left(S\right)$.
We write $\mathbf{RefS}$ to denote the category of reffinite space;
this is the category whose objects are reffinite spaces and whose
morphisms are continuous maps that respect the definite structure,
in the sense that for any reffinite sets $R,S$ and any morphism $\psi:R\rightarrow S$
of reffinite sets:
\begin{align}
\psi\left(\mathbb{N}_{0}\left(R\right)\right) & \subseteq\mathbb{N}_{0}\left(S\right)\\
\psi\left(R^{\prime}\right) & \subseteq S^{\prime}
\end{align}
\end{defn}
\begin{rem}
The name ``reffinite'' is based on the formula \emph{reflexive}
+ \emph{finite}. Profinite objects are the projective limits of things,
while ind-finite objects are the injective limits of things, and these
two types of objects are categorically dual to one another. Since
objects which are self-dual are called reflexive, it seems reasonable
to unite profinite and ind-finite as ``reffinite''.
\end{rem}
\begin{defn}
Given a reffinite space $S$, the \textbf{rising topology }on $S$
is the topology generated by the following sets:

I. Sets which are open in the profinite topology on $S$. (We call
these \textbf{profinite-open sets}.)

\vphantom{}

II. Singletons $\left\{ s\right\} $, where $s$ is any element of
$\mathbb{N}_{0}\left(S\right)$. (We call these sets \textbf{point
sets}.)

\vphantom{}

III. Sets of the form $\left\{ \mathfrak{s}\right\} \cup\left\{ \pi_{j,\infty}\left(\mathfrak{s}\right):j\geq i\right\} $,
where $\mathfrak{s}$ is any element of $S^{\prime}$ and $i$ is
any element of $I$. (We call these sets \textbf{dust sets}.)

Furthermore, given a topological space $T$, a function $\chi:S\rightarrow T$
is said to be \textbf{rising-continuous }if it is continuous as a
map from $S$ (equipped with its rising topology) to $T$. That is,
$\chi^{-1}\left(U\right)$ is open in $S$'s rising topology for all
open sets in $T$. We write $C^{\textrm{rise}}\left(S,T\right)$ to
denote the space of rising-continuous functions from $S$ to $T$.
\end{defn}
\begin{rem}
Note that point sets and dust sets are both open and closed in the
rising topology.
\end{rem}
The following propositions are then easily proven:
\begin{prop}
Let $R$ be a reffinite set. Then $\mathbb{N}_{0}\left(R\right)$
is dense subset of $R$ in both the profinite and rising topology.
\end{prop}
\begin{prop}
Let $T$ be a topological space, let $S$ be a reffinite set, and
let $\chi:S\rightarrow T$ be rising-continuous. Then, for any reffinite
set $R$, and any continuous map $\psi:R\rightarrow S$ (i.e., continuous
w.r.t. the spaces' profinite topologies) satisfying:
\begin{align}
\psi\left(\mathbb{N}_{0}\left(R\right)\right) & =\mathbb{N}_{0}\left(S\right)\\
\psi\left(R^{\prime}\right) & \subseteq S^{\prime}
\end{align}
the composite $\chi\circ\psi:R\rightarrow T$ is rising-continuous.
\end{prop}
The standard example of such maps are as follows:
\begin{defn}
\label{def:space-change map}Let $p$ and $q$ be distinct primes.
Define the function $\psi_{q,p}:\mathbb{Z}_{p}\rightarrow\mathbb{Z}_{q}$
by:
\begin{equation}
\psi_{q,p}\left(\sum_{n=0}^{\infty}d_{n}p^{n}\right)\overset{\textrm{def}}{=}\sum_{n=0}^{\infty}d_{n}q^{n}\label{eq:definition of psi_p,q}
\end{equation}
I call $\psi_{q,p}$ the \textbf{$\left(p,q\right)$-adic space-change
map}.
\end{defn}
\begin{prop}
\label{prop:Characterization of psi_p,q}Let $p$ and $q$ be primes,
possibly non-distinct.

\vphantom{}

I. $\psi_{q,p}:\mathbb{Z}_{p}\rightarrow\mathbb{Z}_{q}$ is uniformly
continuous for all $p$ and $q$; in particular, it is Hölder continuous
with exponent $\log_{p}q=\frac{\ln q}{\ln p}$, with:
\begin{equation}
\left|\psi_{q,p}\left(\mathfrak{a}\right)-\psi_{q,p}\left(\mathfrak{b}\right)\right|_{q}\leq\left|\mathfrak{a}-\mathfrak{b}\right|_{p}^{\log_{p}q}\label{eq:Holder estimate for psi_p,q}
\end{equation}
This inequality holds with equality for all $\mathfrak{a},\mathfrak{b}$
if and only if $p\leq q$.

\vphantom{}

II. $\psi_{q,p}$ is injective if and only if $p\leq q$, and is surjective
if and only if $p\geq q$.

\vphantom{}

III. $\psi_{q,p}$ is the unique rising-continuous function satisfying
the functional equations:
\begin{equation}
\psi_{q,p}\left(p\mathfrak{z}+j\right)=j+q\psi_{q,p}\left(\mathfrak{z}\right),\textrm{ }\forall\mathfrak{z}\in\mathbb{Z}_{p},\textrm{ }\forall j\in\left\{ 0,\ldots,p-1\right\} \label{eq:Functoinal equations of psi_p,q}
\end{equation}
\end{prop}
\begin{rem}
Precomposing a $p$-adic F-series $X$ with $\psi_{p,q}:\mathbb{Z}_{q}\rightarrow\mathbb{Z}_{p}$
then gives a $q$-adic F-series $X\circ\psi_{p,q}$.
\end{rem}
Importantly, the existence of the Fourier transform of an F-series
remains unchanged when we precompose the F-series by a space-change
map.
\begin{example}
Given a $p$-adic F-series $X$ characterized by:
\begin{equation}
X\left(p\mathfrak{z}+k\right)=a_{k}X\left(\mathfrak{z}\right)+b_{k},\textrm{ }\forall k\in\left\{ 0,\ldots,p-1\right\} ,\textrm{ }\forall\mathfrak{z}\in\mathbb{Z}_{p}
\end{equation}
let $Y$ denote $X\circ\psi_{p,q}$. Then:
\begin{align*}
Y\left(q\mathfrak{z}+j\right) & =X\left(\psi_{p,q}\left(q\mathfrak{z}+j\right)\right)\\
 & =X\left(p\psi_{p,q}\left(\mathfrak{z}\right)+j\right)\\
 & =a_{j}X\left(\psi_{p,q}\left(\mathfrak{z}\right)\right)+b_{j}\\
 & =a_{j}Y\left(\mathfrak{z}\right)+b_{j}
\end{align*}
for all $j\in\left\{ 0,\ldots,p-1\right\} $. Thus, we can use the
formula from \textbf{Theorem \ref{thm:The-breakdown-variety} }to
compute Fourier transforms for $X$ and $Y$ over their respective
spaces.
\end{example}
This suggests one might be able to view F-series as functors out of,
say, the category of light profinite abelian groups, in the manner
of Clausen \& Scholze's condensed mathematics, due to the passing
resemblance that this set-up has to theirs \cite{condensed}. This
might present an avenue for realizing F-series as curves, in the manner
discussed in final subsections of \cite{first blog paper}.

Finally, for completeness' sake (pun intended), it is worth discussing
the adèlic aspects of frames and, more generally, completions of locally
convex topological rings. Indeed despite their apparent novelty, there
is a sense in which frames are a natural construction, especially
from a number theoretic perspective. To see this, however, we must
take up the language of \textbf{locally convex topological spaces}.
This is an extensive area of functional analysis, one that can be
applied to a host of different algebraic structures (algebras, rings,
vector spaces, modules, groups, etc.), united by the same underlying
idea, which we recall for the reader's edification in the case where
the underlying structure is a ring.
\begin{defn}
Let $\mathcal{X}$ be a topological ring (not necessarily unital or
commutative), let $I$ be an indexing set, and let $\mathscr{P}\overset{\textrm{def}}{=}\left\{ p_{i}\right\} _{i\in I}$
be a collection of seminorms on $\mathcal{X}$. For any $x\in\mathcal{X}$,
any real number $\epsilon>0$, and any $i\in I$, the open ball centered
at $x$ of radius $\epsilon$ with respect to $p_{i}$ is defined
by the set:
\begin{equation}
B_{i}\left(x;\epsilon\right)\overset{\textrm{def}}{=}\left\{ y\in\mathcal{X}:p_{i}\left(y-x\right)<\epsilon\right\} 
\end{equation}
The \textbf{locally convex topology }(LCT)\textbf{ on $\mathcal{X}$
induced/generated by $\mathscr{P}$} is the topology generated by
arbitrary unions and finite intersections of the $B_{i}\left(x;\epsilon\right)$.
A sequence $\left\{ x_{n}\right\} _{n\geq1}$ in $\mathcal{X}$ is
said to \textbf{converge }to $x\in\mathcal{X}$ with respect to $\mathscr{P}$
(i.e., in the topology induced by $\mathscr{P}$) precisely when:
\begin{equation}
\lim_{n\rightarrow\infty}p_{i}\left(x-x_{n}\right)\overset{\mathbb{R}}{=}0,\textrm{ }\forall i\in I
\end{equation}
A \textbf{locally convex topological ring }(LCTR) is a topological
ring equipped with an LCT.

Furthermore, a sequence $\left\{ x_{n}\right\} _{n\geq1}$ in $\mathcal{X}$
is said to be \textbf{Cauchy }with respect to $\mathscr{P}$ whenever,
for every $\epsilon>0$ and every $i\in I$, there is an integer $N_{i}\geq0$
so that:
\begin{equation}
p_{i}\left(x_{m}-x_{n}\right)<\epsilon,\textrm{ }\forall m,n\geq N_{i}
\end{equation}
A locally convex topology is said to be \textbf{complete }if every
Cauchy sequence converges. A LCTR is said to be \textbf{complete}
if its LCT is complete. Every locally convex topological ring has
a completion which is unique up to isomorphism; this completion is
defined in the usual way as the collection of equivalence classes
of Cauchy sequences.

An \textbf{isomorphism }of LCTRs is an isomorphism of the underlying
rings which is also continuous with respect to the LCT.
\end{defn}
For this we need to recall how to take the maximum of a collection
of finitely many seminorms.
\begin{defn}
Let $p_{1},\ldots,p_{N}\in\mathscr{S}\left(\mathcal{X}\right)$. Then,
I write $\max\left\{ p_{1},\ldots,p_{N}\right\} $ to denote the map
$\mathcal{X}\rightarrow\left[0,\infty\right)$ defined by:
\[
\max\left\{ p_{1},\ldots,p_{N}\right\} \left(x\right)\overset{\textrm{def}}{=}\max_{1\leq n\leq N}p_{n}\left(x\right),\textrm{ }\forall x\in\mathcal{X}
\]
\end{defn}
\begin{prop}
If $\rho_{1},\ldots,\rho_{N}\in\mathscr{S}\left(\mathcal{X}\right)$,
then $\max\left\{ \rho_{1},\ldots,\rho_{N}\right\} \in\mathscr{S}\left(\mathcal{X}\right)$.
Also, the maximum of finitely many ring norms is a ring norm.
\end{prop}
Proof: An easy exercise.
Q.E.D.

\vphantom{}

We also have the following fundamental observation:
\begin{prop}
Let $\mathcal{X}$ be a ring, and let $\mathscr{P}=\left\{ \nu_{1},\ldots,\nu_{d}\right\} $
be a set of finitely many norms on $\mathcal{X}$. Then, the following
topologies on $\mathcal{X}$ are equal:

I. The metric topology induced by the norm:
\begin{equation}
\left\Vert x\right\Vert \overset{\textrm{def}}{=}\max_{1\leq i\leq d}\nu_{i}\left(x\right),\textrm{ }\forall x\in\mathcal{X}
\end{equation}

II. The locally convex topology induced by $\mathscr{P}$.

\vphantom{}
\end{prop}
Proof: A sequence in $\mathcal{X}$ converges in (I) if and only if
it converges in (II). The finiteness of $\mathscr{P}$ is essential
to this argument. More specifically, identify $\mathcal{X}$ with
its image in the cartesian product $\mathcal{X}^{d}$ under the diagonal
embedding $x\mapsto\left(x,\ldots,x\right)$. Then, equip the $i$th
coordinate of this embedding with the norm $\nu_{i}$ and give the
entire embedding the product topology. Then, the product topology,
(I), and (II) are all equivalent.

Q.E.D.

\vphantom{}The number theoretic connection comes by way of the \textbf{Weak
Approximation Theorem }of algebraic number theory.
\begin{thm}[The Weak Approximation Theorem \cite{Neukirch}]
Let $K$ be any field, let $N\in\mathbb{N}_{1}$, and let $\left|\cdot\right|_{1},\ldots,\left|\cdot\right|_{N}$
be inequivalent, non-trivial absolute values on $K$. Choose any $N$
elements $a_{1},\ldots,a_{N}\in K$. Then, for any $\epsilon>0$,
there exists an $a\in K$ so that $\left|a-a_{n}\right|_{n}<\epsilon$
for all $n\in\left\{ 1,\ldots,N\right\} $.
\end{thm}
Before we use this, however, we need to establish the following:
\begin{cor}
Let $R$ be a ring, and let $\nu_{1},\ldots,\nu_{d}$ be finitely
many inequivalent, non-trivial norms on $R$, and let $\alpha$ be
their maximum:
\begin{equation}
\alpha\left(r\right)\overset{\textrm{def}}{=}\max_{1\leq n\leq d}\nu_{n}\left(r\right),\textrm{ }\forall r\in R
\end{equation}
Finally, let $\hat{R}$ denote the metric completion of $R$ with
respect to $\alpha$.

In this situation, $\hat{R}$ is isometric to $\prod_{n=1}^{N}R_{n}$,
the metric space given by taking the cartesian product of the $R_{n}$s,
where for each $n$, $R_{n}$ is the completion of $K$ with respect
to $\nu_{n}$. The norm on the product is:
\begin{equation}
\left\Vert \mathbf{z}\right\Vert \overset{\textrm{def}}{=}\max_{1\leq n\leq d}\nu_{n}\left(z_{i}\right)
\end{equation}
\end{cor}
Proof: An element of $\prod_{n=1}^{d}R_{n}$ is a tuple of the form
$\left(r_{1},\ldots,r_{d}\right)$ where, for each $n\in\left\{ 1,\ldots,d\right\} $,
$r_{n}\in R_{n}$. Now, since $R$ is dense in each $R_{n}$, for
each $n\in\left\{ 1,\ldots,d\right\} $, there exists a sequence $\left\{ a_{m,n}\right\} _{m\geq1}\in K$
so that $\nu_{n}\left(r_{n}-a_{m,n}\right)<1/m$ for all $m\geq1$
and all $n\in\left\{ 1,\ldots,d\right\} $. By the \textbf{Weak Approximation
Theorem}, for each $m$, there is a $b_{m}\in K$ so that $\nu_{n}\left(a_{m,n}-b_{m}\right)<1/m$
for all $n\in\left\{ 1,\ldots,d\right\} $. Hence, for all $m$ and
$n$, we have:
\[
\nu_{n}\left(r_{n}-b_{m}\right)\leq\nu_{n}\left(r_{n}-a_{m,n}\right)+\nu_{n}\left(a_{m,n}-b_{m}\right)<\frac{1}{2m}
\]
Thus, the sequence $\left\{ b_{m}\right\} _{m\geq1}$ is a sequence
in $R$ which, for each $n$, converges in $R_{n}$ to $r_{n}$ as
$m\rightarrow\infty$.

Now, consider $R$ as a metric space with $\alpha$ as its metric,
consider the diagonal embedding: $\phi:R\rightarrow R^{d}$ defined
by:
\begin{equation}
\phi\left(r\right)\overset{\textrm{def}}{=}\underbrace{\left(r,\ldots,r\right)}_{N\textrm{ times}}
\end{equation}
Note that we can make $R^{d}$ a metric space by defining $\delta:R^{d}\times R^{d}\rightarrow\mathbb{R}$
by:
\begin{equation}
\delta\left(\mathbf{r},\mathbf{s}\right)\overset{\textrm{def}}{=}\max_{1\leq n\leq d}\nu_{n}\left(r_{n}-s_{n}\right)
\end{equation}
where $\mathbf{r}=\left(r_{1},\ldots,r_{d}\right)$ and $\mathbf{s}=\left(s_{1},\ldots,s_{d}\right)$
are elements of $R^{d}$. Let $\Delta\left(R^{d}\right)$ denote the
image of $\phi$ in $R^{d}$. As constructed, $\phi$ is an isometry
from $\left(R,\alpha\right)$ onto $\left(\Delta\left(R^{d}\right),\delta\right)$.
Thus, the completion of $\left(R,\alpha\right)$ will be isometric
to the completion of $\left(\Delta\left(R^{d}\right),\delta\right)$.
Now, as we did above, pick:
\begin{equation}
\mathbf{r}\overset{\textrm{def}}{=}\left(r_{1},\ldots,r_{d}\right)\in\prod_{n=1}^{d}R_{n}
\end{equation}
Then, there exists a sequence $\left\{ b_{m}\right\} _{m\geq1}$ in
$R$ so that: 
\begin{equation}
\nu_{n}\left(r_{n}-b_{m}\right)<\frac{1}{2m},\textrm{ }\forall m\geq1,\forall n\in\left\{ 1,\ldots,d\right\} 
\end{equation}
Consequently, letting $\mathbf{b}_{m}$ be the element of $\Delta\left(R^{d}\right)$
whose every entry is $b_{m}$ (ex., $\mathbf{b}_{1}$ is a tuple of
$d$ consecutive $b_{1}$s), observe that for any $m_{1},m_{2}\geq1$:
\begin{align*}
\delta\left(\mathbf{b}_{m_{1}},\mathbf{b}_{m_{2}}\right) & =\max_{1\leq n\leq d}\nu_{n}\left(b_{m_{1}}-b_{m_{2}}\right)\\
 & \leq\max_{1\leq n\leq d}\nu_{n}\left(b_{m_{1}}-z_{n}\right)+\max_{1\leq n\leq d}\nu_{n}\left(z_{n}-b_{m_{2}}\right)\\
 & <\frac{1}{2m_{1}}+\frac{1}{2m_{2}}
\end{align*}
Thus, for any $\epsilon>0$, choose $\min\left\{ m_{1},m_{2}\right\} \geq\frac{1}{\epsilon}$.
Then:
\begin{equation}
\delta\left(\mathbf{b}_{m_{1}},\mathbf{b}_{m_{2}}\right)<\frac{\epsilon}{2}+\frac{\epsilon}{2}<\epsilon
\end{equation}
which shows that $\left\{ \mathbf{b}_{m}\right\} _{m\geq1}$ is Cauchy
in $\left(\Delta\left(R^{d}\right),\delta\right)$. Since $\delta$
extends to a metric on $\prod_{n=1}^{d}R_{n}$, any completion of
$\left(\Delta\left(R^{d}\right),\delta\right)$ is necessarily contained
in any completion of $\prod_{n=1}^{d}R_{n}$. In particular, since
$\mathbf{r}\in\prod_{n=1}^{d}R_{n}$, and:
\begin{equation}
\delta\left(\mathbf{r},\mathbf{b}_{m}\right)=\max_{1\leq n\leq d}\nu_{n}\left(z_{n}-b_{m}\right)<\frac{1}{2m}
\end{equation}
we conclude that $\mathbf{b}_{m}$ converges to $\mathbf{r}$ in $\prod_{n=1}^{d}R_{n}$.

Since $\mathbf{r}$ was arbitrary, we have that there is an isometric
copy of $\prod_{n=1}^{d}R_{n}$ contained within the completion of
$\left(\Delta\left(R^{d}\right),\delta\right)$, and hence, contained
within the completion of $\left(R,\alpha\right)$.

To finish, we just need to show that every element of $\left(\Delta\left(R^{d}\right),\delta\right)$'s
completion is contained in $\prod_{n=1}^{d}R_{n}$. For this, we first
note that since $\prod_{n=1}^{d}R_{n}$ is a cartesian product of
finitely many complete metric spaces (the $R_{n}$s), it itself is
a complete metric space under $\delta$. Since $\left(\Delta\left(R^{d}\right),\delta\right)$
embeds in $\prod_{n=1}^{d}R_{n}$, we then have that the completion
of $\left(\Delta\left(R^{d}\right),\delta\right)$ is contained in
$\prod_{n=1}^{d}R_{n}$, and we are done.

Q.E.D.

\vphantom{}A simple modification of this argument lets us work with
completions of the ring of functions $X\rightarrow K$, in a manner
that is clearly frame-theoretic.
\begin{cor}
Let $X$ be a set, let $K$ be a field, and let $R$ be the ring of
functions $X\rightarrow K$ under the standard pointwise operations.
Let $d\in\mathbb{N}_{1}$, let $\nu_{1},\ldots,\nu_{d}$ be inequivalent,
non-trivial norms values on $K$, and for each $n\in\left\{ 1,\ldots,d\right\} $,
let $\left\Vert \cdot\right\Vert _{n}$ be the supremum norm:
\begin{equation}
\left\Vert f\right\Vert _{n}\overset{\textrm{def}}{=}\sup_{x\in X}\nu_{n}\left(f\left(x\right)\right)
\end{equation}
Then:
\begin{equation}
\alpha\left(f\right)\overset{\textrm{def}}{=}\max_{1\leq n\leq d}\left\Vert f\right\Vert _{n}
\end{equation}
induces a norm on $R$, and the completion of $R$ with respect to
$\alpha$ is then a complete metric ring which is ring isometric to:
\begin{equation}
\left(\prod_{n=1}^{d}K_{n}\right)^{X}
\end{equation}
the ring of functions from $X$ to $\prod_{n=1}^{d}K_{n}$.
\end{cor}
Because of this, we see that frames and, more generally, locally convex
topologies, provide a natural method of keeping track of sequences
of objects (numbers, functions, etc.) whose convergence occurs in
potentially multiple different completions of a given ring. Although
I originally devised frames as a purely ad hoc device to simplify
discussion of convergence behavior of this type \cite{My Dissertation},
it is clear that we can turn this ad hoc aspect of the construction
on its head to capture arithmetic information. In particular, we can
consider projective systems of frames, where a frame $\mathcal{F}$
comes before another frame $\mathcal{G}$ if every $\mathcal{G}$-convergent
sequence is also $\mathcal{F}$-convergent. By considering the inverse
limit of this system, one ought to get something very adèle-adjacent,
if not fully adèlic, possibly with the requirement that we impose
a restriction condition in order to guarantee that the limit frame
generates adélic completions. In the spirit of \textbf{Example \ref{exa:Let-,-where}},
this suggests that perhaps the ``correct'' way to deal with pointwise
convergence of F-series whose parameters are all scalars is by considering
the image of the F-series under the diagonal embedding into an appropriately
constructed frame-theoretic version of the adèle ring of a global
field, though questions remain as to how to interpret an F-series
in this way when, at a given $\mathfrak{z}$, there might exist places
of the underlying global field for which $X\left(\mathfrak{z}\right)$
is either non-convergent, or divergent to $\infty$ in absolute value.
I discuss this a bit in \cite{my frames paper}, suggesting that projective
spaces might be used to deal with divergences, though I don't think
this technical detail will be resolved that easily.

On the other side of things\textemdash domain-wise\textemdash in viewing
F-series as functors out of profinite abelian groups, just as we can
realize $X$'s image as being adèle-valued by way of the diagonal
embedding, we might be able to realize F-series as functions out of
rings of finite adèles (direct products of all the possible completions
of $\mathcal{O}_{K}$ with respect to every non-trivial place of $K$)
by defining, the image of $X$ on $\mathcal{O}_{K_{\ell}}$, for every
place $\ell$ of $K$, by the composite $X\circ\psi_{p,\ell}$, where
$\psi_{p,\ell}:\mathcal{O}_{K_{\ell}}\rightarrow\mathbb{Z}_{p}$ is
the space-change map. In this way, we would have $X$ as being from
the adèles to the adèles, capturing all possible arithmetic information
one could get. Figuring out these details is just one of the many
avenues of future investigations my work opens onto.

\subsection{\label{subsec:Ascent-from-Fractional}Ascent from Fractional Ideals
of Polynomial Rings\protect 
}At the heart of our work is the idea of a $p$-adic F-series, an object
$X$ satisfying equations of the shape:
\begin{equation}
X\left(p\mathfrak{z}+k\right)=a_{k}X\left(\mathfrak{z}\right)+b_{k},\textrm{ }\forall k\in\left\{ 0,\ldots,p-1\right\} \label{eq:objects-1}
\end{equation}
for objects $a_{k}$ and $b_{k}$, and where $\mathfrak{z}$ is a
$p$-adic variable. Obviously, the nature and extent of $X$'s existence
depend what kinds of objects we are talking about.

In my dissertation and subsequent research, even in the case when
the $a_{k}$s and $b_{k}$s were elements of a field, I frequently
treated them as indeterminates, and found that all of the analysis
worked out without issue. I am pleased to report that this informal
procedure can be made rigorous.

To begin with, if we treat the $a_{k}$s and $b_{k}$s as formal indeterminates,
evaluating (\ref{eq:objects-1}) at $\mathfrak{z}=0$ and setting
$k=0$ gives us:
\begin{equation}
X\left(0\right)=a_{0}X\left(0\right)+b_{0}\label{eq:initial initial condition}
\end{equation}
which we can then solve to obtain:
\begin{equation}
X\left(0\right)=\frac{b_{0}}{1-a_{0}}\label{eq:initial condition-1}
\end{equation}
Setting $\mathfrak{z}=0$ in (\ref{eq:objects-1}) but letting $k$
vary, we would then get:
\begin{equation}
X\left(k\right)=a_{k}X\left(0\right)+b_{k}=\frac{a_{k}b_{0}}{1-a_{0}}+b_{k}
\end{equation}
Writing:
\begin{equation}
X\left(pk_{2}+k_{1}\right)=a_{k_{1}}X\left(k_{2}\right)+b_{k_{1}},\textrm{ }\forall k_{1},k_{2}\in\left\{ 0,\ldots,p-1\right\} 
\end{equation}
we then see that the value of $X$ on $\left\{ 0,\ldots,p-1\right\} $
determine the values of $X$ on $\left\{ p,\ldots,p^{2}-1\right\} $.
More generally, for any $n\geq1$, the values of $X$ on $\left\{ p^{n-1},\ldots,p^{n}-1\right\} $
then completely determine the values of $X$ on $\left\{ p^{n},\ldots,p^{n+1}-1\right\} $,
and thus, that the existence of the right-hand side of (\ref{eq:initial condition-1})
is sufficient to guarantee that $X\left(\mathfrak{z}\right)$ is uniquely
defined for all $\mathfrak{z}\in\mathbb{N}_{0}$.

Moreover, note that for any $\mathfrak{z}\in\mathbb{N}_{0}$, $X\left(\mathfrak{z}\right)$
will be a rational function in the $a_{k}$s and $b_{k}$s with coefficients
in $\mathbb{Z}$, and that, crucially, $\left(1-a_{0}\right)X\left(\mathfrak{z}\right)$
will always be an element of $\mathbb{Z}\left[a_{0},\ldots,a_{p-1},b_{0},\ldots,b_{p-1}\right]$.
This shows that $X\mid_{\mathbb{N}_{0}}$ takes values in a fractional
ideal of $\mathbb{Z}\left[a_{0},\ldots,a_{p-1},b_{0},\ldots,b_{p-1}\right]$.
For technical reasons, however, rather than use fractional ideals,
we will instead need to use localizations.

As is classical in commutative algebra and algebraic geometry, given
a ring of polynomials, we can impose relations on the polynomials'
variables by modding out by an ideal. Viewing $X$ as module-valued
over $\mathbb{N}_{0}$ allows us to impose nearly any relations or
evaluations on the $a_{k}$s and $b_{k}$s as we wish by simply taking
an appropriate quotient of $X$'s target space. The one caveat to
note is that, in order to get (\ref{eq:initial condition-1}) from
(\ref{eq:initial initial condition}) it must be that the relation
$a_{0}=1$ does not hold, else (\ref{eq:initial condition-1}) becomes
division by zero. That being said, note that if $a_{0}=1$, then (\ref{eq:initial initial condition})
forces the relation $b_{0}=0$ to hold. In this case, we are then
free to choose any value for $X\left(0\right)$, after which the recursive
structure of (\ref{eq:objects-1}) will guarantee that (\ref{eq:objects-1})
will have a unique solution on $\mathbb{N}_{0}$ for each chosen initial
condition.

In analogy with and generalization of the notion of an affine vector
space, we also have the notion of an affine subring.
\begin{defn}
Let $R$ be a commutative unital ring. . Then, we call a set $A\subseteq M$
an \textbf{affine subring} of $R$ whenever, for any $b\in A$, the
set $\left\{ a-b:a\in A\right\} $ is a subring of $M$.
\end{defn}
\begin{rem}
Just like with affine vector spaces, given $r\in R$ and a subring
$S$ of $R$, the coset $r+S$ is then an affine subring.
\end{rem}
Our set-up is then as follows:
\begin{prop}
\label{prop:module set up}Let $p$ be an integer $\geq2$, let $K$
be a global field (if $\textrm{char}K>0$, assume $\textrm{char}K$
is co-prime to $p$), let $a_{0},\ldots,a_{p-1}$ and $b_{0},\ldots,b_{p-1}$
be indeterminates, let $R$ denote the polynomial ring $\mathcal{O}_{K}\left[a_{0},\ldots,a_{p-1},b_{0},\ldots,b_{p-1}\right]$,
and let $\mathcal{R}$ denote the quotient:
\begin{equation}
\mathcal{R}\overset{\textrm{def}}{=}\frac{R\left[x\right]}{\left\langle \left(1-a_{0}\right)x-1\right\rangle }
\end{equation}
that is, the localization of $R$ away from the prime ideal $\left\langle 1-a_{0}\right\rangle $.
Then, the functional equation:
\begin{equation}
X\left(n\right)=a_{\left[n\right]_{p}}X\left(\theta_{p}\left(n\right)\right)+b_{\left[n\right]_{p}},\textrm{ }\forall n\in\mathbb{N}_{0}\label{eq:eq}
\end{equation}
has a unique solution $X:\mathbb{N}_{0}\rightarrow\mathcal{R}$.

Next, let $I$ be an ideal of $R$. Then:

I. The solution $X:\mathbb{N}_{0}\rightarrow\mathcal{R}$ of (\ref{eq:eq})
induces a unique solution $X:\mathbb{N}_{0}\rightarrow\mathcal{R}/I\mathcal{R}$
if and only if $\left\langle 1-a_{0}\right\rangle \nsubseteq I$;
here $X:\mathbb{N}_{0}\rightarrow\mathcal{R}/I\mathcal{R}$ is just
the image of $X:\mathbb{N}_{0}\rightarrow\mathcal{R}$ under the canonical
projection $\mathcal{R}\rightarrow\mathcal{R}/I\mathcal{R}$.

II. If $\left\langle 1-a_{0}\right\rangle \subseteq I$ and $\left\langle b_{0}\right\rangle \nsubseteq I$,
the set of all $X:\mathbb{N}_{0}\rightarrow\mathcal{R}/I\mathcal{R}$
satisfying (\ref{eq:eq}) is empty.

III. If $\left\langle 1-a_{0},b_{0}\right\rangle \subseteq I$, the
set of all $X:\mathbb{N}_{0}\rightarrow\mathcal{R}/I\mathcal{R}$
satisfying (\ref{eq:eq}) is an affine submodule of $\mathcal{R}/I\mathcal{R}$
of dimension $1$.
\end{prop}
Proof: The only part needing explication is (III). So, letting $I$
be as given, let $X,Y:\mathbb{N}_{0}\rightarrow\mathcal{R}/I\mathcal{R}$
be solutions of (\ref{eq:eq}). Then, for $Z=X-Y$:
\begin{align*}
Z\left(n\right) & =X\left(n\right)-Y\left(n\right)\\
 & =a_{\left[n\right]_{p}}X\left(\theta_{p}\left(n\right)\right)+b_{\left[n\right]_{p}}-\left(a_{\left[n\right]_{p}}Y\left(\theta_{p}\left(n\right)\right)+b_{\left[n\right]_{p}}\right)\\
 & =a_{\left[n\right]_{p}}\left(X\left(\theta_{p}\left(n\right)\right)-Y\left(\theta_{p}\left(n\right)\right)\right)\\
 & =a_{\left[n\right]_{p}}Z\left(\theta_{p}\left(n\right)\right)
\end{align*}
where all equalities are in $\mathcal{R}/I\mathcal{R}$. For $n\in\left\{ 0,\ldots,p-1\right\} $,
we obtain:
\begin{equation}
Z\left(n\right)=a_{n}Z\left(0\right)
\end{equation}
Moreover, since $a_{0}=1$ in $\mathcal{R}/I\mathcal{R}$, we see
that $Z\left(0\right)=a_{0}Z\left(0\right)=Z\left(0\right)$, which
leaves $Z\left(0\right)$ undetermined. However, once $Z\left(0\right)$
is chosen, that choice completely determines the value of $Z\left(n\right)$
for all $n\geq0$. So, given any two functions $Z_{1},Z_{2}:\mathbb{N}_{0}\rightarrow\mathcal{R}/I\mathcal{R}$
so that:
\begin{equation}
Z_{j}\left(n\right)=a_{\left[n\right]_{p}}Z_{j}\left(n\right),\textrm{ }\forall n\geq0
\end{equation}
we see that, for any $r_{1},r_{2}\in R$:
\begin{equation}
\left(r_{1}Z_{1}+r_{2}Z_{2}\right)\left(n\right)=a_{\left[n\right]_{p}}\left(r_{1}Z_{1}+r_{2}Z_{2}\right)\left(n\right),\textrm{ }\forall n\geq0
\end{equation}
Hence, the set of $f:\mathbb{N}_{0}\rightarrow\mathcal{R}/I\mathcal{R}$
satisfying the equation: 
\begin{equation}
f\left(n\right)=a_{\left[n\right]_{p}}f\left(n\right),\textrm{ }\forall n\geq0
\end{equation}
is an $R$-subring of $\mathcal{R}/I\mathcal{R}$ of dimension $1$.
Since $Z$ was the difference of two $\mathcal{R}/I\mathcal{R}$-valued
solutions of (\ref{eq:eq}), we conclude that the set of $\mathcal{R}/I\mathcal{R}$-valued
solutions of (\ref{eq:eq}) is a one-dimensional affine subring of
$\mathcal{R}/I\mathcal{R}$.

Q.E.D.

\vphantom{}In light of this, we will use the following terminology
to describe our ideals.
\begin{defn}
Given $R$ as in \textbf{Proposition \ref{prop:module set up}}, we
say an ideal $I\subseteq R$ is \textbf{admissible }if either:

I. $\left\langle 1-a_{0}\right\rangle \nsubseteq I$.

II. $\left\langle 1-a_{0},b_{0}\right\rangle \subseteq I$

We say $I$ is a \textbf{unique solution ideal }when (I) holds, and
say $I$ is a \textbf{non-unique solution ideal }when (II) holds.
\end{defn}
With \textbf{Proposition \ref{prop:module set up}} now established,
we can freely consider the quotient of $\mathcal{R}$ by any admissible
ideal. This formalizes the idea of treating the $a_{k}$s and $b_{k}$s
as indeterminates. With that being done, the next order of business
is to consider how to extend $X$ from a function on $\mathbb{N}_{0}$
to one on $\mathbb{Z}_{p}$ or a subset thereof. For this, we nest
equation (\ref{eq:eq}) $N$ times to obtain:
\begin{equation}
X\left(m\right)\overset{\mathcal{R}}{=}X\left(\theta_{p}^{\circ N}\left(m\right)\right)\prod_{n=0}^{N-1}a_{\left[\theta_{p}^{\circ n}\left(m\right)\right]_{p}}+\sum_{n=0}^{N-1}b_{\left[\theta_{p}^{\circ n}\left(m\right)\right]_{p}}\prod_{k=0}^{n-1}a_{\left[\theta_{p}^{\circ k}\left(m\right)\right]_{p}}
\end{equation}
for all integers $m\geq0$ and where, as indicated, the equality occurs
in $\mathcal{R}$. Letting $\mathfrak{z}\in\mathbb{Z}_{p}$ be arbitrary,
set $m=\left[\mathfrak{z}\right]_{p^{N}}$. Since $\theta_{p}^{\circ N}\left(m\right)=0$
for all $m\in\left\{ 0,\ldots,p^{N}-1\right\} $, with the help of
\textbf{Proposition \ref{prop:Product identity for kappa_X}},\textbf{
}we get:
\begin{equation}
X\left(\left[\mathfrak{z}\right]_{p^{N}}\right)\overset{\mathcal{R}}{=}X\left(0\right)a_{0}^{N}\kappa\left(\left[\mathfrak{z}\right]_{p^{N}}\right)+\sum_{n=0}^{N-1}b_{\left[\theta_{p}^{\circ n}\left(\mathfrak{z}\right)\right]_{p}}a_{0}^{n}\kappa\left(\left[\mathfrak{z}\right]_{p^{n}}\right),\textrm{ }\forall N\geq0,\textrm{ }\forall\mathfrak{z}\in\mathbb{Z}_{p}\label{eq:X of z mod p to the N-1}
\end{equation}
where $\kappa$ denotes:
\begin{equation}
\kappa\left(m\right)\overset{\textrm{def}}{=}\prod_{k=1}^{p-1}\left(\frac{a_{k}}{a_{0}}\right)^{\#_{p:k}\left(m\right)},\textrm{ }\forall m\in\mathbb{N}_{0}
\end{equation}

\begin{rem}
Technically, it is an abuse of notation to write $\kappa$ separately,
as $\kappa\left(m\right)\notin\mathcal{R}$. Rather, the correct function
to define is the M-function $\left\{ a_{0}^{n}\kappa\left(\left[\mathfrak{z}\right]_{p^{n}}\right)\right\} _{n\geq0}$.
The presence of the factor $a_{0}^{n}$ guarantees that $a_{0}^{n}\kappa\left(\left[\mathfrak{z}\right]_{p^{n}}\right)\in\mathcal{R}$
for all $n$ and $\mathfrak{z}$. Indeed, we have: 
\begin{equation}
a_{0}^{n}\kappa\left(\left[\mathfrak{z}\right]_{p^{n}}\right)\in R=\mathcal{O}_{K}\left[a_{0},\ldots,a_{p-1},b_{0},\ldots,b_{p-1}\right]
\end{equation}
for all $n$ and $\mathfrak{z}$.
\end{rem}
The obvious thing to do is to take the limit of (\ref{eq:X of z mod p to the N-1})
as $N\rightarrow\infty$. Not only would this extend/interpolate $X$
to $\mathfrak{z}\in\mathbb{Z}_{p}$ by defining $X\left(\mathfrak{z}\right)$
to be $\lim_{N\rightarrow\infty}X\left(\left[\mathfrak{z}\right]_{p^{N}}\right)$
for all $\mathfrak{z}$ for which the limit exists, it would also
give us the F-series representation of $X$ for free. In order to
do this, however, we need to appeal to some notion of convergence.
We will do this by invoking a frame and taking a completion. However,
as this involves working with functions valued in a fractional ideal,
rather than a ring, we need to modify our notion of frames slightly.
Thankfully, this is easily done using quotient frames to pass from
a frame on $R$ to a frame on $\mathcal{R}$, and we will maintain
the compatibility with quotients provided we restrict ourselves to
unique solution ideals. The following example illustrates why:
\begin{example}
Since:
\begin{equation}
\mathcal{R}=\frac{R\left[x\right]}{\left\langle \left(1-a_{0}\right)x-1\right\rangle }
\end{equation}
the $\mathcal{F}_{\mathfrak{z}}$s of an $\mathcal{R}$-frame $\mathcal{F}$
will be a composite of the evaluate-at-$\mathfrak{z}$-map $\mathcal{R}^{D\left(\mathcal{F}\right)}\rightarrow\mathcal{R}$
and an absolute value $\left|\cdot\right|_{\mathcal{F}\left(\mathfrak{z}\right)}\in\mathscr{V}\left(\mathcal{R}\right)$.
Here, $\left|\cdot\right|_{\mathcal{F}\left(\mathfrak{z}\right)}$
is going to be the quotient absolute value generated by some absolute
value $v\in\mathscr{V}\left(R\left[x\right]\right)$:
\begin{equation}
\left|r\right|_{\mathcal{F}\left(\mathfrak{z}\right)}=\inf_{s\in\left\langle \left(1-a_{0}\right)x-1\right\rangle }v\left(r+s\right)
\end{equation}
for any representative $r$ of an equivalence class in $\mathcal{R}$.
Thus, given any ideal $I\subseteq R$, we have that: 
\begin{equation}
\mathcal{R}/I\mathcal{R}=\frac{R\left[x\right]}{\left\langle \left(1-a_{0}\right)x-1,I\right\rangle _{R\left[x\right]}}
\end{equation}
where $\left\langle \left(1-a_{0}\right)x-1,I\right\rangle _{R\left[x\right]}$
is the ideal generated over $R\left[x\right]$ by $\left(1-a_{0}\right)x-1$
and all of the elements of $I$. As such:
\begin{equation}
\left|r\right|_{\left(\mathcal{F}/I\right)\left(\mathfrak{z}\right)}=\inf_{s\in\left\langle \left(1-a_{0}\right)x-1,I\right\rangle _{R\left[x\right]}}v\left(r+s\right)\label{eq:F mod I for R script}
\end{equation}
for any representative $r$ of an equivalence class in $\mathcal{R}/I\mathcal{R}$.
If $I$ is not a unique solution ideal, then $\left\langle 1-a_{0}\right\rangle \subseteq I$,
in which case:
\begin{equation}
-1=\left(1-a_{0}\right)x-1+\left(1-a_{0}\right)x\in\left\langle \left(1-a_{0}\right)x-1,I\right\rangle _{R\left[x\right]}
\end{equation}
and hence, $\left\langle \left(1-a_{0}\right)x-1,I\right\rangle _{R\left[x\right]}=R\left[x\right]$,
which forces (\ref{eq:F mod I for R script}) to vanish for all $r$
(simply choose $s=-r$).
\end{example}
Thus, we can justify taking the limit as $N\rightarrow\infty$ in
(\ref{eq:X of z mod p to the N-1}) by picking an $\mathcal{R}$-frame
$\mathcal{F}$ on $\mathbb{Z}_{p}$ so that everything converges.
The theorem below summarizes both what happens and what we will need.
\begin{thm}
\label{thm:variety frames}Let everything be as given in \textbf{Proposition
\ref{prop:module set up}}. Let $\mathcal{F}$ be a reffinite, evaluative,
multiplicative $\mathcal{R}$-frame on $\mathbb{Z}_{p}$ so that the
limit:
\begin{equation}
X\left(\mathfrak{z}\right)\overset{\textrm{def}}{=}\lim_{N\rightarrow\infty}\sum_{n=0}^{N-1}b_{\left[\theta_{p}^{\circ n}\left(\mathfrak{z}\right)\right]_{p}}a_{0}^{n}\kappa\left(\left[\mathfrak{z}\right]_{p^{n}}\right)\label{eq:X variety frame}
\end{equation}
is $\mathcal{F}$-convergent. $X$ is then the unique $\mathcal{F}$-rising-continuous
solution $X\in C\left(\mathcal{F}\right)$ to the functional equation:
\begin{equation}
X\left(\mathfrak{z}\right)\overset{\mathcal{F}\left(\mathfrak{z}\right)}{=}a_{\left[\mathfrak{z}\right]_{p}}X\left(\theta_{p}\left(\mathfrak{z}\right)\right)+b_{\left[\mathfrak{z}\right]_{p}},\textrm{ }\forall\mathfrak{z}\in D\left(\mathcal{F}\right)\label{eq:rising continuation}
\end{equation}
where $\mathcal{F}\left(\mathfrak{z}\right)$ denotes the completion
of $\mathcal{R}$ with respect to the absolute value $\left|\cdot\right|_{\mathcal{F}\left(\mathfrak{z}\right)}$
invoked by $\mathcal{F}$ at $\mathfrak{z}$. Moreover, for any unique
solution ideal $I\subseteq R$, $X$ descends to an $\left(\mathcal{F}/I\mathcal{R}\right)$-rising-continuous
function $X\in C\left(\mathcal{F}/I\mathcal{R}\right)$ (given by
the limit of (\ref{eq:X variety frame}) with respect to $\mathcal{F}/I\mathcal{R}$)
so that (\ref{eq:rising continuation}) holds with equality in $\left(\mathcal{F}/I\mathcal{R}\right)\left(\mathfrak{z}\right)$,
the completion of $\mathcal{R}/I\mathcal{R}$ with respect to the
absolute value $\left|\cdot\right|_{\left(\mathcal{F}/I\mathcal{R}\right)\left(\mathfrak{z}\right)}$
invoked by $\mathcal{F}/I\mathcal{R}$ at $\mathfrak{z}$.
\end{thm}
Proof: Take limits of (\ref{eq:X of z mod p to the N-1}).

Q.E.D.
\begin{rem}
The reason why $R$ is the ring of polynomials over $\mathcal{O}_{K}$
rather than over $\mathbb{Z}$ is because choosing $\mathcal{O}_{K}$
allows us to use absolute values that come from $\mathcal{O}_{K}$.
Any absolute value on $R$ is necessarily an extension of an absolute
value on $\mathcal{O}_{K}$, and since $\mathcal{R}$ is a quotient
of $R\left[x\right]$, any absolute value on $\mathcal{R}$ is going
to be induced by an absolute value on $R\left[x\right]$. Thus, $\mathcal{R}$
admits absolute values derived from any of the absolute values of
$\mathcal{O}_{K}$. If we restricted our attention solely to $\mathbb{Z}$,
we would only be able to use absolute values coming from $\mathbb{Z}$.
\end{rem}
To use \textbf{Theorem \ref{thm:variety frames}}, we simply need
to construct a frame satisfying its hypotheses. This is done easily
enough. For the sake of convenience, we give our construction a name.
\begin{defn}[The Digit Construction]
\label{def:digit construction-1}Let $D_{0}$ denote $\mathbb{N}_{0}$,
and let $D_{1}$ denote the set of all $p$-adic integers containing
infinitely many $1$s digits. Then, by induction, for each $k\in\left\{ 2,\ldots,p-1\right\} $,
define $D_{k}$ to be the set of all $p$-adic integers $\mathfrak{z}$
so that:

\textbullet{} $\mathfrak{z}$ has infinitely many $k$s digits.

\textbullet{} $\mathfrak{z}$ is \emph{not} contained in $D_{1},\ldots,D_{k-1}$.

Note that the $D_{k}$s then form a partition of $\mathbb{Z}_{p}$
into pair-wise disjoint sets. We call this partition the \textbf{digit
domain}.

Let $R$ be an arbitrary commutative ring. A \textbf{digital $R$-frame
}is an $R$-frame $\mathcal{F}$ so that there are absolute values
$v_{0},\ldots,v_{p-1}\in\mathscr{V}\left(R\right)$ for which:
\begin{equation}
\mathcal{F}\left(\mathfrak{z}\right)=\begin{cases}
v_{0} & \textrm{if }\mathfrak{z}\in D_{0}\cap D\left(\mathcal{F}\right)\\
v_{1} & \textrm{if }\mathfrak{z}\in D_{1}\cap D\left(\mathcal{F}\right)\\
\vdots & \vdots\\
v_{p-1} & \textrm{if }\mathfrak{z}\in D_{p-1}\cap D\left(\mathcal{F}\right)
\end{cases},\textrm{ }\forall\mathfrak{z}\in D\left(\mathcal{F}\right)
\end{equation}
That is: 
\begin{equation}
\mathcal{F}_{\mathfrak{z}}\left(f\right)=\sum_{k=0}^{p-1}\left[\mathfrak{z}\in D_{k}\cap D\left(\mathcal{F}\right)\right]\left|f\left(\mathfrak{z}\right)\right|_{v_{k}},\textrm{ }\forall f\in C\left(\mathcal{F}\right)
\end{equation}
\end{defn}
We then have the following:
\begin{prop}
Let $\mathcal{F}$ be the digital $\mathcal{R}$-frame so that $\mathbb{N}_{0}\subseteq D\left(\mathcal{F}\right)$
and:
\begin{equation}
\left|a_{0}\right|_{\mathcal{F}\left(\mathfrak{z}\right)}\limsup_{n\rightarrow\infty}\left|\kappa\left(\left[\mathfrak{z}\right]_{p^{n}}\right)\right|_{\mathcal{F}\left(\mathfrak{z}\right)}^{1/n}<1,\textrm{ }\forall\mathfrak{z}\in D\left(\mathcal{F}\right)
\end{equation}
Then, $\mathcal{F}$ satisfies the hypothesis of \textbf{Theorem \ref{thm:variety frames}}.
\end{prop}
Proof: Use the root test for series convergence. Moreover, note that
$\mathcal{F}$ will then be shift-invariant thanks to $\kappa$ (use
\textbf{Proposition \ref{prop:Kappa shift equation}}).
Q.E.D.

\vphantom{}

Applying this to the case of a family of finitely many F-series $X_{1},\ldots,X_{d}$,
we get the following:
\begin{thm}
\label{thm:variety frame general}Let $d,p$ be positive integers,
with $d\geq1$ and $p\geq2$, let $K$ be a global field (if $\textrm{char}K>0$,
then assume $\textrm{char}K$ is co-prime to $p$). For each $j\in\left\{ 1,\ldots,d\right\} $
and $k\in\left\{ 0,\ldots,p-1\right\} $ let $a_{j,k}$ and $b_{j,k}$
be indeterminates, and let $R_{d}\left(K\right)$ (a.k.a., $R_{d}$,
$R$) be the ring of polynomials in the $a_{j,k}$s and $b_{j,k}$s
with coefficients in $\mathcal{O}_{K}$, so that $R$ is then a free
$\mathcal{O}_{K}$-algebra in $2dp$ indeterminates. Then, let:
\begin{equation}
\mathcal{R}_{d}\left(K\right)\overset{\textrm{def}}{=}\frac{R\left[x_{1},\ldots,x_{d}\right]}{\left\langle \left(1-a_{1,0}\right)x_{1}-1,\ldots,\left(1-a_{d,0}\right)x_{d}-1\right\rangle }
\end{equation}
(a.k.a., $\mathcal{R}_{d}$, $\mathcal{R}$) be the localization of
$R$ away from the prime ideal $\left\langle 1-a_{1,0},\ldots,1-a_{d,0}\right\rangle $.

Writing:
\begin{equation}
\kappa_{j}\left(n\right)\overset{\textrm{def}}{=}\prod_{k=1}^{p-1}\left(\frac{a_{j,k}}{a_{j,0}}\right)^{\#_{p:k}\left(n\right)}
\end{equation}
let $\mathcal{F}$ be a digital $\mathcal{R}$-frame on $\mathbb{Z}_{p}$
so that:
\begin{equation}
\left|a_{j,0}\right|_{\mathcal{F}\left(\mathfrak{z}\right)}\limsup_{n\rightarrow\infty}\prod_{k=1}^{p-1}\left|\frac{a_{j,k}}{a_{j,0}}\right|_{\mathcal{F}\left(\mathfrak{z}\right)}^{\#_{p:k}\left(\left[\mathfrak{z}\right]_{p^{n}}\right)}<1,\textrm{ }\forall\mathfrak{z}\in D\left(\mathcal{F}\right),\textrm{ }\forall j\in\left\{ 1,\ldots,d\right\} \label{eq:root condition for individual j}
\end{equation}
(For example, given a choice of the absolute values for $\mathcal{F}\left(\mathfrak{z}\right)$,
define $\mathcal{F}$'s domain to be the set of all $\mathfrak{z}\in\mathbb{Z}_{p}$
so that (\ref{eq:root condition for individual j}) holds true.) Then:
\begin{equation}
X_{j}\left(\mathfrak{z}\right)\overset{\textrm{def}}{=}\lim_{N\rightarrow\infty}\sum_{n=0}^{N-1}b_{j,\left[\theta_{p}^{\circ n}\left(\mathfrak{z}\right)\right]_{p}}a_{j,0}^{n}\kappa_{j}\left(\left[\mathfrak{z}\right]_{p^{n}}\right)\label{eq:X variety frame-1}
\end{equation}
is $\mathcal{F}$-convergent for all $j\in\left\{ 1,\ldots,d\right\} $.
Each $X_{j}$ is then the unique $\mathcal{F}$-rising-continuous
solution $X_{j}\in C\left(\mathcal{F}\right)$ to the functional equation:
\begin{equation}
X_{j}\left(\mathfrak{z}\right)\overset{\mathcal{F}\left(\mathfrak{z}\right)}{=}a_{j,\left[\mathfrak{z}\right]_{p}}X_{j}\left(\theta_{p}\left(\mathfrak{z}\right)\right)+b_{j,\left[\mathfrak{z}\right]_{p}},\textrm{ }\forall\mathfrak{z}\in D\left(\mathcal{F}\right)\label{eq:rising continuation-1}
\end{equation}
Moreover, for any ideal $I\subseteq R$ so that $\left\langle 1-a_{j,0}\right\rangle \nsubseteq I$
for all $j\in\left\{ 1,\ldots,d\right\} $ (we call such an ideal
a \textbf{unique solution ideal for the $X_{j}$s}), the $X_{j}$s
descend to $\mathcal{F}/I\mathcal{R}$-rising-continuous functions
$X_{j}\in C\left(\mathcal{F}/I\mathcal{R}\right)$ (given by the limit
of (\ref{eq:X variety frame-1}) with respect to $\mathcal{F}/I\mathcal{R}$)
so that (\ref{eq:rising continuation-1}) holds with equality in $\left(\mathcal{F}/I\mathcal{R}\right)\left(\mathfrak{z}\right)$.
\end{thm}
Proof: Same as above.

Q.E.D.
\begin{rem}
From here on out, in order to save space, I will commit a small abuse
of notation and write $\mathcal{F}/I$ instead of the more appropriate
notation $\mathcal{F}/I\mathcal{R}$.
\end{rem}
With the above, we can finally rigorously identify the ambient space
in which Section \ref{subsec:The-Central-Computation}'s computations
occurred.

In the paper's introduction, we initially defined an F-series as a
function from $\mathbb{Z}_{p}$ to a ring of formal power series generated
by the F-series' parameters. Thus, if $X$ was characterized by the
equations:
\begin{equation}
X\left(p\mathfrak{z}+k\right)=a_{k}X\left(\mathfrak{z}\right)+b_{k}
\end{equation}
for all $k\in\left\{ 0,\ldots,p-1\right\} $, we could view $X$ as
a map $\mathbb{Z}_{p}\rightarrow\mathcal{O}_{K}\left[\left[a_{0},\ldots,a_{p-1},b_{0},\ldots b_{p-1}\right]\right]$,
or $\mathbb{Z}_{p}\rightarrow\mathcal{O}_{K}\left[\left[a_{0},\ldots,a_{p-1}\right]\right]\left[b_{0},\ldots b_{p-1}\right]$,
seeing as any $b_{k}$ never occurs in $X$ with an exponent $\geq2$.
While this set-up works fine as a provisional definition, it has poor
compatibility with quotients by ideals. This is due to the following
fact about quotients of power series rings over a Dedekind domain:
\begin{fact}
\label{fact:Let--be}Let $R$ be a Dedekind domain which is not a
field, and let $a,b\in R$ with $b\neq0$, so that $\left\langle a\right\rangle +\left\langle b\right\rangle =R$
(the ideals are co-prime). Then, the ring $R\left[\left[a/b\right]\right]$
is defined by the quotient:
\begin{equation}
R\left[\left[\frac{a}{b}\right]\right]\overset{\textrm{def}}{=}\frac{R\left[\left[x\right]\right]}{\left\langle bx-a\right\rangle }
\end{equation}
where $x$ is an indeterminate. Then:
\begin{equation}
R\left[\left[\frac{a}{b}\right]\right]\cong\begin{cases}
\left\{ 0\right\}  & \textrm{if }a\in R^{\times}\\
R & \textrm{if }a=0\\
\prod_{\mathfrak{p}\mid\left\langle a\right\rangle }R_{\mathfrak{p}} & \textrm{if }a\in R\backslash\left(\left\{ 0\right\} \cup R^{\times}\right)
\end{cases}
\end{equation}
where the product is taken over all prime ideals $\mathfrak{p}$ of
$R$ which divide $\left\langle a\right\rangle $, and where $R_{\mathfrak{p}}$
is the $\mathfrak{p}$-adic completion of $R$.
\end{fact}
Because of this, we see that evaluating the indeterminate of a formal
power series ring over a Dedekind domain $R$ at an element $a/b$
of the field of fractions of $R$ ignores denominators. Moreover,
whenever the evaluation homomorphism on such a ring spits out something
other than $R$ or $0$, the ring we get is a product of completions
of $R$ with respect to the primes that divide $a$. As such, if the
image of $R\left[\left[x\right]\right]$ under an evaluation homomorphism
is non-trivial (i.e., neither $0$ nor $R$), then the image is going
to naturally be a non-archimedean space. Because of this, if we had
an $R\left[\left[x\right]\right]$-frame $\mathcal{F}$ and wanted
to consider the quotient frame induced by an ideal $I$ of $R\left[\left[x\right]\right]$,
unless $R\left[\left[x\right]\right]/I$ was either $0$ or $R$,
the seminorms produced by $\mathcal{F}/I$ would have to be non-archimedean,
because $R\left[\left[x\right]\right]/I$ is non-archimedean. This
would then restrict us to having to use non-archimedean seminorms
for $\mathcal{F}$ itself, as the seminorms induced by a quotient
have the same quality (archimedean or non-archimedean) as the seminorms
that induced them.

Thankfully, frames allow us to circumvent this technical dilemma.
Instead of treating our $X$s as valued in a ring of formal power
series, we will treat our $X$s as being element of $C\left(\mathcal{F}\right)$,
where $\mathcal{F}$ is an $\mathcal{R}$-frame as described in \textbf{Theorem
\ref{thm:variety frame general}}. In this case, to use that theorem's
terminology, the $X$s will take values in the locally convex topological
$R$-algebra obtained by completing $\mathcal{R}^{D\left(\mathcal{F}\right)}$
with respect to the locally convex topology induced by the seminorms
$\left\{ \mathcal{F}_{\mathfrak{z}}:\mathfrak{z}\in D\left(\mathcal{F}\right)\right\} $.
As such:
\begin{assumption}
Henceforth, we will be working using the set-up from \textbf{Theorem
\ref{thm:variety frame general}}. I will use $R_{d}$ and $\mathcal{R}_{d}$
when I want to indicate the number of distinct F-series we are considering
at a given moment. Thus, if I say, ``consider the case of $\mathcal{R}_{1}$'',
it means we are working with a single F-series and our indeterminates
will be $a_{0},\ldots,a_{p-1}$ and $b_{0},\ldots b_{p-1}$.

Thus, in general, we will be working with $\mathcal{R}$ and a unique
solution ideal $I$. Our F-series will then take values in $\mathcal{R}/I\mathcal{R}$.
Their Fourier transforms will be elements of $\textrm{Frac}\left(R/I\right)\left(\zeta_{p^{\infty}}\right)$,
the maximal $p$-power cyclotomic extension of the field of fractions
of $R/I$.

Because of the compatibility of our arguments with quotients by unique
solution ideals $I$ of $R$, we will do all of our work ignoring
$I$, and having our F-series take values in $\mathcal{R}$ and their
Fourier transforms take values in $\textrm{Frac}\left(R\right)\left(\zeta_{p^{\infty}}\right)$.
\end{assumption}

\subsection{\label{subsec:Quasi-Integrability-=000026-Degenerate}Quasi-Integrability
\& Degenerate Measures}

FOR THIS SECTION, WE FIX AN INTEGER $p\geq2$ AND A COMMUTATIVE, UNITAL,
INTEGRAL DOMAIN WHOSE CHARACTERISTIC, IF POSITIVE, IS CO-PRIME TO
$p$. WE ALSO FIX AN EMBEDDING OF THE MULTIPLICATIVE GROUP OF $p$-POWER
ROOTS OF UNITY IN $R\left(\zeta_{p^{\infty}}\right)$

Though the interested reader can refer to \cite{2nd blog paper} for
more details, the essential idea behind the Fourier analysis we've
been using is that for profinite abelian groups like $\mathbb{Z}_{p}$,
one can formulate Fourier theory in a very general way using the Schwartz-Bruhat
(SB) function. For Pontryagin duality, all unitary characters $\mathbb{Z}_{p}\rightarrow R\left(\zeta_{p^{\infty}}\right)$
are maps of the form:
\begin{equation}
\mathfrak{z}\in\mathbb{Z}_{p}\mapsto e^{2\pi i\left\{ t\mathfrak{z}\right\} _{p}}\in R\left(\zeta_{p^{\infty}}\right)
\end{equation}
for some $t\in\hat{\mathbb{Z}}_{p}=\mathbb{Z}\left[1/p\right]/\mathbb{Z}$,
where, here, recall that $e^{2\pi i\left\{ t\mathfrak{z}\right\} _{p}}$
denotes the image of the associated $p$-power root of unity under
our chosen embedding into $R\left(\zeta_{p^{\infty}}\right)$. The
reason everything works is due to the fundamental identity:
\begin{equation}
\left[\mathfrak{z}\overset{p^{n}}{\equiv}k\right]=\frac{1}{p^{n}}\sum_{\left|t\right|_{p}\leq p^{n}}e^{2\pi i\left\{ t\left(\mathfrak{z}-k\right)\right\} _{p}},\textrm{ }\forall n\geq0,\textrm{ }\forall k\in\left\{ 0,\ldots,p^{n}-1\right\} ,\textrm{ }\forall\mathfrak{z}\in\mathbb{Z}_{p}\label{eq:fourier}
\end{equation}
which expresses the indicator function of $k+p^{n}\mathbb{Z}_{p}$
as a linear combination of unitary characters and $p$-power roots
of unity. The condition on the characteristic of $R$ is needed in
order to guarantee that we can divide out by $p^{n}$. Since everything
in (\ref{eq:fourier}) occurs at a purely algebraic level, we can
use this identity to formally compute the Fourier transform of any
element of $\mathcal{S}\left(\mathbb{Z}_{p},R\right)$. Indeed, (\ref{eq:fourier})
tells us that:
\[
t\mapsto\frac{\mathbf{1}_{0}\left(p^{n}t\right)}{p^{n}}e^{-2\pi ikt}
\]
is the Fourier transform of $\left[\mathfrak{z}\overset{p^{n}}{\equiv}k\right]$.
As such, by linearity, the Fourier transform of any locally constant
function $f$ mod $p^{N}$ is $\hat{f}:\hat{\mathbb{Z}}_{p}\rightarrow\textrm{Frac}\left(R\right)\left(\zeta_{p^{\infty}}\right)$
given by:
\begin{equation}
\hat{f}\left(t\right)=\frac{\mathbf{1}_{0}\left(p^{N}t\right)}{p^{N}}\sum_{n=0}^{p^{N}-1}f\left(n\right)e^{-2\pi int}
\end{equation}
which vanishes for all $\left|t\right|_{p}>p^{N}$. Since $f$ is
locally constant mod $p^{N}$, the above sum is precisely the expected
Fourier integral:
\begin{equation}
\hat{f}\left(t\right)=\int_{\mathbb{Z}_{p}}f\left(\mathfrak{z}\right)e^{-2\pi i\left\{ t\mathfrak{z}\right\} _{p}}d\mathfrak{z},\textrm{ }\forall t\in\hat{\mathbb{Z}}_{p},\textrm{ }\forall f\in\mathcal{S}\left(\mathbb{Z}_{p},R\right)
\end{equation}
Since every SB function is a linear combination of finitely many indicator
functions of the form $\left[\mathfrak{z}\overset{p^{n}}{\equiv}k\right]$
for various values of $n$ and $k$, this shows that Fourier analysis
for functions in $\mathcal{S}\left(\mathbb{Z}_{p},R\right)$ be done
on a purely formal/algebraic basis, without any need to appeal to
analytical questions like topology or convergence. , Furthermore,
Fourier inversion always holds over $\mathcal{S}\left(\mathbb{Z}_{p},R\right)$:
\begin{equation}
f\left(\mathfrak{z}\right)=\sum_{t\in\hat{\mathbb{Z}}_{p}}\hat{f}\left(t\right)e^{2\pi i\left\{ t\mathfrak{z}\right\} _{p}},\textrm{ }\forall\mathfrak{z}\in\mathbb{Z}_{p},\textrm{ }\forall f\in\mathcal{S}\left(\mathbb{Z}_{p},R\right)
\end{equation}
simply because $\hat{f}$ has finite support in $\hat{\mathbb{Z}}_{p}$.

Just as local fields (resp. local rings) arise from taking metric
completions of global fields (resp., their rings of integers) with
respect to a given absolute value, we can get richer spaces of functions
out of the $R$-algebra $\mathcal{S}\left(\mathbb{Z}_{p},R\right)$
by considering its completions with respect to the norms/absolute
values on $\mathcal{S}\left(\mathbb{Z}_{p},R\right)$ induced by norms/absolute
values on $R$. If, for example, $R=\mathbb{Q}$, then given an absolute
value $\left|\cdot\right|_{q}$ on $\mathbb{Q}$, we get a norm:
\begin{equation}
\left\Vert f\right\Vert _{p,q}\overset{\textrm{def}}{=}\sup_{\mathfrak{z}\in\mathbb{Z}_{p}}\left|f\left(\mathfrak{z}\right)\right|_{q},\textrm{ }\forall f\in\mathcal{S}\left(\mathbb{Z}_{p},\mathbb{Q}\right)
\end{equation}
on $\mathcal{S}\left(\mathbb{Z}_{p},\mathbb{Q}\right)$. The completion
of $\mathcal{S}\left(\mathbb{Z}_{p},\mathbb{Q}\right)$ with respect
to this norm is then the Banach space $C\left(\mathbb{Z}_{p},\mathbb{Q}_{q}\right)$
of continuous functions $\mathbb{Z}_{p}\rightarrow\mathbb{Q}_{q}$,
where $\mathbb{Q}_{q}=\mathbb{R}$ in the case $\left|\cdot\right|_{q}$
is the real absolute value. In this way, we can obtain a theory of
Fourier analysis on $C\left(\mathbb{Z}_{p},\mathbb{Q}_{q}\right)$
simply by understanding how Fourier analysis on $\mathcal{S}\left(\mathbb{Z}_{p},\mathbb{Q}\right)$
behaves with respect to the completion: we compute things in $C\left(\mathbb{Z}_{p},\mathbb{Q}_{q}\right)$
by taking limits of the corresponding computations performed over
sequences in $\mathcal{S}\left(\mathbb{Z}_{p},\mathbb{Q}\right)$.
In the same way, we can go from Fourier analysis on $\mathcal{S}\left(\mathbb{Z}_{p},R\right)$
to Fourier analysis on $\mathcal{S}\left(\mathbb{Z}_{p},\textrm{Frac}\left(R_{\ell}\right)\right)$
where $\textrm{Frac}\left(R_{\ell}\right)$ is the field of fractions
of the completion of $R_{\ell}$ with respect to some absolute value
$\ell$.
\begin{rem}
Note that even without requiring convergence properties to vary from
point to point\textemdash as is the case with F-series\textemdash this
contruction (which is key to the success of my universal Fourier analysis)
is already deeply intertwined with the ideas of completions of rings
with respect to (semi)norms (and, more generally, with locally convex
topologies). This suggests that the Schwartz-Bruhat functions are
playing a functorial role in all this. Once again, this also makes
sense: distributions are simply the duals of SB functions, so it is
natural that a functorial realization of spaces of distributions would
dualize to a functorial realization of spaces of SB functions.
\end{rem}
With frames, we can perform the above procedures with even greater
freedom: we can allow the chosen completions to vary with respect
to $\mathfrak{z}$. Given an $R$-frame $\mathcal{F}$, the seminorms
produced by $\mathcal{F}$ induce a locally convex topology on $\mathcal{S}\left(\mathbb{Z}_{p},R\right)$,
and $C\left(\mathcal{F}\right)$ will be contained in the completion
of $\mathcal{S}\left(\mathbb{Z}_{p},R\right)$ with respect to this
topology; indeed, $C\left(\mathcal{F}\right)$ will be equal to the
completion precisely when $D\left(\mathcal{F}\right)=\mathbb{Z}_{p}$.
Just like with the case of going from $\mathcal{S}\left(\mathbb{Z}_{p},\mathbb{Q}\right)$
to $C\left(\mathbb{Z}_{p},\mathbb{Q}_{q}\right)$, all of our Fourier
analysis will follow by taking limits in $\mathcal{S}\left(\mathbb{Z}_{p},R\right)$.

It's natural to ask why one would bother using such a complicated
set-up. For an analyst like me, the answer to that is quite simple:
it allows us to meaningfully integrate functions like F-series which,
in general, fail to be integrable in the classical sense. This fact
is of independent interest because it shows that $\left(p,q\right)$-adic
analysis is much more nuanced than was previously thought. As I explain
in detail in \cite{2nd blog paper} (see also \cite{My Dissertation}),
in the case where $R=\mathbb{Q}_{q}$ (or a metrically complete extension
thereof) for an integer $q\geq2$ which is co-prime to $p$, a $\left(p,q\right)$-adic
function $f:\mathbb{Z}_{p}\rightarrow R$ is Haar-integrable \emph{\uline{if
and only if}} $f$ is continuous. (This result was originally discovered
by W.H. Schikhof in his PhD dissertation \cite{Schikhof's Thesis};
see also the final chapter of \cite{van Rooij - Non-Archmedean Functional Analysis},
though the latter is out of print.) It's easy to show that, with the
help of a frame to make sense of their convergence, many F-series
can be used to give representations of $\left(p,q\right)$-adic functions
which fail to be integrable in Schikhof's sense of the term, but which
are, as I would say, ``quasi-integrable''.
\begin{example}
Fix co-prime integers $p,q$ with $p\notin\left\{ -1,0,1\right\} $
and $q$ odd, and consider the $2$-adic F-series:
\begin{equation}
X_{q,p}\left(\mathfrak{z}\right)\overset{\textrm{def}}{=}\sum_{n=0}^{\infty}\frac{q^{\#_{2:1}\left(\left[\mathfrak{z}\right]_{2^{n}}\right)}}{p^{n}}
\end{equation}
Fix a prime divisor $\ell$ of $q$, and let $\mathcal{F}_{2,\ell}$
be the $\mathbb{Q}$-frame on $\mathbb{Z}_{2}$ with the seminorms
defined by:
\begin{equation}
\mathcal{F}_{2,\ell;\mathfrak{z}}\left(f\right)\overset{\textrm{def}}{=}\begin{cases}
\left|f\left(\mathfrak{z}\right)\right|_{\mathbb{R}} & \forall\mathfrak{z}\in\mathbb{N}_{0}\\
\left|f\left(\mathfrak{z}\right)\right|_{\ell} & \forall\mathfrak{z}\in\mathbb{Z}_{2}\backslash\mathbb{N}_{0}
\end{cases},\textrm{ }\forall f:\mathbb{Z}_{2}\rightarrow\mathbb{Q}
\end{equation}
and then extended naturally to $C\left(\mathcal{F}_{2,\ell}\right)$.
It's an elementary exercise in series convergence to show that $X_{q,p}\left(\mathfrak{z}\right)$
converges with respect to $\mathcal{F}_{2,\ell}$, and that the sum
of the F-series defines a function $X_{q,p}:\mathbb{Z}_{2}\rightarrow\mathbb{Z}_{\ell}$
satisfying the functional equation:
\begin{align}
X_{q,p}\left(2\mathfrak{z}\right) & =1+\frac{1}{p}X_{q,p}\left(\mathfrak{z}\right)\\
X_{q,p}\left(2\mathfrak{z}+1\right) & =1+\frac{q}{p}X_{q,p}\left(\mathfrak{z}\right)
\end{align}
for all $\mathfrak{z}\in\mathbb{Z}_{2}$. (Moreover, $X_{q,p}$ is
the unique rising-continuous $\left(2,\ell\right)$-adic function
satisfying these functional equations.) With this, we get:
\begin{align}
X_{q,p}\left(0\right) & =\frac{p}{p-1}\\
X_{q,p}\left(1\right) & =1+\frac{q}{p-1}
\end{align}
and:
\begin{align*}
X_{q,p}\left(2^{n}\right) & =\frac{1}{p^{n}}X_{q,p}\left(1\right)+\sum_{m=0}^{n-1}\frac{1}{p^{m}}\\
 & =\frac{1}{p^{n}}\left(1+\frac{q}{p-1}\right)+\frac{1-\frac{1}{p^{n}}}{1-\frac{1}{p}}\\
 & =\frac{p}{p-1}+\frac{1}{p^{n}}\frac{q-1}{p-1}
\end{align*}
If $X_{q,p}:\mathbb{Z}_{2}\rightarrow\mathbb{Z}_{\ell}$ were continuous,
the convergence of $2^{n}$ to $0$ in $\mathbb{Z}_{2}$ as $n\rightarrow\infty$
would force $X_{q,p}\left(2^{n}\right)$ to converge $\ell$-adically
to $X_{q,p}\left(0\right)$ as $n\rightarrow\infty$. However, since
$p$ and $q$ are co-prime and $\ell$ is a prime divisor of $q$,
we have that $p\in\mathbb{Z}_{\ell}^{\times}$, and as such, $\lim_{n\rightarrow\infty}X_{q,p}\left(2^{n}\right)$
does not exist in $\mathbb{Z}_{\ell}$. Thus, $X_{q,p}$ fails to
be continuous at $0$, and is therefore not an element of $C\left(\mathbb{Z}_{2},\mathbb{Q}_{\ell}\right)$.
Because of this, $X_{q,p}$ fails to be $\left(2,\ell\right)$-adically
integrable, in the sense that:
\[
\lim_{N\rightarrow\infty}\frac{1}{2^{N}}\sum_{n=0}^{2^{N}-1}X_{q,p}\left(n\right)
\]
fails to converge in $\mathbb{Q}_{\ell}$. (This can be shown directly
using the functional equations of $X_{q,p}$ to solve for the partial
sum recursively.) However, using \textbf{Theorem \ref{thm:The-breakdown-variety}},
we can compute $X_{q,p}$'s Fourier transform, and in particular,
compute its integral. Here:
\[
\hat{X}_{q,p}\left(t\right)=
\]
\begin{equation}
\hat{X}_{q,p}\left(t\right)=\begin{cases}
\min\left\{ 0,v_{2}\left(t\right)\right\} \prod_{n=0}^{-v_{2}\left(t\right)-1}\frac{1+qe^{-2\pi i2^{n}t}}{2p} & \textrm{if }\frac{q+1}{2p}=1\\
\frac{1}{1-\frac{q+1}{2p}}\prod_{n=0}^{-v_{2}\left(t\right)-1}\frac{1+qe^{-2\pi i2^{n}t}}{2p} & \textrm{if }\frac{q+1}{2p}\neq1
\end{cases},\textrm{ }\forall t\in\hat{\mathbb{Z}}_{2}
\end{equation}
where, recall, $v_{2}\left(t\right)=+\infty$ when $t=0$ and $v_{2}\left(t\right)<0$
for $t\in\hat{\mathbb{Z}}_{2}\backslash\left\{ 0\right\} $. Since
the integral of a function on $\mathbb{Z}_{2}$ is the function's
$0$th Fourier coefficient, we get:
\begin{equation}
\int_{\mathbb{Z}_{2}}X_{q,p}\left(\mathfrak{z}\right)d\mathfrak{z}=\begin{cases}
0 & \textrm{if }\frac{q+1}{2p}=1\\
\frac{1}{1-\frac{q+1}{2p}} & \textrm{if }\frac{q+1}{2p}\neq1
\end{cases}
\end{equation}
a computation that cannot be done using Schikhof's theory of integration
(the Monna-Springer integral).
\end{example}
To that end, in my dissertation, I introduced the following definition,
though the version given here is much more general, thanks to the
maturity of my theory of frames in the interim.
\begin{defn}
Let $\mathcal{F}$ be an $R$-frame on $\mathbb{Z}_{p}$, and let
$X\in C\left(\mathcal{F}\right)$. We say $X$ is \textbf{$\mathcal{F}$-quasi-integrable
}/ \textbf{quasi-integrable with respect to $\mathcal{F}$ }whenever
there exists a function $\hat{X}:\hat{\mathbb{Z}}_{p}\rightarrow\textrm{Frac}\left(R\right)\left(\zeta_{p^{\infty}}\right)$
so that:
\begin{equation}
\lim_{N\rightarrow\infty}\sum_{\left|t\right|_{p}\leq p^{N}}\hat{X}\left(t\right)e^{2\pi i\left\{ t\mathfrak{z}\right\} _{p}}\overset{\mathcal{F}}{=}X\left(\mathfrak{z}\right)\label{eq:q.i.}
\end{equation}
More generally, we say a function is \textbf{quasi-integrable} when
it is quasi-integrable with respect to some frame.
\end{defn}
As the set of quasi-integrable functions with respect to a given frame
is rather messy to deal with, it helps that we can make sense of such
functions using the more orthodox constructions of distributions and
measures. The key player here is the Parseval-Plancherel identity:
\begin{equation}
\int_{\mathbb{Z}_{p}}f\left(\mathfrak{z}\right)g\left(\mathfrak{z}\right)d\mathfrak{z}=\sum_{t\in\hat{\mathbb{Z}}_{p}}\hat{f}\left(t\right)\hat{g}\left(-t\right),\textrm{ }\forall f,g\in\mathcal{S}\left(\mathbb{Z}_{p},\textrm{Frac}\left(R\right)\left(\zeta_{p^{\infty}}\right)\right)
\end{equation}
which we use to obtain a distribution $X\left(\mathfrak{z}\right)d\mathfrak{z}$
from a quasi-integrable function $X$ by way of the formula:
\begin{equation}
\int_{\mathbb{Z}_{p}}f\left(\mathfrak{z}\right)X\left(\mathfrak{z}\right)d\mathfrak{z}=\sum_{t\in\hat{\mathbb{Z}}_{p}}\hat{f}\left(t\right)\hat{X}\left(-t\right),\textrm{ }\forall f\in\mathcal{S}\left(\mathbb{Z}_{p},\textrm{Frac}\left(R\right)\left(\zeta_{p^{\infty}}\right)\right)
\end{equation}
If there is an absolute value $\left|\cdot\right|_{\ell}\in\mathscr{V}\left(R\left(\zeta_{p^{\infty}}\right)\right)$
so that $\sup_{t\in\hat{\mathbb{Z}}_{p}}\left|\hat{X}\left(t\right)\right|_{\ell}<\infty$,
then the above formula extends to define $X\left(\mathfrak{z}\right)d\mathfrak{z}$
as an element of the continuous dual of the Wiener algebra $W\left(\mathbb{Z}_{p},\overline{\overline{\textrm{Frac}\left(R_{\ell}\right)}}\right)$
(and hence, as a measure), where $\overline{\overline{\textrm{Frac}\left(R_{\ell}\right)}}$
is the metric completion of the algebraic closure of the field of
fractions of the completion of $R$ with respect to $\ell$. This
is advantageous, because it allows us to do Fourier analysis with
a much larger class of $\left(p,q\right)$-adic functions than are
available to us through the classical theory (viz. \cite{Schikhof's Thesis,van Rooij - Non-Archmedean Functional Analysis}).
Not only that, but, as I demonstrate in \cite{2nd blog paper}, realizing
quasi-integrable functions as distributions/measures is actually quite
natural: if $X$ is quasi-integrable, with a Fourier transform $\hat{X}$,
and if we write $d\mu$ to denote the distribution (or measure) whose
Fourier-Stieltjes transform is $\hat{X}$, then (\ref{eq:q.i.}) is
formally equivalent to the statement that $X\left(\mathfrak{z}\right)$
is the \textbf{Radon-Nikodym derivative} of the distribution/measure
$d\mu$ with respect to $d\mathfrak{z}$. Indeed, if $d\mu$ is a
real or complex-valued measure on $\mathbb{Z}_{p}$ with Fourier-Stieltjes
transform $\hat{X}$, then if the limit in (\ref{eq:q.i.}) converges
in $\mathbb{C}$, say, for almost every $\mathfrak{z}\in\mathbb{Z}_{p}$,
that limit is precisely the Radon-Nikodym derivative of $d\mu$ in
the classical sense of the term (ex: \cite{Folland - real analysis})!
So, independent of its applications to number theory, my discovery
of the phenomenon of quasi-integrability in \cite{My Dissertation}
shows that there's a great deal of unexplored structure in measures
valued in non-archimedean spaces ripe for exploration.

However, as the saying goes, all magic comes with a price. Consequently,
it should come as no surprise that quasi-integrability has a price
of its own, one that shows a significant departure from classical
Fourier analysis. To see how this manifests, let us return to our
work from \textbf{Section \ref{subsec:The-Central-Computation}}.
\begin{example}
Now that we have a rigorous understanding of quasi-integrability,
I can specify exactly what I meant by computing a Fourier transform
for the $X_{\mathbf{m}}$s as per \textbf{Section \ref{subsec:The-Central-Computation}}.
Here, we work with the set-up of \textbf{Theorem \ref{thm:variety frame general}}.
Our starting point is \textbf{Proposition \ref{prop:X_3-hat formal solution}}.
Recall, this said that if $\hat{X}_{\mathbf{m}}$ was known for all
$\mathbf{m}<\mathbf{n}$ for a given $\mathbf{n}$, attempting to
use the formal equation (\ref{eq:formal solution}): 
\begin{equation}
\hat{\chi}\left(t\right)=\alpha_{\mathbf{n}}\left(t\right)\hat{\chi}\left(pt\right)+\sum_{\mathbf{m}<\mathbf{n}}\alpha_{\mathbf{m},\mathbf{n}}\left(t\right)\hat{X}_{\mathbf{m}}\left(pt\right)
\end{equation}
to solve for the unknown function $\hat{\chi}$ as a means of obtaining
a Fourier transform of $\hat{X}_{\mathbf{n}}$ could result in one
of three outcomes:

I. $\alpha_{\mathbf{n}}\left(0\right)\neq1$, in which case the solution
$\hat{\chi}$ was uniquely determined by the $\alpha_{\mathbf{m},\mathbf{n}}$s
and by the $\hat{X}_{\mathbf{m}}$s for $\mathbf{m}<\mathbf{n}$.

II. $\alpha_{\mathbf{n}}\left(0\right)=1$ and $\sum_{\mathbf{m}<\mathbf{n}}\alpha_{\mathbf{m},\mathbf{n}}\left(0\right)\hat{X}_{\mathbf{m}}\left(0\right)\neq0$,
in which case no solution $\hat{\chi}$ exists.

III. $\alpha_{\mathbf{n}}\left(0\right)=1$ and $\sum_{\mathbf{m}<\mathbf{n}}\alpha_{\mathbf{m},\mathbf{n}}\left(0\right)\hat{X}_{\mathbf{m}}\left(0\right)=0$,
in which case infinitely many solutions $\hat{\chi}$ exist, and all
of these lay in a $1$-dimensional affine subspace of the space of
functions $\hat{\mathbb{Z}}_{p}\rightarrow\textrm{Frac}\left(\mathcal{R}\right)\left(\zeta_{p^{\infty}}\right)$.
In particular, for any two solutions $\hat{\chi}$, $\hat{\chi}^{\prime}$,
we had that: 
\begin{equation}
\hat{\chi}\left(t\right)-\hat{\chi}^{\prime}\left(t\right)=\left(\hat{\chi}\left(0\right)-\hat{\chi}^{\prime}\left(0\right)\right)\hat{A}_{\mathbf{n}}\left(t\right)
\end{equation}
where:
\begin{equation}
\hat{A}_{\mathbf{n}}\left(t\right)=\prod_{n=0}^{-v_{p}\left(t\right)-1}\alpha_{\mathbf{n}}\left(p^{n}t\right)
\end{equation}

So, considering the third case, where $\alpha_{\mathbf{n}}\left(0\right)=1$
and $\sum_{\mathbf{m}<\mathbf{n}}\alpha_{\mathbf{m},\mathbf{n}}\left(0\right)\hat{X}_{\mathbf{m}}\left(0\right)=0$,
let $\hat{\chi}$ and $\hat{\chi}^{\prime}$ be two distinct solutions
to (\ref{eq:formal solution}). Then, letting $c$ denote the quantity
$\hat{\chi}\left(0\right)-\hat{\chi}^{\prime}\left(0\right)$, we
have:
\begin{align*}
\sum_{\left|t\right|_{p}\leq p^{N}}\left(\hat{\chi}\left(t\right)-\hat{\chi}^{\prime}\left(t\right)\right)e^{2\pi i\left\{ t\mathfrak{z}\right\} _{p}} & =c\sum_{\left|t\right|_{p}\leq p^{N}}\hat{A}_{\mathbf{n}}\left(t\right)e^{2\pi i\left\{ t\mathfrak{z}\right\} _{p}}\\
\left(\mathbf{Proposition}\textrm{ }\mathbf{\ref{prop:A_X sum}}\right); & =r_{\mathbf{n},0}^{N}\kappa_{\mathbf{n}}\left(\left[\mathfrak{z}\right]_{p^{N}}\right)+\underbrace{\left(1-\alpha_{\mathbf{n}}\left(0\right)\right)}_{0}\sum_{n=0}^{N-1}r_{\mathbf{n},0}^{n}\kappa_{\mathbf{n}}\left(\left[\mathfrak{z}\right]_{p^{n}}\right)
\end{align*}
Now, suppose there is an $\mathcal{R}$-frame $\mathcal{F}$ on $\mathbb{Z}_{p}$
so that both of the limits:
\begin{align}
\lim_{N\rightarrow\infty}\sum_{\left|t\right|_{p}\leq p^{N}}\hat{\chi}\left(t\right)e^{2\pi i\left\{ t\mathfrak{z}\right\} _{p}}\\
\lim_{N\rightarrow\infty}\sum_{\left|t\right|_{p}\leq p^{N}}\hat{\chi}^{\prime}\left(t\right)e^{2\pi i\left\{ t\mathfrak{z}\right\} _{p}}
\end{align}
converge with respect to $\mathcal{F}$ to functions $\chi\left(\mathfrak{z}\right)$
and $\chi^{\prime}\left(\mathfrak{z}\right)$, respectively. If, in
addition:
\begin{equation}
\lim_{N\rightarrow\infty}r_{\mathbf{n},0}^{N}\kappa_{\mathbf{n}}\left(\left[\mathfrak{z}\right]_{p^{N}}\right)\overset{\mathcal{F}}{=}0
\end{equation}
we then have: 
\begin{equation}
\chi\left(\mathfrak{z}\right)-\chi^{\prime}\left(\mathfrak{z}\right)\overset{\mathcal{F}}{=}\lim_{N\rightarrow\infty}\sum_{\left|t\right|_{p}\leq p^{N}}\left(\hat{\chi}\left(t\right)-\hat{\chi}^{\prime}\left(t\right)\right)e^{2\pi i\left\{ t\mathfrak{z}\right\} _{p}}\overset{\mathcal{F}}{=}\lim_{N\rightarrow\infty}r_{\mathbf{n},0}^{N}\kappa_{\mathbf{n}}\left(\left[\mathfrak{z}\right]_{p^{N}}\right)\overset{\mathcal{F}}{=}0
\end{equation}
This shows that $\chi$ and $\chi^{\prime}$ are the same function,
even though their respective Fourier series were generated by two
completely different formulae ($\hat{\chi}$ and $\hat{\chi}^{\prime}$,
respectively). In this way, we see that \emph{both} $\hat{\chi}$
and $\hat{\chi}^{\prime}$ can be called ``Fourier transforms''
of $\chi$, in the sense that they both generate Fourier series that
sum to $\chi$ with respect to $\mathcal{F}$. However, the distributions
$d\mu$ and $d\mu^{\prime}$ induced by $\hat{\chi}$ and $\hat{\chi}^{\prime}$,
respectively via the Parseval-Plancherel construction are NOT equal!
Thus, we can realize $\chi$ as a distribution in \emph{two completely
different ways}: one, by identifying it with $d\mu$, the distribution
with $\hat{\chi}$ as its Fourier-Stieltjes transform; the other,
by identifying it withe $d\mu^{\prime}$, the measure with $\hat{\chi}^{\prime}$
as its Fourier-Stieltjes transform. In fact, for this frame particular,
\emph{every }element of the $1$-dimensional affine subspace of solutions
of (\ref{eq:formal solution}) will be a Fourier transform of $\chi$,
giving us a $1$-dimensional affine space of distributions that may
be used as realizations of $\chi$.

In summary, the quasi-integrable function $\chi$ has infinitely many
possible Fourier transforms, which may be distinguished from one another
by considering the different distributions they induced via the Parseval-Plancherel
construction.
\end{example}
Unfortunately, things get worse.
\begin{example}
\label{exa:even worse}Fix $\mathcal{F}$ and $\mathbf{n}$ so that
$\alpha_{\mathbf{n}}\left(0\right)=1$ and: 
\begin{equation}
\lim_{N\rightarrow\infty}r_{\mathbf{n},0}^{N}\kappa_{\mathbf{n}}\left(\left[\mathfrak{z}\right]_{p^{N}}\right)\overset{\mathcal{F}}{=}0
\end{equation}
Then, by \textbf{Proposition \ref{prop:A_X sum}},\textbf{ }$\alpha_{\mathbf{n}}\left(0\right)=1$
implies: 
\begin{equation}
\sum_{\left|t\right|_{p}\leq p^{N}}\hat{A}_{\mathbf{n}}\left(t\right)e^{2\pi i\left\{ t\mathfrak{z}\right\} _{p}}=r_{\mathbf{n},0}^{N}\kappa_{\mathbf{n}}\left(\left[\mathfrak{z}\right]_{p^{N}}\right),\textrm{ }\forall N,\mathfrak{z}
\end{equation}
and so:
\begin{equation}
\lim_{N\rightarrow\infty}\sum_{\left|t\right|_{p}\leq p^{N}}\hat{A}_{\mathbf{n}}\left(t\right)e^{2\pi i\left\{ t\mathfrak{z}\right\} _{p}}\overset{\mathcal{F}}{=}0
\end{equation}
Now, let $\chi\in C\left(\mathcal{F}\right)$ be \emph{any} $\mathcal{F}$-quasi-integrable
function, and let $\hat{\chi}$ be a Fourier transform of $\chi$.
Then, for any constant $c$, $\hat{\chi}+c\hat{A}_{\mathbf{n}}$ will
\emph{also }be a Fourier transform of $\chi$:
\begin{equation}
\lim_{N\rightarrow\infty}\sum_{\left|t\right|_{p}\leq p^{N}}\left(\hat{\chi}\left(t\right)+c\hat{A}_{\mathbf{n}}\left(t\right)\right)e^{2\pi i\left\{ t\mathfrak{z}\right\} _{p}}\overset{\mathcal{F}}{=}\lim_{N\rightarrow\infty}\sum_{\left|t\right|_{p}\leq p^{N}}\hat{\chi}\left(t\right)e^{2\pi i\left\{ t\mathfrak{z}\right\} _{p}}\overset{\mathcal{F}}{=}\chi\left(\mathfrak{z}\right)
\end{equation}
Thus, the non-uniqueness of the Fourier transform of a quasi-integrable
function is not a mere quirk of the $\hat{X}_{\mathbf{n}}$s, but
rather a defining feature of the theory as a whole. In embracing frames
to broaden our Fourier-analytic horizons, we have lost the uniqueness
of a function's Fourier transform, in stark contrast to classical
Fourier analysis.
\end{example}
I first made this disturbing observation in my doctoral dissertation.
I expended a great deal of effort trying to prove it away before I
gave in and accepted this state of affairs. This led me to coin the
following terminology:
\begin{defn}
\label{def:degen}Let $R$ be an integral domain, and let $\mathcal{F}$
be an $R$-frame on $\mathbb{Z}_{p}$. I say a distribution $d\mu\in\mathcal{S}\left(\mathbb{Z}_{p},\textrm{Frac}\left(R\right)\right)^{\prime}$
is \textbf{$\mathcal{F}$-degenerate }(or simply \textbf{degenerate},\textbf{
}for short) whenever:
\begin{equation}
\lim_{N\rightarrow\infty}\sum_{\left|t\right|_{p}\leq p^{N}}\hat{\mu}\left(t\right)e^{2\pi i\left\{ t\mathfrak{z}\right\} _{p}}\overset{\mathcal{F}}{=}0
\end{equation}
where $\hat{\mu}$ is the Fourier-Stieltjes transform\textbf{ }of
$d\mu$.
\end{defn}
\begin{notation}
Given a function $\hat{\phi}:\hat{\mathbb{Z}}_{p}\rightarrow\textrm{Frac}\left(R\right)\left(\zeta_{p^{\infty}}\right)$,
we write:
\begin{equation}
\left\Vert \hat{\phi}\right\Vert _{p,\ell}\overset{\textrm{def}}{=}\sup_{t\in\hat{\mathbb{Z}}_{p}}\left|\hat{\phi}\left(t\right)\right|_{\ell}
\end{equation}
\end{notation}
\begin{rem}
Give a distribution $d\mu\in\mathcal{S}\left(\mathbb{Z}_{p},\textrm{Frac}\left(R\right)\right)^{\prime}$,
if there is an absolute value $\ell\in\mathscr{V}\left(R\right)$
so that $\left\Vert \hat{\mu}\right\Vert _{p,\ell}<\infty$, then
$d\mu$ extends to a $\left(p,\ell\right)$-adic measure $d\mu\in W\left(\mathbb{Z}_{p},\overline{\overline{\textrm{Frac}\left(R_{\ell}\right)}}\right)^{\prime}$
defined by the Parseval-Plancherel construction. Hence, we can apply
\textbf{Definition \ref{def:degen}} to measures, in addition to distributions.
As such, in a minor abuse of notation, I will use the term ``degenerate
measure'', even when talking distributions.
\end{rem}
Though I failed to get rid of the non-uniqueness, my consolation prize
was my realization that this non-uniqueness is relatively manageable.
\begin{prop}
Let $f\in C\left(\mathcal{F}\right)$ be $\mathcal{F}$-quasi-integrable,
and let $\hat{f}$ be any Fourier transform of $f$. Then, for any
$\mathcal{F}$-degenerate measure $d\mu$, $\hat{f}+\hat{\mu}$ is
a Fourier transform of $f$ as well.
\end{prop}
Proof:
\begin{align*}
\lim_{N\rightarrow\infty}\sum_{\left|t\right|_{p}\leq p^{N}}\left(\hat{f}\left(t\right)+\hat{\mu}\left(t\right)\right)e^{2\pi i\left\{ t\mathfrak{z}\right\} _{p}} & \overset{\mathcal{F}}{=}\lim_{N\rightarrow\infty}\sum_{\left|t\right|_{p}\leq p^{N}}\hat{f}\left(t\right)e^{2\pi i\left\{ t\mathfrak{z}\right\} _{p}}+\lim_{N\rightarrow\infty}\sum_{\left|t\right|_{p}\leq p^{N}}\hat{\mu}\left(t\right)e^{2\pi i\left\{ t\mathfrak{z}\right\} _{p}}\\
 & \overset{\mathcal{F}}{=}f\left(\mathfrak{z}\right)+0
\end{align*}

Q.E.D.

\vphantom{}

Because of this, we see the following:
\begin{thm}
Let $\mathcal{F}$ be an $R$-frame on $\mathbb{Z}_{p}$. Then, the
set of functions $\hat{\mathbb{Z}}_{p}\rightarrow\textrm{Frac}\left(R\right)\left(\zeta_{p^{\infty}}\right)$
which are Fourier transforms of $\mathcal{F}$-quasi-integrable functions
is affine subspace of the $\textrm{Frac}\left(R\right)$-vector space
of functions $\hat{\mathbb{Z}}_{p}\rightarrow\textrm{Frac}\left(R\right)\left(\zeta_{p^{\infty}}\right)$,
with the difference of any two Fourier transforms being an $\mathcal{F}$-degenerate
measure.

Consequently, if we let $\textrm{Degen}\left(\mathcal{F}\right)$
denote the $\textrm{Frac}\left(R\right)$-vector space of functions
$\hat{\mathbb{Z}}_{p}\rightarrow\textrm{Frac}\left(R\right)\left(\zeta_{p^{\infty}}\right)$
which are Fourier-Stieltjes transforms of $\mathcal{F}$-degenerate
measures, we have that the Fourier transform of an $\mathcal{F}$-quasi-integrable
function is a unique element of the quotient space:
\begin{equation}
\left\{ \hat{f}:\hat{\mathbb{Z}}_{p}\rightarrow\textrm{Frac}\left(R\right)\left(\zeta_{p^{\infty}}\right)\right\} /\textrm{Degen}\left(\mathcal{F}\right)
\end{equation}
\end{thm}

\section{\label{sec:Making-Things-Precise}Making Things Precise}

Now that we have established the conceptual backdrop needed to understand
what was happening in \textbf{Section \ref{subsec:The-Central-Computation}},
we can proceed with the task of justifying everything that was done.
The first two subsections of this section are devoted to establishing
notation and technical results needed to make the estimates and asymptotics
that \textbf{Section \ref{subsec:Assumption--=000026}} will use in
its proofs of the results that justify the work of \textbf{Section
\ref{subsec:The-Central-Computation}}.

\subsection{\label{subsec:M-functions-=000026-Asymptotics}Asymptotics and M-functions}
\begin{notation}
We use Vinogradov notation. Given sequences $x_{n}$ and $y_{n}$
of non-negative real numbers, we write:
\begin{equation}
x_{n}\ll y_{n}\textrm{ as }n\rightarrow\infty
\end{equation}
to mean there exists a real constant $C>0$ so that:
\begin{equation}
x_{n}\leq Cy_{n}
\end{equation}
holds for all sufficiently large $n$.

When working with non-negative real-valued functions $f_{n}\left(\mathfrak{z}\right),g_{n}\left(\mathfrak{z}\right)$,
we write:
\begin{equation}
f_{n}\left(\mathfrak{z}\right)\ll_{\mathfrak{z}}g_{n}\left(\mathfrak{z}\right)\textrm{ as }n\rightarrow\infty,\textrm{ }\forall\mathfrak{z}\in Z
\end{equation}
to indicate that, for each $\mathfrak{z}$ in the set $Z$, there
is a real constant $C_{\mathfrak{z}}>0$, depending on $\mathfrak{z}$,
so that:
\begin{equation}
f_{n}\left(\mathfrak{z}\right)\leq C_{\mathfrak{z}}g_{n}\left(\mathfrak{z}\right)
\end{equation}
holds for all sufficiently large $n$.
\end{notation}
As we will see, the key to making the limits in \textbf{Section \ref{subsec:The-Central-Computation}}
rigorous lies in proving through an inductive-recursive procedure
that the decay rate of the tail of the Fourier series of $X_{I}$
for $I\subseteq\mathcal{J}$ gets passed on to the decay rate of the
tail of the Fourier series of $X_{J}$, where $J\subseteq\mathcal{J}$
contains $I$ as a proper subset. To do this, we need to be able to
show, for example, that for real-valued M-functions $\left\{ E_{n}\right\} _{n\geq0},\left\{ F_{n}\right\} _{n\geq0}$
and a $\mathfrak{z}\in\mathbb{Z}_{p}$ so that:
\begin{equation}
\lim_{n\rightarrow\infty}\left(E_{n}\left(\mathfrak{z}\right)+F_{n}\left(\mathfrak{z}\right)\right)\overset{\mathbb{R}}{=}0
\end{equation}
then there exists a real-valued M-function $\left\{ G_{n}\right\} _{n\geq0}$
so that $E_{n}\left(\mathfrak{z}\right)+F_{n}\left(\mathfrak{z}\right)\ll_{\mathfrak{z}}G_{n}\left(\mathfrak{z}\right)$
as $n\rightarrow\infty$. There are several other combinations of
this type that need to be dealt with.

To make things work, we will restrict the $\mathfrak{z}$ under consideration
to belong to a subset of $\mathbb{Z}_{p}$ so that $\lim_{n\rightarrow\infty}M_{n}^{1/n}\left(\mathfrak{z}\right)$
exists in $\mathbb{R}$, for all real-valued M-functions $\left\{ M_{n}\right\} _{n\geq0}$
Thus, if there are non-negative real numbers $\lambda_{E}\left(\mathfrak{z}\right)$
and $\lambda_{F}\left(\mathfrak{z}\right)$ strictly less than $1$
for which:
\begin{align}
\lim_{n\rightarrow\infty}E_{n}^{1/n}\left(\mathfrak{z}\right) & \overset{\mathbb{R}}{=}\lambda_{E}\left(\mathfrak{z}\right)\\
\lim_{n\rightarrow\infty}F_{n}^{1/n}\left(\mathfrak{z}\right) & \overset{\mathbb{R}}{=}\lambda_{F}\left(\mathfrak{z}\right)
\end{align}
We will have a bound:
\begin{equation}
E_{n}\left(\mathfrak{z}\right)+F_{n}\left(\mathfrak{z}\right)\ll_{\mathfrak{z}}\lambda_{E}^{n}\left(\mathfrak{z}\right)+\lambda_{F}^{n}\left(\mathfrak{z}\right)
\end{equation}
for all sufficiently large $n$, and thus, that $E_{n}\left(\mathfrak{z}\right)+F_{n}\left(\mathfrak{z}\right)$
will be bounded by an M-function $\left\{ G_{n}\right\} _{n\geq0}$
for which $\lim_{n\rightarrow\infty}G_{n}^{1/n}\left(\mathfrak{z}\right)=\lambda_{G}\left(\mathfrak{z}\right)$
is a real number less than $1$ but greater than either $\lambda_{E}\left(\mathfrak{z}\right)$
or $\lambda_{F}\left(\mathfrak{z}\right)$. This same argument will
work for all the cases we need.

We begin with a technical lemma:
\begin{lem}
\label{lem:ll lemma}Let $\left\{ x_{n}\right\} _{n\geq0}$ and $\left\{ y_{n}\right\} _{n\geq0}$
be sequences of positive real numbers tending to $0$. Suppose that
the quantities:
\begin{equation}
x\overset{\textrm{def}}{=}\lim_{n\rightarrow\infty}x_{n}^{1/n}
\end{equation}
\begin{equation}
y\overset{\textrm{def}}{=}\lim_{n\rightarrow\infty}y_{n}^{1/n}
\end{equation}
exist in $\mathbb{R}$ and are non-zero.

Now, let $\gamma:\mathbb{N}_{0}\rightarrow\mathbb{N}_{0}$ be any
function satisfying:
\begin{align}
\gamma\left(n\right) & \leq n,\textrm{ }\forall n\geq0\\
\lim_{n\rightarrow\infty}\gamma\left(n\right) & \overset{\mathbb{R}}{=}+\infty\\
\lim_{n\rightarrow\infty}\frac{\gamma\left(n\right)}{n} & \overset{\mathbb{R}}{=}r\in\left(0,1\right]
\end{align}
Then, $x^{\gamma\left(n\right)}\ll y^{n}$ as $n\rightarrow\infty$
implies $x_{\gamma\left(n\right)}\ll y_{n}$ as $n\rightarrow\infty$.
\end{lem}
Proof: Let everything be as given. By way of contradiction, suppose
that $x^{\gamma\left(n\right)}\ll y^{n}$ and yet, for any real constant
$C>0$, that there are infinitely many $n$ for which $x_{\gamma\left(n\right)}>Cy_{n}$.
Fixing such a $C$, we can pick a strictly increasing sequence of
positive integers $\left\{ n_{\ell}\right\} _{\ell\geq0}$ so that
$x_{\gamma\left(n_{\ell}\right)}>Cy_{n_{\ell}}$ holds for all $\ell$.
Taking $\gamma\left(n_{\ell}\right)$th roots from both sides, we
get:
\begin{equation}
x_{\gamma\left(n_{\ell}\right)}^{1/\gamma\left(n_{\ell}\right)}>C^{1/\gamma\left(n_{\ell}\right)}y_{n_{\ell}}^{1/\gamma\left(n_{\ell}\right)},\textrm{ }\forall\ell\geq0\label{eq:ell limit}
\end{equation}
Here:
\begin{equation}
\lim_{\ell\rightarrow\infty}y_{n_{\ell}}^{1/\gamma\left(n_{\ell}\right)}=\lim_{\ell\rightarrow\infty}\left(y_{n_{\ell}}^{1/n_{\ell}}\right)^{n_{\ell}/\gamma\left(n_{\ell}\right)}=y^{1/r}
\end{equation}
Since $C$ is positive and $\gamma\left(n\right)\rightarrow\infty$
as $n\rightarrow\infty$, $C^{1/\gamma\left(n_{\ell}\right)}$ tends
to $1$ as $\ell\rightarrow\infty$, so letting $\ell\rightarrow\infty$
in (\ref{eq:ell limit}), we obtain $x>y^{1/r}$.

Now, if $x=1$, then $y$ must be strictly less than $1$, hence $y^{n}\rightarrow0$,
and so $x^{\gamma\left(n\right)}\ll y^{n}$ forces $x=0$, which is
impossible. Thus, $x$ must be $<1$, in which case $x^{\gamma\left(n\right)}\ll y^{n}$
implies there is a $C^{\prime}>0$ so that $x^{\gamma\left(n\right)}\leq C^{\prime}y^{n}$
for all sufficiently large $n$, say $n\geq N_{C^{\prime}}$. But
then, choosing $n\geq N_{C^{\prime}}$ yields:
\begin{equation}
C^{\prime}y^{n}\geq x^{\gamma\left(n\right)}=\left(\frac{x^{\frac{\gamma\left(n\right)}{n}}}{y}\right)^{n}y^{n}=\left(x^{\frac{\gamma\left(n\right)}{n}-r}\times\frac{x^{r}}{y}\right)^{n}y^{n}
\end{equation}
Since $\gamma\left(n\right)/n\sim r$ as $n\rightarrow\infty$, and
since $\gamma\left(n\right)\leq n$, upon letting $\epsilon>0$ be
arbitrarily small, we can and do choose an $N_{\epsilon}$ larger
than $N_{C^{\prime}}$ so that $-\epsilon<\gamma\left(n\right)-rn<\epsilon$
for all $n\geq N_{\epsilon}$. Since $x$ is non-zero, it must be
between $0$ and $1$, and so:
\begin{equation}
x^{\epsilon/n}<x^{\frac{\gamma\left(n\right)}{n}-r}<x^{-\epsilon/n}
\end{equation}
which gives:
\begin{equation}
C^{\prime}y^{n}>\left(x^{\frac{\gamma\left(n\right)}{n}-r}\times\frac{x^{r}}{y}\right)^{n}y^{n}>x^{\epsilon}\left(\frac{x^{r}}{y}\right)^{n}y^{n}
\end{equation}
In particular:
\begin{equation}
C^{\prime}>x^{\epsilon}\left(\frac{x}{y^{1/r}}\right)^{rn}\label{eq:C lower bound}
\end{equation}
Since $x>y^{1/r}$, $\left(x/y^{1/r}\right)^{r}>1$, and thus, the
lower bound on the far right of (\ref{eq:C lower bound}) tends to
$\infty$ as $n\rightarrow\infty$, which is impossible, as $C^{\prime}<\infty$.
This is the desired contradiction.

Thus, if $x^{\gamma\left(n\right)}\ll y^{n}$ as $n\rightarrow\infty$,
it must be that $x_{\gamma\left(n\right)}\ll y_{n}$ as $n\rightarrow\infty$.

Q.E.D.

\vphantom{}In order to guarantee that $\lim_{n\rightarrow\infty}M_{n}^{1/n}\left(\mathfrak{z}\right)$
exists, we will need to use the notion of the \textbf{density }of
the digits of a $p$-adic integer.
\begin{defn}
Let $\mathfrak{z}\in\mathbb{Z}_{p}$. Then, for each $k\in\left\{ 0,\ldots,p-1\right\} $,
the \textbf{upper }(resp. \textbf{lower}) \textbf{density of $k$s
digits in $\mathfrak{z}$} are defined by:
\begin{equation}
\overline{d_{p:k}}\left(\mathfrak{z}\right)\overset{\textrm{def}}{=}\limsup_{n\rightarrow\infty}\frac{\#_{p:k}\left(\left[\mathfrak{z}\right]_{p^{n}}\right)}{n}
\end{equation}
and:
\begin{equation}
\underline{d_{p:k}}\left(\mathfrak{z}\right)\overset{\textrm{def}}{=}\liminf_{n\rightarrow\infty}\frac{\#_{p:k}\left(\left[\mathfrak{z}\right]_{p^{n}}\right)}{n}
\end{equation}
respectively. The \textbf{(natural) density of $k$s digits in $\mathfrak{z}$}
is defined by the limit:
\begin{equation}
d_{p:k}\left(\mathfrak{z}\right)\overset{\textrm{def}}{=}\lim_{n\rightarrow\infty}\frac{\#_{p:k}\left(\left[\mathfrak{z}\right]_{p^{n}}\right)}{n}
\end{equation}
whenever it exists.

It is well-known that, given a $k\in\left\{ 0,\ldots,p-1\right\} $,
there exist infinitely many $\mathfrak{z}\in\mathbb{Z}_{p}$ so that
$d_{p:k}\left(\mathfrak{z}\right)$ does not exist. As such, we write
$\textrm{NiceDig}\left(\mathbb{Z}_{p}\right)$\textbf{ }to denote
the set of all $\mathfrak{z}\in\mathbb{Z}_{p}$ so that $d_{p:k}\left(\mathfrak{z}\right)$
exists for all $k\in\left\{ 0,\ldots,p-1\right\} $. (As the symbol
suggests, this is the set of $p$-adic integers whose digits are ``nice''.)
\end{defn}
The key properties of the digits of $\mathfrak{z}$s in $\textrm{NiceDig}\left(\mathbb{Z}_{p}\right)$
are chronicled below.
\begin{prop}
\label{prop:nice digits}Let $\mathfrak{z}\in\textrm{NiceDig}\left(\mathbb{Z}_{p}\right)$.
Then:

I. $\theta_{p}\left(\mathfrak{z}\right)\in\textrm{NiceDig}\left(\mathbb{Z}_{p}\right)$,
and $d_{p:k}\left(\theta_{p}\left(\mathfrak{z}\right)\right)=d_{p:k}\left(\mathfrak{z}\right)$
for all $k\in\left\{ 0,\ldots,p-1\right\} $.

II. Let $\gamma:\mathbb{N}_{0}\rightarrow\mathbb{N}_{0}$ be any function
satisfying:
\begin{align}
\gamma\left(n\right) & \leq n,\textrm{ }\forall n\geq0\\
\lim_{n\rightarrow\infty}\gamma\left(n\right) & \overset{\mathbb{R}}{=}+\infty\\
\lim_{n\rightarrow\infty}\frac{\gamma\left(n\right)}{n} & \overset{\mathbb{R}}{=}L
\end{align}
Then: 
\begin{equation}
\lim_{n\rightarrow\infty}\frac{1}{n}\#_{p:k}\left(\left[\theta_{p}^{\circ\gamma\left(n\right)}\left(\mathfrak{z}\right)\right]_{p^{n-\gamma\left(n\right)}}\right)=\left(1-L\right)d_{p:k}\left(\mathfrak{z}\right)
\end{equation}
\end{prop}
Proof: For both (I) and (II), we use \textbf{Proposition \ref{prop:fundamental functional equations-1}}
to deal with $\theta_{p}$.

I. Letting $k$ be arbitrary, we have:
\begin{align*}
\lim_{n\rightarrow\infty}\frac{\#_{p:k}\left(\left[\theta_{p}\left(\mathfrak{z}\right)\right]_{p^{n}}\right)}{n} & =\lim_{n\rightarrow\infty}\left(\frac{\#_{p:k}\left(\left[\mathfrak{z}\right]_{p^{n+1}}\right)}{n}-\underbrace{\frac{\#_{p:k}\left(\left[\mathfrak{z}\right]_{p}\right)}{n}}_{\rightarrow0}\right)\\
 & =\lim_{n\rightarrow\infty}\overbrace{\frac{n+1}{n}}^{\rightarrow1}\frac{\#_{p:k}\left(\left[\mathfrak{z}\right]_{p^{n+1}}\right)}{n+1}\\
 & =d_{p:k}\left(\mathfrak{z}\right)
\end{align*}
Thus, $d_{p:k}\left(\theta_{p}\left(\mathfrak{z}\right)\right)$ exists
and is equal to $d_{p:k}\left(\mathfrak{z}\right)$ for all $k$,
which proves (I).

II. 
\begin{align*}
\lim_{n\rightarrow\infty}\frac{1}{n}\#_{p:k}\left(\left[\theta_{p}^{\circ\gamma\left(n\right)}\left(\mathfrak{z}\right)\right]_{p^{n-\gamma\left(n\right)}}\right) & =\lim_{n\rightarrow\infty}\frac{\#_{p:k}\left(\left[\mathfrak{z}\right]_{p^{n}}\right)-\#_{p:k}\left(\left[\mathfrak{z}\right]_{p^{\gamma\left(n\right)}}\right)}{n}\\
 & =d_{p:k}\left(\mathfrak{z}\right)-\lim_{n\rightarrow\infty}\frac{\gamma\left(n\right)}{n}\frac{\#_{p:k}\left(\left[\mathfrak{z}\right]_{p^{\gamma\left(n\right)}}\right)}{\gamma\left(n\right)}\\
 & =d_{p:k}\left(\mathfrak{z}\right)-Ld_{p:k}\left(\mathfrak{z}\right)\\
 & =\left(1-L\right)d_{p:k}\left(\mathfrak{z}\right)
\end{align*}

Q.E.D.

\vphantom{}

To make things easier to read, we need some notation, as well as a
term for the qualitative control we will seek to establish on the
decay of the tails of the Fourier series of the $X_{J}$s.
\begin{notation}
Let $\left\{ M_{n}\right\} _{n\geq0}\in\textrm{MFunc}\left(\mathbb{Z}_{p},K\right)$,
let $\alpha\in\mathscr{V}\left(K\right)$ be an absolute value on
$K$. Then, we write $\alpha\left(M_{n}\right)$ denote the function
$\mathfrak{z}\mapsto\alpha\left(M_{n}\left(\mathfrak{z}\right)\right)$.
Note that since absolute values are multiplicative, we have that $\left\{ \alpha\left(M_{n}\right)\right\} _{n\geq0}$
is then a real-valued $p$-adic M-function.
\end{notation}
Proof: Absolute values are multiplicative and positive-definite.

Q.E.D.

\vphantom{}

\begin{defn}
\label{def:nice decay}Let $\left\{ E_{n}\right\} _{n\geq0}$ be a
real-valued $p$-adic M-function. Given $\mathfrak{z}\in\mathbb{Z}_{p}$,
we say $\left\{ E_{n}\right\} _{n\geq0}$ \textbf{decays summably
at $\mathfrak{z}$ / has summable decay at $\mathfrak{z}$ }whenever
$\sum_{n=0}^{\infty}E_{n}\left(\mathfrak{z}\right)$ converges in
$\mathbb{R}$.

More generally, given an absolute value $\alpha\in\mathscr{V}\left(K\right)$
and a $p$-adic $K$-valued M-function $\left\{ M_{n}\right\} _{n\geq0}$,
we say $\left\{ M_{n}\right\} _{n\geq0}$ \textbf{decays $\alpha$-summably
at $\mathfrak{z}$ / has $\alpha$-summable decay at $\mathfrak{z}$
}whenever the real-valued M-function $\left\{ \alpha\left(M_{n}\right)\right\} _{n\geq0}$
decays summably $\mathfrak{z}$.
\end{defn}
\begin{prop}
\label{prop:root test condition}Let $\left\{ E_{n}\right\} _{n\geq0}$
be a real-valued $p$-adic M-function with multipliers $e_{0},\ldots,e_{p-1}$.
Then $E_{n}$ decays summably at $\mathfrak{z}$ whenever:
\begin{equation}
\ln e_{0}+\sum_{k:e_{k}>e_{0}}\overline{d_{p:k}}\left(\mathfrak{z}\right)\left|\ln\frac{e_{k}}{e_{0}}\right|<\sum_{k:e_{k}\leq e_{0}}\underline{d_{p:k}}\left(\mathfrak{z}\right)\left|\ln\frac{e_{k}}{e_{0}}\right|\label{eq:decay root test condition}
\end{equation}
\end{prop}
Proof: Since:
\begin{equation}
\limsup_{n\rightarrow\infty}E_{n}^{1/n}\left(\mathfrak{z}\right)=e_{0}\prod_{k=1}^{p-1}\left(\frac{e_{k}}{e_{0}}\right)^{\underline{d_{p:k}}\left(\mathfrak{z}\right)\left[e_{k}\leq e_{0}\right]+\overline{d_{p:k}}\left(\mathfrak{z}\right)\left[e_{k}>e_{0}\right]}
\end{equation}
the \textbf{Root Test} for series convergence guarantees that:
\begin{equation}
\sum_{n=0}^{\infty}E_{n}\left(\mathfrak{z}\right)<\infty
\end{equation}
will occur whenever:
\begin{equation}
e_{0}\prod_{k=1}^{p-1}\left(\frac{e_{k}}{e_{0}}\right)^{\underline{d_{p:k}}\left(\mathfrak{z}\right)\left[e_{k}\leq e_{0}\right]+\overline{d_{p:k}}\left(\mathfrak{z}\right)\left[e_{k}>e_{0}\right]}<1
\end{equation}
i.e.:
\begin{equation}
\ln e_{0}+\sum_{k:e_{k}>e_{0}}\overline{d_{p:k}}\left(\mathfrak{z}\right)\left|\ln\frac{e_{k}}{e_{0}}\right|<\sum_{k:e_{k}\leq e_{0}}\underline{d_{p:k}}\left(\mathfrak{z}\right)\left|\ln\frac{e_{k}}{e_{0}}\right|
\end{equation}

Q.E.D.
\vphantom{}Next, we have a notation to simplify the statements of
our next few propositions.
\begin{defn}
Let $\left\{ a_{n}\right\} _{n\geq0}$ and $\left\{ b_{n}\right\} _{n\geq0}$
be sequences in $\mathbb{N}_{0}$, and let $\mathfrak{z}\in\mathbb{Z}_{p}$.
Given real-valued $p$-adic M-functions $\left\{ E_{n}\right\} _{n\geq0}$
and $\left\{ F_{n}\right\} _{n\geq0}$, we write: 
\begin{equation}
E_{n}\preceq_{a_{n},b_{n},\mathfrak{z}}F_{n}\label{eq:ordering def}
\end{equation}
whenever the following conditions hold:

I. $\lim_{n\rightarrow\infty}E_{a_{n}}^{1/n}\left(\theta_{p}^{\circ b_{n}}\left(\mathfrak{z}\right)\right)$
converges in $\mathbb{R}$ to an element of $\left(0,1\right)$.

II. $\lim_{n\rightarrow\infty}F_{n}^{1/n}\left(\mathfrak{z}\right)$
converges in $\mathbb{R}$ to an element of $\left(0,1\right)$.

III. $E_{a_{n}}\left(\theta_{p}^{\circ b_{n}}\left(\mathfrak{z}\right)\right)\ll_{\mathfrak{z}}F_{n}\left(\mathfrak{z}\right)$
as $n\rightarrow\infty$; that is, there is a real constant $C>0$
so that $E_{a_{n}}\left(\theta_{p}^{\circ b_{n}}\left(\mathfrak{z}\right)\right)\leq CF_{n}\left(\mathfrak{z}\right)$
holds for all sufficiently large $n$.
\end{defn}
Our first proposition shows that for suitable $\gamma$, the decay
of $E_{\gamma\left(n\right)}\left(\mathfrak{z}\right)$ will be bounded
from above by the decay of an M-function $\left\{ F_{n}\right\} _{n\geq0}$.
\begin{prop}
\label{prop:gamma prop}Fix $\mathfrak{z}\in\textrm{NiceDig}\left(\mathbb{Z}_{p}\right)$,
and let $\left\{ E_{n}\right\} _{n\geq0}$ be a real-valued $p$-adic
M-function that has summable decay at $\mathfrak{z}$. Next, let $\gamma:\mathbb{N}_{0}\rightarrow\mathbb{N}_{0}$
be any function satisfying:
\begin{align}
\gamma\left(n\right) & \leq n,\textrm{ }\forall n\geq0\\
\lim_{n\rightarrow\infty}\gamma\left(n\right) & \overset{\mathbb{R}}{=}+\infty\\
0<\lim_{n\rightarrow\infty}\frac{\gamma\left(n\right)}{n} & \leq1
\end{align}
where the limit $\lim_{n\rightarrow\infty}\gamma\left(n\right)/n$
exists in $\mathbb{R}$.

If $\lim_{n\rightarrow\infty}E_{n}^{1/n}\left(\mathfrak{z}\right)<1$,
then there is a real-valued $p$-adic M-function $\left\{ F_{n}\right\} _{n\geq0}$
so that $E_{n}\preceq_{\gamma\left(n\right),0,\mathfrak{z}}F_{n}$.
\end{prop}
Proof: Since $\mathfrak{z}$ has nice $p$-adic digits:
\begin{equation}
d_{p:k}\left(\mathfrak{z}\right)=\lim_{n\rightarrow\infty}\frac{\#_{p:k}\left(\left[\mathfrak{z}\right]_{p^{n}}\right)}{n}
\end{equation}
exists and is in $\left[0,1\right]$ for all $k\in\left\{ 0,\ldots,p-1\right\} $.
Hence:
\begin{equation}
\lim_{n\rightarrow\infty}E_{n}^{1/n}\left(\mathfrak{z}\right)=\lim_{n\rightarrow\infty}e_{0}\prod_{k=1}^{p-1}\left(e_{k}/e_{0}\right)^{\#_{p:k}\left(\left[\mathfrak{z}\right]_{p^{n}}\right)/n}\overset{\mathbb{R}}{=}e_{0}\prod_{k=1}^{p-1}\left(e_{k}/e_{0}\right)^{d_{p:k}\left(\mathfrak{z}\right)}
\end{equation}
and so, we have the asymptotic equivalence:
\begin{equation}
E_{n}\left(\mathfrak{z}\right)\sim e_{0}^{n}\prod_{k=1}^{p-1}\left(e_{k}/e_{0}\right)^{nd_{p:k}\left(\mathfrak{z}\right)}\textrm{ as }n\rightarrow\infty
\end{equation}
Letting $f_{0},\ldots,f_{p-1}$ be undetermined constants in $\left(0,1\right)$,
set:
\begin{equation}
F_{n}\left(\mathfrak{z}\right)=\prod_{m=0}^{n-1}f_{\left[\theta_{p}^{n}\left(\mathfrak{z}\right)\right]_{p}}
\end{equation}
Then:
\begin{equation}
\lim_{n\rightarrow\infty}F_{n}^{1/n}\left(\mathfrak{z}\right)\overset{\mathbb{R}}{=}f_{0}\prod_{k=1}^{p-1}\left(f_{k}/f_{0}\right)^{d_{p:k}\left(\mathfrak{z}\right)}\textrm{ as }n\rightarrow\infty
\end{equation}
As such, by \textbf{Lemma \ref{lem:ll lemma}}, we have that $E_{\gamma\left(n\right)}\left(\mathfrak{z}\right)\ll F_{n}\left(\mathfrak{z}\right)\textrm{ as }n\rightarrow\infty$
will occur provided:
\begin{equation}
e_{0}^{\gamma\left(n\right)}\prod_{k=1}^{p-1}\left(e_{k}/e_{0}\right)^{\gamma\left(n\right)d_{p:k}\left(\mathfrak{z}\right)}\ll_{\mathfrak{z}}f_{0}^{n}\prod_{k=1}^{p-1}\left(f_{k}/f_{0}\right)^{nd_{p:k}\left(\mathfrak{z}\right)}\textrm{ as }n\rightarrow\infty\label{eq:proviso-1}
\end{equation}
Here:

\begin{align}
\lim_{n\rightarrow\infty}\left(e_{0}^{\gamma\left(n\right)}\prod_{k=1}^{p-1}\left(e_{k}/e_{0}\right)^{\gamma\left(n\right)d_{p:k}\left(\mathfrak{z}\right)}\right)^{1/n} & \overset{\mathbb{R}}{=}\left(e_{0}\prod_{k=1}^{p-1}\left(e_{k}/e_{0}\right)^{d_{p:k}\left(\mathfrak{z}\right)}\right)^{L}\\
\lim_{n\rightarrow\infty}\left(f_{0}^{n}\prod_{k=1}^{p-1}\left(f_{k}/f_{0}\right)^{nd_{p:k}\left(\mathfrak{z}\right)}\right)^{1/n} & \overset{\mathbb{R}}{=}f_{0}\prod_{k=1}^{p-1}\left(f_{k}/f_{0}\right)^{d_{p:k}\left(\mathfrak{z}\right)}
\end{align}
Thus, (\ref{eq:proviso-1}), will occur provided:
\begin{equation}
\left(e_{0}\prod_{k=1}^{p-1}\left(e_{k}/e_{0}\right)^{d_{p:k}\left(\mathfrak{z}\right)}\right)^{Ln}\ll_{\mathfrak{z}}\left(f_{0}\prod_{k=1}^{p-1}\left(f_{k}/f_{0}\right)^{d_{p:k}\left(\mathfrak{z}\right)}\right)^{n}\textrm{ as }n\rightarrow\infty
\end{equation}
This, in turn, occurs provided:
\begin{equation}
\left(\underbrace{e_{0}\prod_{k=1}^{p-1}\left(e_{k}/e_{0}\right)^{d_{p:k}\left(\mathfrak{z}\right)}}_{=\lim_{n\rightarrow\infty}E_{n}^{1/n}\left(\mathfrak{z}\right)<1}\right)^{L}<f_{0}\prod_{k=1}^{p-1}\left(f_{k}/f_{0}\right)^{d_{p:k}\left(\mathfrak{z}\right)}\label{eq:proviso-2}
\end{equation}
Taking logarithms, and choosing the $f_{k}$s to be real numbers satisfying:
\begin{equation}
0<f_{0}\prod_{k=1}^{p-1}\left(f_{k}/f_{0}\right)^{d_{p:k}\left(\mathfrak{z}\right)}<1
\end{equation}
(\ref{eq:proviso-2}) then becomes:
\begin{equation}
-L\left|\ln\left(e_{0}\prod_{k=1}^{p-1}\left(e_{k}/e_{0}\right)^{d_{p:k}\left(\mathfrak{z}\right)}\right)\right|<-\left|\ln\left(f_{0}\prod_{k=1}^{p-1}\left(f_{k}/f_{0}\right)^{d_{p:k}\left(\mathfrak{z}\right)}\right)\right|
\end{equation}
and hence:
\begin{equation}
L>\left|\frac{\ln\left(f_{0}\prod_{k=1}^{p-1}\left(f_{k}/f_{0}\right)^{d_{p:k}\left(\mathfrak{z}\right)}\right)}{\ln\left(e_{0}\prod_{k=1}^{p-1}\left(e_{k}/e_{0}\right)^{d_{p:k}\left(\mathfrak{z}\right)}\right)}\right|
\end{equation}
Since $0<L\leq1$, we can always make this work by choosing the $f_{k}$s
so that $f_{0}\prod_{k=1}^{p-1}\left(f_{k}/f_{0}\right)^{d_{p:k}\left(\mathfrak{z}\right)}$
is a number in $\left(0,1\right)$ which is only slightly larger than
$e_{0}\prod_{k=1}^{p-1}\left(e_{k}/e_{0}\right)^{d_{p:k}\left(\mathfrak{z}\right)}$.

Q.E.D.

\vphantom{}Next, we deal with getting an M-function upper bound for
$E_{n}\left(\theta_{p}^{\circ n-\gamma\left(n\right)}\left(\mathfrak{z}\right)\right)$.
\begin{prop}
\label{prop:gamma theta prop}Fix $\mathfrak{z}\in\textrm{NiceDig}\left(\mathbb{Z}_{p}\right)$,
and let $\left\{ E_{n}\right\} _{n\geq0}$ be a real-valued $p$-adic
M-function that has summable decay at $\mathfrak{z}$. Next, let $\gamma:\mathbb{N}_{0}\rightarrow\mathbb{N}_{0}$
be any function satisfying:

\begin{align}
\gamma\left(n\right) & \leq n,\textrm{ }\forall n\geq0\\
\lim_{n\rightarrow\infty}\gamma\left(n\right) & \overset{\mathbb{R}}{=}+\infty\\
0\leq\lim_{n\rightarrow\infty}\frac{\gamma\left(n\right)}{n} & <1
\end{align}
where the limit $L\overset{\textrm{def}}{=}\lim_{n\rightarrow\infty}\gamma\left(n\right)/n$
exists in $\mathbb{R}$.

If $\lim_{n\rightarrow\infty}E_{n}^{1/n}\left(\mathfrak{z}\right)<1$,
then there is a real-valued $p$-adic M-function $\left\{ F_{n}\right\} _{n\geq0}$
so that $E_{n}\preceq_{n-\gamma\left(n\right),n,\mathfrak{z}}F_{n}$.
\end{prop}
Proof: Since $\mathfrak{z}$ has nice digits, and since $\gamma$
satisfies the hypothesis of \textbf{Proposition \ref{prop:nice digits}},
we have:
\begin{equation}
\lim_{n\rightarrow\infty}\frac{1}{n}\#_{p:k}\left(\left[\theta_{p}^{\circ\gamma\left(n\right)}\left(\mathfrak{z}\right)\right]_{p^{n-\gamma\left(n\right)}}\right)=\left(1-L\right)d_{p:k}\left(\mathfrak{z}\right)
\end{equation}
Thus:
\begin{align*}
E_{n-\gamma\left(n\right)}^{1/n}\left(\theta_{p}^{\circ\gamma\left(n\right)}\left(\mathfrak{z}\right)\right) & =\lim_{n\rightarrow\infty}\left(e_{0}^{n-\gamma\left(n\right)}\prod_{k=1}^{p-1}\left(e_{k}/e_{0}\right)^{\#_{p:k}\left(\left[\theta_{p}^{\circ\gamma\left(n\right)}\left(\mathfrak{z}\right)\right]_{p^{n-\gamma\left(n\right)}}\right)}\right)^{1/n}\\
 & =\left(e_{0}\prod_{k=1}^{p-1}\left(e_{k}/e_{0}\right)^{d_{p:k}\left(\mathfrak{z}\right)}\right)^{1-L}
\end{align*}
So, letting $f_{0},\ldots,f_{p-1}$ be undetermined constants in $\left(0,1\right)$,
set:
\begin{equation}
F_{n}\left(\mathfrak{z}\right)=\prod_{m=0}^{n-1}f_{\left[\theta_{p}^{n}\left(\mathfrak{z}\right)\right]_{p}}
\end{equation}
Then:
\begin{equation}
F_{n}\left(\mathfrak{z}\right)\sim f_{0}^{n}\prod_{k=1}^{p-1}\left(f_{k}/f_{0}\right)^{nd_{p:k}\left(\mathfrak{z}\right)}\textrm{ as }n\rightarrow\infty
\end{equation}
As such, $E_{n-\gamma\left(n\right)}\left(\theta_{p}^{\circ\gamma\left(n\right)}\left(\mathfrak{z}\right)\right)\ll_{\mathfrak{z}}F_{n}\left(\mathfrak{z}\right)\textrm{ as }n\rightarrow\infty$
will occur provided:
\begin{equation}
\left(\underbrace{e_{0}\prod_{k=1}^{p-1}\left(e_{k}/e_{0}\right)^{d_{p:k}\left(\mathfrak{z}\right)}}_{=\lim_{n\rightarrow\infty}E_{n}^{1/n}\left(\mathfrak{z}\right)<1}\right)^{1-L}<f_{0}\prod_{k=1}^{p-1}\left(f_{k}/f_{0}\right)^{d_{p:k}\left(\mathfrak{z}\right)}\label{eq:proviso}
\end{equation}

Taking logarithms, and letting the $f_{k}$s be such so that:
\begin{equation}
0<f_{0}\prod_{k=1}^{p-1}\left(f_{k}/f_{0}\right)^{d_{p:k}\left(\mathfrak{z}\right)}<1
\end{equation}
we have:
\begin{equation}
-\left(1-L\right)\left|\ln\left(e_{0}\prod_{k=1}^{p-1}\left(e_{k}/e_{0}\right)^{d_{p:k}\left(\mathfrak{z}\right)}\right)\right|<-\left|\ln\left(f_{0}\prod_{k=1}^{p-1}\left(f_{k}/f_{0}\right)^{d_{p:k}\left(\mathfrak{z}\right)}\right)\right|
\end{equation}
hence:
\begin{equation}
L<1-\left|\frac{\ln\left(f_{0}\prod_{k=1}^{p-1}\left(f_{k}/f_{0}\right)^{d_{p:k}\left(\mathfrak{z}\right)}\right)}{\ln\left(e_{0}\prod_{k=1}^{p-1}\left(e_{k}/e_{0}\right)^{d_{p:k}\left(\mathfrak{z}\right)}\right)}\right|
\end{equation}
Since $0\leq L<1$, we can always make this work by choosing the $f_{k}$s
so that $f_{0}\prod_{k=1}^{p-1}\left(f_{k}/f_{0}\right)^{d_{p:k}\left(\mathfrak{z}\right)}$
is a number in $\left(0,1\right)$ which is only slightly larger than
$e_{0}\prod_{k=1}^{p-1}\left(e_{k}/e_{0}\right)^{d_{p:k}\left(\mathfrak{z}\right)}$.

Q.E.D.

\vphantom{}Next, we deal with the sum of two M-functions:
\begin{prop}
\label{prop:sum of M-functions bound}Fix $\mathfrak{z}\in\textrm{NiceDig}\left(\mathbb{Z}_{p}\right)$,
and let $\left\{ E_{n}\right\} _{n\geq0}$ and $\left\{ F_{n}\right\} _{n\geq0}$
be real-valued $p$-adic M-functions. If: 
\begin{equation}
\max\left\{ \lim_{n\rightarrow\infty}E_{n}^{1/n}\left(\mathfrak{z}\right),\lim_{n\rightarrow\infty}F_{n}^{1/n}\left(\mathfrak{z}\right)\right\} <1
\end{equation}
Then, there is a real-valued $p$-adic M-function $\left\{ G_{n}\right\} _{n\geq0}$
so that:

I. $E_{n}\left(\mathfrak{z}\right)+F_{n}\left(\mathfrak{z}\right)\ll_{\mathfrak{z}}G_{n}\left(\mathfrak{z}\right)\textrm{ as }n\rightarrow\infty$.

II. $\lim_{n\rightarrow\infty}G_{n}^{1/n}\left(\mathfrak{z}\right)$$<1$.
\end{prop}
Proof: Since $\mathfrak{z}$ has nice digits, we can write:
\begin{align}
E_{n}\left(\mathfrak{z}\right) & \sim e_{0}^{n}\prod_{k=1}^{p-1}\left(e_{k}/e_{0}\right)^{nd_{p:k}\left(\mathfrak{z}\right)}\textrm{ as }n\rightarrow\infty\\
F_{n}\left(\mathfrak{z}\right) & \sim f_{0}^{n}\prod_{k=1}^{p-1}\left(f_{k}/f_{0}\right)^{nd_{p:k}\left(\mathfrak{z}\right)}\textrm{ as }n\rightarrow\infty
\end{align}
Hence, letting:
\begin{equation}
G_{n}\left(\mathfrak{z}\right)=g_{0}^{n}\prod_{k=1}^{p-1}\left(\frac{g_{k}}{g_{0}}\right)^{\#_{p:k}\left(\left[\mathfrak{z}\right]_{p^{n}}\right)}
\end{equation}
we need only choose the $g_{k}$s to be small enough so that: 
\begin{equation}
g_{0}\prod_{k=1}^{p-1}\left(g_{k}/g_{0}\right)^{d_{p:k}\left(\mathfrak{z}\right)}<1
\end{equation}
yet large enough so that:
\begin{equation}
g_{0}\prod_{k=1}^{p-1}\left(\frac{g_{k}}{g_{0}}\right)^{d_{p:k}\left(\mathfrak{z}\right)}\geq\max\left\{ e_{0}\prod_{k=1}^{p-1}\left(\frac{e_{k}}{e_{0}}\right)^{d_{p:k}\left(\mathfrak{z}\right)},f_{0}\prod_{k=1}^{p-1}\left(\frac{f_{k}}{f_{0}}\right)^{d_{p:k}\left(\mathfrak{z}\right)}\right\} 
\end{equation}
which can always be done, thanks to the hypothesis on $E_{n}$ and
$F_{n}$.

Q.E.D.

\vphantom{}Almost done, we now deal with the maximum of two M-functions:
\begin{prop}
\label{prop:max of two M functions}Fix $\mathfrak{z}\in\textrm{NiceDig}\left(\mathbb{Z}_{p}\right)$,
and let $\left\{ E_{n}\right\} _{n\geq0}$ and $\left\{ F_{n}\right\} _{n\geq0}$
be real-valued $p$-adic M-functions. If:
\begin{equation}
\max\left\{ \lim_{n\rightarrow\infty}E_{n}^{1/n}\left(\mathfrak{z}\right),\lim_{n\rightarrow\infty}F_{n}^{1/n}\left(\mathfrak{z}\right)\right\} <1
\end{equation}
Then, there is a real-valued $p$-adic M-function $\left\{ G_{n}\right\} _{n\geq0}$
so that:

I. $\max\left\{ E_{n}\left(\mathfrak{z}\right),F_{n}\left(\mathfrak{z}\right)\right\} \ll_{\mathfrak{z}}G_{n}\left(\mathfrak{z}\right)\textrm{ as }n\rightarrow\infty$.

II. $\lim_{n\rightarrow\infty}G_{n}^{1/n}\left(\mathfrak{z}\right)<1$.
\end{prop}
Proof: Observe that for all $\mathfrak{z}$ and all $n$, we have
the bound:
\begin{equation}
\max\left\{ E_{n}\left(\mathfrak{z}\right),F_{n}\left(\mathfrak{z}\right)\right\} \leq\left(\max\left\{ e_{0},f_{0}\right\} \right)^{n}\prod_{k=1}^{p-1}\left(\max\left\{ \frac{e_{k}}{e_{0}},\frac{f_{k}}{f_{0}}\right\} \right)^{\#_{p:k}\left(\left[\mathfrak{z}\right]_{p^{n}}\right)}
\end{equation}
The upper bound is a real-valued M-function, by \textbf{Proposition
\ref{prop:M function characterization}}. So, we choose it to be our
$G_{n}$. The conditions on the multipliers of $E_{n}$ and $F_{n}$
then pass to the multipliers of $G_{n}$.

\vphantom{}

The last ingredient we need is a technical lemma that will be used
in the next subsection in conjunction with \textbf{Proposition \ref{prop:Kappa shift equation}}.
\begin{prop}
\label{prop:quotient prop}Let $a:\mathbb{N}_{0}\rightarrow\left(0,\infty\right)$
be a sequence of positive real numbers with: 
\begin{equation}
\lim_{n\rightarrow\infty}a\left(n\right)\overset{\mathbb{R}}{=}0
\end{equation}
Then, there exists a non-decreasing $\gamma:\mathbb{N}_{0}\rightarrow\mathbb{N}_{0}$
with: 
\begin{align}
\gamma\left(n\right) & \leq n,\textrm{ }\forall n\geq0\\
\lim_{n\rightarrow\infty}\gamma\left(n\right) & \overset{\mathbb{R}}{=}+\infty
\end{align}
so that:
\begin{equation}
\lim_{n\rightarrow\infty}\frac{a\left(n\right)}{a\left(\gamma\left(n\right)\right)}\overset{\mathbb{R}}{=}0
\end{equation}
\end{prop}
Proof: Set:
\begin{equation}
c_{n}\overset{\textrm{def}}{=}\max\left\{ a\left(m\right):m\geq m\right\} 
\end{equation}
Note that $c_{n}$ is a non-increasing sequence tending to $0$. Next,
let $I\overset{\textrm{def}}{=}\left\{ c_{n}:n\geq0\right\} $, and
enumerate the elements of $I$ in decreasing order as $i_{1}>i_{2}>\ldots$.
Also, let $i_{0}$ denote $+\infty$. Then, for each $n$, let $x_{n}$
be the largest integer $\geq n$ for which $a\left(x_{n}\right)=i_{n}$,
and set $x_{0}=0$. By construction, we have that $x_{1}<x_{2}<\cdots$.
Lastly, for any $\epsilon>0$, let $N_{\epsilon}$ be the smallest
non-negative integer so that $a\left(n\right)<\epsilon$ for all $n\geq N_{\epsilon}$.

We then define $\gamma:\mathbb{N}_{0}\rightarrow\mathbb{N}_{0}$ like
so: 
\begin{equation}
\gamma\left(n\right)\overset{\textrm{def}}{=}\sum_{k=0}^{\infty}x_{k}\left[N_{i_{k}/2^{k}}\leq n<N_{i_{k+1}/2^{k+1}}\right]
\end{equation}
where $N_{i_{0}}=0$, because $a\left(n\right)<i_{0}=\infty$ occurs
for all $n\geq0$.

So, for example, we have $\gamma\left(0\right)=\gamma\left(1\right)=\cdots=\gamma\left(n\right)=0$
until we reach $n=N_{i_{1}/2}$, the smallest $n$ past which $a_{n}$
will always be less than $i_{1}/2$. Then, we have $\gamma\left(N_{i_{1}/2}\right)=\gamma\left(N_{i_{1}/2}+1\right)=\cdots=\gamma\left(n\right)=x_{1}$
until we reach $n=N_{i_{2}/4}$, the smallest $n$ past which $a\left(n\right)$
will always be less than $i_{2}/4$; then $\gamma\left(N_{i_{2}/4}\right)=\gamma\left(N_{i_{2}/4}+1\right)=\cdots=\gamma\left(n\right)=x_{2}$,
until... (and so on and so forth).

As constructed, we have that:
\begin{equation}
N_{i_{k}/2^{k}}\leq n<N_{i_{k+1}/2^{k+1}}\Rightarrow\begin{cases}
a\left(n\right)\leq\frac{i_{k}}{2^{k}}\\
\gamma\left(n\right)=x_{k}
\end{cases},\textrm{ }\forall k\geq0
\end{equation}
and so:
\begin{equation}
N_{i_{k}/2^{k}}\leq n<N_{i_{k+1}/2^{k+1}}\Rightarrow a\left(\gamma\left(n\right)\right)=a\left(x_{k}\right)=i_{k}
\end{equation}
Taking limsups, we get:
\begin{align*}
\limsup_{n\rightarrow\infty}\frac{a\left(n\right)}{a\left(\gamma\left(n\right)\right)} & \leq\limsup_{k\rightarrow\infty}\max_{N_{i_{k}/2^{k}}\leq n<N_{i_{k+1}/2^{k+1}}}\frac{a\left(n\right)}{a\left(\gamma\left(n\right)\right)}\\
 & \leq\limsup_{k\rightarrow\infty}\frac{i_{k}/2^{k}}{i_{k}}\\
 & \overset{\mathbb{R}}{=}0
\end{align*}
which proves that:
\begin{equation}
\lim_{n\rightarrow\infty}\frac{a\left(n\right)}{a\left(\gamma\left(n\right)\right)}\overset{\mathbb{R}}{=}0
\end{equation}

Q.E.D.

\subsection{\label{subsec:Frame-Theoretic-Asymptotics}Frame-Theoretic Asymptotics
and the Main LimitLemma}

In this section, we will prove a key lemma (\textbf{Lemma \ref{lem:main limit lemma}}
on page \pageref{lem:main limit lemma}) needed to justify limits
of the type taken in \textbf{Assumption \pageref{assu:main limit lemma}}.
In order to prove this result, we need to modify standard asymptotic
notation to be compatible with frames.
\begin{defn}
Let $X$ be a set, let $R$ be a Dedekind domain, and let $\mathcal{F}$
be an $R$-frame on $X$.

I. Given $f,g\in C\left(\mathcal{F}\right)$, we write:
\begin{equation}
f\left(x\right)\overset{\mathcal{F}}{\ll}g\left(x\right)
\end{equation}
to mean that for each $x\in D\left(\mathcal{F}\right)$, there is
a positive real constant $C_{x}$ so that:
\begin{equation}
\mathcal{F}_{x}\left(f\right)\leq C_{x}\mathcal{F}_{x}\left(g\right)
\end{equation}

II. Given sequences $\left\{ f_{N}\right\} _{N\geq0}$ and $\left\{ g_{N}\right\} _{N\geq0}$
in $C\left(\mathcal{F}\right)$. We write:

\begin{equation}
f_{N}\left(x\right)\overset{\mathcal{F}}{\ll}g_{N}\left(x\right)\textrm{ as }N\rightarrow\infty
\end{equation}
and also:
\begin{equation}
f_{N}\left(x\right)\overset{\mathcal{F}}{=}O\left(g_{N}\left(x\right)\right)\textrm{ as }N\rightarrow\infty
\end{equation}
to mean that, for each $x\in D\left(\mathcal{F}\right)$, there is
a positive real constant $C_{x}$ and an integer $N_{x}\geq0$ so
that:
\begin{equation}
\mathcal{F}_{x}\left(f_{N}\right)\leq C_{x}\mathcal{F}_{x}\left(g_{N}\right),\textrm{ }\forall N\geq N_{x}
\end{equation}

III. Given $d\geq1$ and $f_{1},\ldots,f_{d},g_{1},\ldots,g_{d}\in C\left(\mathcal{F}\right)$,
we write:
\begin{equation}
\max\left\{ f_{1}\left(x\right),\ldots,f_{d}\left(x\right)\right\} \overset{\mathcal{F}}{\ll}\max\left\{ g_{1}\left(x\right),\ldots,g_{d}\left(x\right)\right\} 
\end{equation}
whenever:
\begin{equation}
\max\left\{ \mathcal{F}_{x}\left(f_{1}\right),\ldots,\mathcal{F}_{x}\left(f_{d}\right)\right\} \ll\max\left\{ \mathcal{F}_{x}\left(g_{1}\right),\ldots,\mathcal{F}_{x}\left(g_{d}\right)\right\} 
\end{equation}
and likewise for minimums.

Also, given sequences $\left\{ f_{1,n}\right\} _{n\geq0},\ldots,\left\{ f_{d,n}\right\} _{n\geq0}\in C\left(\mathcal{F}\right)$,
we write:
\begin{equation}
\lim_{n\rightarrow\infty}\max\left\{ f_{1,n}\left(x\right),\ldots,f_{d,n}\left(x\right)\right\} \overset{\mathcal{F}}{=}0
\end{equation}
to mean that:
\begin{equation}
\lim_{n\rightarrow\infty}\max\left\{ \mathcal{F}_{x}\left(f_{1,n}\right),\ldots,\mathcal{F}_{x}\left(f_{d,n}\right)\right\} \overset{\mathbb{R}}{=}0,\textrm{ }\forall x\in D\left(\mathcal{F}\right)
\end{equation}
and likewise for minimums.

IV. Given an M-function $\left\{ M_{n}\right\} _{n\geq0}\in\textrm{MFunc}\left(\mathbb{Z}_{p},R\right)$,
we say that $\left\{ M_{n}\right\} _{n\geq0}$ is \textbf{strongly
$\mathcal{F}$-summable }whenever:
\begin{equation}
\sum_{n=0}^{\infty}M_{n}\left(\mathfrak{z}\right)
\end{equation}
is $\mathcal{F}$-convergent, with:
\begin{equation}
\sum_{n=0}^{\infty}\mathcal{F}_{\mathfrak{z}}\left(M_{n}\right)<\infty
\end{equation}
for all $\mathfrak{z}\in D\left(\mathcal{F}\right)$ for which $\mathcal{F}_{\mathfrak{z}}$
is archimedean.

V. Let $\left\{ g_{n}\right\} _{n\geq0}$ be a sequence in $C\left(\mathcal{F}\right)$.
Then, I say the \textbf{$g_{n}$s have $\left(\mathcal{F},M\right)$-type
decay }whenever, for each $\mathfrak{z}\in D\left(\mathcal{F}\right)$,
there is an M-function $\left\{ E_{\mathfrak{z},n}\right\} _{n\geq0}$
(depending on $\mathfrak{z}$!) so that:
\begin{equation}
\lim_{n\rightarrow\infty}\mathcal{F}_{\mathfrak{z}}\left(E_{\mathfrak{z},n}\right)\overset{\mathbb{R}}{=}0,\textrm{ }\forall\mathfrak{z}\in D\left(\mathcal{F}\right)
\end{equation}
and, for all $m,n\geq0$:
\begin{equation}
g_{n}\circ\theta_{p}^{\circ m}\overset{\mathcal{F}}{\ll}E_{\mathfrak{z},n}\circ\theta_{p}^{\circ m}\label{eq:f kappa type definition}
\end{equation}
where the constant of proportionality depends only on $\mathfrak{z}$.

More generally, I say a sequence $\left\{ f_{n}\right\} _{n\geq0}$
in $C\left(\mathcal{F}\right)$ are\textbf{ $\left(\mathcal{F},M\right)$-convergent
}whenever:

i. The $f_{n}$s $\mathcal{F}$-converge to some $f\in C\left(\mathcal{F}\right)$;

ii. The function $g_{n}=f-f_{n}$ has $\left(\mathcal{F},M\right)$-type
decay.

and I then say that the $f_{n}$s \textbf{$\left(\mathcal{F},M\right)$-converge
}to $f$.
\end{defn}
\begin{rem}
For a given $\mathfrak{z}\in D\left(\mathcal{F}\right)$, if $\mathcal{F}_{\mathfrak{z}}\left(M_{n}\right)$
has summable decay at $\mathfrak{z}$, then $\sum_{n=0}^{\infty}M_{n}\left(\mathfrak{z}\right)$
converges in $\mathcal{F}\left(\mathfrak{z}\right)$. However, the
converse is not always true: the convergence of $\sum_{n=0}^{\infty}M_{n}\left(\mathfrak{z}\right)$
in $\mathcal{F}\left(\mathfrak{z}\right)$ need not imply that $\mathcal{F}_{\mathfrak{z}}\left(M_{n}\right)$
has summable decay at $\mathfrak{z}$. We use the slightly more general
condition of strong $\mathcal{F}$-summability because it increases
the range of applicability of \textbf{Lemma \ref{lem:main limit lemma}},
which is of independent interest in its own right.
\end{rem}
\begin{rem}
(\ref{eq:f kappa type definition}) means that for each $\mathfrak{z}\in D\left(\mathcal{F}\right)$,
there is a real constant $C_{\mathfrak{z}}>0$ so that:
\begin{equation}
\mathcal{F}_{\mathfrak{z}}\left(g_{n}\circ\theta_{p}^{\circ m}\right)\leq C_{\mathfrak{z}}\mathcal{F}_{\mathfrak{z}}\left(E_{\mathfrak{z},n}\circ\theta_{p}^{\circ m}\right),\textrm{ }\forall m,n\in\mathbb{N}_{0}\label{eq:f kappa type elaborated}
\end{equation}
When $\mathcal{F}$ is evaluative and $E_{\mathfrak{z},n}$ is of
the form:
\begin{equation}
E_{\mathfrak{z},n}\left(\mathfrak{y}\right)=s_{\mathfrak{z},0}^{n}\prod_{k=1}^{p-1}\left(\frac{s_{\mathfrak{z},k}}{s_{\mathfrak{z},0}}\right)^{\#_{p:k}\left(\left[\mathfrak{y}\right]_{p^{n}}\right)}
\end{equation}
for constants $s_{\mathfrak{z},0},\ldots,s_{\mathfrak{z},p-1}\in R$,
(\ref{eq:f kappa type elaborated}) reduces to:
\begin{equation}
\left|g_{n}\left(\theta_{p}^{\circ m}\left(\mathfrak{z}\right)\right)\right|_{\mathcal{F}\left(\mathfrak{z}\right)}\leq C_{\mathfrak{z}}\left|s_{\mathfrak{z},0}^{n}\prod_{k=1}^{p-1}\left(\frac{s_{\mathfrak{z},k}}{s_{\mathfrak{z},0}}\right)^{\#_{p:k}\left(\left[\theta_{p}^{\circ m}\left(\mathfrak{z}\right)\right]_{p^{n}}\right)}\right|_{\mathcal{F}\left(\mathfrak{z}\right)},\textrm{ }\forall m,n\in\mathbb{N}_{0},\textrm{ }\forall\mathfrak{z}\in D\left(\mathcal{F}\right)
\end{equation}
The decay requirement on the $E_{\mathfrak{z},n}$s is then:
\begin{equation}
\lim_{n\rightarrow\infty}\left|s_{\mathfrak{z},0}^{n}\prod_{k=1}^{p-1}\left(\frac{s_{\mathfrak{z},k}}{s_{\mathfrak{z},0}}\right)^{\#_{p:k}\left(\left[\mathfrak{z}\right]_{p^{n}}\right)}\right|_{\mathcal{F}\left(\mathfrak{z}\right)}\overset{\mathbb{R}}{=}0,\textrm{ }\forall\mathfrak{z}\in D\left(\mathcal{F}\right)
\end{equation}
\end{rem}
\begin{rem}
To get a bit ahead of ourselves, if only for the sake of keeping our
eye on the prize, our eventual goal will be to show through induction/recursion
that, given $\mathbf{m}$, if $X_{\mathbf{k}}$'s Fourier series is
$\left(\mathcal{F},M\right)$-convergent for all $\mathbf{k}<\mathbf{m}$,
then $X_{\mathbf{m}}$'s Fourier series will be $\left(\mathcal{F},M\right)$-convergent,
as well. Not only will this then prove the quasi-integrability of
$X_{\mathbf{m}}$ (by showing that $X_{\mathbf{m}}$'s Fourier series
does, in fact, converge $\mathcal{F}$-wise to $X_{\mathbf{m}}$),
it will also then feed in to the next step of the inductive argument,
thanks to \textbf{Lemma \ref{lem:main limit lemma}}.
\end{rem}
At long last, we can finally begin to give justifications for the
work done in Section \textbf{\ref{subsec:The-Central-Computation}}.
Recall that, in that section, our principal assumption (\textbf{Assumption
\ref{assu:main limit lemma}}) was that the limit:

\begin{equation}
\lim_{N\rightarrow\infty}\sum_{n=0}^{N-2}r_{\mathbf{n},0}^{n}\kappa_{\mathbf{n}}\left(\left[\mathfrak{z}\right]_{p^{n}}\right)r_{\mathbf{m},\mathbf{n},\left[\theta_{p}^{\circ n}\left(\mathfrak{z}\right)\right]_{p}}\left(\tilde{X}_{\mathbf{m},N-n-1}\left(\theta_{p}^{\circ n+1}\left(\mathfrak{z}\right)\right)-\hat{X}_{\mathbf{m}}\left(0\right)\right)\label{eq:limit}
\end{equation}
converged to:
\begin{equation}
\sum_{n=0}^{\infty}r_{\mathbf{n},0}^{n}\kappa_{\mathbf{n}}\left(\left[\mathfrak{z}\right]_{p^{n}}\right)r_{\mathbf{m},\mathbf{n},\left[\theta_{p}^{\circ n}\left(\mathfrak{z}\right)\right]_{p}}\left(X_{\mathbf{m}}\left(\theta_{p}^{\circ n+1}\left(\mathfrak{z}\right)\right)-\hat{X}_{\mathbf{m}}\left(0\right)\right)
\end{equation}
Generalizing the notation slightly to clean up the presentation, the
anatomical layout of the limit (\ref{eq:limit}) is as follows:
\begin{equation}
\lim_{N\rightarrow\infty}\sum_{n=0}^{N-2}M_{n}\left(\mathfrak{z}\right)\rho\left(\left[\theta_{p}^{\circ n}\left(\mathfrak{z}\right)\right]_{p}\right)\tilde{X}_{N-n-1}\left(\theta_{p}^{\circ n+1}\left(\mathfrak{z}\right)\right)\label{eq:anatomy}
\end{equation}
where:
\begin{itemize}
\item (\ref{eq:anatomy}) is being considered with respect to a reffinite
$R$-frame $\mathcal{F}$ on $\mathbb{Z}_{p}$, where $R$ is a Dedekind
domain;
\item $\left\{ M_{n}\right\} _{n\geq0}$ is an $R$-valued $p$-adic M-function;
\item $\rho$ is a function from $\mathbb{Z}_{p}/p\mathbb{Z}_{p}$ to $R$
(the former being identified with both $\mathbb{Z}/p\mathbb{Z}$ and
the set $\left\{ 0,\ldots,p-1\right\} $ equipped with the operation
of addition mod $p$);
\item $X\in C\left(\mathcal{F}\right)$ is $\mathcal{F}$-quasi-integrable,
and: 
\begin{equation}
\tilde{X}_{N-n-1}\left(\theta_{p}^{\circ n+1}\left(\mathfrak{z}\right)\right)=\sum_{\left|t\right|_{p}\leq p^{N-n-1}}\hat{X}\left(t\right)e^{2\pi i\left\{ t\theta_{p}^{\circ n+1}\left(\mathfrak{z}\right)\right\} _{p}}
\end{equation}
where $\hat{X}$ is \emph{any }Fourier transform of $X$ with respect
to $\mathcal{F}$.
\end{itemize}
Viewed in this way, we want to show is that the limit of (\ref{eq:anatomy})
as $N\rightarrow\infty$ is $\mathcal{F}$-convergent to:
\begin{equation}
\sum_{n=0}^{\infty}M_{n}\left(\left[\mathfrak{z}\right]_{p^{n}}\right)\rho\left(\left[\theta_{p}^{\circ n}\left(\mathfrak{z}\right)\right]_{p}\right)X\left(\theta_{p}^{\circ n+1}\left(\mathfrak{z}\right)\right)
\end{equation}
Our next lemma gives us sufficient conditions on $\left\{ M_{n}\right\} _{n\geq0}$,
$\mathcal{F}$, and $X$ in order to guarantee that we get the convergence
we seek.
\begin{lem}[Main Limit Lemma]
\label{lem:main limit lemma}Let $R$ be a Dedekind domain, let $p$
be an integer $\geq2$, let $\rho:\mathbb{Z}/p\mathbb{Z}\rightarrow R$
be a function, let $\left\{ M_{n}\right\} _{n\geq0}$ be a $p$-adic
$R$-valued M-function, let $\mathcal{F}$ be an evaluative, multiplicative,
reffinite $R$-frame on $\mathbb{Z}_{p}$, let $X\in C\left(\mathcal{F}\right)$
be $\mathcal{F}$-quasi-integrable, and let $\hat{X}$ be any Fourier
transform of $X$ with respect to $\mathcal{F}$. Fix an integer $n_{0}\geq1$,
and let $\left\{ F_{N}\right\} _{N\geq n_{0}}$ be the sequence in
$\mathcal{S}\left(\mathbb{Z}_{p},R\right)$ defined by:
\begin{equation}
F_{N}\left(\mathfrak{z}\right)\overset{\textrm{def}}{=}\sum_{n=0}^{N-1}\rho\left(\left[\theta_{p}^{\circ n}\left(\mathfrak{z}\right)\right]_{p}\right)M_{n}\left(\mathfrak{z}\right)\tilde{X}_{N-n-1}\left(\theta_{p}^{\circ n+1}\left(\mathfrak{z}\right)\right)
\end{equation}
Then:
\begin{equation}
\lim_{N\rightarrow\infty}F_{N}\left(\mathfrak{z}\right)\overset{\mathcal{F}}{=}\sum_{n=0}^{\infty}\rho\left(\left[\theta_{p}^{\circ n}\left(\mathfrak{z}\right)\right]_{p}\right)M_{n}\left(\mathfrak{z}\right)X\left(\theta_{p}^{\circ n+1}\left(\mathfrak{z}\right)\right)
\end{equation}
occurs, and the limit is $\mathcal{F}$-rising-continuous, whenever
the following conditions are satisfied:

I. $\left\{ M_{n}\right\} _{n\geq0}$ is strongly $\mathcal{F}$-summable.
\label{eq:frame-wise decay}

II. For each $\mathfrak{z}\in D\left(\mathcal{F}\right)$:
\begin{equation}
\sup_{n\geq0}\mathcal{F}_{\mathfrak{z}}\left(X\circ\theta_{p}^{\circ n}\right)<\infty\label{eq:shift bound}
\end{equation}

III. The sequence $\left\{ X-\tilde{X}_{N}\right\} _{N\geq0}$ has
$\left(\mathcal{F},M\right)$-type decay.
\end{lem}
\begin{rem}
Recall that we write:
\begin{equation}
\Delta_{N}^{\left(m\right)}\left\{ X\right\} \left(\mathfrak{z}\right)\overset{\textrm{def}}{=}X\left(\theta_{p}^{\circ m}\left(\mathfrak{z}\right)\right)-\tilde{X}_{N}\left(\theta_{p}^{\circ m}\left(\mathfrak{z}\right)\right)
\end{equation}
for all $m,N\geq0$ and all $\mathfrak{z}$.
\end{rem}
Proof: Let everything be as given, and suppose (I), (II), and (III)
all hold. The trick to the proof is to add and subtract $X\left(\theta_{p}^{\circ n+1}\left(\mathfrak{z}\right)\right)$
to the summand of $F_{N}\left(\mathfrak{z}\right)$:
\begin{align*}
F_{N}\left(\mathfrak{z}\right) & =\overbrace{\sum_{n=0}^{N-1}\rho\left(\left[\theta_{p}^{\circ n}\left(\mathfrak{z}\right)\right]_{p}\right)M_{n}\left(\mathfrak{z}\right)X\left(\theta_{p}^{\circ n+1}\left(\mathfrak{z}\right)\right)}^{\textrm{call this I}_{N}\left(\mathfrak{z}\right)}\\
 & +\underbrace{\sum_{n=0}^{N-1}\rho\left(\left[\theta_{p}^{\circ n}\left(\mathfrak{z}\right)\right]_{p}\right)M_{n}\left(\mathfrak{z}\right)\Delta_{N-n-1}^{\left(n+1\right)}\left\{ X\right\} \left(\mathfrak{z}\right)}_{\textrm{call this II}_{N}\left(\mathfrak{z}\right)}
\end{align*}

\begin{claim}
$\lim_{N\rightarrow\infty}\textrm{I}_{N}\left(\mathfrak{z}\right)$
is $\mathcal{F}$-convergent. Moreover, letting $\textrm{I}_{\infty}\left(\mathfrak{z}\right)$
denote the function to which the limit $\mathcal{F}$-converges, we
have that:
\begin{equation}
\textrm{I}_{\infty}\left(\mathfrak{z}\right)\overset{\mathcal{F}}{=}\sum_{n=0}^{\infty}\rho\left(\left[\theta_{p}^{\circ n}\left(\mathfrak{z}\right)\right]_{p}\right)M_{n}\left(\mathfrak{z}\right)X\left(\theta_{p}^{\circ n+1}\left(\mathfrak{z}\right)\right)
\end{equation}
and that $\textrm{I}_{\infty}\left(\mathfrak{z}\right)$ is rising-continuous
with respect to $\mathcal{F}$.

Proof of Claim: Fix $\mathfrak{z}\in\mathbb{Z}_{p}$. Since $\rho$'s
domain is finite, for each $\mathfrak{z}\in\mathbb{Z}_{p}$, there
exists a real constant $\left\Vert \rho\right\Vert _{\mathfrak{z}}>0$
depending only on $\mathfrak{z}$ and $\mathcal{F}$ so that:
\begin{equation}
\sup_{n\geq0}\mathcal{F}_{\mathfrak{z}}\left(\rho\left(\left[\theta_{p}^{\circ n}\left(\cdot\right)\right]_{p}\right)\right)\leq\left\Vert \rho\right\Vert _{\mathfrak{z}}
\end{equation}
Using this in conjunction with (I) and (II), we have that the $n$th
term of $\textrm{I}_{N}\left(\mathfrak{z}\right)$ satisfies:
\begin{equation}
\rho\left(\left[\theta_{p}^{\circ n}\left(\mathfrak{z}\right)\right]_{p}\right)M_{n}\left(\mathfrak{z}\right)X\left(\theta_{p}^{\circ n+1}\left(\mathfrak{z}\right)\right)\overset{\mathcal{F}}{\ll}M_{n}\left(\mathfrak{z}\right)
\end{equation}
for some constant depending only on $\rho$ and $\mathfrak{z}$. The
strong $\mathcal{F}$-convergence of: 
\begin{equation}
\sum_{n=0}^{\infty}M_{n}\left(\mathfrak{z}\right)\label{eq:geo}
\end{equation}
given by (I) tells us that the upper bound on the $n$th term of $\textrm{I}_{N}\left(\mathfrak{z}\right)$
tends to zero with respect to $\mathcal{F}_{\mathfrak{z}}$ for each
$\mathfrak{z}\in D\left(\mathcal{F}\right)$. When $\mathcal{F}_{\mathfrak{z}}$
is non-archimedean, this is sufficient to guarantee the convergence
of $\textrm{I}_{N}\left(\mathfrak{z}\right)$ to $\textrm{I}_{\infty}\left(\mathfrak{z}\right)$
in $\mathcal{F}\left(\mathfrak{z}\right)$ as $N\rightarrow\infty$.

On the other hand, if $\mathcal{F}_{\mathfrak{z}}$ is archimedean,
(\ref{eq:frame-wise decay})'s requirement that (\ref{eq:geo}) converge
absolutely with respect to $\mathcal{F}_{\mathfrak{z}}$ then shows
that $\textrm{I}_{N}\left(\mathfrak{z}\right)$ is dominated by an
absolutely convergent series, and thus that $\textrm{I}_{N}\left(\mathfrak{z}\right)$
converges to $\textrm{I}_{\infty}\left(\mathfrak{z}\right)$ in $\mathcal{F}\left(\mathfrak{z}\right)$
as $N\rightarrow\infty$.

Since $\mathcal{F}_{\mathfrak{z}}$ is either archimedean or non-archimedean
for each $\mathfrak{z}\in D\left(\mathcal{F}\right)$, and since $\mathfrak{z}$
was arbitrary, this proves the $\mathcal{F}$-convergence of $\textrm{I}_{N}$
to $\textrm{I}_{\infty}$.

This proves the Claim. $\checkmark$
\end{claim}
So, all that remains is to show the $\mathcal{F}$-convergence of
$\lim_{N\rightarrow\infty}\textrm{II}_{N}\left(\mathfrak{z}\right)$
to $0$. To do this, fix $\mathfrak{z}\in\mathbb{Z}_{p}$, let $N$
arbitrary and large, and let $\gamma:\mathbb{N}_{0}\rightarrow\mathbb{N}_{0}$
any non-decreasing function satisfying the properties: 
\begin{align}
\gamma\left(n\right) & \leq n,\textrm{ }\forall n\geq0\\
\lim_{n\rightarrow\infty}\gamma\left(n\right) & \overset{\mathbb{R}}{=}+\infty
\end{align}
We then split the $n$-sum in $\textrm{II}_{N}$ into two pieces:
\begin{align*}
\textrm{II}_{N}\left(\mathfrak{z}\right) & =\overbrace{\sum_{n=0}^{\gamma\left(N\right)-1}\rho\left(\left[\theta_{p}^{\circ n}\left(\mathfrak{z}\right)\right]_{p}\right)M_{n}\left(\mathfrak{z}\right)\Delta_{N-n-1}^{\left(n+1\right)}\left\{ X\right\} \left(\mathfrak{z}\right)}^{\textrm{Call this }\textrm{III}_{N}\left(\mathfrak{z}\right)}\\
 & +\underbrace{\sum_{n=\gamma\left(N\right)}^{N-1}\rho\left(\left[\theta_{p}^{\circ n}\left(\mathfrak{z}\right)\right]_{p}\right)M_{n}\left(\mathfrak{z}\right)\Delta_{N-n-1}^{\left(n+1\right)}\left\{ X\right\} \left(\mathfrak{z}\right)}_{\textrm{Call this }\textrm{IV}_{N}\left(\mathfrak{z}\right)}
\end{align*}
To see what will be done, let $A_{n}$ denote $\rho\left(\left[\theta_{p}^{\circ n}\left(\mathfrak{z}\right)\right]_{p}\right)M_{n}\left(\mathfrak{z}\right)$
and let $B_{N,n}$ denote $\Delta_{N-n-1}^{\left(n+1\right)}\left\{ X\right\} \left(\mathfrak{z}\right)$.
As we saw in our investigation of $\textrm{I}_{\infty}$, hypothesis
(I) guarantees that $A_{n}$ tends to zero $\mathcal{F}$-wise as
$n\rightarrow\infty$. $B_{N,n}$'s behavior, however, is more complicated.
Heuristically, since $n$ is relatively small in $\textrm{III}_{N}$
(that is, $n<\gamma\left(N\right)$), $A_{n}$ will be bounded in
$\textrm{III}_{N}$ while $B_{N,n}$ will be small. On the other hand,
in $\textrm{IV}_{N}$, where $n$ is ``large'' ($n\geq\gamma\left(N\right)$),
$A_{n}$ will be small and $B_{N,n}$ will be bounded. Together, both
$\textrm{III}_{N}$ and $\textrm{IV}_{N}$ end up being small, and
therefore tend to $0$ as $N\rightarrow\infty$. The complications
arise from unravelling $B_{N,n}$.

In all of this, note that we can ignore $\rho\left(\left[\theta_{p}^{\circ n}\left(\mathfrak{z}\right)\right]_{p}\right)$,
as it will be bounded by $\left\Vert \rho\right\Vert _{\mathfrak{z}}$
uniformly with respect to $n$. So, we will.

\vphantom{}\textbullet{} \textbf{Bounding $\textrm{III}_{N}$}: By
hypothesis (I), $\textrm{III}_{N}$'s $A_{n}$ is bounded uniformly
with respect to $n$. As for $B_{N,n}$, we have:
\begin{equation}
B_{N,n}=\Delta_{N-n-1}^{\left(n+1\right)}\left\{ X\right\} \left(\mathfrak{z}\right)=\tilde{X}_{N-n-1}\left(\theta_{p}^{\circ n+1}\left(\mathfrak{z}\right)\right)-X\left(\theta_{p}^{\circ n+1}\left(\mathfrak{z}\right)\right)
\end{equation}
By (III), we know there is a real-valued $p$-adic M-function $\left\{ E_{\mathfrak{z},n}\right\} _{n\geq0}$
so that:
\begin{equation}
\max_{0\leq n\leq\gamma\left(N\right)-1}\mathcal{F}_{\mathfrak{z}}\left(\Delta_{N-n-1}^{\left(n+1\right)}\left\{ X\right\} \right)\ll\max_{n\leq\gamma\left(N\right)-1}\mathcal{F}_{\mathfrak{z}}\left(E_{\mathfrak{z},N-n-1}\circ\theta_{p}^{\circ n+1}\right)\textrm{ as }N\rightarrow\infty\label{eq:by III}
\end{equation}
where the constant of proportionality depends on $\mathcal{F}$ and
$\mathfrak{z}$, and is \emph{independent} of both $n$ and $N$.
\begin{claim}
\ 
\begin{equation}
\lim_{N\rightarrow\infty}\max_{n\leq\gamma\left(N\right)-1}E_{\mathfrak{z},N-n-1}\left(\theta_{p}^{\circ n+1}\left(\mathfrak{z}\right)\right)\overset{\mathbb{R}}{=}0\label{eq:kappa estimate for III}
\end{equation}

Proof of Claim: Since $\left\{ E_{\mathfrak{z},m}\right\} _{m\geq0}$
is an M-function, \textbf{Proposition \ref{prop:Kappa shift equation}}
allows us to write:
\begin{equation}
\mathcal{F}_{\mathfrak{z}}\left(E_{\mathfrak{z},N-n-1}\circ\theta_{p}^{\circ n+1}\right)=\mathcal{F}_{\mathfrak{z}}\left(\frac{E_{\mathfrak{z},N-n-1+n+1}}{E_{\mathfrak{z},n+1}}\right)\overset{!}{=}\frac{\mathcal{F}_{\mathfrak{z}}\left(E_{\mathfrak{z},N}\right)}{\mathcal{F}_{\mathfrak{z}}\left(E_{\mathfrak{z},n+1}\right)}\label{eq:kappa shift}
\end{equation}
The equality (!) holds because $\mathcal{F}$ is multiplicative and
evaluative.

As $\lim_{n\rightarrow\infty}\mathcal{F}_{\mathfrak{z}}\left(E_{\mathfrak{z},n}\right)\overset{\mathbb{R}}{=}0$,
since, in $\textrm{III}_{N}$, $n$ is restricted to the set $\left\{ 0,\ldots,\gamma\left(N\right)-1\right\} $,
we have that the value of $\mathcal{F}_{\mathfrak{z}}\left(E_{\mathfrak{z},n+1}\right)$
will be minimized by setting $n$ to be as large as possible; in this
case, $n=\gamma\left(N\right)-1$. Using (\ref{eq:kappa shift}) on
the right-hand side of (\ref{eq:by III}) then gives us the upper
bound:
\begin{equation}
\max_{0\leq n\leq\gamma\left(N\right)-1}\mathcal{F}_{\mathfrak{z}}\left(\Delta_{N-n-1}^{\left(n+1\right)}\left\{ X\right\} \right)\ll\sup_{n\leq\gamma\left(N\right)-1}\mathcal{F}_{\mathfrak{z}}\left(\frac{E_{\mathfrak{z},N}}{E_{\mathfrak{z},n+1}}\right)=\frac{\mathcal{F}_{\mathfrak{z}}\left(E_{\mathfrak{z},N}\right)}{\mathcal{F}_{\mathfrak{z}}\left(E_{\mathfrak{z},\gamma\left(N\right)}\right)}\label{eq:III estimate}
\end{equation}
which holds for all sufficiently large $N$. Note the crucial use
of the fact that the constant of proportionality is independent of
$n$ and $N$.

Finally, to conclude, note that since the $E_{\mathfrak{z},N}$s are
M-functions, there are no $N$ and $\mathfrak{y}$ for which $E_{\mathfrak{z},N}\left(\mathfrak{y}\right)=0$.
Since $\mathcal{F}$ is multiplicative and evaluative, this guarantees
that $\mathcal{F}_{\mathfrak{z}}\left(E_{\mathfrak{z},N}\right)\neq0$
for any $N$. Combined with the decay $\lim_{N\rightarrow\infty}\mathcal{F}_{\mathfrak{z}}\left(E_{\mathfrak{z},N}\right)\overset{\mathbb{R}}{=}0$
required by (III), we see that $\left\{ \mathcal{F}_{\mathfrak{z}}\left(E_{\mathfrak{z},n}\right)\right\} _{n\geq0}$
satisfies the hypotheses of \textbf{Proposition \ref{prop:quotient prop}.}
to choose $\gamma\left(N\right)$ so that (\ref{eq:III estimate})
tends to $0$ as $N\rightarrow\infty$.

This proves the Claim. $\checkmark$
\end{claim}
Thus, we see that the $n$th term of $\textrm{III}_{N}\left(\mathfrak{z}\right)$
is bounded frame-wise by:
\begin{equation}
\underbrace{\rho\left(\left[\theta_{p}^{\circ n}\left(\mathfrak{z}\right)\right]_{p}\right)M_{n}\left(\mathfrak{z}\right)}_{\textrm{bounded w.r.t. }n}\underbrace{\Delta_{N-n-1}^{\left(n+1\right)}\left\{ X\right\} \left(\mathfrak{z}\right)}_{\rightarrow0\textrm{ as }N\rightarrow\infty}
\end{equation}
which then decays to $0$. In particular, this shows that:
\begin{equation}
\lim_{N\rightarrow\infty}\max_{n\leq\gamma\left(N\right)-1}\rho\left(\left[\theta_{p}^{\circ n}\left(\mathfrak{z}\right)\right]_{p}\right)M_{n}\left(\mathfrak{z}\right)\Delta_{N-n-1}^{\left(n+1\right)}\left\{ X\right\} \left(\mathfrak{z}\right)\overset{\mathbb{R}}{=}0
\end{equation}
In the case where $\mathcal{F}_{\mathfrak{z}}$ is non-archimedean,
this proves the $n$th term of $\textrm{III}_{N}\left(\mathfrak{z}\right)$
tends to $0$ as $N\rightarrow\infty$ uniformly with respect to $n$,
thus making the whole series vanish as $N\rightarrow\infty$. Meanwhile,
in the archimedean case, we can write:
\begin{equation}
\mathcal{F}_{\mathfrak{z}}\left(\textrm{III}_{N}\right)\ll_{\mathfrak{z}}\underbrace{\mathcal{F}_{\mathfrak{z}}\left(\sum_{m=0}^{\gamma\left(N\right)-1}\rho\left(\left[\theta_{p}^{\circ m}\left(\cdot\right)\right]_{p}\right)M_{n}\left(\cdot\right)\right)}_{\textrm{bounded as }N\rightarrow\infty}\times\underbrace{\max_{0\leq n\leq\gamma\left(N\right)-1}\mathcal{F}_{\mathfrak{z}}\left(\Delta_{N-n-1}^{\left(n+1\right)}\right)\left\{ X\right\} }_{\textrm{vanishes as }N\rightarrow\infty}
\end{equation}
which shows the $\mathcal{F}$-wise vanishing of $\textrm{III}_{N}$
as $N\rightarrow\infty$.

\vphantom{}\textbf{\textbullet{} Bounding $\textrm{IV}_{N}$}: By
(III), we know there is a real-valued $p$-adic M-function $\left\{ E_{\mathfrak{z},n}\right\} _{n\geq0}$
so that: 
\begin{equation}
\mathcal{F}_{\mathfrak{z}}\left(\Delta_{N-n-1}^{\left(n+1\right)}\left\{ X\right\} \right)\ll_{\mathfrak{z}}\mathcal{F}_{\mathfrak{z}}\left(E_{\mathfrak{z},N-n+1}\circ\theta_{p}^{\circ n+1}\right)
\end{equation}
for all $N$ and all $n\in\left\{ \gamma\left(N\right),\ldots,N-n_{0}\right\} $,
where, as indicated, the constant of proportionality depends only
on $\mathfrak{z}$. Once again, since $\mathcal{F}_{\mathfrak{z}}$
is multiplicative and evaluative, we get (\ref{eq:kappa shift}),
and so, taking the supremum over $n\in\left\{ \gamma\left(N\right),\ldots,N-n_{0}\right\} $,
we have:
\begin{equation}
\sup_{\gamma\left(N\right)\leq n<N}\mathcal{F}_{\mathfrak{z}}\left(\Delta_{N-n-1}^{\left(n+1\right)}\left\{ X\right\} \right)\ll\sup_{\gamma\left(N\right)\leq n\leq N-1}\frac{\mathcal{F}_{\mathfrak{z}}\left(E_{\mathfrak{z},N}\right)}{\mathcal{F}_{\mathfrak{z}}\left(E_{\mathfrak{z},n+1}\right)}\overset{!}{=}\frac{\mathcal{F}_{\mathfrak{z}}\left(E_{\mathfrak{z},N}\right)}{\mathcal{F}_{\mathfrak{z}}\left(E_{\mathfrak{z},N}\right)}=1\label{eq:IV bound}
\end{equation}
Here, the equality marked (!) follows from the fact that since (III)
forces $\mathcal{F}_{\mathfrak{z}}\left(E_{\mathfrak{z},n+1}\right)$
to tend to $0$ as $n\rightarrow\infty$, we can minimize $\mathcal{F}_{\mathfrak{z}}\left(E_{\mathfrak{z},n+1}\right)$
by maximizing $N$.

Using (\ref{eq:IV bound}), we have:
\begin{equation}
\max_{\gamma\left(N\right)\leq n\leq N-1}\rho\left(\left[\theta_{p}^{\circ n}\left(\mathfrak{z}\right)\right]_{p}\right)M_{n}\left(\mathfrak{z}\right)\Delta_{N-n-1}^{\left(n+1\right)}\left\{ X\right\} \ll_{\mathfrak{z}}M_{\gamma\left(N\right)}\left(\mathfrak{z}\right)
\end{equation}
Thus, if $\mathcal{F}_{\mathfrak{z}}$ is non-archimedean, we get
that the $n$th term of $\textrm{IV}_{N}\left(\mathfrak{z}\right)$
tends to $0$ as $N\rightarrow\infty$, and does so uniformly with
respect to $n$, which then guarantees the convergence of $\textrm{IV}_{N}\left(\mathfrak{z}\right)$
to $0$. Alternatively, if $\mathcal{F}_{\mathfrak{z}}$ is archimedean,
we can write:
\begin{equation}
\mathcal{F}_{\mathfrak{z}}\left(\textrm{IV}_{N}\right)\ll_{\mathfrak{z}}\sum_{n=\gamma\left(N\right)}^{N-1}\mathcal{F}_{\mathfrak{z}}\left(M_{n}\right)\leq\sum_{n=\gamma\left(N\right)}^{\infty}\mathcal{F}_{\mathfrak{z}}\left(M_{n}\right)<\infty
\end{equation}
for all $N$, using hypothesis (I). As the upper bound vanishes as
$N\rightarrow\infty$, this then proves the convergence of $\textrm{IV}_{N}\left(\mathfrak{z}\right)$
to $0$ in the archimedean case.

Since $\mathfrak{z}\in D\left(\mathcal{F}\right)$ was arbitrary,
this shows that $\textrm{II}_{N}$ converges to $0$ with respect
to $\mathcal{F}$ as $N\rightarrow\infty$, and thus $F_{N}\left(\mathfrak{z}\right)$
$\mathcal{F}$-converges to $F\left(\mathfrak{z}\right)$ as $N\rightarrow\infty$.

Q.E.D.

\subsection{\label{subsec:Assumption--=000026}Assumption \ref{assu:main limit lemma}
\& the Inductive Step}

With \textbf{Lemma \ref{lem:main limit lemma} }now in hand, we know
the conditions required in order to show that the limits of the type
considered in \textbf{Assumption \ref{assu:main limit lemma}} hold
rigorously. In order to prove the quasi-integrability of $X_{\mathbf{n}}$
and the other $X_{\mathbf{m}}$s, we need to figure out the conditions
needed in order to apply \textbf{Lemma \ref{lem:main limit lemma}},
and then work out the base case of the inductive argument that will
show that the conditions being satisfied for all $X_{\mathbf{k}}$
for all $\mathbf{k}<\mathbf{m}$ for a given $\mathbf{m}$ imply the
necessary conditions are satisfied by $X_{\mathbf{m}}$ as well.

So, returning to \textbf{Assumption \ref{assu:main limit lemma}},
there, we claimed that
\begin{equation}
\tilde{f}_{\mathbf{n},N}\left(\mathfrak{z}\right)=\sum_{n=0}^{N-2}r_{\mathbf{n},0}^{n}\kappa_{\mathbf{n}}\left(\left[\mathfrak{z}\right]_{p^{n}}\right)\sum_{\mathbf{m}<\mathbf{n}}r_{\mathbf{m},\left[\theta_{p}^{\circ n}\left(\mathfrak{z}\right)\right]_{p}}\tilde{X}_{\mathbf{m},N-n-1}\left(\theta_{p}^{\circ n+1}\left(\mathfrak{z}\right)\right)
\end{equation}
tends to:
\begin{equation}
\sum_{n=0}^{\infty}r_{\mathbf{n},0}^{n}\kappa_{\mathbf{n}}\left(\left[\mathfrak{z}\right]_{p^{n}}\right)\sum_{\mathbf{m}<\mathbf{n}}r_{\mathbf{m},\left[\theta_{p}^{\circ n}\left(\mathfrak{z}\right)\right]_{p}}X_{\mathbf{m}}\left(\theta_{p}^{\circ n+1}\left(\mathfrak{z}\right)\right)
\end{equation}
as $N\rightarrow\infty$. Here, the convergence will be frame-wise.
Looking at \textbf{Lemma \ref{lem:main limit lemma}}'s\textbf{ }hypotheses,
we deduce a collection of conditions sufficient to apply our lemma.
\begin{prop}
\textbf{\label{prop:dealing with main limit assumption}Assumption
\ref{assu:main limit lemma}} will hold with $\mathcal{F}$-convergence
whenever there is a multiplicative, evaluative, reffinite $\mathcal{R}$-frame
on $\mathbb{Z}_{p}$ so that:

\vphantom{}

I. The M-function $\left(r_{\mathbf{n},0},\kappa_{\mathbf{n}}\right)$
is strongly $\mathcal{F}$-summable.

\vphantom{}

II. For each $\mathfrak{z}\in D\left(\mathcal{F}\right)$, $\max_{1\leq j\leq d}\sup_{n\geq0}\mathcal{F}_{\mathfrak{z}}\left(X_{j}\circ\theta_{p}^{\circ n}\right)<\infty$.

\vphantom{}

III. For each $\mathbf{m}<\mathbf{n}$, $X_{\mathbf{m}}$ is $\mathcal{F}$-quasi-integrable
with a Fourier transform $\hat{X}_{\mathbf{m}}$ so that $\Delta_{N}^{\left(0\right)}\left\{ X_{\mathbf{m}}\right\} $
has $\left(\mathcal{F},M\right)$-type decay with respect to $N$.
\end{prop}
Proof: The only thing that isn't immediately obvious is (II). Applying
\textbf{Lemma \ref{lem:main limit lemma} }directly to justify \textbf{Assumption
\ref{assu:main limit lemma}} would require:
\begin{equation}
\max_{\mathbf{m}<\mathbf{n}}\sup_{n\geq0}\mathcal{F}_{\mathfrak{z}}\left(X_{\mathbf{m}}\circ\theta_{p}^{\circ n}\right)<\infty
\end{equation}
 However, by definition:
\begin{equation}
X_{\mathbf{m}}=\prod_{j=1}^{d}X_{j}^{m_{j}}
\end{equation}
where the $X_{j}$s are degree $1$ F-series, each of which is $\mathcal{F}$-compatible.
Since\emph{ $\mathcal{F}_{\mathfrak{z}}$ }is multiplicative, and
since each $X_{j}$ is $\mathcal{F}$-compatible, we have:
\begin{equation}
\mathcal{F}_{\mathfrak{z}}\left(X_{\mathbf{m}}\circ\theta_{p}^{\circ n}\right)=\prod_{j=1}^{d}\left(\mathcal{F}_{\mathfrak{z}}\left(X_{j}\circ\theta_{p}^{\circ n}\right)\right)^{m_{j}}
\end{equation}
and so:
\begin{align*}
\max_{\mathbf{m}<\mathbf{n}}\sup_{n\geq0}\mathcal{F}_{\mathfrak{z}}\left(X_{\mathbf{m}}\circ\theta_{p}^{\circ n}\right) & =\max_{\mathbf{m}<\mathbf{n}}\sup_{n\geq0}\prod_{j=1}^{d}\left(\mathcal{F}_{\mathfrak{z}}\left(X_{j}\circ\theta_{p}^{\circ n}\right)\right)^{m_{j}}\\
 & =\max_{\mathbf{m}<\mathbf{n}}\sup_{n\geq0}\left(\max_{1\leq j\leq d}\mathcal{F}_{\mathfrak{z}}\left(X_{j}\circ\theta_{p}^{\circ n}\right)\right)^{\Sigma\left(\mathbf{m}\right)}\\
 & \leq\max_{\mathbf{m}<\mathbf{n}}\max_{1\leq j\leq d}\left(\sup_{n\geq0}\mathcal{F}_{\mathfrak{z}}\left(X_{j}\circ\theta_{p}^{\circ n}\right)\right)^{\Sigma\left(\mathbf{m}\right)}\\
 & \leq\left(\max_{1\leq j\leq d}\sup_{n\geq0}\mathcal{F}_{\mathfrak{z}}\left(X_{j}\circ\theta_{p}^{\circ n}\right)\right)^{\Sigma\left(\mathbf{n}\right)}
\end{align*}
As $\Sigma\left(\mathbf{n}\right)$ is finite, we see that $\max_{1\leq j\leq d}\sup_{n\geq0}\mathcal{F}_{\mathfrak{z}}\left(X_{j}\circ\theta_{p}^{\circ n}\right)<\infty$
is then sufficient to give us what we need.

Q.E.D.

\vphantom{}Of the three conditions given by \textbf{Proposition \ref{prop:dealing with main limit assumption}},
(I) and (II) do not require demonstrating, or even \emph{assuming
}the quasi-integrability of $X_{\mathbf{m}}$ for any $\mathbf{m}<\mathbf{n}$.
The quasi-integrability of the $X_{\mathbf{m}}$s is only relevant
to (III). As such, our proof of the quasi-integrability of the $X_{\mathbf{m}}$s
will proceed by strong induction with (III). For the inductive step,
we will assume that (III) holds for all proper subsets $\mathbf{m}$
of $\mathbf{n}$ and then use this to show that (III) holds for $X_{\mathbf{n}}$
as well.

To do this, we need to compute $\Delta_{N}^{\left(0\right)}\left\{ X_{\mathbf{n}}\right\} $.
Since $X_{\mathbf{n}}=f_{\mathbf{n}}-g_{\mathbf{n}}$, we tackle $f_{\mathbf{n}}$
and $g_{\mathbf{n}}$ separately. For $f_{\mathbf{n}}$, we pull out
the $n=0$th through $n=N-2$th terms from its series representation
to write:

\begin{align*}
\Delta_{N}^{\left(0\right)}\left\{ f_{\mathbf{n}}\right\} \left(\mathfrak{z}\right) & =f_{\mathbf{n}}\left(\mathfrak{z}\right)-\tilde{f}_{\mathbf{n},N}\left(\mathfrak{z}\right)\\
 & =\overbrace{\sum_{n=0}^{\infty}r_{\mathbf{n},0}^{n}\kappa_{\mathbf{n}}\left(\left[\mathfrak{z}\right]_{p^{n}}\right)\sum_{\mathbf{m}<\mathbf{n}}r_{\mathbf{m},\left[\theta_{p}^{\circ n}\left(\mathfrak{z}\right)\right]_{p}}\left(X_{\mathbf{m}}\left(\theta_{p}^{\circ n+1}\left(\mathfrak{z}\right)\right)-\hat{X}_{\mathbf{m}}\left(0\right)\right)}^{\textrm{pull out terms from here to cancel \ensuremath{\hat{X}_{\mathbf{m}}\left(0\right)} in the line below}}\\
 & -\sum_{n=0}^{N-2}r_{\mathbf{n},0}^{n}\kappa_{\mathbf{n}}\left(\left[\mathfrak{z}\right]_{p^{n}}\right)\sum_{\mathbf{m}<\mathbf{n}}r_{\mathbf{m},\left[\theta_{p}^{\circ n}\left(\mathfrak{z}\right)\right]_{p}}\left(\tilde{X}_{\mathbf{m},N-n-1}\left(\theta_{p}^{\circ n+1}\left(\mathfrak{z}\right)\right)-\hat{X}_{\mathbf{m}}\left(0\right)\right)\\
 & =\sum_{n=N-1}^{\infty}r_{\mathbf{n},0}^{n}\kappa_{\mathbf{n}}\left(\left[\mathfrak{z}\right]_{p^{n}}\right)\sum_{\mathbf{m}<\mathbf{n}}r_{\mathbf{m},\left[\theta_{p}^{\circ n}\left(\mathfrak{z}\right)\right]_{p}}\left(X_{\mathbf{m}}\left(\theta_{p}^{\circ n+1}\left(\mathfrak{z}\right)\right)-\hat{X}_{\mathbf{m}}\left(0\right)\right)\\
 & +\sum_{n=0}^{N-2}r_{\mathbf{n},0}^{n}\kappa_{\mathbf{n}}\left(\left[\mathfrak{z}\right]_{p^{n}}\right)\sum_{\mathbf{m}<\mathbf{n}}r_{\mathbf{m},\left[\theta_{p}^{\circ n}\left(\mathfrak{z}\right)\right]_{p}}\underbrace{\left(X_{\mathbf{m}}\left(\theta_{p}^{\circ n+1}\left(\mathfrak{z}\right)\right)-\tilde{X}_{\mathbf{m},N-n-1}\left(\theta_{p}^{\circ n+1}\left(\mathfrak{z}\right)\right)\right)}_{\Delta_{N-n-1}^{\left(n+1\right)}\left\{ X_{\mathbf{m}}\right\} \left(\mathfrak{z}\right)}
\end{align*}
We then compute $\Delta_{N}^{\left(0\right)}\left\{ g_{\mathbf{n}}\right\} $
using \textbf{Proposition} \textbf{\ref{prop:iterative triangle}}.
Combining this with $\Delta_{N}^{\left(0\right)}\left\{ f_{\mathbf{n}}\right\} $
yields:
\begin{align*}
\Delta_{N}^{\left(0\right)}\left\{ X_{\mathbf{n}}\right\} \left(\mathfrak{z}\right) & =\Delta_{N}^{\left(0\right)}\left\{ f_{\mathbf{n}}\right\} \left(\mathfrak{z}\right)-\Delta_{N}^{\left(0\right)}\left\{ g_{\mathbf{n}}\right\} \left(\mathfrak{z}\right)\\
 & =\sum_{n=0}^{N-2}r_{\mathbf{n},0}^{n}\kappa_{\mathbf{n}}\left(\left[\mathfrak{z}\right]_{p^{n}}\right)\sum_{\mathbf{m}<\mathbf{n}}r_{\mathbf{m},\left[\theta_{p}^{\circ n}\left(\mathfrak{z}\right)\right]_{p}}\Delta_{N-n-1}^{\left(n+1\right)}\left\{ X_{J}\right\} \left(\mathfrak{z}\right)\\
 & +\sum_{n=N-1}^{\infty}r_{\mathbf{n},0}^{n}\kappa_{\mathbf{n}}\left(\left[\mathfrak{z}\right]_{p^{n}}\right)\sum_{\mathbf{m}<\mathbf{n}}r_{\mathbf{m},\left[\theta_{p}^{\circ n}\left(\mathfrak{z}\right)\right]_{p}}\left(X_{\mathbf{m}}\left(\theta_{p}^{\circ n+1}\left(\mathfrak{z}\right)\right)-\hat{X}_{\mathbf{m}}\left(0\right)\right)\\
 & -\left(\Delta_{0}^{\left(N\right)}\left\{ g_{\mathbf{n}}\right\} \left(\mathfrak{z}\right)+\left[\alpha_{\mathbf{n}}\left(0\right)=1\right]N\beta_{\mathbf{n}}\left(0\right)\right)r_{\mathbf{n},0}^{N}\kappa_{\mathbf{n}}\left(\left[\mathfrak{z}\right]_{p^{N}}\right)
\end{align*}
Now, we take estimates. The various pieces of import are labeled below
using mostly roman numerals

\begin{align}
\Delta_{N}^{\left(0\right)}\left\{ X_{\mathbf{n}}\right\} \left(\mathfrak{z}\right) & \overset{\mathcal{F}}{\ll}\underbrace{\sum_{n=0}^{N-2}\overbrace{r_{\mathbf{n},0}^{n}\kappa_{\mathbf{n}}\left(\left[\mathfrak{z}\right]_{p^{n}}\right)}^{\textrm{ia}_{N}\left(\mathfrak{z}\right)}\overbrace{\sum_{\mathbf{m}<\mathbf{n}}r_{\mathbf{m},\left[\theta_{p}^{\circ n}\left(\mathfrak{z}\right)\right]_{p}}\Delta_{N-n-1}^{\left(n+1\right)}\left\{ X_{\mathbf{m}}\right\} \left(\mathfrak{z}\right)}^{\textrm{ib}_{N}\left(\mathfrak{z}\right)}}_{(\textrm{i}_{N}\left(\mathfrak{z}\right))}\label{eq:key bound}\\
 & +\underbrace{\max_{\mathbf{m}<\mathbf{n}}\max_{m\geq N}\left(X_{\mathbf{m}}\left(\theta_{p}^{\circ m}\left(\mathfrak{z}\right)\right)-\hat{X}_{\mathbf{m}}\left(0\right)\right)}_{\textrm{ii}_{N}\left(\mathfrak{z}\right)}\times\underbrace{\sum_{n=N-1}^{\infty}r_{\mathbf{n},0}^{n}\kappa_{\mathbf{n}}\left(\left[\mathfrak{z}\right]_{p^{n}}\right)}_{\textrm{iii}_{N}\left(\mathfrak{z}\right)}\nonumber \\
 & +\underbrace{\left(\Delta_{0}^{\left(N\right)}\left\{ g_{\mathbf{n}}\right\} \left(\mathfrak{z}\right)+\left[\alpha_{\mathbf{n}}\left(0\right)=1\right]N\beta_{\mathbf{n}}\left(0\right)\right)}_{\textrm{iv}_{N}\left(\mathfrak{z}\right)}\times\underbrace{r_{\mathbf{n},0}^{N}\kappa_{\mathbf{n}}\left(\left[\mathfrak{z}\right]_{p^{N}}\right)}_{\textrm{v}_{N}\left(\mathfrak{z}\right)}\nonumber 
\end{align}

Now, by strong induction, \emph{suppose that the hypotheses of }\textbf{\emph{Proposition
\ref{prop:dealing with main limit assumption}}}\emph{ hold}. As  $\textrm{i}_{N}\left(\mathfrak{z}\right)$
is the most complicated of the terms to deal with, we'll save it for
last. As for the others:
\begin{itemize}
\item $\textrm{ii}_{N}\left(\mathfrak{z}\right)$: Because $\hat{X}_{\mathbf{m}}\left(0\right)$
is a constant and because $\mathbf{n}$ has finitely many subsets,
we have:
\begin{align*}
\overbrace{\max_{\mathbf{m}<\mathbf{n}}\max_{m\geq N}\left(X_{\mathbf{m}}\left(\theta_{p}^{\circ m}\left(\mathfrak{z}\right)\right)-\hat{X}_{\mathbf{m}}\left(0\right)\right)}^{\textrm{ii}_{N}\left(\mathfrak{z}\right)} & \overset{\mathcal{F}}{\ll}\max_{\mathbf{m}<\mathbf{n}}\max_{m\geq N}X_{\mathbf{m}}\left(\theta_{p}^{\circ m}\left(\mathfrak{z}\right)\right)\\
\left(\mathcal{F}_{\mathfrak{z}}\textrm{ is a seminorm};X_{\mathbf{m}}=\prod_{j=1}^{d}X_{j}^{m_{j}}\right); & \overset{\mathcal{F}}{\ll}\max_{1\leq j\leq d}\max_{m\geq N}X_{j}^{m_{j}}\left(\theta_{p}^{\circ m}\left(\mathfrak{z}\right)\right)\\
 & \overset{\mathcal{F}}{\ll}\max_{1\leq j\leq d}\max_{m\geq N}X_{j}\left(\theta_{p}^{\circ m}\left(\mathfrak{z}\right)\right)
\end{align*}
the bottom line of which is bounded by hypothesis (II) of \textbf{Proposition
\ref{prop:dealing with main limit assumption}}.
\item $\textrm{iii}_{N}\left(\mathfrak{z}\right)$: By hypothesis (I) of
\textbf{Proposition \ref{prop:dealing with main limit assumption}},
$\left(r_{\mathbf{n},0},\kappa_{\mathbf{n}}\right)$ is strongly $\mathcal{F}$-summable,
which gives us the estimate;
\begin{equation}
\textrm{iii}_{N}\left(\mathfrak{z}\right)=\sum_{n=N-1}^{\infty}r_{\mathbf{n},0}^{n}\kappa_{\mathbf{n}}\left(\left[\mathfrak{z}\right]_{p^{n}}\right)\overset{\mathcal{F}}{\ll}r_{\mathbf{n},0}^{N}\kappa_{\mathbf{n}}\left(\left[\mathfrak{z}\right]_{p^{N}}\right)
\end{equation}
the upper bound of which tends to $0$.
\item $\textrm{iv}_{N}\left(\mathfrak{z}\right)$: We have:
\begin{equation}
\textrm{iv}_{N}\left(\mathfrak{z}\right)=g_{\mathbf{n}}\left(\theta_{p}^{\circ N}\left(\mathfrak{z}\right)\right)+\left[\alpha_{\mathbf{n}}\left(0\right)=1\right]N\beta_{\mathbf{n}}\left(0\right)
\end{equation}
When $\mathcal{F}_{\mathfrak{z}}$ is non-archimedean, $\max_{N\geq0}\mathcal{F}_{\mathfrak{z}}\left(N\right)<\infty$,
which gives us:
\begin{equation}
\textrm{iv}_{N}\left(\mathfrak{z}\right)\overset{\mathcal{F}}{\ll}\max_{n\geq0}g_{\mathbf{n}}\left(\theta_{p}^{\circ n}\left(\mathfrak{z}\right)\right)\textrm{ as }N\rightarrow\infty
\end{equation}
Otherwise, when $\mathcal{F}_{\mathfrak{z}}$ is archimedean, we have:
\begin{equation}
\textrm{iv}_{N}\left(\mathfrak{z}\right)\overset{\mathcal{F}}{\ll}\max\left\{ N\left[\alpha_{\mathbf{n}}\left(0\right)=1\right],g_{\mathbf{n}}\left(\theta_{p}^{\circ N}\left(\mathfrak{z}\right)\right)\right\} 
\end{equation}
\item $\textrm{v}_{N}\left(\mathfrak{z}\right)$: As $\textrm{v}_{N}\left(\mathfrak{z}\right)=r_{\mathbf{n},0}^{N}\kappa_{\mathbf{n}}\left(\left[\mathfrak{z}\right]_{p^{N}}\right)$,
the assumed validity of hypothesis (I) of \textbf{Proposition \ref{prop:dealing with main limit assumption}}
gives us the desired decay to $0$ as $N\rightarrow\infty$.
\end{itemize}
Putting all of this together, we see that (\ref{eq:key bound}) reduces
to:

\begin{align}
\Delta_{N}^{\left(0\right)}\left\{ X_{\mathbf{n}}\right\} \left(\mathfrak{z}\right) & \overset{\mathcal{F}}{\ll}\underbrace{\sum_{n=0}^{N-2}\overbrace{r_{\mathbf{n},0}^{n}\kappa_{\mathbf{n}}\left(\left[\mathfrak{z}\right]_{p^{n}}\right)}^{\textrm{(ia)}}\overbrace{\sum_{\mathbf{m}<\mathbf{n}}r_{\mathbf{m},\left[\theta_{p}^{\circ n}\left(\mathfrak{z}\right)\right]_{p}}\Delta_{N-n-1}^{\left(n+1\right)}\left\{ X_{\mathbf{m}}\right\} \left(\mathfrak{z}\right)}^{\textrm{(ib)}}}_{\textrm{(i)}}\label{eq:key bound, simplified}\\
 & +\max\left\{ N\left[\alpha_{\mathbf{n}}\left(0\right)=1\right],g_{\mathbf{n}}\left(\theta_{p}^{\circ N}\left(\mathfrak{z}\right)\right)\right\} r_{\mathbf{n},0}^{N}\kappa_{\mathbf{n}}\left(\left[\mathfrak{z}\right]_{p^{N}}\right)\nonumber 
\end{align}

Conveniently, we can deal with both the lower line \emph{and }the
base case of our inductive argument (namely, when $\Sigma\left(\mathbf{n}\right)=1$)
in one fell swoop. This is done by establishing a variant of \textbf{Theorem
\ref{thm:The-breakdown-variety}}. Here, we will draw from the algebraic
set-up in \textbf{Section \ref{subsec:Ascent-from-Fractional}} so
as to ensure that all of the arguments to come will work by descent
from an $\mathcal{R}$-frame $\mathcal{F}$ to the quotient frame
$\mathcal{F}/I$ induced by an ideal $I$ of a polynomial ring.
\begin{thm}
\label{thm:mini thesis}Let $K$ be a global field, let $p$ be an
integer $\geq2$, let $a_{0},\ldots,a_{p-1}$ and $b_{0},\ldots,b_{p-1}$
be indeterminates, let $R$ denote the polynomial ring $\mathcal{O}_{K}\left[a_{0},\ldots,a_{p-1},b_{0},\ldots,b_{p-1}\right]$,
and let $\mathcal{R}$ denote the localization of $R$ away from the
prime ideal $\left\langle 1-a_{0}\right\rangle $. Let $\mathcal{F}$
be a reffinite $\mathcal{R}$-frame on $\mathbb{Z}_{p}$, let $\kappa_{X}:\mathbb{N}_{0}\rightarrow\textrm{Frac}R$
(where $\textrm{Frac}R$ is the field of fractions of $R$) be defined
by:
\begin{equation}
\kappa_{X}\left(n\right)\overset{\textrm{def}}{=}\prod_{j=1}^{p-1}\left(\frac{a_{j}}{a_{0}}\right)^{\#_{p:j}\left(n\right)},\textrm{ }\forall n\in\mathbb{N}_{0}
\end{equation}
Finally, let $I$ be any ideal of $R$ so that $\left\langle 1-a_{0}\right\rangle \nsubseteq I$
and:
\begin{equation}
\limsup_{n\rightarrow\infty}\left|a_{0}^{n}\kappa_{X}\left(\left[\mathfrak{z}\right]_{p^{n}}\right)\right|_{\left(\mathcal{F}/\mathfrak{I}\right)\left(\mathfrak{z}\right)}^{1/n}<1,\textrm{ }\forall\mathfrak{z}\in D\left(\mathcal{F}\right)\label{eq:root test condition}
\end{equation}
then the following hold:

I. The foliated $p$-adic F-series:

\begin{equation}
X\left(\mathfrak{z}\right)=\sum_{n=0}^{\infty}\left(\sum_{j=0}^{p-1}\hat{b}_{j}\left(\varepsilon_{n}\left(\mathfrak{z}\right)\right)^{j}\right)a_{0}^{n}\prod_{k=1}^{p-1}\left(a_{k}/a_{0}\right)^{\#_{k}\left(\left[\mathfrak{z}\right]_{p^{n}}\right)}\label{eq:Y}
\end{equation}
is both compatible and quasi-integrable with respect to $\mathcal{F}/I$,
where:
\begin{equation}
\hat{b}_{j}=\sum_{k=0}^{p-1}b_{k}e^{2\pi ijk/p}
\end{equation}
The formula for $\hat{X}$ given in the statement of \textbf{Theorem
\ref{thm:The-breakdown-variety}} is then a Fourier transform of $X$.
Moreover:
\begin{equation}
X\left(\mathfrak{z}\right)=a_{\left[\mathfrak{z}\right]}X\left(\theta_{p}\left(\mathfrak{z}\right)\right)+b_{\left[\mathfrak{z}\right]}\label{eq:func}
\end{equation}
holds for all $\mathfrak{z}\in D\left(\mathcal{F}/I\right)$. In particular,
we have:
\begin{equation}
\Delta_{N}^{\left(0\right)}\left\{ X\right\} \left(\mathfrak{z}\right)\overset{\mathcal{F}/I}{\ll}a_{0}^{N}\kappa_{X}\left(\left[\mathfrak{z}\right]_{p^{N}}\right)\textrm{ as }N\rightarrow\infty\label{eq:delta_N of X decay}
\end{equation}

\vphantom{}

II.
\begin{equation}
\sup_{N\geq0}\left(\mathcal{F}/I\right)_{\mathfrak{z}}\left(X\circ\theta_{p}^{\circ N}\right)<1,\textrm{ }\forall\mathfrak{z}\in D\left(\mathcal{F}\right)\label{eq:shift sup boundedness}
\end{equation}
Note that in (\ref{eq:root test condition}), the value of the limsup
is dependent on $\mathfrak{z}$, and thus might not be $<1$ uniformly
with respect to $\mathfrak{z}\in D\left(\mathcal{F}\right)$. Likewise,
in (\ref{eq:shift sup boundedness}), the bound is not uniform with
respect to $\mathfrak{z}$, either.
\end{thm}
Proof: Let everything be as given.

(\ref{eq:root test condition}) implies (\ref{eq:func}) by \textbf{Theorem
\ref{thm:variety frame general}}, which also gives us the compatibility
of the F-series (\ref{eq:Y}).

Next, since $\mathcal{F}$ is multiplicative, evaluative, and reffinite,
so is $\mathcal{F}/I$, and we have:
\begin{equation}
\left(\mathcal{F}/I\right)_{\mathfrak{z}}\left(X\circ\theta_{p}^{\circ N}\right)=\left|X\left(\theta_{p}^{\circ N}\left(\mathfrak{z}\right)\right)\right|_{\left(\mathcal{F}/I\right)\left(\theta_{p}^{\circ N}\left(\mathfrak{z}\right)\right)}=\left|X\left(\theta_{p}^{\circ N}\left(\mathfrak{z}\right)\right)\right|_{\left(\mathcal{F}/I\right)\left(\mathfrak{z}\right)}
\end{equation}
for all $N\geq0$ and all $\mathfrak{z}\in D\left(\mathcal{F}\right)$,
and so, by \textbf{Proposition \ref{prop:Kappa shift equation}}:
\begin{align*}
\left(\mathcal{F}/I\right)_{\mathfrak{z}}\left(X\circ\theta_{p}^{\circ N}\right) & \ll\sum_{n=0}^{\infty}\left|a_{0}^{n}\kappa_{X}\left(\left[\theta_{p}^{\circ N}\left(\mathfrak{z}\right)\right]_{p^{n}}\right)\right|_{\left(\mathcal{F}/I\right)\left(\mathfrak{z}\right)}\\
 & =\sum_{n=0}^{\infty}\left|a_{0}^{n}\frac{\kappa_{X}\left(\left[\mathfrak{z}\right]_{p^{N+n}}\right)}{\kappa_{X}\left(\left[\mathfrak{z}\right]_{p^{N}}\right)}\right|_{\left(\mathcal{F}/I\right)\left(\mathfrak{z}\right)}\\
\left(\times\frac{a_{0}^{N}}{a_{0}^{N}}\right); & =\sum_{n=0}^{\infty}\left|\frac{a_{0}^{N+n}\kappa_{X}\left(\left[\mathfrak{z}\right]_{p^{N+n}}\right)}{a_{0}^{N}\kappa_{X}\left(\left[\mathfrak{z}\right]_{p^{N}}\right)}\right|_{\left(\mathcal{F}/I\right)\left(\mathfrak{z}\right)}\\
\left(\textrm{re-index}\right); & =\frac{\sum_{n=N}^{\infty}\left|a_{0}^{n}\kappa_{X}\left(\left[\mathfrak{z}\right]_{p^{n}}\right)\right|_{\left(\mathcal{F}/I\right)\left(\mathfrak{z}\right)}}{\left|a_{0}^{N}\kappa_{X}\left(\left[\mathfrak{z}\right]_{p^{N}}\right)\right|_{\left(\mathcal{F}/I\right)\left(\mathfrak{z}\right)}}
\end{align*}
By the hypothesis (\ref{eq:root test condition}), the \textbf{Root
Test} shows that the series in the numerator converges in $\mathbb{R}$.
Moreover, we get the estimate:
\begin{equation}
\sum_{n=N}^{\infty}\left|a_{0}^{n}\kappa_{X}\left(\left[\mathfrak{z}\right]_{p^{n}}\right)\right|_{\left(\mathcal{F}/I\right)\left(\mathfrak{z}\right)}\ll_{\mathfrak{z}}\left|a_{0}^{N}\kappa_{X}\left(\left[\mathfrak{z}\right]_{p^{N}}\right)\right|_{\left(\mathcal{F}/I\right)\left(\mathfrak{z}\right)}
\end{equation}
where, as indicated, the constant of proportionality depends on $\mathfrak{z}$,
though not on $N$. Thus:
\begin{equation}
\left(\mathcal{F}/I\right)_{\mathfrak{z}}\left(X\circ\theta_{p}^{\circ N}\right)\ll_{\mathfrak{z}}\frac{\sum_{n=N}^{\infty}\left|a_{0}^{n}\kappa_{X}\left(\left[\mathfrak{z}\right]_{p^{n}}\right)\right|_{\left(\mathcal{F}/I\right)\left(\mathfrak{z}\right)}}{\left|a_{0}^{N}\kappa_{X}\left(\left[\mathfrak{z}\right]_{p^{N}}\right)\right|_{\left(\mathcal{F}/I\right)\left(\mathfrak{z}\right)}}\ll_{\mathfrak{z}}\frac{\left|a_{0}^{N}\kappa_{X}\left(\left[\mathfrak{z}\right]_{p^{N}}\right)\right|_{\left(\mathcal{F}/I\right)\left(\mathfrak{z}\right)}}{\left|a_{0}^{N}\kappa_{X}\left(\left[\mathfrak{z}\right]_{p^{N}}\right)\right|_{\left(\mathcal{F}/I\right)\left(\mathfrak{z}\right)}}=1
\end{equation}
Taking the supremum over $N\geq0$ then gives (\ref{eq:shift sup boundedness}).

Finally, for quasi-integrability, like in \textbf{Proposition \ref{prop:X}},
without assuming the \emph{truth }of \textbf{Theorem \ref{thm:The-breakdown-variety}},
we can and very much will \emph{use} the formula it gives for $\hat{X}\left(t\right)$.
Applying \textbf{Proposition \ref{prop:X}} and \textbf{Proposition}
\textbf{\ref{prop:iterative triangle}}, we get:
\begin{equation}
\Delta_{N}^{\left(0\right)}\left\{ X\right\} \left(\mathfrak{z}\right)=\left(\Delta_{0}^{\left(N\right)}\left\{ X\right\} \left(\mathfrak{z}\right)+\left[\alpha_{X}\left(0\right)=1\right]N\beta_{X}\left(0\right)\right)a_{0}^{N}\kappa_{X}\left(\left[\mathfrak{z}\right]_{p^{N}}\right)
\end{equation}
where:
\begin{equation}
\Delta_{0}^{\left(N\right)}\left\{ X\right\} \left(\mathfrak{z}\right)=X\left(\theta_{p}^{\circ N}\left(\mathfrak{z}\right)\right)-\tilde{X}_{0}\left(\theta_{p}^{\circ N}\left(\mathfrak{z}\right)\right)=X\left(\theta_{p}^{\circ N}\left(\mathfrak{z}\right)\right)-\hat{X}\left(0\right)
\end{equation}
By (\ref{eq:shift sup boundedness}), this is bounded frame-wise uniformly
with respect to $N$:
\begin{equation}
\Delta_{0}^{\left(N\right)}\left\{ X\right\} \left(\mathfrak{z}\right)+\left[\alpha_{X}\left(0\right)=1\right]N\beta_{X}\left(0\right)\overset{\mathcal{F}/I}{\ll}N\textrm{ as }N\rightarrow\infty
\end{equation}
Hence:
\begin{equation}
\Delta_{N}^{\left(0\right)}\left\{ X\right\} \left(\mathfrak{z}\right)\overset{\mathcal{F}/I}{\ll}Na_{0}^{N}\kappa_{X}\left(\left[\mathfrak{z}\right]_{p^{N}}\right)
\end{equation}
By (\ref{eq:shift sup boundedness}), for each $\mathfrak{z}\in D\left(\mathcal{F}\right)$,
$\left(\mathcal{F}/I\right)_{\mathfrak{z}}\left(a_{0}^{N}\kappa_{X}\left(\left[\cdot\right]_{p^{N}}\right)\right)$
tends to $0$ at a rate of $O\left(r^{N}\right)$ for some $r\in\left(0,1\right)$.
This overwhelms $\left(\mathcal{F}/I\right)_{\mathfrak{z}}\left(N\right)$
(which is bounded if $\mathfrak{z}$ is non-archimedean and is at
most $N$ if $\mathfrak{z}$ is archimedean), and shows that:
\begin{equation}
\Delta_{N}^{\left(0\right)}\left\{ X\right\} \left(\mathfrak{z}\right)\overset{\mathcal{F}/I}{\ll}Na_{0}^{N}\kappa_{X}\left(\left[\mathfrak{z}\right]_{p^{N}}\right)\overset{\mathcal{F}/I}{\ll}a_{0}^{N}\kappa_{X}\left(\left[\mathfrak{z}\right]_{p^{N}}\right)
\end{equation}
which proves (\ref{eq:delta_N of X decay}). Taking the limit as $N\rightarrow\infty$,
we get $\lim_{N\rightarrow\infty}\Delta_{N}^{\left(0\right)}\left\{ X\right\} \left(\mathfrak{z}\right)\overset{\mathcal{F}/I}{=}0$,
which proves the quasi-integrability of $X$, as we have just shown
that the Fourier series generated by our choice for $\hat{X}$ converges
$\mathcal{F}/I$-wise to $X$.

Q.E.D.

\vphantom{}
\begin{prop}
\label{prop:dealing with g}Let $R$, $\mathcal{R}$, and an ideal
$I\subseteq R$ be as given in \textbf{Theorem \ref{thm:variety frame general}},
and let $\mathcal{F}$ be an evaluative, multiplicative, reffinite
$\mathcal{R}$-frame on $\mathbb{Z}_{p}$. If:
\begin{equation}
\limsup_{N\rightarrow\infty}\left(\left(\mathcal{F}/I\right)_{\mathfrak{z}}\left(r_{\mathbf{n},0}^{N}\kappa_{\mathbf{n}}\left(\left[\cdot\right]_{p^{N}}\right)\right)\right)^{1/N}<1,\textrm{ }\forall\mathfrak{z}\in D\left(\mathcal{F}\right)
\end{equation}
then $g_{\mathbf{n}}$ is $\mathcal{F}/I$-quasi-integrable and:
\begin{equation}
\max\left\{ N\left[\alpha_{\mathbf{n}}\left(0\right)=1\right],g_{\mathbf{n}}\left(\theta_{p}^{\circ N}\left(\mathfrak{z}\right)\right)\right\} r_{\mathbf{n},0}^{N}\kappa_{\mathbf{n}}\left(\left[\mathfrak{z}\right]_{p^{N}}\right)\overset{\mathcal{F}/I}{\ll}r_{\mathbf{n},0}^{N}\kappa_{\mathbf{n}}\left(\left[\mathfrak{z}\right]_{p^{N}}\right)\textrm{ as }N\rightarrow\infty\label{eq:decay of g and N}
\end{equation}
\end{prop}
Proof: Apply \textbf{Theorem \ref{thm:mini thesis}} to $g_{\mathbf{n}}$.

Q.E.D.

\vphantom{}

In our discussion of (\ref{eq:key bound, simplified}), the inductive
hypothesis requires us to show that the $\left(\mathcal{F},M\right)$-type
decay of $\Delta_{N}^{\left(0\right)}\left\{ X_{\mathbf{m}}\right\} $
for all $\mathbf{m}<\mathbf{n}$ guarantees the same type of decay
will occur for $\Delta_{N}^{\left(0\right)}\left\{ X_{\mathbf{n}}\right\} $.
Since we have already shown that the lower line of (\ref{eq:key bound, simplified})
has $\left(\mathcal{F},M\right)$-type decay, we can get the $\left(\mathcal{F},M\right)$-type
decay of $\Delta_{N}^{\left(0\right)}\left\{ X_{\mathbf{n}}\right\} $
by establishing the $\left(\mathcal{F},M\right)$-type decay of $\textrm{i}_{N}$
from (\ref{eq:key bound, simplified}) and then applying \textbf{Proposition
\ref{prop:sum of M-functions bound}} to get a single M-function as
an upper bound for the sum of two M-functions.

Thankfully, dealing with $\textrm{i}_{N}$ is pretty simple, primarily
because we've already done it. Indeed, choosing the term corresponding
to a single $\mathbf{m}$ in $\textrm{i}_{N}$, the term we get is
in precisely the same form as the function:
\begin{equation}
\textrm{II}_{N}\left(\mathfrak{z}\right)=\sum_{n=0}^{N-1}\rho\left(\left[\theta_{p}^{\circ n}\left(\mathfrak{z}\right)\right]_{p}\right)M_{n}\left(\mathcal{\mathfrak{z}}\right)\Delta_{N-n-1}^{\left(n+1\right)}\left\{ X\right\} \left(\mathfrak{z}\right)
\end{equation}
we worked with in our proof of \textbf{Lemma \ref{lem:main limit lemma}}.
Recall, there, we showed that:

\begin{align*}
\textrm{II}_{N}\left(\mathfrak{z}\right) & =\overbrace{\sum_{n=0}^{\gamma\left(N\right)-1}\rho\left(\left[\theta_{p}^{\circ n}\left(\mathfrak{z}\right)\right]_{p}\right)r^{n}\kappa\left(\left[\mathfrak{z}\right]_{p^{n}}\right)\Delta_{N-n-1}^{\left(n+1\right)}\left\{ X\right\} \left(\mathfrak{z}\right)}^{\textrm{Call this }\textrm{III}_{N}\left(\mathfrak{z}\right)}\\
 & +\underbrace{\sum_{n=\gamma\left(N\right)}^{N-1}\rho\left(\left[\theta_{p}^{\circ n}\left(\mathfrak{z}\right)\right]_{p}\right)r^{n}\kappa\left(\left[\mathfrak{z}\right]_{p^{n}}\right)\Delta_{N-n-1}^{\left(n+1\right)}\left\{ X\right\} \left(\mathfrak{z}\right)}_{\textrm{Call this }\textrm{IV}_{N}\left(\mathfrak{z}\right)}\\
 & \ll_{\mathfrak{z}}\frac{E_{\mathfrak{z},N}\left(\mathfrak{z}\right)}{E_{\mathfrak{z},\gamma\left(N\right)}\left(\mathfrak{z}\right)}+r^{\gamma\left(N\right)}\kappa\left(\left[\mathfrak{z}\right]_{p^{\gamma\left(N\right)}}\right)\\
 & =E_{\mathfrak{z},N-\gamma\left(N\right)}\left(\theta_{p}^{\circ\gamma\left(N\right)}\left(\mathfrak{z}\right)\right)+\underbrace{r^{\gamma\left(N\right)}\kappa\left(\left[\mathfrak{z}\right]_{p^{\gamma\left(N\right)}}\right)}_{\textrm{call this }F_{\mathfrak{z},\gamma\left(N\right)}\left(\mathfrak{z}\right)}\\
 & =E_{\mathfrak{z},N-\gamma\left(N\right)}\left(\theta_{p}^{\circ\gamma\left(N\right)}\left(\mathfrak{z}\right)\right)+F_{\mathfrak{z},\gamma\left(N\right)}\left(\mathfrak{z}\right)
\end{align*}
where $\gamma$ was any non-decreasing function $\mathbb{N}_{0}\rightarrow\mathbb{N}_{0}$
satisfying:
\begin{align}
\gamma\left(n\right) & \leq n,\textrm{ }\forall n\geq0\\
\lim_{n\rightarrow\infty}\gamma\left(n\right) & \overset{\mathbb{R}}{=}+\infty
\end{align}
In order to apply \textbf{Propositions \ref{prop:gamma prop}}, \textbf{\ref{prop:gamma theta prop}},
and \textbf{\ref{prop:sum of M-functions bound}}, we need for the
multipliers of $E_{\mathfrak{z},N-\gamma\left(N\right)}\left(\theta_{p}^{\circ\gamma\left(N\right)}\left(\mathfrak{z}\right)\right)$
and $F_{\mathfrak{z},\gamma\left(N\right)}\left(\mathfrak{z}\right)$
to satisfy the properties required by those propositions' hypotheses.
Helpfully, those propositions also show that if said properties are
satisfied by the M-functions that we feed into them, then the M-functions
those propositions give us as our upper bounds will also satisfy those
properties. This is the key to making our inductive step work; not
only does the quasi-integrability of $X_{\mathbf{m}}$ for $\mathbf{m}<\mathbf{n}$
imply the quasi-integrability of $X_{\mathbf{n}}$, the decay properties
of the M-function bounds on $\Delta_{N}^{\left(0\right)}\left\{ X_{\mathbf{m}}\right\} $
for $\mathbf{m}<\mathbf{n}$ are also passed on to $\Delta_{N}^{\left(0\right)}\left\{ X_{\mathbf{n}}\right\} $.
The trick that makes this all work is restricting our frame's domain
to $\textrm{NiceDig}\left(\mathbb{Z}_{p}\right)$, so that we can
apply \textbf{Propositions \ref{prop:gamma prop}}, \textbf{\ref{prop:gamma theta prop}}.
\begin{notation}
Given an $\mathcal{R}$-frame $\mathcal{F}$, let $\left(r,\kappa\right)$
be an $\mathcal{R}$-valued $p$-adic M-function with multipliers
$r_{0},\ldots,r_{p-1}$, with $r_{0}=r$. Then, for any $\mathfrak{z}\in D\left(\mathcal{F}\right)$,
we write $\mathcal{F}_{\mathfrak{z}}\left(r,\kappa\right)$ to denote
the real-valued $p$-adic M-function whose multipliers are $\mathcal{F}_{\mathfrak{z}}\left(r_{0}\right),\ldots,\mathcal{F}_{\mathfrak{z}}\left(r_{p-1}\right)$.
\end{notation}
Here, then, is the inductive step:
\begin{thm}[The Inductive Step]
\label{thm:inductive step}Let $R$, $\mathcal{R}$, and an ideal
$I\subseteq R$ be as given in \textbf{Theorem \ref{thm:variety frame general}},
and let $\mathcal{F}$ be an evaluative, multiplicative, reffinite
$\mathcal{R}$-frame on $\mathbb{Z}_{p}$ with $D\left(\mathcal{F}\right)\subseteq\textrm{NiceDig}\left(\mathbb{Z}_{p}\right)$.
Consider the hypotheses:

I. 
\begin{equation}
\lim_{N\rightarrow\infty}\left(\left(\mathcal{F}/I\right)_{\mathfrak{z}}\left(r_{\mathbf{n},0}^{N}\kappa_{\mathbf{n}}\left(\left[\cdot\right]_{p^{N}}\right)\right)\right)^{1/N}<1,\textrm{ }\forall\mathfrak{z}\in D\left(\mathcal{F}\right)
\end{equation}

\vphantom{}

II. 
\begin{equation}
\max_{j\in\mathbf{n}}\sup_{n\geq0}\left(\mathcal{F}/I\right)_{\mathfrak{z}}\left(X_{j}\circ\theta_{p}^{\circ n}\right)<\infty,\textrm{ }\forall\mathfrak{z}\in D\left(\mathcal{F}\right)
\end{equation}

\vphantom{}

III. For each $\mathbf{m}$ with $\mathbf{0}<\mathbf{m}<\mathbf{n}$,
$X_{\mathbf{m}}$ is $\mathcal{F}/I$-quasi-integrable with a Fourier
transform $\hat{X}_{\mathbf{m}}$ so that $\Delta_{N}^{\left(0\right)}\left\{ X_{\mathbf{m}}\right\} $
has $\left(\mathcal{F}/I,M\right)$-type decay, with:
\begin{equation}
\Delta_{N}^{\left(0\right)}\left\{ X_{\mathbf{m}}\right\} \overset{\mathcal{F}/I}{\ll_{\mathbf{m}}}E_{\mathbf{m},\mathfrak{z},n}\left(\mathfrak{z}\right)\label{eq:X_J delta bound}
\end{equation}
where, as indicated, the constant of proportionality depends on $\mathbf{m}$,
and where, for each $\mathfrak{z}\in D\left(\mathcal{F}\right)$,
$\left\{ E_{\mathbf{m},\mathfrak{z},n}\right\} _{n\geq0}\in\textrm{MFunc}\left(\mathbb{Z}_{p},\mathbb{R}\right)$
so that $\lim_{n\rightarrow\infty}E_{\mathbf{m},\mathfrak{z},n}^{1/n}\left(\mathfrak{z}\right)\overset{\mathbb{R}}{<}1$
for all non-empty $\mathbf{m}<\mathbf{n}$.

\vphantom{}If (I), (II), and (III) are satisfied, then $X_{\mathbf{n}}\left(\mathfrak{z}\right)$
is $\mathcal{F}/I$-quasi-integrable with a Fourier transform $\hat{X}_{\mathbf{n}}$
so that $\Delta_{N}^{\left(0\right)}\left\{ X_{\mathbf{n}}\right\} $
has $\left(\mathcal{F}/I,M\right)$-type decay, with:
\begin{equation}
\Delta_{N}^{\left(0\right)}\left\{ X_{\mathbf{n}}\right\} \overset{\mathcal{F}/I}{\ll}E_{\mathbf{n},\mathfrak{z},n}\left(\mathfrak{z}\right)
\end{equation}
where $\left\{ E_{\mathbf{n},\mathfrak{z},n}\right\} _{n\geq0}\in\textrm{MFunc}\left(\mathbb{Z}_{p},\mathbb{R}\right)$
satisfies $\lim_{n\rightarrow\infty}E_{\mathbf{n},\mathfrak{z},n}^{1/n}\left(\mathfrak{z}\right)\overset{\mathbb{R}}{<}1$.
\end{thm}
Proof: Let everything be as given and fix $\mathfrak{z}\in D\left(\mathcal{F}\right)$.
Because the arguments apply verbatim when considering descent through
the quotient $\mathcal{F}\rightarrow\mathcal{F}/I$, we will give
the proof writing $\mathcal{F}$ rather than $\mathcal{F}/I$.

By \textbf{Proposition \ref{prop:dealing with g}}, hypothesis (I)
guarantees that the lower line of the key estimate (\ref{eq:key bound, simplified})
tends to $0$ frame-wise as $N\rightarrow\infty$. Thus, all that
remains to be done is to show that the top line of that estimate:
\begin{equation}
\textrm{i}_{N}\left(\mathfrak{z}\right)=\sum_{n=0}^{N-2}r_{\mathbf{n},0}^{n}\kappa_{\mathbf{n}}\left(\left[\mathfrak{z}\right]_{p^{n}}\right)\sum_{\mathbf{m}<\mathbf{n}}r_{\mathbf{m},\left[\theta_{p}^{\circ n}\left(\mathfrak{z}\right)\right]_{p}}\Delta_{N-n-1}^{\left(n+1\right)}\left\{ X_{\mathbf{m}}\right\} \left(\mathfrak{z}\right)
\end{equation}
tends to $0$ with respect to $\mathcal{F}$ as $N\rightarrow\infty$.

For brevity, let $\left(\rho_{\mathfrak{z}},\kappa_{\mathfrak{z}}\right)$
denote the M-function $\mathcal{F}_{\mathfrak{z}}\left(r_{\mathbf{n},0},\kappa_{\mathbf{n}}\right)$.
To complete the proof, we just need to show that $\textrm{i}_{N}\left(\mathfrak{z}\right)$
can be bounded frame-wise by a real-valued M-function satisfying the
root test condition (\ref{eq:decay root test condition}).

Proceeding like in \textbf{Lemma \ref{lem:main limit lemma}}, we
split $i_{N}$ at $n=\gamma\left(N\right)-1$:

\begin{align*}
\textrm{i}_{N}\left(\mathfrak{z}\right) & =\sum_{n=0}^{\gamma\left(N\right)-1}r_{\mathbf{n},0}^{n}\kappa_{\mathbf{n}}\left(\left[\mathfrak{z}\right]_{p^{n}}\right)\sum_{\mathbf{m}<\mathbf{n}}r_{\mathbf{m},\left[\theta_{p}^{\circ n}\left(\mathfrak{z}\right)\right]_{p}}\Delta_{N-n-1}^{\left(n+1\right)}\left\{ X_{\mathbf{m}}\right\} \left(\mathfrak{z}\right)\\
 & +\sum_{n=\gamma\left(N\right)}^{N-2}r_{\mathbf{n},0}^{n}\kappa_{\mathbf{n}}\left(\left[\mathfrak{z}\right]_{p^{n}}\right)\sum_{\mathbf{m}<\mathbf{n}}r_{\mathbf{m},\left[\theta_{p}^{\circ n}\left(\mathfrak{z}\right)\right]_{p}}\Delta_{N-n-1}^{\left(n+1\right)}\left\{ X_{\mathbf{m}}\right\} \left(\mathfrak{z}\right)
\end{align*}
Using (\ref{eq:X_J delta bound}) and noting the uniform boundedness
of the $r_{\mathbf{m},\left[\theta_{p}^{\circ n}\left(\mathfrak{z}\right)\right]_{p}}$s
with respect to $\mathbf{m}<\mathbf{n}$ and $n$, we have:
\begin{equation}
\mathcal{F}_{\mathfrak{z}}\left(\max_{n\in S}\sum_{\mathbf{m}<\mathbf{n}}r_{\mathbf{m},\left[\theta_{p}^{\circ n}\left(\cdot\right)\right]_{p}}\Delta_{N-n-1}^{\left(n+1\right)}\left\{ X_{\mathbf{m}}\right\} \right)\ll_{\mathfrak{z},\mathbf{n}}\max_{n\in S}\max_{\mathbf{m}<\mathbf{n}}E_{\mathbf{m},\mathfrak{z},N-n-1}\left(\theta_{p}^{\circ n+1}\left(\mathfrak{z}\right)\right)
\end{equation}
for any set $S$. By \textbf{Proposition \ref{prop:max of two M functions}},
for each $\mathfrak{z}$, we can choose an M-function $E_{\mathfrak{z},n}$
so that:
\begin{equation}
\max_{\mathbf{m}<\mathbf{n}}E_{\mathbf{m},\mathfrak{z},n}\left(\mathfrak{y}\right)\leq E_{\mathfrak{z},n}\left(\mathfrak{y}\right),\textrm{ }\forall\mathfrak{y}\in\mathbb{Z}_{p},\textrm{ }\forall n\geq0
\end{equation}
and so that $\lim_{n\rightarrow\infty}E_{\mathfrak{z},n}^{1/n}\left(\mathfrak{z}\right)<1$.
Then:
\begin{align*}
\mathcal{F}_{\mathfrak{z}}\left(\max_{n\in S}\sum_{\mathbf{m}<\mathbf{n}}r_{\mathbf{m},\left[\theta_{p}^{\circ n}\left(\cdot\right)\right]_{p}}\Delta_{N-n-1}^{\left(n+1\right)}\left\{ X_{\mathbf{m}}\right\} \right) & \ll_{\mathfrak{z},\mathbf{n}}\max_{n\in S}E_{\mathbf{m},\mathfrak{z},N-n-1}\left(\theta_{p}^{\circ n+1}\left(\mathfrak{z}\right)\right)\\
\left(\mathbf{Proposition}\textrm{ }\mathbf{\ref{prop:Kappa shift equation}}\right); & =\max_{n\in S}\frac{E_{\mathfrak{z},N}\left(\mathfrak{z}\right)}{E_{\mathfrak{z},n+1}\left(\mathfrak{z}\right)}
\end{align*}
Thus:
\begin{align*}
\mathcal{F}_{\mathfrak{z}}\left(\textrm{i}_{N}\right) & \overset{\mathcal{F}}{\ll}_{\mathfrak{z},\mathbf{n}}\overbrace{\sum_{n=0}^{\gamma\left(N\right)-1}\rho_{\mathfrak{z}}^{n}\kappa_{\mathfrak{z}}\left(\left[\mathfrak{z}\right]_{p^{n}}\right)\max_{0\leq n\leq\gamma\left(N\right)-1}\frac{E_{\mathfrak{z},N}\left(\mathfrak{z}\right)}{E_{\mathfrak{z},n+1}\left(\mathfrak{z}\right)}}^{\textrm{Call this }\textrm{III}_{N}\left(\mathfrak{z}\right)}\\
 & +\underbrace{\sum_{n=\gamma\left(N\right)}^{N-2}\rho_{\mathfrak{z}}^{n}\kappa_{\mathfrak{z}}\left(\left[\mathfrak{z}\right]_{p^{n}}\right)\max_{\mathbf{m}<\mathbf{n}}\max_{\gamma\left(N\right)\leq n\leq N-2}\frac{E_{\mathfrak{z},N}\left(\mathfrak{z}\right)}{E_{\mathfrak{z},n+1}\left(\mathfrak{z}\right)}}_{\textrm{Call this }\textrm{IV}_{N}\left(\mathfrak{z}\right)}
\end{align*}
Since $\left(r_{\mathbf{n},0},\kappa_{\mathbf{n}}\right)$ is strongly
$\mathcal{F}$-summable, we have that $\rho_{\mathfrak{z}}^{n}\kappa_{\mathfrak{z}}\left(\left[\mathfrak{z}\right]_{p^{n}}\right)$
will be small with respect to $\mathcal{F}$ when $n$ is large. Note
that as $n\rightarrow\infty$, $1/E_{\mathfrak{z},n+1}\left(\mathfrak{z}\right)$
will be maximized when $n$ is large, and so, we get the bounds:
\begin{align}
\mathcal{F}_{\mathfrak{z}}\left(\textrm{III}_{N}\right) & \overset{\mathcal{F}}{\ll}_{\mathfrak{z},\mathbf{n}}E_{\mathfrak{z},N-\gamma\left(N\right)}\left(\theta_{p}^{\circ\gamma\left(N\right)}\left(\mathfrak{z}\right)\right)\\
\mathcal{F}_{\mathfrak{z}}\left(\textrm{IV}_{N}\right) & \overset{\mathcal{F}}{\ll}_{\mathfrak{z},\mathbf{n}}\rho_{\mathfrak{z}}^{\gamma\left(N\right)}\kappa\left(\left[\mathfrak{z}\right]_{p^{\gamma\left(N\right)}}\right)\frac{E_{\mathfrak{z},N}\left(\mathfrak{z}\right)}{E_{\mathfrak{z},N-1}\left(\mathfrak{z}\right)}\overset{\mathcal{F}}{\ll}_{\mathfrak{z},\mathbf{n}}\rho_{\mathfrak{z}}^{\gamma\left(N\right)}\kappa\left(\left[\mathfrak{z}\right]_{p^{\gamma\left(N\right)}}\right)
\end{align}
where we have used the fact that $\frac{E_{\mathfrak{z},N}\left(\mathfrak{z}\right)}{E_{\mathfrak{z},N-1}\left(\mathfrak{z}\right)}$
is going to be equal to one of the multipliers of $\left\{ E_{\mathfrak{z},n}\right\} _{n\geq0}$,
and hence, that it will be bounded independently of $N$ and $\mathfrak{z}$.

Writing:
\begin{equation}
F_{\mathfrak{z},n}\left(\mathfrak{z}\right)\overset{\textrm{def}}{=}\rho_{\mathfrak{z}}^{n}\kappa_{\mathfrak{z}}\left(\left[\mathfrak{z}\right]_{p^{n}}\right)
\end{equation}
we have that $\lim_{n\rightarrow\infty}F_{\mathfrak{z},n}^{1/n}\left(\mathfrak{z}\right)<1$.
$\textrm{i}_{N}$ then satisfies the estimate:
\begin{equation}
\mathcal{F}_{\mathfrak{z}}\left(\textrm{i}_{N}\right)\ll_{\mathfrak{z},\mathbf{n}}E_{\mathfrak{z},N-\gamma\left(N\right)}\left(\theta_{p}^{\circ\gamma\left(N\right)}\left(\mathfrak{z}\right)\right)+F_{\mathfrak{z},\gamma\left(N\right)}\left(\mathfrak{z}\right)\label{eq:III.1}
\end{equation}

Thus, the desired frame-wise decay of $\textrm{i}_{N}$ will follow
once we apply \textbf{Propositions \ref{prop:gamma prop}},\textbf{
\ref{prop:gamma theta prop}}, and \textbf{\ref{prop:sum of M-functions bound}}
to (\ref{eq:III.1}). For this, we just need to be able to choose
a $\gamma$ so that the hypotheses of those propositions are satisfied.
This is where the niceness of $\mathfrak{z}$'s digits plays a key
role. We have three cases.

\vphantom{}\textbf{Case 1}: Suppose $\mathfrak{z}\in\mathbb{N}_{0}$.
Since $E_{\mathfrak{z},N}$ was constructed so as to guarantee that
$\lim_{N\rightarrow\infty}E_{\mathfrak{z},N}^{1/N}\left(\mathfrak{z}\right)<1$
we have:
\begin{equation}
1>\lim_{N\rightarrow\infty}E_{\mathfrak{z},N}^{1/N}\left(\mathfrak{z}\right)=e_{\mathfrak{z},0}\lim_{N\rightarrow\infty}\prod_{\ell=1}^{p-1}\left(\frac{e_{\mathfrak{z},\ell}}{e_{\mathfrak{z},0}}\right)^{\#_{p:\ell}\left(\left[\mathfrak{z}\right]_{p^{N}}\right)/N}=e_{\mathfrak{z},0}
\end{equation}
Consequently, for all sufficiently large $N$:
\begin{equation}
E_{\mathfrak{z},N}\left(\mathfrak{z}\right)=e_{\mathfrak{z},0}^{N}\prod_{\ell=1}^{p-1}\left(\frac{e_{\mathfrak{z},\ell}}{e_{\mathfrak{z},0}}\right)^{\#_{p:\ell}\left(\mathfrak{z}\right)}
\end{equation}
Next, since $\left(r_{\mathbf{n},0},\kappa_{\mathbf{n}}\right)$ is
strongly $\mathcal{F}$-summable, we have that $\left\{ F_{\mathfrak{z},n}\right\} _{n\geq0}=\left(\rho_{\mathfrak{z}},\kappa_{\mathfrak{z}}\right)=\mathcal{F}_{\mathfrak{z}}\left(r_{\mathcal{\mathbf{n}},0},\kappa_{\mathbf{n}}\right)$
satisfies:
\begin{equation}
\lim_{n\rightarrow\infty}\rho_{\mathfrak{z}}^{n}\kappa_{\mathfrak{z}}\left(\left[\mathfrak{z}\right]_{p^{n}}\right)\overset{\mathbb{R}}{=}0
\end{equation}
Letting $\rho_{\mathfrak{z},0},\ldots,\rho_{\mathfrak{z},p-1}$ be
the multipliers of $\left(\rho_{\mathfrak{z}},\kappa_{\mathfrak{z}}\right)$,
the same argument given for $E_{\mathfrak{z},n}$ above shows that
$\rho_{\mathfrak{z},k}<1$.

Thus, for all sufficiently large $N$, the upper bound (\ref{eq:III.1})
then becomes:
\begin{align*}
\mathcal{F}_{\mathfrak{z}}\left(\textrm{i}_{N}\right) & \ll_{\mathfrak{z},\mathbf{n}}E_{\mathfrak{z},N-\gamma\left(N\right)}\left(\theta_{p}^{\circ\gamma\left(N\right)}\left(\mathfrak{z}\right)\right)+F_{\mathfrak{z},\gamma\left(N\right)}\left(\mathfrak{z}\right)\\
 & =\frac{E_{\mathfrak{z},N}\left(\mathfrak{z}\right)}{E_{\mathfrak{z},\gamma\left(N\right)}\left(\mathfrak{z}\right)}+F_{\mathfrak{z},\gamma\left(N\right)}\left(\mathfrak{z}\right)\\
 & \ll e_{\mathfrak{z},k}^{N-\gamma\left(N\right)}+\rho_{\mathfrak{z},k}^{\gamma\left(N\right)}
\end{align*}
So, choosing $\gamma\left(N\right)\overset{\textrm{def}}{=}\left\lfloor N/2\right\rfloor $,
we have that:
\begin{align}
\mathcal{F}_{\mathfrak{z}}\left(\textrm{i}_{N}\right) & \ll_{\mathfrak{z},\mathbf{n}}\left(\max\left\{ \sqrt{e_{\mathfrak{z},k}},\sqrt{\rho_{\mathfrak{z},k}}\right\} \right)^{N}
\end{align}
which shows that at such a $\mathfrak{z}$, $\mathcal{F}_{\mathfrak{z}}\left(\textrm{i}_{N}\right)$
is bounded by the M-function with $\max\left\{ \sqrt{e_{\mathfrak{z},k}},\sqrt{\rho_{\mathfrak{z},k}}\right\} $
as its zero multiplier and with all other multipliers being $1$.

\vphantom{}\textbf{Case 2}:\textbf{ }Suppose $\mathfrak{z}\in\textrm{NiceDig}\left(\mathbb{Z}_{p}\right)\backslash\mathbb{N}_{0}$.
Fix $k\in\left\{ 0,\ldots,p-1\right\} $, and suppose that $d_{p:k}\left(\mathfrak{z}\right)=1$.
Then $d_{p:\ell}\left(\mathfrak{z}\right)=0$ for all $\ell\neq k$.
Since $E_{\mathfrak{z},N}$ was constructed so as to guarantee that
$\lim_{N\rightarrow\infty}E_{\mathfrak{z},N}^{1/N}\left(\mathfrak{z}\right)<1$
we have:
\begin{equation}
1>\lim_{N\rightarrow\infty}E_{\mathfrak{z},N}^{1/N}\left(\mathfrak{z}\right)=e_{\mathfrak{z},0}\prod_{\ell=1}^{p-1}\left(\frac{e_{\mathfrak{z},\ell}}{e_{\mathfrak{z},0}}\right)^{d_{p:\ell}\left(\mathfrak{z}\right)}=e_{\mathfrak{z},0}\frac{e_{\mathfrak{z},k}}{e_{\mathfrak{z},0}}=e_{\mathfrak{z},k}
\end{equation}
Consequently, there exist real constants $C_{\mathfrak{z}},C_{\mathfrak{z}}^{\prime}>0$
so that:
\begin{equation}
C_{\mathfrak{z}}^{\prime}e_{\mathfrak{z},k}^{N}\leq E_{\mathfrak{z},N}\left(\mathfrak{z}\right)\leq C_{\mathfrak{z}}e_{\mathfrak{z},k}^{N}
\end{equation}
for all sufficiently large $N$.

Next, since $\left(r_{\mathbf{n},0},\kappa_{\mathbf{n}}\right)$ is
strongly $\mathcal{F}$-summable, we have that $\left\{ F_{\mathfrak{z},n}\right\} _{n\geq0}=\left(\rho_{\mathfrak{z}},\kappa_{\mathfrak{z}}\right)=\mathcal{F}_{\mathfrak{z}}\left(r_{\mathbf{n},0},\kappa_{\mathbf{n}}\right)$
satisfies:
\begin{equation}
\lim_{n\rightarrow\infty}\rho_{\mathfrak{z}}^{n}\kappa_{\mathfrak{z}}\left(\left[\mathfrak{z}\right]_{p^{n}}\right)\overset{\mathbb{R}}{=}0
\end{equation}
Letting $\rho_{\mathfrak{z},0},\ldots,\rho_{\mathfrak{z},p-1}$ be
the multipliers of $\left(\rho_{\mathfrak{z}},\kappa_{\mathfrak{z}}\right)$,
the same argument given for $E_{\mathfrak{z},n}$ above shows that
$\rho_{\mathfrak{z},k}<1$.

Thus, for all sufficiently large $N$, the upper bound (\ref{eq:III.1})
then becomes:
\begin{align*}
\mathcal{F}_{\mathfrak{z}}\left(\textrm{i}_{N}\right) & \ll_{\mathfrak{z},\mathbf{n}}E_{\mathfrak{z},N-\gamma\left(N\right)}\left(\theta_{p}^{\circ\gamma\left(N\right)}\left(\mathfrak{z}\right)\right)+F_{\mathfrak{z},\gamma\left(N\right)}\left(\mathfrak{z}\right)\\
 & =\frac{E_{\mathfrak{z},N}\left(\mathfrak{z}\right)}{E_{\mathfrak{z},\gamma\left(N\right)}\left(\mathfrak{z}\right)}+F_{\mathfrak{z},\gamma\left(N\right)}\left(\mathfrak{z}\right)\\
 & \ll e_{\mathfrak{z},k}^{N-\gamma\left(N\right)}+\rho_{\mathfrak{z},k}^{\gamma\left(N\right)}
\end{align*}
again, choosing $\gamma\left(N\right)\overset{\textrm{def}}{=}\left\lfloor N/2\right\rfloor $,
we have that:
\begin{align}
\mathcal{F}_{\mathfrak{z}}\left(\textrm{i}_{N}\right) & \ll_{\mathfrak{z},\mathbf{n}}\left(\max\left\{ \sqrt{e_{\mathfrak{z},k}},\sqrt{\rho_{\mathfrak{z},k}}\right\} \right)^{N}
\end{align}
which shows that at such a $\mathfrak{z}$, $\mathcal{F}_{\mathfrak{z}}\left(\textrm{i}_{N}\right)$
is bounded by the M-function with $\max\left\{ \sqrt{e_{\mathfrak{z},k}},\sqrt{\rho_{\mathfrak{z},k}}\right\} $
as its zero multiplier and with all other multipliers being $1$.

\vphantom{}

\textbf{Case 3}: Suppose $\mathfrak{z}\in\textrm{NiceDig}\left(\mathbb{Z}_{p}\right)\backslash\mathbb{N}_{0}$
and that there is no $k\in\left\{ 0,\ldots,p-1\right\} $ so that
all but finitely many of $\mathfrak{z}$'s digits are equal to $k$.
Then, since $\mathfrak{z}$ has nice digits, $d_{p:k}\left(\mathfrak{z}\right)$
exists for all $k$ and is not equal to $1$ for any $k$. Thus, for
our situation, let $\delta\in\left(0,1\right)$ be slightly larger
than $\max_{0\leq k\leq p-1}d_{p:k}\left(\mathfrak{z}\right)$, and
set $\gamma\left(N\right)=\left\lfloor \delta N\right\rfloor $. Then,
$\gamma\left(N\right)$, $E_{\mathfrak{z},N-\gamma\left(N\right)}\left(\theta_{p}^{\circ\gamma\left(N\right)}\left(\mathfrak{z}\right)\right)$,
and $F_{\mathfrak{z},\gamma\left(N\right)}\left(\mathfrak{z}\right)$
satisfy the hypotheses of \textbf{Propositions \ref{prop:gamma prop}},
\textbf{\ref{prop:gamma theta prop}}, and \textbf{\ref{prop:sum of M-functions bound}}.

\vphantom{}

Thus, in either case, (\ref{eq:III.1}) is bounded from above by an
M-function satisfying the necessary decay conditions, and we are done.

Q.E.D.

\vphantom{}

We now have everything we need in order to prove our main result.
This will be done in two steps. First, we will deal with $\mathcal{R}_{1}$;
this entails showing that for a single F-series $X$, the Fourier
transform ends up inducing an isomorphism out of the polynomial algebra
$K\left[X\left(\mathfrak{z}\right)d\mathfrak{z}\right]$ generated
over $K$ by the measure $X\left(\mathfrak{z}\right)d\mathfrak{z}$
with pointwise multiplication. For this case\textemdash which, for
obvious reasons, I call the \textbf{Power Case}\textemdash we will
use \textbf{Theorem} \textbf{\ref{thm:mini thesis}} to establish
the quasi-integrability of $X$ as the base case, and then apply \textbf{Theorem
\ref{thm:inductive step}} for the inductive step. Once that is done,
we will be able to deal with the general case ($\mathcal{R}_{d}$,
for $d\geq2$) like so: given a product $X_{1}^{e_{1}}\cdots X_{d}^{e_{d}}$
of $d$ distinct F-series where the $e_{j}$s are positive integer
exponents, we will set $Y_{j}=X_{j}^{e_{j}}$ for all $j\in\left\{ 1,\ldots,d\right\} $
and then apply \textbf{Theorem \ref{thm:inductive step}}, using the
\textbf{Power Case }as the base case to establish the quasi-integrability
of each $Y_{j}$.
\begin{thm}[The Power Case]
\label{thm:powers of X}Let $p\geq2$, let $K$ be a global field
(if $\textrm{char}K>0$, then assume $\textrm{char}K$ is co-prime
to $p$). For each $k\in\left\{ 0,\ldots,p-1\right\} $, let $a_{k}$
and $b_{k}$ be indeterminates, and let $R_{1}\left(K\right)$ (a.k.a.,
$R_{1}$, $R$) be the ring of polynomials in the $a_{k}$s and $b_{k}$s
with coefficients in $\mathcal{O}_{K}$, so that $R$ is then a free
$\mathcal{O}_{K}$-algebra in $2p$ indeterminates. Then, let:
\begin{equation}
\mathcal{R}_{1}\left(K\right)\overset{\textrm{def}}{=}\frac{R\left[x\right]}{\left\langle \left(1-a_{0}\right)x-1\right\rangle }
\end{equation}
(a.k.a., $\mathcal{R}_{1}$, $\mathcal{R}$) be the localization of
$R$ away from the prime ideal $\left\langle 1-a_{0}\right\rangle $.

Now, let $I$ be a unique solution ideal of $R$ ($\left\langle 1-a_{0}\right\rangle \nsubseteq I$),
and let $\mathcal{F}$ be any multiplicative, evaluative, reffinite
$\mathcal{R}$-frame on $\mathbb{Z}_{p}$ with $D\left(\mathcal{F}\right)\subseteq\textrm{NiceDig}\left(\mathbb{Z}_{p}\right)$.
If:
\begin{equation}
\lim_{n\rightarrow\infty}\left|a_{0}^{n}\kappa_{1}\left(\left[\mathfrak{z}\right]_{p^{n}}\right)\right|_{\mathcal{F}\left(\mathfrak{z}\right)}^{1/n}<1,\textrm{ }\forall\mathfrak{z}\in D\left(\mathcal{F}\right)\label{eq:root test condition, power case}
\end{equation}
where:
\begin{equation}
\kappa_{1}\left(m\right)=\prod_{k=1}^{p-1}\left(\frac{a_{k}}{a_{0}}\right)^{\#_{p:k}\left(m\right)}
\end{equation}
then:

\vphantom{}I. There exists a unique rising-continuous function $X\in C\left(\mathcal{F}/I\right)$
so that:
\begin{equation}
X\left(\mathfrak{z}\right)\overset{\left(\mathcal{F}/I\right)\left(\mathfrak{z}\right)}{=}a_{\left[\mathfrak{z}\right]_{p}}X\left(\theta_{p}\left(\mathfrak{z}\right)\right)+b_{\left[\mathfrak{z}\right]_{p}},\textrm{ }\forall\mathfrak{z}\in D\left(\mathcal{F}/I\right)\label{eq:functional equations-1}
\end{equation}
where $\left(\mathcal{F}/I\right)\left(\mathfrak{z}\right)$ is the
completion of $\mathcal{R}/I\mathcal{R}$ with respect to the absolute
value $\mathcal{F}/I$ associates to $\mathfrak{z}$.

\vphantom{}II. For any integer $n\geq0$, the function $X^{n}\in C\left(\mathcal{F}/I\right)$
with:
\begin{equation}
\mathfrak{z}\in\mathbb{Z}_{p}\mapsto X^{n}\left(\mathfrak{z}\right)\in\left(\mathcal{F}/I\right)\left(\mathfrak{z}\right)
\end{equation}
is $\mathcal{F}/I$-quasi-integrable, where $X^{n_{}}\left(\mathfrak{z}\right)$
is the pointwise product of $n$ copies of $X$, with $X^{n}$ being
$1$ identically when $n=0$.

\vphantom{}III. For each integer $n\geq0$, define $\hat{X}^{*n}:\hat{\mathbb{Z}}_{p}\rightarrow\textrm{Frac}\left(\mathcal{R}/I\mathcal{R}\right)\left(\zeta_{p^{\infty}}\right)$
as follows. Let $\hat{X}^{*0}\left(t\right)$ denote the function
$\mathbf{1}_{0}\left(t\right)$, and let $\hat{X}^{*1}\left(t\right)$
denote $\hat{X}\left(t\right)$, given by:
\begin{equation}
\hat{X}\left(t\right)=\begin{cases}
0 & \textrm{if }t=0\\
\left(\beta_{1}\left(0\right)v_{p}\left(t\right)+\gamma_{1}\left(\frac{t\left|t\right|_{p}}{p}\right)\right)\hat{A}_{1}\left(t\right) & \textrm{if }t\neq0
\end{cases},\textrm{ }\forall t\in\hat{\mathbb{Z}}_{p}
\end{equation}
if $\alpha_{1}\left(0\right)=1$ occurs in $\mathcal{R}/I\mathcal{R}$,
and given by:
\begin{equation}
\hat{X}\left(t\right)=\begin{cases}
\frac{\beta_{1}\left(0\right)}{1-\alpha_{1}\left(0\right)} & \textrm{if }t=0\\
\left(\frac{\beta_{1}\left(0\right)}{1-\alpha_{1}\left(0\right)}+\gamma_{1}\left(\frac{t\left|t\right|_{p}}{p}\right)\right)\hat{A}_{1}\left(t\right) & \textrm{if }t\neq0
\end{cases},\textrm{ }\forall t\in\hat{\mathbb{Z}}_{p}
\end{equation}
when $\alpha_{1}\left(0\right)\neq1$ in $\mathcal{R}/I\mathcal{R}$.

Then, by recursion, for each $n\geq2$, let $\hat{X}^{*n}\left(t\right)$
denote $\hat{f}_{n}\left(t\right)-\hat{g}_{n}\left(t\right)$, where:
\begin{equation}
\hat{f}_{n}\left(t\right)\overset{\textrm{def}}{=}\sum_{k=0}^{n-1}\binom{n}{k}\sum_{m=0}^{-v_{p}\left(t\right)-2}\hat{X}^{*k}\left(p^{m+1}t\right)\alpha_{k,n}\left(p^{m}t\right)\prod_{h=0}^{m-1}\alpha_{n}\left(p^{h}t\right)\label{eq:f-hat, power case}
\end{equation}
where:
\begin{equation}
\alpha_{k,n}\left(t\right)=\frac{1}{p}\sum_{j=0}^{p-1}a_{j}^{k}b_{j}^{n-k}e^{-2\pi ijt}
\end{equation}
 and:

\begin{equation}
\hat{g}_{n}\left(t\right)=\begin{cases}
0 & \textrm{if }t=0\\
\left(\beta_{n}\left(0\right)v_{p}\left(t\right)+\gamma_{n}\left(\frac{t\left|t\right|_{p}}{p}\right)\right)\hat{A}_{n}\left(t\right) & \textrm{if }t\neq0
\end{cases},\textrm{ }\forall t\in\hat{\mathbb{Z}}_{p}\label{eq:g_script J hat , alpha equals 1-1-1}
\end{equation}
when $\alpha_{n}\left(0\right)=1$ and:
\begin{equation}
\hat{g}_{n}\left(t\right)=\begin{cases}
\frac{\beta_{n}\left(0\right)}{1-\alpha_{n}\left(0\right)} & \textrm{if }t=0\\
\left(\frac{\beta_{n}\left(0\right)}{1-\alpha_{n}\left(0\right)}+\gamma_{n}\left(\frac{t\left|t\right|_{p}}{p}\right)\right)\hat{A}_{n}\left(t\right) & \textrm{if }t\neq0
\end{cases},\textrm{ }\forall t\in\hat{\mathbb{Z}}_{p}\label{eq:g_script J hat , alpha not equal to 1-1-1}
\end{equation}
when $\alpha_{n}\left(0\right)\neq1$, where, for all $n\in\mathbb{N}_{0}$,
$\alpha_{n},\beta_{n},\gamma_{n},\hat{A}_{n}:\hat{\mathbb{Z}}_{p}\rightarrow\textrm{Frac}\left(\mathcal{R}/I\mathcal{R}\right)\left(\zeta_{p^{\infty}}\right)$
are given by:
\begin{equation}
\alpha_{n}\left(t\right)\overset{\textrm{def}}{=}\frac{1}{p}\sum_{j=0}^{p-1}a_{j}^{n}e^{-2\pi ijt}
\end{equation}
\begin{equation}
\beta_{n}\left(t\right)\overset{\textrm{def}}{=}\frac{1}{p}\sum_{j=0}^{p-1}c_{n,j}e^{-2\pi ijt}
\end{equation}
\begin{equation}
\gamma_{n}\left(t\right)\overset{\textrm{def}}{=}\frac{\beta_{n}\left(t\right)}{\alpha_{n}\left(t\right)}
\end{equation}
\begin{equation}
\hat{A}_{n}\left(t\right)\overset{\textrm{def}}{=}\prod_{m=0}^{-v_{p}\left(t\right)-1}\alpha_{n}\left(p^{m}t\right)
\end{equation}
where:
\begin{equation}
c_{n,j}\overset{\textrm{def}}{=}-\sum_{k=0}^{n-1}\binom{n}{k}a_{j}^{k}b_{j}^{n-k}\hat{X}^{*k}\left(0\right),\textrm{ }\forall n\geq1,\textrm{ }
\end{equation}
In all of the above, all sums are defined to be $0$ when the lower
limit of summation is greater than the upper limit of summation. Likewise,
all products are defined to be $1$ when the lower limit of multiplication
is greater than the upper limit of multiplication.

As defined, $\hat{X}^{*n}$ is then a Fourier transform of $X^{n}$
with respect to $\mathcal{F}$.

\vphantom{}IV. For all $n\geq0$, $\hat{X}^{*n}$ will be the Fourier-Stieltjes
transform of a distribution in $\mathcal{S}\left(\mathbb{Z}_{p},\mathcal{R}/I\mathcal{R}\right)^{\prime}$,
which we identify with the symbol $X^{n}\left(\mathfrak{z}\right)d\mathfrak{z}$.
For any integers $m,n\geq0$, we then define a pointwise product:
\begin{equation}
\left(X^{m}\left(\mathfrak{z}\right)d\mathfrak{z}\right)\left(X^{n}\left(\mathfrak{z}\right)d\mathfrak{z}\right)\overset{\textrm{def}}{=}X^{m+n}\left(\mathfrak{z}\right)d\mathfrak{z}
\end{equation}
This then shows that the distribution $X\left(\mathcal{\mathfrak{z}}\right)d\mathfrak{z}$
then generates a free algebra $K\left[X\left(\mathfrak{z}\right)d\mathfrak{z}\right]$
over $K$ with respect to pointwise multiplication. Moreover, the
Fourier transform induces an isomorphism from this algebra to the
$K$-algebra generated by $\hat{X}$ under convolution.

\textbf{WARNING}: In the convolution algebra generated by $\hat{X}$,
given any integer $n\geq0$, we \emph{define }the convolution of $n$
copies of $\hat{X}$ to be $\hat{X}^{*n}$. It will generally not
be the case that computing the actual convolution directly using the
standard definition:
\begin{equation}
\left(\hat{f}*\hat{g}\right)\left(t\right)\overset{\textrm{def}}{=}\sum_{s\in\hat{\mathbb{Z}}_{p}}\hat{f}\left(t-s\right)\hat{g}\left(s\right)
\end{equation}
will yield, for instance, $\hat{X}*\hat{X}=\hat{X}^{*2}$, because
the sum defining the convolution will usually be either divergent
or non-convergent.
\end{thm}
Proof: As stated earlier, the proof is by using \textbf{Theorem \ref{thm:mini thesis}}
as the base case and then applying \textbf{Theorem \ref{thm:inductive step}}
for $\mathcal{R}_{1}$. This then justifies \textbf{Assumption \ref{assu:main limit lemma}},
which shows that the computations from \textbf{Theorem \ref{thm:formal quasi-integrability of X script J}}
hold true. Here, since $\mathbb{N}_{0}^{d}=\mathbb{N}_{0}$, we have
that:
\begin{align}
r_{m,n,k} & =\binom{n}{m}a_{k}^{m}b_{k}^{n-m}\\
r_{n,k} & =a_{k}^{n}
\end{align}
from which we obtain the formulae for $\hat{f}_{n}$, $\hat{g}_{n}$,
and the rest.

Q.E.D.

\vphantom{}

We now have our main result:
\begin{cor*}
\label{main restated}Let $d,p$ be positive integers, with $d\geq1$
and $p\geq2$, let $K$ be a global field (if $\textrm{char}K>0$,
then assume $\textrm{char}K$ is co-prime to $p$). For each $j\in\left\{ 1,\ldots,d\right\} $
and $k\in\left\{ 0,\ldots,p-1\right\} $ let $a_{j,k}$ and $b_{j,k}$
be indeterminates, and let $R_{d}\left(K\right)$ (a.k.a., $R_{d}$,
$R$) be the ring of polynomials in the $a_{j,k}$s and $b_{j,k}$s
with coefficients in $\mathcal{O}_{K}$, so that $R$ is then a free
$\mathcal{O}_{K}$-algebra in $2dp$ indeterminates. Then, let:
\begin{equation}
\mathcal{R}_{d}\left(K\right)\overset{\textrm{def}}{=}\frac{R\left[x_{1},\ldots,x_{d}\right]}{\left\langle \left(1-a_{1,0}\right)x_{1}-1,\ldots,\left(1-a_{d,0}\right)x_{d}-1\right\rangle }
\end{equation}
(a.k.a., $\mathcal{R}_{d}$, $\mathcal{R}$) be the localization of
$R$ away from the prime ideal $\left\langle 1-a_{1,0},\ldots,1-a_{d,0}\right\rangle $.

Writing:
\begin{equation}
\kappa_{j}\left(n\right)\overset{\textrm{def}}{=}\prod_{k=1}^{p-1}\left(\frac{a_{j,k}}{a_{j,0}}\right)^{\#_{p:k}\left(n\right)}
\end{equation}
let $\mathcal{F}$ be a digital $\mathcal{R}$-frame on $\mathbb{Z}_{p}$
so that:
\begin{equation}
\left|a_{j,0}\right|_{\mathcal{F}\left(\mathfrak{z}\right)}\limsup_{n\rightarrow\infty}\prod_{k=1}^{p-1}\left|\frac{a_{j,k}}{a_{j,0}}\right|_{\mathcal{F}\left(\mathfrak{z}\right)}^{\#_{p:k}\left(\left[\mathfrak{z}\right]_{p^{n}}\right)}<1,\textrm{ }\forall\mathfrak{z}\in D\left(\mathcal{F}\right),\textrm{ }\forall j\in\left\{ 1,\ldots,d\right\} \label{eq:root condition for individual j-1}
\end{equation}
Let $I$ be a unique solution ideal of $R$ ($\left\langle 1-a_{j,0}\right\rangle \cap I=\left\{ 0\right\} $
for all $j$). Then:

\vphantom{}I. For each $j\in\left\{ 1,\ldots,d\right\} $, there
exists a unique rising-continuous function $X_{j}\in C\left(\mathcal{F}/I\right)$
so that:
\begin{equation}
X_{j}\left(\mathfrak{z}\right)\overset{\left(\mathcal{F}/I\right)\left(\mathfrak{z}\right)}{=}a_{j,\left[\mathfrak{z}\right]_{p}}X_{j}\left(\theta_{p}\left(\mathfrak{z}\right)\right)+b_{j,\left[\mathfrak{z}\right]_{p}},\textrm{ }\forall\mathfrak{z}\in D\left(\mathcal{F}/I\right)\label{eq:functional equations}
\end{equation}
where $\left(\mathcal{F}/I\right)\left(\mathfrak{z}\right)$ is the
completion of $\mathcal{R}/I\mathcal{R}$ with respect to the absolute
value $\mathcal{F}/I$ associates to $\mathfrak{z}$.

\vphantom{}

II. For any $\mathbf{n}\in\mathbb{N}_{0}^{d}$, the function $X_{\mathbf{n}}=\prod_{j=1}^{d}X_{j}^{n_{j}}\in C\left(\mathcal{F}/I\right)$
is $\mathcal{F}/I$-quasi-integrable. Moreover, we denote a Fourier
transform of this function by:
\begin{equation}
\hat{X}_{\mathbf{n}}=\bigstar_{j=1}^{d}\hat{X}_{j}^{*n_{j}}\left(t\right)
\end{equation}
where $\bigstar$ denotes formal convolution, with $\hat{X}_{j}^{*n_{j}}\left(t\right)$
being defined as the convolution identity element $\mathbf{1}_{0}$
whenever $n_{j}=0$. Here, $\bigstar_{j=1}^{d}\hat{X}_{j}^{*n_{j}}$
is a function $\hat{\mathbb{Z}}_{p}\rightarrow\textrm{Frac}\left(\mathcal{R}/I\mathcal{R}\right)\left(\zeta_{p^{\infty}}\right)$.

Using \textbf{Theorem \ref{thm:powers of X}} to compute $\hat{X}_{j}$
for all $j\in\left\{ 1,\ldots,d\right\} $, we then have the following
formula for $\hat{X}_{\mathbf{n}}$:

\begin{equation}
\hat{X}_{\mathbf{n}}\left(t\right)=\hat{f}_{\mathbf{n}}\left(t\right)-\hat{g}_{\mathbf{n}}\left(t\right)
\end{equation}
where:
\begin{equation}
\hat{f}_{\mathbf{n}}\left(t\right)=\sum_{\mathbf{m}<\mathbf{n}}\sum_{n=0}^{-v_{p}\left(t\right)-2}\left(\prod_{m=0}^{n-1}\alpha_{\mathbf{n}}\left(p^{m}t\right)\right)\alpha_{\mathbf{m},\mathbf{n}}\left(p^{n}t\right)\hat{X}_{\mathbf{m}}\left(p^{n+1}t\right)\label{eq:f_script J hat-1}
\end{equation}
where the sum is taken over all $\mathbf{m}=\left(m_{1},\ldots,m_{d}\right)$
so that $m_{j}\leq n_{j}$ for all $j\in\left\{ 1,\ldots,d\right\} $
and $m_{j}<n_{j}$ for at least one $j\in\left\{ 1,\ldots,d\right\} $,
and where:
\begin{equation}
\hat{g}_{\mathbf{n}}\left(t\right)=\begin{cases}
0 & \textrm{if }t=0\\
\left(\beta_{\mathbf{n}}\left(0\right)v_{p}\left(t\right)+\gamma_{\mathbf{n}}\left(\frac{t\left|t\right|_{p}}{p}\right)\right)\hat{A}_{\mathbf{n}}\left(t\right) & \textrm{if }t\neq0
\end{cases},\textrm{ }\forall t\in\hat{\mathbb{Z}}_{p}\label{eq:g_script J hat , alpha equals 1-1}
\end{equation}
if $\alpha_{\mathbf{n}}\left(0\right)=1$ and:
\begin{equation}
\hat{g}_{\mathbf{n}}\left(t\right)=\begin{cases}
\frac{\beta_{\mathbf{n}}\left(0\right)}{1-\alpha_{\mathbf{n}}\left(0\right)} & \textrm{if }t=0\\
\left(\frac{\beta_{\mathbf{n}}\left(0\right)}{1-\alpha_{\mathbf{n}}\left(0\right)}+\gamma_{\mathbf{n}}\left(\frac{t\left|t\right|_{p}}{p}\right)\right)\hat{A}_{\mathbf{n}}\left(t\right) & \textrm{if }t\neq0
\end{cases},\textrm{ }\forall t\in\hat{\mathbb{Z}}_{p}\label{eq:g_script J hat , alpha not equal to 1-1}
\end{equation}
if $\alpha_{\mathbf{n}}\left(0\right)\neq1$, where, for all $\mathbf{m},\mathbf{n}\in\mathbb{N}_{0}^{d}$
with $m_{j}\leq n_{j}$ for all $j\in\left\{ 1,\ldots,d\right\} $,
$\alpha_{\mathbf{m},\mathbf{n}},\beta_{\mathbf{n}},\gamma_{\mathbf{n}},\hat{A}_{\mathbf{n}}:\hat{\mathbb{Z}}_{p}\rightarrow\textrm{Frac}\left(\mathcal{R}/I\mathcal{R}\right)\left(\zeta_{p^{\infty}}\right)$
are given by:
\begin{equation}
\alpha_{\mathbf{m},\mathbf{n}}\left(t\right)\overset{\textrm{def}}{=}\frac{1}{p}\sum_{k=0}^{p-1}r_{\mathbf{m},\mathbf{n},k}e^{-2\pi ikt}
\end{equation}
\begin{equation}
\beta_{\mathbf{n}}\left(t\right)\overset{\textrm{def}}{=}\frac{1}{p}\sum_{k=0}^{p-1}c_{\mathbf{n},k}e^{-2\pi ikt}
\end{equation}
\begin{equation}
\gamma_{\mathbf{n}}\left(t\right)\overset{\textrm{def}}{=}\frac{\beta_{\mathbf{n}}\left(t\right)}{\alpha_{\mathbf{n}}\left(t\right)}
\end{equation}
\begin{equation}
\hat{A}_{\mathbf{n}}\left(t\right)\overset{\textrm{def}}{=}\prod_{n=0}^{-v_{p}\left(t\right)-1}\alpha_{\mathbf{n}}\left(p^{n}t\right)
\end{equation}
where:
\begin{equation}
c_{\mathbf{n},k}=-\sum_{\mathbf{m}<\mathbf{n}}r_{\mathbf{m},k}\hat{X}_{\mathbf{m}}\left(0\right)
\end{equation}
\begin{equation}
r_{\mathbf{m},\mathbf{n},k}\overset{\textrm{def}}{=}\binom{\mathbf{n}}{\mathbf{m}}a_{\mathbf{m},k}b_{\mathbf{n}-\mathbf{m},k}
\end{equation}
\begin{equation}
r_{\mathbf{n},k}\overset{\textrm{def}}{=}r_{\mathbf{n},\mathbf{n},k}
\end{equation}
\begin{equation}
a_{\mathbf{m},k}\overset{\textrm{def}}{=}\prod_{j=1}^{d}a_{j,k}^{m_{j}}
\end{equation}
\begin{equation}
b_{\mathbf{n}-\mathbf{m},k}\overset{\textrm{def}}{=}\prod_{j=1}^{d}b_{j,k}^{n_{j}-m_{j}}
\end{equation}
\begin{equation}
\alpha_{\mathbf{n}}\left(t\right)\overset{\textrm{def}}{=}\alpha_{\mathbf{n},\mathbf{n}}\left(t\right)
\end{equation}

\vphantom{}III. For all $j\in\left\{ 1,\ldots,d\right\} $, $\hat{X}_{j}$
will be the Fourier-Stieltjes transform of a distribution in $\mathcal{S}\left(\mathbb{Z}_{p},\textrm{Frac}\left(\mathcal{R}/I\mathcal{R}\right)\right)^{\prime}$,
which we identify with the symbol $X_{j}\left(\mathfrak{z}\right)d\mathfrak{z}$.
Then, for any $\mathbf{n}\in\mathbb{N}_{0}^{d}$, we define a pointwise
product:
\begin{equation}
X_{\mathbf{n}}\left(\mathfrak{z}\right)d\mathfrak{z}\overset{\textrm{def}}{=}\prod_{j=1}^{d}\left(X_{j}\left(\mathfrak{z}\right)d\mathfrak{z}\right)^{n_{j}}\overset{\textrm{def}}{=}\left(\prod_{j=1}^{d}X_{j}^{n_{j}}\left(\mathfrak{z}\right)\right)d\mathfrak{z}
\end{equation}
where $\left(\prod_{j=1}^{d}X_{j}^{n_{j}}\left(\mathfrak{z}\right)\right)d\mathfrak{z}$
is the distribution whose Fourier-Stieltjes transform is $\bigstar_{j=1}^{d}\hat{X}_{j}^{*n_{j}}$,
as given by (II). Consequently, we get a polynomial algebra $K\left[X_{1}\left(\mathfrak{z}\right)d\mathfrak{z},\ldots,X_{d}\left(\mathfrak{z}\right)d\mathfrak{z}\right]$.
Moreover, the Fourier transform induces an isomorphism from this algebra
to the $K$-algebra generated by $\hat{X}_{1},\ldots,\hat{X}_{d}$
under convolution.
\end{cor*}
Proof: As stated earlier, the proof is by using \textbf{Theorem \ref{thm:powers of X}}
as the base case and then applying \textbf{Theorem \ref{thm:inductive step}}
for $\mathcal{R}_{d}$. This then justifies \textbf{Assumption \ref{assu:main limit lemma}},
which shows that the computations from \textbf{Theorem \ref{thm:formal quasi-integrability of X script J}}
hold true.

Q.E.D.
\begin{rem}
The restriction of $D\left(\mathcal{F}\right)$ to a subset of $\textrm{NiceDig}\left(\mathbb{Z}_{p}\right)$
is \emph{probably} not necessary, but I have not investigated it in
detail, as the proofs here are already intricate enough. The main
point of sensitivity is \textbf{Proposition \ref{prop:nice digits}}
and the way it interacts with \textbf{Propositions \ref{prop:gamma prop}}
and \textbf{\ref{prop:gamma theta prop}}.
\end{rem}
While the $X_{\mathbf{n}}$s will always be distributions, the places
of $\mathcal{R}$ (and, in the quotient case, $\mathcal{R}/I\mathcal{R}$)
for which $X_{\mathbf{n}}$ becomes a measure captures information
about the norms of $\alpha_{\mathbf{n}}$ and the $\alpha_{\mathbf{m},\mathbf{n}}$s
with respect to the supremum norm induced by the chosen place of $\mathcal{R}$.
\begin{rem}
For what follows, given an ideal $I\subseteq R$, we write $\left\Vert \cdot\right\Vert _{p,\ell}$
to denote the norm:
\begin{equation}
\left\Vert \hat{f}\right\Vert _{p,\ell}\overset{\textrm{def}}{=}\sup_{t\in\hat{\mathbb{Z}}_{p}}\left|\hat{f}\left(t\right)\right|_{\ell}
\end{equation}
where $\ell$ is a place of $\mathcal{R}/I\mathcal{R}$ and $\left|\cdot\right|_{\ell}$
is any associated absolute value; $\hat{f}$ here is a function from
$\hat{\mathbb{Z}}_{p}\rightarrow\textrm{Frac}\left(\mathcal{R}/I\mathcal{R}\right)\left(\zeta_{p^{\infty}}\right)$.
\end{rem}
\begin{cor}
\label{thm:X_Js as measures}Fix integers $p\geq2$ and $d\geq1$,
and a unique solution ideal $I\subseteq R_{d}$. Using \textbf{Theorem
\ref{cor:main result}}, fix a choice of Fourier transforms $\hat{X}_{\mathbf{n}}$
of $X_{\mathbf{n}}$ for all $\mathbf{n}\in\mathbb{N}_{0}^{d}$, and
use these to define $\hat{f}_{\mathbf{n}}$ and $\hat{g}_{\mathbf{n}}$
for all $\mathbf{n}$ with $\Sigma\left(\mathbf{n}\right)\geq2$.
Let $\ell$ be a non-trivial place of $\mathcal{R}_{d}/I\mathcal{R}_{d}$;
if $\ell$ is non-archimedean, we require the residue field of the
completion $\textrm{Frac}\left(\mathcal{R}_{d}/I\mathcal{R}_{d}\right)_{\ell}$
to have characteristic co-prime to $p$. Finally, let $\mathbf{n}\in\mathbb{N}_{0}^{d}\backslash\left\{ \mathbf{0}\right\} $.

I. If $\ell$ is non-archimedean and $\max_{\mathbf{m}\leq\mathbf{n}}\left\Vert \alpha_{\mathbf{m},\mathbf{n}}\right\Vert _{p,\ell}\leq1$,
then the formula:
\begin{align*}
W\left(\mathbb{Z}_{p},\overline{\overline{\textrm{Frac}\left(\mathcal{R}_{d}/I\mathcal{R}_{d}\right)_{\ell}}}\right) & \rightarrow\overline{\overline{\textrm{Frac}\left(\mathcal{R}_{d}/I\mathcal{R}_{d}\right)_{\ell}}}\\
\phi & \mapsto\sum_{t\in\hat{\mathbb{Z}}_{p}}\hat{\phi}\left(t\right)\hat{X}_{\mathbf{n}}\left(-t\right)
\end{align*}
defines an element of $W\left(\mathbb{Z}_{p},\overline{\overline{\textrm{Frac}\left(\mathcal{R}_{d}/I\mathcal{R}_{d}\right)_{\ell}}}\right)^{\prime}$.
This shows that $\hat{X}_{\mathbf{n}}$ is then the Fourier-Stieltjes
transform of $\left(p,\ell\right)$-adic measure, and as such, we
may identify this measure with $X_{\mathbf{n}}\left(\mathfrak{z}\right)d\mathfrak{z}$
and write:
\begin{equation}
\int_{\mathbb{Z}_{p}}\phi\left(\mathfrak{z}\right)X_{\mathbf{n}}\left(\mathfrak{z}\right)d\mathfrak{z}\overset{\textrm{def}}{=}\sum_{t\in\hat{\mathbb{Z}}_{p}}\hat{\phi}\left(t\right)\hat{X}_{\mathbf{n}}\left(-t\right),\textrm{ }\forall\phi\in W\left(\mathbb{Z}_{p},\overline{\overline{\textrm{Frac}\left(\mathcal{R}_{d}/I\mathcal{R}_{d}\right)_{\ell}}}\right)
\end{equation}
Moreover, the operator norm of this measure is:
\begin{equation}
\left\Vert X_{\mathbf{n}}\left(\mathfrak{z}\right)d\mathfrak{z}\right\Vert =\left\Vert \hat{X}_{\mathbf{n}}\right\Vert _{p,\ell}\leq1
\end{equation}

II. If $\ell$ is archimedean and:
\begin{equation}
\max_{\mathbf{m}:\mathbf{0}<\mathbf{m}\leq\mathbf{n}}\left\{ \left\Vert \alpha_{\mathbf{m},\mathbf{n}}\right\Vert _{p,\ell}\right\} <1\label{eq:measure hypothesis}
\end{equation}
the Parseval-Plancherel construction given in (I) makes $X_{\mathbf{n}}\left(\mathfrak{z}\right)d\mathfrak{z}$
into an element of $W\left(\mathbb{Z}_{p},\overline{\overline{\textrm{Frac}\left(\mathcal{R}_{d}/I\mathcal{R}_{d}\right)_{\ell}}}\right)^{\prime}$
if and only if $\alpha_{\mathbf{n}}\left(0\right)\neq1$ in $\mathcal{R}_{d}/I\mathcal{R}_{d}$.
When $\alpha_{\mathbf{n}}\left(0\right)=1$ in $\mathcal{R}_{d}/I\mathcal{R}_{d}$,
$X_{\mathbf{n}}\left(\mathfrak{z}\right)d\mathfrak{z}$ will merely
be a distribution in $\mathcal{S}\left(\mathbb{Z}_{p},\mathcal{R}_{d}/I\mathcal{R}_{d}\right)^{\prime}$.
\end{cor}
Proof: For brevity, let $\mathbb{K}$ denote $\textrm{Frac}\left(\mathcal{R}/I\mathcal{R}\right)$.
By construction, the continuous dual of $W\left(\mathbb{Z}_{p},\overline{\overline{\mathbb{K}_{\ell}}}\right)$
is then isometrically isomorphic to $B\left(\hat{\mathbb{Z}}_{p},\overline{\overline{\mathbb{K}_{\ell}}}\right)$
(see any of \cite{2nd blog paper,My Dissertation,Schikhof's Thesis,van Rooij - Non-Archmedean Functional Analysis}
for details), the Banach space of functions $\hat{\mathbb{Z}}_{p}\rightarrow B\left(\hat{\mathbb{Z}}_{p},\overline{\overline{\mathbb{K}_{\ell}}}\right)$
with:
\begin{equation}
\lim_{\left|t\right|_{p}\rightarrow\infty}\left|\hat{f}\left(t\right)\right|_{\ell}\overset{\mathbb{R}}{=}0
\end{equation}
Thus, to show that $X_{\mathbf{n}}\left(\mathfrak{z}\right)d\mathfrak{z}$
is an element of $W\left(\hat{\mathbb{Z}}_{p},\overline{\overline{\mathbb{K}_{\ell}}}\right)^{\prime}$,
it suffices to show that that $\hat{X}_{\mathbf{n}}$ is $\ell$-adically
bounded over $\hat{\mathbb{Z}}_{p}$. We do this by induction on $\Sigma\left(\mathbf{n}\right)$.

For the base case, letting $j\in\left\{ 1,\ldots,d\right\} $, \textbf{Theorem
\ref{thm:mini thesis}} gives us:

\begin{equation}
\hat{X}_{j}\left(t\right)=\begin{cases}
0 & \textrm{if }t=0\\
\left(\beta_{j}\left(0\right)v_{p}\left(t\right)+\gamma_{j}\left(\frac{t\left|t\right|_{p}}{p}\right)\right)\hat{A}_{j}\left(t\right) & \textrm{if }t\neq0
\end{cases},\textrm{ }\forall t\in\hat{\mathbb{Z}}_{p}
\end{equation}
when $\alpha_{j}\left(0\right)=1$, and:
\begin{equation}
\hat{X}_{j}\left(t\right)=\begin{cases}
\frac{\beta_{j}\left(0\right)}{1-\alpha_{j}\left(0\right)} & \textrm{if }t=0\\
\frac{\beta_{j}\left(0\right)\hat{A}_{j}\left(t\right)}{1-\alpha_{j}\left(0\right)}+\gamma_{j}\left(\frac{t\left|t\right|_{p}}{p}\right)\hat{A}_{j}\left(t\right) & \textrm{if }t\neq0
\end{cases},\textrm{ }\forall t\in\hat{\mathbb{Z}}_{p}
\end{equation}
when $\alpha_{j}\left(0\right)\neq1$ as Fourier transforms of $X_{j}$.

Here:
\begin{equation}
\hat{A}_{j}\left(t\right)=\prod_{n=0}^{-v_{p}\left(t\right)-1}\alpha_{j}\left(p^{n}t\right)=\prod_{n=0}^{-v_{p}\left(t\right)-1}\left(\frac{1}{p}\sum_{k=0}^{p-1}a_{j,k}e^{-2\pi ikp^{n}t}\right)
\end{equation}
and:
\begin{equation}
\gamma_{j}\left(t\right)=\frac{\beta_{j}\left(t\right)}{\alpha_{j}\left(t\right)}=\frac{\frac{1}{p}\sum_{k=0}^{p-1}b_{j,k}e^{-2\pi ikt}}{\frac{1}{p}\sum_{k=0}^{p-1}a_{j,k}e^{-2\pi ikt}}
\end{equation}
Here, $t\left|t\right|_{p}/p$ is a rational number of the form $k/p$
for all $t\in\hat{\mathbb{Z}}_{p}$, for some $k\in\left\{ 0,\ldots,p-1\right\} $,
with $k=0$ if and only if $t=0$. Hence, for each $t$:
\begin{equation}
\gamma_{j}\left(\frac{t\left|t\right|_{p}}{p}\right)\in\mathbb{K}\left(e^{2\pi it}\right)
\end{equation}
Likewise $\hat{A}_{j}\left(t\right)\in\mathbb{K}\left(e^{2\pi it}\right)$.
Thus, for each $t$, the quantity $\hat{X}_{j}\left(t\right)$ is
contained in the cyclotomic extension $\mathbb{K}\left(e^{2\pi it}\right)$
of $\mathbb{K}$.

So, letting $\ell$ be a place of $\mathbb{K}$, have that $\left|\cdot\right|_{\ell}$
admits a unique extension to $\mathbb{K}\left(e^{2\pi it}\right)$
for each $t$, and so we may consider the extension of $\left|\cdot\right|_{\ell}$
to the maximal $p$-power cyclotomic extension of $\mathbb{K}$.

Since $\gamma_{j}\left(t\left|t\right|_{p}/p\right)$ takes on only
finitely many values as $t$ varies, our formula for $\hat{X}_{j}\left(t\right)$
gives:
\begin{equation}
\left|\hat{X}_{j}\left(t\right)\right|_{\ell}\ll\begin{cases}
\left|\min\left\{ 0,v_{p}\left(t\right)\right\} \right|_{\ell}\left|\hat{A}_{j}\left(t\right)\right|_{\ell} & \textrm{if }R_{j}\left(0\right)=1\\
\left|\hat{A}_{j}\left(t\right)\right|_{\ell} & \textrm{if }R_{j}\left(0\right)\neq1
\end{cases}
\end{equation}
where $v_{p}\left(0\right)=+\infty$ and $v_{p}\left(t\right)<0$
for all $t\in\hat{\mathbb{Z}}_{p}\backslash\left\{ 0\right\} $. Every
$\ell$-adically bounded function $\hat{\mathbb{Z}}_{p}\rightarrow\overline{\overline{\mathbb{K}_{\ell}}}$
induces a measure on $W\left(\mathbb{Z}_{p},\overline{\overline{\mathbb{K}_{\ell}}}\right)$
by the Parseval-Plancherel construction given in \textbf{Theorem \ref{thm:The-breakdown-variety}}
(see \cite{2nd blog paper} for details). Hence, letting:
\begin{equation}
\textrm{m}_{p,j}\left(t\right)\overset{\textrm{def}}{=}\begin{cases}
\min\left\{ 0,v_{p}\left(t\right)\right\}  & \textrm{if }R_{j}\left(0\right)=1\\
1 & \textrm{if }R_{j}\left(0\right)\neq1
\end{cases}
\end{equation}
we have:
\begin{equation}
\sup_{t\in\hat{\mathbb{Z}}_{p}}\left|\textrm{m}_{p,j}\left(t\right)\hat{A}_{j}\left(t\right)\right|_{\ell}<\infty
\end{equation}
implies $X_{j}$ is a measure in $W\left(\mathbb{Z}_{p},\overline{\overline{\mathbb{K}_{\ell}}}\right)^{\prime}$
via the aforementioned the Parseval-Plancherel construction. Now,
we have two cases, based on the quality of $\ell$.

\vphantom{}

I. Suppose $\ell$ is non-archimedean, and that:
\begin{equation}
\max_{\mathbf{m}\leq\mathbf{n}}\left\Vert \alpha_{\mathbf{m},\mathbf{n}}\right\Vert _{p,\ell}\leq1
\end{equation}
Here, the non-archimedean quality of $\ell$ gives us $\left|\textrm{m}_{p,j}\left(t\right)\right|_{\ell}\leq1$
for all $j$. So, for any $\mathbf{m},\mathbf{n}$ with $\mathbf{m}\leq\mathbf{n}$:
\begin{equation}
\left|\hat{A}_{\mathbf{m},\mathbf{n}}\left(t\right)\right|_{\ell}=\prod_{n=0}^{-v_{p}\left(t\right)-1}\left|\alpha_{\mathbf{m},\mathbf{n}}\left(p^{n}t\right)\right|_{\ell}\leq\left\Vert \alpha_{\mathbf{m},\mathbf{n}}\right\Vert _{p,\ell}^{\max\left\{ 0,-v_{p}\left(t\right)\right\} },\textrm{ }\forall t\in\hat{\mathbb{Z}}_{p}
\end{equation}
as such, using \textbf{Theorem \ref{cor:main result}}'s formula for
$\hat{X}_{j}$ for any given $j\in\left\{ 1,\ldots,d\right\} $, we
have: 
\begin{equation}
\left|\hat{X}_{j}\left(t\right)\right|_{\ell}\ll_{j,\ell}\left\Vert \alpha_{j}\right\Vert _{p,\ell}^{\max\left\{ 0,-v_{p}\left(t\right)\right\} },\textrm{ }\forall t\in\hat{\mathbb{Z}}_{p}\label{eq:base case estimate}
\end{equation}
as well as the trivial bound on $\hat{X}_{\mathbf{0}}=\mathbf{1}_{0}$,
where, as indicated, the constant of proportionality depends only
on $j$ and $\ell$. Thus, $\left\Vert \alpha_{\mathbf{n}}\right\Vert _{p,\ell}\leq1$
is sufficient to prove that $\hat{X}_{\mathbf{n}}\left(t\right)$
is the Fourier transform of a $\left(p,\ell\right)$-adic measure
in the base case where $\Sigma\left(\mathbf{n}\right)=1$.

Moving to the inductive step, suppose that $\max_{\mathbf{m}\leq\mathbf{n}}\left\Vert \alpha_{\mathbf{m}}\right\Vert _{p,\ell}\leq1$
and: 
\begin{equation}
\max_{\mathbf{m}<\mathbf{n}}\left\Vert \hat{X}_{\mathbf{m}}\right\Vert _{p,\ell}<\infty
\end{equation}
Using:

\begin{equation}
\hat{f}_{\mathbf{n}}\left(t\right)=\sum_{\mathbf{m}<\mathbf{n}}\sum_{n=0}^{-v_{p}\left(t\right)-2}\left(\prod_{m=0}^{n-1}\alpha_{\mathbf{n}}\left(p^{m}t\right)\right)\alpha_{\mathbf{m},\mathbf{n}}\left(p^{n}t\right)\hat{X}_{\mathbf{m}}\left(p^{n+1}t\right)
\end{equation}
we apply the $\ell$-adic ultrametric inequality to get:
\begin{equation}
\left|\hat{f}_{\mathbf{m}}\left(t\right)\right|_{\ell}\leq\max_{\mathbf{m}<\mathbf{n}}\max_{0\leq n\leq-v_{p}\left(t\right)-2}\left|\left(\prod_{m=0}^{n-1}\alpha_{\mathbf{n}}\left(p^{m}t\right)\right)\alpha_{\mathbf{m},\mathbf{n}}\left(p^{n}t\right)\hat{X}_{\mathbf{m}}\left(p^{n+1}t\right)\right|_{\ell}
\end{equation}
Meanwhile, the formula for $\hat{g}_{\mathbf{n}}\left(t\right)$ and
the $\ell$-adic boundedness of $v_{p}\left(t\right)$ on $\hat{\mathbb{Z}}_{p}\backslash\left\{ 0\right\} $
gives:
\begin{equation}
\left|\hat{g}_{\mathbf{n}}\left(t\right)\right|_{\ell}\ll_{\mathbf{n},\ell}\left|\hat{A}_{\mathbf{n}}\left(t\right)\right|_{\ell},\textrm{ }\forall t\in\hat{\mathbb{Z}}_{p}
\end{equation}
where, as indicated, the constant of proportionality depends only
on $\mathbf{n}$ and $\ell$. Then:
\begin{align*}
\left|\hat{X}_{\mathbf{n}}\left(t\right)\right|_{\ell} & \leq\max\left\{ \left|\hat{g}_{\mathbf{n}}\left(t\right)\right|_{\ell},\left|\hat{f}_{\mathbf{n}}\left(t\right)\right|_{\ell}\right\} \\
 & \ll_{\mathbf{n},\ell}\max\left\{ \left|\hat{A}_{\mathbf{n}}\left(t\right)\right|_{\ell},\max_{I\subset\mathbf{n}}\max_{0\leq n\leq-v_{p}\left(t\right)-2}\left|\left(\prod_{m=0}^{n-1}\alpha_{\mathbf{n}}\left(p^{m}t\right)\right)\alpha_{\mathbf{m},\mathbf{n}}\left(p^{n}t\right)\hat{X}_{\mathbf{m}}\left(p^{n+1}t\right)\right|_{\ell}\right\} \\
 & \leq\max\left\{ \left|\hat{A}_{\mathbf{n}}\left(t\right)\right|_{\ell},\underbrace{\left(\max_{\mathbf{m}<\mathbf{n}}\left\Vert \hat{X}_{\mathbf{m}}\right\Vert _{p,\ell}\right)}_{<\infty}\underbrace{\left(\max_{\mathbf{m}^{\prime}<\mathbf{n}}\left\Vert \alpha_{\mathbf{m}^{\prime},\mathbf{n}}\right\Vert _{p,\ell}\right)}_{\infty}\left(\max\left\{ 1,\left\Vert \alpha_{\mathbf{n}}\right\Vert _{p,\ell}\right\} \right)^{\max\left\{ 0,-v_{p}\left(t\right)-2\right\} }\right\} \\
 & \ll_{\mathbf{n},\ell}\max\left\{ \left|\hat{A}_{\mathbf{n}}\left(t\right)\right|_{\ell},\left(\max\left\{ 1,\left\Vert \alpha_{\mathbf{n}}\right\Vert _{p,\ell}\right\} \right)^{\max\left\{ 0,-v_{p}\left(t\right)-2\right\} }\right\} 
\end{align*}
Here, note that $\max_{\mathbf{m}^{\prime}<\mathbf{n}}\left\Vert \alpha_{\mathbf{m}^{\prime},\mathbf{n}}\right\Vert _{p,\ell}<\infty$
independent of any hypothesis, since each $\alpha_{\mathbf{m}^{\prime},\mathbf{n}}$
is a degree $p-1$ polynomial in $e^{-2\pi it}$. Meanwhile:
\begin{equation}
\left|\hat{A}_{\mathbf{n}}\left(t\right)\right|_{\ell}=\left|\prod_{m=0}^{-v_{p}\left(t\right)-1}\alpha_{\mathbf{n}}\left(p^{m}t\right)\right|_{\ell}\leq\left(\max\left\{ 1,\left\Vert \alpha_{\mathbf{n}}\right\Vert _{p,\ell}\right\} \right)^{\max\left\{ 0,-v_{p}\left(t\right)-1\right\} }
\end{equation}
and so:
\begin{equation}
\left|\hat{X}_{\mathbf{n}}\left(t\right)\right|_{\ell}\ll_{\mathbf{n},\ell}\left(\max\left\{ 1,\left\Vert \alpha_{\mathbf{n}}\right\Vert _{p,\ell}\right\} \right)^{\max\left\{ 0,-v_{p}\left(t\right)-1\right\} }\leq1,\textrm{ }\forall t\in\hat{\mathbb{Z}}_{p}
\end{equation}
Hence, the hypothesis $\left\Vert \alpha_{\mathbf{n}}\right\Vert _{p,\ell}\leq1$
is then sufficient to show that $\hat{X}_{\mathbf{n}}$ is the Fourier-Stieltjes
transform of a $\left(p,\ell\right)$-adic measure, $X_{\mathbf{n}}\left(\mathfrak{z}\right)d\mathfrak{z}$.
Since the Fourier-Stieltjes transform is an isometric isomorphism
from $W\left(\mathbb{Z}_{p},\overline{\overline{\mathbb{K}_{\ell}}}\right)^{\prime}$
onto the space of bounded functions $\hat{\mathbb{Z}}_{p}\rightarrow\overline{\overline{\mathbb{K}_{\ell}}}$,
we then get the estimate:
\begin{equation}
\left\Vert X_{\mathbf{n}}\left(\mathfrak{z}\right)d\mathfrak{z}\right\Vert =\left\Vert \hat{X}_{\mathbf{n}}\right\Vert _{p,\ell}\leq1
\end{equation}
as desired.

\vphantom{}

II. Suppose $\ell$ is archimedean. Like before, we get the estimate:
\begin{equation}
\left|\hat{A}_{\mathbf{m},\mathbf{n}}\left(t\right)\right|_{\ell}\leq\left\Vert \alpha_{\mathbf{m},\mathbf{n}}\right\Vert _{p,\ell}^{\max\left\{ 0,-v_{p}\left(t\right)\right\} },\textrm{ }\forall t\in\hat{\mathbb{Z}}_{p}
\end{equation}
$\alpha_{\mathbf{n}}\left(0\right)\neq1$, we do not need to contend
with a factor of $v_{p}\left(t\right)$ in $\hat{g}_{\mathbf{n}}$,
and so we once again get the estimate:
\begin{equation}
\left|\hat{g}_{\mathbf{n}}\left(t\right)\right|_{\ell}\ll_{\mathbf{n},\ell}\left|\hat{A}_{\mathbf{n}}\left(t\right)\right|_{\ell}\leq\left\Vert \alpha_{\mathbf{n}}\right\Vert _{p,\ell}^{\max\left\{ 0,-v_{p}\left(t\right)\right\} },\textrm{ }\forall t\in\hat{\mathbb{Z}}_{p}
\end{equation}
As for $\hat{f}_{\mathbf{n}}$, using the ordinary triangle inequality
gives:
\begin{equation}
\left|\hat{f}_{\mathbf{n}}\left(t\right)\right|_{\ell}\ll_{\mathbf{n}}\sum_{\mathbf{m}<\mathbf{n}}\sum_{n=0}^{-v_{p}\left(t\right)-2}\left\Vert \alpha_{\mathbf{m},\mathbf{n}}\right\Vert _{p,\ell}^{n}\left|\hat{X}_{\mathbf{m}}\left(p^{n+1}t\right)\right|_{\ell}
\end{equation}
Since $\hat{X}_{\mathbf{0}}=\mathbf{1}_{0}$, $\mathbf{m}=\mathbf{0}$
term of the upper bound is bounded by:
\begin{equation}
\sum_{n=0}^{-v_{p}\left(t\right)-2}\left\Vert \alpha_{\mathbf{0},\mathbf{n}}\right\Vert _{p,\ell}^{n}\leq\sum_{n=0}^{\infty}\left\Vert \alpha_{\mathbf{0},\mathbf{n}}\right\Vert _{p,\ell}^{n},\textrm{ }\forall t\in\hat{\mathbb{Z}}_{p}
\end{equation}
Hence, our assumption that $\left\Vert \alpha_{\mathbf{0},\mathbf{n}}\right\Vert _{p,\ell}<1$
guarantees the bound on the $\mathbf{m}=\mathbf{0}$th term. As such,
it suffices to deal with the terms for which $\mathbf{m}>\mathbf{0}$.

Now, for $j\in\left\{ 1,\ldots,d\right\} $, we have:
\begin{equation}
\left|\hat{X}_{j}\left(t\right)\right|_{\ell}\ll_{j,\ell}\left(-v_{p}\left(t\right)\right)^{\left[\alpha_{j}\left(0\right)=1\right]}\left\Vert \alpha_{j}\right\Vert _{p,\ell}^{\max\left\{ 0,-v_{p}\left(t\right)\right\} },\textrm{ }\forall t\in\hat{\mathbb{Z}}_{p}
\end{equation}
Thus, for $\Sigma\left(\mathbf{n}\right)=2$, letting $\mathbf{e}_{j}$
denote the element of $\mathbb{N}_{0}^{d}$ with a $1$ in the $j$the
entry and a $0$ in all other entries, we have:
\begin{align}
\left|\hat{f}_{\mathbf{n}}\left(t\right)\right|_{\ell} & \ll_{\mathbf{n}}\sum_{j\in\left\{ 1,\ldots,d\right\} :n_{j}\geq1}\sum_{n=0}^{-v_{p}\left(t\right)-2}\left\Vert \alpha_{\mathbf{e}_{j},\mathbf{n}}\right\Vert _{p,\ell}^{n}\left|\hat{X}_{j}\left(p^{n+1}t\right)\right|_{\ell}\nonumber \\
 & \ll\sum_{j\in\left\{ 1,\ldots,d\right\} :n_{j}\geq1}\sum_{n=0}^{-v_{p}\left(t\right)-2}\left\Vert \alpha_{\mathbf{e}_{j},\mathbf{n}}\right\Vert _{p,\ell}^{n}\left(-v_{p}\left(t\right)-n-1\right)^{\left[\alpha_{j}\left(0\right)=1\right]}\left\Vert \alpha_{j}\right\Vert _{p,\ell}^{\max\left\{ 0,-v_{p}\left(t\right)-n-1\right\} }\nonumber \\
\left(m=-v_{p}\left(t\right)-n\right); & =\sum_{j\in\left\{ 1,\ldots,d\right\} :n_{j}\geq1}\sum_{m=1}^{-v_{p}\left(t\right)-1}\left\Vert \alpha_{\mathbf{e}_{j},\mathbf{n}}\right\Vert _{p,\ell}^{-v_{p}\left(t\right)-m-1}m^{\left[\alpha_{j}\left(0\right)=1\right]}\left\Vert \alpha_{j}\right\Vert _{p,\ell}^{m+1}\nonumber \\
 & \ll_{\mathbf{n}}\sum_{j\in\left\{ 1,\ldots,d\right\} :n_{j}\geq1}\left\Vert \alpha_{\mathbf{e}_{j},\mathbf{n}}\right\Vert _{p,\ell}^{-v_{p}\left(t\right)}\sum_{m=1}^{-v_{p}\left(t\right)-1}m^{\left[\alpha_{j}\left(0\right)=1\right]}\left(\frac{\left\Vert \alpha_{j}\right\Vert _{p,\ell}}{\left\Vert \alpha_{\mathbf{e}_{j},\mathbf{n}}\right\Vert _{p,\ell}}\right)^{m}\label{eq:this}
\end{align}

Here, we get several possibilities:

\textbullet{} $\alpha_{j}\left(0\right)\neq1$, in which case:

\begin{equation}
\sum_{m=1}^{-v_{p}\left(t\right)-1}m^{\left[\alpha_{j}\left(0\right)=1\right]}\left(\frac{\left\Vert \alpha_{j}\right\Vert _{p,\ell}}{\left\Vert \alpha_{\mathbf{e}_{j},\mathbf{n}}\right\Vert _{p,\ell}}\right)^{m}\ll1+\left(\frac{\left\Vert \alpha_{j}\right\Vert _{p,\ell}}{\left\Vert \alpha_{\mathbf{e}_{j},\mathbf{n}}\right\Vert _{p,\ell}}\right)^{-v_{p}\left(t\right)}
\end{equation}

\textbullet{} $\alpha_{j}\left(0\right)=1$ and $\left\Vert \alpha_{j}\right\Vert _{p,\ell}/\left\Vert \alpha_{\mathbf{e}_{j},\mathbf{n}}\right\Vert _{p,\ell}\neq1$,
in which case:
\begin{equation}
\sum_{m=1}^{-v_{p}\left(t\right)-1}m^{\left[\alpha_{j}\left(0\right)=1\right]}\left(\frac{\left\Vert \alpha_{j}\right\Vert _{p,\ell}}{\left\Vert \alpha_{\mathbf{e}_{j},\mathbf{n}}\right\Vert _{p,\ell}}\right)^{m}\ll1+\left|v_{p}\left(t\right)\right|\left(\frac{\left\Vert \alpha_{j}\right\Vert _{p,\ell}}{\left\Vert \alpha_{\mathbf{e}_{j},\mathbf{n}}\right\Vert _{p,\ell}}\right)^{-v_{p}\left(t\right)}
\end{equation}

\textbullet{} $\alpha_{j}\left(0\right)=1$ and $\left\Vert \alpha_{j}\right\Vert _{p,\ell}/\left\Vert \alpha_{\mathbf{e}_{j},\mathbf{n}}\right\Vert _{p,\ell}=1$,
in which case:
\begin{equation}
\sum_{m=1}^{-v_{p}\left(t\right)-1}m^{\left[\alpha_{j}\left(0\right)=1\right]}\left(\frac{\left\Vert \alpha_{j}\right\Vert _{p,\ell}}{\left\Vert \alpha_{\mathbf{e}_{j},\mathbf{n}}\right\Vert _{p,\ell}}\right)^{m}\ll\left|v_{p}\left(t\right)\right|^{2}
\end{equation}
Using this with (\ref{eq:this}), we see that $\left|\hat{f}_{\mathbf{n}}\left(t\right)\right|_{\ell}$
will remain bounded as $\left|t\right|_{p}\rightarrow\infty$ provided
that:

\begin{equation}
\max\left\{ \left\Vert \alpha_{j}\right\Vert _{p,\ell},\left\Vert \alpha_{\mathbf{e}_{j},\mathbf{n}}\right\Vert _{p,\ell}\right\} <1,\textrm{ }\forall j:\mathbf{e}_{j}\leq\mathbf{n}
\end{equation}
In particular, if this condition is satisfied, we will have:
\begin{equation}
\left|\hat{f}_{\mathbf{n}}\left(t\right)\right|_{\ell}\ll_{\mathbf{n}}\left(\max\left\{ 1,\max_{j\in\left\{ 1,\ldots,d\right\} :n_{j}\geq1}\left\{ \left\Vert \alpha_{j}\right\Vert _{p,\ell},\left\Vert \alpha_{\mathbf{e}_{j},\mathbf{n}}\right\Vert _{p,\ell}\right\} \right\} \right)^{\max\left\{ 0,-v_{p}\left(t\right)\right\} }
\end{equation}
 By induction on $\left|\mathbf{n}\right|$, we have that:
\begin{equation}
\left|\hat{f}_{\mathbf{n}}\left(t\right)\right|_{\ell}\ll_{\mathbf{n}}\left(\max\left\{ 1,\max_{\mathbf{m}:\mathbf{0}<\mathbf{m}\leq\mathbf{n}}\left\{ \left\Vert \alpha_{\mathbf{m}}\right\Vert _{p,\ell},\left\Vert \alpha_{\mathbf{m},\mathbf{n}}\right\Vert _{p,\ell}\right\} \right\} \right)^{\max\left\{ 0,-v_{p}\left(t\right)\right\} }
\end{equation}
will occur provided that:
\begin{equation}
\max_{\mathbf{m}:\mathbf{0}<\mathbf{m}\leq\mathbf{n}}\max\left\{ \left\Vert \alpha_{\mathbf{m}}\right\Vert _{p,\ell},\left\Vert \alpha_{\mathbf{m},\mathbf{n}}\right\Vert _{p,\ell}\right\} <1
\end{equation}
Meanwhile, the formula for $\hat{g}_{\mathbf{n}}\left(t\right)$ gives:
\begin{equation}
\left|\hat{g}_{\mathbf{n}}\left(t\right)\right|_{\ell}\ll_{\mathbf{n},\ell}\left|v_{p}\left(t\right)\right|^{\left[\alpha_{\mathbf{n}}\left(0\right)=1\right]}\left|\hat{A}_{\mathbf{n}}\left(t\right)\right|_{\ell}\leq\left|v_{p}\left(t\right)\right|^{\left[\alpha_{\mathbf{n}}\left(0\right)=1\right]}\left\Vert \alpha_{\mathbf{n}}\right\Vert _{p,\ell}^{-v_{p}\left(t\right)},\textrm{ }\forall t\in\hat{\mathbb{Z}}_{p}\backslash\left\{ 0\right\} 
\end{equation}
Hence, as hypothesized:
\begin{equation}
\max_{\mathbf{m}:\mathbf{0}<\mathbf{m}\leq\mathbf{n}}\max\left\{ \left\Vert \alpha_{\mathbf{m}}\right\Vert _{p,\ell},\left\Vert \alpha_{\mathbf{m},\mathbf{n}}\right\Vert _{p,\ell}\right\} <1
\end{equation}
is then sufficient to guarantee that both $\hat{f}_{\mathbf{n}}$
and $\hat{g}_{\mathbf{n}}$ are bounded with respect to $\left|\cdot\right|_{\ell}$,
and thus, that $\hat{X}_{\mathbf{n}}$ is the Fourier-Stieltjes transform
of an element of $W\left(\mathbb{Z}_{p},\overline{\overline{\mathbb{K}_{\ell}}}\right)^{\prime}$.

Finally, note that if there is an $\mathbf{m}\leq\mathbf{n}$ with
$\Sigma\left(\mathbf{m}\right)\geq1$ so that either $\alpha_{\mathbf{m},\mathbf{n}}\left(0\right)=1$
or $\alpha_{\mathbf{m}}\left(0\right)=1$, it will follow that $\max_{\mathbf{m}:\mathbf{0}<\mathbf{m}\leq\mathbf{n}}\max\left\{ \left\Vert \alpha_{\mathbf{m}}\right\Vert _{p,\ell},\left\Vert \alpha_{\mathbf{m},\mathbf{n}}\right\Vert _{p,\ell}\right\} =1$,
which contradicts the hypothesis (\ref{eq:measure hypothesis}). In
this case, $\hat{X}_{\mathbf{n}}$ will no longer be bounded with
respect to $t$, and thus will \emph{not} be the Fourier-Stieltjes
transform of an element of $W\left(\mathbb{Z}_{p},\overline{\overline{\mathbb{K}_{\ell}}}\right)^{\prime}$.

Q.E.D.
\begin{rem}
One can do fun things with this result. Keeping the notation as from
the above theorem, suppose $\ell=q$ is non-archimedean, and that
we are working in the power case, with $X^{n}$ being a $\left(p,q\right)$-adic
measure for all $n\geq0$. For simplicity, let's suppose that we can
realize $X$ as a $q$-adic valued function on $\mathbb{Z}_{p}$.
An concrete example of this would be:
\begin{equation}
X\left(\mathfrak{z}\right)=\sum_{n=0}^{\infty}a_{0}^{n}\prod_{j=1}^{p-1}a_{j}^{\#_{p:j}\left(\left[\mathfrak{z}\right]_{p^{n}}\right)}
\end{equation}
where the $a_{j}$s are elements of $\mathbb{Q}^{\times}$ so that
$\left|a_{j}\right|_{q}<1$ for all $j\in\left\{ 1,\ldots,p-1\right\} $.
This F-series is then compatible with respect to the frame which assigns
the real topology to $\mathfrak{z}\in\mathbb{N}_{0}$ and assigns
the $q$-adic topology to $\mathbb{Z}_{p}\backslash\mathbb{N}_{0}$.

Next, letting $\mathbb{C}_{q}$ be the metric completion of the algebraic
closure of $\mathbb{Q}_{q}$, let:
\begin{equation}
U_{q}\overset{\textrm{def}}{=}\left\{ \mathfrak{z}\in\mathbb{C}_{q}:\sum_{n=0}^{\infty}\frac{\mathfrak{z}^{n}}{n!}\textrm{ converges in }\mathbb{C}_{q}\right\} 
\end{equation}
so that $U_{q}$ is the disk of convergence of the $q$-adic exponential
function $\exp_{q}:U_{q}\rightarrow\mathbb{C}_{q}$.

For each $n\in\mathbb{N}_{0}$, $X^{n}\left(\mathfrak{z}\right)d\mathfrak{z}$
is a $\left(p,q\right)$-adic measure, with:
\begin{equation}
\int_{\mathbb{Z}_{p}}\phi\left(\mathfrak{z}\right)X^{n}\left(\mathfrak{z}\right)d\mathfrak{z}=\sum_{t\in\hat{\mathbb{Z}}_{p}}\hat{\phi}\left(t\right)\hat{X}^{*n}\left(-t\right)
\end{equation}
for all continuous $\phi:\mathbb{Z}_{p}\rightarrow\mathbb{C}_{q}$.
Here, the theorem shows that:
\begin{equation}
\left|\int_{\mathbb{Z}_{p}}\phi\left(\mathfrak{z}\right)X^{n}\left(\mathfrak{z}\right)d\mathfrak{z}\right|_{q}\leq\left\Vert \phi\right\Vert _{p,q}\underbrace{\left\Vert \hat{X}^{*n}\right\Vert _{p,q}}_{\leq1}\leq\left\Vert \phi\right\Vert _{p,q}
\end{equation}
As such, letting $\mathfrak{y}\in U_{q}$, the expression:
\begin{equation}
\int_{\mathbb{Z}_{p}}\phi\left(\mathfrak{z}\right)\sum_{n=0}^{\infty}\frac{\mathfrak{y}^{n}}{n!}X^{n}\left(\mathfrak{z}\right)d\mathfrak{z}
\end{equation}
converges absolutely in $\mathbb{C}_{q}$, and therefore defines $\left(p,q\right)$-adic
measure:
\begin{equation}
\int_{\mathbb{Z}_{p}}\phi\left(\mathfrak{z}\right)\exp_{q}\left(\mathfrak{y}X\left(\mathfrak{z}\right)\right)d\mathfrak{z}\overset{\textrm{def}}{=}\sum_{t\in\hat{\mathbb{Z}}_{p}}\hat{\phi}\left(t\right)\sum_{n=0}^{\infty}\frac{\hat{X}^{*n}\left(-t\right)}{n!}\mathfrak{y}^{n}\label{eq:moment generating function of X}
\end{equation}
In this way, the map: 
\begin{align}
U_{q} & \rightarrow W\left(\mathbb{Z}_{p},\mathbb{C}_{q}\right)^{\prime}\\
\mathfrak{y} & \mapsto\exp_{q}\left(\mathfrak{y}X\left(\mathfrak{z}\right)\right)d\mathfrak{z}
\end{align}
is then analytic, and, for fixed $\phi\in W\left(\mathbb{Z}_{p},\mathbb{C}_{q}\right)$,
(\ref{eq:moment generating function of X}) is then an analytic function
$U_{q}\rightarrow\mathbb{C}_{q}$. $\exp_{q}\left(\mathfrak{y}X\left(\mathfrak{z}\right)\right)d\mathfrak{z}$
is the \textbf{moment generating measure} of $X$.

Similarly, in the archimedean case, if we can realize $X$ as a uniformly
continuous function $X:\mathbb{Z}_{p}\rightarrow\mathbb{C}$, we then
have that: 
\begin{equation}
\phi_{X}\left(s\right)\overset{\textrm{def}}{=}\int_{\mathbb{Z}_{p}}e^{sX\left(\mathfrak{z}\right)}d\mathfrak{z}
\end{equation}
is an entire function of $s$, and is the moment generating function
of $X$, because:
\begin{equation}
\int_{\mathbb{Z}_{p}}e^{sX\left(\mathfrak{z}\right)}d\mathfrak{z}=\sum_{n=0}^{\infty}\left(\int_{\mathbb{Z}_{p}}X^{n}\left(\mathfrak{z}\right)d\mathfrak{z}\right)\frac{s^{n}}{n!}=\sum_{n=0}^{\infty}\hat{X}^{*n}\left(0\right)\frac{s^{n}}{n!}
\end{equation}
If $X$ is merely in $L^{1}\left(\mathbb{Z}_{p},\mathbb{C}\right)$
(absolutely integrable with respect to the real-valued Haar probability
measure), the moment generating function need no longer be entire.
For a concrete example, consider:
\begin{align}
X\left(2\mathfrak{z}\right) & =\frac{X\left(\mathfrak{z}\right)}{1+i}\\
X\left(2\mathfrak{z}+1\right) & =X\left(\mathfrak{z}\right)+i
\end{align}
Though the restriction of $X$ to $\mathbb{N}_{0}$ is well-defined,
the extension of $X$ to $\mathbb{Z}_{2}$ is convergent only for
those $\mathfrak{z}$ whose $1$s digits do not have density $1$
among $\mathfrak{z}$'s $2$-adic digits. Here, $\phi_{X}$ satisfies:
\begin{equation}
\phi_{X}\left(s\right)=\frac{1}{2}\int_{\mathbb{Z}_{2}}e^{\frac{s}{1+i}X\left(\mathfrak{z}\right)}d\mathfrak{z}+\frac{1}{2}\int_{\mathbb{Z}_{2}}e^{s\left(X\left(\mathfrak{z}\right)+i\right)}d\mathfrak{z}=\phi_{X}\left(\frac{s}{1+i}\right)+\frac{e^{si}}{2}\phi_{X}\left(s\right)
\end{equation}
from which we obtain:
\begin{equation}
\phi_{X}\left(s\right)=\phi_{X}\left(\frac{s}{\left(1+i\right)^{N}}\right)\prod_{n=0}^{N-1}\left(2-e^{\frac{is}{\left(1+i\right)^{n}}}\right)^{-1},\textrm{ }\forall N\geq1
\end{equation}
It is easy to show that the limit of the partial product as $N\rightarrow\infty$
converges uniformly on all compact subsets of $\mathbb{C}$ avoiding
the points: 
\begin{equation}
s=\left(1+i\right)^{n}\left(2\pi k-i\ln2\right),\textrm{ }\forall n\in\mathbb{N}_{0},\textrm{ }\forall k\in\mathbb{Z}\label{eq:poles-1}
\end{equation}
Using the fact that $\phi_{X}\left(0\right)=1$, it follows that:
\begin{equation}
\phi_{X}\left(s\right)=\prod_{n=0}^{\infty}\left(2-e^{\frac{is}{\left(1+i\right)^{n}}}\right)^{-1}
\end{equation}
is meromorphic on $\mathbb{C}$, as its reciprocal is an entire function.
The poles are all simple, and are located at (\ref{eq:poles-1}).
As $-i\ln2$ is the pole closest to the origin, we then have the asymptotic
equivalence:
\[
\hat{X}^{*n}\left(0\right)\sim C_{0}\left(i\ln2\right)^{n}
\]
for some complex constant $C_{0}$. More generally:
\begin{equation}
\phi_{X}\left(s,t\right)\overset{\textrm{def}}{=}\int_{\mathbb{Z}_{p}}e^{sX\left(\mathfrak{z}\right)}e^{-2\pi i\left\{ t\mathfrak{z}\right\} _{p}}d\mathfrak{z},\textrm{ }\forall t\in\hat{\mathbb{Z}}_{p}
\end{equation}
is the exponential generating function whose $n$th coefficient is
$\hat{X}^{*n}\left(t\right)$:
\begin{equation}
\phi_{X}\left(s,t\right)=\sum_{n=0}^{\infty}\hat{X}^{*n}\left(t\right)\frac{s^{n}}{n!}
\end{equation}
\end{rem}

\section{\label{sec:encoding}Measures, Varieties, \& Arithmetic Dynamics}
\begin{rem}
For everything we do in this section, we can replace each instance
of $\textrm{Frac}\left(\mathcal{R}_{d}\right)\left(\zeta_{p^{\infty}}\right)$
with $\textrm{Frac}\left(\mathcal{R}_{d}/I\mathcal{R}_{d}\right)\left(\zeta_{p^{\infty}}\right)$,
where $I\subseteq R_{d}$ is any unique solution ideal. As such, we
will omit reference to $I$, except when stating major results.
\end{rem}
\textbf{Section \ref{subsec:Formal-Solutions}} establishes the paper's
third main result, \textbf{Theorem \ref{thm:formal solutions}} on
page \pageref{thm:formal solutions}, thereby furnishing a more comprehensive
resolution of \textbf{Questions \ref{que:key question (easy part)}
}and\textbf{ \ref{que:key question (hard part)}} from page \pageref{que:key question (easy part)}
of \textbf{Section \ref{subsec:The-Central-Computation}}.

\textbf{Section \ref{subsec:Varieties-=000026-Measures}} uses the
work of \textbf{Section \ref{subsec:Formal-Solutions}} to illuminate
the ways in which F-series and degenerate measures can be used to
encode information about algebraic varieties. Finally, \textbf{Section
\ref{subsec:Arithmetic-Dynamics,-Varieties,}} gives some examples
and conjectural explorations for future research in using F-series
to better understand certain types of arithmetic dynamical systems.

\subsection{\label{subsec:Formal-Solutions}Formal Solutions \& Breakdown Varieties}

In \textbf{Example \ref{exa:proceeding formally}} from page \textbf{\pageref{exa:proceeding formally}}
of \textbf{Section \ref{subsec:The-Central-Computation}}, we treated
$X_{\mathbf{n}}$ as if it was integrable and assumed that it satisfied
the natural transformation laws under changes of variables induced
by its functional equations. Under those assumptions, we had the following:

\begin{align*}
\hat{X}_{\mathbf{n}}\left(t\right) & =\int_{\mathbb{Z}_{p}}e^{-2\pi i\left\{ t\mathfrak{z}\right\} _{p}}X_{\mathbf{n}}\left(\mathfrak{z}\right)d\mathfrak{z}\\
 & =\sum_{k=0}^{p-1}\int_{\mathbb{Z}_{p}}\left[\mathfrak{z}\overset{p}{\equiv}k\right]e^{-2\pi i\left\{ t\mathfrak{z}\right\} _{p}}X_{\mathbf{n}}\left(\mathfrak{z}\right)d\mathfrak{z}\\
 & =\sum_{k=0}^{p-1}\int_{p\mathbb{Z}_{p}+k}e^{-2\pi i\left\{ t\mathfrak{z}\right\} _{p}}X_{n}\left(\mathfrak{z}\right)d\mathfrak{z}\\
\left(\textrm{change of variable}\right); & =\frac{1}{p}\sum_{k=0}^{p-1}\int_{\mathbb{Z}_{p}}e^{-2\pi i\left\{ t\left(p\mathfrak{z}+k\right)\right\} _{p}}X_{\mathbf{n}}\left(p\mathfrak{z}+k\right)d\mathfrak{z}\\
\left(\textrm{apply }(\ref{eq:Z functional equation})\right); & =\frac{1}{p}\sum_{k=0}^{p-1}\int_{\mathbb{Z}_{p}}e^{-2\pi i\left\{ pt\mathfrak{z}\right\} _{p}}e^{-2\pi ikt}\left(\sum_{\mathbf{m}\leq\mathbf{n}}r_{\mathbf{m},\mathbf{n},k}X_{\mathbf{m}}\left(\mathfrak{z}\right)\right)d\mathfrak{z}\\
 & =\sum_{\mathbf{m}\leq\mathbf{e}}\underbrace{\left(\frac{1}{p}\sum_{k=0}^{p-1}r_{\mathbf{m},\mathbf{n},k}e^{-2\pi ikt}\right)}_{\alpha_{\mathbf{m},\mathbf{n}}\left(t\right)}\underbrace{\int_{\mathbb{Z}_{p}}e^{-2\pi i\left\{ pt\mathfrak{z}\right\} _{p}}X_{\mathbf{m}}\left(\mathfrak{z}\right)d\mathfrak{z}}_{\hat{X}_{\mathbf{m}}\left(pt\right)}\\
 & =\sum_{\mathbf{m}\leq\mathbf{n}}\alpha_{\mathbf{m},\mathbf{n}}\left(t\right)\hat{X}_{\mathbf{m}}\left(pt\right)
\end{align*}
We showed in \textbf{Proposition \ref{prop:X_3-hat formal solution}}
that the functional equation:
\begin{equation}
\hat{X}_{\mathbf{n}}\left(t\right)=\sum_{\mathbf{m}\leq\mathbf{n}}\alpha_{\mathbf{m},\mathbf{n}}\left(t\right)\hat{X}_{\mathbf{m}}\left(pt\right)\label{eq:alpha formal equation}
\end{equation}
could be used to recursively solve for $\hat{X}_{\mathbf{n}}$ uniquely
provided $\alpha_{\mathbf{n}}\left(0\right)\neq1$.

By using truncations (precomposing $X_{\mathbf{n}}$ with the projection
map $\mathbb{Z}_{p}\rightarrow\mathbb{Z}/p^{N}\mathbb{Z}$, we saw
a second way to obtain this functional equation, namely, as the formal
limit as $N\rightarrow\infty$ of the functional equations satisfied
by the Fourier transform of $X_{\mathbf{n},N}$, the $N$th truncation
of $X_{\mathbf{n}}$:
\begin{equation}
\hat{X}_{\mathbf{n},N}\left(t\right)=\sum_{\mathbf{m}\leq\mathbf{n}}\alpha_{\mathbf{m},\mathbf{n}}\left(t\right)\hat{X}_{\mathbf{m},N-1}\left(pt\right)
\end{equation}
If the Fourier transforms of the $\hat{X}_{\mathbf{m}}$s exist in
the classical sense\textemdash say, there is a metrically complete,
algebraically closed valued field $\mathbb{K}$ with characteristic
either co-prime to $p$ (if the characteristic is positive) and/or
with a residue field of characteristic co-prime to $p$ (if $\mathbb{K}$
is non-archimedean) so that $X_{\mathbf{m}}\in W\left(\mathbb{Z}_{p},\mathbb{K}\right)$
for all $\mathbf{m}\leq\mathbf{n}$\textemdash then the limit $\lim_{N\rightarrow\infty}\hat{X}_{\mathbf{m},N}\left(t\right)$
converges in $\mathbb{K}$ to $\hat{X}_{\mathbf{m}}\left(t\right)$,
the Fourier integral of $X_{\mathbf{m}}$:
\begin{equation}
\lim_{N\rightarrow\infty}\hat{X}_{\mathbf{m},N}\left(t\right)\overset{\mathbb{K}}{=}\int_{\mathbb{Z}_{p}}X_{\mathbf{m}}\left(\mathfrak{z}\right)e^{-2\pi i\left\{ t\mathfrak{z}\right\} _{p}}d\mathfrak{z},\textrm{ }\forall\mathbf{m}\leq\mathbf{n}
\end{equation}
where $d\mathfrak{z}$ is the $\mathbb{K}$-valued Haar probability
measure on $\mathbb{Z}_{p}$. Here, the existence of the integral
follows by classical abstract harmonic analysis and Pontryagin duality.
Moreover, since $X_{\mathbf{m}}$ is in the Wiener algebra, the Fourier
series generated by $\hat{X}_{\mathbf{m}}$ will sum to $X_{\mathbf{m}}$
for all $\mathbf{m}\leq\mathbf{n}$.

Thus, in this case, provided $\alpha_{\mathbf{n}}\left(0\right)\neq1$
and all the $\hat{X}_{\mathbf{m}}$s are already known, the equation
(\ref{eq:alpha formal equation}) can be solved recursively to obtain
a closed-form expression for the Fourier transform of $X_{\mathbf{n}}$.
This is much nicer than the methods we've used so far, because it
involves no guesswork and can be done completely systematically. The
purpose of this subsection is to show that, even when $X_{\mathbf{n}}$
is merely quasi-integrable, so long as $\alpha_{\mathbf{n}}\left(0\right)\neq1$,
this same formal method of using (\ref{eq:alpha formal equation})
to compute $\hat{X}_{\mathbf{n}}$ does, in fact, ``work'', in the
sense that the answer obtained is indeed a Fourier transform of $\hat{X}_{\mathbf{n}}$
with respect to an appropriately chosen frame.

We begin with some definitions.
\begin{defn}
Viewing $\mathbb{N}_{0}^{d}$ as a partially ordered set with respect
to the relations $\leq$ and $<$, given $\mathbf{m},\mathbf{n}$
with $\mathbf{m}<\mathbf{n}$, we write $\left[\mathbf{m},\mathbf{n}\right]$
to denote the set of all $k\in\mathbb{N}_{0}^{d}$ so that $\mathbf{m}\leq\mathbf{k}\leq\mathbf{n}$
and call $\left[\mathbf{m},\mathbf{n}\right]$ the \textbf{closed
interval from $\mathbf{m}$ to $n$}. We similarly co-opt standard
interval notation to denote open intervals and half-open intervals
in the obvious way (Ex: $\left[\mathbf{m},\mathbf{n}\right)$ is the
set of all $\mathbf{k}$ for which $\mathbf{m}\leq\mathbf{k}<\mathbf{n}$).
We say an interval $I\subseteq\mathbb{N}_{0}^{d}$ is \textbf{bounded
}whenever $\sup_{\mathbf{k}\in I}\Sigma\left(\mathbf{k}\right)<\infty$.

Let $\textrm{Map}\left(\hat{\mathbb{Z}}_{p},\mathcal{R}_{d}\left(\zeta_{p^{\infty}}\right)\right)$
denote the space of all functions $\hat{\mathbb{Z}}_{p}\rightarrow\mathcal{R}_{d}\left(\zeta_{p^{\infty}}\right)$.
Given an $\mathbf{n}\in\mathbb{N}_{0}^{d}$, the \textbf{$\mathbf{n}$th
formal equation }/ \textbf{$\mathbf{n}$th functional equation }is:
\begin{equation}
\hat{Y}_{\mathbf{n}}\left(t\right)=\sum_{\mathbf{m}\in\left[\mathbf{0},\mathbf{n}\right]\cap I}\alpha_{\mathbf{m},\mathbf{n}}\left(t\right)\hat{Y}_{\mathbf{m}}\left(pt\right)\label{eq:nth formal equation}
\end{equation}
for undetermined functions $\left\{ \hat{Y}_{\mathbf{m}}\right\} _{\mathbf{m}\leq\mathbf{n}}$
in $\textrm{Map}\left(\hat{\mathbb{Z}}_{p},\mathcal{R}_{d}\left(\zeta_{p^{\infty}}\right)\right)$.
Given a bounded interval $I\subseteq\mathbb{N}_{0}^{d}$, the \textbf{$I$th
formal system / $I$th functional equation }is the system of equations:
\begin{equation}
\left\{ \hat{Y}_{\mathbf{n}}\left(t\right)=\sum_{\mathbf{m}\in\left[\mathbf{0},\mathbf{n}\right]\cap I}\alpha_{\mathbf{m},\mathbf{n}}\left(t\right)\hat{Y}_{\mathbf{m}}\left(pt\right),\textrm{ }\forall\mathbf{n}\in I\right\} \label{eq:Ith formal system}
\end{equation}
for undetermined functions $\left\{ \hat{Y}_{\mathbf{n}}\right\} _{\mathbf{n}\in I}$
in $\textrm{Map}\left(\hat{\mathbb{Z}}_{p},\mathcal{R}_{d}\left(\zeta_{p^{\infty}}\right)\right)$.
Note that for any $\mathbf{n}$, the $\left[\mathbf{n},\mathbf{n}\right]$th
formal system is then precisely the $\mathbf{n}$th formal equation.

Given an interval $I$, an \textbf{$I$th formal solution} is a collection
of functions $\left\{ \hat{\chi}_{\mathbf{n}}\right\} _{\mathbf{n}\in I}$
in $\textrm{Map}\left(\hat{\mathbb{Z}}_{p},\mathcal{R}_{d}\left(\zeta_{p^{\infty}}\right)\right)$
so that setting $\hat{Y}_{\mathbf{n}}=\hat{\chi}_{\mathbf{n}}$ for
all $\mathbf{n}\in I$ makes (\ref{eq:Ith formal system}) hold true
for all $t\in\hat{\mathbb{Z}}_{p}$.

Given an interval $J\subseteq\mathbb{N}_{0}^{d}$, a \textbf{$J$th
initial condition} is a collection of functions $\left\{ \hat{\chi}_{\mathbf{n}}\right\} _{\mathbf{n}\in J}$
in $\textrm{Map}\left(\hat{\mathbb{Z}}_{p},\mathcal{R}_{d}\left(\zeta_{p^{\infty}}\right)\right)$.
The \textbf{$I$th formal system (}a.k.a. \textbf{$I$th functional
equation) relative to an initial condition }$\left\{ \hat{\chi}_{\mathbf{n}}\right\} _{\mathbf{n}\in J}$\textbf{
}is:
\begin{equation}
\left\{ \hat{Y}_{\mathbf{n}}\left(t\right)=\sum_{\mathbf{m}\in\left[\mathbf{0},\mathbf{n}\right]\cap\left(I\backslash J\right)}\alpha_{\mathbf{m},\mathbf{n}}\left(t\right)\hat{Y}_{\mathbf{m}}\left(pt\right)+\sum_{\mathbf{m}\in\left[\mathbf{0},\mathbf{n}\right]\cap J}\alpha_{\mathbf{m},\mathbf{n}}\left(t\right)\hat{\chi}_{\mathbf{m}}\left(pt\right),\textrm{ }\forall\mathbf{n}\in I\right\} \label{eq:relative system}
\end{equation}
for undetermined functions $\left\{ \hat{Y}_{\mathbf{n}}\right\} _{\mathbf{n}\in I\backslash J}$
in $\textrm{Map}\left(\hat{\mathbb{Z}}_{p},\mathcal{R}_{d}\left(\zeta_{p^{\infty}}\right)\right)$.
An \textbf{$I$th formal solution relative to an initial condition
}$\left\{ \hat{\chi}_{\mathbf{n}}\right\} _{\mathbf{n}\in J}$\textbf{
}is a collection of functions $\left\{ \hat{\chi}_{\mathbf{n}}\right\} _{\mathbf{n}\in I\backslash J}$
in $\textrm{Map}\left(\hat{\mathbb{Z}}_{p},\mathcal{R}_{d}\left(\zeta_{p^{\infty}}\right)\right)$
so that $\hat{Y}_{\mathbf{n}}=\hat{\chi}_{\mathbf{n}}$ for all $\mathbf{n}\in I\backslash J$
satisfies (\ref{eq:relative system}) for all $t$.

For all of the notions presented so far, given a unique solution ideal
$I\subseteq R_{d}$, we will add the qualifier \textbf{mod $I$ }to
indicate that the equalities occur in $\left(\mathcal{R}_{d}/I\mathcal{R}_{d}\right)\left(\zeta_{p^{\infty}}\right)$.
\end{defn}
Now, as we observed in \textbf{Proposition \ref{prop:X_3-hat formal solution}}
from \textbf{Section \ref{subsec:The-Central-Computation}}, given
any $\mathbf{n}>\mathbf{0}$ and any initial condition $\left\{ \hat{\chi}_{\mathbf{m}}\right\} _{\mathbf{m}\in\left[\mathbf{0},\mathbf{n}\right)}$,
the $\mathbf{n}$th functional equation relative to this initial condition
has a unique solution if and only if $\alpha_{\mathbf{n}}\left(0\right)\neq1$.
When $\alpha_{\mathbf{n}}\left(0\right)=1$, the equation breaks down,
and we have one of two possibilities: either:
\begin{equation}
\sum_{\mathbf{m}<\mathbf{n}}\alpha_{\mathbf{m},\mathbf{n}}\left(0\right)\hat{\chi}_{\mathbf{m}}\left(0\right)=0
\end{equation}
in which case there are infinitely many $\mathbf{n}$th solutions
$\hat{Y}_{\mathbf{n}}\left(t\right)$, one for each initial value
$\hat{Y}_{\mathbf{n}}\left(0\right)\in\mathcal{R}_{d}\left(\zeta_{p^{\infty}}\right)$,
or the above equality fails to hold, in which case there is no $\mathbf{n}$th
solution.

This motivates the following definition:
\begin{defn}
Fix $d\geq1$ and a unique solution ideal $I\subseteq R_{d}$. Given
any $\mathbf{n}\in\mathbb{N}_{0}^{d}\backslash\left\{ \mathbf{0}\right\} $,
we write $\textrm{Break}_{I}\left(\mathbf{n}\right)$ to denote the
locus of $a_{j,k}$s and $b_{j,k}$s in $K$ so that $\alpha_{\mathbf{n}}\left(0\right)=1$
holds true in $\mathcal{R}_{d}/I\mathcal{R}_{d}$. We call $\textrm{Break}_{I}\left(\mathbf{n}\right)$
the \textbf{$\mathbf{n}$th breakdown variety relative to }$I$. We
also call $\textrm{Break}_{I}\left(\mathbf{n}\right)$ the \textbf{breakdown
variety of $X_{\mathbf{n}}$ relative to $I$},\textbf{ }and write
it as $\textrm{Break}_{I}\left(X_{\mathbf{n}}\right)$

We then write $\textrm{CoBreak}_{I}\left(\mathbf{n}\right)$ to denote
the locus of all $a_{j,k}$s and $b_{j,k}$s in $K$ so that the equation
$\alpha_{\mathbf{n}}\left(0\right)\neq1$ holds true in $\mathcal{R}_{d}/I\mathcal{R}_{d}$.
We call $\textrm{CoBreak}_{I}\left(\mathbf{n}\right)$ the \textbf{$\mathbf{n}$th
co-breakdown variety relative to }$I$. We also call $\textrm{CoBreak}_{I}\left(\mathbf{n}\right)$
the \textbf{co-breakdown variety of $X_{\mathbf{n}}$ relative to
$I$}, and write it as $\textrm{CoBreak}_{I}\left(X_{\mathbf{n}}\right)$

Next, we write: 
\begin{equation}
\textrm{BI}_{d}\left(I\right)\overset{\textrm{def}}{=}\left\{ \mathbf{n}\in\mathbb{N}_{0}^{d}\backslash\left\{ \mathbf{0}\right\} :\left\langle 1-\alpha_{\mathbf{n}}\left(0\right)\right\rangle \subseteq I\right\} 
\end{equation}
and: 
\begin{equation}
\textrm{CBI}_{d}\left(I\right)\overset{\textrm{def}}{=}\left\{ \mathbf{n}\in\mathbb{N}_{0}^{d}\backslash\left\{ \mathbf{0}\right\} :\left\langle 1-\alpha_{\mathbf{n}}\left(0\right)\right\rangle \nsubseteq I\right\} 
\end{equation}
We call these the sets of \textbf{breakdown indices }and \textbf{co-breakdown
indices }associated to $I$.
\end{defn}
\begin{example}
When $d=1$, the $n\in\mathbb{N}_{0}^{1}\backslash\left\{ 0\right\} $th
breakdown variety relative to $\left\langle 0\right\rangle $ (a.k.a.,
the $n$th breakdown variety of $X$ relative to $0$) is the locus
of all $a_{0},\ldots,a_{p-1}\in K$ so that:
\begin{equation}
\sum_{k=0}^{p-1}a_{k}^{n}=p
\end{equation}
Thus, for example, the breakdown variety of $X^{2}$ relative to $\left\langle 0\right\rangle $
is the locus of points in $K$ lying on the $\left(p-1\right)$-sphere
of radius $p$.
\end{example}
\begin{defn}
We \textbf{breakdown indices }$\textrm{BI}_{d}\left(I\right)$ as
the set of all $\mathbf{n}\in\mathbb{N}_{0}^{d}$ with $\Sigma\left(\mathbf{n}\right)\geq2$
for which $\left\langle 1-\alpha_{\mathbf{n}}\left(0\right)\right\rangle \subseteq I$.
Let $\textrm{CBI}_{d}\left(I\right)$, the \textbf{co-breakdown indices
}be the set of all $\mathbf{n}\in\mathbb{N}_{0}^{d}$ with $\Sigma\left(\mathbf{n}\right)\geq2$
for which $\left\langle 1-\alpha_{\mathbf{n}}\left(0\right)\right\rangle \notin I$.
We $\mathbf{n}\in\mathbb{N}_{0}^{d}$ is \textbf{trivial }when $\Sigma\left(\mathbf{n}\right)\leq1$.
\end{defn}
Relative to an initial condition $\left\{ \hat{\chi}_{\mathbf{m}}\right\} _{\mathbf{m}\in\left[\mathbf{0},\mathbf{n}\right)}$,
the $\mathbf{n}$th functional equation can be written as:
\begin{equation}
\hat{Y}_{\mathbf{n}}\left(t\right)-\alpha_{\mathbf{n}}\left(t\right)\hat{Y}_{\mathbf{n}}\left(pt\right)=\sum_{\mathbf{m}<\mathbf{n}}\alpha_{\mathbf{m},\mathbf{n}}\left(t\right)\hat{\chi}_{\mathbf{m}}\left(pt\right)\label{eq:nth conditioned}
\end{equation}
Letting $dY_{\mathbf{n}}$ and $d\chi_{\mathbf{m}}$ denote the distributions
with $\hat{Y}_{\mathbf{n}}$ and $\hat{\chi}_{\mathbf{m}}$ as their
Fourier transforms, we can use (\ref{eq:nth conditioned}) to compute
a functional equation giving the action of $dY_{\mathbf{n}}$ on an
SB function $\phi$ in terms of the actions of the $d\chi_{\mathbf{m}}$s.
\begin{prop}
\label{prop:distribution action of nth solution}Let $\hat{Y}_{\mathbf{n}}$
be any solution of the $\mathbf{n}$th functional equation relative
to an initial condition $\left\{ \hat{\chi}_{\mathbf{m}}\right\} _{\mathbf{m}\in\left[0,\mathbf{n}\right)}$,
and let $dY_{\mathbf{n}}$ and $d\chi_{\mathbf{m}}$ be the distributions
on $\mathcal{S}\left(\mathbb{Z}_{p},\mathcal{R}_{d}\left(\zeta_{p^{\infty}}\right)\right)$
with $\hat{Y}_{\mathbf{n}}$ and the $\hat{\chi}_{\mathbf{m}}$s as
their Fourier-Stieltjes transforms. Then, for all $\phi\in\mathcal{S}\left(\mathbb{Z}_{p},\mathcal{R}_{d}\left(\zeta_{p^{\infty}}\right)\right)$:
\begin{equation}
\int_{\mathbb{Z}_{p}}\left(\phi\left(\mathfrak{z}\right)-\frac{1}{p}\sum_{k=0}^{p-1}r_{\mathbf{n},k}\phi\left(p\mathfrak{z}+k\right)\right)dY_{\mathbf{n}}\left(\mathfrak{z}\right)=\int_{\mathbb{Z}_{p}}\left(\frac{1}{p}\sum_{k=0}^{p-1}r_{\mathbf{m},\mathbf{n},k}\phi\left(p\mathfrak{z}+k\right)\right)d\chi_{\mathbf{m}}\left(\mathfrak{z}\right)\label{eq:distribution action of nth solution}
\end{equation}
\end{prop}
Proof: Multiplying (\ref{eq:nth conditioned}) by $\hat{\phi}\left(-t\right)$
and summing over all $t\in\hat{\mathbb{Z}}_{p}$, where $\phi\in\mathcal{S}\left(\mathbb{Z}_{p},\mathcal{R}_{d}\left(\zeta_{p^{\infty}}\right)\right)$,
and using:
\begin{equation}
\sum_{t\in\hat{\mathbb{Z}}_{p}}\hat{\phi}\left(-t\right)\hat{\mu}\left(t\right)=\int_{\mathbb{Z}_{p}}\phi\left(\mathfrak{z}\right)d\mu\left(\mathfrak{z}\right),\textrm{ }\forall d\mu\in\mathcal{S}\left(\mathbb{Z}_{p},\textrm{Frac}\left(\mathcal{R}_{d}\right)\left(\zeta_{p^{\infty}}\right)\right)^{\prime}
\end{equation}
we get: 
\begin{equation}
\int_{\mathbb{Z}_{p}}\phi\left(\mathfrak{z}\right)dY_{\mathbf{n}}\left(\mathfrak{z}\right)-\sum_{t\in\hat{\mathbb{Z}}_{p}}\hat{\phi}\left(-t\right)\alpha_{\mathbf{n}}\left(t\right)\hat{Y}_{\mathbf{n}}\left(pt\right)=\sum_{\mathbf{m}<\mathbf{n}}\sum_{t\in\hat{\mathbb{Z}}_{p}}\hat{\phi}\left(-t\right)\alpha_{\mathbf{m},\mathbf{n}}\left(t\right)\hat{\chi}_{\mathbf{m}}\left(pt\right)
\end{equation}
Here, using the Fourier-Stieltjes transform of $dY_{\mathbf{n}}$,
we get:
\begin{align*}
\sum_{t\in\hat{\mathbb{Z}}_{p}}\hat{\phi}\left(-t\right)\alpha_{\mathbf{n}}\left(t\right)\hat{Y}_{\mathbf{n}}\left(pt\right) & =\int_{\mathbb{Z}_{p}}\left(\sum_{t\in\hat{\mathbb{Z}}_{p}}\hat{\phi}\left(-t\right)\alpha_{\mathbf{n}}\left(t\right)e^{-2\pi i\left\{ pt\mathfrak{z}\right\} _{p}}\right)dY_{\mathbf{n}}\left(\mathfrak{z}\right)\\
\left(\alpha_{\mathbf{n}}\left(t\right)=\frac{1}{p}\sum_{k=0}^{p-1}r_{\mathbf{n},k}e^{-2\pi ikt}\right); & =\int_{\mathbb{Z}_{p}}\left(\frac{1}{p}\sum_{k=0}^{p-1}r_{\mathbf{n},k}\sum_{t\in\hat{\mathbb{Z}}_{p}}\hat{\phi}\left(-t\right)e^{2\pi i\left\{ -t\left(p\mathfrak{z}+k\right)\right\} _{p}}\right)dY_{\mathbf{n}}\left(\mathfrak{z}\right)\\
 & =\int_{\mathbb{Z}_{p}}\left(\frac{1}{p}\sum_{k=0}^{p-1}r_{\mathbf{n},k}\phi\left(p\mathfrak{z}+k\right)\right)dY_{\mathbf{n}}\left(\mathfrak{z}\right)
\end{align*}
Likewise, for $d\chi_{\mathbf{m}}$:
\begin{equation}
\sum_{t\in\hat{\mathbb{Z}}_{p}}\hat{\phi}\left(-t\right)\alpha_{\mathbf{m},\mathbf{n}}\left(t\right)\hat{\chi}_{\mathbf{m}}\left(pt\right)=\int_{\mathbb{Z}_{p}}\left(\frac{1}{p}\sum_{k=0}^{p-1}r_{\mathbf{m},\mathbf{n},k}\phi\left(p\mathfrak{z}+k\right)\right)d\chi_{\mathbf{m}}\left(\mathfrak{z}\right)
\end{equation}
Putting all this together, (\ref{eq:nth conditioned}) becomes: 
\begin{equation}
\int_{\mathbb{Z}_{p}}\left(\phi\left(\mathfrak{z}\right)-\frac{1}{p}\sum_{k=0}^{p-1}r_{\mathbf{n},k}\phi\left(p\mathfrak{z}+k\right)\right)dY_{\mathbf{n}}\left(\mathfrak{z}\right)=\int_{\mathbb{Z}_{p}}\sum_{\mathbf{m}<\mathbf{n}}\left(\frac{1}{p}\sum_{k=0}^{p-1}r_{\mathbf{m},\mathbf{n},k}\phi\left(p\mathfrak{z}+k\right)\right)d\chi_{\mathbf{m}}\left(\mathfrak{z}\right)
\end{equation}

Q.E.D.

\vphantom{}The formula that we got motivates the following definition:
\begin{defn}
For any $\mathbf{m},\mathbf{n}\in\mathbb{N}_{0}^{d}$ with $\mathbf{m}\leq\mathbf{n}$,
define $L_{\mathbf{m},\mathbf{n}}:\mathcal{S}\left(\mathbb{Z}_{p},\textrm{Frac}\left(\mathcal{R}_{d}\right)\right)\rightarrow\mathcal{S}\left(\mathbb{Z}_{p},\textrm{Frac}\left(\mathcal{R}_{d}\right)\right)$
by:
\begin{equation}
L_{\mathbf{m},\mathbf{n}}\left\{ \phi\right\} \left(\mathfrak{z}\right)\overset{\textrm{def}}{=}\frac{1}{p}\sum_{k=0}^{p-1}r_{\mathbf{m},\mathbf{n},k}\phi\left(p\mathfrak{z}+k\right),\textrm{ }\forall\phi\in\mathcal{S}\left(\mathbb{Z}_{p},\textrm{Frac}\left(\mathcal{R}_{d}\right)\right)
\end{equation}
We then write $L_{\mathbf{n}}$ to denote $L_{\mathbf{n},\mathbf{n}}$.
\end{defn}
The following is then obvious:
\begin{prop}
\label{prop:properties of L_m,n}$L_{\mathbf{m},\mathbf{n}}$ is a
linear operator. Moreover, letting $v$ be any absolute value on $\mathcal{R}_{d}$,
and letting $\overline{\overline{\mathbb{K}_{v}}}$ be the metric
completion of the algebraic closure of the metric completion $\textrm{Frac}\left(\mathcal{R}_{d}\right)$
with respect to $v$, we have that $L_{\mathbf{m},\mathbf{n}}$ uniquely
extends to a \emph{continuous }linear operator $L:W\left(\mathbb{Z}_{p},\overline{\overline{\mathbb{K}_{v}}}\right)\rightarrow W\left(\mathbb{Z}_{p},\overline{\overline{\mathbb{K}_{v}}}\right)$.
We have the operator norm of $L_{\mathbf{m},\mathbf{n}}$ over $W\left(\mathbb{Z}_{p},\overline{\overline{\mathbb{K}_{v}}}\right)$
satisfies the bound:
\begin{equation}
\left\Vert L_{\mathbf{m},\mathbf{n}}\right\Vert _{p,v}\leq\begin{cases}
\max_{0\leq k<p}\left|\frac{r_{\mathbf{m},\mathbf{n},k}}{p}\right|_{v} & \textrm{if }v<\infty\\
\sum_{k=0}^{p-1}\left|\frac{r_{\mathbf{m},\mathbf{n},k}}{p}\right|_{v} & \textrm{else}
\end{cases}
\end{equation}
where $v<\infty$ means $v$ is a non-archimedean place.
\end{prop}
Proof: Triangle inequality in the archimedean case; ultrametric inequality
in the non-archimedean case.

Q.E.D.

\vphantom{}With the $L_{\mathbf{m},\mathbf{n}}$s in hand, we can
write the result of \textbf{Proposition \ref{prop:distribution action of nth solution}}
as:
\begin{equation}
\int_{\mathbb{Z}_{p}}\left(1-L_{\mathbf{n}}\right)\left\{ \phi\right\} \left(\mathfrak{z}\right)dY_{\mathbf{n}}\left(\mathfrak{z}\right)=\int_{\mathbb{Z}_{p}}\sum_{\mathbf{m}<\mathbf{n}}L_{\mathbf{m},\mathbf{n}}\left\{ \phi\right\} \left(\mathfrak{z}\right)d\chi_{\mathbf{m}}\left(\mathfrak{z}\right)\label{eq:1 - L_n first equation}
\end{equation}
The left-hand side is the image of $\phi$ under the composite $dY_{\mathbf{n}}\circ\left(1-L_{\mathbf{n}}\right):\mathcal{S}\left(\mathbb{Z}_{p},\textrm{Frac}\left(\mathcal{R}_{d}\right)\right)\rightarrow\textrm{Frac}\left(\mathcal{R}_{d}\right)$.
As such, observe that if $1-L_{\mathbf{n}}:\mathcal{S}\left(\mathbb{Z}_{p},\textrm{Frac}\left(\mathcal{R}_{d}\right)\right)\rightarrow\mathcal{S}\left(\mathbb{Z}_{p},\textrm{Frac}\left(\mathcal{R}_{d}\right)\right)$
was \emph{invertible}, by setting $\phi$ equal to $\left(1-L_{\mathbf{n}}\right)^{-1}\left\{ \psi\right\} $,
where $\psi\in\mathcal{S}\left(\mathbb{Z}_{p},\textrm{Frac}\left(\mathcal{R}_{d}\right)\right)$
was arbitrary, we would get:
\begin{equation}
\int_{\mathbb{Z}_{p}}\psi\left(\mathfrak{z}\right)dY_{\mathbf{n}}\left(\mathfrak{z}\right)=\int_{\mathbb{Z}_{p}}\sum_{\mathbf{m}<\mathbf{n}}\left(L_{\mathbf{m},\mathbf{n}}\circ\left(1-L_{\mathbf{n}}\right)^{-1}\right)\left\{ \psi\right\} \left(\mathfrak{z}\right)d\chi_{\mathbf{m}}\left(\mathfrak{z}\right),\textrm{ }\forall\psi\in\mathcal{S}\left(\mathbb{Z}_{p},\textrm{Frac}\left(\mathcal{R}_{d}\right)\right)\label{eq:inverting dY_n}
\end{equation}
The right-hand side is then a closed-form expression of $dY_{\mathbf{n}}$
in terms of the $d\chi_{\mathbf{m}}$s. In particular, we see that:
\begin{prop}
\label{prop:The-th-formal}The $\mathbf{n}$th formal solution $dY_{\mathbf{n}}$
relative to the initial condition $\left\{ \hat{\chi}_{\mathbf{m}}\right\} _{\mathbf{m}\in\left[0,\mathbf{n}\right)}$
is uniquely determined by the initial condition whenever $1-L_{\mathbf{n}}$
is invertible on $\mathcal{S}\left(\mathbb{Z}_{p},\textrm{Frac}\left(\mathcal{R}_{d}\right)\right)$.
\end{prop}
Unsurprisingly, it turns out that the invertibility holds if and only
if $\alpha_{\mathbf{n}}\left(0\right)\neq1$. To show this, we first
compute a closed-form expression for the image of the indicator function
$\left[\mathfrak{z}\overset{p^{n}}{\equiv}k\right]$ under $m$ iterations
of $1-L_{\mathbf{m},\mathbf{n}}$.
\begin{prop}
\label{prop:1 minus L_r in terms of shift}Let $\mathbf{n}\in\mathbb{N}_{0}^{d}$,
and let $m,n\geq0$, and $k\in\left\{ 0,\ldots,p^{n}-1\right\} $.
Then:
\begin{equation}
L_{\mathbf{n}}^{m}\left\{ \left[\cdot\overset{p^{n}}{\equiv}k\right]\right\} \left(\mathfrak{z}\right)=\begin{cases}
\frac{r_{\mathbf{n},0}^{m}\kappa_{\mathbf{n}}\left(\left[k\right]_{p^{m}}\right)}{p^{m}}\left[\mathfrak{z}\overset{p^{n-m}}{\equiv}\theta_{p}^{\circ m}\left(k\right)\right] & \textrm{if }m<n\\
\frac{r_{\mathbf{n},0}^{n}\kappa_{\mathbf{n}}\left(k\right)}{p^{n}}\left(\frac{1}{p}\sum_{j=0}^{p-1}r_{\mathbf{n},j}\right)^{m-n} & \textrm{if }m\geq n
\end{cases}\label{eq:1 minus L_r in terms of shift}
\end{equation}
where
\begin{equation}
\kappa_{\mathbf{n}}\left(n\right)=\prod_{j=1}^{p-1}\left(\frac{r_{\mathbf{n},j}}{r_{\mathbf{n},0}}\right)^{\#_{p:j}\left(n\right)}\label{eq:definition of kappa 3}
\end{equation}
\end{prop}
Proof: For brevity, we will omit the subscript $\mathbf{n}$s. Using
the definition of $L$, we have:
\begin{equation}
L\left\{ \left[\cdot\overset{p^{n}}{\equiv}k\right]\right\} \left(\mathfrak{z}\right)=\frac{1}{p}\sum_{j=0}^{p-1}r_{j}\left[p\mathfrak{z}+j\overset{p^{n}}{\equiv}k\right]
\end{equation}
By induction, we obtain
\begin{equation}
L\left\{ \left[\cdot\overset{p^{n}}{\equiv}k\right]\right\} \left(\mathfrak{z}\right)=\left(\prod_{j=0}^{m-1}\frac{r_{\left[\theta_{p}^{\circ j}\left(k\right)\right]_{p}}}{p}\right)\left[\mathfrak{z}\overset{p^{n-m}}{\equiv}\theta_{p}^{\circ m}\left(k\right)\right]
\end{equation}
By \textbf{Proposition \ref{prop:1 minus L_r in terms of shift}},
we have:

\begin{equation}
\prod_{j=0}^{m-1}\frac{r_{\left[\theta_{p}^{\circ j}\left(k\right)\right]_{p}}}{p}=\frac{1}{p^{m}}r_{0}^{m}\prod_{j=1}^{p-1}\left(r_{j}/r_{0}\right)^{\#_{p:j}\left(\left[k\right]_{p^{m}}\right)}=\frac{r_{0}^{m}\kappa\left(\left[k\right]_{p^{m}}\right)}{p^{m}}
\end{equation}
Hence:
\begin{equation}
L^{m}\left\{ \left[\cdot\overset{p^{n}}{\equiv}k\right]\right\} \left(\mathfrak{z}\right)=\left(\prod_{j=0}^{m-1}\frac{r_{\left[\theta_{p}^{\circ j}\left(k\right)\right]_{p}}}{p}\right)\left[\mathfrak{z}\overset{p^{n-m}}{\equiv}\theta_{p}^{\circ m}\left(k\right)\right]=\frac{r_{0}^{m}\kappa\left(\left[k\right]_{p^{m}}\right)}{p^{m}}\left[\mathfrak{z}\overset{p^{n-m}}{\equiv}\theta_{p}^{\circ m}\left(k\right)\right]
\end{equation}

Now, note that when $m=n$, we have: 
\begin{equation}
L^{n}\left\{ \left[\cdot\overset{p^{n}}{\equiv}k\right]\right\} \left(\mathfrak{z}\right)=\frac{r_{0}^{n}\kappa\left(\left[k\right]_{p^{n}}\right)}{p^{n}}\left[\mathfrak{z}\overset{1}{\equiv}\theta_{p}^{\circ n}\left(k\right)\right]
\end{equation}
Since $\mathfrak{z}$ and $\theta_{p}^{\circ n}\left(k\right)$ are
both $p$-adic integers, the congruence $\mathfrak{z}\overset{1}{\equiv}\theta_{p}^{\circ n}\left(k\right)$
is always true, and so: 
\begin{align*}
L^{n}\left\{ \left[\cdot\overset{p^{n}}{\equiv}k\right]\right\} \left(\mathfrak{z}\right) & =\frac{r_{0}^{n}\kappa\left(\left[k\right]_{p^{n}}\right)}{p^{n}}\\
\left(0\leq k\leq p^{n}-1\right); & =\frac{r_{0}^{n}\kappa\left(k\right)}{p^{n}}
\end{align*}
is a constant function. Also, for the constant function $1$, we note
that:
\begin{equation}
L\left\{ 1\right\} \left(\mathfrak{z}\right)=\frac{1}{p}\sum_{j=0}^{p-1}r_{j}
\end{equation}
and so:
\begin{equation}
L^{m}\left\{ 1\right\} \left(\mathfrak{z}\right)=\left(\frac{1}{p}\sum_{j=0}^{p-1}r_{j}\right)^{m},\textrm{ }\forall m\geq1
\end{equation}
Consequently, when $m\geq n$, we have:
\begin{align*}
L^{m}\left\{ \left[\cdot\overset{p^{n}}{\equiv}k\right]\right\} \left(\mathfrak{z}\right) & =L^{m-n}\left\{ L^{n}\left\{ \left[\cdot\overset{p^{n}}{\equiv}k\right]\right\} \right\} \left(\mathfrak{z}\right)\\
 & =L^{m-n}\left\{ \frac{r_{0}^{n}\kappa\left(k\right)}{p^{n}}\right\} \left(\mathfrak{z}\right)\\
 & =\frac{r_{0}^{n}\kappa\left(k\right)}{p^{n}}\left(\frac{1}{p}\sum_{j=0}^{p-1}r_{j}\right)^{m-n}
\end{align*}

Q.E.D.

\vphantom{}Proving the invertibility of $L_{\mathbf{n}}$ is by a
straightforward, but unusual, argument. 
\begin{lem}
\label{lem:Inversion of L_nu}Let $\mathbf{n}\in\mathbb{N}_{0}^{d}\backslash\left\{ \mathbf{0}\right\} $,
and let $I$ be a unique solution ideal in $R_{d}$. Then, $1-L_{\mathbf{n}}:\mathcal{S}\left(\mathbb{Z}_{p},\textrm{Frac}\left(\mathcal{R}_{d}/I\mathcal{R}_{d}\right)\right)\rightarrow\mathcal{S}\left(\mathbb{Z}_{p},\textrm{Frac}\left(\mathcal{R}_{d}/I\mathcal{R}_{d}\right)\right)$
is an invertible linear transformation if and only if $\left\langle 1-\alpha_{\mathbf{n}}\left(0\right)\right\rangle \cap I=\left\{ 0\right\} $.
When the inverse exists, we have that:
\begin{equation}
\frac{1}{1-\alpha_{\mathbf{n}}\left(0\right)}\frac{r_{\mathbf{n},0}^{n}\kappa_{\mathbf{n}}\left(k\right)}{p^{n}}+\sum_{m=0}^{n-1}\frac{r_{\mathbf{n},0}^{m}\kappa_{\mathbf{n}}\left(\left[k\right]_{p^{m}}\right)}{p^{m}}\left[\mathfrak{z}\overset{p^{n-m}}{\equiv}\theta_{p}^{\circ m}\left(k\right)\right]
\end{equation}
is the image of $\left[\mathfrak{z}\overset{p^{n}}{\equiv}k\right]$
under $\left(1-L\right)^{-1}$.

Moreover, letting $\overline{\overline{\textrm{Frac}\left(\mathcal{R}_{d}/I\mathcal{R}_{d}\right)_{\mathfrak{q}}}}$
be the metric completion of the algebraic closure of the metric completion
of $\textrm{Frac}\left(\mathcal{R}_{d}/I\mathcal{R}_{d}\right)$ with
respect to an absolute value $\left|\cdot\right|_{\mathfrak{q}}$
on $\mathcal{R}_{d}/I\mathcal{R}_{d}$, it then follows that $1-L_{\mathbf{n}}$
extends to be invertible on $W\left(\mathbb{Z}_{p},\overline{\overline{\textrm{Frac}\left(\mathcal{R}_{d}/I\mathcal{R}_{d}\right)_{\mathfrak{q}}}}\right)$
whenever:

I. If $\mathfrak{q}$ is non-archimedean:
\begin{equation}
\max_{0\leq k<p}\left|\frac{r_{\mathbf{m},\mathbf{n},k}}{p}\right|_{\mathfrak{q}}\leq1
\end{equation}

II. If $\mathfrak{q}$ is archimedean:
\begin{equation}
\sum_{k=0}^{p-1}\left|\frac{r_{\mathbf{m},\mathbf{n},k}}{p}\right|_{\mathfrak{q}}<1
\end{equation}
\end{lem}
Proof: Again, we drop the subscript $\mathbf{n}$s except where necessary
for clarity. As everyone knows, one can always write:
\begin{equation}
\left(1-L\right)^{-1}=\frac{1}{1-L}=\sum_{n=0}^{\infty}L^{n}
\end{equation}
provided that there is a complete normed space $\mathcal{X}$ so that
the series on the right converges in $\mathcal{B}\left(\mathcal{X}\right)$,
the Banach algebra of continuous linear operators $\mathcal{X}\rightarrow\mathcal{X}$,
with operator composition as its multiplication map.

So, working in $\mathcal{S}\left(\mathbb{Z}_{p},\textrm{Frac}\left(\mathcal{R}_{d}\right)\right)$\textemdash \emph{note},
we ARE NOT invoking the quotient by $I$ yet!\textemdash since $\mathcal{R}_{d}\left(\zeta_{p^{\infty}}\right)$
is the maximal $p$-power cyclotomic extension of the localization
of $R_{d}$ away from $\left\langle 1-a_{1,0},\ldots,1-a_{d,0}\right\rangle $,
where $R_{d}$ is a ring of polynomials with coefficients in $\mathcal{O}_{K}$,
we have that: 
\begin{equation}
\alpha_{\mathbf{n}}\left(0\right)=\frac{1}{p}\sum_{k=0}^{p-1}r_{\mathbf{n},k}=\frac{1}{p}\sum_{k=0}^{p-1}\prod_{j=1}^{d}a_{j,k}^{n_{j}}
\end{equation}
is \emph{not} a unit of $\mathcal{R}_{d}$. Indeed, if it was a unit,
it would give us a relation among the $a_{j,k}$s, which is impossible,
because $R_{d}$ is a free $\mathcal{O}_{K}$-algebra in the $a_{j,k}$s
and $b_{j,k}$s.

Next, since $R_{d}$ is a Dedekind domain, so is its localization
$\mathcal{R}_{d}$. Thus, the fact that $\alpha_{\mathbf{n}}\left(0\right)\notin\mathcal{R}_{d}^{\times}$
guarantees that the ideal $\left\langle \alpha_{\mathbf{n}}\left(0\right)\right\rangle $
is contained in some non-zero prime ideal $\mathfrak{q}\subseteq\mathcal{R}_{d}$.
As such, considering the $\mathfrak{p}$-adic absolute value $\left|\cdot\right|_{\mathfrak{q}}$
on $\mathcal{R}_{d}$, we have that $\left|\alpha_{\mathbf{n}}\left(0\right)\right|_{\mathfrak{q}}<1$.

Setting $m\geq n$, our formula from \textbf{Proposition} \textbf{\ref{prop:1 minus L_r in terms of shift}}
the previous proposition gives:
\begin{align*}
\left\Vert L_{\mathbf{n}}^{m}\left\{ \left[\cdot\overset{p^{n}}{\equiv}k\right]\right\} \right\Vert _{p,\mathfrak{q}} & =\left|\frac{r_{\mathbf{n},0}^{n}\kappa_{\mathbf{n}}\left(k\right)}{p^{n}}\left(\frac{1}{p}\sum_{j=0}^{p-1}r_{\mathbf{n},j}\right)^{m-n}\right|_{\mathfrak{q}}\\
\left(\frac{1}{p}\sum_{j=0}^{p-1}r_{\mathbf{n},j}=\alpha_{\mathbf{n}}\left(0\right)\right); & \ll_{n}\left|\alpha_{\mathbf{n}}\left(0\right)\right|_{\mathfrak{q}}^{m}
\end{align*}
where the constant of proportionality depends on $n$. Thus:
\begin{equation}
\sum_{m=n}^{\infty}\left\Vert L_{\mathbf{n}}^{m}\left\{ \left[\cdot\overset{p^{n}}{\equiv}k\right]\right\} \right\Vert _{p,\mathfrak{q}}\ll_{n}\sum_{m=n}^{\infty}\left|\alpha_{\mathbf{n}}\left(0\right)\right|_{\mathfrak{q}}^{m}<1
\end{equation}
because $\left|\alpha_{\mathbf{n}}\left(0\right)\right|_{\mathfrak{q}}<1$.
This shows that series:
\begin{equation}
\sum_{m=0}^{\infty}L_{\mathbf{n}}^{m}\left\{ \left[\cdot\overset{p^{n}}{\equiv}k\right]\right\} \left(\mathfrak{z}\right)
\end{equation}
converges in $\textrm{Frac}\left(\mathcal{R}_{d}\right)_{\mathfrak{q}}$,
the completion of $\textrm{Frac}\left(\mathcal{R}_{d}\right)$ with
respect to $\left|\cdot\right|_{\mathfrak{q}}$, albeit pointwise
with respect to $n$. Since the $\left[\mathfrak{z}\overset{p^{n}}{\equiv}k\right]$s
form a basis for $\mathcal{S}\left(\mathbb{Z}_{p},\textrm{Frac}\left(\mathcal{R}_{d}\right)_{\mathfrak{q}}\right)$,
this shows that $1-L_{\mathbf{n}}$ is invertible on $\mathcal{S}\left(\mathbb{Z}_{p},\textrm{Frac}\left(\mathcal{R}_{d}\right)_{\mathfrak{q}}\right)$.
Moreover, using \textbf{Proposition \ref{prop:1 minus L_r in terms of shift}}
again, we have that:
\begin{align*}
\sum_{m=0}^{\infty}L_{\mathbf{n}}^{m}\left\{ \left[\cdot\overset{p^{n}}{\equiv}k\right]\right\} \left(\mathfrak{z}\right) & =\sum_{m=0}^{n-1}L_{\mathbf{n}}^{m}\left\{ \left[\cdot\overset{p^{n}}{\equiv}k\right]\right\} \left(\mathfrak{z}\right)+\sum_{m=n}^{\infty}\frac{r_{\mathbf{n},0}^{n}\kappa_{\mathbf{n}}\left(k\right)}{p^{n}}\left(\alpha_{\mathbf{n}}\left(0\right)\right)^{m-n}\\
 & =\frac{1}{1-\alpha_{\mathbf{n}}\left(0\right)}\frac{r_{\mathbf{n},0}^{n}\kappa_{\mathbf{n}}\left(k\right)}{p^{n}}+\sum_{m=0}^{n-1}\frac{r_{\mathbf{n},0}^{m}\kappa_{\mathbf{n}}\left(\left[k\right]_{p^{m}}\right)}{p^{m}}\left[\mathfrak{z}\overset{p^{n-m}}{\equiv}\theta_{p}^{\circ m}\left(k\right)\right]
\end{align*}

So, let:
\begin{equation}
S_{n,k}\left(\mathfrak{z}\right)\overset{\textrm{def}}{=}\frac{r_{\mathbf{n},0}^{n}\kappa_{\mathbf{n}}\left(k\right)/p^{n}}{1-\alpha_{\mathbf{n}}\left(0\right)}+\sum_{m=0}^{n-1}\frac{r_{\mathbf{n},0}^{m}\kappa_{\mathbf{n}}\left(\left[k\right]_{p^{m}}\right)}{p^{m}}\left[\mathfrak{z}\overset{p^{n-m}}{\equiv}\theta_{p}^{\circ m}\left(k\right)\right]
\end{equation}
Observe that $S_{n,k}$ takes values in $\textrm{Frac}\left(\mathcal{R}_{d}\right)$.
Hence, the linear map $T:\mathcal{S}\left(\mathbb{Z}_{p},\textrm{Frac}\left(\mathcal{R}_{d}\right)\right)\rightarrow\mathcal{S}\left(\mathbb{Z}_{p},\textrm{Frac}\left(\mathcal{R}_{d}\right)\right)$
defined by:
\begin{equation}
T\left\{ \left[\cdot\overset{p^{n}}{\equiv}k\right]\right\} =S_{n,k}
\end{equation}
is actually the inverse of $1-L_{\mathbf{n}}$ over $\textrm{Frac}\left(\mathcal{R}_{d}\right)$.

Indeed, we don't even need the absolute value $\left|\cdot\right|_{\mathfrak{q}}$
to demonstrate this! Applying $1-L_{\mathbf{n}}$ to $S_{n,k}$ directly
and using \textbf{Proposition \ref{prop:1 minus L_r in terms of shift}},
we get:

\begin{align*}
\left(1-L_{\mathbf{n}}\right)\left\{ S_{n,k}\right\} \left(\mathfrak{z}\right) & =\frac{r_{\mathbf{n},0}^{n}\kappa_{\mathbf{n}}\left(k\right)/p^{n}}{1-\alpha_{\mathbf{n}}\left(0\right)}+\sum_{m=0}^{n-1}\frac{r_{\mathbf{n},0}^{m}\kappa_{\mathbf{n}}\left(\left[k\right]_{p^{m}}\right)}{p^{m}}\left[\mathfrak{z}\overset{p^{n-m}}{\equiv}\theta_{p}^{\circ m}\left(k\right)\right]\\
 & -\alpha_{\mathbf{n}}\left(0\right)\frac{r_{\mathbf{n},0}^{n}\kappa_{\mathbf{n}}\left(k\right)/p^{n}}{1-\alpha_{\mathbf{n}}\left(0\right)}\\
 & -\sum_{m=0}^{n-1}\frac{r_{\mathbf{n},0}^{m}\kappa_{\mathbf{n}}\left(\left[k\right]_{p^{m}}\right)}{p^{m}}\frac{r_{\mathbf{n},0}\kappa_{\mathbf{n}}\left(\left[\theta_{p}^{\circ m}\left(k\right)\right]_{p}\right)}{p}\left[\mathfrak{z}\overset{p^{n-m-1}}{\equiv}\theta_{p}^{\circ m+1}\left(k\right)\right]\\
\left(\mathbf{Proposition}\textrm{ }\ref{prop:Kappa shift equation}\right); & =\frac{r_{\mathbf{n},0}^{n}\kappa_{\mathbf{n}}\left(k\right)}{p^{n}}+\overbrace{\sum_{m=0}^{n-1}\frac{r_{\mathbf{n},0}^{m}\kappa_{\mathbf{n}}\left(\left[k\right]_{p^{m}}\right)}{p^{m}}\left[\mathfrak{z}\overset{p^{n-m}}{\equiv}\theta_{p}^{\circ m}\left(k\right)\right]}^{\textrm{cancel}}\\
 & -\underbrace{\sum_{m=1}^{n}\frac{r_{\mathbf{n},0}^{m}\kappa_{\mathbf{n}}\left(\left[k\right]_{p^{m}}\right)}{p^{m}}\left[\mathfrak{z}\overset{p^{n-m}}{\equiv}\theta_{p}^{\circ m}\left(k\right)\right]}_{\textrm{cancel}}\\
 & =\frac{r_{\mathbf{n},0}^{n}\kappa_{\mathbf{n}}\left(k\right)}{p^{n}}-\frac{r_{\mathbf{n},0}^{n}\kappa_{\mathbf{n}}\left(k\right)}{p^{n}}+\left[\mathfrak{z}\overset{p^{n}}{\equiv}k\right]\\
 & =\left[\mathfrak{z}\overset{p^{n}}{\equiv}k\right]
\end{align*}
and we needed no absolute values or topologies to do any of that!
This shows that $S_{n,k}$ truly is an element of $\mathcal{S}\left(\mathbb{Z}_{p},\textrm{Frac}\left(\mathcal{R}_{d}\right)\right)$
that $1-L_{\mathbf{n}}$ sends to $\left[\mathfrak{z}\overset{p^{n}}{\equiv}k\right]$.
Moreover $S_{n,k}$ is the \emph{unique }function in $\mathcal{S}\left(\mathbb{Z}_{p},\textrm{Frac}\left(\mathcal{R}_{d}\right)\right)$
with this property, for if there was another such function, then it
would then be an element of $\mathcal{S}\left(\mathbb{Z}_{p},\textrm{Frac}\left(\mathcal{R}_{d}\right)_{\mathfrak{q}}\right)$
which $1-L_{\mathbf{n}}$ sends to $\left[\mathfrak{z}\overset{p^{n}}{\equiv}k\right]$,
and thus must be equal to $S_{n,k}$, because of the invertibility
of $1-L_{\mathbf{n}}$ over $\mathcal{S}\left(\mathbb{Z}_{p},\textrm{Frac}\left(\mathcal{R}_{d}\right)_{\mathfrak{q}}\right)$.

Now, if we consider the quotient by an ideal $I\subseteq R_{d}$,
note that $S_{n,k}$ is a perfectly well-defined element of $\mathcal{S}\left(\mathbb{Z}_{p},\textrm{Frac}\left(\mathcal{R}_{d}/I\mathcal{R}_{d}\right)\right)$
if and only if $\left\langle 1-\alpha_{\mathbf{n}}\left(0\right)\right\rangle \nsubseteq I$,
from which the invertibility of $1-L_{\mathbf{n}}$ follows (or doesn't,
if $\left\langle 1-\alpha_{\mathbf{n}}\left(0\right)\right\rangle \subseteq I$).

Finally, for the matter of extending to the Wiener algebra, it is
known that for any algebraically closed Banach ring $B$, the set:
\begin{equation}
\mathcal{B}\overset{\textrm{def}}{=}\left\{ \left[\mathfrak{z}\overset{p^{\lambda_{p}\left(n\right)}}{\equiv}n\right]:n\geq0\right\} 
\end{equation}
is a basis for the Banach space $W\left(\mathbb{Z}_{p},B\right)$\textemdash the
\textbf{van der Put basis} \cite{2nd blog paper,My Dissertation,Ultrametric Calculus,Robert's Book}.
Every $f\in W\left(\mathbb{Z}_{p},B\right)$ admits a unique \textbf{van
der Put series }representation with respect to this basis:
\begin{equation}
f\left(\mathfrak{z}\right)=\sum_{n=0}^{\infty}c_{n}\left(f\right)\left[\mathfrak{z}\overset{p^{\lambda_{p}\left(n\right)}}{\equiv}n\right]
\end{equation}
where the coefficients $c_{n}\left(f\right)$ are called the \textbf{van
der Put }coefficients of $f$. Such series necessarily converge in
$B$ uniformly with respect to $\mathfrak{z}\in\mathbb{Z}_{p}$. We
will prove that the place $v$ we selected guarantees that $\left\Vert S_{\lambda_{p}\left(n\right),n}\right\Vert _{p,v}\leq1$
holds for all $n\geq0$. Because of this, we can then define $\left(1-L_{\mathbf{n}}\right)^{-1}$
for every element of $C\left(\mathbb{Z}_{p},B\right)$ by:
\begin{equation}
\left(1-L_{\mathbf{n}}\right)^{-1}\left\{ f\right\} \left(\mathfrak{z}\right)\overset{\textrm{def}}{=}\sum_{n=0}^{\infty}c_{n}\left(f\right)S_{\lambda_{p}\left(n\right),n}\left(\mathfrak{z}\right)
\end{equation}

\begin{claim}
Let $\mathfrak{q}$ be as given in the statement of the lemma. Then,
$\left\Vert S_{\lambda_{p}\left(n\right),n}\right\Vert _{p,\mathfrak{q}}\leq1$
occurs for all $n\geq0$. The equality is strict if $\mathfrak{q}$
is archimedean

Proof of Claim: Letting $n\geq0$ and $k\in\left\{ 0,\ldots,p^{n}-1\right\} $,
taking supremum norms:
\begin{align*}
\sum_{m=0}^{n-1}\left\Vert \left(1-L_{\mathbf{n}}\right)^{m}\left\{ \left[\cdot\overset{p^{n}}{\equiv}k\right]\right\} \right\Vert _{p,\mathfrak{q}} & \leq\begin{cases}
\max_{0\leq m<n}\left(\max_{0\leq j<p}\left|\frac{r_{\mathbf{n},j}}{p}\right|_{\mathfrak{q}}\right)^{m} & \textrm{if }v<\infty\\
\sum_{m=0}^{n-1}\left(\sum_{j=0}^{p-1}\left|\frac{r_{\mathbf{n},j}}{p}\right|_{\mathfrak{q}}\right)^{m} & \textrm{else}
\end{cases}\\
 & \leq\begin{cases}
1 & \textrm{if }\mathfrak{q}<\infty\\
\sum_{m=0}^{n-1}d^{m} & \textrm{else}
\end{cases}
\end{align*}
where 
\begin{equation}
d\overset{\textrm{def}}{=}\sum_{j=0}^{p-1}\left|\frac{r_{\mathbf{n},j}}{p}\right|_{\mathfrak{q}}
\end{equation}
is strictly less then $1$ when $\mathfrak{q}$ is archimedean. Since
$0<d<1$, we have:
\begin{equation}
\sum_{m=0}^{n-1}d^{m}<\sum_{m=0}^{\infty}d^{m}=\frac{1}{1-d}
\end{equation}
a bound independent of $n$, and so:
\begin{equation}
\sup_{n\geq0}\sup_{0\leq k<p^{n}}\sum_{m=0}^{n-1}\left\Vert \left(1-L_{\mathbf{n}}\right)^{m}\left\{ \left[\cdot\overset{p^{n}}{\equiv}k\right]\right\} \right\Vert _{p,\mathfrak{q}}<\infty
\end{equation}
Thus:
\begin{align*}
\left\Vert S_{n,k}\right\Vert _{p,\mathfrak{q}} & \leq\left|\frac{r_{\mathbf{n},0}^{n}\kappa_{\mathbf{n}}\left(k\right)}{p^{n}}\frac{1}{1-\alpha_{\mathbf{n}}\left(0\right)}\right|_{\mathfrak{q}}+\underbrace{\sum_{m=0}^{n-1}\left\Vert \left(1-L_{\mathbf{n}}\right)^{m}\left\{ \left[\cdot\overset{p^{n}}{\equiv}k\right]\right\} \right\Vert _{p,\mathfrak{q}}}_{O\left(1\right)\textrm{ with respect to }n,k}\\
\left(\left\langle 1-\alpha_{\mathbf{n}}\left(0\right)\right\rangle \nsubseteq I\right); & \ll\left|\frac{r_{\mathbf{n},0}^{n}\kappa_{\mathbf{n}}\left(k\right)}{p^{n}}\right|_{\mathfrak{q}}
\end{align*}
Here, if $\mathfrak{q}$ is non-archimedean, our assumption that:
\begin{equation}
\max_{0\leq j<p}\left|\frac{r_{\mathbf{n},j}}{p}\right|_{\mathfrak{q}}\leq1
\end{equation}
then shows that:
\begin{equation}
\left\Vert S_{\lambda_{p}\left(n\right),n}\right\Vert _{p,\mathfrak{q}}\leq1,\textrm{ }\forall n\geq0
\end{equation}
If, on the other hand, $\mathfrak{q}$ is archimedean, our assumption
that:
\begin{equation}
\sum_{j=0}^{p-1}\left|\frac{r_{\mathbf{n},j}}{p}\right|_{\mathfrak{q}}<1
\end{equation}
forces:
\begin{equation}
\max_{0\leq j<p}\left|\frac{r_{\mathbf{n},j}}{p}\right|_{\mathfrak{q}}<1
\end{equation}
and hence: 
\begin{equation}
\left\Vert S_{\lambda_{p}\left(n\right),n}\right\Vert _{p,\mathfrak{q}}\ll\prod_{j=0}^{p-1}\left|\frac{r_{\mathbf{n},j}}{p}\right|^{\#_{p:j}\left(n\right)}<1
\end{equation}
This proves the Claim. $\checkmark$
\end{claim}
Finally we conclude, working in our two cases separately. Let $B$
denote $\overline{\overline{\textrm{Frac}\left(\mathcal{R}_{d}/I\mathcal{R}_{d}\right)_{\mathfrak{q}}}}$,
where $\mathfrak{q}$ is as given in the statement of the lemma.

I. First, if $\mathfrak{q}$ is a non-archimedean absolute value on
$B$, then for any $f\in W\left(\mathbb{Z}_{p},B\right)$:
\begin{equation}
\left\Vert c_{n}\left(f\right)S_{\lambda_{p}\left(n\right),n}\right\Vert _{p,\mathfrak{q}}=\left|c_{n}\left(f\right)\right|_{\mathfrak{q}}\left\Vert S_{\lambda_{p}\left(n\right),n}\right\Vert _{p,\mathfrak{q}}\leq\left|c_{n}\left(f\right)\right|_{\mathfrak{q}}
\end{equation}
Since $f$ is continuous, $\left|c_{n}\left(f\right)\right|_{\mathfrak{q}}\rightarrow0$
in $\mathbb{R}$ as $n\rightarrow\infty$. Consequently:
\begin{equation}
\left(1-L_{\mathbf{n}}\right)^{-1}\left\{ f\right\} \left(\mathfrak{z}\right)\overset{\textrm{def}}{=}\sum_{n=0}^{\infty}c_{n}\left(f\right)S_{\lambda_{p}\left(n\right),n}\left(\mathfrak{z}\right)
\end{equation}
is a sum of functions in the non-archimedean Banach ring $W\left(\mathbb{Z}_{p},B\right)$
so that the norm of the $n$th term tends to $0$ in $\mathbb{R}$
as $n\rightarrow\infty$. Because this Banach ring is non-archimedean,
this decay condition then guarantees the convergence of $\left(1-L_{\mathbf{n}}\right)^{-1}\left\{ f\right\} \left(\mathfrak{z}\right)$
in $W\left(\mathbb{Z}_{p},B\right)$. Moreover, by the ultrametric
inequality of this Banach space's norm, we have that:
\begin{equation}
\left\Vert \left(1-L_{\mathbf{n}}\right)^{-1}\left\{ f\right\} \right\Vert _{p,\mathfrak{q}}\leq\sup_{n\geq0}\left(\left|c_{n}\left(f\right)\right|_{\mathfrak{q}}\left\Vert S_{\lambda_{p}\left(n\right),n}\right\Vert _{p,\mathfrak{q}}\right)\leq\sup_{n\geq0}\left|c_{n}\left(f\right)\right|_{\mathfrak{q}}
\end{equation}
Since the $c_{n}$s decay to $0$ $\mathfrak{q}$-adically, there
are distinct $m,n$ so that $\left|c_{m}\left(f\right)\right|_{\mathfrak{q}}\neq\left|c_{n}\left(f\right)\right|_{\mathfrak{q}}$.
Consequently, the ultrametric inequality tells us that:
\begin{equation}
\left\Vert f\right\Vert _{p,\mathfrak{q}}=\left\Vert \sum_{n=0}^{\infty}c_{n}\left(f\right)\left[\mathfrak{z}\overset{p^{\lambda_{p}\left(n\right)}}{\equiv}n\right]\right\Vert _{p,\mathfrak{q}}=\sup_{n\geq0}\left|c_{n}\left(f\right)\right|_{\mathfrak{q}}
\end{equation}
and hence:
\begin{equation}
\left\Vert \left(1-L_{\mathbf{n}}\right)^{-1}\left\{ f\right\} \right\Vert _{p,\mathfrak{q}}\leq\sup_{n\geq0}\left|c_{n}\left(f\right)\right|_{\mathfrak{q}}=\left\Vert f\right\Vert _{p,\mathfrak{q}}
\end{equation}
which proves that $\left(1-L_{\mathbf{n}}\right)^{-1}$ is a continuous
linear operator on $W\left(\mathbb{Z}_{p},B\right)$, and hence, that
$1-L_{\mathbf{n}}$ is invertible on $W\left(\mathbb{Z}_{p},B\right)$.

\vphantom{}II. Suppose $\mathfrak{q}$ is archimedean. Then, $f\in W\left(\mathbb{Z}_{p},B\right)$
implies that $f$'s integral converges with respect to $\mathfrak{q}$;
this is:
\begin{equation}
\int_{\mathbb{Z}_{p}}f\left(\mathfrak{z}\right)d\mathfrak{z}=\sum_{n=0}^{\infty}\frac{c_{n}\left(f\right)}{p^{\lambda_{p}\left(n\right)}}
\end{equation}
Moreover, as shown in \cite{2nd blog paper}, we have:
\begin{equation}
\hat{f}\left(t\right)=\sum_{n=\frac{1}{p}\left|t\right|_{p}}^{\infty}\frac{c_{n}\left(f\right)}{p^{\lambda_{p}\left(n\right)}}e^{-2\pi int},\textrm{ }\forall t\in\hat{\mathbb{Z}}_{p}
\end{equation}
Noting that:
\begin{equation}
e^{2\pi i\left\{ t\mathfrak{z}\right\} _{p}}=\sum_{m=0}^{\max\left\{ 0,\left|t\right|_{p}-1\right\} }e^{2\pi imt}\left[\mathfrak{z}\overset{\left|t\right|_{p}}{\equiv}m\right]
\end{equation}
we get:
\begin{align*}
\left(1-L_{\mathbf{n}}\right)^{-1}\left\{ e^{2\pi i\left\{ t\cdot\right\} _{p}}\right\} \left(\mathfrak{z}\right) & =\sum_{m=0}^{\left|t\right|_{p}-1}e^{2\pi imt}\frac{r_{\mathbf{n},0}\kappa_{\mathbf{n}}\left(\left[m\right]_{p}\right)}{p}\left[\mathfrak{z}\overset{\left|t\right|_{p}/p}{\equiv}\theta_{p}\left(m\right)\right]\\
 & =\frac{r_{\mathbf{n},0}}{p}\sum_{m=0}^{\left|t\right|_{p}-1}e^{2\pi imt}\kappa_{\mathbf{n}}\left(\left[m\right]_{p}\right)\left[\mathfrak{z}\overset{\left|t\right|_{p}/p}{\equiv}\theta_{p}\left(m\right)\right]
\end{align*}
Now:
\begin{align*}
\left(1-L_{\mathbf{n}}\right)^{-1}\left\{ f-\hat{f}\left(0\right)\right\} \left(\mathfrak{z}\right) & =\left(1-L_{\mathbf{n}}\right)^{-1}\left(\sum_{t\in\hat{\mathbb{Z}}_{p}\backslash\left\{ 0\right\} }\hat{f}\left(t\right)e^{2\pi i\left\{ t\mathfrak{z}\right\} _{p}}\right)\\
 & =\frac{r_{\mathbf{n},0}}{p}\sum_{t\in\hat{\mathbb{Z}}_{p}\backslash\left\{ 0\right\} }\hat{f}\left(t\right)\sum_{m=0}^{\left|t\right|_{p}-1}e^{2\pi imt}\kappa_{\mathbf{n}}\left(\left[m\right]_{p}\right)\left[\mathfrak{z}\overset{\left|t\right|_{p}/p}{\equiv}\theta_{p}\left(m\right)\right]\\
 & =\frac{r_{\mathbf{n},0}}{p}\sum_{n=1}^{\infty}\sum_{m=0}^{p^{n}-1}\kappa_{\mathbf{n}}\left(\left[m\right]_{p}\right)\left[\mathfrak{z}\overset{p^{n-1}}{\equiv}\theta_{p}\left(m\right)\right]\sum_{\left|t\right|_{p}=p^{n}}\hat{f}\left(t\right)e^{2\pi imt}
\end{align*}
Using the summation identity:
\begin{equation}
\sum_{m=0}^{p^{n}-1}g\left(m\right)=\sum_{j=0}^{p-1}\sum_{m=0}^{p^{n-1}-1}g\left(pm+j\right)
\end{equation}
we get:
\begin{align*}
\left(1-L_{\mathbf{n}}\right)^{-1}\left\{ f-\hat{f}\left(0\right)\right\} \left(\mathfrak{z}\right) & =\frac{r_{\mathbf{n},0}}{p}\sum_{n=1}^{\infty}\sum_{j=0}^{p-1}\sum_{m=0}^{p^{n-1}-1}\kappa_{\mathbf{n}}\left(j\right)\left[\mathfrak{z}\overset{p^{n-1}}{\equiv}m\right]\sum_{\left|t\right|_{p}=p^{n}}\hat{f}\left(t\right)e^{2\pi i\left(pm+j\right)t}\\
 & =\frac{r_{\mathbf{n},0}}{p}\sum_{j=0}^{p-1}\kappa_{\mathbf{n}}\left(j\right)\sum_{n=1}^{\infty}\sum_{\left|t\right|_{p}=p^{n}}\hat{f}\left(t\right)e^{2\pi i\left(p\left[\mathfrak{z}\right]_{p^{n-1}}+j\right)t}
\end{align*}
Applying the archimedean absolute value $\left|\cdot\right|_{\mathfrak{q}}$,
we get: 
\begin{equation}
\left|\left(1-L_{\mathbf{n}}\right)^{-1}\left\{ f-\hat{f}\left(0\right)\right\} \left(\mathfrak{z}\right)\right|_{\mathfrak{q}}\ll_{\mathbf{n}}\sum_{n=1}^{\infty}\sum_{\left|t\right|_{p}=p^{n}}\left|\hat{f}\left(t\right)\right|_{\mathfrak{q}}\leq\sum_{t\in\hat{\mathbb{Z}}_{p}}\left|\hat{f}\left(t\right)\right|_{\mathfrak{q}}=\left\Vert f\right\Vert _{W\left(\mathbb{Z}_{p},B\right)}<\infty
\end{equation}
where the equality at the right occurs because of the definition of
$W\left(\mathbb{Z}_{p},B\right)$ when $\mathfrak{q}$ is archimedean.
This shows that $\left(1-L_{\mathbf{n}}\right)^{-1}\left\{ f\right\} \left(\mathfrak{z}\right)-\left(1-L_{\mathbf{n}}\right)^{-1}\left\{ \hat{f}\left(0\right)\right\} \left(\mathfrak{z}\right)$
converges with respect to $\left|\cdot\right|_{\mathfrak{q}}$ uniformly
in $\mathfrak{z}\in\mathbb{Z}_{p}$. Since:
\begin{equation}
\left(1-L_{\mathbf{n}}\right)^{-1}\left\{ \hat{f}\left(0\right)\right\} \left(\mathfrak{z}\right)=\hat{f}\left(0\right)S_{0,0}\left(\mathfrak{z}\right)=\frac{\hat{f}\left(0\right)}{1-\alpha_{\mathbf{n}}\left(0\right)}
\end{equation}
which is a constant, we see that $\left(1-L_{\mathbf{n}}\right)^{-1}\left\{ f\right\} :\mathbb{Z}_{p}\rightarrow B$
is continuous.

Now, let $g$ denote $\left(1-L_{\mathbf{n}}\right)^{-1}\left\{ f\right\} $.
Then:
\begin{equation}
g\left(\mathfrak{z}\right)=\frac{\hat{f}\left(0\right)}{1-\alpha_{\mathbf{n}}\left(0\right)}+\frac{r_{\mathbf{n},0}}{p}\sum_{j=0}^{p-1}\kappa_{\mathbf{n}}\left(j\right)\sum_{n=1}^{\infty}\sum_{\left|s\right|_{p}=p^{n}}\hat{f}\left(s\right)e^{2\pi i\left(p\left[\mathfrak{z}\right]_{p^{n-1}}+j\right)s}
\end{equation}
where, as we just showed, the series converges uniformly in $\mathfrak{z}$.
As such, we can integrate it term-by-term. The main integral is:
\begin{equation}
\int_{\mathbb{Z}_{p}}e^{2\pi ip\left[\mathfrak{z}\right]_{p^{n-1}}s}e^{-2\pi i\left\{ t\mathfrak{z}\right\} _{p}}d\mathfrak{z}
\end{equation}
which we deal with by noting that since $s$ has $\left|s\right|_{p}=p^{n}$,
$ps$ has $\left|ps\right|_{p}=p^{n-1}$, and so:
\begin{equation}
p\left[\mathfrak{z}\right]_{p^{n-1}}s\overset{1}{\equiv}\left\{ ps\right\} _{p},\textrm{ }\forall s:\left|s\right|_{p}=p^{n}
\end{equation}
Hence:
\begin{equation}
\int_{\mathbb{Z}_{p}}e^{2\pi ip\left[\mathfrak{z}\right]_{p^{n-1}}s}e^{-2\pi i\left\{ t\mathfrak{z}\right\} _{p}}d\mathfrak{z}=\int_{\mathbb{Z}_{p}}e^{2\pi i\left\{ \left(ps-t\right)\mathfrak{z}\right\} _{p}}d\mathfrak{z}=\mathbf{1}_{0}\left(ps-t\right)=\left[t\overset{1}{\equiv}ps\right]
\end{equation}
which will be $1$ whenever $ps\overset{1}{\equiv}t$ and will be
$0$ otherwise. We then end up with:
\begin{equation}
\hat{g}\left(t\right)=\frac{\hat{f}\left(0\right)\mathbf{1}_{0}\left(t\right)}{1-\alpha_{\mathbf{n}}\left(0\right)}+\frac{r_{\mathbf{n},0}}{p}\sum_{j=0}^{p-1}\kappa_{\mathbf{n}}\left(j\right)\sum_{n=1}^{\infty}\sum_{\left|s\right|_{p}=p^{n}}\hat{f}\left(s\right)e^{2\pi ijs}\left[t\overset{1}{\equiv}ps\right]
\end{equation}
Applying \textbf{Proposition \ref{prop:adjoint}}, we have:
\begin{align*}
\sum_{\left|s\right|_{p}\leq p^{n}}\hat{f}\left(s\right)e^{2\pi ijs}\left[t\overset{1}{\equiv}ps\right] & =\sum_{\left|s\right|_{p}\leq p^{n-1}}\sum_{\left|r\right|_{p}\leq p}\hat{f}\left(\frac{s}{p}+r\right)e^{2\pi ij\left(\frac{s}{p}+r\right)}\left[t\overset{1}{\equiv}s\right]\\
 & =\mathbf{1}_{0}\left(p^{n-1}t\right)\sum_{\left|r\right|_{p}\leq p}\hat{f}\left(\frac{t}{p}+r\right)e^{2\pi ij\left(\frac{t}{p}+r\right)}
\end{align*}
Hence:
\begin{align*}
\sum_{\left|s\right|_{p}=p^{n}}\hat{f}\left(s\right)e^{2\pi ijs}\left[t\overset{1}{\equiv}ps\right] & =\sum_{\left|s\right|_{p}\leq p^{n}}\hat{f}\left(s\right)e^{2\pi ijs}\left[t\overset{1}{\equiv}ps\right]-\sum_{\left|s\right|_{p}\leq p^{n-1}}\hat{f}\left(s\right)e^{2\pi ijs}\left[t\overset{1}{\equiv}ps\right]\\
 & =\left(\mathbf{1}_{0}\left(p^{n-1}t\right)-\mathbf{1}_{0}\left(p^{n-2}t\right)\right)\sum_{\left|r\right|_{p}\leq p}\hat{f}\left(\frac{t}{p}+r\right)e^{2\pi ij\left(\frac{t}{p}+r\right)}
\end{align*}
for all $n\geq2$, and we get:
\begin{align*}
\hat{g}\left(t\right) & =\frac{\hat{f}\left(0\right)\mathbf{1}_{0}\left(t\right)}{1-\alpha_{\mathbf{n}}\left(0\right)}+\frac{r_{\mathbf{n},0}}{p}\sum_{j=0}^{p-1}\kappa_{\mathbf{n}}\left(j\right)\sum_{\left|s\right|_{p}=p}\hat{f}\left(s\right)e^{2\pi ijs}\underbrace{\left[t\overset{1}{\equiv}ps\right]}_{\mathbf{1}_{0}\left(t\right)}\\
 & +\frac{r_{\mathbf{n},0}}{p}\sum_{j=0}^{p-1}\kappa_{\mathbf{n}}\left(j\right)\sum_{n=2}^{\infty}\sum_{\left|s\right|_{p}=p^{n}}\hat{f}\left(s\right)e^{2\pi ijs}\left[t\overset{1}{\equiv}ps\right]\\
 & =\mathbf{1}_{0}\left(t\right)\underbrace{\left(\frac{\hat{f}\left(0\right)}{1-\alpha_{\mathbf{n}}\left(0\right)}+\frac{r_{\mathbf{n},0}}{p}\sum_{j=0}^{p-1}\kappa_{\mathbf{n}}\left(j\right)\sum_{\left|s\right|_{p}=p}\hat{f}\left(s\right)e^{2\pi ijs}\right)}_{=C,\textrm{ a constant independent of }t}\\
 & +\frac{r_{\mathbf{n},0}}{p}\sum_{j=0}^{p-1}\kappa_{\mathbf{n}}\left(j\right)\underbrace{\sum_{n=0}^{\infty}\left(\mathbf{1}_{0}\left(p^{n+1}t\right)-\mathbf{1}_{0}\left(p^{n}t\right)\right)}_{\textrm{telescoping}}\sum_{\left|r\right|_{p}\leq p}\hat{f}\left(\frac{t}{p}+r\right)e^{2\pi ij\left(\frac{t}{p}+r\right)}\\
 & =C\mathbf{1}_{0}\left(t\right)+\frac{r_{\mathbf{n},0}}{p}\sum_{j=0}^{p-1}\kappa_{\mathbf{n}}\left(j\right)\left(\underbrace{\lim_{n\rightarrow\infty}\mathbf{1}_{0}\left(p^{n}t\right)}_{1,\textrm{ }\forall t}-\mathbf{1}_{0}\left(t\right)\right)\sum_{\left|r\right|_{p}\leq p}\hat{f}\left(\frac{t}{p}+r\right)e^{2\pi ij\left(\frac{t}{p}+r\right)}\\
 & =C\mathbf{1}_{0}\left(t\right)+\left(1-\mathbf{1}_{0}\left(t\right)\right)\frac{r_{\mathbf{n},0}}{p}\sum_{j=0}^{p-1}\kappa_{\mathbf{n}}\left(j\right)\sum_{\left|r\right|_{p}\leq p}\hat{f}\left(\frac{t}{p}+r\right)e^{2\pi ij\left(\frac{t}{p}+r\right)}
\end{align*}
and so:
\begin{align*}
\sum_{t\in\hat{\mathbb{Z}}_{p}}\left|\hat{g}\left(t\right)\right|_{\mathfrak{q}} & \ll\sum_{t\in\hat{\mathbb{Z}}_{p}}\left|\sum_{\left|r\right|_{p}\leq p}\hat{f}\left(\frac{t}{p}+r\right)e^{2\pi ij\left(\frac{t}{p}+r\right)}\right|_{\mathfrak{q}}\\
 & \leq\sum_{\left|r\right|_{p}\leq p}\sum_{t\in\hat{\mathbb{Z}}_{p}}\left|\hat{f}\left(\frac{t}{p}+r\right)\right|_{\mathfrak{q}}\\
 & \ll\sum_{\left|r\right|_{p}\leq p}\left\Vert f\right\Vert _{W\left(\mathbb{Z}_{p},B\right)}\\
 & =p\left\Vert f\right\Vert _{W\left(\mathbb{Z}_{p},B\right)}
\end{align*}
and so $\left(1-L_{\mathbf{n}}\right)^{-1}\left\{ f\right\} =g\in W\left(\mathbb{Z}_{p},B\right)$.
Moreover, the above estimate then shows that $\left(1-L_{\mathbf{n}}\right)^{-1}$
is continuous with respect to the norm of $W\left(\mathbb{Z}_{p},B\right)$,
as desired.

Q.E.D.

\vphantom{}

We now have:
\begin{thm}
\label{thm:inversion method}Let $\mathbf{n}\in\mathbb{N}_{0}^{d}\backslash\left\{ \mathbf{0}\right\} $,
and let $I$ be a unique solution ideal in $R_{d}$. Then, the $\mathbf{n}$th
formal solution $dY_{\mathbf{n}}\in\mathcal{S}\left(\mathbb{Z}_{p},\mathcal{R}_{d}/I\mathcal{R}_{d}\right)^{\prime}$
relative to the initial condition $\left\{ \hat{\chi}_{\mathbf{m}}\right\} _{\mathbf{m}\in\left[0,\mathbf{n}\right)}$
is uniquely determined by the initial condition whenever $\left\langle 1-\alpha_{\mathbf{n}}\left(0\right)\right\rangle \nsubseteq I$.

Moreover, letting $\left|\cdot\right|_{\mathfrak{q}}$ be an absolute
value on $\mathcal{R}_{d}/I\mathcal{R}_{d}$, if $\hat{\chi}_{\mathbf{m}}\in W\left(\mathbb{Z}_{p},\overline{\overline{\textrm{Frac}\left(\mathcal{R}_{d}/I\mathcal{R}_{d}\right)_{\mathfrak{q}}}}\right)^{\prime}$
for all $\mathbf{m}\in\left[0,\mathbf{n}\right)$, then $dY_{\mathbf{n}}\in W\left(\mathbb{Z}_{p},\overline{\overline{\textrm{Frac}\left(\mathcal{R}_{d}/I\mathcal{R}_{d}\right)_{\mathfrak{q}}}}\right)^{\prime}$
whenever:

I. If $\mathfrak{q}$ is non-archimedean:
\begin{equation}
\max_{0\leq k<p}\left|\frac{r_{\mathbf{m},\mathbf{n},k}}{p}\right|_{\mathfrak{q}}\leq1
\end{equation}

II. If $\mathfrak{q}$ is archimedean:
\begin{equation}
\sum_{k=0}^{p-1}\left|\frac{r_{\mathbf{m},\mathbf{n},k}}{p}\right|_{\mathfrak{q}}<1
\end{equation}
\end{thm}
Proof: The first part of the statement follows from \textbf{Lemma
\ref{lem:Inversion of L_nu}} and \textbf{Proposition \ref{prop:The-th-formal}}.
As for being a measure on the wiener algebra, this follows from the
hypothesis on the initial condition and the invertibility of $1-L_{\mathbf{n}}$
over the Wiener algebra when (I) and (II) are satisfied, as guaranteed
by \textbf{Lemma \ref{lem:Inversion of L_nu}}.
Q.E.D.

\vphantom{}Finally, we will show that the Fourier transform of $\hat{Y}_{\mathbf{n}}$
we have obtained is in fact a Fourier transform of $X_{\mathbf{n}}$.
We do this with the help of the $p$-adic Dirichlet Kernel.
\begin{defn}
Let $R$ be a commutative unital ring whose characteristic, if positive,
is co-prime to $p$. Then, the \textbf{$p$-adic Dirichlet Kernel},
$\left\{ D_{N}\right\} _{N\geq0}$ is the sequence of functions in
$\mathcal{S}\left(\mathbb{Z}_{p},R\right)$ defined by:
\begin{equation}
D_{N}\left(\mathfrak{z}\right)\overset{\textrm{def}}{=}p^{N}\left[\mathfrak{z}\overset{p^{N}}{\equiv}0\right]
\end{equation}
\end{defn}
Just like in classical Fourier analysis, convolving a function with
$D_{N}$ outputs the $N$th partial sum of the function's Fourier
series. Given any $\phi\in\mathcal{S}\left(\mathbb{Z}_{p},R\right)$,
we have:
\begin{equation}
\left(\phi*D_{N}\right)\left(\mathfrak{z}\right)=\int_{\mathbb{Z}_{p}}\phi\left(\mathfrak{y}\right)D_{N}\left(\mathfrak{z}-\mathfrak{y}\right)d\mathfrak{y}=\sum_{\left|t\right|_{p}\leq p^{N}}\hat{\phi}\left(t\right)e^{2\pi i\left\{ t\mathfrak{z}\right\} _{p}}
\end{equation}
In particular, when convolving $D_{N}$ with the distribution $dY_{\mathbf{n}}$,
we have:
\begin{equation}
\left(dY_{\mathbf{n}}*D_{N}\right)\left(\mathfrak{z}\right)=\int_{\mathbb{Z}_{p}}D_{N}\left(\mathfrak{z}-\mathfrak{y}\right)dY_{\mathbf{n}}\left(\mathfrak{y}\right)=\sum_{\left|t\right|_{p}\leq p^{N}}\hat{Y}_{\mathbf{n}}\left(t\right)e^{2\pi i\left\{ t\mathfrak{z}\right\} _{p}}
\end{equation}
By inverting $1-L_{\mathbf{n}}$ and using our formula (\ref{eq:inverting dY_n})
for the action of $dY_{\mathbf{n}}$ on a given $\psi$, we have:
\begin{equation}
\int_{\mathbb{Z}_{p}}\psi\left(\mathfrak{y}\right)dY_{\mathbf{n}}\left(\mathfrak{y}\right)=\int_{\mathbb{Z}_{p}}\sum_{\mathbf{m}<\mathbf{n}}\left(L_{\mathbf{m},\mathbf{n}}\circ\left(1-L_{\mathbf{n}}\right)^{-1}\right)\left\{ \psi\right\} \left(\mathfrak{y}\right)d\chi_{\mathbf{m}}\left(\mathfrak{y}\right)
\end{equation}
So, fixing $\mathfrak{z}\in\mathbb{Z}_{p}$, set $\psi\left(\mathfrak{y}\right)=D_{N}\left(\mathfrak{z}-\mathfrak{y}\right)$.
Then:
\begin{align}
\sum_{\left|t\right|_{p}\leq p^{N}}\hat{Y}_{\mathbf{n}}\left(t\right)e^{2\pi i\left\{ t\mathfrak{z}\right\} _{p}} & =\int_{\mathbb{Z}_{p}}\sum_{\mathbf{m}<\mathbf{n}}\left(L_{\mathbf{m},\mathbf{n}}\circ\left(1-L_{\mathbf{n}}\right)^{-1}\right)\left\{ D_{N}\left(\mathfrak{z}-\cdot\right)\right\} \left(\mathfrak{y}\right)d\chi_{\mathbf{m}}\left(\mathfrak{y}\right)\nonumber \\
 & =\sum_{\mathbf{m}<\mathbf{n}}\sum_{k=0}^{p-1}\frac{r_{\mathbf{m},\mathbf{n},k}}{p}\int_{\mathbb{Z}_{p}}\left(1-L_{\mathbf{n}}\right)^{-1}\left\{ D_{N}\left(\mathfrak{z}-\cdot\right)\right\} \left(p\mathfrak{y}+k\right)d\chi_{\mathbf{m}}\left(\mathfrak{y}\right)\label{eq:tells us}
\end{align}
Noting that:
\begin{equation}
D_{N}\left(\mathfrak{z}-\mathfrak{y}\right)=p^{N}\left[\mathfrak{z}\overset{p^{N}}{\equiv}\mathfrak{y}\right]=p^{N}\left[\mathfrak{y}\overset{p^{N}}{\equiv}\left[\mathfrak{z}\right]_{p^{N}}\right]
\end{equation}
we can use the formula from \textbf{Lemma \ref{lem:Inversion of L_nu}
}to directly compute $\left(1-L_{\mathbf{n}}\right)^{-1}\left\{ D_{N}\left(\mathfrak{z}-\cdot\right)\right\} \left(p\mathfrak{y}+k\right)$.
This is:
\begin{equation}
\left(1-L_{\mathbf{n}}\right)^{-1}\left\{ D_{N}\left(\mathfrak{z}-\cdot\right)\right\} \left(p\mathfrak{y}+k\right)=\frac{r_{\mathbf{n},0}^{N}\kappa_{\mathbf{n}}\left(\left[\mathfrak{z}\right]_{p^{N}}\right)/p^{N}}{1-\alpha_{\mathbf{n}}\left(0\right)}+\sum_{n=0}^{N-1}\frac{r_{\mathbf{n},0}^{n}\kappa_{\mathbf{n}}\left(\left[\mathfrak{z}\right]_{p^{n}}\right)}{p^{n}}\left[\mathfrak{z}\overset{p^{N-n}}{\equiv}\theta_{p}^{\circ n}\left(\mathfrak{z}\right)\right]\label{eq:L_nu inverse of D_N of py plus j}
\end{equation}
We can simplify this further by noting that $0\leq n\leq N-1$ allows
us to write:
\begin{equation}
\left[p\mathfrak{y}+k\overset{p^{N-n}}{\equiv}\theta_{p}^{\circ n}\left(\mathfrak{z}\right)\right]=\left[\theta_{p}^{\circ n}\left(\mathfrak{z}\right)\overset{p}{\equiv}k\right]\left[p\mathfrak{y}+k\overset{p^{N-n}}{\equiv}\theta_{p}^{\circ n}\left(\mathfrak{z}\right)\right]
\end{equation}
Seeing as the bracket $\left[\theta_{p}^{\circ n}\left(\mathfrak{z}\right)\overset{p}{\equiv}k\right]$
guarantees that $\theta_{p}^{\circ n}\left(\mathfrak{z}\right)-k$
is a multiple of $p$, the expression on the far right can be written
as:
\begin{equation}
\left[p\mathfrak{y}+k\overset{p^{N-n}}{\equiv}\theta_{p}^{\circ n}\left(\mathfrak{z}\right)\right]=\left[\mathfrak{y}\overset{p^{N-n-1}}{\equiv}\frac{\theta_{p}^{\circ n}\left(\mathfrak{z}\right)-k}{p}\right]=\underbrace{\left[\mathfrak{y}\overset{p^{N-n-1}}{\equiv}\theta_{p}^{\circ n+1}\left(\mathfrak{z}\right)\right]}_{D_{N-n-1}\left(\mathfrak{y}-\theta_{p}^{\circ n+1}\left(\mathfrak{z}\right)\right)/p^{N-n-1}}
\end{equation}
and so:
\begin{equation}
\left[p\mathfrak{y}+k\overset{p^{N-n}}{\equiv}\theta_{p}^{\circ n}\left(\mathfrak{z}\right)\right]=\frac{\left[\theta_{p}^{\circ n}\left(\mathfrak{z}\right)\overset{p}{\equiv}k\right]}{p^{N-n-1}}D_{N-n-1}\left(\mathfrak{y}-\theta_{p}^{\circ n+1}\left(\mathfrak{z}\right)\right)
\end{equation}
Plugging this into (\ref{eq:L_nu inverse of D_N of py plus j}) gives:
\begin{align*}
\left(1-L_{\mathbf{n}}\right)^{-1}\left\{ D_{N}\left(\cdot-\mathfrak{z}\right)\right\} \left(p\mathfrak{y}+k\right) & =\frac{r_{\mathbf{n},0}^{N}\kappa_{\mathbf{n}}\left(\left[\mathfrak{z}\right]_{p^{N}}\right)/p^{N}}{1-\alpha_{\mathbf{n}}\left(0\right)}\\
 & +\sum_{n=0}^{N-1}\frac{r_{\mathbf{n},0}^{n}\kappa_{\mathbf{n}}\left(\left[\mathfrak{z}\right]_{p^{n}}\right)}{p^{n-N}}\frac{\left[\theta_{p}^{\circ n}\left(\mathfrak{z}\right)\overset{p}{\equiv}k\right]}{p^{N-n-1}}D_{N-n-1}\left(\mathfrak{y}-\theta_{p}^{\circ n+1}\left(\mathfrak{z}\right)\right)
\end{align*}
which simplifies to:

\begin{align}
\left(1-L_{\mathbf{n}}\right)^{-1}\left\{ D_{N}\left(\cdot-\mathfrak{z}\right)\right\} \left(p\mathfrak{y}+k\right) & =\frac{r_{\mathbf{n},0}^{N}\kappa_{\mathbf{n}}\left(\left[\mathfrak{z}\right]_{p^{N}}\right)/p^{N}}{1-\alpha_{\mathbf{n}}\left(0\right)}\label{eq:L_nu inverse of D_N of py plus j, simplified}\\
 & +p\sum_{n=0}^{N-1}r_{\mathbf{n},0}^{n}\kappa_{\mathbf{n}}\left(\left[\mathfrak{z}\right]_{p^{n}}\right)\left[\theta_{p}^{\circ n}\left(\mathfrak{z}\right)\overset{p}{\equiv}k\right]D_{N-n-1}\left(\mathfrak{y}-\theta_{p}^{\circ n+1}\left(\mathfrak{z}\right)\right)\nonumber 
\end{align}

Equation (\ref{eq:tells us}) tells us that if we integrate this against
$d\chi_{\mathbf{m}}\left(\mathfrak{y}\right)$, multiply the result
by $r_{\mathbf{m},\mathbf{n},\mathbf{k}}/p$, and then sum over all
$k\in\left\{ 0,\ldots,p-1\right\} $ and $\mathbf{m}<\mathbf{n}$,
we will then obtain: 
\begin{equation}
\sum_{\left|t\right|_{p}\leq p^{N}}\hat{Y}_{\mathbf{n}}\left(t\right)e^{2\pi i\left\{ t\mathfrak{z}\right\} _{p}}
\end{equation}
Since the only part of (\ref{eq:L_nu inverse of D_N of py plus j, simplified})
that depends on $\mathfrak{y}$ is the Dirichlet Kernel, we have that:
\begin{equation}
\int_{\mathbb{Z}_{p}}\left(1-L_{\mathbf{n}}\right)^{-1}\left\{ D_{N}\left(\cdot-\mathfrak{z}\right)\right\} \left(p\mathfrak{y}+k\right)d\chi_{\mathbf{m}}\left(\mathfrak{y}\right)
\end{equation}
is:
\begin{align*}
 & =\frac{r_{\mathbf{n},0}^{N}\kappa_{\mathbf{n}}\left(\left[\mathfrak{z}\right]_{p^{N}}\right)/p^{N}}{1-\alpha_{\mathbf{n}}\left(0\right)}\overbrace{\int_{\mathbb{Z}_{p}}d\chi_{\mathbf{m}}\left(\mathfrak{y}\right)}^{\hat{\chi}_{\mathbf{m}}\left(0\right)}\\
 & +p\sum_{n=0}^{N-1}r_{\mathbf{n},0}^{n}\kappa_{\mathbf{n}}\left(\left[\mathfrak{z}\right]_{p^{n}}\right)\left[\theta_{p}^{\circ n}\left(\mathfrak{z}\right)\overset{p}{\equiv}k\right]\underbrace{\int_{\mathbb{Z}_{p}}D_{N-n-1}\left(\mathfrak{y}-\theta_{p}^{\circ n+1}\left(\mathfrak{z}\right)\right)d\chi_{\mathbf{m}}\left(\mathfrak{y}\right)}_{=\tilde{\chi}_{\mathbf{m},N-n-1}\left(\theta_{p}^{\circ n+1}\left(\mathfrak{z}\right)\right)}
\end{align*}
Multiplying by $r_{\mathbf{m},\mathbf{n},k}/p$ and summing over all
$k\in\left\{ 0,\ldots,p-1\right\} $ and $\mathbf{m}<\mathbf{n}$
gives:
\begin{prop}
\label{prop:L_n inverse of D_N of py plus k}We have:
\begin{align*}
\sum_{\left|t\right|_{p}\leq p^{N}}\hat{Y}_{\mathbf{n}}\left(t\right)e^{2\pi i\left\{ t\mathfrak{z}\right\} _{p}} & =\frac{r_{\mathbf{n},0}^{N}\kappa_{\mathbf{n}}\left(\left[\mathfrak{z}\right]_{p^{N}}\right)/p^{N}}{1-\alpha_{\mathbf{n}}\left(0\right)}\sum_{\mathbf{m}<\mathbf{n}}r_{\mathbf{m},\mathbf{n},k}\hat{\chi}_{\mathbf{m}}\left(0\right)\\
 & +\sum_{n=0}^{N-1}r_{\mathbf{n},0}^{n}\kappa_{\mathbf{n}}\left(\left[\mathfrak{z}\right]_{p^{n}}\right)\sum_{\mathbf{m}<\mathbf{n}}r_{\mathbf{m},\mathbf{n},\left[\theta_{p}^{\circ n}\left(\mathfrak{z}\right)\right]_{p}}\tilde{\chi}_{\mathbf{m},N-n-1}\left(\theta_{p}^{\circ n+1}\left(\mathfrak{z}\right)\right)
\end{align*}
\end{prop}
Letting $N\rightarrow\infty$, the sum on the lower line is precisely
the same main limit we dealt with in earlier sections. As such, we
immediately obtain:
\begin{thm}
\label{thm:formal solutions}Let $\mathcal{F}$ be a frame as given
in \textbf{Corollary \ref{cor:main result}}. Let $I\subseteq R_{d}$
be a unique solution ideal, let $\mathbf{n}\in\textrm{CBI}_{d}\left(I\right)$,
for each $\mathbf{m}<\mathbf{n}$, let $\hat{X}_{\mathbf{m}}$ be
a Fourier transform of $X_{\mathbf{m}}$ with respect to $\mathcal{F}$,
and let $dY_{\mathbf{n}}$ be the distribution obtained by solving
the $\mathbf{n}$th formal equation relative to the initial condition
$\left\{ \hat{X}_{\mathbf{m}}\right\} _{\mathbf{m}\in\left[0,\mathbf{n}\right)}$
in the manner described in \textbf{Theorem \ref{thm:inversion method}}.
Then, $\hat{Y}_{\mathbf{n}}$ is a Fourier transform of $X_{\mathbf{n}}$
with respect to $\mathcal{F}$.
\end{thm}
Proof: We already did this in the course of our paper-long proof of
\textbf{Corollary \ref{cor:main result}}.

Q.E.D.

As a corollary, we then get:
\begin{cor}
Fix a unique solution ideal $I\subseteq R_{d}$, let $\mathcal{F}$
be an $\mathcal{R}_{d}$-frame as given in \textbf{Theorem \ref{thm:variety frame general}}.
Then, $\textrm{BI}_{d}\left(I\right)$ is the set of $\mathbf{n}\in\mathbb{N}_{0}^{d}\backslash\left\{ \mathbf{0}\right\} $
so that, given any choice of Fourier transforms $\left\{ \hat{X}_{\mathbf{m}}\right\} _{\mathbf{m}\in\left[\mathbf{0},\mathbf{n}\right)}$,
there is no Fourier transform of $X_{\mathbf{n}}$ satisfying the
$\mathbf{n}$th functional equation with initial condition $\left\{ \hat{X}_{\mathbf{m}}\right\} _{\mathbf{m}\in\left[\mathbf{0},\mathbf{n}\right)}$.
\end{cor}
In summary, given $\mathbf{n}$, provided we've already solved for
$\hat{X}_{\mathbf{m}}$ for all $\mathbf{m}<\mathbf{n}$, we can use
the formal equation to solve for $\hat{X}_{\mathbf{n}}$ if and only
if $\mathbf{n}\in\textrm{CBI}_{d}\left(I\right)$. When $\mathbf{n}\in\textrm{BI}_{d}\left(I\right)$,
the $1-L_{\mathbf{n}}$ inversion method breaks down, and we must
appeal to \textbf{Corollary \ref{cor:main result}} to obtain a Fourier
transform of $\hat{X}_{\mathbf{n}}$. This is just one of the many
ways that the breakdown variety plays a fundamental role in the Fourier-theoretic
behavior of the $X_{\mathbf{n}}$s.

\subsection{\label{subsec:Varieties-=000026-Measures}Varieties \& Distributions}

As mentioned, the purpose of this subsection is to illustrate some
of the basic aspects in which the distributions and measures induced
by F-series can ``encode'' algebraic varieties.
\begin{example}
\label{exa:first encoding example}The impetus for the idea is our
old friend, \textbf{Proposition \ref{prop:A_X sum}}. Recall, this
gave us the formula:
\begin{equation}
\sum_{\left|t\right|_{p}\leq p^{N}}\hat{A}_{X}\left(t\right)e^{2\pi i\left\{ t\mathfrak{z}\right\} _{p}}=a_{0}^{N}\kappa_{X}\left(\left[\mathfrak{z}\right]_{p^{N}}\right)+\left(1-\alpha_{X}\left(0\right)\right)\sum_{n=0}^{N-1}a_{0}^{n}\kappa_{X}\left(\left[\mathfrak{z}\right]_{p^{n}}\right)
\end{equation}
Choose a frame $\mathcal{F}$ so that $\lim_{n\rightarrow\infty}a_{0}^{n}\kappa_{X}\left(\left[\mathfrak{z}\right]_{p^{n}}\right)$
tends to $0$ nicely enough, we can take the limit as $N\rightarrow\infty$
and get the following:
\begin{equation}
\sum_{t\in\hat{\mathbb{Z}}_{p}}\hat{A}_{X}\left(t\right)e^{2\pi i\left\{ t\mathfrak{z}\right\} _{p}}=\left(1-\alpha_{X}\left(0\right)\right)\sum_{n=0}^{\infty}a_{0}^{n}\kappa_{X}\left(\left[\mathfrak{z}\right]_{p^{n}}\right)
\end{equation}
Using $\hat{A}_{X}$, we can then realize the F-series on the right-hand
side as a $\textrm{Frac}\left(\mathcal{R}_{1}\right)$-valued distribution
$dA_{X}\left(\mathfrak{z}\right)$ by defining its Fourier-Stieltjes
transform to be $\hat{A}_{X}\left(t\right)$. Given any unique solution
ideal $I\subseteq R_{1}$, the distribution $dA_{X}/I$ induced by
$dA_{X}$ in $\mathcal{S}\left(\mathbb{Z}_{p},\mathcal{R}_{1}/I\mathcal{R}_{1}\right)^{\prime}$
will be $\mathcal{F}$-degenerate if and only if $\alpha_{X}\left(0\right)=1$;
thus, its degeneracy detects whether or not the equality:
\begin{equation}
\frac{1}{p}\sum_{k=0}^{p-1}a_{k}=1\label{eq:alpha break}
\end{equation}
holds in $R_{1}/I$. In particular, $dA_{X}/I$ will be $\mathcal{F}$-degenerate
if and only if $\left\langle \alpha_{X}\left(0\right)-1\right\rangle \subseteq I$.
In this way, we can think of $dA_{X}$ as a map which, to every unique
solution ideal $I$ of $R$, outputs a distribution in $\mathcal{S}\left(\mathbb{Z}_{p},\textrm{Frac}\left(\mathcal{R}_{1}/I\mathcal{R}_{1}\right)\right)^{\prime}$.
\end{example}
\begin{rem}
As an aside, I have yet to develop the proper framework for treating
this map; is it a functor (say, with the natural projections $\mathcal{R}\rightarrow\mathcal{R}/I\mathcal{R}$
as morphisms)? A sheaf? Some kind of profinite (?) object indexed
by all the inclusions on the ideals of $R$ that we can quotient out
$\mathcal{R}$ by? Something else? I'm not quite sure, myself, as
this material has already reached the upper bound of the level of
abstraction that I'm currently comfortable with. Figuring out these
details is a natural and necessary direction for future research.
\end{rem}
That being said, \textbf{Example \ref{exa:first encoding example}}
has more to it than might be apparent at first glance.
\begin{example}
\label{exa:measures}Building on \textbf{Example \ref{exa:first encoding example}},
observe that for any absolute value $\ell$ on $\textrm{Frac}\left(\mathcal{R}_{1}/I\mathcal{R}_{1}\right)$,
we can consider the field $\mathbb{K}_{I,\ell}$ given by the metric
completion of the algebraic closure of the metric completion of $\textrm{Frac}\left(\mathcal{R}_{1}/I\mathcal{R}_{1}\right)$
with respect to $\ell$. In this set-up, the distributions in $\mathcal{S}\left(\mathbb{Z}_{p},\textrm{Frac}\left(\mathcal{R}_{1}/I\mathcal{R}_{1}\right)\right)^{\prime}$
which extend to measures in $W\left(\mathbb{Z}_{p},\mathbb{K}_{I,\rho}\right)^{\prime}$
are precisely those distributions $d\mu$ whose Fourier-Stieltjes
transforms $\hat{\mu}:\hat{\mathbb{Z}}_{p}\rightarrow\mathbb{K}_{I,\rho}$
which are $\ell$-adically bounded.

By construction:
\begin{equation}
\left\Vert \hat{A}_{X}\right\Vert _{p,\ell}=\sup_{t\in\hat{\mathbb{Z}}_{p}}\left|\hat{A}_{X}\left(t\right)\right|_{\ell}=\sup_{t\in\hat{\mathbb{Z}}_{p}}\prod_{n=0}^{-v_{p}\left(t\right)-1}\left|\alpha_{X}\left(p^{n}t\right)\right|_{\ell}
\end{equation}
Thus, for example, the boundedness of $\left\Vert \hat{A}_{X}\right\Vert _{p,\ell}$
will occur if $\left\Vert \alpha_{X}\right\Vert _{p,\ell}\leq1$.
Since we can always get absolute values on $\textrm{Frac}\left(\mathcal{R}_{1}/I\mathcal{R}_{1}\right)$
by evaluating all the indeterminates of $\mathcal{R}_{1}/I\mathcal{R}_{1}$
at elements of $K$ and then letting $\ell$ be an absolute value
on $K$, this shows that the algebras in which $dA_{X}$ induces a
measure relates to the absolute values of the left-hand side of (\ref{eq:alpha break})
under an evaluation map.
\end{example}
\begin{example}
\label{exa:morphisms}A natural way of extending the algebraic varieties
encodable in $dA_{X}$ is with the help of a morphism between varieties.
For example. By considering the change-of-variables $a_{j}=x_{j}^{d}$
for all $j\in\left\{ 0,\ldots,p-1\right\} $ for some $d\geq1$ we
will then have that $dA_{X}$ is $\mathcal{F}$-degenerate if and
only if:
\begin{equation}
\frac{1}{p}\sum_{j=0}^{p-1}x_{j}^{d}=1
\end{equation}
This leads to an interesting dilemma. As we've seen, the breakdown
variety of $\left(dA_{X}\right)^{d}$ is:
\begin{equation}
\frac{1}{p}\sum_{j=0}^{p-1}a_{j}^{d}=1\label{eq:d sphere}
\end{equation}
As such, both the image of $dA_{X}$ under the change of variables
$a_{j}=x_{j}^{d}$ and the $d$-fold pointwise product of $dA_{X}$
with itself ``encode'' the variety (\ref{eq:d sphere}). The main
difference between these two encodings is that while the first will
be degenerate if and only if we lie on a point of the variety, the
second need not be degenerate at all. Rather, as \textbf{Theorem \ref{thm:formal solutions}}
showed, the way in which $\left(dA_{X}\right)^{n}$ will manifest
sensitivity to (\ref{eq:d sphere}) is in whether or not a Fourier
transform of $\left(dA_{X}\right)^{n}$ is uniquely determined by
a given choice of Fourier transforms for $\left(dA_{X}\right)^{k}$
for $k\in\left\{ 0,\ldots,n-1\right\} $.
\end{example}
At the time of writing this paper, it is still not entirely clear
as to \emph{how }or to what extent the relationship between an F-series-induced
distribution/measure and its breakdown variety encodes the properties
of the breakdown variety. Ideally, one would like to establish a ``dictionary''
relating the distributions' properties to those of their breakdown
varieties. Our next few results are the first steps toward elucidating
this dictionary in the case where the distributions encode their breakdown
variety through their degeneracy. That being said, other than \textbf{Theorem
\ref{thm:formal solutions}}, the question is much less clear when
it comes to F-series which are not degenerate distributions when $I$
contains their breakdown variety. Figuring out how to characterize
the properties of an F-series induced distribution depend on the relation
between $I$ and the breakdown variety is currently a major open question.
On the plus side, at least we have something to work with in the degenerate
case.

For simplicity's sake, we need a new piece of terminology.
\begin{defn}
Let $R$ be a commutative, unital, integral domain whose characteristic
is either $0$ or co-prime to $p$, and let $\left\{ M_{n}\right\} _{n\geq0}$
be a $p$-adic $\textrm{Frac}\left(R\right)$-valued M-function, with
multipliers $r_{0},\ldots,r_{p-1}$. Then, we write $\alpha_{M},\hat{M}:\hat{\mathbb{Z}}_{p}\rightarrow\textrm{Frac}\left(R\right)\left(\zeta_{p^{\infty}}\right)$
to denote the functions:
\begin{align}
\alpha_{M}\left(t\right) & \overset{\textrm{def}}{=}\frac{1}{p}\sum_{k=0}^{p-1}r_{k}e^{-2\pi ikt}\\
\hat{M}\left(t\right) & \overset{\textrm{def}}{=}\prod_{n=0}^{-v_{p}\left(t\right)-1}\alpha_{M}\left(p^{n}t\right)
\end{align}
and write $dM\in\mathcal{S}\left(\mathbb{Z}_{p},\textrm{Frac}\left(R\right)\right)^{\prime}$
to denote the distribution whose Fourier-Stieltjes transform is $\hat{M}$.
We call $\hat{M}$ \textbf{the Fourier(-Stieltjes) transform of $\left\{ M_{n}\right\} _{n\geq0}$
}and call $dM$ \textbf{the distribution induced by $\left\{ M_{n}\right\} _{n\geq0}$}.

Next, we write $L_{M}:\mathcal{S}\left(\mathbb{Z}_{p},\textrm{Frac}\left(R\right)\right)\rightarrow\mathcal{S}\left(\mathbb{Z}_{p},\textrm{Frac}\left(R\right)\right)$
to denote the $R$-linear operator:
\begin{equation}
L_{M}\left\{ \phi\right\} \left(\mathfrak{z}\right)\overset{\textrm{def}}{=}\frac{1}{p}\sum_{k=0}^{p-1}a_{k}\phi\left(p\mathfrak{z}+k\right),\textrm{ }\forall\phi\in\mathcal{S}\left(\mathbb{Z}_{p},\textrm{Frac}\left(R\right)\right),\textrm{ }\forall\mathfrak{z}\in\mathbb{Z}_{p}
\end{equation}
We call $L_{M}$ the \textbf{sorting operator induced by }$\left\{ M_{n}\right\} _{n\geq0}$.

Letting $\ell$ be a (possibly trivial) place of $R$, and letting
$\overline{\overline{\textrm{Frac}\left(R_{\ell}\right)}}$ be the
metric completion of the algebraic closure of the field of fractions
of the metric completion of $R$ with respect to $\ell$, we say $dM$
is a $\left(p,\ell\right)$-adic measure whenever it extends to an
element of $W\left(\mathbb{Z}_{p},\overline{\overline{\textrm{Frac}\left(R_{\ell}\right)}}\right)^{\prime}$,
in which case, we denote this extension by $dM$ and call it the \textbf{$\left(p,\ell\right)$-adic
measure induced by $\left\{ M_{n}\right\} _{n\geq0}$}. Likewise,
we call $L_{M}$ the \textbf{$\left(p,\ell\right)$-adic sorting operator
induced by $\left\{ M_{n}\right\} _{n\geq0}$} whenever $L_{M}$ admits
an extension to a continuous linear operator $W\left(\mathbb{Z}_{p},\overline{\overline{\textrm{Frac}\left(R_{\ell}\right)}}\right)\rightarrow W\left(\mathbb{Z}_{p},\overline{\overline{\textrm{Frac}\left(R_{\ell}\right)}}\right)$.
\end{defn}
\begin{defn}
Given a linear operator $L:\mathcal{S}\left(\mathbb{Z}_{p},\textrm{Frac}\left(R\right)\right)\rightarrow\mathcal{S}\left(\mathbb{Z}_{p},\textrm{Frac}\left(R\right)\right)$,
recall that the \textbf{pullback }of $L$, denoted $L^{*}:\mathcal{S}\left(\mathbb{Z}_{p},\textrm{Frac}\left(R\right)\right)^{\prime}\rightarrow\mathcal{S}\left(\mathbb{Z}_{p},\textrm{Frac}\left(R\right)\right)^{\prime}$,
is the linear map given by the rule:
\begin{equation}
d\mu\mapsto L^{*}\left\{ d\mu\right\} ,\textrm{ }\forall d\mu\in\mathcal{S}\left(\mathbb{Z}_{p},\textrm{Frac}\left(R\right)\right)^{\prime}
\end{equation}
where $L^{*}\left\{ d\mu\right\} $ is the distribution defined by:
\begin{equation}
\int_{\mathbb{Z}_{p}}\phi\left(\mathfrak{z}\right)L^{*}\left\{ d\mu\right\} \left(\mathfrak{z}\right)\overset{\textrm{def}}{=}\int_{\mathbb{Z}_{p}}L\left\{ \phi\right\} \left(\mathfrak{z}\right)d\mu\left(\mathfrak{z}\right),\textrm{ }\forall\phi\in\mathcal{S}\left(\mathbb{Z}_{p},\textrm{Frac}\left(R\right)\right)
\end{equation}
\end{defn}
\begin{prop}
\label{prop:L_A and dA}We have:

\begin{equation}
\int_{\mathbb{Z}_{p}}\left(1-L_{M}\right)\left\{ \phi\right\} \left(\mathfrak{z}\right)dM\left(\mathfrak{z}\right)=\left(1-\alpha_{M}\left(0\right)\right)\hat{\phi}\left(0\right),\textrm{ }\forall\phi\in\mathcal{S}\left(\mathbb{Z}_{p},\textrm{Frac}\left(R\right)\right)\label{eq:1 - L dA}
\end{equation}
This identity also holds for all $\phi\in W\left(\mathbb{Z}_{p},\overline{\overline{\textrm{Frac}\left(R_{\ell}\right)}}\right)$
whenever both $dM$ and $L_{M}$ admit $\left(p,\ell\right)$-adic
extensions.

From this, we see that $\textrm{Im}\left(1-L_{M}\right)\subseteq\textrm{Ker}\left(dM\right)$
if and only if $\alpha_{M}\left(0\right)=1$. Equivalently, $dM$
is which is to say that $dM$ is an eigenvector (technically, eigendistribution)
of $L_{M}^{*}$ corresponding to the eigenvalue $1$ if and only if
$\alpha_{M}\left(0\right)=1$.
\end{prop}
Proof: Letting $\phi$ be arbitrary, we have:

\begin{align*}
\int_{\mathbb{Z}_{p}}\phi\left(\mathfrak{z}\right)dM\left(\mathfrak{z}\right) & =\sum_{t\in\hat{\mathbb{Z}}_{p}}\phi\left(t\right)\hat{M}\left(-t\right)\\
 & =\hat{\phi}\left(0\right)\hat{M}\left(0\right)+\sum_{t\in\hat{\mathbb{Z}}_{p}\backslash\left\{ 0\right\} }\phi\left(t\right)\hat{M}\left(-t\right)\\
\left(\hat{M}\left(t\right)=\alpha_{M}\left(t\right)\hat{M}\left(pt\right),\textrm{ }\forall t\neq0\right); & =\hat{\phi}\left(0\right)\hat{M}\left(0\right)+\underbrace{\sum_{t\in\hat{\mathbb{Z}}_{p}\backslash\left\{ 0\right\} }\phi\left(t\right)\alpha_{M}\left(-t\right)\hat{M}\left(-pt\right)}_{\textrm{add \& subtract }t=0\textrm{ term}}\\
 & =\left(1-\alpha_{M}\left(0\right)\right)\hat{\phi}\left(0\right)\hat{M}\left(0\right)+\sum_{t\in\hat{\mathbb{Z}}_{p}}\phi\left(t\right)\alpha_{M}\left(-t\right)\hat{M}\left(-pt\right)
\end{align*}
Here, writing out $\alpha_{M}$ and writing $\hat{M}$ by its defining
Fourier integral, we get:
\begin{align*}
\sum_{t\in\hat{\mathbb{Z}}_{p}}\phi\left(t\right)\alpha_{M}\left(-t\right)\hat{M}\left(-pt\right) & =\sum_{t\in\hat{\mathbb{Z}}_{p}}\phi\left(t\right)\left(\frac{1}{p}\sum_{k=0}^{p-1}r_{k}e^{2\pi ikt}\right)\left(\int_{\mathbb{Z}_{p}}e^{2\pi i\left\{ pt\mathfrak{z}\right\} _{p}}dM\left(\mathfrak{z}\right)\right)\\
\left(\textrm{swap }\sum\textrm{ \& }\int\right); & \overset{!}{=}\int_{\mathbb{Z}_{p}}\left(\frac{1}{p}\sum_{k=0}^{p-1}r_{k}\sum_{t\in\hat{\mathbb{Z}}_{p}}\phi\left(t\right)e^{2\pi i\left\{ t\left(p\mathfrak{z}+k\right)\right\} _{p}}\right)dM\left(\mathfrak{z}\right)\\
 & \overset{!!}{=}\int_{\mathbb{Z}_{p}}\left(\frac{1}{p}\sum_{k=0}^{p-1}a_{k}\phi\left(p\mathfrak{z}+k\right)\right)dM\left(\mathfrak{z}\right)\\
 & =\int_{\mathbb{Z}_{p}}L_{M}\left\{ \phi\right\} \left(\mathfrak{z}\right)dM\left(\mathfrak{z}\right)
\end{align*}
Here, the interchange of sum and integral in the step marked (!) is
valid because $\phi$ is an SB function; thus, its Fourier transform
has finite support, and so the $t$-sum reduces to a finite sum. As
for the step marked (!!), again, because $\phi$ is an SB function,
its Fourier series reduces to a finite and necessarily convergent
sum. Thus, we see that:
\begin{align*}
\int_{\mathbb{Z}_{p}}\phi\left(\mathfrak{z}\right)dM\left(\mathfrak{z}\right) & =\left(1-\alpha_{M}\left(0\right)\right)\hat{\phi}\left(0\right)\underbrace{\hat{M}\left(0\right)}_{1}+\sum_{t\in\hat{\mathbb{Z}}_{p}}\phi\left(t\right)\alpha_{M}\left(-t\right)\hat{M}\left(-pt\right)\\
 & =\left(1-\alpha_{M}\left(0\right)\right)\hat{\phi}\left(0\right)+\int_{\mathbb{Z}_{p}}L_{M}\left\{ \phi\right\} \left(\mathfrak{z}\right)dM\left(\mathfrak{z}\right)
\end{align*}
Moving the integral on the far right to the left-hand side then yields
the result.

The equivalent statement involving the images and kernels then follow
from considering the composite of $1-L_{M}$ and $dM$ and using (\ref{eq:1 - L dA}).
As for the eigendistribution claim, note that:
\begin{equation}
\left(1-L_{M}\right)^{*}\left\{ dM\right\} \left(\mathfrak{z}\right)=\left(1-\alpha_{M}\left(0\right)\right)d\mathfrak{z}
\end{equation}
where $d\mathfrak{z}$ is the Haar distribution. Since the pullback
of the identity map is the identity map, we can re-write the above
as:
\begin{equation}
L_{M}^{*}\left\{ dM\right\} \left(\mathfrak{z}\right)=dM\left(\mathfrak{z}\right)-\left(1-\alpha_{M}\left(0\right)\right)d\mathfrak{z}
\end{equation}
and thus, we see that $dM$ is an eigendistribution of $L_{M}^{*}$
if and only if $\alpha_{M}\left(0\right)=1$.

Q.E.D.

\vphantom{}

Using \textbf{Proposition \ref{prop:A_X sum}}, for an appropriate
frame (or, if we just take everything as occurring in $\mathbb{Z}\left[\left[r_{0},\ldots,r_{p-1}\right]\right]$),
we have:
\begin{equation}
\sum_{t\in\hat{\mathbb{Z}}_{p}}\hat{M}\left(t\right)e^{2\pi i\left\{ t\mathfrak{z}\right\} _{p}}\overset{\mathcal{F}}{=}\left(1-\alpha_{M}\left(0\right)\right)\sum_{n=0}^{\infty}M_{n}\left(\mathfrak{z}\right)
\end{equation}
Letting $M$ denote the F-series on the right, provided we do not
lie in the breakdown variety of $M$ (i.e., $\alpha_{M}\left(0\right)\neq1$),
$M$ is characterized by:
\begin{equation}
M\left(2\mathfrak{z}+k\right)=r_{k}M\left(\mathfrak{z}\right)+1-\alpha_{M}\left(0\right),\textrm{ }\forall\mathfrak{z}\in\mathbb{Z}_{p},\textrm{ }\forall k\in\left\{ 0,\ldots,p-1\right\} 
\end{equation}
Proceeding formally, fixing an ideal $I$ to quotient by, we then
have that the breakdown variety of $M^{n}$ is:
\begin{equation}
\sum_{k=0}^{p-1}r_{k}^{n}=p\label{eq:nth breakdown variety of M}
\end{equation}
where the equality holds mod $I$. By \textbf{Theorem \ref{thm:formal solutions}},
since $\hat{M}$ is already chosen, we have that formally solving
for $\hat{M}^{*n}$ will produce a unique solution for all $n\geq2$
until we reach $n_{1}$, the first integer $n\geq2$ for which (\ref{eq:nth breakdown variety of M})
holds true, at which point we will need to appeal to \textbf{Corollary
\ref{cor:main result} }to obtain a formula for $\hat{M}^{*n_{1}}$.
However, we will then be able to use $\left\{ \hat{M}^{*n}\right\} _{0\leq n\leq n_{1}}$
(with $\hat{}$ $\hat{M}^{*0}=\mathbf{1}_{0}$) as an initial condition
to formally solve for $\hat{M}^{*n}$ for all $n>n_{1}$ until we
hit $n_{2}$, the second integer $\geq2$ for which (\ref{eq:nth breakdown variety of M})
holds true, and so on and so forth. That being said, note that for
fixed $r_{0},\ldots,r_{p-1}\in K$, there can be at most $p$ values
of $n$ for which (\ref{eq:nth breakdown variety of M}) holds true.

Multiplication of different F-series yields a kind of ``combination''
of the parameters, and hence, of their breakdown varieties. For example,
given M-functions $M_{1,n},\ldots,M_{d,n}$ with multipliers $r_{1,k},\ldots,r_{d,k}$,
we have that the breakdown variety of $M_{1,n}^{e_{1}}\cdots M_{d,n}^{e_{d}}$
(for integers $e_{1},\ldots,e_{d}\geq0$), is:
\begin{equation}
\sum_{k=0}^{p-1}r_{1,k}^{e_{1}}\cdots r_{d,k}^{e_{d}}=p
\end{equation}

As a final observation, we can use the notion of the tensor product
of distributions/measures. Given $d\mu_{1}\left(\mathfrak{z}_{1}\right)$
and $d\mu_{2}\left(\mathfrak{z}_{2}\right)$ that act on functions
of $\mathfrak{x}_{1}\in\mathbb{Z}_{p_{1}}$ and $\mathfrak{x}_{2}\in\mathbb{Z}_{p_{2}}$,
respectively, recall that we define $\left(d\mu_{1}\otimes d\mu_{2}\right)\left(\mathfrak{z}_{1},\mathfrak{z}_{2}\right)$
as the distribution/measure given by:
\begin{equation}
\phi\mapsto\int_{\mathbb{Z}_{p_{1}}\times\mathbb{Z}_{p_{2}}}\phi\left(\mathfrak{z}_{1},\mathfrak{z}_{2}\right)d\mu_{1}\left(\mathfrak{z}_{1}\right)d\mu_{2}\left(\mathfrak{z}_{2}\right)
\end{equation}
We then have that the Fourier(-Stieltjes) transform distributes across
tensor products, with $\widehat{d\mu_{1}\otimes d\mu_{2}}$ being
a function of $\left(t_{1},t_{2}\right)\in\hat{\mathbb{Z}}_{p_{1}}\times\hat{\mathbb{Z}}_{p_{2}}$
given by:
\begin{align*}
\widehat{d\mu_{1}\otimes d\mu_{2}} & =\int_{\mathbb{Z}_{p_{1}}\times\mathbb{Z}_{p_{2}}}e^{-2\pi i\left\{ t_{1}\mathfrak{z}_{1}+t_{2}\mathfrak{z}_{2}\right\} _{p_{1}\times p_{2}}}d\mu_{1}\left(\mathfrak{z}_{1}\right)d\mu_{2}\left(\mathfrak{z}_{2}\right)\\
 & =\int_{\mathbb{Z}_{p_{1}}\times\mathbb{Z}_{p_{2}}}e^{-2\pi i\left\{ t_{1}\mathfrak{z}_{1}\right\} _{p_{1}}}e^{-2\pi i\left\{ t_{2}\mathfrak{z}_{2}\right\} _{p_{2}}}d\mu_{1}\left(\mathfrak{z}_{1}\right)d\mu_{2}\left(\mathfrak{z}_{2}\right)\\
 & =\int_{\mathbb{Z}_{p_{1}}}e^{-2\pi i\left\{ t_{1}\mathfrak{z}_{1}\right\} _{p_{1}}}d\mu_{1}\left(\mathfrak{z}_{1}\right)\int_{\mathbb{Z}_{p_{2}}}e^{-2\pi i\left\{ t_{2}\mathfrak{z}_{2}\right\} _{p_{2}}}d\mu_{2}\left(\mathfrak{z}_{2}\right)\\
 & =\hat{\mu}_{1}\left(t_{1}\right)\hat{\mu}_{2}\left(t_{2}\right)
\end{align*}
In the equality:
\begin{equation}
e^{-2\pi i\left\{ t_{1}\mathfrak{z}_{1}+t_{2}\mathfrak{z}_{2}\right\} _{p_{1}\times p_{2}}}=e^{-2\pi i\left\{ t_{1}\mathfrak{z}_{1}\right\} _{p_{1}}}e^{-2\pi i\left\{ t_{2}\mathfrak{z}_{2}\right\} _{p_{2}}}
\end{equation}
the expression on the left is a character of the group $\mathbb{Z}_{p_{1}}\times\mathbb{Z}_{p_{2}}$,
and standard Pontryagin duality tells us that a character of a direct
product of locally compact abelian groups is the direct sum of characters
of the individual groups; in this case, this equality bears witness
to the isomorphism:
\begin{equation}
\widehat{\mathbb{Z}_{p_{1}}\times\mathbb{Z}_{p_{2}}}\cong\hat{\mathbb{Z}}_{p_{1}}\oplus\hat{\mathbb{Z}}_{p_{2}}
\end{equation}
Because of this, we have:
\begin{prop}
\label{prop:tensors}Fix $d\geq2$, integers $p_{1},\ldots,p_{d}\geq2$,
and for each $m\in\left\{ 1,\ldots,d\right\} $, let $\left\{ M_{m,n}\right\} _{n\geq0}$
be a $p_{m}$-adic M-function with Fourier transform $\hat{M}_{m}$
and sorting operator $L_{m}$ and $\alpha$ function $\alpha_{m}$.
Then, for any $\phi\in\mathcal{S}\left(\prod_{m=1}^{d}\mathbb{Z}_{p_{m}},\textrm{Frac}\left(R\right)\right)$:
\begin{equation}
\int_{\prod_{m=1}^{d}\mathbb{Z}_{p_{m}}}\left(\prod_{m=1}^{d}\left(1-L_{m}\right)\right)\left\{ \phi\right\} \left(\mathfrak{z}_{1},\ldots,\mathfrak{z}_{d}\right)dM_{1}\left(\mathfrak{z}_{1}\right)\cdots dM_{d}\left(\mathfrak{z}_{d}\right)=\hat{\phi}\left(\mathbf{0}\right)\prod_{m=1}^{d}\left(1-\alpha_{m}\left(0\right)\right)
\end{equation}
where:
\begin{equation}
\hat{\phi}\left(\mathbf{0}\right)=\int_{\prod_{m=1}^{d}\mathbb{Z}_{p_{m}}}\phi\left(\mathbf{z}\right)d\mathbf{z}
\end{equation}
where $\mathbf{z}=\left(\mathfrak{z}_{1},\ldots,\mathfrak{z}_{d}\right)$
and $d\mathbf{z}=d\mathfrak{z}_{1}\cdots d\mathfrak{z}_{d}$. That
is to say, take the image of $\phi$ under:
\begin{equation}
\bigotimes_{m=1}^{d}\left(1-L_{m}\right):\bigotimes_{m=1}^{d}\mathcal{S}\left(\mathbb{Z}_{p_{m}},\textrm{Frac}\left(R\right)\right)\rightarrow\bigotimes_{m=1}^{d}\mathcal{S}\left(\mathbb{Z}_{p_{m}},\textrm{Frac}\left(R\right)\right)
\end{equation}
and then apply $\bigotimes_{m=1}^{d}dM_{m}$ to that image. Then,
the result will be the product of the integral of $\phi$ and $\prod_{m=1}^{d}\left(1-\alpha_{m}\left(0\right)\right)$.
That is, if $V_{1},\ldots,V_{d}$ are the breakdown varieties of the
distributions $dM_{1},\ldots,dM_{d}$, then $V_{1}\cup\cdots\cup V_{d}$
is the breakdown variety of the distribution $dM_{1}\otimes\cdots\otimes dM_{d}$.
Furthermore, the following are equivalent:

\vphantom{}

I. $\alpha_{m}\left(0\right)=1$ for some $m\in\left\{ 1,\ldots,d\right\} $,
and for any $\textrm{Frac}\left(R\right)$-frame $\mathcal{F}$ on
$\mathbb{Z}_{p}$ so that: 
\begin{equation}
\lim_{N\rightarrow\infty}\sum_{\left|t\right|_{p}\leq p^{N}}\hat{M}_{m}\left(t\right)e^{2\pi i\left\{ t\mathfrak{z}\right\} _{p}}
\end{equation}
is $\mathcal{F}$-convergent, $dM_{m}$ is an $\mathcal{F}$-degenerate
distribution. If, in addition, $\left\Vert \hat{M}_{m}\right\Vert _{p,\ell}<\infty$,
then $dM_{m}$ is an $\mathcal{F}$-degenerate measure in $W\left(\mathbb{Z}_{p},\overline{\overline{\textrm{Frac}\left(R_{\ell}\right)}}\right)^{\prime}$.

\vphantom{}

III. The image of $\bigotimes_{m=1}^{d}\mathcal{S}\left(\mathbb{Z}_{p_{m}},\textrm{Frac}\left(R\right)\right)$
under $\bigotimes_{m=1}^{d}\textrm{Im}\left(1-L_{m}\right)$ is contained
in the kernel of $\bigotimes_{m=1}^{d}\left(dM_{m}\mid_{\textrm{Im}\left(1-L_{m}\right)}\right)$.
\end{prop}
Thus, we see that tensor products of measures can be used to encode
unions of algebraic varieties. This appears to only be scratching
the surface of what can be done. Hopefully, future work (by myself,
and/or others) will help clarify and shed further light on this interesting
new correspondence.

\subsection{\label{subsec:Arithmetic-Dynamics,-Varieties,}Arithmetic Dynamics,
Measures, and Varieties}

In ancient Rome, Janus, the god of doorways (and the namesake of the
month of January) was commonly depicted as a man with two heads exhibiting
symmetry with respect to a reflection across a sagittal plane, so
as to illustrate how the god's gaze fell on a doorway's traveler regardless
of whether they were coming or going, entering or leaving. This same
Januarian duality is at work in this, the final section of my paper.
As readers familiar with my work will know, my original impetus for
all of this was my investigations into the infamous Collatz Conjecture
and the behavior of related arithmetic dynamical systems. In the process,
however, I seem to have stumbled my way into arithmetic geometry.
In that respect, just like Janus, this paper's work looks in two directions
simultaneously. One: how might we use the unexpectedly robust algebraic
structure of the algebras of measures induced by F-series to study
algebraic varieties? Two: in what ways does the relation between F-series
and algebraic varieties have bearing on the dynamics of Collatz-type
maps that generated F-series in the first place?

For an odd number $q\geq1$, we can generalize the Collatz map by
considering the \textbf{shortened $qx+1$ map }$T_{q}:\mathbb{Z}\rightarrow\mathbb{Z}$
defined by:
\begin{equation}
T_{q}\left(n\right)\overset{\textrm{def}}{=}\begin{cases}
\frac{n}{2} & \textrm{if }n\overset{2}{\equiv}0\\
\frac{qn+1}{2} & \textrm{if }n\overset{2}{\equiv}1
\end{cases}
\end{equation}
(I call $n/2$ the \textbf{even branch }of $T_{q}$ and call $\left(qn+1\right)/2$
the \textbf{odd branch}).\textbf{ }Note that $T_{q}$ then has a unique
extension $T_{q}:\mathbb{Z}_{2}\rightarrow\mathbb{Z}_{2}$. To $T_{q}$,
we associate the function, called the \textbf{numen},\textbf{ }$\chi_{q}:\mathbb{Z}_{2}\rightarrow\mathbb{Z}_{q}$
characterized by:
\begin{align}
\chi_{q}\left(2\mathfrak{z}\right) & =\frac{\chi_{q}\left(\mathfrak{z}\right)}{2}\\
\chi_{q}\left(2\mathfrak{z}+1\right) & =\frac{q\chi_{q}\left(\mathfrak{z}\right)+1}{2}
\end{align}
$\chi_{q}$ admits an F-series representation:
\begin{equation}
\chi_{q}\left(\mathfrak{z}\right)=-\frac{1}{q-1}+\frac{1}{2\left(q-1\right)}\sum_{n=0}^{\infty}\frac{q^{\#_{1}\left(\left[\mathfrak{z}\right]_{2^{n}}\right)}}{2^{n}}\label{eq:chi_q}
\end{equation}
which converges with respect to the frame $\mathcal{F}_{2,q}$ which
associates the real topology to $\mathfrak{z}\in\mathbb{N}_{0}$ and
associates the $q$-adic topology to $\mathfrak{z}\in\mathbb{Z}_{2}^{\prime}$.
As proven in \cite{My Dissertation,first blog paper}, $\chi_{q}$
satisfies a \textbf{Correspondence Principle} linking its value distribution
to the dynamics of $T_{q}$:
\begin{thm}[The Correspondence Principle for $T_{q}$\cite{first blog paper,My Dissertation}]

I. Let $x\in\mathbb{Z}_{2}$. Then, $x$ is a periodic point of $T_{q}$
if and only if\emph{ there exists a $\mathfrak{z}\in\mathbb{Z}_{2}\cap\mathbb{Q}$
so that $\chi_{q}\left(\mathfrak{z}\right)=x$.}

II. If there is a $\mathfrak{z}\in\mathbb{Z}_{2}\backslash\mathbb{Q}$
so that $\chi_{q}\left(\mathfrak{z}\right)\in\mathbb{Z}$, then $\chi_{q}\left(\mathfrak{z}\right)$
is a divergent point of $T_{q}$ (i.e., the iterates of $\chi_{q}\left(\mathfrak{z}\right)$
under $T_{q}$ tend to either $+\infty$ or $-\infty$)
\end{thm}
\begin{rem}
I have not yet been able to establish a full converse of (II); i.e.,
I have not been able to prove that any divergent points of $T_{q}$
in $\mathbb{Z}$ are the image of $\chi_{q}$ at some irrational $2$-adic
integer.
\end{rem}
More generally, given any hydra map $H:\mathbb{Z}\rightarrow\mathbb{Z}$
(the name I have given to the class of Collatz-type maps I am studying)
with $p$ branches ($T_{q}$ has $2$ branches), $H$ has a numen
$\chi_{H}:\mathbb{N}_{0}\rightarrow\mathbb{Q}$, and if there is a
$\mathbb{Q}$-frame $\mathcal{F}$ on $\mathbb{Z}_{p}$ with, say,
$D\left(\mathcal{F}\right)=\mathbb{Z}_{p}$ so that $\chi_{H}$ extends
to an element of $C\left(\mathcal{F}\right)$, $\chi_{H}$ will satisfy
a corresponding Correspondence Principle relative to the dynamics
of $H$. (The frame $\mathcal{F}$ doesn't need to have $D\left(\mathcal{F}\right)=\mathbb{Z}_{p}$
for this to hold in the case of periodic points, but I digress. See
\cite{my frames paper} for details.) As discussed in \cite{first blog paper},
the Correspondence Principle even extends to hydra maps on the ring
of integers of any global field, though my dissertation \cite{My Dissertation}
only considered the number field case; thankfully, the arguments given
there verbatim over function fields. For reference, given a global
field $K$, a hydra map on $\mathcal{O}_{K}$ consists of the datum:
\begin{itemize}
\item A non-trivial proper ideal $I\subset\mathcal{O}_{K}$.
\item For each $j\in\mathcal{O}_{K}/I$, an invertible endomorphism $r_{j}\in\left(\textrm{End}K\right)^{\times}$
and an element $c_{j}\in K$ and an invertible affine-linear map $H_{j}:K\rightarrow K$
(called the \textbf{$j$th branch}) given by $H_{j}\left(z\right)=r_{j}\left(z\right)+c$.
For any $z\in\mathcal{O}_{K}$, we require $H_{j}\left(z\right)\in\mathcal{O}_{K}$
if and only if $z=j$ in $\mathcal{O}_{K}/I$.
\item A map (called the \textbf{hydra})\textbf{ }$H:\mathcal{O}_{K}\rightarrow\mathcal{O}_{K}$
defined by $H\left(z\right)=H_{\left[z\right]_{I}}\left(z\right)$
for all $z\in\mathcal{O}_{K}$, where $\left[z\right]_{I}$ is the
image of $z$ under the projection $\mathcal{O}_{K}\rightarrow\mathcal{O}_{K}/I$.
\end{itemize}
Returning to the case of $\chi_{q}$, we see that $p=2$, $a_{0}=1/2$,
$a_{1}=q/2$, $b_{0}=0$, $b_{1}=1/2$. To treat $q$ as an indeterminate,
we start from: 
\begin{equation}
R_{1}=\mathbb{Z}\left[a_{0},a_{1},b_{0},b_{1}\right]
\end{equation}
\begin{equation}
\mathcal{R}_{1}=\mathbb{Z}\left[\frac{1}{1-a_{0}}\right]\left[a_{0},a_{1},b_{0},b_{1}\right]
\end{equation}
corresponding to:
\begin{align}
X\left(2\mathfrak{z}\right) & =a_{0}X\left(\mathfrak{z}\right)+b_{0}\\
X\left(2\mathfrak{z}\right) & =a_{1}X\left(\mathfrak{z}\right)+b_{1}
\end{align}
To get to $\chi_{q}$, we take the quotient:
\begin{equation}
\frac{\mathcal{R}_{1}}{\left\langle a_{0}-\frac{1}{2},b_{0},b_{1}-\frac{1}{2}\right\rangle }\cong\mathbb{Z}\left[\frac{1}{2},a_{1}\right]
\end{equation}
and then apply the homomorphism which sends $a_{1}$ to $q/2$, which
gives:
\begin{equation}
\chi_{q}:\mathbb{N}_{0}\rightarrow\mathbb{Z}\left[\frac{1}{2},q\right]
\end{equation}
We can get a frame $\mathcal{F}$ on $\mathbb{Z}\left[1/2,q\right]$
by setting: 
\begin{equation}
\mathcal{F}_{\mathfrak{z}}\left(f\right)=\begin{cases}
\left|f\left(\mathfrak{z}\right)\right|_{\infty} & \textrm{if }\mathfrak{z}\in\mathbb{N}_{0}\\
\left|f\left(\mathfrak{z}\right)\right|_{q} & \textrm{if }\mathfrak{z}\in\mathbb{Z}_{2}^{\prime}
\end{cases},\textrm{ }\forall f:\mathbb{Z}_{2}\rightarrow\mathbb{Z}\left[\frac{1}{2},q\right]
\end{equation}
where $\left|\cdot\right|_{\infty}$ is any archimedean absolute value
on $\mathbb{Z}\left[1/2,q\right]$ and $\left|\cdot\right|_{q}$ is
any absolute value generated by the prime ideal $\left\langle q\right\rangle $
of $\mathbb{Z}\left[\frac{1}{2},q\right]$. The ring $R=\mathbb{Z}\left[\frac{1}{2},q\right]$
will then be the source of our ideals $I$ corresponding to relations
on $\mathbb{Z}\left[\frac{1}{2},q\right]$. A unique solution ideal,
in this case, is any ideal $I\subseteq\mathbb{Z}\left[1/2,q\right]$
so that $\left\langle 1/2\right\rangle \nsubseteq I$. Since $\left\langle 1/2\right\rangle $
generates $\mathbb{Z}\left[1/2,q\right]$, this shows that any proper
ideal of $\mathbb{Z}\left[1/2,q\right]$ will be a unique solution
ideal.

Any Fourier transform of $\chi_{q}$ ends up taking values in $\textrm{Frac}\left(R\right)\left(\zeta_{2^{\infty}}\right)$,
which is $\mathbb{Q}\left(q,\zeta_{2^{\infty}}\right)$. The breakdown
variety of $\chi_{q}$ is:
\begin{equation}
\frac{q+1}{4}=1
\end{equation}
and we have the formulae:
\begin{equation}
\hat{\chi}_{q}\left(t\right)=\begin{cases}
\begin{cases}
0 & \textrm{if }t=0\\
\left(\frac{1}{4}v_{2}\left(t\right)+\frac{1}{2}\right)\hat{A}_{3}\left(t\right) & \textrm{if }t\neq0
\end{cases} & \textrm{if }q=3\\
\\
\begin{cases}
-\frac{1}{q-3} & \textrm{if }t=0\\
-\frac{2\hat{A}_{q}\left(t\right)}{\left(q-1\right)\left(q-3\right)} & \textrm{if }t\neq0
\end{cases} & \textrm{else}
\end{cases}\label{eq:poles}
\end{equation}
where:
\begin{equation}
\hat{A}_{q}\left(t\right)=\prod_{n=0}^{-v_{2}\left(t\right)-1}\frac{1+qe^{-2\pi i2^{n}t}}{4}
\end{equation}
It is very suspicious that, for values of $q$ not in the breakdown
variety of $\chi_{q}$, the $q$-poles of $\hat{\chi}_{q}\left(t\right)$
just so happen to be $q=1$ and $q=3$. When $q=1$, the map $T_{1}$
is easily shown to eventually iterate all positive integers to $1$.
When $q\geq5$, it is conjectured that $T_{q}$ has infinitely many
divergent points in $T_{q}$. In fact, it is conjectured that these
divergent points form a set of density $1$ in $\mathbb{N}_{0}$.
Meanwhile, we have a pole at $q=3$, which is precisely the value
of $q$ for which it is conjectured that $T_{q}$ eventually iterates
every positive integer to $1$.

So: what is going on here? More generally, considering $R_{1}=K\left[a_{0},\ldots,a_{p-1},b_{0},\ldots,b_{p-1}\right]$
and a degree $1$ F-series $X$, whose parameters are the $a_{j}$s
and $b_{j}$s, let $\hat{X}:\hat{\mathbb{Z}}_{p}\rightarrow\textrm{Frac}\left(R_{1}\right)\left(\zeta_{p^{\infty}}\right)$
be the formula for a Fourier transform of $X$ in the situation where
$\alpha_{X}\left(0\right)\neq1$:
\begin{equation}
\hat{X}\left(t\right)=\begin{cases}
\frac{\beta_{X}\left(0\right)}{1-\alpha_{X}\left(0\right)} & \textrm{if }t=0\\
\frac{\beta_{X}\left(0\right)\hat{A}_{X}\left(t\right)}{1-\alpha_{X}\left(0\right)}+\gamma_{X}\left(\frac{t\left|t\right|_{p}}{p}\right)\hat{A}_{X}\left(t\right) & \textrm{if }t\neq0
\end{cases},\textrm{ }\forall t\in\hat{\mathbb{Z}}_{p}\label{eq:X hat}
\end{equation}
(As an aside, note that unless we impose relations on $R_{1}$, the
equality $\alpha_{X}\left(0\right)=1$ will never occur.) Then, letting
$\mathbf{p}=\left(a_{0},\ldots,a_{p-1},b_{0},\ldots,b_{p-1}\right)$
vary in $K^{2p}$, let $e_{\mathbf{p}}$ denote the ring homomorphism
out of $R_{1}$ which evaluates the $a_{j}$s and $b_{j}$s at the
values specified by $\mathbf{p}$, and extend $e_{\mathbf{p}}$ to
a field homomorphism out of $\textrm{Frac}\left(R_{1}\right)\left(\zeta_{p^{\infty}}\right)$.
For a given unique solution ideal $I\subseteq R_{1}$, as $\mathbf{p}$
varies among those points of $K^{2p}$ satisfying all the relations
encoded by $I$, there will exist values of $\mathbf{p}$ for which
the image of $\hat{X}\left(t\right)$ under $e_{\mathbf{p}}$ is singular
for some $t$. Let's call these values of $\mathbf{p}$ the \textbf{poles
of $\hat{X}$ relative to $I$}.

Finally, fixing $I$, given any $\mathbf{p}$ which satisfies the
relations in $I$, let $H_{0},\ldots,H_{p-1}$ denote the affine linear
maps $K\rightarrow K$ given by:

\begin{equation}
H_{j}\left(x\right)\overset{\textrm{def}}{=}a_{j}x+b_{j},\textrm{ }\forall j\in\left\{ 0,\ldots,p-1\right\} ,\textrm{ }\forall x\in K
\end{equation}
where the $a_{j}$s and $b_{j}$s are the values specified by $\mathbf{p}$.
\begin{question}
\label{que:pole dynamics}Letting everything be as given in the above
discussion, suppose we can realize the $H_{j}$s as the branches of
a Hydra map $H:\mathcal{O}_{K}\rightarrow\mathcal{O}_{K}$. How, then,
do the dynamics of $H$ relate to the values of $\mathbf{p}$ which
are poles of $\hat{X}$?
\end{question}
Things get even more interesting when we consider Fourier transforms
of $\chi_{q}^{n}$ for $n\geq2$. The breakdown variety of $\chi_{q}^{n}$
is:
\begin{equation}
q^{n}+1=2^{n+1}
\end{equation}
This can be found by using the formal solution method, where we have:
\begin{align*}
\hat{\chi}_{q}^{*n}\left(t\right) & =\int_{\mathbb{Z}_{2}}\chi_{q}^{n}\left(\mathfrak{z}\right)e^{-2\pi i\left\{ t\mathfrak{z}\right\} _{2}}d\mathfrak{z}\\
 & =\frac{1}{2}\int_{\mathbb{Z}_{2}}\left(\frac{1}{2}\chi_{q}\left(\mathfrak{z}\right)\right)^{n}e^{-2\pi i\left\{ t2\mathfrak{z}\right\} _{2}}d\mathfrak{z}+\frac{1}{2}\int_{\mathbb{Z}_{2}}\left(\frac{q\chi_{q}\left(\mathfrak{z}\right)+1}{2}\right)^{n}e^{-2\pi i\left\{ t\left(2\mathfrak{z}+1\right)\right\} _{2}}d\mathfrak{z}\\
 & =\frac{1}{2^{n+1}}\hat{\chi}_{q}^{*n}\left(2t\right)+\frac{1}{2}\sum_{k=0}^{n}\binom{n}{k}\frac{q^{k}e^{-2\pi it}}{2^{n}}\int_{\mathbb{Z}_{2}}\chi_{q}^{k}\left(\mathfrak{z}\right)e^{-2\pi i\left\{ 2t\mathfrak{z}\right\} _{2}}d\mathfrak{z}\\
 & =\frac{1}{2^{n+1}}\hat{\chi}_{q}^{*n}\left(2t\right)+\frac{e^{-2\pi it}}{2^{n+1}}\sum_{k=0}^{n}\binom{n}{k}q^{k}\hat{\chi}_{q}^{*k}\left(2t\right)
\end{align*}
and hence:
\begin{equation}
\hat{\chi}_{q}^{*n}\left(t\right)=\underbrace{\frac{1+q^{n}e^{-2\pi it}}{2^{n+1}}}_{\alpha_{n}\left(t\right)}\hat{\chi}_{q}^{*n}\left(2t\right)+\frac{e^{-2\pi it}}{2^{n+1}}\sum_{k=0}^{n-1}\binom{n}{k}q^{k}\hat{\chi}_{q}^{*k}\left(2t\right)\label{eq:formal equation for Chi_q hat star n}
\end{equation}
and we get the breakdown variety of $\chi_{q}^{n}$ by setting $t=0$.
By \textbf{Theorem \ref{thm:formal solutions}}, if we have already
found Fourier transforms for $\hat{\chi}_{q},\ldots,\hat{\chi}_{q}^{*n-1}$,
solving this recursively for $\hat{\chi}_{q}^{*n}$ will give us a
Fourier transform of $\chi_{q}^{n}$ provided we stay away from any
ideals of $R=\mathbb{Z}\left[\frac{1}{2},q\right]$ for which $q^{n}+1=2^{n+1}$.

While every Hydra map induces an F-series by way of its numen, I have
yet to determine if the reverse also holds: that is, if every F-series
can be viewed as the numen of a Hydra map acting on some discrete
ring. That being said, in a case like this, where our F-series began
as the numen of a Hydra map, we can, in fact, go back and forth between
Hydras and numens.
\begin{example}
\label{exa:The-breakdown-variety}The breakdown variety of $\chi_{q}^{2}$
is $q^{2}=7$. Intriguingly, we realize $T_{\sqrt{7}}$ as a Hydra-like
map on $\mathbb{Z}\left[\sqrt{7}\right]$ (which is the ring of integers
of $\mathbb{Q}\left(\sqrt{7}\right)$), with the obvious choice of
branches:
\begin{align}
x & \mapsto\frac{x}{2}\\
x & \mapsto\frac{x\sqrt{7}+1}{2}
\end{align}
The problem is, since: 
\begin{equation}
2=9-7=\left(3-\sqrt{7}\right)\left(3+\sqrt{7}\right)
\end{equation}
$2$ is no longer prime in $\mathbb{Z}\left[\sqrt{7}\right]$. As
such, it no longer makes sense to ask ``what is the value of $3+\sqrt{7}$
mod $2$'', and that is important, because in order for $T_{\sqrt{7}}$
to be a Hydra map on $\mathbb{Z}\left[\sqrt{7}\right]$, we need to
be able to determine which branch of the map we apply based on the
value of the input modulo some prime ideal of $\mathbb{Z}\left[\sqrt{7}\right]$.

A way to salvage this is to consider how the branches of $T_{\sqrt{7}}$
transform elements of $\mathbb{Z}\left[\sqrt{7}\right]$. In order
to get the behavior we want, given $x\in\mathbb{Z}\left[\sqrt{7}\right]$,
we need there to be a unique branch of $T_{\sqrt{7}}$ which sends
$x$ to another element of $\mathbb{Z}\left[\sqrt{7}\right]$. So,
letting $x=a+b\sqrt{7}$ for $a,b\in\mathbb{Z}$, for the even branch,
set $y=x/2$. Then, $y$ satisfies:
\begin{equation}
y^{2}-ay+\frac{a^{2}-7b^{2}}{4}=0
\end{equation}
Here, $y\in\mathbb{Z}\left[\sqrt{7}\right]$ if and only if this polynomial
has integer coefficients. This occurs if and only if:
\begin{equation}
a^{2}\overset{4}{\equiv}3b^{2}
\end{equation}
which occurs if and only if $a\overset{2}{\equiv}b\overset{2}{\equiv}0$.

On the other hand, for the odd branch, setting $y=\left(x\sqrt{7}+1\right)/2$,
we get:
\begin{align*}
y & =\frac{\left(a+b\sqrt{7}\right)\sqrt{7}+1}{2}\\
 & \Updownarrow\\
0 & =y^{2}-\left(7b+1\right)y+\frac{49b^{2}+14b-7a^{2}+1}{4}
\end{align*}
So, for $y$ to be in $\mathbb{Z}\left[\sqrt{7}\right]$, we need:
\begin{equation}
49b^{2}+14b-7a^{2}+1\overset{4}{\equiv}0
\end{equation}
Since the only perfect squares mod $4$ are $0$ and $1$, this yields:
\begin{equation}
\begin{cases}
2b\overset{4}{\equiv}3 & \textrm{if }a^{2}\overset{4}{\equiv}0,\textrm{ }b^{2}\overset{4}{\equiv}0\\
b\overset{2}{\equiv}1 & \textrm{if }a^{2}\overset{4}{\equiv}1,\textrm{ }b^{2}\overset{4}{\equiv}0\\
b\overset{2}{\equiv}1 & \textrm{if }a^{2}\overset{4}{\equiv}0,\textrm{ }b^{2}\overset{4}{\equiv}1\\
2b\overset{4}{\equiv}1 & \textrm{if }a^{2}\overset{4}{\equiv}1,\textrm{ }b^{2}\overset{4}{\equiv}1
\end{cases}
\end{equation}
The only case where this is possible is when $b\overset{2}{\equiv}1$
and $a^{2}\overset{4}{\equiv}0$, which is the same as $a\overset{2}{\equiv}0$
and $b\overset{2}{\equiv}1$.

In conclusion, given any $a,b\in\mathbb{Z}$, if $a$ is even, we
can apply $T_{\sqrt{7}}$ to $x=a+b\sqrt{7}$ by choosing the branch
corresponding to the parity of $b$, and the result will again be
an element of $\mathbb{Z}\left[\sqrt{7}\right]$. That is:
\begin{equation}
T_{\sqrt{7}}\left(a+b\sqrt{7}\right)=\begin{cases}
\frac{a+b\sqrt{7}}{2} & \textrm{if }a\overset{2}{\equiv}0\textrm{ \& }b\overset{2}{\equiv}0\\
\frac{\left(a+b\sqrt{7}\right)\sqrt{7}+1}{2} & \textrm{if }a\overset{2}{\equiv}0\textrm{ \& }b\overset{2}{\equiv}1
\end{cases}
\end{equation}
If $a$ is odd, however, then neither branch will output an element
of $\mathbb{Z}\left[\sqrt{7}\right]$. Thus, the forward orbit of
a given $x\in\mathbb{Z}\left[\sqrt{7}\right]$ will eventually stop
if it ever sends the rational component of $x$ to an odd integer.

Running $T_{\sqrt{7}}$ on a computer, an associate of mine found
that the only cycles with $\max\left\{ \left|a\right|,\left|b\right|\right\} <100$
are $\left\{ 0\right\} $ and: 
\begin{equation}
2+\sqrt{7}\rightarrow4+\sqrt{7}\rightarrow4+2\sqrt{7}
\end{equation}
which seems suspiciously Collatz-like.
\end{example}
More generally, letting $q_{n}$ denote the unique positive real root
of $q^{n}+1=2^{n+1}$, we have that the extent to which we can realize
$T_{q_{n}}$ as a map out of the ring of integers of the number field
$\mathbb{Q}\left(q_{n}\right)$ will require understanding both the
factorization of $\left\langle 2\right\rangle $ in $\mathcal{O}_{\mathbb{Q}\left(q_{n}\right)}$
and which points, if any, of in $\mathcal{O}_{\mathbb{Q}\left(q_{n}\right)}$
that the branches of $T_{q_{n}}$ send to non-integer elements of
$\mathbb{Q}\left(q_{n}\right)$. In particular, note that $T_{q_{n}}$
will map $\mathcal{O}_{\mathbb{Q}\left(q_{n}\right)}\rightarrow\mathcal{O}_{\mathbb{Q}\left(q_{n}\right)}$
whenever $\left\langle 2\right\rangle $ is inert in $\mathcal{O}_{\mathbb{Q}\left(q_{n}\right)}$.

This then leads us to a natural generalization of \textbf{Question
}\ref{que:pole dynamics}:
\begin{question}
Let everything be as given in the discussion leading up to \textbf{Question
\ref{que:pole dynamics}}, and let $n\geq2$. This time, let $\mathbf{p}$
be a point in $\overline{K}^{2p}$ (where $\overline{K}$ is the algebraic
closure of $K$) which is either in the breakdown variety of $X^{n}$
relative to $I$ or which is a pole of $\hat{X}^{*n}$ relative to
$I$. Let $K\left(\mathbf{p}\right)$ denote the algebraic extension
of $K$ given by adjoining all the entries of $\mathbf{p}$, and let
the $H_{j}$s be the affine linear maps on $K\left(\mathbf{p}\right)$
given by setting the $a_{j}$s and $b_{j}$s to be the entries of
$\mathbf{p}$. As \textbf{Example \ref{exa:The-breakdown-variety}
}showed, we might not be able to stitch the $H_{j}$s together to
form a hydra map $H$ $\mathcal{O}_{K\left(\mathbf{p}\right)}$ (that
is, there might be elements of $\mathcal{O}_{K\left(\mathbf{p}\right)}$
that none of the $H_{j}$s map to elements of $\mathcal{O}_{K\left(\mathbf{p}\right)}$)
but, if we can, as in \textbf{Example \ref{exa:The-breakdown-variety}},
define $H$ on \emph{some }of $\mathcal{O}_{K\left(\mathbf{p}\right)}$,
how might the dynamics of $H$ relate to $\mathbf{p}$? Likewise,
what happens if we do this for values of $\mathbf{p}$ which do not
lie in $X^{n}$'s breakdown variety relative to $I$, and/or which
are not poles of $\hat{X}^{*n}$ relative to $I$?
\end{question}
Also:
\begin{question}
What is the significance, if any, to poles of $\hat{X}^{*n}$ relative
to $I$ which are not also points of $X^{n}$'s breakdown variety
relative to $I$?
\end{question}
And, of course:
\begin{question}
Is there any formalism to help state these questions in a more straightforward
manner?
\end{question}
Another advantage of treating the parameters of an F-series as formal
indeterminates is that it allows us to introduce a tool totally absent
from classical $\left(p,q\right)$-adic analysis: differentiation
and, more generally, derivations. The formalism presented in this
paper allows for F-series and their Fourier transforms to be differentiated
with respect to their parameters without having to worry about issues
of compatibility. (Indeed, the root test conditions given in \textbf{Corollary
\ref{cor:main result}} ensure that the F-series in question will
remain $\mathcal{F}$-convergent after differentiating with respect
to any of their parameters.
\begin{example}
Consider the $2$-adic F-series $X_{b,a}:\mathbb{Z}_{2}\rightarrow\mathbb{Z}\left[\left[a,b\right]\right]$
given by:
\begin{equation}
X_{b,a}\left(\mathfrak{z}\right)=\sum_{n=0}^{\infty}a^{n}b^{\#_{2:1}\left(\left[\mathfrak{z}\right]_{2^{n}}\right)}
\end{equation}
where $a$ and $b$ are indeterminates. This satisfies the functional
equations:
\begin{align}
X_{b,a}\left(2\mathfrak{z}\right) & =1+aX_{b,a}\left(\mathfrak{z}\right)\\
X_{b,a}\left(2\mathfrak{z}+1\right) & =1+abX_{b,a}\left(\mathfrak{z}\right)
\end{align}
and has a Fourier transform:

\begin{equation}
\hat{X}_{b,a}\left(t\right)=\frac{\hat{A}_{b,a}\left(t\right)}{1-\frac{a\left(b+1\right)}{2}}=\frac{1}{1-\frac{a\left(b+1\right)}{2}}\prod_{n=0}^{-v_{2}\left(t\right)-1}\left(a+abe^{-2\pi i2^{n}t}\right),\textrm{ }\forall t\in\hat{\mathbb{Z}}_{2}\label{eq:X_b,a hat co-breakdown}
\end{equation}
The breakdown variety is:
\begin{equation}
a+ab=2
\end{equation}
Given any unique solution ideal $I\subseteq\mathbb{Z}\left[a,b\right]$
(here, $\left\langle 1-a\right\rangle \nsubseteq I$) we have that
the image of (\ref{eq:X_b,a hat co-breakdown}) under the projections
mod $I$ is again a Fourier transform of $X_{b,a}$ whenever $\left\langle a+ab-2\right\rangle \nsubseteq I$.
When $\left\langle a+ab-2\right\rangle \subseteq I$, we instead must
use the image of the formula:
\begin{equation}
\hat{X}_{b,a}\left(t\right)=\min\left\{ v_{2}\left(t\right),0\right\} \hat{A}_{b,a}\left(t\right),\textrm{ }\forall t\in\hat{\mathbb{Z}}_{2}\label{eq:X_b,a hat breakdown}
\end{equation}
under the projection mod $I$, where, recall, $v_{2}\left(0\right)=+\infty$
and $v_{2}\left(t\right)<0$ for all $t\in\hat{\mathbb{Z}}_{2}\backslash\left\{ 0\right\} $.

Regardless of whether or not $I$ contains the ideal $\left\langle a+ab-2\right\rangle $
generated by the breakdown variety, we can formally differentiate
both $X_{b,a}$ and $\hat{X}_{b,a}$ with respect to $a$ and/or $b$,
and the results will be well-behaved; ex:
\begin{equation}
\sum_{t\in\hat{\mathbb{Z}}_{2}}\frac{\partial\hat{X}_{b,a}}{\partial b}\left(t\right)e^{2\pi i\left\{ t\mathfrak{z}\right\} _{2}}\overset{\mathcal{F}}{=}\frac{\partial}{\partial b}\left(\sum_{t\in\hat{\mathbb{Z}}_{2}}\hat{X}_{b,a}\left(t\right)e^{2\pi i\left\{ t\mathfrak{z}\right\} _{2}}\right)=\frac{\partial X_{b,a}}{\partial b}\left(\mathfrak{z}\right)
\end{equation}
and likewise for $\partial/\partial a$.

For example, we can deduce functional equations for the derivatives
of $X_{b,a}$ by formally differentiating the functional equations.
This can be verified by performing the same differentiation to the
F-series of $X_{a,b}$ and using the identity:
\begin{equation}
\#_{2:1}\left(\left[2\mathfrak{z}+j\right]_{2^{n}}\right)=\#_{2:1}\left(j\right)+\#_{2:1}\left(\left[\mathfrak{z}\right]_{2^{n-1}}\right),\textrm{ }\forall n\geq1,\textrm{ }\forall j\in\left\{ 0,1\right\} ,\textrm{ }\forall\mathfrak{z}\in\mathbb{Z}_{2}
\end{equation}
Doing so, we get:
\begin{align}
\frac{\partial X_{b,a}}{\partial a}\left(2\mathfrak{z}\right) & =X_{b,a}\left(\mathfrak{z}\right)+a\frac{\partial X_{b,a}}{\partial a}\left(\mathfrak{z}\right)\\
\frac{\partial X_{b,a}}{\partial a}\left(2\mathfrak{z}+1\right) & =bX_{b,a}\left(\mathfrak{z}\right)+ab\frac{\partial X_{b,a}}{\partial a}\left(\mathfrak{z}\right)
\end{align}
and:
\begin{align}
\frac{\partial X_{b,a}}{\partial b}\left(2\mathfrak{z}\right) & =a\frac{\partial X_{b,a}}{\partial b}\left(\mathfrak{z}\right)\\
\frac{\partial X_{b,a}}{\partial b}\left(2\mathfrak{z}+1\right) & =aX_{b,a}\left(\mathfrak{z}\right)+ab\frac{\partial X_{b,a}}{\partial b}\left(\mathfrak{z}\right)
\end{align}
Because of this, we see that Fourier transforms of derivatives of
$X_{b,a}$ can then be expressed recursively in terms of Fourier transforms
of $X_{a,b}$ itself.
\end{example}
\begin{rem}
I intend to investigate the behavior of F-series under differential
operators like these in a later paper.
\end{rem}
The above discussions show how we can use my algebraic formalism to
enhance our study of F-series, and by extension, of the dynamics of
hydra type maps. However, both with a mind toward discovering applications
and toward better understanding the relations between the measures
induced by F-series and the dynamics of the associated hydra maps
(if any), it's also worth exploring how the structures presented in
this paper might be employed to study algebraic varieties. As algebraic
geometry is, for me, as near to the opposite of a forté as one can
get, there's not that much that I can say, though what I do have to
say makes it clear that there's definitely much more going on here
that merits further exploration in its own right.

To that end, I figure I might as well end this paper with what I hope
will be a way forward for future work. Simply put: I believe we should
try to develop a way to treat $\chi_{q}$ and other numens as geometric
objects\textemdash specifically, as \emph{curves}\textemdash so as
to enable us to define and compute algebraic invariants of these curves.
If it could be shown that the dynamical properties of hydra maps end
up being encoded in the geometry of the curves defined by their numens,
proving results about the hydras' dynamics might reduce to a matter
of comparing their numens' invariants. I first discussed this possibility
in \cite{first blog paper}, though I have not been able to progress
on it very much, due to the pertinent material being situated deep
inside the areas of algebraic geometry that, at least at the time
of writing, I haplessly struggle to confront.

Perhaps it is just my own naïveté speaking, but I can't help but see
parallels that suggest the research direction I am proposing here
is worth exploring. Indeed, at a purely formal level, there is already
a well-established precedent for treating objects like $\chi_{q}$
as curves: the so-called \textbf{de Rham curves} \cite{de rham}.
\begin{defn}
Let $\left(M,d\right)$ be a complete metric space, let $p$ be an
integer $\geq2$, and let $H_{0},\ldots,H_{p-1}:M\rightarrow M$ be
contraction mappings; that is, for each $j\in\left\{ 0,\ldots,p-1\right\} $,
there is a real constant $0<\lambda_{j}<1$ so that:
\begin{equation}
d\left(H_{j}\left(x\right),H_{j}\left(y\right)\right)\leq\lambda_{j}d\left(x,y\right),\textrm{ }\forall x,y\in M
\end{equation}
It can be shown that for any infinite sequence of integers $\left\{ a_{n}\right\} _{n\geq0}$
in $\left\{ 0,1,\ldots,p-1\right\} $, the infinite composition:
\begin{equation}
H_{a_{0}}\circ H_{a_{1}}\circ\cdots
\end{equation}
will be a constant function, mapping any point in $M$ to an output
that depends solely on the sequence $\left\{ a_{n}\right\} _{n\geq0}$;
this constant will be the unique fixed point of $H_{a_{0}}\circ H_{a_{1}}\circ\cdots$.
Since:
\[
\sum_{n=0}^{\infty}\frac{1}{p^{n}}\overset{\mathbb{R}}{=}\frac{p}{p-1}
\]
we have that every real number $x\in\left[0,1\right]$ can be represented
as:
\begin{equation}
\frac{p-1}{p}\sum_{n=0}^{\infty}\frac{a_{n}\left(x\right)}{p^{n}}\label{eq:digit representation}
\end{equation}
for some sequence of integers $a_{n}\left(x\right)\in\left\{ 0,1,\ldots,p-1\right\} $.
Then, the \textbf{de Rham curve} \textbf{generated by the $H_{j}$}s
is the map $c:\left[0,1\right]\rightarrow M$ which, to each $x\in\left[0,1\right]$,
outputs the fixed point of the composition:
\begin{equation}
H_{a_{0}}\circ H_{a_{1}}\circ\cdots
\end{equation}
\end{defn}
\begin{rem}
Due to their construction, de Rham curves tend to exhibit fractal,
self-similar behavior. In general, these curves will not be continuous,
however their continuity can be guaranteed in the following manner.
Since each $H_{j}$ is a contraction map, by the \textbf{Banach Fixed-Point
Theorem}, each $H_{j}$ has a unique fixed point $x_{j}\in M$. Because
of this, the de Rham curve generated by the $H_{j}$s will be continuous
whenever the $H_{j}$s satisfy the property:
\begin{equation}
H_{j}\left(x_{n}\right)=H_{j+1}\left(x_{0}\right),\textrm{ }\textrm{for }j\in\left\{ 0,\ldots,p-2\right\} \label{eq:de Rham continuity condition}
\end{equation}
This condition is necessary in order to deal with the fact that not
every real number in $\left[0,1\right]$ has a \emph{unique} representation
in of the form (\ref{eq:digit representation}), the classic example
being the equality of $1$ and $0.999\ldots$ in base $10$.
\end{rem}
Oddly enough, de Rham's the construction is formally identical to
how I originally obtained $\chi_{q}$: by indexing arbitrary compositions
of sequences of the maps $x\mapsto x/2$ and $x\mapsto\left(qx+1\right)/2$
on $\mathbb{Q}$ using sequences of $0$s and $1$s, and then identifying
these sequences with the digits of $2$-adic integers.
\begin{example}
Let $K$ be a global field, and let $X$ be a $p$-adic foliated F-series
satisfying:
\begin{equation}
X\left(p\mathfrak{z}+j\right)=a_{j}X\left(\mathfrak{z}\right)+b_{j},\textrm{ }\forall\mathfrak{z}\in\mathbb{Z}_{p},\textrm{ }\forall j\in\left\{ 0,\ldots,p-1\right\} 
\end{equation}
where the $a_{j}$s and $b_{j}$s are elements of $K$, with $a_{j}\notin\mathcal{O}_{K}^{\times}\cup\left\{ 0\right\} $
for all $j$. The condition on the $a_{j}$s guarantees that for each
$j$, there exists a place $\ell_{j}$ of $K$ so that $\left|a_{j}\right|_{\ell_{j}}<1$.
This makes the map $H_{j}:K_{\ell_{j}}\rightarrow K_{\ell_{j}}$ defined
by:
\begin{equation}
H_{j}\left(z\right)\overset{\textrm{def}}{=}a_{j}z+b_{j}
\end{equation}
into a contraction on $K_{\ell_{j}}$, with:
\begin{equation}
\left|H_{j}\left(\mathfrak{x}\right)-H_{j}\left(\mathfrak{y}\right)\right|_{\ell_{j}}\leq\left|a_{j}\right|_{\ell_{j}}\left|\mathfrak{x}-\mathfrak{y}\right|_{\ell_{j}},\textrm{ }\forall\mathfrak{x},\mathfrak{y}\in K_{\ell_{j}}
\end{equation}
It can be shown that $X\left(\mathfrak{z}\right)$ is the function
which, to every: 
\begin{equation}
\mathfrak{z}=\sum_{n=0}^{\infty}d_{n}\left(\mathfrak{z}\right)p^{n}\in\mathbb{Z}_{p}
\end{equation}
(where $d_{n}\left(\mathfrak{z}\right)$ is the $n$th $p$-adic digit
of $\mathfrak{z}$), outputs the fixed point of the composition sequence:
\begin{equation}
H_{d_{0}\left(\mathfrak{z}\right)}\circ H_{d_{1}\left(\mathfrak{z}\right)}\circ\cdots
\end{equation}
If $\mathfrak{z}\in\mathbb{Z}_{p}^{\prime}$, then there is a $j\in\left\{ 1,\ldots,p-1\right\} $
so that $d_{n}\left(\mathfrak{z}\right)=j$ for infinitely many $n$.
In this case, the composition sequence will converge to a limit in
$K_{\ell_{j}}$. On the other hand, if $\mathfrak{z}\in\mathbb{N}_{0}$,
then we have that $d_{n}\left(\mathfrak{z}\right)=0$ for all sufficiently
large $n$, and thus, that the composition sequence will converge
to a limit in $K_{\ell_{0}}$. We can then use the absolute values
to construct a digital frame compatible with $X$ in the manner discussed
in \textbf{Definition \ref{def:digit construction-1}} back on page
\pageref{def:digit construction-1}.
\end{example}
There are only three main distinctions between my construction and
de Rham's:
\begin{itemize}
\item Instead of using $\left[0,1\right]$ as the domain of my curve (the
F-series, $X$), I use $\mathbb{Z}_{p}$.
\item Instead of working with contraction mappings on a complete metric
space, I work with maps that are contractions in a \emph{completion
of }an underlying global field.
\item The continuity condition (\ref{eq:de Rham continuity condition})
is not satisfied. (Indeed, for $\chi_{q}$, the branches' fixed points
are $x=0$ and $x=-1/\left(q-2\right)$. $\left(qx+1\right)/2$ sends
$0$ to $1/2$, while $x/2$ sends $-1/\left(q-2\right)$ to $-1/\left(2\left(q-2\right)\right)$.)
\end{itemize}
While these departures from de Rham's hypotheses might seem discouraging
at first, if we look closer, interactions with this paper's work emerge.
In both my constructions and de Rham's, the underlying idea is to
build a map from the set of all sequences in the symbols $\left\{ 0,1,\ldots,p-1\right\} $
to the set of composition sequences of the maps $H_{0},\ldots,H_{p-1}$.
We just choose two different realizations of this sequence space;
indeed, they differ only by a reciprocal.

In the $2$-adic case, for instance, consider the map $\eta_{2}:\mathbb{Z}_{2}\rightarrow\left[0,1\right]$
given by:
\begin{equation}
\eta_{2}\left(\sum_{n=0}^{\infty}d_{n}\left(\mathfrak{z}\right)2^{n}\right)\overset{\textrm{def}}{=}\sum_{n=0}^{\infty}\frac{d_{n}\left(\mathfrak{z}\right)}{2^{n+1}}
\end{equation}
This a standard embedding of $\mathbb{Z}_{2}$ in the reals; see \cite{Robert's Book}.
Notably, this map is surjective, but not injective. Nevertheless,
if we choose an inverse of $\eta_{2}$, we can realize $\chi_{q}$
as a function $C_{q}$ parameterized over $\left[0,1\right]$ instead
of $\mathbb{Z}_{2}$:
\begin{align}
C_{q}\left(\frac{t}{2}\right) & =\frac{1}{2}C_{q}\left(t\right)\\
C_{q}\left(\frac{t+1}{2}\right) & =\frac{qC_{q}\left(t\right)+1}{2}
\end{align}
for all $t\in\left[0,1\right]$. Formally: 
\begin{equation}
C_{q}\left(t\right)\overset{\textrm{def}}{=}\left(\chi_{q}\circ\eta_{2}^{-1}\right)\left(t\right),\textrm{ }\forall t\in\left[0,1\right]\label{eq:Definition of C_q}
\end{equation}
We even recover an archimedean analogue of rising-continuity: $C_{q}$
is easily seen to be \textbf{left-continuous}:
\begin{equation}
\lim_{n\rightarrow\infty}C_{q}\left(t_{n}\right)\overset{\mathbb{Z}_{q}}{=}C_{q}\left(t\right)
\end{equation}
for all $t\in\left[0,1\right]$ and any strictly increasing sequence
$\left\{ t_{n}\right\} _{n\geq1}$ in $\left[0,1\right]$ converging
to $t$ \cite{first blog paper}.

Note that this is an archimedean incarnation of the functoriality
we saw in our (admittedly brief) discussion of reffinite spaces and
space-change maps, back on page \pageref{def:space-change map} at
the end of \textbf{Section \ref{subsec:Frame-Theory-=000026}}. There,
recall, we wrote $\psi_{q,p}:\mathbb{Z}_{p}\rightarrow\mathbb{Z}_{q}$
to denote the map:
\begin{equation}
\psi_{q,p}\left(\sum_{n=0}^{\infty}d_{n}p^{n}\right)=\sum_{n=0}^{\infty}d_{n}q^{n}
\end{equation}
Here, $p$ and $q$ are numbers. However, if we take a page from F-series
and view $p$ and $q$ merely as formal indeterminates, all $\psi_{q,p}$
is doing is changing series in one formal indeterminate to another.
$\eta_{2}$ does this very same thing, changing the symbol $2^{n}$
into the symbol $1/2^{n}$. In this way, it can be argued that de
Rham's choice of $\left[0,1\right]$ as the space over which to parameterize
his curves was somewhat arbitrary. Just as polynomials are a way of
stringing together a finite sequence as a single algebraic object,
so too are power series a means of stringing together an infinite
sequence as a single algebraic object. That we can go from $\chi_{q}$
to $C_{q}$ by precomposing with a chosen inverse of $\eta_{2}$ is
in the spirit of the functoriality mentioned above. The only difference
is that $\left[0,1\right]$ is not a profinite space. I mean, at a
formal level, it is\textemdash all of its elements are representable
as power series in a single symbol with coefficients drawn from a
finite set\textemdash even though its topology is completely different
from that of a profinite space. This suggests that there might be
a way of giving a unified treatment of both the archimedean and non-archimedean
spaces over which we can pull back a given F-series.

From what vanishingly little of Scholze and Clausen's theory of condensed
mathematics I can currently comprehend, I am aware that one of that
theory's goals is to give a non-archimedean foundation for archimedean
spaces like the real or complex numbers \cite{condensed}. Though
the details of their theory escape me, even \emph{I} can make out
the striking thematic parallels between what they are doing and what
I would like to do. This paper shows that there is Fourier-analytic
structure which is compatible with the interpretation of F-series
and the associated distributions/measures as functors out of the category
of reffinite abelian groups with space-change maps as the morphisms.
At the risk of being presumptuous, I believe that the work I have
presented here and in my other publications offers an entirely new
paradigm for dealing with the dynamics of Collatz-type maps, one whose
pursuit might contribute to our understanding of non-archimedean geometry
as much as it might benefit from it. I just hope that others will
be willing to join me on this expedition, and to take this work where
I cannot.

\section*{Statements and Declarations}

The author states that there are no conflicts of interest, and no
data to be made available. This research did not receive any specific
grant from funding agencies in the public, commercial, or not-for-profit
sectors.

\section*{Acknowledgements}

This paper draws from and builds on my PhD dissertation \cite{My Dissertation}done
at the University of Southern California under the obliging supervision
of Professors Sheldon Kamienny and Nicolai Haydn. Thanks must also
be given to Jeffery Lagarias, Steven J. Miller, Alex Kontorovich,
Andrei Khrennikov, K.R. Matthews, Susan Montgomery, Amy Young and
all the helpful staff of the USC Mathematics Department, and all the
other kindly strangers I've met along the way.

\vphantom{}

\end{document}